\documentclass[11pt,a4paper,makeidx]{amsart}
\usepackage[latin1]{inputenc}        
\usepackage[breaklinks]{hyperref}
\usepackage{appendix}
\usepackage{graphicx}
\usepackage{float}
\usepackage{mathrsfs}
\usepackage{html}
\usepackage{htmllist}
\let\mathcal\mathscr
\usepackage{color}              
\usepackage{amssymb}
\usepackage{amsmath}
\usepackage{amsfonts}
\usepackage{latexsym}
\usepackage[T1]{fontenc}
\usepackage[latin1]{inputenc}
\usepackage[all,ps]{xy}
\RequirePackage{amsthm}
\RequirePackage{amssymb}
\usepackage{latexsym}
\usepackage{amscd}
\usepackage[latin1]{inputenc}
\usepackage[T1]{fontenc}
\usepackage{color}
\RequirePackage{amsthm}
\RequirePackage{amssymb}
\usepackage[all,ps]{xy}
\usepackage{latexsym}

\def\ra{\rightarrow}

\def\la{\leftarrow}
\def\rras{\rightrightarrows}

\def\mpo{\mapsto}

\def\bsh{\backslash}

\def\sbs{\subset}

\def\ts{\times}

\def\uts{\underline{\times}}

\def\a{\alpha}
\def\b{\beta}
\def\d{\delta}
\def\D{\Delta}
\def\e{\epsilon}
\def\g{\gamma}
\def\G{\Gamma}

\def\lb{\lambda}

\def\o{\omega}
\def\om{\omega}
\def\Om{\Omega}
\def\s{\sigma}
\def\Si{\Sigma}

\def\z{\zeta}
\def\r{\rho}

\def\cA{{\mathcal A}}
\def\cB{{\mathcal B}}
\def\cC{{\mathcal C}}
\def\cD{{\mathcal D}}
\def\cE{{\mathcal E}}
\def\cF{{\mathcal F}}
\def\cG{{\mathcal G}}
\def\cH{{\mathcal H}}

\def\cJ{{\mathcal J}}

\def\cM{{\mathcal M}}

\def\cO{{\mathcal O}}
\def\cP{{\mathcal P}}

\def\cS{{\mathcal S}}
\def\cT{{\mathcal T}}
\def\cU{{\mathcal U}}
\def\cV{{\mathcal V}}

\def\cX{{\mathcal X}}
\def\cY{{\mathcal Y}}
\def\cZ{{\mathcal Z}}

\def\tf{\tilde{F}}

\def\th{\tilde{H}}

\def\tt{\tilde{T}}

\def\ba{\mathbb A}
\def\bc{\mathbb C}
\def\be{\mathbb E}
\def\bF{\mathbb F}
\def\bg{\mathbb G}

\def\bn{\mathbb N}
\def\bp{\mathbb P}
\def\bq{\mathbb Q}
\def\br{\mathbb R}

\def\bz{\mathbb Z}
\def\hbz{{\hat{\mathbb Z}}}

\def\op{^{\mathrm{op}}}
\def\cop{^{\mathrm{co-op}}}

\def\do{D_{0}}

\def\pa{\partial}

\def\build#1_#2^#3{\mathrel{
\mathop{\kern 0pt#1}\limits_{#2}^{#3}}}

\def\mygg{\mathfrak{G}}

\def\Aut{\mathrm{ Aut}}
\def\Out{\mathrm{ Out}}
\def\obs{\mathrm{ obs}}

\def\fAut{\mathfrak{Aut}}
\def\fHom{\mathfrak{Hom}}

\def\hol{\mathrm{hol}}
\def\red{{\mathrm{red}}}
\def\Vect{{\mathrm{Vect}}}

\def\Ext{\mathrm{ Ext}}
\def\ext{\underline{\Ext}}
\def\hExt{\be\mathrm{xt}}
\def\hext{\underline{\hExt}}
\def\uom{\underline{\om}}

\def\Hom{\mathrm{Hom}}
\def\uHom{\underline{\mathrm{Hom}}}
\def\Home{\mathrm{Hom}_{\mathrm{Ens}}}
\def\id{\mathrm{id}}

\def\rB{\mathrm{B}}
\def\PSL{\mathrm{PSL}}
\def\uZ{\underline{Z}}
\def\uF{\underline{F}}
\def\uV{\underline{V}}

\def\ul{{\underline{\ell}}}
\def\uG{\underline{\Gamma}}
\def\uuG{\underline{\underline{\Gamma}}}
\def\rH{\mathrm{H}}
\def\htg{\mathrm{H}_{\mathrm{Giraud}}^2}
\def\rP{\mathrm{P}}
\def\rp{\mathrm{pt}}
\def\rS{\mathrm{S}}
\def\rK{\mathrm{K}}
\def\rE{\mathrm{E}}
\def\rF{\mathrm{F}}
\def\rD{\mathrm{D}}
\def\rT{\mathrm{T}}
\def\rO{\mathrm{O}}
\def\rI{\mathrm{I}}
\def\et{\acute{\mathrm{E}}\mathrm{t}}
\def\het{\hat{\mathrm{Et}}}
\def\Ens{{\mathrm{Ens}}}
\def\Top{{\mathrm{Top}}}
\def\Grpd{\mathrm{Grpd}}
\def\cts{{\mathrm{cts}}}
\def\fin{{\mathrm{fin}}}
\def\pro{{\mathrm{pro}}}

\def\gP{\mathfrak{P}}
\def\gp{\mathfrak{p}}
\def\gq{\mathfrak{q}}

\def\gQ{\mathfrak{Q}}
\def\gG{\mathfrak{G}}
\def\gF{\mathfrak{F}}
\def\gC{\mathfrak{C}}

\swapnumbers

\newtheorem{thm}{Theorem}[subsection]

\newtheorem{fact}[thm]{Fact}

\newtheorem{cor}[thm]{Corollary}
\newtheorem{prop}[thm]{Proposition}
\newtheorem{lem}[thm]{Lemma}
\newtheorem{flem}[thm]{Further Lemma}

\theoremstyle{definition}
\newtheorem{defn}[thm]{Definition}

\newtheorem{factdef}[thm]{Fact/Definition}

\newtheorem{rmkdef}[thm]{Remark/Definition}
\newtheorem{altdef}[thm]{Alternative Definition}
\newtheorem{newnot}[thm]{New Notation}
\newtheorem{summary}[thm]{Summary}

\newtheorem{flav}[thm]{Flavours}
\newtheorem{example}[thm]{Example}
\newtheorem{setup}[thm]{Set Up}
\newtheorem{setuprev}[thm]{Set Up/Revision}
\newtheorem{claim}[thm]{Claim}
\newtheorem{subclaim}[thm]{Sub-Claim}
\newtheorem{sublem}[thm]{Sub-Lemma}

\theoremstyle{remark}
\newtheorem{ex}[thm]{Example}

\newtheorem{notation}[thm]{Notation}

\newtheorem{rmk}[thm]{Remark}
\newtheorem{warning}[thm]{Warning}

\newtheorem{scholion}[thm]{Scholion}
\newtheorem{Cscholion}[thm]{Continuation of the Scholion}
\newtheorem{DefPi2}[thm]{Defining $\pi^\pro_2$}
\newtheorem{DefPos}[thm]{Definition 
and existence of p\underline{ointed} Postnikov sequences}
\newtheorem{UniPos}[thm]{Uniqueness  
of pointed and un-pointed Postnikov sequences}
\newtheorem{Pro2G}[thm]{The pro-2-Galois correspondence}
\newtheorem{ProF}[thm]{Unconditional pro-finite theory}

\title[Elementary topology of champs]{Elementary topology of champs}
\author[Michael McQuillan]{Michael McQuillan}
\date{\today}

\begin{document}
\begin{abstract}
Broadly speaking the present is a homotopy
complement to the book of Giraud, \cite{giraud}, 
albeit in a couple of different ways. In the
first place there is a representability theorem
for maps to a topological champ
(the translation stack will
be  eschewed
\ref{faq5.1})
and whence an extremely convenient global atlas,
{\it i.e.} the path space, which permits an
immediate importation of the familiar definitions
of homotopy groups and covering spaces as
encountered in elementary text books, \ref{faqR}.
In the second place,
it provides the adjoint to Giraud's co-homology,
{\it i.e.} the homotopy 2-group $\Pi_2$,
by way of the 2-Galois theory of
covering champs. 
In the sufficiently path connected case this is achieved by
much the same construction employed in constructing
1-covers, {\it i.e.} quotients of the path space by
a groupoid, \ref{faqII}. In the general case, 
so inter alia the pro-finite theory appropriate for
algebraic geometry, 
the development parallels the axiomatic
Galois theory of \cite[expose\'e V]{sga1},
\ref{faqIII}.
The resulting
explicit description of the homotopy 2-type
can be applied to prove theorems in algebraic
geometry: 
optimal generalisations to $\Pi_2$
(by a very different
method, which even gives improvements
to the original case) of the Lefschetz theorems (over a
locally Noetherian base) of \cite{sga2}, \ref{faq4.5}-\ref{faq4.7},
and a counterexample to the extension from 
co-homology to homotopy of the smooth base
change theorem, \ref{faq4.3}. These limited goals 
are achieved, 
albeit arguably at the price of obscuring the
higher categorical structure,
without leaving the 2-category of groupoids.
\end{abstract}
\maketitle
\numberwithin{equation}{section}
\newpage
\tableofcontents
\newpage
\section{F.A.Q.}\label{S:0}
\setcounter{subsection}{-1}
\subsection{What's this about, and why should I care ?}
\subsubsection{Is this one of those manuscripts that
I need to know 25 alternative definitions of
 \texorpdfstring{$n$}{n}-Category 
and 16 generalisations of the \'etale site in order to
get started ?}\label{faq1}

No, it's elementary. As such the higher category
pre-requisites are no more complicated than 
2-categories, 2-functors, their transformations,
and, occasionally, their modifications, all 
understood weakly rather than strictly. The slightly
more demanding pre-requisite is that one should be
familiar with how one uses 2-categories to construct
what might be termed ``exotic gluings'' of spaces,
be they topological, differential, symplectic, algebraic
or whatever, {\it i.e.} in an arbitrary site. This is
Grothendieck's theory of champs, \ref{faq5.1}, and, unfortunately,
the literature tends towards an obsession with 
algebraic geometry. Notable exceptions are Grothendieck
himself, ``il s'est av\'er\'e impossible de faire
de la descente dans la cat\'egorie des pr\'eschemas,
m\^{e}me dans des cas particuliers, sans avoir
d\'evelopp\'e au pr\'ealable avec assez de soin
le langage de la descente dans les cat\'egories 
g\'en\'erales'', \cite[Expos\'e VI]{sga1}, and
the book of Giraud, \cite{giraud}, ``Cohomologie
non ab\'elienne''. What's involved, however, is, 
in the first instance, just
extending the definition of sheaf from a functor
with values in sets together with some gluing conditions,
to that of a 2-sheaf, {\it i.e.} a (weak) 2-functor
with values in groupoids with some further gluing
conditions, \ref{SS:TwoSheaf}, to which some geometric
conditions such as the existence of an atlas should
be added. Thus although \ref{SS:TwoSheaf} is rather
terse it shouldn't be impossible to become acquainted
with the language as one goes along, albeit having
some examples in mind such as
quotients by foliations, dynamical systems, the fibre
$O\ts_K S$ for $O$ an orbifold supported on a knot
$K$ in the 2 or 3 sphere, $S$, 
may be necessary
for this to work.
Alternatively, one could take
the first few chapters of \cite{L-MB}, and replace
``sch\'emas'' by ``my favourite category''. 

\subsubsection{But I know 26 definitions of
 \texorpdfstring{$n$}{n}-Category
and 17 generalisations of the \'etale site,
are you saying this isn't for me ?}
\label{faq2}
One has to distinguish two different phenomenon.
The first is that
the most common case of champs
that one tends to encounter (even in algebraic
geometry) are orbifolds, which have the particular
feature that they are spaces almost everywhere. 
Consequently, the richness (even under the separation
hypothesis of \S.\ref{S:I}-\ref{S:II}) 
of the generalisation and the clear introduction/motivation
for higher category theory that it
provides get obscured. The second is the very
considerable advances in higher category theory,
notable \cite{lurie}.
Consequently, in higher category terms,
the value tends to be limited to the
space like description, albeit in very
low degrees, of things which have a comparable,
and more general, descriptions in topos 
theoretic terms.

\subsubsection{Now I'm confused,
isn't champ just a pseudonym for
orbifold, and to imagine any similar 
structure
which isn't a space almost everywhere is just nonsense.
What do you intend to do ? Distinguish between a
space and a group acting trivially on a space ?}
\label{faq3}
Basically, that is exactly the distinction that
one should make. For example the action of 
a group on a (closed) point, $\G\rras\rp$, 
is invariably trivial. On the other hand if
one expresses a space, $X$, by way of an open
cover $U=(\coprod_\a U_\a)\ra X$ and the
gluing $R:=(\coprod_{\a\b} U_\a\cap U_\b)=U\ts_X U\rras U$
of the same, then the rule for giving a map
$f:T\ra X$ expressed in terms of the equivalence
relation $R\rras U$ is 
best seen by viewing the later as a category,
{\it i.e.} there is a
functor $V\ts_T V\ra R$ for some open cover
$V=\coprod_\a V_\a$ of $T$, and the set
$\Hom(T,X)$ can be identified with
\begin{equation}\label{feq1}
\text{functors $V\ts_T V\ra R$ for some
cover $V$}/(\text{equivalence over a common refinement})
\end{equation}
and if one applies the prescription \eqref{feq1}
to the category $\G\rras\rp$ the resulting set
is the set of $\G$-torsors over $T$. Moreover
the (2-category) of equivalence relations is 
rather particular since (by the simple expedient
of taking the quotient) it's equivalent to the
1-category of sets, a.k.a. categories in which
every arrow is an identity, so the general prescription
is not to define $\Hom(T,-)$ as a set via \ref{feq1}
but as a category with objects
\begin{equation}\label{feq2}
\begin{split}
& \text{functors $F:V\ts_T V\ra R$ for some
cover $V$ and}\\
&\text{arrows natural transformations
$\xi:F_W\ra G_W$}
\end{split}
\end{equation}
for $W$ a common refinement of the covers on which
$F$ and $G$ are a priori defined, which in the
particular case of $\G\rras\rp$ leads to
Grothendieck's definition of $\rB_\G$, {\it i.e.}
the 
2-functor 
(understood weakly unless one wants to get into
large category issues \ref{faq5.3})
from spaces, $T$, to groupoids (a.k.a. a category
in which every arrow is invertible) which to $T$
associates the category whose objects are $\G$-torsors
over $T$, with arrows torsor maps.
In any case, the 
basic relation between 
champs and spaces is 
that between groupoids and sets.

\subsubsection{That's all very clever, but I'm a 
topologist, I hate 
French mathematics, and I
can mimic that sort of thing
using \texorpdfstring{$\rK(\G,1)$}{K(G,1)}'s.}
\label{faq4}
By construction there's certainly a map
$\rK(\G,1)\ra \rB_\G$, which, unsurprisingly,
is a weak homotopy equivalence, but there is
no non-trivial map in the other direction.
Nevertheless, 
the (pointed) homotopy category is an
incredibly rich category, and the addition
that Grothendieck's $\rB_\G$ makes to the
topologist's arsenal is limited. As such, the
objection not only has some substance, but,
clarifies what is arguably the principle
utility, if not the nature, of the construction,
{\it i.e.} the generally valid definition
of $\rK(\G,1)$ 
and related homotopy constructions
in whatever category.

\subsubsection{So you can make symplectic,
or holomorphic, or even characteristic
\texorpdfstring{$p$}{p} versions of Postnikov and Whitehead
towers using this stuff ?}
\label{faq5}
The number of stages, $n$, 
that one can do in the Whitehead, respectively
Postnikov,
tower is a function of having
the right definition of $n$, respectively $n+1$, sheaf
and,
\ref{SS:TwoSheaf}, champs
are the case $n=2$. More precisely, the Whitehead
tower of a 
(path connected)
space, $X$, with homotopy groups, $\pi_n$,
starts from $X_0=X$, 
and continues as $X_{n}\ra X_{n-1}$ a fibration
in $\rK(\pi_n, n-1)$'s such that $X_n$ has no
homotopy in degrees at most $n$. In particular
$X_1\ra X_0$ can (provided it exists) be taken to be the universal
cover. Nevertheless the tower is a priori only
defined up to homotopy, so anything homotopic
to the universal cover is allowed. Consequently, from the
homotopy point of view, the reason why one can
define say a holomorphic structure on $X_1$ as
soon as this exists on $X_0$ is that $\rK(\pi_1, 0)$
has a canonical realisation- 
or,  following the terminology of \cite{prst},
``canonical modelizer''-
as the discrete set
with elements $\pi_1$. Certainly, therefore, \ref{faq3},
one gets a map $X\ra \rB_{\pi_1}$, {\it i.e.}
the first step of the Postnikov tower, but,
and this is the subject of \S.\ref{S:II}-\ref{S:III}
one also gets the 2nd Whitehead stage, 
or better universal 2-cover, $X_2\ra X_1$
as an \'etale fibration with fibre $\rB_{\pi_2}$.

\subsubsection{The usual model for the 2nd Whitehead stage of
\texorpdfstring{$\bp^1_\bc$}{P\_1} is 
\texorpdfstring{$\rS^3$}{S\_3} or \texorpdfstring{$\bc^2\bsh\{0\}$}{C\_2-0} if
you want some holomorphicity, but you're claiming
there's a model which is an \'etale fibration,
whereas the only \'etale fibration over \texorpdfstring{$\bp^1_\bc$}{P\_1} 
is the identity. You on drugs ?}
\label{faq6}
Indeed if $f:X\ra \bp^1_\bc$ is an \'etale fibration,
then $f$ is an isomorphism iff $X$ is a space. A champ,
however, is not a space, and just as one can construct
$\bp^n_\bc$ as the quotient of the action
\begin{equation}\label{feq3}
\bg_m \ts (\bc^{n+1}\bsh\{0\}){\build\rras_{(\lb,x)\mpo x}
^{(\lb,x)\mpo \lb x}} \bc^{n+1}\bsh\{0\}
\end{equation}
it's universal 2-cover, $\cP^n_2$, is the quotient of the action
\begin{equation}\label{feq4}
\bg_a \ts (\bc^{n+1}\bsh\{0\}){\build\rras_{(\ell,x)\mpo x}
^{(\ell,x)\mpo \exp(\ell) x}} \bc^{n+1}\bsh\{0\}
\end{equation}
and the map $\cP^n_2\ra \bp^n_\bc$ is both a local
homeomorphism and a fibration, \ref{defn:fib}. In
particular, $\cP^n_2$ is complete K\"ahler, but it does not 
satisfy the conclusions of Hodge theory, since
by the Leray spectral sequence
\begin{equation}\label{feq5}
\rH^p (\cP^n_2, \bz) =\begin{cases}
\bz,\, p=0\, \text{or}\, 2n+1\\
0,\, \text{otherwise}\end{cases}
\,\text{but}\,\,\,
\rH^p (\cP^n_2, \Omega^q) =\begin{cases}
\bc,\, p=q\, \text{or}\, q+1\\
0,\, \text{otherwise}\end{cases}
\end{equation}
so, $\cP^n_2$ is very far not just from a space
but even from $\bp^n_\bc\ts\rB_\bz$.

\subsubsection{So if the first Whitehead stage is
the universal cover, does, by analogy, this second stage admit
some ``group object in the category of categories''
action, or even Galois theory ?}
\label{faq7}
Indeed just as a (connected)
homotopy 1-type is the same thing
as a group, 
equivalence classes of 
``group objects in the category of categories'', 
{\it i.e.} so called, {\it cf.} \ref{faq8}, 
2-groups, are \cite[Theorem 43]{baez} in 1-to-1
correspondence with (connected) 
topological 2-types, \cite{WhiteheadMaclane},
{\it i.e.} a group $\pi_1$, a $\pi_1$-module $\pi_2$,
and the {\it Postnikov class} $k_3\in\rH^3(\pi_1,\pi_2)$,
\S.\ref{SS:II.1}. Notwithstanding the subtlety, \ref{fact:group1},
that this correspondence isn't a faithful functor from
the homotopy category  to 2-types, it certainly
shouldn't surprise that 2-groups have the same
relation to the universal 2-cover/2nd stage in
the Whitehead tower as groups do to the
universal cover/1st  stage. All of which
is just the ``2 piece'' of a larger story
about $n$-groups, $n$-topoi, and $n$-Galois
theory, $n\leq\infty$, amongst which the
case $n=\infty$, \cite[A1]{lurieA}
is a bit easier than $n$ finite since
it supposes the maximal amount of local
connectivity,  
whereas the latter requires hypothesis
of the form locally $n-1$-connected and
semi-locally $n$-connected, \cite[2.15]{new},
\ref{faq3.1}.
A priori one might think that
this only applies to CW-complexes, wherein
the story began, \cite{toen}, but properly
understood, \cite[3.5]{new}, it also covers
the general pro-finite case which is valid
more or less unconditionally- indeed it may 
even be formulated without supposing locally $0$-connected,
\ref{faq3.7}- because all higher pro-finite
$\varprojlim^{(i)}$'s vanish, \cite[Th\'eor\`eme 7.1]{jensen}.
Consequently the obvious realm of the $n<\infty$ case
is point set topology wherein it constitutes one
of our recurring themes, {\it i.e.} the ``Huerwicz border''
for which one has a homotopy interpretation of
sheaf co-homology, \ref{faq1.7} \& \ref{faq3.7},
and which can actually be used to prove theorems
in point set topology, \ref{fact:395}-\ref{fact:396},
albeit we never go beyond the case $n=2$, where our
treatment closely parallels \cite{twogalois} which
in turn follows \cite{giraud} rather than \cite{lurie}.
Irrespectively however of how one arrives to
the conclusion 
the points is that
2-groups are the 
``canonical modelizer'' of homotopy 2-types
envisaged by \cite{prst}. They are properly
speaking objects of (pro)-discrete group theory,
\ref{faq2.1}, and in such terms
admit a minimal, if not quite
canonical since one has to choose
a co-cycle representing the Postnikov class, 
description which is what affords the
holomorphic (symplectic, differential or
whatever) theory of the initial stages
of the Whitehead tower, \ref{faq5}.
\subsubsection{I'm not buying that. Fibre products
in 2-categories are only  well
defined up to equivalence, so  the group axioms
are devoid of sense in a 2-category, and this whole
2-group thing is a lot of sad rubbish.}
\label{faq8}
This is 
my paraphrasing of a substantive objection by 
Nick Shepherd Barron. The most
superficial way  to deal with it is
to confine the definition of 2-group to 2-categories
admitting strict fibre products, \ref{def:FibreS},
rather than the more general definition, \ref{def:Fibre2}.
This is pathetic but  agrees with
current mathematical practice, and is wholly 
adequate for the purposes of this manuscript.
A better answer to the question is that the
right definition of 2-group is: 
\begin{equation}\label{feq6}
\text{a 2-category with 1 0-cell, and all
1/2-cells invertible.}
\end{equation}
Irrespective of a somewhat  subtle difference between
2-group equivalence and bi-category equivalence,
\eqref{eq:groupTom}, if one accepts \eqref{feq6} is
the response to the objection, then one is logically
obliged to accept that the right definition of
group is not that employed for the last 200 years, but
\begin{equation}\label{feq7}
\text{a groupoid with 1 object.}
\end{equation}
Now this is very much ``parvus error in principio
magnus est in fine'', \cite{thomas}, territory, but,
on balance, this manuscript 
(and 
one could
cite a lot of supporting evidence in the
theory of higher topoi)
tends to support the
view that \eqref{feq7} is the ``right'' definition
of group, {\it e.g.} if we were to express a (pointed path
connected) space, $X_*$, as a quotient of its
path space $P\ra X$ via the equivalence relation
$R:=P\ts_X P\rras P$ then $\pi_0(R_*)$-
pointed in $*\ts *$ would be not
only a very natural definition of $\pi_1(X_*)$
in the path connected context, \ref{SS:I.5},
but, 
a priori,
the groupoid $\pi_0(R_*)\rras \pi_0(P)=\rp$ naturally
presents itself as an objects in a category
of ``canonical modelizers'' of 1-types, and,
indeed if one didn't already know 
Grothendieck's definition,
\ref{faq3}, of $\rB_{\pi_1}$, one would 
rediscover it from asking his ``canonical modelizers''
question about path connected spaces. Nevertheless,
I'm by no means certain that \ref{feq6} is the answer
to the objection, which is, in any case, something that
has to be taken much more seriously than just appealing
to strict fibre products.
\subsubsection{Granted you might be able to define
the universal 2-cover of a space as some kind of
champ, but I bet you that the 2-cover of a champ
is a champ of champs, or some 3, or 4, category
widget that you don't know how to define.}
\label{faq9}
No. The 2-category of champs in the \'etale
(understood in the sense of local homeomorphism
in the classical topology and in the usual sense
algebraically) topology is closed for all operations
that pertain to the homotopy 2-type. The easiest
illustration is that a topological champ which
is locally simply connected and semi-locally
2 connected admits a universal 2-cover, \ref{fact:visit5}
or \ref{fact:396}. Of course algebraically one
should work pro-finitely, 
\ref{prop:376},
while the general topological
situation, \ref{fact:396}, is pro-discrete. The
intervention of pro-objects is,
however, the nature of the 1-category of algebraic
spaces, respectively of general topological spaces,
and not only occurs already at the level of classical
Galois theory, but, in the topological case 
(or, indeed algebraic if one doesn't suppose
locally Noetherian)
there's
even an issue, \ref{def:391}-\ref{def:392},
at the level of $\pi_0$, {\it i.e.} a universal
map to a discrete space need not exist.
\subsubsection{Having some closure like that is
good, but I only care about manifolds, and 
orbifolds, so all this champ this champ that
looks like overkill since I presumably only
need some sub-2-category of champ to do what
I want.}
\label{faq10}
Not really. The important distinction is
between separated and non-separated, or
better, in order to have sufficient generality
for manifolds and orbifolds ``an \'etale fibration
over separated'' albeit, \ref{claim:sep2}, such
fibrations are quite definitely separated in
a valuative sense and it's purely a matter of
convention to distinguish them. This said
as soon as one admits an interest in orbifolds,
then, \ref{eg:eg++1}, the only substantive
restriction on how far the universal 2-cover
is from a space is $\pi_2$, so, as the dimension
increases pretty much anything could happen.
Similarly, if one's point of departure is
group actions on normal algebraic varieties,
then pretty much any normal separated champ could
intervene in the description of the \'etale 2-homotopy
of the quotient, while pretty much any (separated) thing could,
{\it cf.} \ref{faq4.1}, intervene as a
closed sub-champ. Modulo, however, what we've
said about \'etale fibrations over separated,
there are large areas such as 3-manifolds or
(classical) algebraic geometry where one may
not have an a priori interest in the non-separated case.
\subsubsection{The general separated case must
be where this stops. After all
the non-separated case is the perogative of
``non-commutative geometry''. I think you should
keep your neb out of stuff you don't understand.} 
\label{faq11}
There
are several distinctions to be made. The first is
that the constructions of ``non-commutative geometry''
are not per se constructions about ``exotic'' classifying
champ such as $[\br/\bq]$, or more generally transverse
dynamics, but rather about the orbits themselves from
which, occasionally, one gets transverse information but
this is per accidens. Secondly examples like $[\br/\bq]$,
the ``non-commutative torus'', and so forth, are 
separated under any reasonable understanding of Grothendieck's
valuative criteria in the topological case, so they
tend to have extremely good properties, and admit very
easy spectral sequences for calculating (functorially
with respect to the ideas) whatever one's favourite 
co-homology theory is. Probably, 
albeit {\it cf.} \ref{faq1.8},
the first natural
example which fails to satisfy even a valuative criterion
for separation would, \cite{malgrange}, be something like the
classifying champ of the general 
holomorphic foliation on $\bp^2_\bc$, and, irrespectively,
since it's locally 2-connected this example 
has a well defined homotopy 2-group, 
{\it i.e.} not just a pro-2-group,
and
admits, \ref{fact:396}, a universal 2-cover.
\subsubsection{It might well be that this language
permits a conformal, or positive characteristic, or
whatever, definition of homotopy 2-types, but, at
least pro-finitely, so does Artin-Mazur, and it works
in all degrees, so you're just engaging in a large
scale waste of time.}
\label{faq12}
It is certainly true that Artin-Mazur, \cite{artin-mazur},
provides a definition of the pro-finite homotopy groups
of any site, 
and it's inconceivable that this definition won't
agree
with any other.
Nevertheless, to achieve it's goal
\cite{artin-mazur} has recourse to geometric realisation,
and whence the practical implementation of the theory
is usually far from straightforward, {\it e.g.} if
one wants to prove anything about $\pi_1$ one 
invariably uses the definition of
\cite[expose\'e V]{sga1}, which is indeed equivalent,
to the definition of \cite{artin-mazur}, but
this is a theorem- {\it op. cit.} \S 10 for
$\pi_1$ and \cite[3.5]{new} in general.
Our own
mantra, however, is that the correct definitions
of homotopy groups are, already for topological
spaces, Galois theoretic and the theorem is the
other way round, {\it i.e.} given sufficient 
path connectedness the Galois definition may
be realised by loops and spheres, \ref{faq1.7},
which would also appear to be the philosophy
behind \cite[A.1]{lurie}.
In particular once one dispenses with geometric
realisation in a context such as algebraic geometry
to which it is alien, one has the definitions in
a form which is intrinsic to the category being
studied, and it is in this form, \ref{faq4.5},
that one uses them, 
albeit here only in degree at most 2,
to prove the Lefschetz theorem. Similarly,
even the (trivial) GAGA theorem, 
\ref{prop:415}, 
wasn't know for the Artin-Mazur groups
prior to \cite[3.5]{new}
without hypothesis such as geometrically
unibranch.
\renewcommand{\thesubsection}{\arabic{section}.\Roman{subsection}}
\subsection{The representability theorem and its consequences.}
\label{faqR}
\subsubsection{What is the representability theorem ?}
\label{faq1.1}
It's the theorem, \ref{factdef:global}, that under 
not far from optimal, \ref{rmk:path},  hypothesis on a
topological space $Y$, {\it i.e.}
\eqref{eq:path0} and compactness, the champ/2-sheaf,
\ref{SS:TwoSheaf}, of maps, $\Hom(Y, \cX)$ to a
separated topological ({\it i.e.} represented by
a groupoid in the Grothendieck
topology of local homeomorphism) champ is again a
separated topological champ. The key to the theorem
is to prove it locally, \ref{fact:path3}, by way of
constructing a ``tubular neighbourhood'' of the graph,
which, as it happens, has limited sense in an arbitrary
topological space, so, ironically, \cite{Mp}, the algebraic,
or any other situation where one can bring infinitesimal
methods to bear is actually somewhat easier, \ref{rmk:path}.
\subsubsection{Could you give a more precise statement of
the theorem ?}
\label{faq1.2}
No. Whether here, or,
elsewhere in the introduction, you're supposed to click on the above
hyperlinks.
%
%
\subsubsection{What's all that pointing for ? Are you
trying to drag the subject back to the dark ages ?}
\label{faq1.4}
The representability theorem, \ref{faq1.1}, trivially
implies a pointed representability theorem, \ref{fact:point1},
but the critical additional fact about pointed maps,
$\Hom_*(Y_*, \cX_*)$, from a pointed
space, $Y_*$, to a pointed champ, $\cX_*$, is, {\it op. cit.},
that it is itself a space. In particular, therefore, the
path, $\rP\cX_*$, and loop spaces, $\Om\cX_*$, \ref{defn:point2},
of a pointed champ are spaces, which via the usual adjunction,
\ref{cor:point}, permits, 
modulo the correct definition of fibration, \ref{defn:fib},
a more or less instantaneous, and in a
wholly elementary way, extension of the homotopy theory of
spaces to 
(locally path connected)
topological champs. It should also be noted that
if one doesn't point, then one gets the wrong answer,
{\it e.g.} $\pi_0(\Hom(\rS^1, \rB_\G))$ is conjugacy
classes of $\G$ rather than $\G=\pi_0(\Hom_*(\rS_*^1, \rB_{\G,*}))$ itself.
As such, while I have some
sympathy with the attempts of \cite{bangor} to eliminate
pointing, its practicality is evident- on
a deeper level, I defer to \cite{vlad}.

\subsubsection{The
passage from the spaces to champs 
can't be that instantaneous. Otherwise
what's all that junk about the universal
cover in aid of ?}
\label{faq1.5}
Not really. The point is that while the
idea in any text book presentation of
``a locally path connected and semi-locally
1-connected space admits a universal cover''
is on the money, the exposition is insufficiently
functorial, and what should really be proved,
{\it cf.} \ref{faq8},
is that if $R_0:=\rP\cX_*\ts \rP\cX_*\rras\rP\cX_*$
is the path groupoid and $R_1\rras\rP\cX_*$ the
connected component of the identity, then the
quotient $[\rP\cX_*/R_1]$ is the universal cover,
\ref{fact:cover2}, and, as it happens, a certain
number of messy homotopy statements should be
replaced by statements about functors and
natural transformations between these groupoids,
\ref{factdef:cover1}. If this were the case,
equivalently if the use of the path space as
an element of a Grothendieck topology was more
diffused- it's systematic but informal in \cite{BottTu}-
then there'd be relatively little to do beyond
noting the extra complication, \eqref{eq:cover3},
that occurs in the path groupoid of a champ
rather than a space.
\subsubsection{Fair enough. Even the special case
of how to construct the universal cover of an
orbifold via loops attracts a lot of attention,
so, having done everything in such generality does
it cast any light on the developability of orbifolds,
i.e. when is their universal cover a space ?}
\label{faq1.6}
Yes. An orbifold is developable iff every tear
drop can be homotoped to a sphere, \ref{fact:tear2}.
This is a simple consequence of the relative
homotopy sequence, \eqref{eq:tear3}, but, even in
the case of an orbifold, $\mu:\cO\ra M$, supported
on a manifold, $M$, in a knot, $K$, the relative
homotopy sequence has to be applied to the embedding
$\mu^{-1}(x)\hookrightarrow\cO$, $x\in K$, where the fibre,
$\mu^{-1}(x)$ is never an orbifold, but it is a champ,
{\it cf.} \ref{faq3}. There's also a 
group of tear drops, but, arising from 
relative homotopy, it's non-commutative, and
one should, \ref{rmk:tear1}, be careful about
its structure.
\subsubsection{What's all that 1-Galois, Huerwicz stuff
at the end of this chapter in aid of ? Indeed looking
ahead, what's the chapter itself in aid of ? The results
of the pro-discrete \'etale theory are much more general,
and completely dispense with all this loop-sphere rubbish.}
\label{faq1.7}
As it happens, Huerwicz is usually stated as a relation
between homotopy groups and singular homology, which,
even for locally path connected spaces, is, subtly
different to the relation, \ref{cor:one2}-\ref{fact:one2},
between homotopy and sheaf co-homology,
{\it cf.} \cite[7.1]{lurie}. That is, however,
highly tangential to the real question, 
{\it cf.} \ref{faq12}, of why bother
about the sphere
way of doing things, when the \'etale way provides much
more general results, {\it i.e.} exactly the same
theorems but with path-connected replaced by connected.
The
best way to illustrate 
what can be gained from the sphere approach is
by considering some examples, to wit:

a) The ``Theorem of the Sphere'' for 3-manifolds,
respectively 2-orbifolds, {\it i.e.} the assertion
that if $\pi_2$, respectively the tear drop group
of \ref{fact:tear2}, is non-trivial then there is
an embedded sphere, respectively american football.

b) Idem, for complex algebraic Fano orbifolds but
where one 
(anticipates) that $\pi_2$ is always non-trivial
and one
asks for a holomorphically (whence
algebraic) embedded sphere, tear drop, or 
american football.

c) Idem as (b) but for any 1-dimensional algebraic
orbifold in any characteristic.

Now (a) is
the starting point of 3-manifold theory, and, 
to prove it one uses 
that the homotopy classes are
represented by spheres. Plainly
this appeal to the existence of
homotopically
non-trivial spheres is essential
to the demonstration, and the 
existence of such things only
follows from the more generally 
valid \'etale definition of \ref{fact:396}
a postiori,
{\it i.e.} the two definitions of $\pi_2$ are
equivalent by the 2-Galois correspondence,
\ref{prop:cor1}, for locally path connected
spaces. The same circle of ideas should
by \ref{fact:tear2} work for 3-orbifolds,
a.k.a. Thurston's orbifold conjecture,
but {\it op. cit.} doesn't give a generically
embedded 
american football which is what one
needs to get started 
(and this sort of thing, {\it i.e.} the
failure of homotopy classes to be generically
embedded may be the extent of
the difficulty in extending Smale's
$h$-cobordism theorem to champs)
whereas the current
logic of the demonstration via Thurston's
orbifold theorem is torturous. Apart from
the case of manifolds, which is a consequence
of Mori's Frobenius trick,
the only
other known
case 
of (b)
is that of algebraic surfaces with quotient
singularities, \cite{JamesSean}, 
which, b.t.w., is a very long enumeration by cases,
and, arguably,  the question is best addressed algebraically.
Nevertheless under (the much) stronger hypothesis, \cite{SiuYau},
of positive sectional curvature one can do it analytically,
and, of course, the starting point 
of {\it op. cit.} is that classes in
$\pi_2$ are represented by spheres. 
The nature of (c) is somewhat different,
since the difficult case, \ref{lem:IV117},
of orbifolds with moduli $\bp^1$,
which currently has no algebraic proof over $\bc$,
needs Van Kampen, and a ``nice'' topology in
which to employ it. The former is, however,
a theorem of \'etale homotopy, \cite[Expos\'e IX.5.1]{sga1},
and, the ``nice'' topology wherein it is to be employed is
just the classical topology. As such, the use
of loops and spheres in (c) is only a matter of convenience,
albeit a useful one. A (tangential) curiosity
here is that the specialisation theorems \cite[Expos\'e X.3.8]{sga1}
fail, \ref{rmk:FeitThomson}, absolutely (Feit-Thomson)
to do the said difficult case, \ref{lem:IV117}, 
in characteristic 2, while the classification of
finite simple groups doesn't leave much room for
manoeuvre in characteristic 3 either, so, as far as I know,
(c) is open for tame orbifolds in characteristics 2 or 3.
\subsubsection{If the whole loop-sphere thing is
so useful why didn't you develop it in the general
non-separated case, or, at least under some valuative
criteria such as paths have unique limits, or even that
locally the champ is the classifier of a discrete group
action ?}
\label{faq1.8}
One should first distinguish the sub-question of
whether there are any better ``separation hypothesis''
than those we have employed, which, {\it cf.}
\ref{faq10}, are really designed for 
conformal homotopy of spaces and orbifolds.
As such, a case which
we've not covered, is where the lack of separation
is the result of a global group action, {\it e.g.}
$[\br/\bq]$, so called ``non-commutative torus'',
{\it etc.}. On the other hand, everything
extends to such global
cases trivially, so the better question is:
are there some natural examples with good
separation properties that don't come from
global group actions, and, to date, I don't know any.
Indeed, one has a very rich set, \cite{ncmt}, of examples in
the form of the classifying champs of foliated
algebraic surfaces, wherein to imagine that
one would find an example which is locally but
not globally the classifying champ of a discrete
group action is just silly. Plausibly, these
examples are valuative separated- it almost
follows from the uniformity of
uniformisation of \cite{ncmt} and I've tried quite
hard to prove it, but \cite{malgrange} is very  much in the
opposite direction. Consequently, 
modulo trivia,
it seems to be
the case that one should either suppose separated
or arbitrary. This said there is a good case
for developing the sphere theory  arbitrarily,
{\it e.g.} 
it's utility along the lines of (a)-(c)
of \ref{faq1.6} in the aforesaid foliation
example should be pretty apparent. On the
other hand, in such generality the path
space would fail to be separated, and since
we never ever employ coverings, or groupoids (or indeed even the
word space) whether it be an \'etale atlas
or the path space, and the resulting \'etale,
respectively path, groupoid unless these objects
are themselves separated. Thus, for example,
we regard non-separated spaces as champs, rather
than spaces, and, while, perhaps, such generality
should have been attempted, it would have been
wholly independent of our use of loops and spheres
via the path groupoid in \S.\ref{S:I}-\ref{S:II} anyway.

\subsection{2-Galois theory via loops and spheres.}
\label{faqII}
\subsubsection{The first section of this chapter 
looks to me like a lengthy rant in group co-homology.
What's going on ?}
\label{faq2.1}
There's plenty of literature about 2-groups,
fundamental 2-group,
crossed modules, coherence theorems, {\it etc.}
all of which is the opposite of help for 
doing \'etale homotopy.
The basic thing to keep in mind is, \ref{faq5},
that size matters inversely, and
what we're dealing with here is an
aspect of discrete group theory.
Consequently,
one would ideally like to define a 2-group
as the 2-type $(\pi_1, \pi_2, k_3)$ of \ref{faq7}.
Unfortunately, this doesn't quite work, 
\ref{fact:group1},
and
what one has to do, \eqref{eq:group1}-\eqref{eq:group2},
is replace $k_3$ by a normalised co-cycle,
$K_3:\pi_1^3\ra\pi_2$. The difference between
such a thing and a crossed module or  strict
2-group is analogous to the difference between
co-homology groups and injective resolutions,
{\it i.e.} the latter are excellent in theory,
but one would never use them in practice. 
More details are in \cite{baez}, albeit we need
this sort of thing for actions of a 2-group
on a groupoid, which is more general than {\it op. cit.},
and the resulting development is most of the
chapter. It's summarised in \ref{sum:group1}, and
essentially the only (guest) appearance of the more common
approach via crossed modules/strict 2-groups is
in the scholion, \ref{warn:LeftRight},  on left vs right 2-group
actions which follows. It may, however, be usefully noted
how one goes from a strict 2-group, $\mygg$, to its
2-type. Specifically, if to fix ideas $\mygg$ were
contained in the automorphisms (be they weak or
strong) of a category $\cC$, then $\pi_1(\mygg)$
is the classes of functors in $\mygg$ modulo
equivalence via natural transformations, 
and $\pi_2(\mygg)$ are the natural transformations
which stabilise the identity. Consequently if
for each $\om\in\pi_1(\mygg)$ we choose a functor
$F_\om$ representing it, there are natural transformations
(unique up to $\pi_2$) $\phi_{\tau,\om}:F_{\tau\om}
\Rightarrow F_\tau F_\om$ and a natural transformation
$S_{\s,\tau,\om}$ (unique given $\phi$) such that
the commutativity of
\begin{equation}\label{feq2.1}
 \xy 
 (0,0)*+{F_{\s\tau\om}}="A";
 (30,0)*+{F_\s F_{\tau\om} }="BB";
 (60,0)*+{F_\s F_\tau F_\om}="B";
 (60,-18)*+{ F_\s F_\tau F_\om}="C";
 (30,-18)*+{F_{\s\tau}F_\om  }="D";
 (0,-18)*+{F_{\s\tau\om}  }="E";
    {\ar@{=>}^{  \phi_{\s, \tau\om}} "BB";"A"};
    {\ar@{=>}^{(F_\s)_*\phi_{\tau\om}} "B";"BB"};
    {\ar@{=}_{} "C";"B"};
    {\ar@{=>}_{F_{\om}^* \phi_{\s,\tau}} "C";"D"};
    {\ar@{=>}_{\phi_{\s\tau,\om}} "D";"E"};
    {\ar@{=>}_{(F_{\s\tau\om})_* (S_{\s,\tau\,\om}) } "A";"E"};
 \endxy 
\end{equation} 
affords a co-cycle $S_{\s\,tau,\om}$, a.k.a. the
Postnikov class $K_3$.
\subsubsection{So if I apply this recipe to the
automorphisms of the universal 2-cover (of which
I wouldn't mind some more explanation b.t.w.) presumably
I get the fundamental 2-group 
\texorpdfstring{$\Pi_2$}{Pi\_2}, or,
equivalently the 2-type 
\texorpdfstring{$(\pi_1, \pi_2, k_3)$}{(pi\_1, pi\_2, k\_3)} ?}
\label{faq2.2}
More or less exactly right. The only error here
is that if $\cX_2\ra\cX$ is the universal cover
of a champ, and one employs
the usual rules for passing from commutative
diagrams to 2-commutative diagrams then the 2-category
in which automorphism is to be understood,
$\et_2(\cX)$, \ref{defn:cor1}, isn't quite the naive
$\underline{\mathrm{Cham}}\mathrm{p}\underline{\mathrm{s}}/\cX$
unless $\cX$ is a space. 
As to the universal 2-cover itself:
under our separation
hypothesis, \ref{faq10} \& \ref{faq1.8}, a (pointed) topological
champ $\cX_*$ has a path space $P:=\rP\cX_*\ra \cX$ and
the Whitehead hierarchy may be written as:
\begin{equation}\label{feq2.2}
\text{$\cX_0:=\cX$ is the classifying champ of
the groupoid $R_0:=P\ts_\cX P\rras P$.}
\end{equation}
whenever $\cX$ is path connected, and if it's
also semi-locally 1-connected
\begin{equation}\label{feq2.3}
\text{the classifying champ, $\cX_1$ of
the groupoid $R_1\rras P$.}
\end{equation}
for $R_1$ the connected component of the identity,
is the universal cover, while, finally if $\cX$
were locally 1-connected and semi-locally 2-connected
\begin{equation}\label{feq2.4}
\text{the classifying champ, $\cX_2$ of
the groupoid $R_2\rras P$.}
\end{equation}
for $R_2$ the universal cover of $R_1$, is the
universal 2-cover. In particular, therefore, the
arrows of $R_2$ which become the identity in
$R_1$ are the $P$-group $S:=P\ts\pi_2$ and $R_2=R_1/S$,
while $\cX_2\ra \cX_1$ is a locally constant 
gerbe in $\rB_{\pi_2}$'s- \ref{faq3}, \ref{faq6}.
\subsubsection{By analogy, therfore, if we
replace group  by 2-group; covering space
by covering champ; and set  by groupoid shouldn't
there be a 2-Galois equivalence between
covering champs and groupoids on which 
\texorpdfstring{$\Pi_2$}{Pi\_2} acts ?}
\label{faq2.3}
Again 
modulo what we've already said 
in \ref{faq2.2}
about the difference 
(which pertains to the 1 and 2 cells)
between
$\et_2(\cX)$ and coverings in
$\underline{\mathrm{Cham}}\mathrm{p}\underline{\mathrm{s}}/\cX$
this is exactly right. In particular
the $0$-cells are the \'etale covers
of which
the generally (in any site, {\it e.g.}
connected but not path connected spaces) valid 
definition is that the fibre over any sufficiently
small open, $U$, is $U\ts \cG$ for $\cG$
a discrete groupoid which, in the topological
and
locally 1-connected case, is equivalent, \ref{cor:sep3}, to the
more succinct property of being an \'etale
fibration of champs. Unsurprisingly, therefore,
the 2-functor
\begin{equation}\label{feq2.5}
\et_2(\cX) \ra \Grpd(\Pi_2)
\end{equation}
affording the 2-Galois equivalence, \ref{prop:cor1},
is the fibre functor which
takes an \'etale fibration $q:\cY\ra\cX_*$ over
a pointed champ and sends it to the ($\Pi_2$-equivariant)
groupoid $q^{-1}(*)$. 
\subsubsection{Neat. Any chance of some examples
to help me develop some intuition for the
2-Galois correspondence ?}
\label{faq2.4}
Indeed there is. Numerous examples are
provided in \S.\ref{SS:II.6}-\ref{SS:II.7},
ranging from the perennial favourite of
group extensions, {\it i.e.} $\et_2(\rB_{\pi_1})$,
to the case of finite $\pi_1$ wherein 
(for $\pi'_1$ a sub-group)
one
encounters classifying champs under 
the action of
extensions
\begin{equation}\label{feq2.6}
0\ra \bg_m^d \ra E \ra \pi'_1\ra 0
\end{equation}
and the extensions themselves can,
\ref{eg:tor2}, be expressed via integration
over spheres. Irrespective of intermediate
examples such as non-associative extensions,
\ref{eq:eg3}, these 2-examples contain 
(given the GAGA theorem \ref{prop:415})
everything
in the specific case of 1-dimensional 
algebraic champs in characteristic zero
addressed in \cite{noohiuni}. In addition
the various alternatives for defining the
Postnikov class such as,
\ref{fact:pos2}, the $d_3^{0,2}$
in the H\"oschild-Serre spectral sequence
should aid the intuition.
\subsection{Pro 2-Galois theory.}
\label{faqIII}
\subsubsection{What's the main difficulty, or
indeed difference, between the loop-sphere
approach and the (pro)-\'etale theory ?}
\label{faq3.1}
In the first place one should distinguish
between pro-finite and pro-discrete. The
latter, \ref{SS:III.9}, requires hypothesis
of the form ``semi-locally-$n$-connected''
to work, 
\ref{faq3.7},
whereas the latter is valid more
or less unconditionally. As such, the 
specifics of the respective cases are
somewhat different, but the root of the
difficulty
in either case is that one doesn't have
a (meaningful) notion of path space, so
the description \eqref{feq2.2}-\eqref{feq2.4}
of the Whitehead hierarchy as groupoids
acting on the same contractible space,
{\it i.e.} the path space,
is invalid.
The most specific manifestation of this in the
pro-finite theory is that in all probability
the universal cover doesn't exist, while in
the pro-discrete theory the game is to 
give conditions- 
``locally-$n-1$-connected and semi-locally-$n$-connected''-
such that the universal $n$-cover exists.
\subsubsection{Surely if there's no path space,
and you don't plan to use geometric realisation,
there's no relation between the loop-sphere
theory, and the pro-\'etale theory ?}
\label{faq3.2}
You'd be surprised. The basic thing about
the formulae \eqref{feq2.2}-\eqref{feq2.4} for
the Whitehead hierarchy in the topological case
is that one has a groupoid, $R_0$, acting on a
space whose quotient is the thing, $\cX_0$, whose covers
we wish to study, and the universal 1, respectively 2,
cover are obtained from the connected component, $R_1$,
of the identity, respectively it's universal cover.
More generally, however, all (connected) \'etale fibrations
of $\cX_0$ can be expressed as a groupoid 
$R\rras P$ acting on the path space where $R$
is a (not necessarily connected) \'etale cover
of components of $R_0$- 
which comes out in the wash,
\ref{fact:cor5}, in the proof of the
2-Galois correspondence, \eqref{feq2.5}.
Now, in practice, particularly
in the pro-finite theory, one doesn't need a
similar fact for all covers 
(in a more or less arbitrary site)
at once, but only for
a finite number, and this can be achieved by replacing
the path space $P$ by a sufficiently fine \'etale
atlas $U$-
\eqref{eq:sep2}. 
Consequently, although 
the path approach is logically independent
from the pro-\'etale approach, there is 
a similarity in the constructions without
which it's somewhat difficult to follow what
is going on.
\subsubsection{That meta-rule of replacing
the path-space by a sufficiently fine \'etale
atlas looks useful, but I don't see how to
make Galois objects from it, and presumably 
you need some sort of 2-Galois objects which
are something else again?}
\label{faq3.3}
Indeed the utility of the path space, or, better
the path fibration, is that it lifts to every
\'etale fibration, {\it i.e.} it comes with
a certain universality, whereas the notion of
``sufficiently fine \'etale atlas'' supposes 
that it covers some a priori thing in the
hierarchy of \'etale 2-(or indeed 1)-covers
of interest. Thus, one needs the right notion
of Galois, which
modulo the distinction between $\et_2(\cX)$
and $\underline{\mathrm{Cham}}\mathrm{p}\underline{\mathrm{s}}/\cX$
of \ref{faq2.2}, is, 
\cite{noohi1},
in the 1-Galois, equivalently representable covers, context 
exactly the definition,
\eqref{eq:plag1}, of \cite[expose\'e V]{sga1}.
The definition of 2-Galois is, however,
rather more involved.
\subsubsection{So when exactly is a map of
champs 2-Galois ?}
\label{faq3.4}
In the first place given any map
$q:\cY\ra\cX$ of champs in $\et_2(\cX)$,
{\it i.e.} proper \'etale covering in the algebraic
setting, there is a unique (up to unique equivalence)
factorisation $\cY\xrightarrow{p}\cY_1\xrightarrow{r}\cX$
into a locally constant gerbe followed by a 
representable cover, \ref{claim:sep1}, which
one can, and should, consider as the Whitehead
hierarchy, \eqref{feq2.2}-\eqref{feq2.4}, in
the particular rather than the universal. As such
there is an obvious necessary condition for the
gerbe $p$ to give rise to a piece of $\pi_2$,
{\it i.e.} $r_*$ should be an isomorphism on
$\pi_1$. This is the quasi-minimality condition 
of \ref{def:funct103}, wherein perhaps 2-minimal
may have been better notation since minimal in the abstract
setting of \cite[expos\'e 195]{fga2} becomes
an isomorphism on $\pi_0$ in practice. Another
obvious condition is that $r$ should be Galois,
and, irrespectively the unicity of the aforesaid
factorisation implies- see \ref{faq2.1}
for notation- a map
\begin{equation}\label{feq3.1}
\pi_1(\Aut_{\et_2(\cX)}(q)) \ra 
\pi_1(\Aut_{\et_2(\cX)}(r))=
\text{Usual Galois group of $\cY_1/\cX$}
\end{equation}
and the final piece in the definition of
2-Galois, \ref{def:noplag2}, is that \eqref{feq3.1}
should be both surjective and admit a section.
\subsubsection{So then you'll have a directed system
\texorpdfstring{$q_i\ra q_j$}{q\_ i -> q\_ j}, giving
rise to a directed system of 2-groups
\texorpdfstring{$\Aut(q_i)$}{Aut(q\_ i)}, the limit
of which is, presumably, the pro-fundamental 2-group?}
\label{faq3.5}
Yes and no. A postiori this works, but a priori it
may well be
wide of the mark since \eqref{feq3.1} will invariably
have a non-trivial kernel, 
\ref{claim:noplag2},
in the pro-finite context,
{\it i.e.} inevitably 2-Galois cells have 
(even modulo equivalence)
too many 
automorphisms. This does not, however, exclude that
\eqref{feq3.1} becomes an isomorphism in the limit,
which, in turn is what has to checked, {\it i.e.}
if $\Pi_2$ is defined as the limit of the automorphisms
of 2-Galois cells, then the resulting $\pi_1(\Pi_2)$
in the general nonsense 2-group sense of \ref{faq2.1},
agrees with the usual 1-Galois sense of $\pi_1$. By
the very definition of weak pro-2-categories this 
is equivalent, \ref{fact:365}, to the
existence of a Postnikov Sequence as defined in
\ref{def:PosSeq}. As such, the critical thing to
do is the (rather large) diagram chase of \ref{SS:III.5}
which shows the existence and uniqueness of Postnikov
sequences. After this, the meta-rule ``path space
<-> sufficiently fine atlas'', {\it cf.} \ref{faq3.2}, leads to a largely
mutatis mutandis demonstration of the pro-finite-2-Galois
correspondence \ref{prop:376} relative to its
loop-sphere variant \ref{prop:cor1}. There are,
however, some tautological bonuses along the way such as 
the Whitehead theorem, \ref{Whitehead}, that if $f:\cY\ra\cX$ 
is an isomorphism on the homotopy groups $\pi_i$, $i\leq 2$,
then
it's already an isomorphism of their fundamental 2-groups $\Pi_2$.
\subsubsection{Isn't all of this in the book of Giraud
``Cohomologie
non ab\'elienne'' ?}
\label{faq3.6}
Surprisingly enough the answer is mostly no, albeit
that for anything which is a cohomological rather
than homotopical proposition the answer is mostly yes,
while if we pose the same question for Lurie's
``Higher topoi'', the answer is, at least implicitly,
yes to everything.
Unquestionably the reason for the lacunae
in \cite{giraud} is largely
historical, {\it e.g.} interest in defining a
non-abelian $\rH^2$, lack of development of 2-groups
and bi-categories at the time of writing, {\it etc.},
and, consequently, \S.\ref{S:II}-\ref{S:III} can,
reasonably, be viewed as the homotopy complement
of \cite{giraud}. In particular, therefore, 
everything about co-homology of locally constant
pro-finite 
sheaves, or, indeed non-abelian links in the sense of {\it op. cit.},
follows (easily) from the 2-Galois 
correspondence, \ref{prop:376}. As such we're once more
in ``parvus error in principio magnus est in fine'' territory,
and the following can reasonably be stressed,

(a) Giraud's basic question of how to continue 
the cohomology sequence of non-abelian groups to
$\rH^2$ already has sense, and is arguably subordinate
to, defining $\rH^1$ of pointed sets. 

(b) Giraud cohomology 
is so tautological 
({\it cf.} \cite[7.2.2.14]{lurie} which is the
same principle, albeit that a topos
tautology isn't quite the same as
a site tautology)
as to pose
the question as to whether it ought not to be
considered logically prior to sheaf co-homology.
The two of course agree
since Giraud cohomology vanishes on injectives, 
\cite[IV.3.4.3]{giraud}, 
but the tautology is perhaps most evident
via hypercoverings
\ref{fact:384}.

(c) The relation between homotopy and co-homology
is best seen as an adjunction, which eventually
becomes the Galois correspondence. It is, therefore,
naturally pro, \cite[A.1.2]{lurie}, \ref{faq3.7},
and the
tautological elegance of Giraud co-homology
passes to the pro-finite, or, more generally,
pro-discrete setting- \ref{def:3811}- which, inter alia,
needs to be compared 
with the pro-sheaf definition
(partial definitions appear to be
a cottage industry) in the abelian case, \ref{defn:3810}.

(d) A useful intermediary in the adjunction
between homotopy and co-homology
is the dual on the latter
{\it i.e.} homology, \ref{schol:386}. As such
homology is (always) tautologically dual
to co-homology, and
should not be confused with a 
theorem (Grothendieck-Poincar\'e-Verdier duality)
equivalently co-homology with compact support.
\subsubsection{Those pro-discrete homotopy theorems
look completely new, but it's a bit odd that
the existence of a universal cover of a locally
connected (but not necessarily path connected)
and semi-locally 1-connected (in the \'etale 
rather than the path sense) wasn't proved in the 30's.
Are you sure it's right ?}
\label{faq3.7}
Point set topology didn't quite die in the 30's,
and I imagine it's in the literature somewhere
but so far I haven't been able to find it- 
even last week's localic \cite[2.15]{new} is a
little different.
Arguably
the trick for proving the theorems about $\pi_1$,
{\it i.e.} a locally connected (point set topology
sense) space has a universal cover iff it's
semi-locally 1-connected in the \'etale sense,
\ref{fact:395}, and its $\pi_2$ variant, \ref{fact:396} is to
realise that in full generality there is a 
Huerwicz theorem for $\pi_0$, \ref{fact:392}, and
that the right definition of $\pi_0$ is as a
pro-set, \ref{def:392}. This in turn sheds light
on the pro-finite theory which ought to be considered
as un-conditional, {\it i.e.} it doesn't even
need local connectedness, \ref{ProFF}.
\subsection{Applications to algebraic geometry.}
\subsubsection{What's the status of the GAGA theorem
for holomorphic champs ?}
\label{faq4.1}
Our present interest is only the 2-category 
$\et_2(\cX^{\hol})$, for $\cX^{\hol}$ the holomorphic
champ associated to an algebraic champs,
$\cX$ of finite type over $\bc$, and in this
case, it easily, \ref{prop:415}, follows
from the comparison of \cite[Expos\'e XVI.4.1]{sga4}
between \'etale and classical co-homology that
$q:\cY\ra\cX^{\hol}$ is algebraic iff it's proper.
If more generally, however, $q:\cY\ra\cX^{\hol}$ were
just a proper map to an algebraic champs, then 
even if $\cX^{\hol}$ were the moduli space of $\cY$,
it
needn't be algebraic, {\it e.g.} glue two affine lines at
infinitely many points then take the classifying
champ of $\pm 1$ acting by interchanging the lines.
As such, without the \'etale covering hypothesis, it's
difficult to avoid supposing that $\cX^{\hol}$ is proper.
Even, however, supposing this the difficulty in 
reducing to the results of \cite{sga4} is that
if $q$ isn't \'etale 
then one may not be able to factor out the stabiliser
of every generic point in a way akin to \ref{faq3.4}
without further hypothesis such as unibranching.
Of course if $\cY$ is normal one has both this
and the local results of \cite{ArtinLocalMonodromy}
so such a $\cY$ is algebraic. Otherwise, in full
generality, the question looks to be open, and
probably requires the extension of the algebraic
space case, \cite{ArtinContraction}, to champs.
In any case, it is not, due to its 
topological character,
reducible, as some have asserted, to algebraisation of coherent sheaves.
\subsubsection{But in positive
characteristic there's no GAGA theorem, and
I want to calculate \texorpdfstring{$\pi_2$}{pi\_2}. What should I do ?}
\label{faq4.2}
The basic tool, \ref{fact:386}, is the short exact sequence
\begin{equation}\label{feq4.1}
0\ra \rH^2_{\mathrm{cts}}(\pi_1, \uZ)
\ra \rH^2(\cX,\uuG) \ra \Hom^{\pi_1}_{\mathrm{cts}}(\pi_2, \uZ)
\xrightarrow{A\mpo A_* K_3} 
\bigl( \rH^3_{\mathrm{cts}}(\pi_1, \uZ), \mathrm{obs} \bigr)
\end{equation}
which holds for any locally constant sheaf $\uZ$.
In the topological case this is the weight 2 term
in the H\"oschild-Serre spectral sequence, and,
although there is no universal cover, the resulting
exact sequence, \ref{feq4.1}, continues to hold.
One also has a proper base change theorem,
\ref{fact:Lef13} for $\Pi_2$, but the smooth
base change theorem can fail. On the plus side,
however, $\pi_2$, unlike $\pi_1$, of a non-proper
variety only depends on the geometric fibre, \ref{lem:l21}.
\subsubsection{You're not seriously asserting that
``prime to p'' smooth proper base change 
from the generic to the special fibre
fails
for \texorpdfstring{$\pi_2$}{pi\_2} ?}
\label{faq4.3}
Remarkably enough, and the validity of such
base change whether for $\pi_1$ or co-homology
notwithstanding, yes. A counterexample is
bi-disc quotients, \ref{fact:sp1}, which if
for example it were given by the action of
$\mathrm{SL}_2(\cO)$ for $\cO$ a real quadratic
ring of integers any sufficiently large $p$ which is inert in $\cO$
will do. This doesn't contradict the corresponding
theorem in co-homology since the failure of 
``prime to p'' specialisation for $\pi_2$ is 
balanced in \ref{feq4.1} by the failure of specialisation for
all of $\pi_1$- ``prime to p'' group co-homology
depends on the ``prime to p'' sub-groups, not
the quotient groups, \cite[III.10.3]{brown}.
Consequently, I don't see how one can get
a ``prime to p'' specialisation theorem for
$\pi_2$ without an extremely strong, but
of Huerwicz type, hypothesis, \ref{fact:sp2},
that specialisation is an isomorphism on the
``prime to p'' group co-homology of
$\pi_1$. This also tends to expose the limits of the utility,
especially (Feit-Thomson) in characteristic 2,
of the ``prime to p'' specialisation theorem
for $\pi_1$. In particular, therefore, \ref{faq1.7},
\ref{rmk:FeitThomson}, it appears to be open
as to whether (tame) hyperbolic triangular
orbifolds in characteristic $2$ or $3$, {\it cf.} \ref{faq1.7}.(c), have
non-trivial $\pi_2$, \ref{lem:IV117}.
\subsubsection{Granted
Lefschetz for \texorpdfstring{$\pi_0$}{pi\_0} and 
\texorpdfstring{$\pi_1$}{pi\_1} 
of champs isn't in SGA, but everything
here has a moduli space, so the translation
must be trivial ?}
\label{faq4.5}
Again ``parvus error in principio magnus
est in fine'', and already for $\pi_0$ the
SGA hypothesis aren't quite right. Specifically,
the hyperplane section is supposed nowhere
a divisor of zero, ``hypoth\`ese peut-\^{e}tre
superflue'', \cite[Expos\'e XIII.4.2]{sga2},
and the variation $X\ra S$ is supposed flat.
In order to see why eliminating such excessive
hypothesis is essential to the nature of the
question, the critical point to realise is that
the (pro-finite)
Lefschetz theorems have a built in gratuitous induction step.
The best known example of this is (up to some
minor regularity hypothesis) if $i:H\hookrightarrow X$
is a hyperplane in a variety of dimension at least 2,
$i_*$ is  a surjection on $\pi_1$ because $i$ is
not only an isomorphism on $\pi_0$, but so is any base
change of it by a finite surjective map $X'\ra X$.
To fully exploit this inductive structure (even over
$\bc$) one has to,

(a) Work with champs, even if the goal was varieties-
{\it cf.} the above example for $\pi_1$, with the
definition of 2-Galois, \ref{faq3.4}, and the analogous
argument for $\pi_2$, \ref{cor:l21}.

(b) Work with necessarily non-flat families defined
by intersecting two (necessarily not always distinct) divisors-
the first appearance of the implied construction is
\eqref{l0p3}.

Item (b) certainly already demands avoiding the
flatness and regularity hypothesis of \cite{sga2}, but
even if one were prepared to assume them, they're not
easily conserved in the generality necessary to do (a),
{\it e.g.} to just copy and paste from \cite{sga2} would
at least
require supposing that all champs are tame, a.k.a. ignore
any ``p''-phenomenon in $\pi_2$.
\subsubsection{I suppose that's true that the SGA strategy
opens a pandora's box of issues about flatness and tameness,
but these are already present for \texorpdfstring{$\pi_0$}{pi\_ 0},
while even admitting 
that you can do it by induction, it has to start somewhere, and 
it's hardly the case that starting in
degree -1 is an option?}
\label{faq4.6}
Juxta modum, {\it i.e.} if there were a Lefschetz theorem
in degree $-1$ it would have a surjectivity corollary for
$\pi_0$ akin to the the more well known example of how to go
from $\pi_0$ to the surjectivity of $\pi_1$ in \ref{faq4.5}.
However it's the corollary that one needs for the induction, and
this, \ref{fact:l01}, is true, to wit:
\begin{equation}
\label{feq4.2}
\text{a hyperplane section of a connected champs
of positive dimension is non-empty.}
\end{equation}
Now we have the starting point of the induction,
leading to a small, but subtle strengthening, \ref{fact:Lef02},
of Lefschetz for $\pi_0$, from the demonstration of
which it's already
largely clear how one proceeds from there to $\Pi_2$-
\ref{prop:Lef21},\ref{cor:LLef21}- via the
corresponding results for $\pi_1$- 
\ref{prop:Lef11}, \ref{prop:l11}-
and, indeed, modulo the right definitions in terms of
higher categories, to all homotopy groups.
\subsubsection{For a trivial induction it does
seem to go on a bit. I suspect there's a catch.}
\label{faq4.7}
Not really. The reasons for the length from
the starting point \eqref{feq4.2}/\ref{fact:l01}
to the concluding \ref{cor:LLef21} may be identified as
follows,

(a) Due to its non-commutative nature $\pi_1$ has to
be treated separately, \S.\ref{SS:L1}.

(b) Modulo the right definitions, the way to organise
the induction for the higher homotopy groups is 
an appropriate use of
Deligne's descent spectral sequence, \cite[5.3.3]{deligne3},
which, arguably, would be easier in the general rather
than the specific way it's employed in \ref{prop:Lef21}.

(c) We have to correct/make up the definitions of
homotopy depth as we go along, \ref{schol:l01}, 
\ref{cschol:Lef11}, \ref{cschol:l11}.

As such, the only substantive issue is (c), wherein
the SGA definition of homotopy depth, $d$, at a closed point.
$x$, of a variety $X$ over a field, $k$, 
is that all local \'etale co-homology $\rH^q_x$
of all (not necessarily commutative if $q=1$)
locally constant sheaves
vanishes for $q<d$- which, for
example, is at least the dimension if $X$
is l.c.i., \ref{ex:lci}. This definition is,
however, 
useless because, \eqref{CschLef1}, unlike algebraic depth
it doesn't imply 
what one actually needs, {\it i.e.}
the rectified vanishing condition,
\cite[Expos\'e XIII.4.3.D\'efinition 2]{sga2},
\begin{equation}\label{feq4.3}
\rH^q_x=0,\,\,\text{whenever:}\,\, q+ \mathrm{Trdeg}_k k(x) <d
\end{equation}
which, nevertheless, is referred to as
rectified homotopy depth
in {\it op. cit.}.
Consequently, we don't follow the SGA usage, and,
modulo some further precision, \ref{def:l11}, on the topology 
being employed, {\it etc.} just call \eqref{feq4.3}
the homotopy depth. Irrespectively, homotopy depth
is a local condition, and what one would like to
know is that (much as for algebraic depth, albeit
without any regularity
assumptions) that on cutting by a hyperplane the
homotopy depth decreases by at most 1. To put the
difficulty in perspective: in as much as it pertains
to $\pi_0$, the ``proof'', \cite[Expos\'e XIII.2.1]{sga2}
is wrong, \ref{cschol:l11}, and the correction in
 \cite[3.1.7]{joinsAndIntersections} needs a non-trivial
trick, while, a priori, what we need to complete the
proof, \ref{prop:Lef21}, for $\Pi_2$ is an analogous
statement for $\pi_1$. Fortunately, however, this is
overkill, and it suffices to know the $\pi_1$
(or, in general $\pi_{q-1}$)
analogue of  \cite[3.1.7]{joinsAndIntersections} 
for the generic hyperplane, \ref{prop:Lef11}, which
in turn can be fitted into the general (supposing
an appropriate higher category definition of $\pi_q$) 
induction so as to bypass the thorny (local) question
of how the homotopy depth behaves on cutting with a
hyperplane. Alternatively if one's ultimate interest
were smooth, or more generally lci, one could profit
from \ref{ex:lci} to just build good local behaviour, upon
cutting with a hyperplane, into the definition of
homotopy depth. 

\subsubsection{Irrespectively of whether you have the definitions
to hand of higher homotopy groups that you want, 
the definition of co-homology is there, and the resulting
Lefschetz theorem that your strategy implies is 
technically better than what can be done via Poincar\'e duality,
so why haven't you done the co-homology theorem beyond
degree 2?}
\label{faq4.8}
Even for cohomology, it
doesn't work without the right higher category definitions.
As has been said, \ref{faq4.5}.(a), 
the gratuitous induction step for $\Pi_2$ is 
not Lefschetz for
$\pi_q$ of spaces, $q\leq 1$,
but of champs. Amongst these
two Lefschetz theorems $\pi_0$  
for champs is identically the
same theorem for spaces since
the moduli map is an isomorphism
on $\pi_0$, but the $\pi_1$
theorem for champs, even those
which are locally constant over
a space, {\it cf.} \eqref{feq5}
even though it isn't algebraic,
is not the $\pi_1$ theorem
for spaces. Similarly, to get
the induction step whether for
$\rH^3$ or $\pi_3$, requires the
$\pi_q$, $q\leq 2$ steps for 3-spaces,
\ref{faq5.5}, which for
$q\leq 1$ ought to reduce to the champ theorem
by a variation of the above moduli
argument, but for $\pi_2$ the 3-space
Lefschetz theorem is more general than
the champ one. Consequently, even in the
presence of \cite[3.5]{new} (which certainly
applies to the \'etale topos of a champ),
and even curbing one's ambition to higher Lefschetz
for champs, or just spaces, rather than $n$-spaces,
the $n$th step of the
induction has to be run in a sufficiently large
$n$-category of $n$-spaces, or, possibly with a
modicum of intelligence $\infty$-topoi. 
This said, everything
in such a $n$-category would be locally constant
and projective (whence a fortiori C\v ech, and devoid
of hypercovering issues) so either cobbling together an {\it ad hoc}
definition in order to push the induction through or
using the $\infty$-topoi definition directly
should work fine, although it will involve
leaving the familiar 2-category terrain of this
manuscript.
\subsection{Miscellaneous: language, notation, sets, and generalisation.}
\subsubsection{Why do you use champ instead of 
the commonly
accepted (mis)translation stack ?}
\label{faq5.1}
I used to, but in thinking about problems such
as (a)-(c) of \ref{faq1.7} I found that it
clouded the mind because (and I suspect that
this was Mumford's, not unreasonable in the
case of $\cM_g$, intention) it created an association
with scheme, which, {\it e.g.} \eqref{feq5}, is
not entirely warranted. Now while there may be a whole
sub-story about how bundle got translated as faisceaux 
then got translated back as sheaf (from 
which the amusing corollary that no-one
has a working translation of gerbe)
the word champ was
coined in French because it correctly reflects the
idea, {\it i.e.} ``faisceaux de niveau 2'', and, functorially
with respect to the ideas, \ref{SS:TwoSheaf}, the right
translation, without any representability
hypothesis on the diagonal or of the existence
of an atlas, is 2-sheaf. Plausibly, therefore,  
under the said representability hypothesis, 2-space,
or maybe 2-manifold,
might be the way to go, but, for the moment following
the French usage 
(particular given the propensity for $n$-stack,
and whence the ridiculous $0$-stack instead
of sheaf)
seems to be the least bad option. This said the
choice of the word champ rather than 2-space has
the following inconvenience: strictly speaking
one usually says CAT-TOP-champ to mean a champ
in the big (albeit the use of genuinely
large categories can usually be avoided
\ref{faq5.3}) site of the topology defined by TOP
on CAT, together with the conditions that there
is an atlas and the diagonal is representable, 
a.k.a. the champ is the 2-sheafification in TOP
of a groupoid in CAT, and just to add confusion
CAT-TOP often gets replaced by the person or persons
who first proposed the definition. Thus, for 
example, if CAT is affine schemes, 
which, rather irritatingly, gives
a slightly different theory than CAT=schemes,
and TOP is
\'etale then one says Deligne-Mumford champ, 
or, more correctly, algebraic Deligne-Mumford champ. The irritating 
difference between affine schemes
and schemes in this example disappears under
the hypothesis that the diagonal is affine,
which although a much weaker hypothesis than separated
doesn't have sense in general. We will, however,
usually have plenty of separation hypothesis to
accommodate such border line phenomenon, and so
we adopt the convention that CAT Deligne-Mumford
champ (or just Deligne-Mumford champ, and even
just champ if the context is clear)  means 
CAT-\'etale ({\it i.e.} locally isomorphic)-champ.
Usually CAT is separated topological spaces or 
separated schemes, \ref{faq5.3}, and in
the rare cases where TOP is not \'etale
it is explictly noted.
\subsubsection{Any other Frenchisms that we
should know about ?}
\label{faq5.2}
Just a
couple. 
The first is separated rather than Hausdorff.
The other is: if $j,k$ belong to an inverse (resp. direct)
system then we say that the system is co-filtered on
the right (respectively filtered on the left) if there 
is a $i$ such that $i\ra j$, and $i\ra k$ (respectively
$j\ra i$ and $k\ra i$). I believe that the English for
these conditions is ``directed'', and I've tried to
employ ``directed'' accordingly. Nevertheless, I find this
usage rather confusing since it seems to me that a partially
ordered set is already directed irrespectively of whether
it's right co-filtered, respectively left filtered.
\subsubsection{Which invites the question of any other
notational peculiarities?}
\label{faq5.2bis}
Again, just a couple. The usage of 
``homotopy depth'' rather than ``rectified
homotopy depth'' as covered post \ref{faq4.7}.(c), 
while it may be usefully emphasised that space, whether it be
algebraic, topological, or whatever, unlike champ, 
always means
separated, \ref{faq5.2}, space unless explicitly
stated otherwise, \ref{faq1.8}.
\subsubsection{You're rather non committal about whether
you're supposing the existence of a universe or not.
What exactly is going on ? Are you some sort of
rigorist ?}
\label{faq5.3}
Mitigated laxist is probably closer to the truth. Nevertheless
while positing a universe may be a convenient
extension of ZFC,
it is not a conservative one, so it obviously shouldn't be
done unless it's actually necessary,
and it isn't, so we don't. 
Specifically, one usually supposes that there is
a universe in order to give sense to something
like $\Hom(\mathrm{Sch}\op, \mathrm{\Ens})$, 
resp. $\Hom_2(\mathrm{Sch}\op, \mathrm{\Grpd})$,
so as
to 
be able to talk about sheaves, 
2-sheaves/champs, \ref{faq5.1}, 
on a large site. In an algebraic geometry
context, however,
one never needs anything as large as $\mathrm{Sch}$, {\it i.e.}
Noetherian rings with residue fields finitely generated over
some fixed base inevitably suffices, and this is a set- any Noetherian
domain embeds in its completion at a maximal ideal- while in
a differential geometry context the (small) category with objects
finite disjoint unions of
$\br^n$'s, 
$n$ not necessarily fixed,
and differentiable maps between them plays a similar
role. Plausibly, however, in the general topology, or non-Noetherian
algebraic geometry, setting, where there are non-void results,
\ref{fact:395}-\ref{fact:396}, \ref{ProFF}, one may need to be
a little more ad-hoc to avoid supposing the existence of a 
universe, {\it e.g.} for the homotopy types of such spaces
or schemes 2-sheaves in the small site of the space, resp.
\'etale site of the scheme is sufficient for all the necessary
2-categories to be defined. Such definitions, {\it i.e.} 
a class of
1, resp. 2,
functors from a small category to a large 1, resp. 2, category
all have perfect sense in NBG set theory as does 
``every essentially surjective fully faithful functor
is an equivalence''
for large categories, and its (large) 2-category variant, \cite[1.5.13]{tom}.
This being at worst (and it's almost certainly overkill) 
what we have to worry about, the quickest way to 
proceed while respecting the maxim
``numquam ponenda est pluralitas sine necessitate''\footnote{Which
could reasonably be applied to the great set theory hoax in general},
\cite[i.27.2.K]{ock}, is to work in NBG set theory.
This is a conservative extension of ZFC, 
(with a helpful 
\htmladdnormallink{wiki}{https://en.wikipedia.org/wiki/Von\_Neumann\%E2\%80\%93Bernays\%E2\%80\%93G\%C3\%B6del\_set\_theory\#Category\_theory})
so it's just
an exercise for the reader who wishes to apply a given
theorem without leaving ZFC, to check that the desired
theorem has sense in ZFC. 
%
%
\subsubsection{It does seem pretty convincing that
the language/calculus of weak 2-categories (with all 2-cells
invertible) is the language/calculus of homotopy 2-types,
but 2 isn't such a big number, are there any more intermediary
steps that I can take before jumping to higher-categories ?}
\label{faq5.5}
Not really. Many of the diagram chases here are already
long and tedious, and
in the higher case one needs to organise them by way 
of the formalism of simplicial sets,
\cite{lurie}, 
\cite{simpson} to make them manageable.
This said, the definition of a 3-space is, in a way that takes
into account \ref{faq8}, absolutely unambiguous, to wit one
starts with a (weak) 2-category
\begin{equation}\label{feq51}
X_2 \rras X_1 \rras X_0
\end{equation}
in which all 1 and 2 cells are in invertible, be it
in affine varieties, disjoint unions of $\br^n$'s, 
\ref{faq5.3}, or whatever
one's favourite poison is and 3-sheafifies it in the obvious
(hint: do the case $\rB_{\rB_Z}$ first, {\it i.e.}
$X_1$, $X_0$ points, $X_2$ an abelian group, {\it cf.}
\ref{fact:eg1} and \cite[IV.3.5.1]{giraud} albeit taking account
of hypercovers \ref{fact:384}) way, which could be
worked through with the explicit goal of
bumping the Lefschetz theorem up to $\pi_3$ of 3-spaces, \ref{faq4.8}.
\subsubsection{I have a C.N.R.S. position and would like
to meet a jeune femme, do you have any theorems that might help ?} 
\label{faq5.6}
No, and I strongly suspect that you have the wrong F.A.Q..
This sort of thing was covered by Greg McShane in his
agony columns for ``Math Actuelle'',  
of which his most pertinent answer was: ``Try a normal
girl, hint: start at the FNAC
\footnote{A French audio-visual chain}, they have their names
on a badge, it gives you a start''. 
Unfortunately, however, Mr. Google is only showing
results for the 
somewhat better known magazine ``Femme Actuelle'', and, 
having asked Greg the current status of his own publication, I never got a
reply.
\subsubsection{So you're not thanking your fellow
Glaswegian
for his contribution to the FAQ ?}
\label{faq5.7}
Actually he's from Wishaw, although I would like to thank
J.-B. Bost, \ref{faq5.1}, A. Epstein, \ref{faq5.3},
J.-V. Periera, \ref{faq1.8}, and Y. Rieck, \ref{faq1.7}.
In a less specific way the motivation for this undertaking
came from inflicting my interest in conformal homotopy on
others during my last couple of visits to N.Y.U., for which
I'd like to thank F.A. Bogomolov, S. Cappell, and M. Gromov.

\numberwithin{equation}{section}
\renewcommand{\thesection}{\Roman{section}}
\renewcommand{\thesubsection}{\Roman{section}.\alph{subsection}}
\newpage

\section{Paths, loops and spheres}\label{S:I}

\subsection{Mapping spaces- local theory}\label{SS:I.1}

Let $\cX$ be a (not necessarily separated) topological champ, and
$a:Y\ra \cX$ a map from 
a separated locally compact CW-complex, or, more
generally, a space which is: 
\begin{equation}\label{eq:path0}
\begin{split}
&\text{locally compact; paracompact; separated; and every
cover, $V_\a$, admits a shrinking}\\
&\text{such that all intersections $V_{\a\b}:=V_\a\cap V_\b$
are connected,}
\end{split}
\end{equation} 
where $\a=\b$ is allowed, and $V'_\a\hookrightarrow V_\a$
is a shrinking if the closure of the former is contained
in the latter.
The important point is that the graph:
\begin{equation}\label{eq:path1}
\begin{CD}
Y@>{\G_a}>> Y\times \cX \\
@. @VV{p}V \\
@. Y
\end{CD}
\end{equation}
fails, in general, to
be an embedding. More precisely if
$\pi:U\ra \cX$ is an \'etale atlas, and
we form the fibre product:
\begin{equation}\label{eq:path2}
\begin{CD}
Z@>{\G_a}>> Y\times U \\
@VV{\rho}V @VV{\mathrm{Id}\times\pi}V \\
Y @>{\G_a}>> Y\times \cX
\end{CD}
\end{equation}
then the upper horizontal arrow
need be no better than net if 
$\cX$ isn't a space. The atlas,
$\pi$,
is, however, \'etale, so, by (\ref{eq:path0}), 
as a cover of $Y$, $\rho$ may
be refined to a locally finite
cover $V=\coprod_{\a} V_\a \ra Y$
by 
relatively compact
open embeddings $V_\a\hookrightarrow Y$
such that the function 
\begin{equation}\label{CorrectionI.1.1}
n_V(y) :=\sharp \{V_\b\ni y\}
\end{equation}
restricted to any $V_\a$ 
is bounded.
Plainly the classifying space of the groupoid 
$V\times_Y V\rras Y$ is just $Y$, and
since everything is separated the
top row of (\ref{eq:path2}) affords
a family of embeddings,
$\G_\a:V_\a\hookrightarrow Y\times U$.
Consequently, if for $R:=U\times_{\cX} U\rras U$
the groupoid implied by the atlas $U$,
and 
$V\times_Y V\rras V$ identified with 
$\coprod_{\a\b} V_{\a\b}:=V_\a\cap V_\b$,
we have 
a diagram of pairs of embeddings,
\begin{equation}\label{eq:path3}
\begin{CD}
(V_{\a\b}\hookrightarrow V_{\a\b}\times R) @>>{t}> (\G_\b:V_\b\hookrightarrow V_\b\times U)\\
@VV{s}V @.\\
(\G_\a:V_\a\hookrightarrow V_\a\times U) 
\end{CD}
\end{equation}
The source, $s$, and sink, $t$, 
are \'etale, so, in a way
depending only on $\cX$, we may
insist a priori that the cover
$U$ is sufficiently fine to guarantee
that $R$ is a topological disjoint
union $\coprod_i R_i$ such that
whether the source or the sink
restricted to each $R_i$ is a
homeomorphism onto an open subset
of $U$. The maps \eqref{eq:path3} arise, however,
from the continuous functor
\begin{equation}\label{Combinatorics1}
\coprod_{\a\ts \b}V_{\a\b}= V\ts_Y V \ra R=\coprod_i R_i
\end{equation}
so by our combinatorial hypothesis, \eqref{eq:path0}, 
for sufficently fine covers $V\ra Y$:
each $V_{\a\b}$ is contained in a unique $R_i$ 
which we denote by $R_{\a\b}$, and:
\begin{fact}\label{claim:path1}
Without loss of generality there is
an open neighbourhood
$R'_{\a\b}=V_{\a\b}\times R_{\a\b}$
of the graph of $V_{\a\b}$ 
in \eqref{eq:path3} such that the source
and sink of op. cit. restrict
to a homeomorphism of $R'_{\a\b}$
with its image, and, better still:
if the diagonal $\cX\xrightarrow{\D} \cX\ts\cX$
is universally closed, the restriction
of $s\ts t$ in \eqref{eq:path3} to $R'_{\a\b}$
is universally closed.
\end{fact}

Applying this we obtain  sets of
arrows,
\begin{equation}\label{eq:path4}
V\times_Y V = \coprod_{\a\b}V_{\a\b}\subseteq R':=\coprod_{\a\b} 
R'_{\a\b} \subseteq Y\times R
\end{equation}
where, for convenience, we can use the
identity map in $R$ rather than \ref{claim:path1}
if $\a=\b$. In addition, 
we have a projection $j:R'\ra V\times_Y V$, and, while there is no difficulty
in arranging that the inverse in $Y\times R$
leaves $R'$ invariant, $R'$ will, in all
likelihood, fail to be a groupoid. We can
remedy this by introducing $R''\subseteq R'$
defined by way of: $f\in R''$ iff there is
an open neighbourhood $F\ni f$ in $R'$ such
that for $f'\in F$
\begin{equation}\label{eq:path5}
\begin{split}
&\text{every arrow},\,
k\in V\times_Y V, \, \text{which itself, or 
its inverse, is compossible, be it on the
left,} \\ &\text{or the right, with}\, j(f'),\,  
\text{can be lifted to- a necessarily unique by \ref{claim:path1}- arrow in,}
\, R'
\end{split}
\end{equation}
The main technical device for constructing
mapping spaces is:
\begin{lem}\label{claim:path2} Notations as above,
then the following hold:

(a) The product groupoid structure on
$V\ts_Y V \times R$ induces an
equivalence relation $R''\rras W''$
with \'etale source and sink, for $W''$
the set of unit arrows 
of the product in $R''$. 

(b) For any shrinking $V'\ra V$, the resulting
groupoid, again denoted $R''\rras W''$, contains
the graph (albeit for $V'$ rather than $V$) of the functor 
\eqref{Combinatorics1}.

(c) If the
diagonal $\cX\xrightarrow{\D}\cX\ts\cX$ is
universally closed, or, 
much weaker $s\ts t$ restricted to any
$R_i$ in \eqref{Combinatorics1} is closed,
there is an open neighbourhood $W$ of the
graph of $V'$ in $V\ts U$ such that the
sliced groupoid,  
$R''\rras W$, is separated.
\end{lem}
\begin{proof}
Recalling that one can define a groupoid
without reference to objects- {\it i.e.}
use identity arrows instead- the verification
that $R''\rras W$ satisfies the groupoid
axioms is routine. 
By definition $R''$ is open in $R'$,
so $(s,t):R''\ra W$ is \'etale.
The map, 
\begin{equation}\label{ConS1}
s\times t: R'\longrightarrow \coprod_{\a\b} (V_\a\times U)\times (V_\b\times U)
\end{equation}
is already injective, so, a fortiori, 
$R''\rras W$ is an equivalence relation,
which does (a).
As to (b): if $R'_r$
is the set of arrows with the composition property
(\ref{eq:path5}) on the right, 
albeit just for $k$ of {\it op. cit.},
then 
$R'_r\supseteq V'\times_Y V'$, so if it contained
an open neighbourhood of the same, then
by symmetry the same holds  whether
for the left composition property, and/or
composition with inverses, and by intersecting
these four possibilities, every point of $V'\ts_Y V'$
would satisfy
the definition of $R''$. As such consider the
set of compossible arrows, $C_{\a\b}\subset C$, defined by the
fibre squares,
\begin{equation}\label{eq:path6}
\begin{CD}
R' @<<< C @<<< C_{\a\b}\\
@VV{s}V @VV{\s}V @VV{\s_{\a\b}}V \\
V'\times U @<{t}<< R' @<<< R'_{\a\b}
\end{CD}
\end{equation}
By base change of the known arrows,
all the new arrows are \'etale. The
set $C_{\a\b}$ is contained in,
\begin{equation}\label{ConSistent}
\coprod_\g R'_{\a\b}\times R'_{\b\g}
\end{equation}
so the cardinality, $m(f)$- 
which is lower semi-continuous-
of the fibre of $\s$ over an arrow $f$
is at most 
$n_{V'}(t(j(f)))$, and the set where this
maximum
is realised is $R'_r$.
Now suppose that $R'_r$ isn't open
at the image of some $f\in V'\ts_Y V'$,
then there is a net of arrows,
$f_i\ra f$ such that
\begin{equation}\label{PathCorrect} 
n_{V'}(t(j(f_i)))> m(f_i)
\end{equation}  
Thus 
$n_i:=n_{V'}(t(j(f_i))) > n:=n_{V'}(t(j(f)))$,
and the possible open sets $V_\g\ni t(j(f_i))$ are,
for $f$ fixed, finite independent of $i$. 
Sub-sequencing the net, we obtain that
the set of such open sets, and whence $N=n_i$, is constant.
We have, therefore, $N$-nets of
arrows, $g_i^p$ in $V'\ts_Y V'$ each with source
the sink of $j(f_i)$
in the same $V'_\b$ and repective sinks in distinct
$V'_{\g^p}$, $1\leq p\leq N$. By construction the 
sources of the $g_i^p$ are convergent, while by
the hypothesis, \eqref{eq:path0}, of local compactness
all limit points of their sinks belong to the opens,
$V_{\g^p}$, in the larger cover which we have shrunk. 
The space, $Y$, is, however, separated, so any
net of arrows whose source and sink have a limit,
has itself a limit. Consequently, in the larger
cover $V$, $n_{V}(t(j(f)))$ is at least $N$, and, similarly,
the set valued function of which it is the cardinality
takes at least all of the values $\g^p$. As such
if we replace $V'$ by $V$ in \eqref{eq:path6}, then
in the resulting \eqref{ConSistent}, all the $\g^p$
occur in the fibre over $f$, which, in turn extends
to a neighbourhood, so on returning to $V'$, 
all $m(f_i)\geq N$ contradicting \eqref{PathCorrect}. 
This not only proves (b), but a little more, {\it viz:}
if $Z$ is the closed complement of $R''$ in $R'$,
both understood for the larger cover $V$, then
the closure, $Y$, of $s(Z)\cup t(Z)$ in $V\ts U$
doesn't meet the graph of the shrunk cover $V'$.
However, as noted in \ref{claim:path1}, $s\ts t$
resticted to $R'$ is closed, so we can take $W$
to be the complement of $Y$ in item (c). 
\end{proof}
At the level of spaces \ref{claim:path2},
under hypothesis (c),  is equivalent
to the following commutative diagram,
\begin{equation}\label{eq:path7}
\begin{CD}
\cX @<<{q_a}< J_a:= [W/R''] @<<{\G}< Y\\
@| @V{p}VV @| \\
\cX @<<{a}< Y @=Y
\end{CD}
\end{equation} 
where $J_a$ is a separated space locally isomorphic
to $V_\a\times U_\a$ for $V_\a\hookrightarrow Y$
open, with $U_\a\ra \cX$ an \'etale neighbourhood
of $a(V_\a)$; and, we deliberately confuse
notation with that of the graph, (\ref{eq:path1}), 
since $J_a$ is just an open neighbourhood of the
graph if $\cX$ were a space. This suggests,
\begin{fact}\label{fact:path3}
Suppose $Y$ is compact; $\cX$ satisfies the
weak separation condition \ref{claim:path2}.(c);
and
let $(S_a,\s:=\G)$ be the pointed space of sections of $p$ in
(\ref{eq:path7})- viewed as a subspace of
$\mathrm{Hom}(Y,J_a)$ in the compact open topology-
with $A:Y\times S\ra \cX$ the implied  space of
morphisms,
then for any pointed germ of deformation, $B:Y\times (T,*)\ra \cX$,
and natural transformation $\eta_*$ such that,
\begin{equation}\label{eq:path8}
\xymatrix{
 Y
   \ar@/^2pc/[rr]_{\quad}^{B_*}="1"
   \ar@/_2pc/[rr]_{A_\s=a}="2"
&& \cX
 \ar@{}"1";"2"|(.2){\,}="7"
  \ar@{}"1";"2"|(.8){\,}="8"
 \ar@{=>}"7" ;"8"^{\eta_*}
} 
\end{equation}
2-commutes there is an open neighbourhood $T'\ni *$, and
a map $b:T'\ra S$ such that,
\begin{equation}\label{eq:path9}
 \xy
 (-18,0)*+{Y\times T'}="L";
 (18,0)*+{\cX}="R";
 (0,16)*+{Y\times S}="T";
    {\ar^{\mathrm{id}\times b} "L";"T"};
    {\ar^{A} "T";"R"};
    {\ar_{B} "L";"R"};
    {\ar@{=>}^{\eta} (0,2);(0,12)}
 \endxy
\end{equation}
commutes for a unique natural transformation $\eta$,
while restricting to (\ref{eq:path8}) over the
base point.
\end{fact}
\begin{proof}
In the first instance, suppose that $B_*=a$.
Now, we're not supposing that $T$ is locally
compact around $*$ so the conditions \eqref{eq:path0}
are not satisfied. Nevertheless, because we're
only interested in the germ of $T$ around $*$,
we may shrink- a precision that
will subsequently be omitted- $T\ni *$ as
necessary, and, 
replacing $Y$ by $Y\ts T$, we can form
$J_B$ in \eqref{eq:path7} in the same way
as we formed $J_a$. More precisely,
at the first stage (\ref{eq:path2}) we have an
\'etale covering of the open neighbourhood
$Y\times T$ of $Y\times *$, so, 
since $Y$ is compact, we may suppose that our covering
of $Y\times T$ is just $V\times T$,
and \ref{claim:path1} holds with the
same $R_{\a\b}$.
The key point where local compactness is
used in the proof of \ref{claim:path1} is
to show that arrows where the composition
property \eqref{eq:path5} fails cannot
accumulate on the graph of the smaller
cover $V'_\a$ which may well fail for
$V'_\a\ts T$, but since we can take $T$
arbitrarily small we only need this for
$V'_\a\ts *$, which holds with the same proof,
so that:
\begin{lem}\label{lem:path1}
For $Y$ compact, 
the formation of $J_a$ in (\ref{eq:path7})
commutes with sufficiently small base change,
{\it i.e.} if $B_*=a$ in 
\ref{fact:path3} and $T'\ni *$ is sufficiently
small there is a commutative diagram,
\begin{equation}\label{eq:path10}
\begin{CD}
\cX @<{q_B}<< J_B @>{p_{T'}}>> Y\times T' \\
@| @VVV @VVV \\
\cX @<{q_a}<< J_a @>{p}>> Y
\end{CD}
\end{equation}
with the second square fibred, and satisfying
the obvious further commutativity on adjoining
the sections $\s=\s_a$ and $\s_B$.
\end{lem}
Applying this to the case in point, then,
at least as far as the existence of $b$
goes, we can suppose that $T=S_{B_*}$
around $\s_{B_*}$. Consequently, we are
largely done if we can establish,
\begin{flem}
Irrespective of the compactness of $Y$,
a natural transformation $\eta_*:a\Rightarrow c$
between 1-morphisms $a,c:Y\ra \cX$ affords
an isomorphism, $i_{\eta}$, between neighbourhoods-
notationally identified with the ambient space-
of $\s_a$ and $\s_c$ in $J_a$, $J_c$ respectively
such that, 
\begin{equation}\label{eq:path11}
 \xy
 (-18,0)*+{J_c}="L";
 (18,0)*+{\cX}="R";
 (0,16)*+{J_a}="T";
    {\ar^{i_{\eta}} "L";"T"};
    {\ar^{q_a} "T";"R"};
    {\ar_{q_c} "L";"R"};
    {\ar@{=>}^{\eta} (0,2);(0,12)}
 \endxy
\end{equation}
commutes
for some unique natural transformation
$\eta$ restricting to $\eta_*$.
\end{flem}
The graph, $\G_\eta$ of $\eta_*$ affords
a commutative diagram,
\begin{equation}\label{eq:path12}
 \xy
 (-20,0)*+{V\times U}="L";
 (20,0)*+{V\times U}="R";
 (0,16)*+{V}="T";
 (0,0)*+{V\times R}="M";
    {\ar_{\G_a} "T";"L"};
    {\ar^{\G_c} "T";"R"};
    {\ar^{s} "M";"L"};
    {\ar_{t} "M";"R"};
    {\ar^{\G_{\eta}} "T";"M"}
 \endxy
\end{equation}
while the cover $W_a$, respectively $W_c$, of \ref{claim:path2} 
of $J_a$, respectively $J_c$, may
be identified to an open neighbourhood of 
$\G_a$, respectively $\G_c$. Since $s$ and $t$ are
\'etale if we
shrink everything as necessary, there is, therefore,
a unique way to extend (\ref{eq:path12}) to a commutative diagram,
\begin{equation}\label{eq:path13}
 \xy
 (-34,0)*+{V\times U\supseteq W_a}="L";
 (34,0)*+{W_c\subseteq V\times U}="R";
 (0,16)*+{W_c}="T";
 (0,0)*+{V\times R}="M";
    {\ar@{=} "T";"R"};
    {\ar_{s} "M";"L"};
    {\ar^{t} "M";"R"};
    {\ar^{\G_{\eta}} "T";"M"}
 \endxy
\end{equation}
By unicity, the composition $s\G_\eta$ glues
to an isomorphism $i_\eta$ between $J_c$ and
$J_a$, allowing us to fill the left diagonal
arrow in (\ref{eq:path13}). Composing this with
the inclusion of $R''_a\rras W_a$ in $R$,
the unicity of  $\G_\eta$ in (\ref{eq:path13}),
now gives a natural transformation between this
and the inclusion of $R''_c\rras W_c$ in $R$.

Applying this with $c=B_*$ yields 
\ref{fact:path3} up to the possible 
lacuna of the uniqueness of $\eta$
if $T$ were smaller than the full
space of sections $S_{B_*}$. To resolve
this, however, we need only consider
the base change variant of (\ref{eq:path12}),
to wit:
\begin{equation}\label{eq:path14}
\begin{CD}
V\times T @<<< V\times * @>>> V\times T\\
@V_{\G_{A(\mathrm{id}\times b)}}VV @V{\G_\eta}VV @V{\G_B}VV \\
V\times T\times U@<s<< V\times T \times R @>t>> V\times T\times U
\end{CD}
\end{equation}
so, again, there is a unique extension
of $\G_\eta$ to a sufficiently small
neighbourhood $V\times T'\supset V\times *$
such that the diagram commutes.
\end{proof}
Let us add some clarification by way of,
\begin{rmk}\label{rmk:path}
The hypothesis of compactness of $Y$ in \ref{fact:path3}
is only necessary if one wants the universality
diagram (\ref{eq:path9}) to hold in an neighbourhood
in $Y\times T$ of the form $Y\times T'$ for $T'\ni *$.
If one is prepared to weaken this to,
\begin{equation}\label{eq:path15}
\text{The universal property (\ref{eq:path9}) holds
in an open neighbourhood of}\, Y\times *\,\, \text{in}\,
Y\times T
\end{equation}
then there is a non-compact alternative with exactly
the same proof.

Some further insight arises from considering
the case where $\cX$ is differentiable. As such,
we have the zero section of the tangent bundle
$0:Y\xleftarrow{\rightarrow}a^* T_{\cX}:p$, 
and, up to appropriate shrinking in a neighbourhood
of $\s$, $J_a$ is 
diffeomorphic to an open neighbourhood of
$0$ in $a^*T_{\cX}$. 

This infinitesimal alternative is particularly
well adapted- \cite{Mp}- to the analogue of \ref{fact:path3}
in algebraic geometry, albeit one has to
replace $a^*T_{\cX}$ by Grothendieck's jets
$a^*\cP_{\cX}$, and, at least in characteristic
$p>0$ compactness becomes indispensable 
because of the wild nature of the \'etale
topology otherwise, cf. \cite[Expos\'e X, 1.9-10]{sga1}.
In a sense, however, the topological case is more
difficult since already topological spaces aren't
Cartesian closed, and even the existence of an
exponential object requires local compactness, which,
in turn needn't be locally compacted, {\it i.e.}
exactly the kind of difficulties encountered in  
the proofs of
\ref{claim:path2} and \ref{fact:path3}.
\end{rmk}

\subsection{Mapping spaces- global theory}\label{SS:I.2}

Plainly we continue to suppose that the space
$Y$ satisfies (\ref{eq:path0}), we do not,
in the first instance, make any global compactness assumption, 
we do, however, suppose throughout this section
suppose that $\cX$ is separated
({\it i.e.} $\cX\xrightarrow{\D}\cX\ts\cX$ proper) and
we choose a
sets worth, $A$, of representatives of
1-morphisms $a:Y\ra \cX$
-every such map can be represented by the
set of
functors of the form \eqref{Combinatorics1}
for some open cover $V\ra Y$.

Now for $a,b\in A$ and $\cX_Y=Y\times X$ consider $J_{ab}$ defined via the Cartesian square,
\begin{equation}\label{eq:global2}
\begin{CD}
J_a @<<< J_{ab}\\
@V{p_a\times q_a}VV     @VVV \\
\cX_Y @<<{p_b\times q_b}< J_b
\end{CD}
\end{equation}
By construction, $p_a\times q_a$ is \'etale, so
this is just a presentation of $\cX_Y$ as a
(necessarily separated) groupoid,
\begin{equation}\label{eq:global3}
\coprod_{ab} J_{ab} \rras J:=(\coprod_a J_a)
\end{equation}
in spaces over $Y$.
In particular  sections, $S_{ab}$, of $p_{ab}$ project
to section of $p_a$, respectively $p_b$,
yielding a groupoid,
\begin{equation}\label{eq:global1}
(s,t): \Sigma:= \coprod_{ab} S_{ab}\rras S:=\coprod_{a} S_{a}
\end{equation}
and we assert,
\begin{claim}\label{claim:global1}
If $Y$ is compact, then the groupoid (\ref{eq:global1}) is separated and \'etale.
\end{claim}
\begin{proof}
Notwithstanding that (\ref{eq:global3}) is separated
and \'etale,
we should be a little cautious, since the
mapping functor $T\mpo \Hom(Y\ts T,X)$ doesn't behave
particularly well with respect to either \'etale
or proper maps. We are, however, considering spaces
of sections, so for $\s_{ab}$ a section of $p_{ab}$ 
projecting to a section $\s_a$ of $p_a$ we have a
diagram,
\begin{equation}\label{eq:global5}
\begin{CD}
Y@>{\s_{ab}}>> J_{ab} \\
@| @VV{s}V \\
Y@>{\s_{a}}>> J_{a}
\end{CD}
\end{equation}
where the horizontal arrows are embeddings, and
$s$ is \'etale. Since $Y$ is compact, and everything
is separated, $s$ extends to a homeomorphism from
a neighbourhood of $\s_{ab}(Y)$ onto its image,
so the source and sink in (\ref{eq:global1})
are \'etale. Finally, let $T_a$, $T_b$ be compacts
in $S_a$, and $S_b$ respectively, with 
$f_n=(f_{an}, f_{bn})$ a net of arrows in $T_{ab}:=T_{a\,s}\times_{t} T_{b}$
converging to $f=(f_a,f_b)\in T_a\times T_b$.
The functions $f_{an}$, $f_{bn}$ take values
in compacts $K_a$, $K_b$ independent of $n$,
so by the separation of (\ref{eq:global3}), the
functions $f_n$ take values in the compact
$K_{ab}:= K_{a\,s}\times_{t} K_{b}$. This latter
set admits a finite cover $W=\coprod_i W_i$ such
that the source and sink of (\ref{eq:global3})
restricted to each $W_i$ are homeomorphisms
with the image. Refining the net as necessary,
we can suppose that there is a relatively compact
open cover $V=\coprod_i V_i\ra Y$ such that
every $f_{an}$, $f_{bn}$ restricted to $V_i$
takes values in compacts $K_{ai}$, $K_{bi}$ 
of $s(W_i)$, respectively $t(W_i)$, independently
of $n$. Consequently the $f_n\mid_{V_i}$
take values in the compact  $K_{abi}=K_{ai\,s}\times_{t} K_{bi}\subseteq K_{ab}$,
where, by construction, the source and sink
restricted to $K_{abi}$ are homeomorphisms
with the image. As such, for each $i$, the net $f_n\mid_{V_i}$
limits on $g_i:V_i\ra K_{abi} \subseteq T_{ab}$.
Since the $f_n$ are continuous, and $J_{ab}$
is separated, the $g_i$ glue to a function
from $Y$ to $T_{ab}$. 
\end{proof} 
Denoting, the 
projection $Y\times [S/\Sigma]\ra\cX$ by $Q$,
we have the following universal property:
\begin{factdef}\label{factdef:global} Let $Y$ be compact, then the 
separated classifying champ $\Hom(Y,\cX)$ defined
by (\ref{eq:global2}) has the following universal
property: for every morphism $F:Y\times T\ra\cX$,
with $T$ a space, there is a morphism $G:T\ra\Hom(Y,\cX)$,
and a natural transformation $\eta$ such that,
\begin{equation}\label{eq:global6}
 \xy
 (-18,0)*+{Y\times T}="L";
 (18,0)*+{\cX_Y}="R";
 (0,16)*+{Y\times\Hom(Y,\cX)}="T";
    {\ar^{\mathrm{id}\times G} "L";"T"};
    {\ar^{\mathrm{id}\times Q } "T";"R"};
    {\ar_{F} "L";"R"};
    {\ar@{=>}^{\eta} (0,2);(0,12)}
 \endxy
\end{equation}
2-commutes. In addition if $(G',\eta')$ is any other
pair satisfying (\ref{eq:global6}) then,
\begin{flalign}\label{eq:globalextra}
\text{there is
a unique natural transformation}\,\, \xi:G \Rightarrow G'
\text{such that}\, \eta'=Q(\mathrm{id}\times\xi)\eta.
\end{flalign}
\end{factdef}
\begin{proof} 
By definition, for $\tau\in T$ there is an $a(\tau)\in A$,
and a natural transformation, $\zeta_{\tau0}:\tau\Rightarrow a(\tau)$. 
By (\ref{eq:path9}), for $W_\tau\ni \tau$ a sufficiently small
\'etale neighbourhood this extends uniquely to a 
commutative diagram,
\begin{equation}\label{eq:global8}
 \xy
 (-18,0)*+{Y\times W_\tau}="L";
 (18,0)*+{Y\times\cX}="R";
 (0,16)*+{Y\times S_{a(\tau)}}="T";
 (36,16)*+{J_{a(\tau)}}="E";  
   {\ar^{\mathrm{id}\times G_\tau} "L";"T"};
    {\ar^{Q} "T";"R"};
    {\ar_{F\mid_{Y\times W_\tau}} "L";"R"};
    {\ar@{=>}^{\zeta_\tau} (0,2);(0,12)};
{\ar^{\mathrm{evaluation}} "T";"E"};
{\ar^{p_a\times q_a} "E";"R"}
 \endxy
\end{equation}
restricting to the initial data over the closed point.
Put $e_\tau$ to be the composite of the leftmost and top
arrows, then the diagram,
\begin{equation}\label{eq:global9}
\xy
(-22,10)*+{Y\times W\, (:=\coprod_\tau W_\tau)} ="L";
(22,10)*+{F^*J \, (:=\coprod_a J_a)}  ="R";
(0,0)*+{Y\times T}="B";
{\ar "L";"B"};
{\ar "R";"B"};
{\ar^{\coprod{e_\tau}} "L";"R"}
\endxy
\end{equation}
is just the refinement of one \'etale cover by another,
so there is a functor of groupoids over $Y$,
\begin{equation}\label{eq:global10}
e: (Y\times (W\times_{T} W) \rras Y\times W) \longrightarrow
(J \times_{\cX_Y} J \rras J)
\end{equation} 
together with a natural transformation $\eta:F\Rightarrow e$,
which may have little, or nothing, to do with 
the $\z_\tau$ in (\ref{eq:global8}). This establishes
(\ref{eq:global6}) a fortiori. 

Now suppose $(G',\eta')$ is another pair satisfying
(\ref{eq:global6}), 
and let $W':= \coprod W'_\tau\ra T$ be a cover
sub-ordinate to $(G')^* S\ra T$ and $W\ra T$,
then for some assignment $a':T\ra A$ we have a
diagram,
\begin{equation}\label{eq:global11}
 \xy
 (-20,0)*+{Y\times W'_\tau}="L";
 (56,0)*+{Y\times\cX}="R";
 (0,10)*+{Y\times S_{a(\tau)}}="T";
 (36,10)*+{J_{a(\tau)}}="E";  
 (0,-10)*+{Y\times S_{a'(\tau)}}="A";
 (36,-10)*+{J_{a'(\tau)}}="B";  
{\ar_{\mathrm{evaluation}} "A";"B"};
{\ar_{p_a'\times q_a'} "B";"R"};
{\ar_{\mathrm{id}\times G'_\tau} "L";"A"};
   {\ar^{\mathrm{id}\times G_\tau} "L";"T"};
    {\ar@{=>}^{\zeta_\tau=\eta'_\tau\eta^{-1}_\tau} (18,-7);(18,7)};
{\ar^{\mathrm{evaluation}} "T";"E"};
{\ar^{p_a\times q_a} "E";"R"}
 \endxy
\end{equation}
The natural transformation $\z$ is, in particular,
a map $\z_\tau:Y\times W'_\tau\ra J_{a(\tau)a'(\tau)}$ so,
tautologically, we get a map $\xi_{\tau}: W'_\tau\ra S_{a(\tau)a'(\tau)}$
such that,
\begin{equation}\label{eq:global12}
 \xy
 (-30,0)*+{Y\times W'_\tau}="L";
(0,10)*+{Y\times S_{a(\tau)}}="T";
(0,-10)*+{Y\times S_{a'(\tau)}}="A";
(30,0)*+{Y\times  S_{a(\tau)a'(\tau)}  }="R";
{\ar_{\mathrm{id}\times G'_\tau} "L";"A"};
   {\ar^{\mathrm{id}\times G_\tau} "L";"T"};
{\ar_{s} "R";"T"};
{\ar^{t} "R";"A"};
{\ar^{\mathrm{id}\times \xi_\tau} "L";"R"}
 \endxy
\end{equation}
commutes. Since $\z$ is a natural transformation,
so is $\xi=\coprod_{\tau} \xi_\tau$, while the
commutativity of (\ref{eq:global12}) is 
equivalent to (\ref{eq:globalextra})- uniqueness
being immediate from the natural injectivity of
arrows in (\ref{eq:global3}) into those of
(\ref{eq:global1}).
\end{proof}
By way of an exercise in the universal property,
\begin{rmk}\label{rmk:global}
For $(G,\eta)$ as in (\ref{eq:global6}), there is
a natural isomorphism between the stabiliser of $G$
and that of $F$.
\end{rmk}

\subsection{Base points}\label{SS:I.3}
Plainly, by a base point of a champ $\cX$
is to be understood a map $*:\rp\ra\cX$.
This will, however, be an embedding iff
$\cX$ is space like around $*$. By way
of notation we make,
\begin{defn}\label{defn:point1}
By $\cX_*$ is to be understood a champ,
or space, with a base point, and 
for $Y$ a compact space satisfying (\ref{eq:path0})
$\Hom_*(Y_*,\cX_*)$- itself pointed in 
the constant map $*$-
is defined, for $\cX$ separated, 
via the fibre square,
\begin{equation}\label{eq:point1}
\begin{CD}
*\times \Hom(Y,\cX) @<<< \Hom_*(Y_*,\cX_*)\\
@V{Q}VV @VVV \\
\cX @<{*}<< \rp
\end{CD}
\end{equation}
\end{defn}
This construction has surprisingly good
properties, but extracting them requires
a little care. In the first place, the fibre
square (\ref{eq:point1}) is, 
\ref{def:Fibre2},
only commutative up to a natural
transformation $\zeta$ between the compositions
of either side of the square, and
we assert:
\begin{claim}\label{fact:point1}
If $Y$ is connected, and $\cX$ is  separated,
then $\Hom_*(Y_*,\cX_*)$ is a separated space,
and the pair $(\Hom_*(Y_*,\cX_*), \zeta)$
has the following universal property:
if $F:Y\times T\ra\cX$ is a map from a space such that
there is a natural transformation  
$\xi:*\Rightarrow F\vert_{*\times T}$ then
there is a unique map $G:T\ra \Hom_*(Y_*,\cX_*)$ along
with
a unique natural transformation $\eta$ such that $\zeta=\eta_*\xi$,
and the following diagram 2-commutes
\begin{equation}\label{eq:pointuni1}
 \xy
 (-18,0)*+{Y\times T}="L";
 (18,0)*+{\cX_Y}="R";
 (0,16)*+{Y\times\Hom_*(Y_*,\cX_*)}="T";
    {\ar^{\mathrm{id}\times G} "L";"T"};
    {\ar^{\mathrm{id}\times Q } "T";"R"};
    {\ar_{F} "L";"R"};
    {\ar@{=>}^{\eta} (0,2);(0,12)}
 \endxy
\end{equation}
\end{claim}
\begin{proof}
The unicity in the asserted universal property implies, 
that every point of $\Hom_*(Y_*,\cX_*)$ has a trivial stabiliser,
cf. \ref{rmk:global}, while
$\Hom_*(Y_*,\cX_*)$ is certainly a separated
champ, so arguing as in \cite[1.1]{km} it's a separated space.
The universal property (\ref{eq:global6}),
and the definition of fibre products afford
a commutative diagram (\ref{eq:pointuni1}),
where necessarily $\z=\eta_*\xi$. 
By 
(\ref{eq:globalextra}) we know that if
$(G',\eta')$ is another such pair then there is a natural
transformation $\g:G\Rightarrow G'$
such that $\eta'=Q(\mathrm{id}\times \g)\eta$,
so $Q(\g_*)$ is identically $1$, and by
\ref{rmk:global} $\g_*$ is also identically $1$.
Since natural transformations of maps
of spaces are trivial, it suffices to prove that if $G=G'$
then $\g=\mathbf{1}$.
Now,
$Y\times\Hom_*(Y_*,\cX_*)$
is represented by a separated \'etale groupoid,
so the  
automorphisms of $\mathrm{id}\times G$ are
a discrete sheaf over $Y\times T$, and whence the
subset, $V$, where $\g=\mathbf{1}$
is open and
closed. By construction, $V$ contains
$*\ts T$, so every 
fibre, $V_t$, is a non-empty
open closed subset of $Y$ containing $*$,
whence $V$ is everything.
\end{proof}
For all the difference that it makes, one 
could suppress $\z$ in \ref{fact:point1}
from the notation. Nevertheless, strictly
speaking there is a mild ambiguity, which
even extends, 
on pointing $\Hom_*(Y_*,\cX_*)$ by the
unique map $\rp\ra \Hom_*(Y_*,\cX_*)$
afforded by $\rp\ra\cX$ via the
universal property (\ref{eq:pointuni1}),
to the universal pointed map implicit in
(\ref{eq:pointuni1}) only being weakly
rather than strictly pointed, albeit that
this can always be remedied by an a priori
change
to a naturally equivalent base point in
$\cX$, which we may often do without further
comment in applying definitions such as:
\begin{defn}\label{defn:point2}
Let $I_0$ be the unit interval pointed
in zero, then we have a path space functor:
\begin{equation}\label{eq:point2}
\cX_* \mpo \rP{\cX_*}:= \Hom_*(I_0,\cX_*) 
\end{equation}
and a loop space functor, $\Omega{\cX_*}$
defined via the  Cartesian square,
\begin{equation}\label{eq:point3}
\begin{CD}
1\times \rP{\cX_*}@<<< \Omega{\cX_*}\\
@V{Q}VV @VVV \\
\cX@<{*}<< \rp
\end{CD}
\end{equation}
where either functor takes values in
pointed spaces, and $\cX$ is supposed
separated.
\end{defn}
The next subtlety is that $\rp\xrightarrow{*}\cX$ need not be
an embedding, so:
\begin{warning}\label{warn:point1}
Neither the top horizontal in (\ref{eq:point1}) nor
that in (\ref{eq:point3}) need be an embedding.
For example if $\cX=\rB_G$ for $G$ a, say
finite, group, then since every torsor
over the interval is trivial,
$\Hom(I,\rB_G)$ has, up to equivalence one point.
By \ref{rmk:global} its stabiliser is $G$,
so in fact $\Hom(I,\rB_G)\xrightarrow{\sim}\rB_G$,
and fails to be contractible for $G\neq 1$, while,
$\rP\rB_{G}\xrightarrow{\sim}\rp$,
and $\Omega\rB_{G}\xrightarrow{\sim} G$.
\end{warning}
While, at first sight, this may be
counter intuitive, it's exactly the
behaviour that one wants, to wit:
\begin{fact}\label{fact:point2}
The path space $\rP{\cX_*}$ is contractible,
and we have a canonical isomorphism: 
\begin{equation}\label{eq:point4}
\Hom_*(\rS^1_*, \cX_*)\xrightarrow{\sim} \Omega{\cX_*}
\end{equation}
\end{fact}  
\begin{proof}
For $I^\lb$ the interval of length $\lb\leq 1$, the series of maps,
\begin{equation}\label{eq:point5}
\lb: I\xrightarrow{x\mpo \lb x} I^\lb \hookrightarrow I,\quad \text{yields}\quad
\Lambda:\rP{\cX_*}\times [0,1] \ra \rP{\cX_*}: a\mpo \vert\lb^* a\vert
\end{equation}
by the universal property (\ref{eq:global6}),
while \ref{fact:point1} applied to $Y$ a point,
shows that $\Lambda_0$ is as expected, {\it i.e.}
the composition:
$\rP{\cX_*} \ra\rp\xrightarrow{*} \rP{\cX_*}$,
proving the obvious contractibility of $\rP{\cX_*}$.

As to the second part, 
to get a map from left to right in (\ref{eq:point4}) is just the
definition of fibre products (of categories)
and universal property (\ref{eq:global6}) applied
to the composition 
$F:I\times \Hom(\rS^1,\cX) \ra \rS^1\times \Hom(\rS^1,\cX)\ra \cX$.
One can, however, usefully note that for $(G,\eta)$ as
in (\ref{eq:global6}) for $Y=I$, but $\z$ as in
\ref{fact:point1} for $Y=\rS^1$, we get natural 
transformations $\eta^{-1}_p\z$ between $*$ and
$Q(p,G)$ which can be quite different at the
end points $p\in\{0,1\}$. This manifests itself
in going the other way round, since by definition,
$\Omega:=\Omega\cX_*$ is a space of maps from
the interval, say: $a:I\times \Omega\ra\cX$
equipped with natural transformations
$\xi_p: *\Rightarrow a\mid_{p\times \Omega}$ 
for $p=0$ or $1$.
We can extract a descent data over $\rS^1\times \Omega$
compatible with how we already
went from left to right:
identify the open interval $I'=(0,1)$ with
$\rS^1\bsh *$ and take a small neighbourhood
$V\supset *\times \Omega$ such that $V\bsh *\times \Omega$
has a topological decomposition $V_0\coprod V_1$,
and we put $\bar{V}_p=V_p\cup  *\times \Omega$.
Since $\cX$ is defined by an \'etale 
groupoid this latter choice can be made
in a way that for some cover $W\ra V$,
each $\xi_p$ extends uniquely to a map, $\tilde{\xi}_p$ of 
$\bar{W}_p:=W\times_V \bar{V}_p$
to the set of arrows of some presentation
$R\rras U$ of $\cX$ with  sink 
$a\mid_{\tilde{p}}$ for $\tilde{p}$ a 
sufficiently small open (half) cover 
of the end points $p\times \Omega$. Consequently, we get a map
$b$ of $W$ to $U$ by taking the source of
$\tilde{\xi}_0\vee_{*\times \Omega} \tilde{\xi}_1$, and- by the unicity of
the extension- a natural transformation
$\tilde{\xi}:=\tilde{\xi}_0\coprod \tilde{\xi}_1$ of 
$b$ restricted to ${W\bsh *\times \Omega}$
with $a$ restricted to ${\tilde{0}\bsh 0\times \Omega}\coprod {\tilde{1}\bsh 1\times \Omega}$.
Since $\tilde{\xi}$ is a natural transformation, the
maps, $a\coprod b: (I'\times \Omega)\coprod W\ra U$,
glue to some $\rS^1\times \Omega\ra \cX$ which
sends $*\times\Omega$
to a natural transformation of $*$,
so we can go from 
right to left in (\ref{eq:point4}) by the universal property
(\ref{eq:pointuni1}). 
\end{proof}
In consequence, even though the non-spatial
maps in \ref{defn:point2} need not be
strictly pointed,
\ref{fact:point2} combined
with \ref{fact:point1} implies,
\begin{cor}\label{cor:point}
The (internal) adjoint
of looping is, as usual, suspension, {\it i.e.} for pointed $T_*$:  
\begin{equation}\label{eq:pointsuspend} 
\Hom_*(T_*,\Omega\cX_*)= \Hom_*(\Si T_*,\cX_*)
\end{equation}
\end{cor}
whence for $p,q\in\bn$:
\begin{equation}\label{eq:point6}
\Omega^q\cX_*=\Hom_*(\rS^1_*,\Omega^{q-1}\cX_*)=\Hom_*(\rS^q_*,\cX_*)\quad
\Hom_*(\rS^p_*,\Omega^{q}\cX_*)= \Hom_*(\rS^{p+q}_*,\cX_*)
\end{equation}
which ($\cX$ separated) shows
all possible definitions of the higher groups are
equivalent, {\it i.e.},
\begin{equation}\label{eq:point7} 
\pi_{p+q}(\cX_*) := \pi_0\left( \Hom_*(\rS^{p+q}_*,\cX_*)\right)=\pi_p(\Omega^{q}\cX_*),\quad
p\in\bn\cup\{0\}, \, q\in\bn
\end{equation}
Where even the $\pi_0$ in (\ref{eq:point7}) is the $\pi_0$ 
of a space by \ref{fact:point1}. 
It, therefore, 
remains to discuss the pointed set, $\pi_0(\cX_*)$ of path connected components, 
pointed in the component of $*$, and the 
related question of concatenation/the
group structure in (\ref{eq:point7}). 
To this end, consider how connectedness
manifests itself
as an equivalence relation on external $\Hom$, to wit:
\begin{equation}\label{eq:point8} 
\begin{split}
\text{for}\, x,y:\rp\ra\cX,\, &\text{define}\,
x\sim  y\, \text{if for 
some natural transformations}\, \xi,\eta, \,\, \text{we have}: 
\\
& x {\build\Rightarrow_{}^{\xi}} \, a(0)\,
 \xrightarrow{a}\, a(1)\, {\build\Rightarrow_{}^{\eta}}\, y,\quad
a:I\ra\cX\, \text{a path}
\end{split}
\end{equation}
while- as in the proof of \ref{fact:point2}-
if we have two paths $a$, $b$ whose sink and source
are related by a natural transform then we can slightly
thicken $a$  about $1$,
and similarly $b$   about $0$ in such
a way that these
thickenings
glue, by a thickening  of
the implied natural transformations. As such
we get a concatenation operation:
\begin{equation}\label{eq:point9} 
x {\build\Rightarrow_{}^{\xi}} \, a(0)\,
 \xrightarrow{a}\, a(1)\, {\build\Rightarrow_{}^{\eta}}\, y\,\, \bigvee\,\,
y {\build\Rightarrow_{}^{\eta'}} \, b(0)\,
 \xrightarrow{a}\, b(1)\, {\build\Rightarrow_{}^{\zeta}}\, z\,\,
\longmapsto\,\,
x {\build\Rightarrow_{}^{\xi'}} \, c(0)\,
 \xrightarrow{a}\, c(1)\, {\build\Rightarrow_{}^{\zeta'}}\, z
\end{equation}
with the not unimportant caveat that $c$ is
only well defined up to natural transformations.
This lack of well definedness of concatenation is neither
here nor there as far as (\ref{eq:point8}) being
an equivalence relation is concerned, but otherwise
it may not be the behaviour that one wants.
Fortunately, \ref{fact:point1} gives an alternative.
Say, for example $a$ belongs to the moduli 
space of paths with fixed end $\z_a:*\Rightarrow a(*)$
the locally constant natural transformation 
appearing in (\ref{eq:pointuni1}), then 
by the universal property (\ref{eq:pointuni1})  we can find some
unique
$ba$ which is naturally equivalent to $c$,
and $\z_a:*\Rightarrow ba(*)$. Consequently,
\begin{fact}\label{fact:point3}
If in the concatenation schema (\ref{eq:point9}),
$x$ is naturally equivalent to $*$, then 
there is a well defined concatenation $ba\in  \rP\cX_*$,
and unique natural transformations,
$\a:(2)^*a\Rightarrow ba\vert_{[0,1/2]}$,
$\b:(2\mathrm{id}-1)^*b\Rightarrow ba\vert_{[1/2,1]}$
such that $\eta'\eta=\a\b^{-1}$.
\end{fact}
\begin{proof} 
Everything is pretty much by definition,
except uniqueness of $\a$ and $\b$, but,
in general these are unique modulo natural
transformations of the glued arrow, so
this follows by (\ref{eq:pointuni1})
\end{proof}
This may seem pedantic,
but it's exactly what's needed
to guarantee, \ref{rmk:fib},
that the path fibration is a
fibration. One should also note
that there is a need to distinguish,
as in the proof of the contractibility
of the path space in (\ref{eq:point5}),
between the base point of $\lb^*a$,
$a\in\rP\cX_*$, $\lb\in [0,1]$ which is a constant,
$\z_a$ say, and that of its moduli $\vert\lb^*a\vert$,
which if $\rP\cX_*\ra \cX_*$ is strictly
pointed, will collapse to the identity as $\lb\ra 0$.
Obviously, there are plenty of variants on this
theme. Consider for example the map,
\begin{equation}\label{point:final}
\mathrm{cat}: \rS^p\ra \rS^p\vee \rS^p
\end{equation}
affording multiplication of spheres,
{\it i.e.} $\s\tau=\mathrm{cat}^*(\s,\tau)$ on
$\Hom(\rS^p,\_ )$, then we have a map,
\begin{equation}\label{point:final1}
\mathrm{cat}^*: \rS^p\times \Hom_*(\rS^p\vee \rS^p ,\cX_*) \ra \cX
\end{equation}
whence by (\ref{eq:pointuni1}), a product $\s\tau\in \Hom_*(\rS^p\vee \rS^p ,\cX_*)$,
albeit $\s\tau$ may now only be uniquely equivalent
to $\mathrm{cat}^*(\s,\tau)$. Fortunately, since
the equivalence is unique, it's functorial, and 
the product so defined remains homotopy associative.

\subsection{Fibrations}\label{SS:I.4}

The finally possibility in (\ref{eq:point7}) has the
advantage that it allows us to quickly extend our
knowledge of homotopy groups of spaces to champs,
modulo some definitions, {\it i.e.}
\begin{defn}\label{defn:fib0}
Let $p:\cE\ra\cB$ be a map of champs, then  
$f:T\ra \cB$ is said to be
weakly liftable 
with ambiguity $\a$,
if there's some $F:T\ra\cE$
together with an equivalence $\a:f\Rightarrow p F$
\end{defn}
Unlike asking for strict liftability, 
weak liftability
is closed under  base change.
Similarly one wants a definition of
fibration that is closed under base change,
which forces:
\begin{defn}\label{defn:fib}
A map $p:\cE\ra\cB$ is said to be a fibration if it
has the following homotopy lifting property for maps from spaces:
if 
for $f:I\times T\ra\cB$, 
$f\vert_{0\times T}$ is liftable to $\tilde{f}_0$ with
ambiguity $\eta$ then $f$ is liftable to some $F$
with ambiguity $\a$ in such a way that there 
exists a natural transformation $\xi:\tilde{f}_0\Rightarrow F_0$
rendering commutative:
\begin{equation}\label{eq:ambig}
 \xy
 (-12,0)*+{p\tilde{f}_0}="L";
 (12,0)*+{pF_0}="R";
 (0,10)*+{f_0}="T";
    {\ar@{=>}_{\eta} "T";"L"};
    {\ar@{=>}^{\a_0 } "T";"R"};
    {\ar@{=>}_{p\xi} "L";"R"}
 \endxy
\end{equation}
\end{defn}
By way of clarification let us make:
\begin{rmk}\label{rmk:fib}
The definition \ref{defn:fib} is the unique variant of
the usual homotopy lifting property that is closed
under base change. 
Consequently if $p:\cE\ra\cB$ is representable, and
$f:T\ra\cB$ is a map from a space, then $f^*p$ is
a fibration in the usual sense, while a small 
diagram chase reveals that if $p:\cE\ra\cB$ is representable,
then it's a fibration iff $f^*p$ is a fibration for
every such $f$. Otherwise, there is something to check,
which in the case of  $\rP\cX_*\ra \cX$ the obstruction
to following the usual proof, {\it i.e.} pull back
the universal map along $I\times I\ra I:(x,y)\mpo x+y/2$,
is overcome by \ref{fact:point3}, 
and we conclude that $\rP\cX_*\ra \cX$ is a
fibration in the sense of \ref{defn:fib}.
\end{rmk}
As such, we don't require the homotopy lifting
property for champs, so, let us make:
\begin{warning}\label{warn:fib}
Supposing 
$\cE$, $\cB$ separated and
$\cB$ path connected, it's plausible that
the fibres of $p$ may not be homotopic. Nevertheless,
if $\cF_*$ is the pointed fibre, the 
homotopy class of the
loop spaces $\Omega^q\cF_*$, $q\in\bn$
only depend on the path connected component of $*\in\cE$.
Better still, $\Omega\cF_*$ is (homotopic to) 
the homotopy fibre of $\Omega\cE_*\ra\Omega\cB_*$
\end{warning}
\begin{proof}
For any path $c:I\ra\cB$,  $c^*p:c^*\cE\ra I$ is
again a fibration. 
Consequently, if all the geometric fibres of $p$
are spaces, then $c^*p$ is a fibration in spaces, and all
fibres must be homotopic as soon as 
$\cB$ is path connected. 
We can apply this to the path fibration $\rP\cB_*$, {\it i.e.}
for $b:\rp\ra \cB$ in
the connected component of $*$ the fibre, 
\begin{equation}\label{eq:fibextra}
\begin{CD}
\Hom(I,\cB)@<<< \rP\cB_{*,b} \\
@V{Q(0,-)\times Q(1,-)}VV @VVV \\
\cB\times\cB  @<{*\times b}<< \rp
\end{CD}
\end{equation}
must be homotopic to
$\Omega\cB_*$, and reversing this gives a
homotopy to $\Omega\cB_b$, so, the loop space
only depends on the (path) connected component
of $*$. Now let $F$ be the homotopy fibre of
$\Omega\cE_*\ra \Omega\cB_*$. 
This is determined by the universal property
(\ref{eq:pointuni1}) applied to the composition
$\rS^1\times \Omega\cE_*\ra\cE\ra\cB$,
so adopting the notation of (\ref{eq:pointuni1})
in the obvious way we get $G:\Omega\cE_*\ra\Omega\cB_*$,
with a uniquely determined natural transformation
$\eta$ satisfying $\z_\cB=\eta_*p(\z_\cE)$. As such,
we can identify $F$ with pairs, $(a,b)$
where $a\in \Omega\cE_*$,  $b:\rS^1\times I\ra 
\Sigma I_1\ra\cB_*$,
and $G(a)=b(0)$, so that the ambiguity in
lifting $b(0)$ is exactly $\eta_a$.
In addition, 
we have maps:
\begin{equation}\label{eq:fib1}
\begin{CD}
F\times \rS^1\times I \supset F\times \rS^1\times 0 @>{a\times s \mpo a(s)}>> \cE\\
@VV{f:=b\times s\times t\mpo b(s,t)}V @. \\
\cB @.
\end{CD}
\end{equation}
to which we can apply \ref{defn:fib}
to get a lift $H$ of the vertical arrow,
along with some ambiguities $\a$ and $\xi$
satisfying (\ref{eq:ambig}) as it applies
in the current situation. The map $H$
affords a distinguished path,
$h=H\vert_{F\times *\times I}$
with source $H(*):=H(*\times 0)$, 
equivalent to $*$ via $\xi_*\z_\cE:*\Rightarrow a(*)\Rightarrow H(*)$.
Similarly, there are loops
$c_t=H\vert_{F\times\rS^1\times t}$
based at $h_t:=h:{F\times t}\ra \cE$. 
Consequently we 
get some $M:F\times\rS^1\times I \ra \cE: t\mpo h_t^{-1}c_t h_t$,
by way of pulling back $H$ along an appropriate
map from $\rS^1\times I$ to itself,
so that the ambiguities for $M$ are
just the pull-back of those already
introduced. In particular there is
a natural transformation
$
\b:* \Rightarrow p(M_1) 
$,
fabricated from an appropriate restriction of $\a\z_{\cB}$.
This determines, 
\eqref{FibreA2},
a fibre square fitting into a 2-commutative diagram: 
\begin{equation}\label{eq:ambig1}
 \xy
 (-18,10)*+{F\times\rS^1}="T";
 (0,0)*+{\cF}="A";
 (18,0)*+{\cE}="B";
 (18,-18)*+{\cB}="C";
 (0,-18)*+{\rp}="D";
  (39,0)*+{\b:*\Rightarrow \, pM_1}="E";
    {\ar^{M_1} "T";"B"};
    {\ar_{m} "T";"A"};
    {\ar_{} "T";"D"};
    {\ar_{} "A";"B"};
    {\ar_{} "A";"D"};
    {\ar^{p} "B";"C"};
    {\ar_{*} "D";"C"};
{\ar@{=>}^{\s} (14,-4);(4,-14)}
 \endxy
\end{equation}
in which, 
for ease of exposition,
the triangles may, 
\ref{def:FibreS},
be taken to be
strictly commutative. Plainly the map $m$
determines a bunch of circles in the fibre,
but what must be shown- since $*$ is not an
embedding- is that they are all
based in equivalent points. To this end,
there is no loss of generality in supposing
that $p$ is strictly pointed, so the
square in (\ref{eq:ambig1}) is strictly
pointed, 
from 
which 
the commutativity of (\ref{eq:ambig}) evaluated
in the base point
affords another 2-commutative diagram,
\begin{equation}\label{eq:ambig2}
 \xy
 (-18,18)*+{F\times *}="T";
 (0,0)*+{\cF}="A";
 (18,0)*+{\cE}="B";
 (18,-18)*+{\cB}="C";
 (0,-18)*+{\rp}="D";
  (49,0)*+{\b_*=\a_0\z_{\cB}(b_0):*\Rightarrow \, pM_1(*)}="E";
    {\ar^{M_1(*)} "T";"B"};
    {\ar_{*} "T";"A"};
    {\ar_{} "T";"D"};
    {\ar_{} "A";"B"};
    {\ar_{} "A";"D"};
    {\ar^{p} "B";"C"};
    {\ar_{*} "D";"C"};
{\ar@{=>}^{\s_*=\mathbf{1}} (14,-4);(4,-14)};
{\ar@{=>}_{\xi_*\z_\cE} (0,3);(1,7)};
{\ar@{=}_{} (-3,-3);(-5,-5)}
 \endxy
\end{equation}
from which the universal property of fibre
products yields a unique equivalence $*\Rightarrow m(*)$
satisfying \eqref{FibreA5},
so that, finally, the universal property (\ref{eq:pointuni1})
gives a map $r:F\ra \Omega\cF_*$.
The verification that this is a homotopy inverse to the natural
inclusion $i: \Omega\cF_*\ra F$
is the same trick, {\it i.e.} use the
condition (\ref{eq:ambig}) to construct
diagrams \`a la (\ref{eq:ambig1}) and (\ref{eq:ambig2})
but now for $M$ pulled back along $i$ 
to a map $\Omega\cF_*\times\rS^1\times I\ra \cE$,
which becomes the top most arrow in the
first diagram, with it's restriction
to $\Omega\cF_*\times * \times I$ being
top most in the second diagram,
to get a homotopy between $ri$ and the identity,
while that for $ir$ already follows
from (\ref{eq:ambig1}) and (\ref{eq:ambig2}).
Consequently,
\begin{equation}\label{eq:fib2}
\Omega\cF_*\ra \Omega\cE_*\ra \Omega\cB_*
\end{equation}
is a fibration in spaces, so we get \ref{warn:fib}
with even the possibility of making un-necessary 
changes of base points in (\ref{eq:fib2}) since
concatenation of loops gives an isomorphism 
between connected components of the loop space.
\end{proof}
Therefore, as promised, since everything in (\ref{eq:fib2}) 
is a space, most of the following is for free
for separated champs.
\begin{fact}\label{fact:fib1}
If $\cF_*$ is the fibre of a pointed fibration $\cE_*\ra\cB_*$; then we
have a long exact sequence
\begin{equation}\label{eq:fib3}
\cdots \ra \pi_1(\cE_*)\ra  \pi_1(\cB_*) \ra \pi_0(\cF_*)\ra \pi_0(\cE_*)\ra \pi_0(\cB_*)
\end{equation}
\end{fact}
\begin{proof}
The only thing left to do is the exactness from
$\pi_1(\cB_*)$ onwards, of which the only remotely
non-trivial part is $\ra\pi_1(\cB_*) \xrightarrow{\d} \pi_0(\cF_*)\ra$. 
The connecting homomorphism $\d$ may usefully be described
using external Hom; identify,
cf. the proof of \ref{fact:point2},
a pointed loop, $\underline{b}$, in $\cB$ with a triple $(b,\b_0,\b_1)$
consisting of $b:I\ra \cB$ and natural transformations
$\b_q$ of the end points $b(q)$ with $*$, $q\in\{0,1\}$,
modulo the action of natural transformations, $\eta$, 
between paths commuting
with both end points.
The end points transformations lift to equivalences
$\tilde{\b}_q:\cF\rightarrow \cF_q$, 
This gives a base point $0_*:\tilde{\xi}_0(*):\rp\ra \cF_0$,
and by base change we have a fibration $b^*p:b^*\cE\ra I$
with very limited possibilities for ambiguity. 
Indeed, 
the homotopy lifting property affords some strict section
$a:I\ra b^*\cE$ over $b$ with $\xi:0_*\Rightarrow a(0)$,
a natural transformation, so
we get some $\d(\underline{b}):\rp\ra\cF$ by way of:
\begin{equation}\label{eq:fibplus}
*\, {\build\ra_{}^{\b_0}}\, 0_* {\build\Rightarrow_{}^{\xi}} \, 
a(0) \, \xrightarrow{a} \, a(1) \, {\build\ra_{}^{\b_1^{-1}}} \, \d(\underline{b}) 
\end{equation}
Given $\underline{b}$, by  (\ref{eq:point8}) the well definedness
of $\d(\underline{b})$ in $\pi_0(\cF_*)$ check
list is as follows: 
if $(a',\xi')$ is another such pair, then
we form some wedge $a\vee a':T=I\vee_{0}I\ra b^*\cE$.
This may involve some more natural
transformation,
and restrictions to either interval which are only
equivalent to $a$, and $a'$, albeit, given the 
natural transformation at the base point, via a
unique natural transformation 
by (\ref{eq:pointuni1}).
Now pull-back $I\times I$ along $T\times I$,
to get, up to some further equivalences, 
a homotopy with end points in
$p^{-1}(1)$ between either interval in the
wedge and the restriction of the homotopy
lifting to $0\times I\subset T\times I$. As such we have
some maps $A$, $A'$ which are uniquely
equivalent to $a$,  $a'$, and
homotopic (as sections of $p^*b$) to a third map $A''$,
where all of $A$, $A'$, $A''$ have 
the same end point, $r$, which in turn is
equivalent to $0_*$ by some $\xi'':O_*\Rightarrow r$.
The change from either $(a,\xi)$, or $(a',\xi')$,
to $(A,\xi'')$, respectively $(A',\xi'')$,
doesn't effect the equivalence class of $\d(\underline{b})$
since the fibres
are now embedded, {\it i.e.} natural transformations
between the end points in $b^*\cE$ are 
natural transformations in $\cF$, 
while going from either of $(A,\xi'')$ or $(A',\xi'')$
to $(A'',\xi'')$ is just a homotopy
of sections of $b^*p$.
There remains the possibility of changing
$\underline{b}$ by either natural transformations
commuting with the end points and/or homotopies,
but this only effects the diagram (\ref{eq:fibplus})
by equivalence of categories and/or homotopies,
which again, doesn't matter since the fibres
remain embedded and the homotopies preserve the
fibre over $1$.

Having verified that $\d$ is well defined, the
fact that it sends $\pi_1(\cE_*)$ to the 
class of $*$ is a diagram chase similar
to the comparison of (\ref{eq:ambig1}) and
(\ref{eq:ambig2}) but easier, since it
doesn't involve the compatability condition (\ref{eq:ambig}).
Indeed, if a based loop
in $\cE_*$ is described by a triple of
a path, $a$, and natural transformations
$\a_q$ of the end points with $*$, then
there is a diagram akin to (\ref{eq:ambig1})
with square pull-back along $b=pa$, and
top most object $I$ mapping to the fibre
along $A$, say. There is also
a second diagram akin to (\ref{eq:ambig2})
where the base point is deduced from,
$q_*=p(\a_q)$, and the strict pointedness of $p$, 
so that, again, the universal
property of fibre products reveals that
$q_*$ is equivalent to $A(q)$.
Conversely if $\d(\underline{b})$ is in the connected component
of $*$ (\ref{eq:point8}) and the
equivalence of fibres gives a path $c:I\ra \cF_{b(1)}$ with 
natural transformations $\eta$, $\z$ at its ends,
and whence a diagram in $b^*\cE$:
\begin{equation}\label{eq:fib4}
0\xleftarrow{b^*p}\, 
 0_* \, {\build\Rightarrow_{}^{\xi}}\,
a(0) \, \xrightarrow{a} \, a(1) \, 
{\build\Rightarrow_{}^{\eta}} \, c(0)\,
 \xrightarrow{c}\,  c(1)\, {\build\Rightarrow_{}^{\zeta}}  
1_* \xrightarrow{b^*p}\, 1
\end{equation}
Under the natural projection to $\underline{b}^*\cE$,
{\it i.e.} the pull-back over $\rS^1$, $0_*$
gets identified to $1_*$, and (\ref{eq:fib4})
becomes a descent data for some 
$\underline{a}:\rS^1\ra \underline{b}^*\cE$
which, up to homotopy, lifts the base strictly, while the fibre
product is 2-commutative, so the projection
of $\underline{a}$ to $\cE$ is a based loop
such that $p\underline{a}$ equivalent to a homotopy of $\underline{b}$,
as required. The only other minor subtlety
in the rest is if a point  $x:\rp\ra\cF$ is
connected to $*$ in $\cE$. In which case
one easily fabricates some based loop
$\underline{b}$ such that $\d(\underline{b})$ is
equivalent to $x$ in $b^*\cE$, but since
the fibres of $b^*\cE$ are embedded,
$x$ must be equivalent to $\d(\underline{b})$
in $\cF$.
\end{proof}

\subsection{The universal 1 and 2 covers}\label{SS:I.5}

Let $\cX_*$ be a pointed (separated) champ, and denote by
$R_0:=\rP\cX_*\times_{\cX}\rP\cX_*\rras \cX$ the separated 
(and almost certainly not \'etale) groupoid in spaces
defined by the universal property of fibre products.
As usual one says that $\cX$ is locally path connected and
semi-locally simply connected if these conditions hold
in some basis of \'etale neighbourhoods of every 
geometric point of $\cX$. The locally path connected
condition implies that $\rP\cX_*\ra\cX$ is open,
which is the same as universally open, which in turn is an
effective descent condition for spaces, so
$\cX=[\rP\cX_*/R_0]$ in the open topology
as soon as $\cX$ is path connected, or, equivalently
should $\cX_*$ be locally contractible,
appropriate slices of this groupoid yield a 
presentation of $\cX$ in the \'etale topology,
and we further assert:
\begin{lem}\label{lem:cover1}
If $\cX$ is  locally $n-1$ connected, $n\in\bn$, and
semi-locally n-connected,
then the graph, $\G_a$ in $I\times\cX$ of any
path has a basis of \'etale neighbourhoods with
the same property.
\end{lem}
\begin{proof}
For simplicity of notation we can, without
loss of generality, suppose that $\cX$ is path connected.
Let $J_a$ be as per (\ref{eq:path7}),
and 
identify $\G_a$ with the image of
the implied section, $\s_a$,
of the projection $p:J_a\ra I$. By hypothesis
we have a basis of coverings
of $\G_a$ in $J_a$ of the form $I_\a\times U_\a$,
where $U_\a\ra\cX$ is a $n-1$-connected, and 
semi-locally $n$-connected \'etale neighbourhood
of $a(I_\a)$. Without loss of generality the
$I_\a$ are finitely many connected open intervals in $I$
arranged by a total ordering such that any $I_\a$
only meets its immediate neighbours, so, 
in particular any $I_{\a\b}:=I_a\cap I_\b$ is either
empty or connected. Now for $I_{\a\b}\neq \emptyset$, $\a<\b$,
choose a point $p\in I_{\a\b}$ together with a
$n-1$-connected, and semi-locally $n$-connected
neighbourhood $V_{\a\b}\ni \s_a(p)$  which is 
contained in $U_\a\cap U_\b$, with
$Z'_{\a}$,  
respectively $Z''_{\b}$
its complement in $U_\a$, respectively $U_\b$.
Shrinking the cover of the interval appropriately
we can suppose that
$Z'_{\a}\ts\bar{I}_{\a\b}$,
respectively $Z''_{\b}\ts\bar{I}_{\a\b}$, miss the
graph, and for $\a<\b<\g$, we put
\begin{equation}\label{fixMayerV}
V_\b := U_\b\ts I_\b \bsh (Z''_{\b}\ts\bar{I}_{\a\b}
\cup Z'_{\b}\ts\bar{I}_{\b\g})
\end{equation}
with suitable notational adjustment if the interval
$I_\b$ contains an end point.
As such, by construction, $V_\a\cap V_\b$ is
$V_{\a\b}\ts I_{\a\b}$, so, in the first instance
this is connected, whence by Van-Kampen and/or
Mayer-Vietoris, the neighbourhood $\cup_\a V_\a$
of the graph
has the required property provided that the $V_\a$
are themselves $n-1$ connected and semi-locally
$n$-connected. On the other hand for $N_{\a\b}$,
$N_{\b\g}$, $N'_{\b}$ small thickenings of
the closed intervals 
$\bar{I}_{\a\b}$, $\bar{I}_{\b\g}$,
$I'_\b:=\bar{I}_\b\bsh I_{\a\b}\cup I_{\a\g}$,
$V_\b$ is up to a small homotopy
\begin{equation}\label{fixMayerVV}
V_{\a\b}\ts N_{\a\b} \cup U_\b\ts N'_\b \cup V_{\b\g}\ts N_{\b\g} 
\end{equation}
so we conclude by repeated application
of Van-Kampen and/or
Mayer-Vietoris.
\end{proof}

We have a space of paths of paths, to wit:
$\Hom(I,\rP\cX_*)$ which by the universal
property (\ref{eq:pointuni1}) we can
equally identify with the fibre product,
\begin{equation}\label{eq:cover1}
\begin{CD}
I\times \Hom(I\times I, \cX) @<<< I\times\Hom(I,\rP\cX_*)\\
@V{(y,A)\mpo (y,A(0,y))}VV @VVV \\
I\times\cX @<{\mathrm{id}\times *}<< I\times \rp 
\end{CD}
\end{equation}
where (we'll attempt to consistently) label points in
$I\times I$ by way of coordinates $(x,y)$. From
this we can construct a further space via the fibre square,
\begin{equation}\label{eq:cover2} 
\begin{CD}
I\times\Hom(I,\rP\cX_*) @<<<  M'\\
@V{(y,A)\mpo (y,A)\times A(1,y)}VV @VVV\\
I\times\Hom(I,\rP\cX_*)\times\cX@<{\mathrm{id}\times A(1,0)}<< I\times\Hom(I,\rP\cX_*)
\end{CD}
\end{equation}
So, explictly, $M'$ consists of a triple $(A,\z,\xi)$ where
$A:I\times I\ra\cX$ is a map, and $\z_{A(y)}:*\Rightarrow A\vert{0\times I}$,
$\xi_{A(y)}:A(1,0)\Rightarrow A\vert{1\times I}$ are natural
transformations, albeit that $\z_{A(y)}$ is exactly
as in (\ref{eq:pointuni1}), so, in a sense it's fixed.
We also dispose of a discrete sheaf, $\mathrm{Aut}(A(1,0))$
over $M'$, and a section $\xi_{A(0)}$ of the same, and
we define $M$ to be the closed subspace of $M'$ where
this section is the identity.
An element of $M$ may, therefore, usefully be
visualised in a diagram,
\begin{equation}\label{eq:cover3}
 \xy
 (0,0)*+{}="A";
 (18,0)*+{}="B";
 (18,-18)*+{}="C";
 (0,-18)*+{}="D";
  (46,-9)*+{=\,\, A(1,y)\,\, {\build\Longleftarrow_{}^{\xi(y)}}\,\, A(1,0),\,\, \xi(0)=\mathbf{1}}="E";
  (-16,-9)*+{*\,\, {\build\Longrightarrow_{}^{\z}} A(0,1)\,\, =}="F";
(9,-9)*+{A}="G";
(9,2)*+{b}="G";
(9,-20)*+{a}="G";
    {\ar@{-}_{} "A";"B"};
    {\ar@{-}_{} "A";"D"};
    {\ar@{-}^{} "B";"C"};
    {\ar@{-}_{} "D";"C"};
 \endxy
\end{equation}
where $a$, and $b$ are the projections $A(x,0)$, and $A(x,1)$
to $\rP\cX_*$. As usual, concatenation along the edges $a$ or
$b$ 
is only homotopy associative, so it
doesn't yield a groupoid structure in $M$, but we 
do have a map to the groupoid $R_0$ by the definition of
fibre products, or, explicitly for $a$, $b$ paths identified
with their moduli in $\rP\cX_*$, the points of $R_0$ are
given by triples $(a,b,\xi)$ where $\xi:a(1)\Rightarrow b(1)$
is a natural transformation, so, 
the map is, $(A,\z,\xi)\mpo (a,b,\xi(1))$-
the transformation $\z$ with $*$ being the same throughout-
and we assert,
\begin{claim}\label{claim:cover1}
If $\cX$ is locally connected, and semi-locally simply
connected then $R_0$ is locally path connected, and
the image $R_1$ of $M\ra R_0$ is an open and closed sub-groupoid.
\end{claim}
\begin{proof} In the above notation the groupoid structure
in $R_0$ is,
\begin{equation}\label{eq:cover4}
(a,b,\xi)\times (b,c,\eta)\ra (a,c,\eta\xi)
\end{equation}
The effect of concatenation in $M$ is an easier
variant of
\ref{fact:point3}: since no 
rescaling in the horizontal- (\ref{eq:cover3})-
direction is involved the usual concatenation
of squares $A$, $B$ 
with horizontal edges $(a,b)$, and $(b,c)$
along $b$ is already a space of paths in $\rP\cX_*$,
so, no extra natural transformations are involved,
and
$BA$ has edges identically $(a,c)$. Whence the role
of the 
transformation with $*$ is trivial
throughout, while if $\xi(y)$, $\eta(y)$ are the
natural transformations at the end points,
we get a concatenated transformation between
$BA(1,0)=A(1,0)$ by taking $\xi(y)$,
respectively $\eta(y)\xi(1)$, where we should,
which since $\eta(0)=\mathbf{1}$ and the
sheaf of transformations is discrete, this
glues to a transformation, written $\eta\xi(y)$,
so $R_1$ is indeed a groupoid. To see that
it's open look at the image of some $A$ of
the form \eqref{eq:cover3},
and take neighbourhoods,
$J_a$, $J_b$ of the form guaranteed by 
\ref{lem:cover1} of the respective graphs of
$a$, $b$ as well as path-connected open
neighbourhoods $U_a$, $U_b$ of the end
points 
$a(1)$, $b(1)$
with $U_{\bullet}\times \bullet\subset J_{\bullet}$,
for $\bullet\in\{a,b\}$. Now in the notation
of \ref{fact:path3} consider the open
neighbourhood, $V$, of the image of $A$ in the set of
arrows in $R_0$ defined by the section
spaces $S_a$, $S_b$ at either end, with
end points in $U_a$ and $U_b$, 
and $\xi\in U_a\ts_\cX U_b$, which
for $U_a$, $U_b$ sufficiently small
we may suppose is homeomorphic to
the same by the source, respectively the sink.
As ever,
$\xi$ extends uniquely to a neighbourhood
of the right horizontal end, so shrinking
$U_a$ and $U_b$ appropriately, for every
pair of paths $c:I\ra U_a$, 
$d:I\ra U_b$ with $c(1)=a(1)$, $d(1)=b(1)$,
we can extend $A$ to a square,
\begin{equation}\label{eq:cover5}
 \xy
 (0,0)*+{}="A";
 (18,0)*+{}="B";
 (18,-18)*+{}="C";
 (0,-18)*+{}="D";
  (22,0)*+{}="E";
  (22,-18)*+{}="F";
(9,-9)*+{A}="G";
(9,2)*+{b}="H";
(9,-20)*+{a}="I";
(20,-20)*+{c}="J";
(20,2)*+{d}="J";
(20,-9)*+{B}="K";
    {\ar@{-}_{} "A";"B"};
    {\ar@{-}_{} "A";"D"};
    {\ar@{-}^{} "B";"C"};
    {\ar@{-}_{} "D";"C"};
    {\ar@{-}_{} "B";"E"};
    {\ar@{-}_{} "E";"F"};
    {\ar@{-}^{} "C";"F"}
 \endxy
\end{equation}
Here we should be careful, since
by \ref{fact:point3} 
the concatenation $BA$ may, 
in principle,
involve
extra natural transformations. It is,
however, \ref{fact:point1},
uniquely unique, and 
close to $a$, respectively $b$
defined a priori
in $V$
%
while for $U_a$, $U_b$ small enough
the natural transformation $\xi_A(y)$ extends
to some $\tilde{\xi}_{BA} (y)$ between $c$ and every horizontal cross
section over $c$ in (\ref{eq:cover5}), 
so, without loss of generality,
the image of \eqref{eq:cover5} in
$R_0$ stays in $V$.
Now let $(a',b',\xi')\in V$;
and choose a path, $c:I\ra U_a$,
from $a(1)$ to $a'(1)$.
By the
definition of $B$ we have an equivalent path
$d$ in $U_b$, which, since $\xi'\in U_a\ts_\cX U_b$,
and we've taken this to be homeomorphic
to its image whether by the source or sink,
is equally a path from $b(1)$ to $b'(1)$.
Now
we can add $a'$ and $b'$ to the diagram (\ref{eq:cover5}),
and 
since $\pi_1(J_a)$ and $\pi_1(J_b)$ have trivial image in $\pi_1(\cX_*)$,
by \ref{cor:point}, 
we can fill $(a')^{-1}ca$,
respectively $bd(b')^{-1}$, to 
triangles mapping to $\cX$; while
$ca$ and $a'$, respectively $db$ and $d'$
are already in $\rP\cX_*$ so these triangles 
concatenate to (\ref{eq:cover5}) without
occasioning any further natural transformations
on any of $BA$, $a'$, or $b'$, so
$(a',b',\xi')\in R_1$, 
and, en passant, 
$V$ is path connected. To conclude that
$R_1$ is also closed let $A_n$ be a net
whose image $(a_n,b_n,\xi_n)$ converges
to $(a',b',\xi')$, and take $V'$ to be
the same form of connected neighbourhood,
but now  defined by
$G_{a'}$, {\it etc.}. Again, shrinking as necessary, we have a
unique extension $\tilde{\xi}'$ of $\xi'$
to an arrow between $U_{a'}$ and $U_{b'}$,
so as soon as $n$ is large enough, $\tilde{\xi}'$
between the ends of $a_n$ and $b_n$ is $\xi_n$.
Fix one such large $n$; join the ends of $a_n$ and $a'$ by some
$c_n$, and those of $b_n$ to $b'$ by its image
$d_n$ under $\tilde{\xi}'$. Again we have a
diagram like (\ref{eq:cover5}) but this time
for $A_n$, $a_n$, {\it etc.} albeit with the
difference that we need to fill the region
between $c_n$, and $d_n$, so, en passant,
define the right vertical edge. This is,
however, only a question of observing that
$\xi_n$ defines a descent data between
small thickenings of $A_n$ and 
$c_n\times I$, and having effected the gluing
via \ref{fact:point3} we get an
extension $\tilde{\xi}_n$
of a conjugate of $\xi_n$  
between a restriction of $c_n$ and
every horizontal cross section in
a small neighbourhood
of the right vertical of $A_n$. Since $R_0\rras \rP\cX_*$
is proper, $\tilde{\xi}_n$ extends
to the entire region bounded by $c_n$
and $d_n$, 
and since as an automorphism of the join 
of $a_n$ and $c_n$ it's identically,
$\mathbf{1}$, by continuity it's value as
as an automorphism of $c_n$ is also identically $\mathbf{1}$.
The composite, $\tilde{\xi'}^{-1}\tilde{\xi}_n$
between $c_n$ and itself is $\mathbf{1}$
at the join with $A_n$, so by connectedness 
it's exactly $\xi'$ where $c_n$ meets $a'$.
Now we can just fill $a_nc_n(a')^{-1}$ 
{\it etc.} and argue as before.    
\end{proof}

Consequently $R_1$ is just some bunch of
connected components of $R_0$, but better still:

\begin{lem}\label{lem:cover2}
The natural map $\Omega\cX_*\ra R_0$ afforded by
the universal property of fibre products is a
homotopy isomorphism, and $R_1$ is the connected
component of the identity arrows, in fact, its
even the connected component of $(*,*)\in R_0$.
\end{lem}
\begin{proof}
Without loss of generality we can
suppose that $\rP\cX_*\ra \cX_*$ is strictly
pointed, and identify $\Omega\cX_*$ with arrows
in $R_0$ with source $*$. As such, $\Omega\cX_*\hookrightarrow R_0$
is embedded, and in the notation of the proof
of \ref{claim:cover1} we can write $\om\in \Omega\cX_*$
as a triple $(*,c,\g)$ where $c\in \rP\cX_*$ and
$\g:*\Rightarrow c(1)$ is a natural transformation.
Now let $(a,b,\xi)\in R_0$ be given, and apply
the concatenation procedure \ref{fact:point3}
to get a unique path of paths $d_{\lb}:=\{(1-\lb\mathrm{id})^*a\}b\in\rP\cX_*$,
$\lb\in I$,
together with natural transformations $\b_\lb:b\Rightarrow b'_{\lb}:=(1/2)^*d_\lb$
and $\a_\lb:\{(1-\lb\mathrm{id})^*a\}
\Rightarrow a':=(1-1/2\mathrm{id})^* d_\lb$
such that $\a_\lb(0)=\b_\lb(1)\g$ for all $\lb$.
Between 
the part of $a$ not so far employed,
{\it i.e.}
$(1-\lb)^*a$, and its moduli,
$\vert(1-\lb)^*a\vert$ there is another
unique family of natural transformation,
$\xi_\lb:\vert(1-\lb)^*a\vert\Rightarrow (1-\lb)^*a$,
so we get a path $\lb\mpo (\vert(1-\lb)^*a\vert, d_{\lb}, \a_\lb(1)\xi_\lb(1))$
in $R_0$.
Since $\rP\cX_*\ra \cX_*$ is strictly pointed,
$\vert(1-\lb)^*a\vert\ra *$ as $\lb\ra 1$,
and the end point of this path is a based loop,
while by (\ref{eq:pointuni1}), 
the natural transformations $\xi_{\lb}$, $\b_\lb$
go to $\mathbf{1}$ as $\lb\ra 0$, so the
starting point of this path is $(a,b,\g)$.
This procedure
is wholly uniform in $R_0$ and affords
maps $r_\lb:R_0\ra R_0$, with $r_1$
a retraction of the inclusion, $i$ of the
loop space, while yielding a homotopy
between the identity and $ir_1$ in the
other direction.
Now, suppose $(a,b,\xi)\in R_1$ is the
image of some $A$ as in (\ref{eq:cover3}),
then applying $r_1$ to such an element
leads to a based loop which is homotopic
to the concatenation of $a^{-1}$ and $a$,
with the identity as the natural transformation,
which in turn homotopes to $*$ 
by way of the concatenation, again 
via the identity in the join, of $(\lb)^*a$
with it's inverse, so $R_1$ is the 
connected component not just of the identity arrows
but of $(*,*)$ by 
\ref{claim:cover1}.
\end{proof}

Describing the above map $i$ as ``the'' natural map is
a slight misnomer since we equally dispose of 
another, to wit: $j(\om)=(c,*, \g^{-1})$ in the
above notation, or, equivalently $i(\om)^{-1}$
where the inverse is understood in $R_0$, and
the proof of \ref{lem:cover2} shows that $i(\om)$
and $j(\om^{-1})$ belong to the same connected 
component of $R_0$; where now $\om^{-1}$ means
homotopy inverse, {\it i.e.}  
take the data \`a la (\ref{eq:point8}) with end 
points $*$ defining $\om$, write it back to front, and take
it's moduli. Now for $\om$ a based loop we dispose of a 
section 
$i_\om$, respectively, $j_\om$,
of the source, 
respectively the sink,
of $R_0$ via concatenation
\begin{equation}\label{eq:cover6}
\rP\cX_*\ra R_0\,\begin{cases} 
i_\om:a\mpo (a, a\omega , \alpha)\\
j_\om:a\mpo (a\omega, a, \alpha^{-1}) = i_\om(a)^{-1}
\end{cases}
\end{equation}
where $a\omega$ is the concatenation of $a$ and $\om$
in the sense of \ref{fact:point3}, {\it i.e.}
in the above notation, gluing along $\z_a\g^{-1}$
at the join with $\alpha$ the unique induced transformation
of op. cit. from $a(1)\Rightarrow (a\om)(1)$.
This results in four possible ways of going from
$R_1$ to the different connected components of
$R_0$, for example using $i_\om$ one has 
a diagram with every square
fibred
\begin{equation}\label{eq:coverFix}
\begin{CD}
\rP\cX_*@<<< R_{\om, 1}@>{c_{\om,1}}>{f\mpo i_\om f}>
R_0^{\om}
@<{c_{1,\om}} <{f\mpo fi_\om}< R_{1, \om}@>>> \rP\cX_*  \\
@V{si_\om}VV @VV{s_\om}V  @. @V{t_\om}VV @VV{ti_\om}V\\
\rP\cX_*@<t<< R_1 @. @. R_1@>s>> \rP\cX_* 
\end{CD}
\end{equation}
where $R_0^{\om}$ is the connected component
of $i_\om(*)$,
and similarly for $j_\om$. 
Such compositions have slightly different
properties. For example, in (\ref{eq:coverFix}),
$s_\om$ is an isomorphism, and $c_{\om,1}$
an injection whereas $t_\om$ is an 
injection and $c_{1,\om}$ is an isomorphism;
while by \ref{lem:cover2} all of the above
are homotopy equivalences. As an example of
how these differences manifest themselves
let's prove
\begin{fact}\label{fact:cover1}
$R_1$ is a normal sub-groupoid of $R_0$.
\end{fact}
\begin{proof}
Normality means
that if $S\in R_1$ stabilises $a\in P\cX_*$,
then for any $f\in R_0$ with source $a$,
 $f^{-1}Sf$ is a stabiliser in $R_1$. 
To check this write, $f=c_{1,\om}(g)$, for
some $g\in R_{1,\om}$, 
so $f=gi_\om(b)$ for some unique path $b$.
As such
$f^{-1}Sf= i_\om(b)^{-1}(g^{-1}sg) i_\omega(b)$;
while conjugation with $i_\om(b)$ maps $R_1$
to a unique connected component of $R_0$,
whence it takes values in $R_1$. 
\end{proof}
Now we're in a position to prove the obvious, where,
for convenience $\Omega\sbs\Omega\cX_*$ is a complete repetition
free list of representatives of homotopy classes
of bases loops
\begin{factdef}\label{factdef:cover1}
Suppose $\cX$ is path connected, and
let $\cX_1:=[P\cX_*/R_1]$ in the open topology,
then $p:\cX_1\ra \cX$ is a representable \'etale
cover- so a postiori $\cX_1$ is equivalent to an
\'etale groupoid- and there is an action, {\it i.e.} a
fibre square
\begin{equation}\label{eq:coverFix2}
 \xy
 (0,0)*+{\pi_1(\cX_*)\uts \cX_1}="A";
 (38,0)*+{\cX_1 }="B";
 (38,-18)*+{\cX}="C";
 (0,-18)*+{\cX_1 }="D";
  (23,-5)*+{}="E";
  (15,-13)*+{}="F";
    {\ar^{q} "A";"B"};
    {\ar_{p} "A";"D"};
    {\ar^{p} "B";"C"};
    {\ar_{p} "D";"C"};
    {\ar@{=>}_{i} "F";"E"};
 \endxy
\end{equation}
where the leftmost vertical is understood to
be trivial, and the second/horizontal projection, $q$,
is thought of as the action, albeit the
natural transformation $i$- 
which can be identified with $i$ of
\eqref{eq:cover6}-
is,  
\ref{def:Fibre2},
part of the data,
and there is an associativity condition given
by a further diagram
\begin{equation}\label{eq:coverplus}
 \xy
 (0,0)*+{\pi_1(\cX_*)\uts \pi_1(\cX_*)\uts \cX_1}="A";
 (38,0)*+{\pi_1(\cX_*)\uts \cX_1 }="B";
 (38,-18)*+{\cX_1}="C";
 (0,-18)*+{\pi_1(\cX_*)\uts \cX_1 }="D";
  (23,-5)*+{}="E";
  (15,-13)*+{}="F";
    {\ar^{} "A";"B"};
    {\ar_{\mathrm{id}\times q} "A";"D"};
    {\ar^{q} "B";"C"};
    {\ar_{q} "D";"C"};
    {\ar@{=>}_{\alpha} "F";"E"};
 \endxy
\end{equation}
such that 
for multiplication in $\pi_1(\cX_*)$
understood to be pre, rather
than post, concatenation we have
for $\om\mpo F_\om$ the functor
afforded by the action \eqref{eq:coverFix2}
the co-cycle condition
\begin{equation}\label{eq:coverFix1}
p(\alpha_{\tau,\om})=i_{(\tau\om)}i_\om^{-1} F_\om^* i_\tau^{-1}
= i_{(\tau\om)}i_\tau^{-1} (F_\tau)_* i_\om^{-1}
\end{equation}
\end{factdef}
\begin{proof} For brevity, let $P=\rP\cX_*$, 
$s,t$ the source and sink of any groupoid that
may appear, and consider
the fibre square,
\begin{equation}\label{eq:cover91}
\begin{CD}
\cX_1@<<< [R_{0\, s}\times_t R_1\rras R_0]\\
@VV{p}V @VVV\\
\cX_0@<<< P
\end{CD}
\end{equation}
Since $R_1$ is a sub-groupoid of $R_0$,
the groupoid in the top right
hand corner is an equivalence relation,
in fact: $F\sim G$ iff $F=hG$ for some
$h\in R_1$. By (\ref{eq:coverFix}), $F$ and $G$ must
belong to the same connected component, say $R_0^\om$
of $R_0$, and we
can write, $F=fi_\om(a)$, $G=gi_\om(a)$ for some
$f,g\in R_1$, and $a\in P$
whence $F\sim G$ iff $s(F)=s(G)$
so slicing the top right hand corner of
(\ref{eq:cover91}) along $\coprod{i_\om}:\coprod_{\om\in\Omega} P\ra R_0$ 
yields an equivalence of this groupoid
with $P\times \pi_1(\cX_*)$.

This trivially proves that $\cX_1\ra \cX$ is a representable
\'etale cover if $\cX$ has a locally contractible atlas. Otherwise,
it's a bit tedious, because it's  a semi-local
proposition, and
one has to worry about the difference between the open
topology and the \'etale topology. Specifically,
let $U\ra \cX$ be a path-connected
and semi-locally simply connected \'etale neighbourhood
of $x\in U$. The local path connectedness of $U$
implies that it may equally be described as the
classifying champ of $R_U:=\rP U_x\times\rP U_x\rras\rP U_x$.
Now choose a path $d$ from $*$ to $x$;
we can
concatenate with elements of $\rP U_x$ under their
natural image in $\Hom(I,\cX)$ to get a functor,
\begin{equation}\label{eq:cover9}
D:R_U\ra R_0: f=(a,b) \mpo (ad, bd)
\end{equation}
This allows us to describe $U_1=U\times_\cX \cX_1$ as
the classifying champ of,
\begin{equation}\label{eq:cover10}
R_{U\, D(t)}\times_s R'_t\times_s R_1\rras R'
\end{equation}
where $R'\subseteq R_0$ is the set of arrows with
source in the image of $D$. Again this is actually
an equivalence relation on $R'$, where $F\sim G$
iff there are arrows $h\in R_1$ and $u\in R_U$
such that $G=hFD(u)^{-1}$. Since the image of
$\pi_1(U)$ is trivial, the image of arrows under $D$
is contained in $R_1$, and $F,G$ 
must belong to the same connected component of $R_0$.
Consequently if we again write $F=fi_\om(c)$,
$f\in R_1$, $\om\in\Omega$, 
$c\in P$
and choose,
$h=f^{-1}\in R_1$,
then we find a projection in the orbit 
of $F$ to arrows of the form $i_\om (D(a))$,
$a\in \rP U_x$. As such if we slice this groupoid by way of
the fibre square:
\begin{equation}\label{eq:cover11}
\begin{CD}
R_{U\, D(t)}\times_s R'_t\times_s R_1@<<< R''\\
@V{s\times t}VV @VVV\\
R'\times R'@<{i_\om(D(a))}<< \coprod_\om \rP U_x\times \coprod_\om \rP U_x
\end{CD}
\end{equation}
then the lower horizontal arrow is at worst
net, and whence the right vertical is an
equivalent topological groupoid. An arrow
between $(a\om, a)$ and $(b\om, b)$ in $R''$
is, perhaps slightly contrary to the
notation in (\ref{eq:cover10}), a pair $(h,u)$ where $u\in R_U$
is an arrow from $a$ to $b$, 
not just $D(u)$ from $D(a)$ to $D(b)$,
and $h=D^*i_\om(b) D(u) D^*j_\om(a)$,
so the righthand
of (\ref{eq:cover11}) is 
equivalent to
$\coprod_{\om\in\Omega} [U/R_U]$,
and by construction this is isomorphic to $\Omega\times U$.
Of which an isomorphism $\cX_1\times\pi_1(\cX_*)\xrightarrow{\sim}\cX_1\times_\cX\cX_1$
is an easier variant- just replace $D$ in (\ref{eq:cover9}) 
by the inclusion $R_1\ra R_0$, so that 
according to this presentation in
groupoids (\ref{eq:coverFix2}) is given by
\begin{equation}\label{eq:coverFix3}
 \xy
 (0,0)*+{\coprod_{\om\in\Omega} R_1}="A";
 (38,0)*+{R_1}="B";
 (38,-18)*+{R_0}="C";
 (0,-18)*+{R_1 }="D";
  (23,-5)*+{}="E";
  (15,-13)*+{}="F";
    {\ar^{(\om, A)\mpo F_\om (A) } "A";"B"};
    {\ar_{(\om, A)\mpo A} "A";"D"};
    {\ar^{\mathrm{inclusion}} "B";"C"};
    {\ar_{\mathrm{inclusion}} "D";"C"};
    {\ar@{=>}_{a\mpo i_\om(a)} "F";"E"};
 \endxy
\end{equation}
where for $\om\in\Omega$, $F_\om$ is the functor,
\begin{equation}\label{eq:coverplus1}
F_\om: R_1 \ra R_1: A\mpo i_\om(t(A)) A j_\om(s(A))= i_\om(t(A)) A i_\om(s(A))^{-1}
\end{equation}
albeit, in general, $F:\pi_1(\cX_*)\ra \Hom_1(R_1,R_1):\om\mpo F_\om$
is not a homomorphism.
Nevertheless, 
for $\tau,\om\in\Omega$, and $(\tau\om)\in\Omega$ the
representative of their product in $\pi_1(\cX_*)$-
so the representative of the concatenation $\om\tau$
in the notation of \ref{fact:point3}- we have 
an a priori $R_0$ valued, but, in fact by \ref{lem:cover2}
$R_1$ valued, continuous
map 
\begin{equation}\label{eq:coverFix11}
\alpha_{\tau,\om}: P\ra R_1: a\mpo i_{(\tau\om)}(a)i_\om(a)^{-1}i_\tau(a\om)^{-1}
\end{equation} 
which is a
natural transformation between
$F_{\tau}F_{\om}$ and $F_{(\tau\om)}$, 
whence (\ref{eq:coverplus}) and (\ref{eq:coverFix1}). 
\end{proof}
One should, therefore, be careful, in understanding
the action of $\pi_1(\cX_*)$ on $\cX_1$, so for
for the moment we confine
ourselves to describing the universal property of $\cX_1$
by way of,
\begin{fact}\label{fact:cover2}
If $\cX$ is path connected,
the champ $\cX_1$ is simply connected, and for any
other (weakly) pointed representable connected
\'etale covering 
$q:\cY_* \ra \cX_*$ 
there is a  pointed map
$r:\cX_{1,*}\ra \cY_*$, {\it i.e.} a pair $(r,\rho)$ 
where $\rho:*\Rightarrow r(*)$,
and a natural transformation, $\eta:p\Rightarrow qr$
such that,
\begin{equation}\label{eq:cover12}
 \xy
 (0,0)*+{qr(*_{\cX_1})}="A";
 (20,0)*+{q(*_{\cY})}="B";
 (20,-18)*+{*_{\cX}}="C";
 (0,-18)*+{p(*_{\cX_1})}="D";
    {\ar@{=>}^{q(\rho)} "B";"A"};
    {\ar@{=>}_{\eta_*} "D";"A"};
    {\ar@{=>}^{y} "C";"B"};
    {\ar@{=>}_{x} "C";"D"};
 \endxy
\end{equation}
commutes, where $x$, $y$, fixed, afford the 
weak pointing of $p$ and $q$ respectively.
Given $\eta$, $\rho$ is unique, and conversely. Better still:
if $(r',\rho',\eta')$ 
is any other  triple such that
(\ref{eq:cover12}) commutes, then there is a
unique natural transformation
$\xi:r\Rightarrow r'$
satisfying $\xi_*\rho=\rho'$, and this requirement
implies $\eta'=q(\xi)\eta$.
\end{fact}
\begin{proof} For brevity, let $P=\rP{\cX_*}$, and
identify $\cX$ and $\cX_1$ with $[P/R_0]$ and
$[P/R_1]$ respectively- so $p:\cX_1\ra\cX$ will
only be strictly pointed if the path fibration is.
By a liberal application of (\ref{eq:pointuni1}),
\ref{rmk:fib}, and (\ref{eq:fibextra}) the 
path space functor
$\rP$ is, up to unique isomorphism, constant on
pointed representable \'etale covers of a path
connected space, so by \ref{lem:cover2}, $\pi_1(\cX_{1*})=\mathbf{1}$.
More generally, if $q:\cY\ra\cX$ is representable
then the diagonal $\cY\ra \cY\times_\cX \cY$
is a closed embedding, and if $q$ is also an
\'etale covering then it's also an open embedding,
so by base change in the (2-commutative) fibre square,
\begin{equation}\label{eq:cover13}
\begin{CD}
R_0=P\times_\cX P @<<< R':=P\times_\cY P\\
@VVV @VVV\\
\cY\times_\cX\cY @<{\D}<< \cY
\end{CD}
\end{equation}
$R'$ is an open and closed sub-groupoid of $R_0$
with $\pi_0(R')\xrightarrow{\sim} \pi_1(\cY_*)$
by \ref{lem:cover2}. As such proving the proposition
is largely a question of being careful about which
map the arrow $P\ra\cY$ implicit in (\ref{eq:cover13})
actually is. To this end, observe that there is a 
fibre square:
\begin{equation}\label{eq:cover14}
\begin{CD}
P\times q^{-1}(*)=q^*P @>>> P\\
@VVV @VVV\\
\cY @>q>> \cX
\end{CD}
\end{equation}
which may well only be 2-commutative
by way of a natural transformation
$\b_q$, say, between the composition
of the upper two arrows and the lower two.
Now, we're identifying $*_{\cX_1}$ with the
module of $*$ in $\rP\cX_*$, so there is
a unique map $\tilde{*}_\cY:\rp\ra q^{-1}(*)$
whose image in $P$ is $*_{\cX_1}$, and,
for good measure, there's even
a unique natural transformation, $\z$, between
the image 
(which needn't be the same thing as a 
point of the fibre)
of $\tilde{*}_\cY$ in $\cY$
and $*_\cY$
%
such that $q(\z)\b_q(\tilde{*}_\cY)$ is $yx^{-1}$
of
(\ref{eq:cover12}). 
Consequently, among the various possibilities for
$P\ra\cY$, the good choice
is the restriction of the left hand vertical in
(\ref{eq:cover14}) to the connected component of $\tilde{*}_\cY$,
since
by (\ref{eq:cover13}), there is a unique
isomorphism of $R_1$ with the connected component
of the identity in $R'$, and this yields a 
strictly commuting triangle
\begin{equation}\label{MoreFixing}
 \xy
 (-18,0)*+{R_1}="L";
 (18,0)*+{R_0}="R";
 (0,16)*+{R'}="T";
    {\ar@{->}^{} "L";"T"};
    {\ar@{->}^{} "T";"R"};
    {\ar@{->}^{} "L";"R"};
 \endxy
\end{equation}
which by \eqref{eq:cover13} can be identified to
a triangle of inclusions amoongst connected
components of $R_0$. To conclude, observe:

(a) All the unicity statements in \ref{fact:cover2}
either follow from 
the connectivity argument of \ref{fact:point1} or
the fact that points of spaces do not have stabilisers.

(b) As such, we only need to prove the existence
statements, which are all trivial for $\cX$, $\cX_1$,
$\cY$ replaced by the equivalent champs $[P/R_0]$,
$[P/R_1]$, and $[P/R']$, {\it e.g.} one can take
$r$ to be the upper left inclusion in 
\eqref{MoreFixing} and $\eta$ trivial by
the strict commutativity of {\it op. cit.}.

(c) The proposition is independent of the equivalence
class of the data, albeit this has to be
be understood 
in the 2-Category $\et_2(\cX)$ of \ref{defn:cor1}
rather than the naive 
$\underline{\mathrm{Cham}}\mathrm{p}\underline{\mathrm{s}}/\cX$.
\end{proof}
As we've observed, $M\rras \rP\cX_*$ is only homotopy
associative. The argument at the beginning of the
proof of \ref{claim:cover1} do however establish,
\begin{factdef}\label{factdef:cover2}
Define $R_2\ra R_1$ as the quotient of $\mu:M\ra R_1$
by the equivalence relation, $m\sim m'$ iff
$m$ and $m'$ belong to the same path connected
component of the same fibre, then in the quotient
topology $R_2\rras \rP\cX_*$
is a groupoid in spaces.
\end{factdef}
In order to make a finer identification of  topological structure of $R_2$, observe:
\begin{lem}\label{lem:cover4}
The map $M\ra R_1$ is a ``fibration'' in spaces,
{\it i.e.} it has the homotopy lifting property
for maps $f:T\ts I$ where $T$ is a direct sum
of connected spaces, so, for example, $T$
locally connected.
\end{lem}
\begin{proof} Basically the same construction, and a little
easier than \ref{claim:cover1} since there's
no further need to use the hypothesis of 
locally path connected and semi-locally
simply connected, {\it i.e.} given $f:T\times I\ra R_1$,
and a lifting of $f_0$ one argues as in
op. cit. to prove that the set where $f$
can be lifted compatibly with the lifting
of $f_0$ is open and closed.
\end{proof}
From which we arrive to the description of the topology of $R_2$
\begin{fact}\label{fact:cover3}
If $\cX_*$ is path connected, locally $1$-connected,
and semi-locally $2$-connected, then 
$R_2\ra R_1$
is the universal cover of $R_1$ with covering
group $\pi_2(\cX_*)$. 
\end{fact}
\begin{proof}
Let $\Omega=\Omega\cX_{1*}$, so, equivalently
the connected component of $*\in \Omega\cX_{*}$,
which, in either case is path connected and
semi-locally simply connected by \ref{lem:cover1}.
As such it has a universal cover, $\tilde{\Omega}\ra \Omega$.
Now identify $\Omega$ with its image under
$i$, 
{\it i.e.} $\om\mpo i_\om(*)$ in \eqref{eq:cover6},
in $R_1$, and let $r=r_\lb:I\times R_1\ra R_1$, for $r_\lb$,
$\lb\in I$, as in the
proof of \ref{lem:cover2}. Observe that the
fibre, $i^*M$ of $M\ra R_1$ is naturally the path
space of $\Omega$ pointed in $*$, and that in $i^*M$
the equivalence relation of \ref{factdef:cover2}
is simply homotopies of paths. The formation of the quotient topology
commutes with closed embeddings, so $i^*M/\sim=i^*R_2=\tilde{\Omega}$.
Now consider the diagram of fibred squares:
\begin{equation}\label{eq:cover15}
\begin{CD}
\tilde{\Omega}@>>> R'@>>> \tilde{\Omega}\\
@VVV @VVV @VVV \\
\Omega@>>i> R_1@>>r_1>\Omega
\end{CD}
\end{equation}
The composite of the top horizontal arrows
is the identity, so the middle vertical is
a connected \'etale cover with covering
group $\pi_2(\cX_*)$, which is the
universal cover of $R_1$ by \ref{claim:cover1}
and
\ref{lem:cover2}.
Now let $r':M\times I \ra R_1$ be the 
pull-back of $r$ along the natural projection
$M\ra R_1$,
then by \ref{lem:cover4} we can
lift $r'$ to some $\tilde{r}':M\times I\ra M$
By the definition of
the quotient topology and fibre products,
this gives a map $\tilde{r}'_0: R_2\ra R'$. In the
other direction, we have the composition
$R'':=i^*M\times R_1\ra R_1$
and the pull-back of $r$ along the same
to $r'':R''\times I \ra R_1$. 
Over
$1$ this can be lifted to $M$ as the
inclusion of $i^*M$, so we get a
lifting $\tilde{r}'':R''\times I\ra M$,
and whence $\tilde{r}''_0:R''\ra R_2$,
which by another application of the homotopy lifting property
descends to a map $\tilde{r}''_0:R'\ra R_2$;
while the fact that $\tilde{r}'_0$ and
$\tilde{r}''_0$ are mutually inverse follows
from the uniqueness of homotopy liftings
modulo homotopy.
\end{proof} 
Consequently we can conclude this section by introducing,
\begin{defn}\label{defn:coverAddon} 
If $\cX_*$ is locally 1-connected, and semi-locally 2-connected
then the classifying champ $\cX_2:=[\rP\cX_*/R_2]\ra \cX$ is
the universal 2-cover, or, similarly constructed one path
component at a time, if $\cX$ is only locally 1-connected, and semi-locally 2-connected.
\end{defn}
Since the properties of the universal 2-cover are rather
less familiar than those of the universal 1-cover, the
discussion of how and why it is universal is postponed 
till \ref{fact:visit5}.

\subsection{Relative homotopy, developability, and tear drops}\label{SS:I.6}

Let $\cA_*\hookrightarrow\cX_*$ be a strictly pointed
embedded sub-champ of a separated champ $\cX$. The
fibre square,
\begin{equation}\label{eq:tear1}
\begin{CD}
\rP\cX_* @<<< \Omega(\cX_*,\cA_*)\\
@VV{p}V @VV{p}V \\
\cX@<<< \cA
\end{CD}
\end{equation}
defines, by base change a fibration, $\Omega(\cX_*,\cA_*)\ra\cA$, 
with fibre $\Omega\cX_*$, and
we make
\begin{defn}\label{defn@tear1}
For $q\in\bn$ the relative homotopy groups
$\pi_q(\cX_*,\cA_*)$ are defined by,
\begin{equation}\label{eq:tear2}
\pi_q(\cX_*,\cA_*) = \pi_{q-1}(\Omega(\cX_*,\cA_*))
\end{equation}
\end{defn}
By \ref{fact:fib1},  we have a long exact sequence
\begin{equation}\label{eq:tear3}
\pi_0(\cA_*)\la \pi_1(\cX_*,\cA_*)\la \pi_1(\cX_*)\la \pi_1(\cA_*)
\la \pi_2(\cX_*,\cA_*) \la \pi_2(\cX_*)\la\cdots
\end{equation} 
Of particular interest is the case where $\cA_*$
is the embedded sub-champ supported on the base
point, {\it i.e.} $\cA_*=\rB_{G(*)}$, for
$G(*)$, or just $G$ is there is no danger of
confusion, is the local monodromy group of $*$.
In this case by \ref{warn:point1}, $\cA_*$ is
weakly equivalent to a $\rK(G,1)$, so the 
relative homotopy groups are just the homotopy
groups for $q\geq 3$, and (\ref{eq:tear3})
becomes an exact sequence,
\begin{equation}\label{eq:tear4}
1 \la \pi_1(\cX_*,\rB_{G(*)})\la \pi_1(\cX_*)\la G(*)
\la \pi_2(\cX_*,\rB_{G(*)}) \la \pi_2(\cX_*)\la 1
\end{equation}
A champ is said to be {\it developable} if it's universal
cover exists, and is a space, so, obviously
\begin{fact}\label{fact:tear1}
A locally connected and semi-locally simply connected
champ $\cX$ is developable iff for every 
point, $G(*)\ra \pi_1(\cX_*)$ is injective, or,
better by (\ref{eq:tear4}) iff $ \pi_2(\cX_*)\xrightarrow{\sim} \pi_2(\cX_*,\rB_{G(*)})$.
\end{fact}
\begin{proof}
The universal cover exists by \ref{fact:cover2}, and
it's probably most useful to observe
that elements of $\pi_2(\cX_*,\rB_{G(*)})$
may schematically be represented by a diagram
\begin{equation}\label{eq:tear5}
 \xy
 (0,0)*+{}="A";
 (18,0)*+{}="B";
 (18,-18)*+{}="C";
 (0,-18)*+{}="D";
  (-9,-9)*+{*\,\, {\build\Longrightarrow_{}^{\z}} }="F";
(9,-9)*+{A}="G";
(9,2)*+{*}="G";
(9,-20)*+{*}="G";
    {\ar@{-}_{} "A";"B"};
    {\ar@{-}_{} "A";"D"};
    {\ar@{-}^{} "B";"C"};
    {\ar@{-}_{} "D";"C"};
 \endxy
\end{equation}
with $A$ and $\z_A$ as in (\ref{eq:cover3}). A priori
this is missing the extra information of the natural
transformation $\xi(y):*\Rightarrow A(1,y)$
in (\ref{eq:cover3}), but this is implicit since the
free face edge lies in $\rB_{G(*)}$, and the map,
$\pi_2(\cX_*,\rB_{G(*)})\ra G(*)$ just sends this
data to $\xi(1)$ by \ref{warn:point1}. This image is, however,
exactly- beginning of the proof of \ref{claim:cover1}- 
the stabiliser of $*$ in $R_1$.

Alternatively,  the group action in
\ref{factdef:cover1} is transitive on geometric points,
so the fibre of $\cX_1\ra\cX$ over
$\rB_{G(*)}$ is of the form, $\coprod_{x\in\G} \rB_{G_1(*)}$
for some discrete set quotient $\G$ of $\pi_1(\cX_*)$  
and sub-group $G_1(*)$ of $G(*)$, whence the fibre over
$*$ is  $\coprod_{x\in\G} G(*)/G_1(*)$ which in turn is 
naturally isomorphic to $\pi_1(\cX_*)$ as a set,
{\it i.e.} $\G$ is the quotient of $\pi_1(\cX_*)$
by the image of $G(*)$ in (\ref{eq:tear4}), which
is isomorphic to $G(*)/G_1(*)$. 
\end{proof}
Plainly, the elements of $\pi_2(\cX_*,\rB_{G(*)})$
are rather close to spheres, and the discrepancy
from such admits an amusing description. Indeed,
let us restore the natural transformation $\xi$ to
the rightmost edge of (\ref{eq:tear5}) by way of:
\begin{equation}\label{eq:tear5plus}
 \xy
 (0,0)*+{}="A";
 (18,0)*+{}="B";
 (18,-18)*+{}="C";
 (0,-18)*+{}="D";
  (-9,-9)*+{*\,\, {\build\Longrightarrow_{}^{\z}} }="F";
(9,-9)*+{A}="G";
(9,2)*+{*}="H";
(9,-20)*+{*}="I";
(22, -16)*+{}="J";
(22, -2)*+{}="K";
    {\ar@{-}_{} "A";"B"};
    {\ar@{-}_{} "A";"D"};
    {\ar@{-}^{} "B";"C"};
    {\ar@{-}_{} "D";"C"};
    {\ar@{=>}_{\xi(y),\, \xi(0)=\mathbf{1}} "J";"K"};
 \endxy
\end{equation}
Now if we view $\cX_*$ as a gerbe over its moduli,
$X_*$, then the composition $A:I\times I\ra \cX\ra X$
certainly sends the boundary in (\ref{eq:tear5plus})
to $*\in X$. Consequently by 
the topological variant of \cite[1.1]{km}, there is
a neighbourhood $V$ of the boundary such that $A|_V$
factors through some embedding $[U/G(*)]\hookrightarrow \cX$ for 
some (not necessarily faithful) action of $G(*)$ 
on an \'etale neighbourhood $U\ra\cX$ of $*$. As such,
the fibre product $V':=V\times_\cX U\ra V$ is a 
$G(*)$-torsor. Plainly, we may suppose that $V$ is just
the points at distance at most $\e$ from the boundary for some $\e>0$,
thus $V$ is homotopic to $\rS^1$ and the isomorphism
class
 of such torsors are $\Hom(\pi_1(\rS^1), G(*))$, so, on choosing
an orientation of $\rS^1$, we have a non-canonical
identification of $\xi(1)$
with the isomorphism class of the torsor. The torsor
is naturally pointed, so let $V_\infty$ be the 
connected component of the point, then $V_\infty\ra V$
is a $\mu_n$ torsor for
$n$ the order of $\xi(1)$. By definition, the
choice of a generator of $\pi_1(\rS^1)$
gives a generator $\g_n\in \mu_n$,
so 
we get a faithful (and even less canonical) 
representation $\rho:\mu_n\ra G(*):\g_n\ra \xi(1)$,
together with a map $A_\infty:V_\infty \ra \cX$ such that
$\rho(\g): A_{\infty}\Rightarrow \g^*A_{\infty}$,
where we reasonably confuse $G(*)$ with arrows 
in $\cX$ via the embedding $[U/G(*)]\hookrightarrow \cX$.
As such if $A_0$ is the restriction of $A$ to the
complement of the boundary, $U_0$, say, then, notwithstanding
the lack of canonicity occasioned by the above choices,
we have
a descent data, 
\begin{equation}\label{eq:tear6}
R'_n{\build\rras_{s}^{t}} U:=U_0\coprod U_\infty 
{\build\longrightarrow_{}^{A_0\coprod A_\infty}} \cX
\end{equation}
for a groupoid representing the {\it tear drop, $\cS_n$, of order $n$},
{\it i.e.} the ``bad-orbifold'' supported on $\rS^2$ with
a unique non-space like point of order $n$. While the
tear drop is not a space, descent data continues to
have its obvious meaning, {\it i.e.} a natural transformation
between $s^*(A_0\coprod A_\infty)$ and $t^*(A_0\coprod A_\infty)$
satisfying the co-cycle condition. As such we certainly
get a map $\cS_n\ra\cX$. One might, however, wonder
how much ambiguity there is in such a map since not only are 
\ref{factdef:global}, and \ref{fact:point1} 
only for spaces,  
but other choices of $\g_n$ and so forth, will
yield a different descent data/groupoid in (\ref{eq:tear6}). 
The answer is however best possible, 
\begin{fact}\label{fact:tear2}
The open and closed subset $\Omega^2_n\cX_*$
of $\Omega\Omega(\cX_*,\cA_*)$ where $\xi_1$ has  order $n$ 
finely represents pointed (in the ``bad point'') tear drops of order
$n$ whose monodromy group $\mu_n$ at the ``bad point'' injects into
$G(*)$, i.e. 
on fixing a presentation of $\cS_n$
there is a pointed map $Q_n:\Omega^2_n\cX_*\times\cS_{n*} \ra \cX$,
and a natural transformation, $\z_n:*\Rightarrow Q_n(\_, *)$
with the said property on the local monodromy
groups such that if
$(F,\xi):T\times \cS_{n*} \ra \cX_*$ is any other
such pair then there is a unique map $G:T\ra \Omega_n^2\cX_*$
and a unique natural transformation, $\eta:F\Rightarrow Q(G\times\mathrm{id})$
for which $\z_n=\eta_*\xi$.
\end{fact}
\begin{proof}
As ever $\xi_1$ is an element of the discrete sheaf of automorphisms
of $*$, so the set $\Omega_g$ 
where this takes the value $g\in G(*)$
is open and closed. By (\ref{eq:pointuni1}),
there is no difficulty in replacing $U_0$ in (\ref{eq:tear6}) by the
complement of $\Omega_g\times I\times I$ in the
boundary, nor in replacing $U_\infty$ by a neighbourhood of
the boundary, so 
just keeping the above choices, {\it i.e.} $\g_n\mpo g$, {\it etc.}
and taking the union over $g$ of order $n$,
we certainly get a map $Q_n$.
Now, we can always fix a map $r:I\times I\ra\cS_n$
by sending a neighbourhood of the boundary to 
a connected $\mu_n$ torsor, and the interior to the
plane. We therefore have the pull-back $r^*Q_n$, 
and the 
restriction, $p$, to $\Omega^2_n\times I\times I$ of
composition of universal maps $I\times I\times\Omega\Omega(\cX_*,\rB_{G(*)})
\ra I\times \Omega(\cX_*,\rB_{G(*)})\ra \cX$,
which by (\ref{eq:pointuni1})
are related by a natural transformation $\rho:p\Rightarrow r^*Q_n$.
We also dispose of a natural transformation,
$\z'':*\ra p\mid\Omega^2_n\times *$, and we put
$\z_n=\rho_*\z''$. Now suppose we have a pair $(F,\xi)$,
then $r^*F$ affords
some unique $G':T\times I\ra\rP\cX_*$ and a unique natural
transformation, $\eta'$ of the pull-back of $F$ with $Q(G'\times\mathrm{id})$
according to the prescriptions of (\ref{eq:pointuni1}) with
$Q$ the universal map of op. cit. By hypothesis, the end
point of the path $G'(t,y)$ is in $\rB_{G(*)}$ for all $(t,y)\in T\times I$,
so $G'$ factors through $\Omega(\cX_*,\rB_{G(*)})$; while
the paths $G'(t,0)$, and $G'(t,1)$ go to the module of
the constant path, {\it i.e.} $*$, so $G':T\times I \ra \Omega(\cX_*,\rB_{G(*)})$
is the same thing as some $G:T\ra \Omega\Omega(\cX_*,\rB_{G(*)})$,
which by definition must factor through $\Omega_n^2\cX_*$; and
similarly $\eta'$ can be identified with 
some $\eta'':r^* F\Rightarrow p(G\times\mathrm{id})$
such that $\z''=\eta''_* (r^*\xi)$. This gives a
natural transformation, $\rho\eta'':r^*F\Rightarrow Q_nr(G\times\mathrm{id})=
r^*(Q_n(G\times\mathrm{id}))$; while $r$ is an 
isomorphism outwith the boundary and $\cX$ is
separated so there is a unique natural transformation
$\eta:F\Rightarrow Q_n(G\times\mathrm{id})$ such
that $r^*\eta= \rho\eta''$ from which $\z_n=r^*(\eta_*\xi)=\eta_*\xi$.
As to uniqueness:
having fixed $r:I\times I\ra \cS_n$,
$G'$, and whence $G$ is certainly unique, so the
only possible lack of unicity is in $\eta$ given
$G$, but this is again excluded by the same connectedness
argument as in the proof \ref{fact:point1}.
\end{proof}
Let us add some comment by way of
\begin{rmk}\label{rmk:tear1}
While \ref{fact:tear2} is an extremely satisfactory description
of the obstruction to developability, {\it i.e.} a locally
path connected and semi-locally simply connected champ is
developable iff every pointed tear drop can be homotoped to
a sphere, it's not a good idea to try and read the group
structure of $\pi_2(\cX_*,\rB_{G(*)})$ from the homotopy
classes of tear drops since there is no canonical map
from $\mu_n\ra \bz/n$. The group structure is best read
from diagrams such as (\ref{eq:tear5plus}) as described at
the start of the proof of \ref{claim:cover1}. In particular,
like any other relative $\pi_2$, $\pi_2(\cX_*,\rB_{G(*)})$
is a priori non-commutative.
\end{rmk}

\subsection{1-Galois theory and the Huerwicz theorem}\label{SS:I.7}

A priori representable \'etale covers of
a champ $\cX$ form a 2-category. There is,
however, an associated 1-category, to wit:

\begin{factdef}\label{factdef:one1}
Let $\cX$ be a path connected separated
topological champ, then the category
$\et_1(\cX)$ has objects representable \'etale covers 
$q:\cY\ra\cX$ and arrows diagrams,
\begin{equation}\label{eq:one1}
 \xy
 (0,0)*+{\cY}="A";
 (28,0)*+{\cY'}="B";
 (14,-12)*+{\cX}="C";
{\ar_{}^{f} "A";"B"};
    {\ar^{q'} "B";"C"};
    {\ar_{q} "A";"C"};
{\ar@{=>}^{\eta} (10,-6);(18,-6)}
\endxy
\end{equation}
modulo the equivalence relation $(f,\eta)\sim (f',\eta')$
iff there is a natural transformation $\xi:f\Rightarrow f'$
such that $\eta'=q'(\xi)\eta$;
with composition the evident composition of
diagrams. In particular, because the objects
of $\et_1(\cX)$ are representable covers,
$\et_1(\cX)$ has fibre products in the sense
of a 1-category, {\it i.e.}
with strict unicity in the universal property.
\end{factdef} 
The definition of  $\et_1(\cX)$ is forced by the 
the universal property \ref{fact:cover2}, it's
also forced by the axiomatic Galois 
theory of a 1-category in the sense of \cite[Expos\'e V.5]{sga1},
\cite{noohi1}, {\it e.g.} 
\begin{lem}\label{lem:one1}
If $\cX$ is path connected locally and globally, and
semi-locally 1-connected, then
the fibre functor $\et_1(\cX)\ra\Ens:q\mpo q^{-1}(*)$
is representable. Indeed, we have a canonical isomorphism:
\begin{equation}\label{eq:one2} 
q^{-1}(*)\xrightarrow{\sim}\Hom_{\et_1}(\cX_1,\cY)
\end{equation}
\end{lem}
\begin{proof}
To give a point in $q^{-1}(*)$ 
is equivalent to giving a map
$*_\cY$, and a natural transformation
$y:*\Rightarrow q(*_\cY)$
as encountered in 
\ref{fact:cover2}, {\it cf.} \ref{fact:cor2},
which
for $\cY'$ the connected components of $*_\cY$,
determines
a weakly pointed map, $\cY'_*\ra\cX_*$, so we can
apply \ref{fact:cover2}, to find a triangle:
\begin{equation}\label{eq:one3}
 \xy
 (0,0)*+{\cX_1}="A";
 (28,0)*+{\cY}="B";
 (14,-12)*+{\cX}="C";
{\ar_{}^{r} "A";"B"};
    {\ar^{q} "B";"C"};
    {\ar_{p} "A";"C"};
{\ar@{=>}^{\eta} (10,-6);(18,-6)}
\endxy
\end{equation}
and the uniqueness statement of \ref{fact:cover2} is
equivalent to this yielding a well defined map
from left to right in (\ref{eq:one2}) of which
the inverse,  
for $p$ strictly pointed,
sends a map $r:\cX\ra\cY$
to the point of the fibre determined by
$r(*):\rp\ra\cY$, and the natural transformation
$\eta_*:* \Rightarrow qr(*)$.
\end{proof}
Now from \ref{factdef:cover1} we have an isomorphism
$\cX_1\times_\cX \cX_1\xrightarrow{\sim}\cX_1\times \pi_1(\cX_*)$, 
given by the fibre square (\ref{eq:coverFix2})
so (cf. (\ref{eq:coverplus1}) if $\om\mpo F_\om$
is the second projection  
of (\ref{eq:coverFix2})
for every
$\om\in \pi_1(\cX_*)$ we have a triangle,
\begin{equation}\label{eq:one3plus}
 \xy
 (0,0)*+{\cX_1}="A";
 (28,0)*+{\cX_1}="B";
 (14,-12)*+{\cX}="C";
{\ar_{}^{F_\om} "A";"B"};
    {\ar^{q} "B";"C"};
    {\ar_{p} "A";"C"};
{\ar@{=>}^{i_\om} (10,-6);(18,-6)}
\endxy
\end{equation}
for $i_\om$ the natural transformation of op. cit.,
which by (\ref{eq:coverFix1})- or, indeed, the
unicity of fibre products of representable maps-
satisfies $i_{\tau\om}=p(\alpha_{\tau,\om})F_\om^*(i_\tau) i_\om$
for $\a$ as in (\ref{eq:coverplus}), and multiplication
of loops pre, rather than the habitual post, concatenation.
Consequently,
\begin{cor}\label{cor:one1}
For $q=p$ the isomorphism (\ref{eq:one2}) is actually
an isomorphism of groups, so for general $q$, (\ref{eq:one2})
defines
a right action of $\pi_1(\cX_*)$ on $q^{-1}(*)$.
\end{cor}
\begin{proof} As we've observed we have a map of
groups, 
$\om\ra (F_\om, i_\om)$,
and the set isomorphism from right to
left of (\ref{eq:one2}) occurring at
the end of the proof of  
\ref{lem:one1}, so it will suffice to
check that their composition is the
identity. This is, however, a formal
consequence of \ref{factdef:cover1},
{\it i.e.} one just takes a point
$*\in p^{-1}(*)$ and adjoins another
fibre square to (\ref{eq:coverFix2})
to obtain
\begin{equation}\label{eq:oneFix2}
 \xy
(-38,0)*+{\pi_1(\cX_*) }="X";
 (-38,-18)*+{\rp}="Y";
(0,0)*+{\pi_1(\cX_*)\uts \cX_1}="A";
 (38,0)*+{\cX_1 }="B";
 (38,-18)*+{\cX}="C";
 (0,-18)*+{\cX_1 }="D";
  (23,-5)*+{}="E";
  (15,-13)*+{}="F";
    {\ar^{\om\ts x\mpo F_\om(x)} "A";"B"};
    {\ar_{\om\ts x\mpo x} "A";"D"};
    {\ar^{p} "B";"C"};
    {\ar_{p} "D";"C"};
    {\ar^{} "X";"A"};
    {\ar_{*} "Y";"D"};
    {\ar^{} "X";"Y"};
    {\ar@{=>}_{i_\om} "F";"E"};
 \endxy
\end{equation}
of which the totality is a fibre square, and
whence the unique lifting of $(F_\om(*), i_\om(*))$
to the top left hand corner is $\om$, or,
more generally an isomorphism thereof if $p$
is not strictly pointed.
\end{proof}

With these pre-liminaries, we get the usual Galois correspondence:
\begin{fact}\label{fact:one1}
If $\cX$ is path connected locally and globally, and
semi-locally 1-connected, then
the fibre functor $\et_1(\cX)\ra\Ens(\pi):q\mpo q^{-1}(*)$
viewed as taking values in sets with $\pi_1(\cX_*)$-action is
an equivalence of categories.
\end{fact}
\begin{proof}
The right hand side of (\ref{eq:one2}) doesn't
depend on the base point, so if two maps yield
the same fibre functor at one point, they yield
the same fibre functor at every point; while 
the objects of $\et_1(\cX)$ are representable
\'etale covers whence the functor is certainly
faithful. To go the other way, we can always
reduce to $\pi_1(\cX_*)$-sets on which the 
action is transitive, so equivalently the
set of cosets $\pi_1(\cX_*)/K$ for some sub-group
$K$. Associated to such, however, by \ref{lem:cover2} is the 
unique open and closed 
sub-groupoid $\iota: R\hookrightarrow R_0$ such that $\iota_*(\pi_0(R))=K$
under the  natural identification of $\pi_0(R_0)$
with $\pi_1(\cX_*)$ of op. cit.;
and the resulting classifying champ $[\rP\cX_*/R]$
is a representable \'etale cover of $\cX$ with fibre
exactly $\pi_1(\cX_*)/K$.
\end{proof}
Let us apply this to obtain an appropriate version
of Huerwicz's theorem. 
A priori singular homology and singular co-homology present
something of an issue for champs, {\it e.g.} by
\ref{factdef:global} it doesn't a priori make sense
to talk about the space of simplices since this
might only be a champ, and if one were to take
its moduli then, \ref{warn:point1}, information can be lost.
The derived functor definition of co-homology
makes perfect sense, however, and coincides
with the C\v{e}ch co-homology 
at the level of $\rH^1$.
By
definition C\v{e}ch 1-co-cycles are descent data,
or slightly more elegantly, for some presentation
$R\rras U$ of $\cX$, a  C\v{e}ch co-cycle with 
values in a group $G$ is just a map $\phi:R\ra G$
such that,
\begin{equation}\label{eq:one4}
U\times G \ni (s(f),x)\xleftarrow{\s} R\times G \xrightarrow{\tau} 
(t(f), \phi(f)x)\in  U\times G
\end{equation}
defines a groupoid. As such, cf. \ref{fact:cover2} and
(\ref{eq:one3plus}), locally constant
sheaves $\underline{G}$ are group objects in $\et(\cX_1)$,
and the $\underline{G}$-C\v{e}ch co-homology,
$\underline{G}$-torsors so that by \ref{fact:one1},
\begin{cor}\label{cor:one2}
If $\cX$ is path connected, locally path connected,
and semi-locally 1-connected, then the category
of locally constant sheaves $\underline{G}$ on $\cX$
is equivalent to the category of groups $G$ with
(right) $\pi_1(\cX_*)$-action, and we have a 
functorial isomorphism,
\begin{equation}\label{eq:oneplus}
\rH^1(\cX, \underline{G})\xrightarrow{\sim}\rH^1(\pi_1(\cX_*), \underline{G})
\end{equation}
\end{cor}
Now using the path fibration we can bump this up to
\begin{fact}\label{fact:one2}
If for $n\in\bn$, $\cX$ is locally and globally $n$-connected,
and semi-locally $n+1$-connected, then for any abelian
group $A$ there are functorial isomorphisms,
\begin{equation}\label{eq:one5}
\rH^q(\cX, A) \xrightarrow{\sim} \begin{cases} 0& \text{if $1\leq q\leq n$},\\
\Hom(\pi_q(\cX_*), A)& \text{if $q=n+1$}
\end{cases}
\end{equation}
\end{fact}
\begin{proof}
Let $P=\rP\cX_*$ be the path space, then we have a complex of spaces,
\begin{equation}\label{eq:one6}
\cdots P^{(3)}:= P\times_\cX P\times_\cX P {\build\ra_{\ra}^{\ra}}
P^{(2)}:= P\times_\cX P\, (=R_1)\rras P^{(1)}:= P
\end{equation}
Since $\cX$ is locally path connected, the path
fibration is open so by \cite[5.3.3]{deligne3},
for any sheaf of abelian groups $\cA$,
there is a spectral sequence
\begin{equation}\label{eq:one7}
\rE_2^{p,q} = \check{\rH}^p( \rH^q(P^{(p)}, \cA)) \Rightarrow \rH^{p+q}(\cX, \cA)
\end{equation}
where the C\v{e}ch co-homology is that defined by the
complex (\ref{eq:one6}). Now $P\ra\cX$ is a fibration,
and $P$ homotopes to $*$ keeping the base point fixed
by \ref{fact:point2}, so each $P^{(p)}$ is homotopic
to $p$-copies $\Omega^{(p)}$ of the loop space $\Omega=\Omega\cX_*$.
By \ref{lem:cover1}, $\Omega$ is locally and globally
$n-1$ connected, and semi-locally $n$-connected; while
C\v{e}ch co-homology of constant sheaves is homotopy
invariant ({\it e.g.} use the Leray spectral sequence
for the projection $f:\cX\ts I\ra \cX$, and $R^qf_*\cF$ 
coincides, by direct calculation, with relative C\v{e}ch,
whence vanishes for $q>1$, for any sheaf $\cF$ independent
of any hypothesis on $\cX$)
so by induction the rows, $1\leq q\leq (n-1)$ will all
be identically zero as soon as $n\geq 2$. Similarly, the first
column is wholly zero for $q>0$, as is the first row
for $p>0$ since by hypothesis $n\geq 1$ so
every $\Omega^{(p)}$ is connected which implies that $d_1^{p,0}$
is the identity for $p$-odd and zero for $p$-even. 
Consequently, 
by induction in $n$, (\ref{eq:one7}) degenerates
at $\rE_2$ for $p+q\leq n+1$, and
we need only check that the kernel of,
\begin{equation}\label{eq:one8}
d_1^{1,n}:\rH^n(\Omega, A) \ra \rH^n(\Omega^2, A)
\end{equation}
is the non-trivial case of (\ref{eq:one5}). The
three maps in (\ref{eq:one6}) from $\Omega^2$ to $\Omega$
are the 2-projections and concatenation, so for
$n=1$, this differential is zero by \ref{cor:one2},
whence in general by induction in $n$ with everything
being functorial by
the functoriality of spectral sequences.
\end{proof}

\newpage

\section{Topological 2-Galois theory}\label{S:II}

\subsection{2-groups}\label{SS:II.1}

Suppose that a {\it topological 2-type} is given,
{\it i.e.} a triple $\underline{\pi}_2:=(\pi_1,\pi_2, k_3)$ where 
$\pi_1$ is a group acting on the left of the
abelian group $\pi_2$ and $k_3$ a class in
$\rH^3(\pi_1,\pi_2)$. Associated to such a 
data is a (weak) 2-group, $\Pi_2$. In detail:
the underlying category of $\Pi_2$ is the
groupoid $(s,t): \pi_1\ltimes\pi_2 \rras \pi_1$ with
source and sink the projection of $\pi_1\ltimes\pi_2\ra \pi_1$,
{\it i.e.} the action on $\pi_1$ is trivial.
The monoidal product, as a map of arrows is just
the group law in $ \pi_1\ltimes\pi_2$,
\begin{equation}\label{eq:group1}
\otimes: \Pi_2\times\Pi_2\ra\Pi_2: (\tau, B)\times (\om, A)\mpo 
(\tau\om, BA^{\tau})
\end{equation}
The identity object is the identity in $\pi_1$,
and similarly inverses are actual inverses in
$\pi_1$, whence, so defined, they're strict
inverses. Finally, to define the associator,
simply choose a normalised co-cycle, $K_3$, 
representing the {\it Postnikov class} $k_3$, so as to get a natural
transformation
\begin{equation}\label{eq:group2}
\kappa: \otimes(\mathrm{id}\times\otimes)\Rightarrow
\otimes(\otimes\times\mathrm{id}):\mathrm{ob}(\Pi_2^3)
\rightarrow \mathrm{Ar}(\Pi_2)  : (\sigma,\tau,\om)\mpo K_3(\sigma,\tau,\om)
\in \mathrm{stab}(\sigma\tau\om) 
\end{equation}
where we profit from the fact that $\pi_2$ is
canonically the stabiliser of any object in
$\Pi_2$ by the definition of the underlying
groupoid. Plainly, 
there is an a priori ambiguity
in the definition, in the form of the
lifting $K_3$ of $k_3$. If, however, $\tilde{K}_3$
were another, and $\tilde{\Pi}_2$ the 2-group
so constructed, then the difference $K_3-\tilde{K}_3$
is the co-boundary of a normalised 2 co-chain,
$c$, which in turn defines a natural transformation,
\begin{equation}\label{eq:group3}
\g:\otimes\Rightarrow \tilde{\otimes}: 
\mathrm{ob}(\Pi_2^2)
\rightarrow \mathrm{Ar}(\Pi_2):
(\tau, \om)\mpo c(\tau, \om)
\end{equation}
so the resulting structures are isomorphic as 
2-groups, {\it i.e.} the pairs 
$I=(\mathrm{id},\g)$, $J=(\mathrm{id},\g^{-1})$
afford (weak) monoidal
functors 
\begin{equation}\label{eq:groupPlus}
I:\Pi_2\rightarrow \tilde{\Pi}_2 
\quad J:\tilde{\Pi}_2 \rightarrow \Pi_2
\end{equation}
such that $IJ$, respectively $JI$ are the identity,
and $I$, respectively $J$ are even unique modulo
automorphisms of $\Pi_2$, respectively $\tilde{\Pi}_2$. 

Now, there is an obvious notion of morphisms
$\underline{\pi}_2\ra \underline{\pi}'_2$ of
topological 2-types, {\it i.e.} maps $f_p: \pi_p\ra \pi_p'$,
$p=1$ or $2$ compatible with the module structure
such that, $f_1^*k_3'=(f_2)_* k_3$. It is not,
however, the case that morphisms of 2-types
can be functorially lifted to morphisms of
2-groups. More precisely,
\begin{fact}\label{fact:group1}\cite[Theorem 43]{baez}
There is a functor $\text{2-Groups}\ra\text{2-types}:\mathfrak{P}\mpo 
\underline{\pi}$
such that,
\begin{enumerate}
\item[(a)] If $\mathfrak{P}$ is a 2-group, and $\Pi_2$ is the
2-group constructed as above from the 2-type of $\mathfrak{P}$
then there is an equivalence of 2-groups $\Pi_2\ra \mathfrak{P}$. 
\item[(b)] The map from 1-homomorphisms of 2-groups,
modulo 2-homomorphisms to maps of 2-types
$$
(\Hom_1(\mathfrak{P}, \mathfrak{P}')/\sim)  \ra \Hom(\underline{\pi},
\underline{\pi}')  
$$
exhibits the former over the latter as a principal
homogeneous space under $\rH^2(\pi_1, \pi_2')$.
\end{enumerate}
\end{fact}
Our ultimate interest is the action of 2-groups on
groupoids, and in light of the objection \ref{faq8}
the only philosophically satisfactory definition
appears to be
\begin{defn}\label{def:groupNick}
An action of a 2-group on a groupoid (or more
generally a category) $\cF$ is a map of
2-groups $\gP\ra\fAut(\cF)$.
\end{defn}
On the other hand, one tends to encounter actions
(particularly transitive ones) as isomorphisms
of certain fibre products (cf. \ref{factdef:cover1}),
and so it is useful to have  
\begin{altdef}\label{altdef:group1}
A left action of a 2-group $\mathfrak{P}$ on a groupoid 
(or more generally a category)
$\cF$ is a
functor, $A:\mathfrak{P}\uts \cF\ra \cF$- 
$\uts$ strict fibre product, \ref{def:FibreS}-
and  natural 
transformations $\a: A(\otimes\uts \mathrm{id}_\cF)\Rightarrow 
A(\mathrm{id}_{\mathfrak{P}}\uts A)$,
$\b:\mathrm{id}_\cF\Rightarrow A(\mathbf{1}\uts\mathrm{id}_\cF)$
such that,
\begin{equation}\label{eq:group4}
 \xy
 (0,0)*+{\gP_1\uts \gP_2\uts\gP_3\uts \cF}="A";
 (60,0)*+{\gP_1\uts\gP_{23}\uts\cF}="B";
 (60,-35)*+{\gP_1\uts\cF}="C";
(0,-35)*+{\gP_1\uts\gP_2\uts\cF}="D";
(50,10)*+{\gP_{12}\uts\gP_3\uts\cF}="E";
 (110,10)*+{\gP_{123}\uts\cF}="F";
 (110,-25)*+{\cF}="G";
(50,-25)*+{\gP_{12}\uts\cF}="H";
(59,-34)*+{}="X";
(100,-26)*+{}="Y";
(60,5)*+{\kappa\Rightarrow}="Z";
(45,-32)*+{\a_{12}\Rightarrow}="T";
(35,-18)*+{\a_{23}\Downarrow}="U";
(70,-7)*+{\a_{(12)3}\Downarrow}="V";
(90,-10)*+{\a_{1(23)}\Downarrow}="W";
(25,-15)*+{\mathbf{1}\Downarrow}="W";
{\ar_{\id_1\uts\otimes_{23} } "A";"B"};
    {\ar^{A_{23} } "B";"C"};
    {\ar_{\id_1\uts A_2 } "D";"C"};
{\ar_{}_{\id_{1}\uts\id_2\uts A_3 } "A";"D"};
    {\ar^{\otimes_{(12)\uts 3} } "E";"F"};
    {\ar^{A_{123} } "F";"G"};
{\ar@{-->}_{A_{12}} "H";"G"}    
{\ar@{-->}_{\id_{12}\uts A_3} "E";"H"}
{\ar@{-->}^{\otimes_{12}\uts\id_{\cF}} "D";"H"}
{\ar^{\otimes_{12}\uts\id_3\uts\id_{\cF}  } "A";"E"};
    {\ar_{ \otimes_{(1)\uts (23)} } "B";"F"};
{\ar_{A_1} "X";"Y"};
\endxy
\end{equation}
2-commutes; with similar diagrams involving $\b$ and the
inverse and identity laws whose exact form we ignore since
in practice our 2-groups will have strict identities and
strict inverses so that the required 2-commutativity of
such will not be an issue.
\end{altdef}
Let us therefore verify
\begin{fact}\label{fact:group2}
There is a $1$ to $1$ correspondence between actions 
in the sense of \ref{def:groupNick} and \ref{altdef:group1}.
\end{fact}
\begin{proof} 
The automorphism group $\gF$ of a groupoid is the
2-group with objects functors from $\cF$ to itself
admitting a weak inverse
and arrows natural transformations between them.
A priori there is an issue about the functoriality
of inversion, but this can be remedied by choosing
an equivalent category $A:\cF_0\hookrightarrow \cF$
with exactly one object for every isomorphism class
in $\cF$ together with a weak inverse $B:\cF\ra\cF_0$
to $A$; so that $F\mpo F_0:= BFA$ yields a strictly
invertible automorphism $F_0$ of $\cF_0$, and whence
an inverse functor $\mathrm{inv}:\gF\ra\gF:F\mpo AF_0^{-1}B$.
The action of a map $\xi:F\Rightarrow G$ of functors
on an arrow $f:x\ra y$ is the arrow $\xi_yF(f):F(x)\ra G(y)$,
which like $\gF$ is strictly associative, {\it i.e.} the
a priori 2-commutative diagram (\ref{eq:group4}) is
actually strictly commutative. Plainly if we have
a map of 2-groups $\gP\ra\gF$, then the composition
with the natural action of $\gF$ on $\cF$ yields
an action of $\gP$ on $\cF$. Conversely, given an
action $A:\gP\times\cF\ra\cF$ to an object $\om$ of
$\gP$ one associates the functor $A_\om:\cF\ra\cF$
which sends an arrow $f$ to $A(\mathbf{1}_\om, f)$;
while an arrow $S:\om\Rightarrow\tau$ goes to the
natural transformation of $A_\om$ and $A_{\tau}$
defined by sending an object $x$ of $\cF$ to the
arrow $A(S,\mathbf{1}_x)$. This procedure defines
a functor $\underline{A}:\gP\ra\gF$, while, 
more or less by definition,
the
natural transformations
$\a$, $\b$ of \ref{altdef:group1} yield  natural transformations 
$\underline{\a}:\underline{A}(\otimes_{\gP})\Rightarrow
\underline{A}\otimes_{\gF}\underline{A}$, 
$\underline{\b}: \mathbf{1}_{\gF}\Rightarrow A(\mathbf{1}_{\gP})$.
The 2-commutative diagram (\ref{eq:group4}) becomes a
2-commutative diagram in which every occurrence of
$A$, respectively $\a$, becomes an occurrence of
$\underline{A}$, respectively $\underline{\a}$,
and similarly for the conditions on inverses and
identities; so that the triple $(\underline{A}, \underline{\a}, \underline{\b})$
is indeed a map of 2-groups; while an inspection
of the above formulae show that the composition
$\underline{A}:\gP\ra\gF$ with the natural action
of $\gF$ on $\cF$ is exactly (not just equivalent to)
the original action $A$, and similarly the 
map $\gF\ra\gF$ obtained by underlining the
natural action is the identity exactly.
\end{proof}

Just as \ref{fact:group1}, actions on groupoids
can be simplified considerably by reduction to
the case where there is a bijection between 
objects of $\cF$ and isomorphism classes of 
objects. Already therefore, as in the proof
of \ref{fact:group2}, automorphisms, $a$, of such
a $\cF$ must be strictly invertible, and should
$a$ send an object $x$ to $y$, then the stabiliser
groups of $x$ and $y$ must be isomorphic.
Consequently such a groupoid $\cF$ is a direct
sum of groupoids $\cF_i$ where every two objects
of $\cF_i$ have the isomorphic stabilisers, $\G_i$,
say, 
while for $i\neq j$, $\G_i$ and $\G_j$ fail to
be isomorphic. Plainly, $\fAut(\cF)$ is isomorphic
to the direct product of $\fAut(\cF_i)$. As
such to determine the automorphisms of a 
groupoid, we're reduced to studying the 
{\it skeletal} case, {\it i.e.} $\cF$ is 
defined by a discrete set of objects, $F$, and a group
$\G$ according to the rule,
\begin{equation}\label{eq:group5}
\Hom(x, y) = \begin{cases} \emptyset& \text{if $x\neq y$},\\
\G& \text{if $x=y$}
\end{cases}
\end{equation}
Invertible functors 
can therefore be written as a pair
$(\om, A)\in \mathrm{Sym}_F \times \Home(F,\Aut(\G))$
which act on arrows by way of
\begin{equation}\label{eq:groupPlus1}
x{\build\longrightarrow_{}^{g}} x \quad \mpo \quad 
\om(x)\xrightarrow{A(\om(x))(g)} \om(x)
\end{equation}
so for $A^\om(x):=A(\om^{-1}(x)$
the left action of $\mathrm{Sym}_F $
on $\Home(F,\Aut(\G))$ they form the group
\begin{equation}\label{eq:group6}
G_1:=\mathrm{Sym}_F \ltimes \Home(F,\Aut(\G))
\end{equation}
Similarly
the natural transformations are $G_2:=\Home(F, \G)$,
which we understand as a left $G_1$-module via
\begin{equation}\label{eq:group8}
(\om, A) \times g\mpo \{x\mpo A(x)(g_{\om^{-1}(x)})\}
\end{equation}
so that the underlying category of $\fAut(\cF)$
is the groupoid,
\begin{equation}\label{eq:group7}
G_1\ltimes G_2\rras G_1
\end{equation}
with source the given functor, and sink
the functor obtained via the 
action of $\Home(F, \G)$ on $\Home(F,\Aut(\G))$
via inner automorphisms. 
As such $\fAut(\cF)$ is
synonymous with the crossed module- \cite[IV.5]{brown},
\cite[\S 10]{noohiuni}- defined, 
for $Z$ the centre of $\G$,by the exact sequence
\begin{equation}\label{eq:group9}
1\ra \G_2=\Home(F,Z)\ra
G_2\xrightarrow{\mathrm{Inn}} 
G_1
\ra  \mathrm{Sym}_F \ltimes \Home(F,\mathrm{Out}(\G)) =\G_1
\ra 1
\end{equation}
To such an exact sequence, \cite[IV.5]{brown}, there is
an associated co-homology class $\obs_\cF\in \rH^3(\G_1, \G_2)$ which,
as the notation suggests, we will view as an obstruction
class, and
\begin{fact}\label{fact:group3}
The 2-type of the automorphism group $\gF$ of
a groupoid, $\cF$, in which all stabiliser groups are
isomorphic to the same group $\G$ is the triple
$\underline{\phi}:=(\G_1, \G_2, \obs_\cF)$ as defined
above for $F$ the moduli, i.e. isomorphism classes of
objects, of $\cF$. 
If $\G$ is abelian the obstruction class
is $0$, otherwise it may well be non-trivial,
e.g. if $F$ is a point then the 2-type $\underline{\phi}$
is the topological 2-type of the classifying space
of the universal fibration in $\rK(\G,1)$'s.
\end{fact}
\begin{proof} Just compare the proof of \ref{fact:group1}
in \cite{baez} with \cite[IV.5]{brown}. The fact that
$\G_p=G_p$, $p=1$ or $2$ implies, op. cit., that the
obstruction is trivial if $\G$ is abelian, otherwise,
already for $F$ a point, the examples of non-triviality are legion,
\cite[IV.6]{brown}. 
\end{proof}
Plainly, therefore, if we have an action of
$\Pi_2$ on $\cF$ we have a map of 2-types
$\underline{\pi}\ra\underline{\phi}$. In
particular $F$ is a left $\pi_1$ set, and
we have a 1 co-cycle
\begin{equation}\label{eq:groupPlus2}
\bar{A}_\om:\pi_1\ra \Home(F, \Out(\G))\quad 
\bar{A}_{\tau\om}=\bar{A}_\tau \bar{A}_\om^\tau
\end{equation}
affording $\pi_1\ra\G_1$, while the action
of $\pi_1$ on $\G_2$ is
\begin{equation}\label{eq:groupPlus3}
\om\ts z\mpo\{ x\mpo A_{\om}(x) (z_{\om^{-1}\cdot x})\}
\end{equation}
which is well defined for any lifting
$A_\om$ of $\bar{A}_\om$ to a family of 
automorphisms of $\G$.

Now suppose that the action of $\Pi_2$ 
on $\cF$ is
{\it transitive, i.e.} $\pi_1$ acts transitively
on $F$, and that the latter is pointed in
$*$. As such, we have the stabiliser
$\pi_1'$ of $*$, the image $\pi_2'$ of
$\pi_2$, and an induced Postnikov class
$k_3'\in\rH^3(\pi_1',\pi_2')$ which we
may view as the restriction $K_3'$ of
the co-cycle (\ref{eq:group2}). Consequently
we have a 2-group $\Pi_2'$ deduced from $\Pi_2$
with 2-type $\underline{\pi}'=(\pi_1',\pi_2',k_3')$
whose action on $\cF$ affords an action
on the one point category $\rB_\G$ supported
at $*$. For this action the co-cycle, $\bar{A}$
of (\ref{eq:groupPlus2}) restricts to a
representation
\begin{equation}\label{eq:groupPlus4}
\bar{A}_\om(*): \pi_1'\ra \Out(\G)_*:=\Home(*,\Out(\G))
\end{equation}
which by (\ref{eq:groupPlus3}) affords, in turn,
a representation on the centre $Z_*:=\Home(*,Z)$,
and we have a particularly simple map of 2-types
$\underline{\pi}'\ra \pi_*=(\Out(\G)_*, Z_*, \mathrm{obs}_*)$.
The resulting action of $\Pi_2'$ on 
$\rB_\G$ will be referred to as the
{\it pointed stabiliser} action, and
we assert
\begin{fact}\label{fact:groupPlus1}
The restriction $2$-functor from pointed transitive
$\Pi_2$ groupoids to actions of sub 2-groups
(i.e. $\Pi_2'$ 
with $K_3'$ induced from $K_3$
for $\pi_1'$ any sub-group of
$\pi_1$, and
$\pi_2'$ any quotient group of $\pi_2$) on
a group, i.e. a
groupoid with $1$-object, defined by
sending a $\Pi_2$-category to its pointed
stabiliser action is an equivalence of
2-categories.
\end{fact}
The proof 
should be the Leray spectral
sequence for $\rB_{\pi_1}\ra\rB_{\pi_1/\pi'_1}$,
but since $\pi'_1$ needn't be normal,
we get 
a series of lemmas in
group co-homology instead. To this end it is
convenient to identify a pointed 
transitive $\pi_1$-set $F_*$ with the
left action on the right co-sets
$\pi_1'/\pi_1$ of the stabiliser. As such
if $\bar{\rho}\in F$ is a set of elements
in $\pi_1$ identified with the
right co-sets then for $X$ any set
the left action of $\pi_1$ on $\Home(F, X)$ is
\begin{equation}\label{eq:groupPlus5}
\pi_1\ts \Home(F,X)\ra  \Home(F,X):x\mpo x^\om=
\{\bar{\rho}\mpo x (\bar{\rho\om})\}
\end{equation}
We also dispose of a (set) map
\begin{equation}\label{eq:groupPlus6}
\pi_1\ra\pi_1': x\mpo x':= x\bar{x}^{-1} 
\end{equation}
along with the possibility to choose $\bar{1}=1$, and we assert
\begin{lem}\label{lem:groupPlus1}
Let 
$\pi_1$ be a group, and
$G$ another group (with, say, trivial $\pi_1$
action, since this is all we need) then 
for $F$ the set of right cosets of
a sub-group $\pi_1'$ with $\pi_1$ acting on the left of
$\Home(F, G)$ via (\ref{eq:groupPlus5})
the natural restriction
\begin{equation}\label{eq:groupPlus7}
\mathrm{res}:
\rH^1(\pi_1, \Home(F, G))\ra 
\rH^1(\pi_1', G)\, (=
\Hom(\pi_1', G)/\mathrm{Inn}_G)
: A\mpo A(*)\vert_{\pi_1'}
\end{equation}
is an isomorphism. Indeed, the inverse is given by inflation 
\begin{equation}\label{eq:groupPlus8}
\mathrm{inf}: B\mpo \{\bar{\rho}\ts \om \mpo 
B((\bar{\rho}\om)')\}
\end{equation}
\end{lem}
\begin{proof} Let $A:\om\mpo A_\om$ be the $\pi_1$ co-cycle of
(\ref{eq:groupPlus7}), and put $a(\bar{\rho})=A_{\bar{\rho}}(*)$,
then $\tilde{A}:\om\mpo aA(a^\om)^{-1}$ is
equivalent. However, by construction,
$\tilde{A}_{\bar{\rho}}(*)=1_G$ for all $\bar{\rho}\in F$,
so $\tilde{A}_{\om}(\bar{\rho})= \tilde{A}_{(\bar{\rho}\om)'}(*)$
by the co-cycle condition.
\end{proof}
In the particular case of $G=\Out(\G)$ we can 
define an automorphism $a:\cF\ra\cF$ of the
category (\ref{eq:group5}) by way of the
functor $(1, A_{\bar{\rho}})$ for some liftings
of the outer automorphisms $\bar{A}_{\bar{\rho}}$
of (\ref{eq:groupPlus4}) to honest automorphisms.
Plainly conjugation by $a$ is a (strict) isomorphism
of $\pi_1$-groupoids, and so
\begin{cor}\label{cor:groupPlus1}
Any (skeletal) pointed transitive $\Pi_2$
groupoid is isomorphic to one in which the
outer co-cycle (\ref{eq:groupPlus2}) satisfies
(\ref{eq:groupPlus8}). 
\end{cor}
Now let us apply similar considerations to
the higher co-homology groups, {\it i.e.}
\begin{lem}\label{lem:groupPlus2}
Let $\pi_1$ be a group, with $F$ the set
of right co-sets of a sub-group $\pi_1'$
such that, for some abelian group $Z$,
$\pi_1$ acts (\ref{eq:groupPlus3}) on $\Home(F,Z)$ by way of
a $\Hom(F,\Aut(Z))$ valued $1$ co-cycle $A$
satisfying (\ref{eq:groupPlus8}) then 
the restriction
map
\begin{equation}\label{eq:groupPlus9}
\mathrm{res}:
\rH^n (\pi_1, \Home(F,Z)) \ra \rH^n(\pi_1', Z):K\mpo K(*)\vert_{\pi_1'}
\end{equation}
is an isomorphism
for all $n\in\bn\cup\{0\}$, with inverse 
the inflation
\begin{equation}\label{eq:groupPlus10} 
\mathrm{inf}:K\mpo \{(x_1,\cdots, x_n)\times \bar{\rho}\mpo
K((\bar{\rho}x_1)', (\overline{\rho x_1}x_2)',
\cdots, (\overline{\rho x_1\cdots x_{n-1}}x_n)')\}
\end{equation}
where, here and elsewhere, formulae such
as (\ref{eq:groupPlus10}) should be read
as functions of as many variables as make
sense, e.g. the inflation of a constant
is the constant function of $\bar{\rho}$.
\end{lem}
\begin{proof} Let $(S^n\subseteq \Home(\pi_1^n\ts F, Z), \pa)$,
$n\in\bn\cup\{0\}$,
be the standard complex for the computation of
the left hand side of (\ref{eq:groupPlus9}) 
via normalised co-chains
and
consider the effect of the homotopy,
\begin{equation}\label{eq:groupPlus11}
h:S^{n+1}\ra S^n: K\mpo\{ (hK): (x_1,\cdots, x_n)\times \bar{\rho}\mpo
K(\bar{\rho}, x_1,\cdots, x_n) (*)\}
\end{equation}
which affords the formula
\begin{equation}\label{eq:groupPlus12}
\pa h + h\pa = \mathrm{id} - P
\end{equation}
where $P$ is the projector $K\mpo PK$ which has values
\begin{equation}\label{eq:groupPlus13}
(x_1,\cdots, x_n)\times \bar{\rho}\mpo
\begin{cases}
K(\bar{\rho}x_1, x_2,\cdots, x_n)(*)
-K(\overline{\rho x_1},  x_2,\cdots, x_n)^{(\overline{\rho}x_1)'}(*)
&\text{if $n\in\bn$}\\
K(*)&\text{if $n=0$}
\end{cases}
\end{equation}
In particular the projector satisfies
\begin{equation}\label{eq:groupPlus14}
(PK)(x_1,\cdots, x_n)(\bar{\rho})=(PK)(\bar{\rho}x_1,\cdots x_n)(*),
\quad\text{and,}\quad (PK)(\bar{x}_1,x_2, \cdots,x_n)(*)=0
\end{equation}
and the co-homology of $(S^n,\pa)$ is the co-homology
of the image of $P$, {\it i.e.} the complex 
\begin{equation}\label{eq:groupPlus15}
C^n=\{K\in\Home(\pi_1^n, Z)\,\vert\, K(x_1,\cdots, x_n)=0\,\,
\text{if},\, x_1=\bar{x}_1,\, \text{or},\, x_i=1,\, i\geq 2\}
\end{equation}
where, critically, the differential 
$d:C^n\ra C^{n+1}:K\mpo dK$
is
\begin{equation}\label{eq:groupPlus16}
K(\bar{x}_1x_2,x_3,\cdots, x_{n+1})^{x'_1}
-K({x}_1x_2,x_3,\cdots, x_{n+1}) +\cdots +
(-1)^{n+1}K(x_1,\cdots, x_n)
\end{equation}
We have therefore found the $p=0$ term in the
sequence of filtered complexes
\begin{equation}\label{eq:groupPlus17} 
C^n=F^0C^n{\build\hookleftarrow_{}^{\mathrm{inf}}}
F^1C^n{\build\hookleftarrow_{}^{\mathrm{inf}}}
F^2C^n\cdots {\build\hookleftarrow_{}^{\mathrm{inf}}}
F^pC^n
\end{equation}
where $F^pC^n$ is the sub-group of $\Home((\pi_1')^p\ts \pi_1^{n-p}, Z)$,
or indeed just $\Home((\pi_1')^n, Z)$ if $p\geq n$
such that
\begin{equation}\label{eq:groupPlus18}
K(x_1, \cdots, x_n)=0,\, \text{if any}\, x_i=1,\, 
1\leq i\leq n,\, 
\text{or,}\,
\bar{x}_{p+1}=x_{p+1}
\end{equation}
while for $p\in\bn$ the inflation maps 
$\mathrm{inf}:F^pC^{n=p+q}\ra F^{p-1}C^{n=(p-1)+(q+1)}$
are given by
\begin{equation}\label{eq:groupPlus19}
\mathrm{inf}(K)(t_1,\cdots,t_{p-1},x_1,\cdots, x_{q+1})
=K(t_1,\cdots, t_{p-1}, x'_1, \bar{x}_1 x_{2}, x_{3},\cdots, x_{q+1})
\end{equation}
and the differentials $d_p:F^pC^{p+q}\ra F^{p}C^{p+q+1}$ are
\begin{equation}\label{eq:groupPlus20}
\begin{split}
 K(t_2,\cdots, t_{p}, x'_{1}, \bar{x}_1x_2,x_3,\cdots, & x_{q+1})^{t_1}
-K(t_1t_2,\cdots, t_p, x'_{1}, \bar{x}_1x_2,x_3,\cdots, x_{q+1}) +\\
K(t_1, t_2t_3,  \cdots, t_p, x'_{1}, \bar{x}_1x_2,x_3,& \cdots,  x_{q+1})
+\cdots +\\ 
(-1)^{p-1}
K(t_1, \cdot, t_{p-1}t_p,  x'_{1}, \bar{x}_1x_2, &x_3,\cdot,  x_{q+1})+
(-1)^p  K(t_1, \cdot, t_{p-1}, t_px'_{1}, \bar{x}_1x_2,x_3,\cdot,x_{q+1}) \\
+ (-1)^{p+1} K( t_1,   \cdot, t_p, x_1x_2, x&_3,\cdot, x_{q+1})
+(-1)^{p+2} K(t_1,\cdot, t_p, x_1, x_2 x_3,\cdot, x_{q+1})
+\cdots+  \\
 (-1)^{p+q} K( t_1,  \cdots,  t_p, x_1, x_2,&\cdots, x_q x_{q+1})
+ (-1)^{p+q+1} K(t_1,\cdots, t_p, x_1, x_2,\cdots, x_q)
\end{split}
\end{equation}
Consequently on understanding $t_1=x'_1$ if $p=0$,
(\ref{eq:groupPlus20}) coincides with (\ref{eq:groupPlus16})
in this case, so that up to a fixed term in the
$x_j$'s we have almost a standard differential
in the $t_i$'s except for the one term where this transgresses
into $x_j$'s, and conversely, from right to left
we have a standard differential in the $x_j$'s
except for the said transgression into the $t_i$'s.
In any case, the composition of $\mathrm{inf}$ followed
by the natural restriction $\mathrm{res}:F^pC^n\ra F^{p+1}C^n$
is the identity, while in the other direction we 
consider the effect of the homotopies
$h_p: F^pC^{p+q+1}\ra F^{p}C^{p+q}$
\begin{equation}\label{eq:groupPlus21}
(h_pK)(t_1,\cdots,t_p,x_1\cdots x_q)=\begin{cases}
(-1)^{p+1}K(t_1,\cdots,t_p,x'_1,\bar{x}_1, x_2,\cdots x_q)
& \text{if $q\geq 1$,}\\
0 &\text{if $q\leq 0$}
\end{cases}
\end{equation}
which yield for some projector $P_p$ of groups
\begin{equation}\label{eq:groupPlus22}
h_pd_p-d_ph_p=\mathrm{id}-P_p
\end{equation} 
and whence a projector of complexes. Indeed
$(-1)^{p}P_p$ is the sum of (group) projectors
\begin{equation}\label{eq:groupPlus23}
\begin{split}
K(t_2,\cdots, t_p, x'_1, (\bar{x}_1x_2)', \overline{x_1x_2},\cdots & x_q)
-K(t_1t_2,\cdots t_p, x'_1, (\bar{x}_1x_2)', \overline{x_1x_2},\cdots x_q)\\
+\cdots + \, 
(-1)^pK(t_1,\cdots, t_{p-1},  t_px&'_1,  (\bar{x}_1x_2)', \overline{x_1x_2},\cdots x_q) +\\
(-1)^{p+1} K(t_1,\cdots, t_p, (x_1x_2)', \overline{x_1x_2}&,\cdots x_q)
+(-1)^{p+2}K(t_1,\cdots t_p, x'_1, \bar{x}_1x_2, \cdots x_q)
\end{split}
\end{equation}
whenever $q\geq 2$,  
while $P_pK(t,x_1)=K(t,x'_1)$ 
for $q=1$, and plainly $P_p$
is the identity
for $q\leq 0$. Furthermore, 
if $\bar{x_1}x_2=\overline{x_1x_2}$,
then the first $p+1$ terms in (\ref{eq:groupPlus23})
vanish 
by (\ref{eq:groupPlus18})
because $(\bar{x_1}x_2)'=1$, while
the last 2-cancel since $(x_1x_2)'=x'_1$. 
Similarly- (\ref{eq:groupPlus18}) again- 
all but the last term 
in (\ref{eq:groupPlus23}) vanish if $K$ is in the
image of inflation, so we obtain
\begin{equation}\label{eq:groupPlus24}
P_p=\mathrm{inf}(\mathrm{res}P_p),\quad 
\mathrm{id}=(\mathrm{res}P_p)\mathrm{inf}
\end{equation}
Consequently (\ref{eq:groupPlus17}) is a chain
of homotopic complexes, and taking $p>n$ yields
(\ref{eq:groupPlus9})-\eqref{eq:groupPlus10}.
\end{proof}
By way of a minor variant of
\ref{lem:groupPlus1} and \ref{lem:groupPlus2}
we also have
\begin{lem}\label{lem:groupPlus3}
Let all the hypothesis be as in \ref{lem:groupPlus2}-
so in particular the form of the action on $\Home(F,Z)$
but with $Z$ non-abelian, then exactly the same
conclusions hold for $n=0,1$.
\end{lem}
Now let us apply these considerations
to the proof of \ref{fact:groupPlus1}
by way of the following assertions
\begin{claim}\label{claim:groupPlus1}
For every action of a sub 
(in the sense of \ref{fact:groupPlus1})
2-group
$\Pi_2'$ on some $\rB_\G$ there is
an inflated pointed transitive action 
of $\Pi_2$ on the pointed skeletal
groupoid $\cF_*$ with objects the right cosets
$F:=\pi_1'\bsh \pi_1$ pointed in the
identity and every automorphism group
$\G$. Furthermore, the restriction of
inflation is the identity, while
the
inflation of restriction is equivalent
to the identity.
\end{claim}
\begin{proof} 
Given the action of $\Pi_2'$ on $\G$,
we have a representation $A'_1:\pi_1'\ra \Out(\G)$
which can be inflated by the formula
(\ref{eq:groupPlus8}) to a co-cycle
$\bar{A}_1:\pi_1\ra\Out(\G)$, and whence,
cf. (\ref{eq:groupPlus2}),
a representation $A_1:\pi_1\ra \G_1=\mathrm{Sym}(F)\ltimes
\Home(F,\Out(\G))$ for $\pi_1$ acting
on the left of $F$. Similarly, we
have a 
representation $A_2':\pi_2'\ra Z$
in the centre of $\G$ which by
(\ref{eq:groupPlus3}) inflates
uniquely to a representation
$A_2:\pi_2 \ra \G_2=\Home(F,Z)$
compatible with $A_1$, and
the actions of $\pi_1$ on $\pi_2$,
respectively $\G_1$ on $\G_2$ by
way of the formula
\begin{equation}\label{eq:groupPlus25}
S\mpo \mathrm{inf}(A'_2)(S)=A_2(S)
=\{\bar{\rho} \mpo A'_2(S^{\bar{\rho}}\vert_{\pi_2'})\}
\end{equation}
By \ref{fact:group1}, or better, 
\cite[Theorem 43]{baez}, the remaining
component of an action is a normalised
2 co-chain $c:\pi_1^2\ra \Home(F,Z)$
such that $(A_2)_* K_3- (A_1)^*\mathrm{obs}_\cF$
is the co-boundary $\pa(c)$.
We already have, however, a
2 co-chain $c':(\pi'_1)^2\ra Z$ such
that $(A'_2)_* K'_3- (A'_1)^*\mathrm{obs}_{\rB_\G}= 
\pa(c')$; while by  
\ref{lem:groupPlus2}
applied to $\G_1$ we can, without
loss of generality, suppose 
that $\mathrm{obs}_\cF=\mathrm{inf} (\mathrm{obs}_{\rB_\G})$.
Similarly, by a double application 
of \ref{lem:groupPlus2}:
$(A_2)_* K_3-\mathrm{inf}((A'_2)_* K_3')=\pa k$ for some
co-chain $k$ whose restriction is identically
zero, and whence $k$ is unique modulo co-boundaries.
As such, if we take $c=\mathrm{inf}(c')+k$,
then we get an action 
(well defined up to equivalence)
of $\Pi_2$ on $\cF$
which restricts identically to the given action of
$\Pi_2'$ on $\rB_\G$.
Conversely, if 
$c$, or equivalently
a pointed
transitive $\Pi_2$ action on $\cF_*$, is given, then
$c-(k + \mathrm{inf}\mathrm{res}(c))$ is a co-cycle
whose restriction is 0, so it's a co-boundary, and
the inflation of restriction is equivalent to
what we started with. 
\end{proof}
To which (\ref{eq:groupTom}) is very pertinent; while
the functoriality of \ref{claim:groupPlus1} is
also critical, so 
\begin{rmk}\label{rmk:groupPlus1}
In the notation of the proof of \ref{claim:groupPlus1},
an explicit solution of the equation
\begin{equation}\label{eq:groupFix1}
(\mathrm{inf}(A'_2))_* K_3=
\mathrm{inf}((A'_2)_* (\mathrm{res}K_3)) + \pa(k)
\end{equation}
is provided by $k=(A'_2)_* (hK_3)$ for 
$(hK_3): F\ts\pi_1^2\ra \pi_2$ defined by
\begin{equation}\label{eq:groupFix2}
(\bar{\rho},\tau,\om)\mpo
K_3(\bar{\rho},\tau,\om) - 
K_3((\bar{\rho}\tau)', \overline{\rho\tau}, \om)
+ K_3( (\bar{\rho}\tau)', (\overline{\rho\tau}\om)', \overline{\rho\tau\om})
\end{equation}
which is a little simpler than the proof of
\ref{lem:groupPlus2} might suggest, so checking
it directly is perhaps simpler than reverse
engineering.
\end{rmk}
In any case,
this proves that the restriction functor
is essentially surjective on 0-cells,
so, \cite[1.5.13]{tom}, \ref{fact:groupPlus1}
will follow if we can establish
\begin{claim}\label{claim:groupPlus2}
Let $\cF_*$, $\cG_*$ be $\Pi_2$ transitive
pointed groupoids then the restriction
functor
\begin{equation}\label{eq:groupPlus26}
\fHom_{\Pi_2}(\cF_*, \cG_*)\ra
\fHom_{\Pi'_2, \Pi''_2}(\mathrm{res}_{\Pi'_2}(\cF), 
\mathrm{res}_{\Pi''_2}(\cG))
\end{equation}
is an equivalence of 1-categories; where
$\mathrm{res}_{\Pi'_2}(\cF)$,
respectively $\mathrm{res}_{\Pi''_2}(\cG)$,
is to be viewed as a $\Pi'_2$ action
on a group $\rB_\G$, respectively
$\Pi''_2$ on $\rB_\D$, such that an
arrow of 2-groups $\Pi'_2\ra \Pi''_2$
is given and the 2-types $\underline{\pi}'\ra \underline{\pi}''$
satisfy $p_1:\pi'_1\hookrightarrow \pi''_1$
is injective, respectively
$p_2:\pi'_2\twoheadrightarrow \pi''_2$
is surjective,
and $(p_2)_* K'_3=(p_1)^*K''_3$.
\end{claim}
\begin{proof}
Quite generally if $\cF$, $\cG$ are
groupoids, or even just categories,
with actions $A$, respectively $B$
by the source and sink of a 
 map of (any) 2-groups $(p, \mu):\gF\ra 
\gG$ then (its philosophically
unsatisfactory nature notwithstanding)
the objects of $\fHom_{\gF,\gG}(\cF,\cG)$
are best described via \ref{altdef:group1}
as pairs $(f,\xi)$ forming a 2-commutative
diagram
\begin{equation}\label{eq:groupPlus27}
 \xy
 (0,0)*+{\gF\uts \cF}="A";
 (38,0)*+{\cF }="B";
 (38,-18)*+{\cG}="C";
 (0,-18)*+{\gG\uts\cG }="D";
  (23,-5)*+{}="E";
  (15,-13)*+{}="F";
    {\ar^{A} "A";"B"};
    {\ar_{p\uts f} "A";"D"};
    {\ar^{f} "B";"C"};
    {\ar_{B} "D";"C"};
    {\ar@{=>}_{\xi} "F";"E"};
 \endxy
\end{equation}
which in turn form a 2-commutative diagram
\begin{equation}\label{eq:groupPlus28}
 \xy
 (0,0)*+{\gF_1\uts \gF_2\uts \cF}="A";
 (60,0)*+{\gF_{12}\uts\cF}="B";
 (60,-35)*+{\gG_{12}\uts\cG}="C";
(0,-35)*+{\gG_1\uts\gG_2\uts\cG}="D";
(50,10)*+{\gF_{1}\uts \cF}="E";
 (110,10)*+{\cF}="F";
 (110,-25)*+{\cG}="G";
(50,-25)*+{\gG_{1}\uts\cG}="H";
(59,-34)*+{}="X";
(100,-26)*+{}="Y";
(60,5)*+{\alpha\Rightarrow}="Z";
(45,-32)*+{\beta\Rightarrow}="T";
(40,-18)*+{\Downarrow\mu\uts f}="U";
(70,-7)*+{\xi_{(1)}\Downarrow}="V";
(90,-10)*+{\xi_{(12)}\Downarrow}="W";
(25,-15)*+{\id_1\uts\xi_2\Downarrow}="W";
{\ar_{\otimes_{12}\uts\id } "A";"B"};
    {\ar^{p_{12}\uts f } "B";"C"};
    {\ar_{ \otimes_{12}\uts\id } "D";"C"};
{\ar_{}_{p_1\uts p_2\uts f } "A";"D"};
    {\ar^{A } "E";"F"};
    {\ar^{f } "F";"G"};
{\ar@{-->}_{B} "H";"G"}    
{\ar@{-->}_{p_{1}\uts f} "E";"H"}
{\ar@{-->}^{\id_1\uts B_2} "D";"H"}
{\ar^{\id_1\uts A_2 } "A";"E"};
    {\ar_{ A } "B";"F"};
{\ar_{B} "X";"Y"};
\endxy
\end{equation}
and 
since we're always able to suppose
that 2-groups, their morphisms, and
actions send identities to identities
strictly, \ref{eq:groupPlus28}
actually forces $\xi\vert_{\mathbf{1}\uts \cF}$
to be the identity. Similarly
morphisms $\z:(f,\xi)\Rightarrow (g,\eta)$
are just natural transformations $\z:f\Rightarrow g$
affording the obvious 2-commutative diagram.

Now in our cases of interest $\mu$ will
always be the identity; while in the
left hand side of 
\ref{eq:groupPlus26} we not only have
$\gF=\gG=\Pi_2$, but, without loss of
generality, we may suppose that
$\pi''_2=\pi'_2$. To profit from the
transitivity of the action:
let $E$, $F$, $G$ be the spaces of right
cosets $\pi'_1\bsh\pi''_1$, $\pi'_1\bsh \pi_1$, 
$\pi''_1\bsh \pi_1$ with pointing in the
identity and left 
$\pi''_1$, respectively
$\pi_1$, action as appropriate. 
In particular, if we identify $E$ with
a set of representatives $e\in\pi''_1\subseteq \pi_1$,
and similarly $G$ with $g\in\pi_1$, then
the set $eg$ represents $F$, and the 
operations of inflation from $\pi'_1$ to
$\pi''_1$ followed by $\pi''_1$ to $\pi_1$,
in the sense of \ref{lem:groupPlus1}-\ref{lem:groupPlus3},
strictly commutes with that from $\pi'_1$
to $\pi''_1$. As such, we'll 
either continue to denote cosets
by a $\bar{\,}$ or, if there is
risk of ambiguity,
write $x\mpo e(x)$,
for the representative of the
coset $\pi_1''x$, {\it etc.}, and extend the notation of
(\ref{eq:groupPlus6}) by way of
\begin{equation}\label{groupPlusPlus1}
\pi_1\ra\pi''_1:x\mpo x'':= xg(x)^{-1},\quad
\pi''_1\ra\pi'_1:x\mpo \dot{x}:= xe(x)^{-1}
\end{equation}
so that $x'=\dot{x}''$. In any case,
we can identify $f$
with a map $x\mpo f(x)$ in $\Home(F,\Hom(\G,\D))$;
$\xi$ with $(\om,x)\mpo \xi_\om(\om\cdot x)$
for $\xi_\om\in \Home(F,\D)$, and 
the 2-commutativity \ref{eq:groupPlus27} yields
\begin{equation}\label{eq:groupPlus29}
f(\om\cdot x) A_\om (\om\cdot x) =
\mathrm{Inn}_{\xi_\om(\om\cdot x)} B_\om (\om\cdot f(x)) f(x)
\end{equation}
for functors $A_\om$, $B_\om$
deduced from the action as in the proof
of \ref{fact:group2}, whence, modulo
conjugation the outer representation
(\ref{eq:groupPlus2}).
Passing to outer representations, we
conclude from \ref{lem:groupPlus1} 
and (\ref{eq:groupPlus29})
that
modulo conjugation in $\D$, $f$ is
the map  
\begin{equation}\label{eq:groupPlusPlus2}
\mathrm{inf}(f(*)):F(=E\cdot G)\ra\Hom(\G,\D):eg\mpo B_{e}^{-1}f(*)
\end{equation}
Similarly, viewing $\a$ as $
(\tau,\om,x)\mpo \a_{\tau,\omega}((\tau\omega).x))$,
(\ref{eq:groupPlus28}) becomes
\begin{equation}\label{eq:groupPlus30}
f(\a_{\tau,\om})\xi_{(\tau\om)}=\xi_\tau\xi^\tau_\om \b_{\tau,\om},
\quad \eta^\tau(x)=B_\tau(f(x))(\eta(\tau^{-1}\cdot x)),
\, \eta\in\Home(F,\D)
\end{equation}
where, notation notwithstanding, the latter formula
in (\ref{eq:groupPlus30}) may not define an 
action of $\pi_1$ on 
$\Home(F,\D)$; more precisely
\begin{equation}\label{eq:groupPlus30Plus}
(\eta^\om)^{\tau}=\mathrm{Inn}_{\b_{\tau,\om}} \z^{\tau\om}
\end{equation} 
from which if
$\d:F\ra \D$ such that
$\mathrm{inf}(f(*))=\mathrm{Inn}_\d f$, 
then on replacing
$(f,\xi)$ by the conjugate pair 
$(\mathrm{inf}(f(*)), \{\om\mpo \d\xi_\om(\d^\om)^{-1}\})$,
the new pair satisfies 
the former equation in (\ref{eq:groupPlus30}),
so, without loss of generality 
we may suppose that $f=\mathrm{inf}(f(*))$.

Now suppose we only have a 
map $(f(*),\xi(*)): 
\rB_\G\ra\rB_\D$ with 
respective actions, then
\begin{equation}\label{eq:groupFix30}
f(*)(\a'_{\tau,\om})\xi'_{(\tau\om)}(*)
=\xi'_\tau(*)(\xi'_\om(*))^\tau \b''_{\tau,\om}
\end{equation}
only holds for $\tau$, $\om$ in $\pi'_1$.
This should first be inflated as a
map of $\Pi''_2$ groupoids $E\uts\rB_\G$
(understood as the inflation of $\rB_\G$)
and $B_\D$. To do this, first consider
\begin{equation}\label{eq:groupFix31}
\xi'':E\ts (\pi''_1)^2\ra \D:
(e,\om)\mpo  B_e(*)^{-1}(\xi'_{\dot{e\tau}}(*) 
\b''_{\dot{e\tau},\overline{e\tau}}
(\b''_{e,\tau})^{-1})
\end{equation}
which for $A'_2:\pi'_2\ra (\text{centre of $\G$})$ verifies
\begin{equation}\label{eq:groupFix32}
(\xi''_\tau (\xi''_\om)^\tau \b''_{\tau,\om} (\xi''_{\tau\om})^{-1})(e)\,=\,
\mathrm{inf}(f(*))(e) (\a'_{\dot{e\tau}, \dot{\overline{e\tau}\om}}
 (A'_2)_* (hK_3)(e,\tau,\om))
\end{equation}
and by (\ref{eq:groupFix2}) the right hand side
of (\ref{eq:groupFix32}) is exactly $\mathrm{inf}(f(*)(\a'')$,
for $\a''$ the inflated associator of the
action of $\Pi''_2$ on $E\times B_\G$. 
To inflate this to a map $\cF\ra\cG$ of $\Pi_2$-groupoids
one just applies the inflation
formula (\ref{eq:groupPlus10})-
which continues to have sense
for non-commutative co-chains-
 to $\xi''$ to get
some $\xi:F\ra \D$
satisfying \ref{eq:groupPlus30}
- there's an $hK_3$ term
in either side of  
\ref{eq:groupPlus30} on similarly inflating
$\a''$ and $\b''$, but they're both central,
so they cancel.

We have, therefore, not only concluded that 
we may, without loss of generality, suppose
that $f=\mathrm{inf}(f(*))$, but that
there
is an object $(\mathrm{inf}(f(*)),\xi)$ in
the left hand side of 
(\ref{eq:groupPlus26}) iff there is an object
$(f(*),\xi'')$ in the right hand side.
At this point, the discussion may be reduced
to our lemmas in group co-homology as follows:
suppose for a given $\om\in\pi_1$ we 
have another solution, $\eta_\om$
of (\ref{eq:groupPlus29}), then
$\eta_\om\xi_\om^{-1}:fA_\om\Rightarrow fA_\om$,
so $\eta_\om(eg)=B_e^{-1}(c_\om(eg))\xi_\om(eg)$ for
some normalised $c_\om\in\Home(F, C)$; where $C$
is the centraliser of $f(\G)$ in $\D$.
Better still by (\ref{eq:groupPlus30})
\begin{equation}\label{eq:groupPlus31}
\tilde{B}:\pi_1\ts F\ra \Aut(C):
(\om,eg) \mpo
B_e(*)
\mathrm{Inn}_{\xi_\om(eg)}B_\om(g)B^{-1}_{\overline{e(g\om)''}}(*) 
\end{equation}
defines a $\Home(F, \Aut(C))$ valued 
co-cycle, and whence an action of
$\pi_1$ on $\Hom(F,C)$
\begin{equation}\label{eq:groupPlus32}
(\om, c)\mpo \, \{x\mpo \tilde{B}_\om(\om\cdot x) (c(\omega^{-1}\cdot x))\}
\end{equation}
satisfying the hypothesis of \ref{lem:groupPlus3};
while given a solution $\xi$ of the
first equation
(\ref{eq:groupPlus30}), we have an isomorphism of sets
\begin{equation}\label{eq:groupPlus33}
\text{co-cycles }\,  \mathrm{Z}^1 (\pi_1\ts F,C) 
\xrightarrow{\sim} \{\text{Solns. of (\ref{eq:groupPlus30})}\}:
c\mpo \{ (eg,\om) \mpo B_e^{-1}(*)(c_\om(eg))\xi_\om(eg) \}
\end{equation}
where the implied $\pi_1$-action on $\Home(F,C)$ is
given by (\ref{eq:groupPlus32}).
Similarly, all morphisms in the categories
envisaged in (\ref{eq:groupPlus26}) are
isomorphisms, which for $f$ fixed, 
and $(c,\xi)\mpo c\cdot\xi$ the
pairing of (\ref{eq:groupPlus33}),
can
be identified with maps $\phi:F\ra C$
such that $\phi:(f,\eta (=c\cdot \xi))\Rightarrow 
(f,\zeta (=d\cdot \xi))$ iff
$c_\om=\phi d_\om (\phi^\om)^{-1}$ 
for the action (\ref{eq:groupPlus32}),
so
that (\ref{eq:groupPlus33}) becomes an
isomorphism
\begin{equation}\label{eq:groupPlus34}
\{\text{objects $(f,\xi)$ with $f$ fixed
modulo isomorphism}\}
\xrightarrow{\sim}
\rH^1(\pi_1, \Home(F,C)) 
\end{equation}
while the automorphisms of 
$(f,\xi)$ are maps $\phi$ such that
for the $\xi$-dependent action (\ref{eq:groupPlus32})
\begin{equation}\label{eq:groupPlus35}
\phi^\om=\phi \quad\text{iff} \quad \phi\in \rH^0(\pi_1, \Home(F,C))
\end{equation}
A similar discussion applies
to the right hand side of (\ref{eq:groupPlus26}),
so choosing $\xi$ in (\ref{eq:groupPlus33})
to be an inflated class, the isomorphism(s)
(\ref{eq:groupPlus34}) and (\ref{eq:groupPlus35})
are compatible with the inflation restriction
isomorphisms of \ref{lem:groupPlus3}. Consequently
the isomorphism classes of objects, and their
automorphisms on either side of
(\ref{eq:groupPlus26}) are mutually inverse
via inflation-restriction; both
sides of
(\ref{eq:groupPlus26})
are groupoids, 
but
not small categories, 
so we conclude by global choice in NBG. 
\end{proof} 
Given the absurd amount of space
that a triviality such as \ref{fact:groupPlus1}
has occupied it's worth making
\begin{summary}\label{sum:group1}
The conclusion of the entire discussion is
that the 2-category of $\Pi_2$-equivariant
(pointed transitive) groupoids is 
(in a way typical of non-abelian co-homology) 
as simple as one could imagine.
Every weak equivalence class of 0-cells maps to a triple
$(\underline{\pi}', \underline{\g'}, \underline{A}')$ consisting of:
a topological 2-type $\underline{\pi}'=(\pi'_1, \pi'_2, k'_3)$, for
$\pi'_1$ a sub-group, respectively a quotient
group $\pi'_2$, of $\pi_1$, respectively $\pi_2$;
the topological 2-type $\underline{\g'}$ of the
universal fibration in $K(\G',1)$'s for some
group $\G'$; and a map $\underline{A}'$ between
these types. Basically, 
\ref{fact:group1} and \ref{def:groupNick},
the weak equivalence class of
0-cells over such a triple ``is'' a principle
homogeneous space under $\rH^2(\pi'_1, \pi_2(\underline{\g'}))$.
Nevertheless there is a subtlety, since this refers
to monoidal functors $\Pi'_2\ra\fAut(\rB_{\G'})$
modulo equivalence, while by (\ref{eq:groupPlus27})
there are more equivalences of $\Pi_2$-equivariant 
groupoids than there are equivalences amongst
monoidal functors. Indeed, in a way which is
typical, cf. \cite[1.5.11]{tom} and
\ref{rmk:corPoint},
of trying to shoe-horn what is really 2-category
theory into 1-category theory by way
of monoidal categories:  
\begin{equation}\label{eq:groupTom}
\text{the monoidal equivalences are
the $\Pi_2$-equivariant equivalences with 
$f=\mathrm{1}$ in (\ref{eq:groupPlus27}).}
\end{equation}
Now, plainly, if one is interested in an
essential surjectivity statement such as
\ref{claim:groupPlus1} then having more
equivalences is the opposite of a problem.
On the other hand, if one wants a nice
linear description of the 0-cells, then
this isn't possible except modulo monoidal
equivalences. More precisely: the $\Pi_2$-equivariant
groupoids modulo equivalence are the quotient
of those modulo monoidal equivalence by
the further action of $\Aut(\G')$ which in the first
instance acts on the 2-types, so the
equivalence classes are given by the action
of sub-groups fixing the 2-type on 
$\rH^2(\pi'_1, \pi_2(\underline{\g'}))$, which
even in a straightforward case case without
a Postnikov invariant or an obstruction class,
so that this is a manifestly linear action,
the action wholly stabilises $0$ and the
quotient is not a linear space.
This is, however, the only subtlety in  
giving a linear algebra description of
the $\Pi_2$-equivariant groupoids. Specifically:
an equivalence class of 1-cells,
determines another triple $(T', T'', \bar{f})$,
where $T'$, respectively $T''$, are the
aforesaid triples of its source, and sink,
which in turn must satisfy $\pi'_1\hookrightarrow\pi''_1$,
$\pi'_2\twoheadrightarrow\pi''_2$; while
$\bar{f}$ is a class in $\rH^1(\G',\G'')$
which takes the $\pi'_1$ action on the
right to the $\pi''_1$ action on the
left, and affords the obvious compatability
conditions on the $\pi_2$'s.
This in turn determines
(up to isomorphism) a group
$C_{\bar{f}}$, to wit: the centraliser
in $\G''$ of the image of $\G'$
under some 
lifting of $\bar{f}$, and the equivalence class of 1-cells 
over $(T', T'', \bar{f})$ is 
either empty, or isomorphic to
$\rH^1(\pi'_1, C_f)$. Here the $\pi'_1$-module
structure $C_f$ 
on the group $C_{\bar{f}}$
depends on a choice of a 1-cell $f$
- or $(f,\xi)$ in the notation
of the proof of 
\ref{claim:groupPlus2}- lying over $(T', T'', \bar{f})$,
and the automorphisms of this 1-cell
are $\rH^0(\pi'_1, C_f)$, which, since
every 2-cell is invertible, is sufficient
to determine everything.
\end{summary}
Unfortunately, we're not quite finished since
we should also add
\begin{scholion}\label{warn:LeftRight} 
(Left vs. right actions)
The various properties of 2-groups
and their actions, such as the so called
pentagon axiom for the associator $K_3$
are most naturally expressed in terms
of left group actions, and whence the
particular choice of definitions of
this chapter. Nevertheless other things,
notably the 2-Galois correspondence, are
most naturally expressed in terms of right
actions. Needless to say this is incredibly
inconvenient. In general terms
cognisant of the objection \ref{faq8}, a right
action of a 2-group, $\Pi_2$, on a category is a
left action of the 2-group $\Pi^{\mathrm{op}}_2$,
{\it i.e.} viewed as a bi-category with 1-object
reverse the direction of the 1-cells but not
the 2-cells, or, again, in terms of
(\ref{eq:group4}) move all the occurrences
of $\gP$ from the left of $\cF$ to the
right. Naturally, therefore, the 
homotopy groups 
of $\Pi^{\mathrm{op}}_2$ are the
opposite group $\pi^{\mathrm{op}}_1$,
and, $\pi_2\op$, {\it i.e.}
$\pi_2$ in 
its natural
left $\pi^{\mathrm{op}}_1$ module structure.
In the usual globular notation
for a 2-category with 1-object,
$*$, the associators change according to 
\begin{equation}\label{eq:LeftRight1}
 \xy
 (0,0)*+{*}="1";
 (30,0)*+{*}="2";
 (60,0)*+{*}="3";
 (90,0)*+{*}="4";
   "1";"2" **\crv{(15,12)} ?(.5)*\dir{>};
   "1";"2" **\crv{(15,-12)} ?(.5)*\dir{>};
   "4";"3" **\crv{(75,12)} ?(.5)*\dir{>};
   "4";"3" **\crv{(75,-12)} ?(.5)*\dir{>};
     (45, 0)*+{\longleftrightarrow}="D";
     (20 , -2 )*+{}="X";
     (20,2 )*+{}="Y";
     (66 , -2 )*+{}="Z";
     (66,2 )*+{}="T";
     (15, 10 )*+{\scriptstyle \s\otimes(\tau\otimes\om) }="F";
     (15,-10 )*+{\scriptstyle (\s\otimes\tau)\otimes\om  }="G";
     (75, 10 )*+{\scriptstyle (\om\op\otimes\tau\op)\otimes\s\op }="H";
     (75,-10 )*+{\scriptstyle \om\op\otimes(\tau\op\otimes\s\op)  }="K";
       {\ar@{=>}^{(\kappa_{\s,\tau,\om})^{-1} } "X";"Y"};
       {\ar@{=>}_{ \kappa\op_{\om\op,\tau\op,\s\op}} "Z";"T"};
\endxy
\end{equation}
where, as usual, qua sets, $\om\op=\om$, {\it etc}.
Nevertheless, as a $\pi_1\op$ group co-chain
with values in $\pi_2\op$,  
$
K^{\mathrm{op}}_3(
\om\op, \tau\op, \s\op)\neq -K_3(\s,\tau,\om)
$, 
in fact this need not even be an opposite co-cycle.
More precisely, if following (\ref{eq:group1}),
we identify the interior of the 2-cell on
the left of (\ref{eq:LeftRight1}) with the
stabiliser of $\s\tau\om$ in $\pi_1\ltimes\pi_2\rras \pi_1$,
then that on the right is the stabiliser
of $(\s\tau\om)\op$ in $\pi_1\op\rtimes\pi_2\rras \pi_1\op$,
where, as (\ref{eq:LeftRight1}) suggests, 
$\pi_2$, respectively $\pi_2\op$, denotes 
not only 
the left action of $\pi_1$, respectively $\pi_1\op$,
but also the right action of
$\pi_1\op$, respectively $\pi_1$. Now there
is a straightforward isomorphism
\begin{equation}\label{eq:LeftRight2}
\pi_1\ltimes\pi_2\xrightarrow{\sim}
\pi_1\rtimes\pi_2\op:(\om, A)\mpo (\om, \om^{-1}._{\text{left action}}A)
\end{equation}
so that in the convention of (\ref{eq:group2}),
on which depends the validity of \ref{fact:group1}, the 
formula for the opposite Postnikov class is
\begin{equation}\label{eq:LeftRight3}
K_3\op:(\pi_1\op)^3\ra\pi_2\op:
(\om\op,\tau\op,\s\op) \mpo
  -(\s\tau\om)^{-1}._{\text{left $\pi_1$ action}}K_3(\s, \tau , \om)
\end{equation}
As such if one employs the obvious maps
\begin{equation}\label{eq:LeftRight4}
\mathrm{op}_1:\pi_1\xrightarrow{\sim}\pi_1\op:
\om\mpo (\om^{-1})\op=(\om\op)^{-1}\quad
\mathrm{op}_2:\pi_2\ra\pi\op_2: A\mpo -A
\end{equation}
then 
$(\mathrm{op}_1)^*K_3\op\neq (\mathrm{op}_2)_* K_3$,
and it's not as obvious as it might be that
(\ref{eq:LeftRight4}) is an isomorphism
of 2-types between $\underline{\pi}$, and
$\underline{\pi}\op$.
The most straightforward way to address this
is to use Maclane's co-herence theorem,
\cite[1.2.15]{tom}, {\it i.e.} choose
an equivalent 2-group, $\gP$, with strict
associativity and strict inverses, or
equivalently, a crossed module. Normally
such strictness is more of a hindrance
than a help, and so far we've only paid
lip service to it in (\ref{eq:group9});
it is however, immediate, that for
a strict 2-groupoid $\gP$
\begin{equation}\label{eq:LeftRight5}
 \xy
 (0,0)*+{x}="1";
 (40,0)*+{y}="2";
 (70,0)*+{x}="3";
 (110,0)*+{y}="4";
   "1";"2" **\crv{(20,20)} ?(.5)*\dir{>};
   "1";"2" **\crv{(20,-20)} ?(.5)*\dir{>};
   "3";"4" **\crv{(90,20)} ?(.5)*\dir{>};
   "3";"4" **\crv{(90,-20)} ?(.5)*\dir{>};
     (55, 0)*+{\build\longleftrightarrow_{}^{\mathrm{op}}}="D";
     (20 , -4 )*+{}="X";
     (20, 4 )*+{}="Y";
     (90 , 6 )*+{}="Z";
     (90, 8 )*+{}="T";
     (90 , -6 )*+{}="M";
     (90, -8 )*+{}="N";
     (20, 12 )*+{\scriptstyle g }="F";
     (20,-12 )*+{\scriptstyle f  }="G";
     (90, 12 )*+{\scriptstyle (g^{-1})\op }="H";
     (90,-12 )*+{\scriptstyle (f^{-1})\op }="K";
       {\ar@{=>}^{\xi } "X";"Y"};
     (90, 3 )*+{\scriptstyle (g^*)^{-1}(f_*)^{-1}(\xi^{-1}) }="L";
(90, -3 )*+{\scriptstyle (f^*)^{-1}(g_*)^{-1}(\xi^{-1}) }="R";
(90, 0 )*+{\scriptstyle \Vert }="S";
       {\ar@{=>}_{} "Z";"T"};
{\ar@{=}_{} "M";"N"};
\endxy
\end{equation}
is a strict isomorphism. Better still if $X:\gP\op\ra\gQ$ is
a morphism of strict 2-groupoids then
\begin{equation}\label{eq:LeftRight6}
\begin{CD}
\gP@>{\sim}>{\mathrm{op}}> \gP\op\\
@V{X\op}VV @VV{X}V\\
\gQ\op@>{\sim}>{\mathrm{op}}> \gQ
\end{CD}
\end{equation} 
strictly commutes. Consequently, if we choose
an equivalence of $\Pi_2$ with some strict
2-group $\gP$, then 
since the automorphism group of a groupoid
is already strict,
the right hand side
of the diagram affords an easy equivalence 
of
2-categories 
\begin{equation}\label{eq:LeftRight7}
\{\text{right $\Pi_2$ groupoids}\}
\xrightarrow{\sim}
\{\text{left $\Pi_2$ groupoids}\}:X\mpo \mathrm{op}^* X
\end{equation}
while the left hand side tells us how to write
it down in practice. Specifically, by
\ref{fact:groupPlus1} and (\ref{eq:LeftRight7})
we only ever have to concern ourselves with
groups $\rB_\G$, on which,
say a 
right action, is given by a 
bunch of functors $R_\om$ with
associators $\rho_{\tau,\om}: R_{(\tau\om)}
=R_{(\tau\om)}\op\Rightarrow R_{\om}R_{\tau}
=R_{\tau}\op R_\om\op$- where the associator
notation is suggested by a right variant of
(\ref{eq:group4}) rather than that for a
left action of $\Pi_2\op$ which would be
$\rho_{\om\op, \tau\op}$, but 
unlike (\ref{eq:LeftRight3}) this is a
purely notational issue. Now since $\rB_\G$
only has one point, the push forwards
in (\ref{eq:LeftRight5}) are all trivial,
so the left action   
given by the composition of the bottom and
left most arrows in (\ref{eq:LeftRight6}) is
\begin{equation}\label{eq:LeftRight8}
L_\om:= R_\om^{-1},\quad 
\lb_{\tau,\om}
:= (R_{(\tau\om)})^{-1}(\rho_{\tau,\om}^{-1})
: L_{(\tau\om)}\Rightarrow L_{\tau}L_{\om}
\end{equation}
As such, there's no specific need to worry
about how $\underline{\pi}$ is isomorphic
to $\underline{\pi}\op$ from the purely
group co-homology perspective. Nevertheless,
it's amusing, so we give the details. In
the first instance, using Maclane's coherence
theorem is akin to working with injective
resolutions, whereas \ref{fact:group1} uses
the standard resolution which
is just C\v{e}ch for the topos
of $\rB_{\pi_1}$. The C\v{e}ch
resolution, however, involves a choice of the
identification of the value of sheaf on a cell
\begin{equation}\label{eq:LeftRight9}
0\xleftarrow{x_1} 1 \xleftarrow{x_2} 2 
\xleftarrow{x_3}
\cdots \xleftarrow{x_{n}} n
\end{equation}
via one of the projections, and this is
conventionally, 
and with no little naturality,
chosen to be the $0$th
projection in the schema (\ref{eq:LeftRight9})
on identifying the sheaf with a
left $\pi_1$-module, {\it i.e.}
the descent datum as an isomorphism 
from source to sink.
If, however, one identifies  
the same sheaf with a right $\pi_1$-module,
{\it i.e.} treat the same descent
datum as an isomorphism from the
sink to the source, then the same
naturality leads one to use the
$n$th projection instead. Now a priori this
isn't a big deal since one just goes from
the complex of left co-chains to that of
right co-chains
by the action of $x_1\cdots x_n$.
It starts getting messy, however, if one 
attempts to identify right co-chains of
$\pi_1$ with left co-chains of $\pi_1\op$
since (\ref{eq:LeftRight9}) becomes
\begin{equation}\label{eq:LeftRight99}
0\xrightarrow{x_1\op} 1 \xrightarrow{x_2\op} 2 
\xrightarrow{x_3\op}
\cdots \xrightarrow{x_{n}\op} n
\end{equation}
so the numbering gets buggered, and to
restore it one has to interchange $i$
with $n-i$. None of this stupidity has
any effect on the topos so if we identify
$\pi_1$ and $\pi_1\op$
by way of 
$\mathrm{op}_1$ of
(\ref{eq:LeftRight4})
we get an automorphism of the standard left co-chain
complex of any left $\pi_1$-module
\begin{equation}\label{eq:LeftRight10}
K\mpo K\cop:=
\{(x_1,\cdots, x_n)\mpo 
(-1)^{(n+1)(n+2)/2}
(x_1\cdots x_n) K(x_n^{-1},\cdots , x_1^{-1})\}
\end{equation}
where, despite appearances, the signs actually
improve by using $\mathrm{co-op}$, rather than
$\mathrm{op}$, which for 2-categories means
invert both 0 and 1-cells, so akin to (\ref{eq:LeftRight4})
$\mathrm{co-op}_1=\mathrm{op}_1$, but 
$\mathrm{co-op}_2=-\mathrm{op}_2=\mathrm{id}$.
Now, evidently, by (\ref{eq:LeftRight3}),
(\ref{eq:LeftRight4})
and (\ref{eq:LeftRight10}),
$(\mathrm{op}_1)^* K_3\op= -K_3\cop$, but
the fact that (\ref{eq:LeftRight10})
is an automorphism does not imply
that the 2-types
$\underline{\pi}$ and $\underline{\pi}\op$
are isomorphic, and, still less, 
a 2-commutative diagram akin to
(\ref{eq:LeftRight6}) in which the
objects are skeletal 2-groups. To do this
one needs a homotopy between $(\mathrm{op}_1)^*$
and $(\mathrm{op}_2)_*$ which commutes with
arbitrary maps of groups and their modules,
so, equivalently, a similarly universal
homotopy between the identity and $\mathrm{co-op}$.
Unsurprisingly such homotopies exist, for example
\begin{equation}\label{eq:LeftRight11}
\pa h + h\pa=\mathbf{1}-(\mathrm{co-op}),\, \text{where:}\,\,
h^n:=\sum_{p=1}^n (-1)^{(n+p+1)(n+p+2)/2} (h_p^n)^*
\end{equation}
and $h_p^n:(\pi_1)^n\ra (\pi_1)^{n+1}$ according to
\begin{equation}\label{eq:LeftRight12}
(x_1,\cdots, x_n)\mpo\begin{cases}
(x_1\cdots x_n, x_n^{-1},\cdots, x_1^{-1}) &\text{if $p=1$},\\
(x_1,\cdots, x_{p-1}, x_p\cdots x_n,  x_n^{-1},\cdots, x_p^{-1}) 
&\text{if $1<p<n$},\\
(x_1,\cdots,x_n, x_n^{-1}) &\text{if $p=n$}
\end{cases}
\end{equation}
Using which, one gets from \ref{fact:group1},
or better 
its source \cite[Theorem 43]{baez}, 
a 2-commutative version of 
(\ref{eq:LeftRight6}), and whence 
another proof of the critical facts (\ref{eq:LeftRight7}) 
and (\ref{eq:LeftRight8}). 
\end{scholion}

\subsection{Fibration and separation}\label{SS:II.2}

Let us begin without any separation hypothesis whatsoever, to wit:

\begin{lem}\label{lem:sep1}
Let $\cY\ra\cX$ be an \'etale map of not necessarily
separated champs then the diagonal $\D_{\cY/\cX}:\cY\ra\cY\times_\cX\cY$
is \'etale.
\end{lem}
\begin{proof}
Let $V\ra \cX$ be an \'etale atlas, put $\cY'=\cY\times_\cX V$,  
and consider the following diagram in which every square is
fibred
\begin{equation}\label{eq:sep1}
\begin{CD}
V@<<< \cY'\times_V\cY'@<<< \cY' \\
@VVV @VVV @VVV\\
\cX@<<< \cY\times_\cX\cY@<{\D_{\cY/\cX}}<< \cY
\end{CD}
\end{equation}
The leftmost vertical is \'etale, whence so are
all the other verticals by base change, and since
the question is local on $\cY$, we may therefore
replace $\cX$ by $V$, and $\cY$ by $\cY'$. If,
however, $U\ra \cY'$ is an atlas for $\cY'$ then
we have a fibre a fibre square
\begin{equation}\label{eq:sep2}
\begin{CD}
U\times_V U @<<{s\times_V t}< R \\
@VVV @VVV\\
\cY'\times_V\cY'@<{\D_{\cY'/V}}<< \cY'
\end{CD}
\end{equation}
for $s,t$ the source and sink of the
groupoid $R\rras U$ representing $\cY'$.
By 
hypothesis $U\ra V$
is \'etale, while
definition the source and sink are
\'etale, so the projection $R\ra V$,
and whence the top
horizontal in (\ref{eq:sep2})
is \'etale.
\end{proof}
Now if $\cY\ra\cX$ is surjective \'etale, then an \'etale
atlas $U\ra\cY$ of $\cY$ also yields an \'etale
atlas of $\cX$, so we can conveniently view
$\cY\ra\cX$ as 
represented by the functor
defined by
the rightmost vertical of the fibre
square
\begin{equation}\label{eq:sep3}
\begin{CD}
\cY @<<{s\times t}< U\times_\cY U= R\\
@V{\D_{\cY/\cX}}VV @VVV\\
\cY\times_\cX\cY@<<< U\times_\cX U=:R_0
\end{CD}
\end{equation}
which by \ref{lem:sep1} is an \'etale map
of groupoids. As such, if we further suppose
that $\cY\ra\cX$ is separated then by definition
the leftmost vertical in (\ref{eq:sep3}) is
proper, so by base change $R\ra R_0$ is too,
whence it's an \'etale covering of its image.
Furthermore, the stabiliser $X\ra U$, respectively
$Y\ra U$ of $\cY$, respectively $\cX$, is the
$U$-group obtained from, say, the fibre square
\begin{equation}\label{eq:sep4}
\begin{CD}
R@<<< X\\
@V{s\times t}VV @VVV\\
U\times U@<{\D_U}<< U
\end{CD}
\end{equation}
so, by yet another base change, $Y\ra X$ is
a proper \'etale map of $U$-groups. Its 
kernel, $Z$, 
is, therefore a finite \'etale covering
of $U$. It's also
a normal
subgroupoid of $R$, so the quotient
$R'=R/Z\rras U$ is an \'etale groupoid in separated spaces-
cf. \cite[7.1-4]{km}- and we assert
\begin{claim}\label{claim:sep1}
({\it cf.} \cite[III.2.1.5]{giraud})
Let $r:\cY\ra\cX$ be a separated surjective \'etale map, then
the classifying champ $\cY':=[U/R']$ affords a
factorisation $\cY\xrightarrow{q}\cY'\xrightarrow{p}\cX$
in separated maps
such that $p$ is a representable \'etale map,
and $q$ is a locally constant gerbe. Furthermore,
such a factorisation into a (not a priori locally
constant) gerbe followed by a (not necessarily
separated) representable \'etale map forming
a 2-commutative triangle
\begin{equation}\label{UniFactor1}
 \xy
 (-18,0)*+{\cY}="L";
 (18,0)*+{\cX}="R";
 (0,16)*+{\cY'}="T";
    {\ar^{q} "L";"T"};
    {\ar^{p} "T";"R"};
    {\ar_{r} "L";"R"};
    {\ar@{=>}^{\eta} (0,2);(0,12)}
 \endxy
\end{equation}
has the following universal property: for any other
such 2-commutative triangle 
\begin{equation}\label{UniFactor2}
 \xy 
 (-18,0)*+{\cY}="L"; 
 (18,0)*+{\cX}="R"; 
 (0,16)*+{\cY''}="T"; 
    {\ar^{y} "L";"T"}; 
    {\ar^{x} "T";"R"}; 
    {\ar_{r} "L";"R"}; 
    {\ar@{=>}^{\xi} (0,2);(0,12)}
 \endxy 
\end{equation} 
one can find a  
dashed arrow, 
$A$,
and fill with $\a$, $\b$ the
other 2-triangles in
\begin{equation}\label{UniFactor3}
 \xy
 (20,-24)*+{\cY'}="A";
 (38,6)*+{\cY}="B";
 (0,0)*+{\cX}="C";
(15,12)*+{\cY''}="D";
{\ar_{}^{q } "B";"A"};
    {\ar^{r } "B";"C"};
    {\ar^{p } "A";"C"};
    {\ar_{x } "D";"C"};
    {\ar_{y } "B";"D"};
{\ar@{-->}_{A} "A";"D"}
{\ar@{=>}_{\b} (10,-8);(14,1)}
{\ar@{=>}^{\a} (24,2);(18,-2)}
\endxy
\end{equation}
such a way that the resulting diagram 2-commutes,
with the following uniqueness property: for any
other filling 
with dashed arrow $A'$, and natural transformations
$\a'$, $\b'$
there is a unique natural transformation
$u:A\Rightarrow A'$ such
that both diagrams resulting from $\a$ and the
filled triangles in \eqref{UniFactor3}
2-commute, {\it i.e.}
\begin{equation}\label{UniFactor4}
(q^*u) \a = \a',\,\,\text{and,}\,\, (x_* u)\b=\b'
\end{equation}
or, what amounts to the same thing, 
the factorisation \eqref{UniFactor1}
is unique 
up to unique equivalence. 
\end{claim}
\begin{proof}
Modulo a notational issue, the top horizontal map in 
(\ref{eq:sep2}) modulo $Z$ is a factorisation of
the same, while the said map modulo $Z$ is $\D_{\cY'/\cX}$
after base change. On the other hand, for any
factorisation $xy$ of a proper map
with $y$ surjective, $x$ is proper, 
so $p$ is separated.
As to its representability:
the fibre of $\cY'$ over $U$ is the 
classifying champ of the
groupoid
\begin{equation}\label{eq:sep5}
R_0\ni f \xleftarrow{\s} R'_t\times_s R_0 \xrightarrow{\tau} f\xi_0 \in R_0
\end{equation}  
for $\xi_0$ the image in $R_0$ of an arrow $\xi\in R'$. In
addition we have a fibre square
\begin{equation}\label{eq:sep6}
\begin{CD}
R'\times R_0 @<<< R'_t\times_s R_0\\
@VVV @VVV\\
R_0\times R_0@<<< R_{0\,t}\times_s R_0
\end{CD}
\end{equation}
where the map on the left is proper since $R\ra R_0$ is
proper, and $R'$ is separated. Consequently the
rightmost map in (\ref{eq:sep6}) is proper. We also
have an isomorphism
\begin{equation}\label{eq:sep7}
C:R_{0\,t}\times_s R_0\xrightarrow{\sim} R_{0\,t}\times_t R_0: (g,f)\mpo (gf,f)
\end{equation}
and a fibre square
\begin{equation}\label{eq:sep8}
\begin{CD}
R_0\times R_0@<<< R_{0\,t}\times_t R_0 \\
@V{t\times t}VV @VVV\\
U\times U @<{\D_U}<< U
\end{CD}
\end{equation}
so the composition of $C$ with the top horizontal in
(\ref{eq:sep6}) is proper, which in turn is
$\s\times\tau$, whence the groupoid of (\ref{eq:sep5}) is
proper. By way of a trivial inspection (\ref{eq:sep5}) is
in fact an equivalence relation, so by 
a topological variant of \cite[1.1]{km},
the classifying champ $[R_0/R'_t\times_s R_0]$ is a 
separated space, and $p$ is representable.

Similarly, the fibre of $q$ over $U$ is the classifying 
champ of 
\begin{equation}\label{eq:sep+1}
R_t\times_s R'\rras R'
\end{equation}
for $R$ acting on the right of $R'$ in the obvious
variation of (\ref{eq:sep5}). In this case, however,
$R\ra R'$ is surjective, so we can slice this along
the identity map $U\ra R'$ to obtain that the
classifying champ of (\ref{eq:sep+1}) is equivalent
to $[U/Z]$. Furthermore, since $Z$ is normal the
isomorphism of $R$-groups
\begin{equation}\label{eq:sep+2}
s^*Z\xrightarrow{\sim} t^*Z: z\mpo fzf^{-1},\quad f\in R
\end{equation}
is a descent datum, whence $Z$ is a locally constant
sheaf on $\cY$. Consequently if $\cY$ is connected
then there is a finite group $ {\vert Z\vert}$ and,
on refining $U$, an isomorphism
\begin{equation}\label{eq:sep+3}
\cY\times_{\cY'} U \xrightarrow{\sim} U\times\rB_{\vert Z\vert}
\end{equation}

Existence established (with $\eta=\mathbf{1}$
by the way)
we turn to the universality
of \eqref{UniFactor1}-\eqref{UniFactor4}. To this
end the first thing to observe is that the 
described property is unchanged on replacing
whether $\cY'$ or $\cY''$ by an equivalent champ,
albeit that this has to be understood in a 2-category
with cells as in \eqref{eq:cor1} rather than the
naive $\underline{\mathrm{Cham}}\mathrm{p}\underline{\mathrm{s}}/\cX$.
As such, by \ref{lem:sep1},
we may without loss of generality suppose
that all champs in \eqref{UniFactor1}-\eqref{UniFactor4}
are the classifying champs of groupoids acting
on the same sufficiently fine \'etale cover $U$,
with the arrows functors, and the 2-cells 
natural transformations between them.
Retaking, therefore, the previous notations
we can express \eqref{UniFactor2} as a factorisation
by way of \'etale maps of groupoids,
\begin{equation}\label{eq:sep10}
R\xrightarrow{y} R'':= U\times_{\cY''} U \xrightarrow{x} R_0
\end{equation}
such that $\xi$ is an equivalence between the
natural projection $r:R\ra R_0$ and $xy$,
with $y$ a gerbe, and $x$ representable.
Consequently, if $Y''$ is the stabaliser of
$R''$ then we have an exact sequence of $U$-groups
\begin{equation}\label{UniSep1}
1\ra Z \ra Y \ra Y''\hookrightarrow X
\end{equation}
The functor $y$ therefore not only factors through $R'$
(as $y'$ say for notational reasons) but
\begin{equation}\label{UniSep2}
\Hom_{R'} (u,v) {\build\rightarrow_{y'}^{\sim}} \Hom_{R''}(yu, yv)
\end{equation}
and similarly $\xi$ factors as $\xi': q\Rightarrow xy'$.
We can, therefore, complete the diagram
\eqref{UniFactor3} with $y'$, the identity
and $\xi$, but to establish the universal
property with these choices of $\cY'$ and
$\cY''$ we also need to go in the other
direction. To this end let $P$ be the 
arrows in $R''$ with sink in the image 
of $y$ then,
\begin{equation}\label{UniSep3}
R_{sy}\ts_t P\rras P: f\xleftarrow{s} (a,f)\xrightarrow{t} y(a)f
\end{equation}
is equivalent to the fibre product 
$\cY\ts_{\cY''} U$ and, by hypothesis,
is a locally constant gerbe over $U$
via the source of $f\in P$ in $R''$.
Consequently for $U$ sufficiently 
fine, $P\ra U$, has a section
\begin{equation}\label{UniSep4}
u\mpo \s(u):=\{ u\xrightarrow{\s_u} y(v)\} \in P\subseteq R''
\end{equation}
As such for $f:u\ra v\in R''$, 
there is by \eqref{UniSep2} a unique
arrow $A(f):u'\ra v'\in R'$ such that
\begin{equation}\label{UniSep5}
\begin{CD}
u@>>f> v \\
@V{\s_u}VV @VV{\s_v}V\\
y'(u') @>{y'(A(f))}>> y'(v')
\end{CD}
\end{equation}
commutes, and whence $f\mpo A(f)$ is a functor.
Better, again by \eqref{UniSep2}, there is a
function $\a:U\ra R'$ such that $y'_* (\a) =(y')^* \s$,
so one can fill the diagram in the other direction
with $A$, $\a$, $(y')^*(\xi')^{-1}x_*(\s)$. 
Finally to get the uniqueness criteria \eqref{UniFactor4}
it suffices to observe that $\a$ affords $\a':\mathbf{1}_{R'}\Rightarrow y'A$,
$\s:\mathbf{1}_{R''}\Rightarrow Ay'$ by construction,
and this is the unique way to make everything 2-commute
since $p$ is injectuve, respectively $q$ is surjective.
\end{proof}
Before progressing let us make
\begin{rmk}\label{rmk:UniFactor} In the 2-category
$\et_2(\cX)$ of \eqref{eq:cor1}, or, more accurately
the obvious variant since our set up here is more
general, the factorisations, \eqref{UniFactor1},
\eqref{UniFactor2} 
may be expressed more cleanly as 1-cells 
$Q:=(q, \eta):r\ra p$,  $Y:=(y, \xi): r\ra x$, 
so that \eqref{UniFactor3} becomes 
\begin{equation}\label{UniFactor5}
 \xy 
 (-18,0)*+{r}="L"; 
 (18,0)*+{p}="R"; 
 (0,16)*+{x}="T"; 
    {\ar^{Y} "L";"T"}; 
    {\ar_{(A,\b)} "R";"T"}; 
    {\ar_{Q} "L";"R"}; 
    {\ar@{=>}^{\a} (0,2);(0,12)}
 \endxy 
\end{equation}
{\it i.e.}, $p$ is a quotient of $r$ unique
up to unique equivalence.
\end{rmk}
As per the beginning of the proof of \ref{claim:sep1}
should we further suppose that $r:\cY\ra\cX$ is proper then
$p$ is proper, and whence
\begin{cor}\label{cor:sep1}
If $r:\cY\ra\cX$ is a proper surjective \'etale map and $\cY$ is connected, then
there is a 
finite group $\vert Z\vert$ and a
unique, up 
to equivalence,
factorisation $\cY\xrightarrow{q}\cY'\xrightarrow{p}\cX$,
where $p$ is a representable \'etale cover, and for
any sufficiently fine \'etale atlas $U\ra \cY'$,
$\cY\times_{\cY'} U \xrightarrow{\sim} U\times\rB_{\vert Z\vert}$.
\end{cor}
Ultimately, we'll suppose  $\cX$ separated. There is,
however, no such hypothesis in any of \ref{lem:sep1},
\ref{claim:sep1}, or \ref{cor:sep1}. As such \ref{cor:sep1}
is what's required for constructing a pro-finite
2-Galois theory of fairly arbitrary champs. This may, however, loose
topological information, so  in general we'll
suppose that $\cY\ra\cX$ is an \'etale fibration 
rather than proper \'etale. Consequently we need a
further lemma  to wit:
\begin{lem}\label{lem:sep2}
Let $p:\cE\ra\cB$ be an \'etale fibration
of not necessarily separated champs and 
$F,G:T\times I\ra\cE$ a pair of maps
together with a natural transformation $\xi_0:F\vert_{T\times 0}
\Rightarrow G\vert_{T\times 0}$ lifting $\eta:pF\Rightarrow pG$ over $0$
then there is a unique natural transformation
$\xi:F\Rightarrow G$ 
lifting $\eta$
such that $\xi\vert_{T\times 0}=\xi_0$.
\end{lem}
\begin{proof}
As in the proof of
\ref{fact:point1} if $\xi$ and $\xi'$ 
extend $\xi_0$ over a subspace of the
form $W\ts I$, then the set, $V$, where 
they coincide is open, closed, and
contains $W\ts 0$. As such every fibre
$V_w$ is a non-emppty connected subset
of $I$, so $V=W\ts I$, and
$\xi$ over $T\ts I$ 
is unique if it exists.

A moderate diagram chase shows, by the universal
property of fibre products, that it is sufficient
to prove the proposition after base changing 
along $pF$, so,
without loss of generality, we can
suppose that $\cB=T\times I$ and that $pF=pG$ is
the identity. As such, $\cE\ra T\times I$ is
surjective, so for a sufficiently fine \'etale
atlas $U\ra \cE$ we can view $F$, $G$ as 
``sections''  of an \'etale functor $P:R\ra R_0$
of groupoids \`a la (\ref{eq:sep3}) albeit that
$R_0\rras U$ is now a separated \'etale equivalence
relation, and, to be precise, $PF$, respectively
$PG$, may be no better than equivalent to the
identity. In any case, since $P$ is \'etale
and every arrow $f\in R$ has a small neighbourhood
$V\ni f$ such that the source and sink are
isomorphisms: 
if we can extend $\xi_0(t)$ to the whole
interval $t\ts I$ then it extends to all
of $V_t\ts I$ for some open neighbourhood
$V_t\ni t$. 
Consequently, by the previous unicity discussion,
such extensions for $s,t\in T$
must glue over $(V_s\cap V_t)\ts I$,
so, without loss of generality, $T$ is a point.
We will therefore be done if we can establish
\begin{claim}\label{claim:sep2}
If $p:\cE\ra I$ is an \'etale fibration of a not necessarily
separated, but connected champ over the interval, then
there is a discrete group $\G$ such that $\cE\xrightarrow{\sim} I\times \rB_\G$
as a champ over $I$.
\end{claim}
Take any two points in the interval, say $0$ and $1$ for
notational convenience, and sets of points $E_0$,
$E_1$ representing the moduli of either fibre, which,
in turn may be identified with $\pi_0(\cE_0)$,
respectively $\pi_0(\cE_1)$. Now given, $x\in\cE_0$,
choose a
section $a_x:I\ra \cE$ with $a(0)$ equivalent to $x$,
so we get a map $A:E_0\ra E_1: x\mpo a_x(1)$, or
more correctly its isomorphism class. Similarly,
we get a map $B:E_1\ra E_0$, and the return map
can be visualised as a diagram
\begin{equation}\label{eq:sep10plus}
x \,{\build\Rightarrow_{}^{\xi}} \,a(0)\ra a(1) \,{\build\Rightarrow_{}^{\eta}}\,
b(1) \ra b(0) \,{\build\Rightarrow_{}^{\zeta}} \,x'
\end{equation}
Although we have no separation hypothesis, whence
no moduli space of paths, we can still form a
wedge $a\vee b:I\vee I\ra \cE$ by gluing (unique up to not
necessarily unique isomorphism) 
the sink of $a$ to the source of $b$ via $\eta$.
Now just as in the proof of \ref{fact:fib1}, the
projection $I\vee I\ra I$ with the wedged point
lying over $1$, homotopes via the wedges of
$[0,s]\vee_s [0,s]$, $s\in (0,1)$ to 
$I\vee I\ra 0$, so by the homotopy lifting
property $x$ of \ref{eq:sep10plus} is, in fact,
connected to $x'$ by a chain of paths and
natural transformations supported in the 
fibre over $0$. Consequently $BA$ is the identity,
and reversing the argument, we deduce that
$\coprod_{x\in E_0} a_x: E_0\times I\ra \cE$
is surjective, and there is no point $e:\rp\Rightarrow\cE$
such that $a_x(p(e))$ is equivalent to
$a_y(p(e))$ for $x\neq y\in E_0$. Now let
$U\ra\cE$ be an \'etale atlas, then $U\ra I$
is \'etale, so $U$ is a disjoint union of
intervals, and for $x\in E_0$ we have a
fibre square of \'etale maps
\begin{equation}\label{eq:sep11}
\begin{CD}
V@>>> U\\
@VVV @VVV\\
I@>{a_x}>> \cE
\end{CD}
\end{equation}
Plainly the top horizontal in (\ref{eq:sep11}) is open,
but it is also closed. Indeed suppose a sequence
of images $\bar{v}_n$ of points $v_n\in V$ converge
to $u\in U$. By the definition of fibre products,
if $u$ is equivalent to a point of the image of
$a_x$ then it is in the image of $V$, so we may
suppose that it's equivalent to some $a_y(p(u))$,
for $x\neq y\in E_0$. The connected component of
$u$ in $U$ is, via $p$, a sub-interval to which
the $\bar{v}_n$ must belong for $n$ sufficiently
large, and so $a_y(p(v_n))$ is equivalent to
$a_x(p(v_n))$, which is nonsense. As we've noted
the image of $V$ is invariant by the groupoid defining
$\cE$, and since $\cE$ is connected, $a_x$ is
not only surjective, but an \'etale atlas, while
$p:\cE\ra I$ expresses $\cE$ as a gerbe over its
moduli. As such the source and sink of 
$I\times_\cE I\ra I$ are the same, and
$\cE$ is the classifying champ $[I/G]$ of
some \'etale $I$-group $G$.

A priori, an $I$-gerbe, or, equivalently, 
the $I$-group $G$ could be complicated.
Nevertheless, every fibre is a $\rB_{G_i}$ for $G_i$
the discrete group which is the fibre of $G$ over $i\in I$.
As such, again, let a pair of points in $I$ be given, and once more
for notational convenience say they're $0$, and $1$.
Now take a closed thickening $J=[-\varepsilon, 1+\varepsilon]$
of $I$, and form the wedge, $T=\rS^1_*\vee_0 J$,
with a pointed circle. We can map this, via $q$ say,
to the unit interval $I$ by first projecting to
$J$ then collapsing $[-\varepsilon,0]$,
respectively $[1,\epsilon]$, to $0$, respectively $1$, and
the identity otherwise. This affords various $T$ points
of $\cE$, to wit: $q^*G$ torsors, and more specifically
the subset obtained by cutting the circle in some point
$t\neq *$, taking the trivial torsor over $T\bsh t$, 
and gluing the ends via $g\in G_0$. An isomorphism
amongst such $q^*G$-torsors is a global section of
$q^*G$, so, by construction, a global section of
$G$, and their isomorphism classes are $G_0$ modulo
the conjugation action of the space $\G$ of global
sections of $G$. To apply the homotopy lifting
property we choose sets, $T_0$, respectively $T_1$,
of torsors representing the conjugacy/isomorphism
classes $G_0/\G$, respectively $G_1/\G$.
For $x\in T_0$, we can homotope $q=q_0$,
by way of
$\rS^1_*\vee_k J$, $k\in I$, together with
appropriate piecewise homotheties of $J$ to
much the same thing, say $q_1:T\ra I$, 
 but now with the roles of $0$ and $1$
reversed,
and
the homotopy lifting property yields 
a map $a_x:T\times I\ra\cE$
whose fibre over $0$ is equivalent to $x$.
By construction, this equivalence must be
given by a global section of $G$, so the
isomorphism class of $a_x(0)$ in $G_0/\G$
is still that of $x$. Furthermore, although
$\rH^1(I,G)$ may be non-trivial, if $f:X\ra I$
is a space over $I$ and $\pi:X\times I$ the
projection, the natural map $\rH^1(X,f^*G)\ra
\rH^1(X\times I, \pi^*f^* G)$ is an
isomorphism, so
the torsors implied by $a_x(k)\vert{T\bsh t}\ra \cE$,
$k\in I$ remain trivial throughout the 
deformation, and we again conclude to a
map $A:G_0/\G\ra G_1/\G:x\mpo a_x(1)$,
and similarly $B:G_1/\G\ra G_0/\G$. To
see that the return map is the identity
is much that same as before, {\it i.e.}
say $y$ is isomorphic to $A(x)$ via
$\g$, form $a_x\vee b_y:T':=T\times (I\vee I)\ra\cE$
by gluing along $\g$ over $1$. As such, we have a
projection of
$T'$ to $J\times(I\vee I)$, and the homotopy, $h_s$, of
$I\vee I\ra I$ to the constant map over $0$ employed
in (\ref{eq:sep10plus}). 
Up to some unimportant behaviour around the
end points of $J$, we may identify it with
$I$, so that inside $J\times(I\vee I)$
we have the wedges of the diagonal $\D\vee\D$,
and we can homotope $J\times(I\vee I)\ra I$
by way of $h_s$ on $\D\vee\D$ and suitable
homotheties of $J$ to $q\times 0:J\times(I\vee I)\ra I$.
Plainly we can arrange that the joins, $*\times (I\vee I)$ of
$\rS^1$ with $J$ map to the same point as
$\D\vee\D$ throughout this homotopy, so that
on applying the homotopy lifting property
we find a homotopy between $x$ and $BA(x)$
through $q^*G$-torsors, and whence, from
what we've already said about the behaviour
of $\rH^1(\underline{\,}, G)$ under homotopy, an
identity of their isomorphism classes as
$q^*G$-torsors, {\it i.e.} $G_0/\G$ by construction.
To conclude: observe that locally around $0$,
$x:T\ra \cE$
is synonymous with a 
$\G$ conjugacy class of sections
$x:V_x\ra G$, where, since
$\G$ is global, the maximal open
connected neighbourhood of $0$, again denoted
$V_x$, where $x$ is defined, depends only on
the conjugacy class. Certainly this set is open,
but it is also closed; since
if $i$ were a limit point, then $A_x(i)$ is, around $i$,
synonymous with a $\G$ 
conjugacy class of
sections $x':V'\ra G$ in a neighbourhood 
$V'\ni i$ which by construction is the same
as the $\G$ conjugacy class of 
$x$ on $V_x\cap V'$, so that $G\xrightarrow{\sim} I\times \G$.
\end{proof}
The immediate application of which is
\begin{cor}\label{cor:sep2}
If $\cY\ra\cX$ is an \'etale fibration of not
necessarily separated champs then the diagonal
$\D_{\cY/\cX}:\cY\ra \cY\times_\cX \cY$ is an
\'etale fibration.
\end{cor}
\begin{proof}
Consider the diagram
\begin{equation}\label{eq:sep12}
\begin{CD}
T @. T\times I @. \\
@V{\tilde{f}_0}VV @V{f}VV @.\\
\cY@>{\D_{\cY/\cX}}>> \cY\times_\cX \cY @>q>> \cY \\
@. @V{p}VV @VVV\\
@. \cY @>>> \cX
\end{CD}
\end{equation} 
together with a given natural transformation
$\eta:f_0\Rightarrow {\D_{\cY/\cX}}( \tilde{f}_0)$.
Put $F=pf$, then by definition of fibre 
products, there is a natural transformation
$p{\D_{\cY/\cX}}( \tilde{f}_0)\Rightarrow \tilde{f}_0$,
which composed with $p\eta$ and inverted yields a natural
transformation $\xi:\tilde{f}_0\Rightarrow F_0$.
As such we get $\a_0={\D_{\cY/\cX}}(\xi)\eta: f_0\Rightarrow 
{\D_{\cY/\cX}}(F_0)$ which by \ref{lem:sep2} is the
restriction to $T\times 0$ of some (unique)
$\a={\D_{\cY/\cX}}(\xi)\eta: f\Rightarrow 
{\D_{\cY/\cX}}(F)$. 
\end{proof}
From which we may deduce a more general variation
of \ref{cor:sep1}
\begin{cor}\label{cor:sep3}
If $r:\cY\ra\cX$ is a surjective \'etale fibration 
of not necessarily separated champs such that
$\cY$ is connected and $\cX$ is locally
1-connected, then
there is a
discrete group $\vert Z\vert$ and a
unique, up to equivalence, 
indeed enjoying exactly the same universal
property as \eqref{UniFactor1}-\eqref{UniFactor4},
factorisation $\cY\xrightarrow{q}\cY'\xrightarrow{p}\cX$,
where $p$ is a representable \'etale cover, and for
any sufficiently fine \'etale atlas $U\ra \cY'$,
$\cY\times_{\cY'} U \xrightarrow{\sim} U\times\rB_{\vert Z\vert}$.
\end{cor}
\begin{proof}
We re-take the notations of \ref{claim:sep1} and thereabouts.
In particular we identify $r$ as the map of groupoids $r':R\ra R_0$
of (\ref{eq:sep3}), with $X/U$, $Y/U$ their respective stabilisers.
Now quite generally if $f:A\ra B$ is an \'etale fibration of spaces
with fibre $F$, 
and $B$
is locally $1$-connected, then every path connected component
of $B$ has a universal cover, so to check that $f$ is a
covering we can suppose that $B$ is simply connected, at which
point the long exact sequence of a fibration implies that
every path connected component of $A$ is isomorphic to $B$.
In particular,
by \ref{cor:sep3} and base change $r'$ is an \'etale covering
of a set of path connected components of $R_0$, so its image
is open and closed. By a further base change,
the kernel $Z$ of $X\ra Y$ is, therefore, 
a locally constant \'etale $U$-group; while
$R\rras U$ is an \'etale groupoid so every arrow $f\in R$
has an open neighbourhood $V\ni f$ such that the source and
sink are isomorphisms from which $V\cap VZ$ is empty, and the
the action of $Z$ on $R$ is discrete. Consequently,
$R'=R/Z$ is a separated space, $R'\rras U$ is an \'etale
groupoid, and, since $R'$ is canonically the image of $r'$,
$R'\ra R_0$ is proper. As such, everything now follows
exactly as in the proof of \ref{claim:sep1}.
\end{proof}
The particular structure of $\cY/\cY'$ merits introducing
\begin{defn}\label{defn:sep1}
If $\cY'$ is connected, then $\cY\ra\cY'$ is said to be
a locally constant gerbe if there is a discrete group
$\G$ such that for some (whence every sufficiently fine)
\'etale atlas $\cY\times_{\cY'} U\xrightarrow{\sim} U\times \rB_{\G}$;
while, in general, $\cY\ra\cY'$ is said to be
a locally constant gerbe if every connected component
is.
\end{defn}
It therfore only remains to clear up a lacuna 
\begin{fact}\label{fact:sep1}
A map $r:\cY\ra \cX$ of (not necessarily separated) 
locally 1-connected
champs is
an \'etale fibration iff it factors as the composition
$\cY\xrightarrow{q}\cY'\xrightarrow{p}\cX$ of a locally
constant gerbe followed by a representable \'etale cover
of a set of path connected components of $\cX$.
\end{fact}
\begin{proof} The only if direction is \ref{cor:sep3}. By
base change, we can suppose that the map we require to
lift is the identity from $T\times I$ to itself; while the
composition of fibrations is a fibration, so we're 
reduced to  finding a section of a locally
constant gerbe $q:\cY\ra T\times I$, given a section over
$T\times 0$. Amongst all \'etale groupoids representing
the interval those corresponding to open covers of $I$
by sub-intervals $I_\a$ 
as encountered in the proof of \ref{lem:cover1} are
co-final, so there is no obstruction to finding a 
section over an interval. As such, we can suppose that
$\cY$ has an \'etale atlas of the form $V\times I$ for
$V\ra T$ an open cover. Consequently $q$ is equivalent
to a map of groupoids $q':R\ra R'\times I$ for 
$R'=V\times_T V\rras V$. By hypothesis $q'$ is an
\'etale cover with fibre $\G$, and a section over
$R'\times 0$, so it has a (groupoid) section everywhere.
\end{proof}

\subsection{The universal 2-cover revisited}\label{SS:II.3}

Let $\cX$ be a separated 
locally path connected champ, then from the 
proof of \ref{cor:sep3}, {\it i.e.} notation of
op. cit. $R'\ra R_0$ is open and closed, every
representable \'etale cover is separated. As such,
\ref{SS:II.2} adds nothing to what has already
been said about the fundamental group of $\cX$,
\ref{SS:I.5}, or 1-Galois theory, \ref{SS:I.7}.
We can, however, employ \ref{cor:sep3} to extend
our knowledge of path spaces, to wit

\begin{fact}\label{fact:visit1}
If $\cY_*\ra\cX_*$ is a (weakly) pointed \'etale fibration,
then there is a lifting $\tilde{Q}:\rP\cX_*\times I\ra  \cY$ of
the universal (pointed) path $Q:\rP\cX_*\times I\ra  \cX$
together with a natural transformation $\tilde{\z}:*\Rightarrow \tilde{Q}(*)$
such that the pair $(\tilde{Q}, \tilde{\z})$ finely represents
pointed paths on $\cY$, {\it i.e.} satisfies the
universal property (\ref{eq:pointuni1}).
\end{fact}
\begin{proof}
The proposition is trivial for representable maps,
so by \ref{cor:sep3} we can, without loss
of generality, suppose that $\cY_*$ is
a strictly pointed locally constant
gerbe over $\cX_*$ with fibre
$\rB_\G$. To construct $\tilde{Q}$ one first
applies the homotopy lifting property to the
deformation retract (\ref{eq:point5}), to get
a lift of $Q_0$, and then applies it to the
universal map itself. 
To 
identify $\tilde{\z}$ and
check the universal property (\ref{eq:pointuni1}),
let a family $(\tilde{F}, \tilde{\xi}):I\times T\ra \cY$ of pointed
intervals be given, with $({F}, {\xi}):I\times T\ra \cX$
it's projection to $\cX$ and $G:T\ra \rP\cX_*$,
$\eta:F\Rightarrow \mathrm{id}\times QG$ as per op. cit. 
Now consider the 2-commutative fibre 
square
\begin{equation}\label{eq:visit1}
\begin{CD}
\cY@<<< \cY' \\
@VVV @VVV \\
\cX @<{\mathrm{id}\times Q}<< I\times\rP\cX_* 
\end{CD}
\end{equation}
The two stage construction of $\tilde{Q}$, and
the proof of \ref{fact:sep1} applied at each 
stage imply that $\cY'\xrightarrow{\sim} (P\times I)\times\rB_\G$.
As such, we may suppose that this isomorphism is
an equality, so that if $\a$ is the natural transformation
between the two sides of the  square
 (\ref{eq:visit1}), then we may identify $\tilde{Q}$
with the natural projection from $P\times I$, and
whence, by the definition of fibre products,
there is a natural transformation $\tilde{\z}:*\Rightarrow \tilde{Q}(*)$
such that $\a_*\pi(\tilde{\z})_*=\z_*$. 
Consequently 
the pair $(\pi\tilde{Q},\pi\tilde{\z})$ also has
the universal property for pointed paths to $\cX$,
so without loss of generality, $(\tilde{Q},\tilde{\z})$
is a strict lifting of $(Q,\z)$, and we require
to prove that the natural transformation,
$\eta$ can be lifted to 
$\tilde{\eta}: \tilde{F}\Rightarrow \mathrm{id}\times \tilde{Q}G$,
given a lifting $\tilde{\eta}_*$ over $*(=0)\times T$;
which as it happens is \ref{lem:sep2}, albeit that
the difficult bit, \ref{claim:sep2}, is a given in
our current situation.
\end{proof}
As soon as we have a pointed path space, we also
have a loop space, 
whence
\begin{cor}\label{cor:visitPlus1}
If $\cX_*\ra\cX'_*$ is a (weakly) pointed \'etale
fibration over a separated champ then for any
pointed space $Y_*$, $\Hom_*(\Sigma Y_*, \cX_*)$ 
is finely represented in the sense of \ref{fact:point1}
by the separated space $\Hom_*(Y_*,\Omega\cX_*)$.
Similarly, $\Hom_*(\rS^1_*, \cX_*)$ is
finely represented by the separated space $\Omega\cX_*$;
while the spaces $M'\cX$, and $M\cX$ are well
defined by the formulae (\ref{eq:cover2}) and
(\ref{eq:cover3}).
In particular the homotopy groups of $\cX$ are
well defined by any of the equivalent formulae
(\ref{eq:point6}). 
\end{cor} 
\begin{proof}
By \ref{fact:visit1}, the path space $\rP\cX_*$
is well defined, so the loop space
is defined exactly as in (\ref{eq:point3}). 
The right hand side of
(\ref{eq:point4});
the left hand side of (\ref{eq:pointsuspend}), along with
the formulae (\ref{eq:cover2}) and (\ref{eq:cover3})
are, therefore, well defined separated spaces irrespectively
of whether $\cX$ is separated or not. It is however,
a formal consequences of the universal property of
the right hand side of (\ref{eq:point4}), respectively
the left hand side of (\ref{eq:pointsuspend}), that it finely
represents the left hand side, respectively the
right hand side, cf. \ref{fact:tear2}.  
\end{proof}
Now in going from the definition of the
space $M$ of homotopies between pointed
paths to the explicit identification,
\ref{claim:cover1}, of the connected
component of the 
identity of the
groupoid $\rP\cX_*\ts_{\cX} \rP\cX_*$,
and whence to the universal 1 and 2 covers,
the separation of $\cX$ is used, but
only in the weaker form
\begin{fact}\label{fact:visitPlus1}
An \'etale fibration $\cX\ra\cX'$ over a
locally path connected
separated champ enjoys the following
separation property: if $a,b:I\ts T \ra \cX$
are maps from the product of a closed interval
with a space $T$, and
$\xi:a\vert [0,1)\ts T \Rightarrow b\vert [0,1)\ts T$
a natural transformation between their 
restrictions to $[0,1)\ts T$ then there
exists a unique natural transformation
$\bar{\xi}:a\Rightarrow b$ which
restricts to $\xi$.
\end{fact}
\begin{proof}
We may suppose that $\cX$ 
is 
path connected, so 
for any base point $*:\rp\ra\cX$, there
exists a path space
$\rP\cX_*$, or better a pair $(\rP\cX_*,\z)$, 
by \ref{fact:visit1}.
The isomorphisms
between $a$ and $b$ come from an \'etale
groupoid, whence the existence of such a
$\xi$ is an open condition, 
so as in the argument immediately
priori to \ref{claim:sep2},
we may suppose that $T$ is a point.
Now choose 
$*$ to be say $a(0)$, then there are
points, $A,B\in\rP\cX_*$ together with
natural transformations
$\a_1:A\Rightarrow a$,
and $\b_1:B\Rightarrow b$, where by
(\ref{eq:pointuni1})  $\a_1$, respectively $\b_1$
is the unique natural transformation
such that $\a_1(0)=\z_A^{-1}$,
respectively $\b_1(0)=\xi(0)\z_B^{-1}$. 
As ever the moduli, $A_t$, of $x\mpo A(tx)$,
affords not just a path in $\rP\cX_*$ from $*$
to $A$, but also natural transformations
$\a_t:A_t\Rightarrow a_t$ uniquely determined
by the initial condition $\a_t(0)=\z_{A_t}^{-1}$ at $*$,
and similarly for $B_t$. 
Now for $t<1$, $A_t$ and $B_t$ are 
equal so $\xi=\b\a^{-1}\vert [0,1)$;
while $\rP\cX_*$ is separated because
$\cX'$ is, so $A=B$, and whence $\bar{\xi}=\b\a^{-1}$.
\end{proof}
The utility of \ref{fact:visitPlus1} is not limited
to constructing the universal 1 and 2 covers, but,
in fact
\begin{fact}\label{fact:visitPlus2}
All of \ref{SS:I.4}-\ref{SS:I.7} hold not just for
separated champs but also for \'etale fibrations
over the same, similarly for such champs
the weak separation condition 
of \ref{claim:path2}.(c) holds.
\end{fact}
\begin{proof}
Apart from the representability theorems \ref{cor:visitPlus1},
the only other thing that we need is separation in
the weak sense of \ref{fact:visitPlus1}.
\end{proof}
There are other ways of deducing 
\ref{fact:visitPlus2}, at least
in the interval \ref{SS:I.5}-\ref{SS:I.7},
from the 2-Galois correspondence
for separated champs. This would,
however, be slightly to miss the
point which is that the definition
of separation, 
unlike the weak or path separation
of \ref{fact:visitPlus1},
is excessive, {\it e.g.}
$\rB_\G$ is separated iff $\G$ is finite,
whereas \ref{fact:visitPlus1} holds iff
$\G$ is discrete. In the same vein
\begin{fact}\label{fact:visit2}
If $\cY_*\ra\cY'_*$ is a locally constant
pointed gerbe with fibre $\rB_\G$ of champs
whose path spaces exist, then
$\Omega\cY_* \ra \Omega\cY'_*$ is a fibration
with fibre $\G$.
\end{fact}
\begin{proof}
Indeed we have
a diagram of fibre squares
\begin{equation}\label{eq:visit2}
\begin{CD}
\rP\cY_* @<<< \Omega\cY'_* @<<< \Omega\cY_*\\
@VVV @VVV @VVV\\
\cY@<<< \rB_\G @<<< *\\
@VVV @VVV @. \\
\cY' @<<< * @.
\end{CD}
\end{equation}
so the arrow in question is a base change of
the fibration $\rp\ra\rB_\G$.
\end{proof}
In similar generality, once the path fibration exists
we have a diagram of fibre squares
\begin{equation}\label{eq:visit3}
\begin{CD}
P\cY_* @<<< P\cY_*\times_\cY P\cY_*@<<< \Omega\cY_*\\
@VVV @VVV @VVV\\
\cY @<<< P\cY_* @<<< *
\end{CD}
\end{equation}
and similarly for $\cY'$, so that the map of groupoids
deduced from the fibre square
\begin{equation}\label{eq:visit4}
\begin{CD}
\cY @<<<  R:=\rP\cY_*\times_\cY \rP\cY_* \\
@V{\D_{\cY/\cY'}}VV @VVV\\
\cY\times_{\cY'}\cY @<<<  R':=\rP\cY'_*\times_{\cY'} \rP\cY'_*
\end{CD}
\end{equation}
is by \ref{cor:sep2} and \ref{cor:sep3} an \'etale 
fibration weakly equivalent to $\Omega\cY_* \ra \Omega\cY'_*$.
Unsurprisingly, therefore
\begin{fact}\label{fact:visit3}
If $\cX'$, is  a separated, locally 1-connected, and
semi-locally 2-connected champ with 
$\cX_*\ra\cX'_*$ a weakly pointed \'etale fibration 
of connected champs
then $\cX$ is also locally 1-connected, and
semi-locally 2-connected. Furthermore: 
their universal 2-covers
$\cX'_2$, $\cX_2$ are 2-connected and isomorphic.
\end{fact}
\begin{proof}
By the long exact sequence of a fibration
and the factorisation lemma \ref{fact:sep1}
$\pi_2(\cX_*)\hookrightarrow \pi_2(\cX'_*)$
is injective, which proves the first part.
The coverings $R_2\ra R_1$, respectively
$R'_2\ra R'_1$, 
of \ref{fact:cover3}, applied to $\cX$, and $\cX'$
respectively,  are just
$\cX_2=\cY$, $\cX_1=\cY'$, respectively $\cX'_2=\cY$,
$\cX'_1=\cY'$
in (\ref{eq:visit4}), which proves the 2-connectedness.
In either case, these coverings are unchanged by
representable \'etale covers, so, again by the factorisation
\ref{fact:sep1}, we can suppose that $\cX\ra\cX'$ is
a locally constant gerbe with fibre $\rB_\G$, whence
by (\ref{eq:visit4}) applied to $\cX=\cY$, and $\cX'=\cY'$,
we get an \'etale fibration $R\ra R'$ with fibre $\G$
such that $R'$ is the fine quotient $R/\G$ for $\G$
identified to a normal sub-group $\rP\cX_*\ts\G$
of the stabiliser. Now plainly, the connected 
components $R_1$, $R'_1$ of the respective identities
must map to each other as \'etale coverings, so
we certainly get a unique homeomorphism $R_2\ra R'_2$
of their universal coverings which sends the identity
to the identity. Now if this map were not an isomorphism
of groupoids, we'd have distinct groupoid
structures on $R_2$ with the same identity  
lifting that on $R_1$. The difference, however, between
such structures is a continuous map $R_{2,\,t}\ts_s R_2\ra\pi_2(\cX)$
which is $\mathbf{1}$ on the identity. 
Whether the source or the sink of $R_2\ra P\cX_*$
is a fibration, whose fibre, $\Om\cX_{2,*}$,
and base are path connected,
whence $R_{2,\,t}\ts_s R_2$ is connected, so
the groupoid
structure on $R_2$ is unique.
\end{proof}

Now 
for things as in \ref{fact:visit2}, then,
as encountered in the above proof, and
just as in \ref{claim:sep1} and \ref{cor:sep3}
but with $U$ replaced by $P=\rP\cY_*$,
the map $R\ra R'$ of (\ref{eq:visit4}) is the quotient by
the kernel of the stabiliser group $\cY\ra \cY'$,
which,
since $P$ is contractible, is the trivial $P$-group $\G$.
On the other hand if $R_1$, $R'_1$ are
the connected components of the identity, then
for any $\g\in\G$ we have the fibre
square
\begin{equation}\label{eq:visit101}
\begin{CD}
R_1\cap (\g R_1\g^{-1}) @>>> \g R_1\g^{-1}\\ 
@VVV @VVV\\
R_1 @>>> R
\end{CD}
\end{equation}
where the right vertical is open and closed, so
the left is as well. Consequently, 
\begin{fact}\label{fact:visit4}
As a $\G$-torsor,
the fibre $R\ts_{R'} R'_1$ over the connected component of
the identity is given by a representation $\rho_Z:\pi_2(\cY')\ra Z$
in the centre $Z$ of $\G$, albeit, for the moment, we
leave open which such representations actually occur.
In any case, 
the conjugation action of $\G$ on the centre is trivial, so
such torsors are classified by a 
subset of $\Hom(\pi_2(\cY'), Z)$; while the particular
case of $\cY=\cX_2$, $\cY'=\cX_1$ implies that
$\pi_2(\cX_*)$ is abelian.
\end{fact}
\begin{proof} 
By \eqref{eq:visit101}, we have 
for each $\g\in \G$ a map 
$f\mpo \g f\g^{-1}$ of the \'etale 
covering $R_1\ra R'_1$, which fixes the
identity arrows, so by the 1-Galois
correspondence, \ref{fact:one1},
all such maps are the identity. As
such the left and right actions
of $\G$ on $R_1$ coincide, so the
1-Galois correspondence implies 
that they arise from a representation
of $\pi_1(R'_1)=\pi_2(\cY'_*)$ in
the centre of $\G$. 
\end{proof}
At which point we come to the universal property of $\cX_2$,
envisaged in \ref{defn:coverAddon},
to wit
\begin{fact}\label{fact:visit5}
If $\cX\ra\cX'$ is an \'etale fibration 
of connected champs with
a separated locally 1-connected,
and semi-locally 2-connected base, then
for any (weakly) pointed
\'etale fibration
$q:\cY_* \ra \cX_*$ from a champ
there is a  pointed map
$r:\cX_{2,*}\ra \cY_*$, {\it i.e.} a pair $(r,\rho)$
where $\rho:*\Rightarrow r(*)$,
and a natural transformation, $\eta:p\Rightarrow qr$
such that,
\begin{equation}\label{eq:visit5}
 \xy
 (0,0)*+{qr(*_{\cX_2})}="A";
 (20,0)*+{q(*_{\cY})}="B";
 (20,-18)*+{*_{\cX}}="C";
 (0,-18)*+{p(*_{\cX_2})}="D";
    {\ar@{=>}_{q(\rho)} "B";"A"};
    {\ar@{=>}_{\eta_*} "D";"A"};
    {\ar@{=>}^{} "C";"B"};
    {\ar@{=>}_{} "C";"D"};
 \endxy
\end{equation}
commutes.
Given $\rho$, $\eta$ is unique; while conversely
given $\eta$ the set of possible $\rho$'s is a
principal homogeneous space under the 
group $\G$ of \ref{fact:visit4} (of the connected
component of $*_\cY$). Similarly,
if $(r',\rho',\eta')$
is any other  triple such that
(\ref{eq:visit5}) commutes, then there is a
unique natural transformation
$\xi:r\Rightarrow r'$
satisfying $\xi_*\rho=\rho'$, and this requirement
implies $\eta'=q(\xi)\eta$.
\end{fact}
\begin{proof}
Let $\cY\ra\cY'\ra\cX$ be the factorisation of $q$ into
a locally constant gerbe followed by a representable
cover guaranteed by \ref{cor:sep3}. This factorisation
is functorial, so to give a map $r:\cX_2\ra\cY$ is the
same thing as giving a pair of 2-commuting maps
\begin{equation}\label{eq:visit6}
\begin{CD}
\cX_2 @>>r> \cY \\
@VVV @VVV \\
\cX_1 @>>r_1> \cY'
\end{CD}
\end{equation}
while geometric fibres don't change under gerbes, so
the question of where the base points go (which depends
only on the representability of $r_1$) is covered
by \ref{fact:cover2}, and we may, therefore, consider $r_1$ as 
a given. Re-taking the notation of \ref{SS:I.5}, at
the level of groupoids we therefore have
\begin{equation}\label{eq:visit7}
\begin{CD}
R_2 @. R \\
@VVV @VVV \\
R_1 @>>r_1> R'
\end{CD}
\end{equation}
where the horizontal arrow identifies $R_1$ with
the connected component of the identity, and $R_2$
is the universal cover of $R_1$. Now, a map of
groupoids must send identities to identities,
so there is at most one possible candidate for
lifting $r_1$, {\it i.e.} the quotient of $R_2$
by the kernel, $K$, of the central representation $\rho_Z$,
and since $K\times P$ is a normal 
(by way of the embedding in the stabiliser) sub-groupoid,
this works, {\it i.e.} 
$R_2\xrightarrow{r} R_2/K\hookrightarrow R \ra R'$ are maps of groupoids.

Now let us turn to the base points. As in the
proof of \ref{fact:cover2}, we can replace 
the base points $*_\cY$, $*_{\cX_2}$ by
equivalent ones, and this has already been
done implicitly in the above. Similarly on
fixing isomorphisms $P\times_{\cY'} \cY\xrightarrow{\sim} P\times\rB_\G$,
and $P\times_{\cX_1} \cX_2\xrightarrow{\sim} P\times\rB_{\pi_2}$-
albeit the latter is somewhat build into
the definitions- we have identifications of
$*_{\cX_2}$, respectively 
$*_\cY$,
with their projections to $\cX_1$, respectively
$\cY'$. Having done this, there is, as we've said,
only one way to define $r$, so, without such 
choices there is only one way to define $r$
up to equivalence. In addition with such 
choices
$\eta$ may
be supposed to be the 
identity, and in any case, up to any
2-commutativity in the factorisation
\ref{cor:sep3}, it's always the
pull-back along the left
hand vertical of op. cit. of the $\eta$ appearing
in (\ref{eq:cover12}). As such, the unicity
of $\eta$ given $\rho$ follows from \ref{fact:cover2},
while the unicity, or otherwise, of $\rho$
is the failure of the stabiliser
of $*_\cY$ to inject via $q$ into the stabiliser
of $q(*_\cY)$, {\it i.e.} the group $\G$ of
\ref{fact:visit4}. That this is equally the
indeterminacy in $r$ follows by \ref{lem:sep2},
or an inspection of (\ref{eq:visit7}) 
in light of \ref{fact:visit4}; while $\eta'=q(\xi)\eta$
follows from the forgetfulness of $q$ and \ref{fact:cover2}.
\end{proof}

\subsection{The action of 
\texorpdfstring{$\Pi_2$}{Pi\_2}}\label{SS:II.4}
For brevity, we'll suppose throughout this section that $\cX$
is a connected, locally 1-connected,
and semi-locally 2-connected champ, 
which admits an \'etale fibration over a
separated champ. We retake
the notations of \ref{SS:I.5}, so that for 
a base point $*':\rp\ra\cX_1$ in the fibre
over $*:\rp\ra\cX$ and
$\om\in\Omega$ we have a triangle
\begin{equation}\label{eq:act1}
 \xy
 (0,0)*+{\pi_2(\cX_{1,*'})}="A";
 (28,0)*+{\pi_2(\cX_{1,\om(*')}) }="B";
 (14,-12)*+{\pi_2(\cX_{*})}="C";
{\ar_{}^{\om} "A";"B"};
    {\ar^{q} "B";"C"};
    {\ar^{p} "C";"A"};
\endxy
\end{equation}
of natural isomorphisms,
and whence the usual (left) action of $\pi_1(\cX_*)$ on
$\pi_2(\cX_*)$ by way of 
$\om\mpo p\om q\in \mathrm{Aut}(\pi_2(\cX_{*}))$. 
Alternatively, since the action of $\om$ is given
by the functor $F_\om$ of (\ref{eq:coverplus1}),
one finds that the composition 
\begin{equation}\label{eq:actplusplus}
\Om_1\xrightarrow{i} R_1 \xrightarrow{F_\om} R_1 \xrightarrow{r} \Om_1
\end{equation}
for $r$ the retract of the inclusion
$i$ of the loop space occurring in
the proof of \ref{lem:cover2} and 
$\Om_1$ the connected component 
of $(*,*)$ (equivalently
$\Om\cX_{1,*'}$) 
that in terms of concatenation of 
paths in the sense of
\ref{fact:point3}, the action deduced from (\ref{eq:act1}) is
the usual conjugation action
\begin{equation}\label{eq:actplusplus1}
\Om_1\ra \Om_1: a\mpo \om^{-1} a\om
\end{equation}
which passes unambiguously to an
automorphism of $\pi_1(\Om_1)$
since this is abelian by \ref{fact:visit4}. The formula
(\ref{eq:actplusplus1}) notwithstanding,
this is a left action because we
define multiplication in $\pi_1$
by left concatenation, {\it i.e.}
the traditional notation \ref{fact:point3}
is mis-leading.
In any case,
we therefore have a skeletal groupoid $\pi_1\ltimes\pi_2\rras \pi_1$,
and even the monoidal product, (\ref{eq:group1}).
As such, although we still haven't identified the
Postnikov class $k_3$, let use denote this
category by $\Pi_2$, and observe
\begin{fact}\label{fact:act1}
There is a lifting of the fibre square (\ref{eq:coverFix2})
to a fibre square
\begin{equation}\label{eq:actFix2}
 \xy
 (0,0)*+{\Pi_2(\cX_*)\uts \cX_2}="A";
 (38,0)*+{\cX_2 }="B";
 (38,-18)*+{\cX}="C";
 (0,-18)*+{\cX_2 }="D";
  (23,-5)*+{}="E";
  (15,-13)*+{}="F";
    {\ar^{q=\om\mpo F_\om} "A";"B"};
    {\ar_{\mathrm{trvial}\,\mathrm{ projection}} "A";"D"};
    {\ar^{p} "B";"C"};
    {\ar_{p} "D";"C"};
    {\ar@{=>}_{i} "F";"E"};
 \endxy
\end{equation}
where $i$ is just the pull-back to $\cX_2$ of
the natural transformation $i$ of op. cit.,
while $F_\om$ are liftings to 
(weak) automorphisms of 
$\cX_2$ of the functors
of (\ref{eq:coverplus1}) defining $q$ 
(the action) in
(\ref{eq:coverFix2}), and we have an
associativity diagram
\begin{equation}\label{eq:actplus}
 \xy
 (0,0)*+{\Pi_2(\cX_*)\uts \Pi_2(\cX_*)\uts \cX_2}="A";
 (38,0)*+{\Pi_2(\cX_*)\uts \cX_2 }="B";
 (38,-18)*+{\cX_2}="C";
 (0,-18)*+{\Pi_2(\cX_*)\uts \cX_2 }="D";
  (23,-5)*+{}="E";
  (15,-13)*+{}="F";
    {\ar_{\quad\quad\otimes\ts\mathrm{id}} "A";"B"};
    {\ar_{\mathrm{id}\times q} "A";"D"};
    {\ar^{q} "B";"C"};
    {\ar_{q} "D";"C"};
    {\ar@{=>}_{\alpha} "F";"E"};
 \endxy
\end{equation}
for some lifting $\a$ of the natural
transformation of (\ref{eq:coverplus}),
satisfying the co-cycle condition
\begin{equation}\label{eq:actFix1}
p(\alpha_{\tau,\om})=i_{(\tau\om)}i_\om^{-1} F_\om^* i_\tau^{-1}
= i_{(\tau\om)}i_\tau^{-1} (F_\tau)_* i_\om^{-1}
\end{equation}
\end{fact}
\begin{proof}
As in the proof of \ref{factdef:cover1}, 
\eqref{eq:cover9}-\eqref{eq:cover10},
the
fibre product $\cX_2\ts_\cX \cX_2$ is represented by the groupoid
\begin{equation}\label{eq:act2}
R_{2\, t}\ts_s R_{0\, t}\ts_s R_2\rras R_0
\end{equation}
with the natural action via the projection of $R_2\ra R_0$.
Every orbit meets some $i_\omega(P)$, 
(\ref{eq:cover6}),
so we get an
equivalent category by slicing along
\begin{equation}\label{eq:act3}
\coprod i_\om: \coprod_{\om\in\Omega} P \ra R_0: a\mpo i_\om(a)
\end{equation}
The arrows between paths $a$, and $b$ in the $\om$th copy
of $P$, is therefore a pair of arrows $A, B\in R_2$ such
that $i_\om(b)= \bar{B}i_\om(a)\bar{A}^{-1}$,
where $\bar{\,}$ denotes 
the image of $R_2$ in $R_1\hookrightarrow R_0$.
As such $A$ is an arrow from $a$ to $b$, and 
$\bar{B}=i_\om(b) \bar{A} j_\om(a)$. By \ref{lem:cover2}
the map $A\mpo\bar{B}$ can certainly be lifted,
but it's a little more useful to note an explicit
lifting as a functor, {\it i.e.} concatenate a square
such as (\ref{eq:cover3}) with a square every horizontal
cross section of which is $\om$, to
obtain a functor from $R_2\ra R_2$, which 
we write $A\mpo F_\om(A)$.
Consequently,
the $\om$th slice of (\ref{eq:act2}) along
(\ref{eq:act3}) is of the form $R_2\times \pi_2(\cX_*)$,
and it remains to check that the multiplication
in $R_2\times \pi_2(\cX_*)$ is the obvious one.
To this end observe that an element $S$ of $\pi_2(\cX_*)$
acts on $A\in R_2$, cf. (\ref{eq:cover3}), by way of
\begin{equation}\label{eq:act5}
 \xy
 (0,0)*+{}="A";
 (18,0)*+{}="B";
 (18,-18)*+{}="C";
 (0,-18)*+{}="D";
  (-20,-9)*+{*\,\, {\build\Longrightarrow_{}^{\z}} }="F";
(11,-9)*+{A}="G";
(0,-9)*+{S}="L";
(9,2)*+{b}="H";
(9,-20)*+{a}="I";
(22, -16)*+{}="J";
(22, -2)*+{}="K";
    (0,-9)*\xycircle(6,9){-}="g";
    {\ar@{-}_{} "A";"B"};
    {\ar@{-}^{} "B";"C"};
    {\ar@{-}_{} "D";"C"};
    {\ar@{=>}_{\xi(y),\, \xi(0)=\mathbf{1}} "J";"K"};
 \endxy
\end{equation}
so that as required we get a lift of 
the diagram (\ref{eq:coverFix3}) to
a fibre square
\begin{equation}\label{eq:actFix3}
 \xy
 (0,0)*+{\coprod_{\om\in\Omega} \pi_2(\cX_*)\ts R_2}="A";
 (60,0)*+{R_2}="B";
 (60,-18)*+{R_0}="C";
 (0,-18)*+{R_2 }="D";
  (34,-5)*+{}="E";
  (26,-13)*+{}="F";
    {\ar^{(\om, S, A)\mpo S F_\om (A) } "A";"B"};
    {\ar_{(\om, S, A)\mpo A} "A";"D"};
    {\ar^{\mathrm{natural}} "B";"C"};
    {\ar_{\mathrm{natural}} "D";"C"};
    {\ar@{=>}_{a\mpo i_\om(a)} "F";"E"};
 \endxy
\end{equation}
which is equivalent to the diagram (\ref{eq:actFix2})
lifting (\ref{eq:coverFix2}). 
To lift the diagram (\ref{eq:coverplus}), with
$\bar{\,}$ image in $R_1$
where appropriate, 
observe that as
maps from $P\ra R_1$ the natural transformations
of (\ref{eq:coverFix11})- albeit now written
$\bar{\alpha}_{\tau,\om}$- satisfy, as a consequence
of the co-cycle condition (\ref{eq:coverFix1})
\begin{equation}\label{eq:act6}
D(\bar{\alpha})(a):=
\bar{\alpha}_{\sigma, (\tau\om)}(a) \bar{F}_{\sigma}(\alpha_{\tau,\om}(a))
[\bar{\alpha}_{(\sigma\tau),\om}(a)\bar{\alpha}_{\sigma,\tau}(a\om)]^{-1}
=\mathbf{1}_{(\sigma\tau\om)a}
\end{equation}
Now, as the notation suggests let $\alpha_{\tau,\om}:P\ra R_2$
be a lifting of $\bar{\alpha}_{\tau,\om}$, then
\begin{equation}\label{eq:act7}
R_2\ra R_2: f\mpo \alpha_{\tau,\om}(s(f))^{-1}F_{(\tau\om)}(f)^{-1}
\alpha_{\tau,\om}(t(f)) F_\tau(F_\om(f))
\end{equation}
takes values in $P\times \pi_2(\cX_*)$, so, in fact,
its the identity, {\it i.e.} any lifting 
$\alpha_{\tau,\om}$ of $\bar{\alpha}_{\tau,\om}$
is a natural transformation between 
$F_\tau F_\om$ and $F_{\tau\om}$, whence 
(\ref{eq:actplus}) and (\ref{eq:actFix1}).
\end{proof}
Before proceeding let us tie this up with
the conjugation action (\ref{eq:act1}) via
\begin{rmk}\label{rmk:act+1}
On identifying $\pi_2=\pi_2(\cX_*)$ 
with the stabiliser of $*$ in $R_2$,
there is a
unique identification of the stabiliser, $S$, of $R_2$
with $P\ts\pi_2$, and whence a map 
$\om\mpo F_\om \vert_S$ from $\Omega$ to $\Aut(\pi_2)$. 
To verify that this is an
action- 
and indeed the action (\ref{eq:act1})-
amounts to checking that the inner
maps $\mathrm{Inn}_{\a_{\tau,\om}}$ are trivial
on $S$. Now these are maps from $P$
to $\Aut(\pi_2)$ so they're certainly constant.
The arrow $\a_{\tau,\om}(*)$ usually
fails to be $\mathbf{1}$, but it is a lifting of
the arrow $\bar{\a}_{\tau,\om}(*)$, and $R_1$
is path connected, so there's a map $\bar{f}:I\ra R_1$
with $\bar{f}(1)$, $\bar{f}(0)$ respectively 
$\bar{\a}_{\tau,\om}(*)$, and the identity on
its source. Consequently there's a lifting
$f:I\ra R_2$ which at $1$ is
${\a}_{\tau,\om}(*)$, 
and an element of $\pi_2$ at $0$,
whence by \ref{fact:visit4} we get a homotopy
between the identity and $\mathrm{Inn}_{\a_{\tau,\om}}$,
so the latter inner functor is indeed $\mathbf{1}$.
\end{rmk} 
Furthermore, as a consequence of (\ref{eq:act6})
we have
\begin{factdef}\label{factdef:act1} 
Notations as above, then $P\ra R_2: a\mpo D(\alpha)(a)$ 
takes values in $P\times\pi_2(\cX_*)$, and so defines
a co-cycle $-K_3(\alpha):\pi_1(\cX_*)^3\ra \pi_2(\cX_*)$,
with the corresponding Postnikov class $k_3(\cX_*)
\in\rH^3(\pi_1(\cX_*), \pi_2(\cX_*))$
being independent of the liftings $\alpha$, and the choices
whether of the sections $i_\om$ of (\ref{eq:cover6})
or the complete repetition free list of homotopy
classes $\Om$. The triple, $(\pi_1(\cX_*), \pi_2(\cX_*), k_3(\cX_*))$,
or just $(\pi_1,\pi_2, k_3)$ if there is no danger
of confusion, with the action (\ref{eq:act1}) is
the topological 2-type of $\cX$, and $\Pi_2(\cX_*)$
in the sense of (\ref{eq:group1}), and (\ref{eq:group2}),
with unicity as per \eqref{eq:group3}-\eqref{eq:groupPlus},
is the homotopy 2-group of $\cX_*$.
\end{factdef}
\begin{proof}
The fact that 
given the sections, and the homotopy classes,
changing the lifting of $\bar{\alpha}$ amounts
to changing $K_3(\alpha)$ by a co-boundary is
clear. Otherwise, labelling by $\pi_1$ rather
than $\Om$, the effect of changing either
homotopy classes or sections is the same,
{\it i.e.} we get another bunch of sections,
$I_\om:P\ra R_0$,
of the source with
values in the same connected component as $i_\om$. 
As such we get functors $\bar{G}_\om:R_1\ra R_1:A\mpo I_\om A I_\om^{-1}$,
which are related to the previous $\bar{F}_\om$ by
$\bar{\d}_\om: \bar{G}_\om\Rightarrow \bar{F}_\om$ for
$\bar{\d}_\om= i_\om (I_\om)^{-1}$, which is a priori
$R_0$ valued, but by \ref{lem:cover2} is actually
a map $\bar{\d}:P\ra R_1$, which since $P$ is
contractible lifts to a map $\d_\om:P\ra R_2$, and
whence a lifting $G_\om: =(\d_\om) F_\om (\d_\om)^{-1}$
of $\bar{G}_\om$ to an endomorphism of $R_2$.
Plainly, replacing $i$ by $I$ in the formula
(\ref{eq:coverFix11}) defines natural transformations
$\bar{\beta}_{\tau,\om}:\bar{G}_\tau\bar{G}_\om\Rightarrow \bar{G}_{(\tau\om)}$,
which, in turn, satisfy 
\begin{equation}\label{eq:act9}
\bar{\beta}_{\tau,\om}=
\bar{\d}_{\tau\om}^{-1}\bar{\alpha}_{\tau,\om}
\bar{\d}_{\tau}\bar{G}_{\tau}(\bar{\d}_\om)
\end{equation}
A lifting of everything on the right hand side of (\ref{eq:act9}),
has been defined, so we can just take this as a definition
of a lifting $\beta_{\tau,\om}$ of the left hand side. As
such (or perhaps, modulo patience and a large piece of paper),
$D(\beta)(\sigma, \tau, \om) = \d_{(\sigma\tau\om)} D(\alpha)(\sigma,\tau,\om) 
\d_{(\sigma\tau\om)}^{-1}$. Now as in
\ref{rmk:act+1}: for any arrow 
$f:a\ra b$ in $R_2$ we can realise it as the
end point of a path of arrows $f_t:a\ra b_t$,
with $f_0=\mathbf{1}_a$ since $R_2$ is path
connected, so the conjugation action of $f=f_1$
on the kernel of $R_2\ra R_1$ is a path of conjugation 
actions on $\pi_2$ starting at the identity, whence
$\pi_2$ is central in
$R_2$ and $D(\beta)=D(\alpha)$,
or if one prefers $D(\beta)=D(\alpha)D(S)$ for $D(S)$
a co-boundary had we chosen a lifting of $\bar{\beta}$ another
than that suggested by (\ref{eq:act9}). 
\end{proof} 
All of which may be applied to conclude
\begin{fact}\label{fact:act2}
The 
second projection, $q$, of 
(\ref{eq:actFix2}) naturally affords
a 
not necessarily faithful, cf. \ref{cor:one1},
action of $\Pi_2$ on $\cX_2$. More precisely, in
the notation of the proof of op. cit., we have a
mapping of 2-groups $F:\Pi_2\ra \mathrm{Aut}(R_2)$ given
by $\om\mpo F_\om$ on $1$-cells, $S\mpo \{F_\om {\build\Rightarrow_{}^{S}} F_\om$\}
on $2$-cells (modulo the canonical identification of $\pi_2\times P$ with the
kernel of $R_2\ra R_1$), together with a natural
transformation $\alpha: F(\otimes) \Rightarrow \otimes(F\times F)$-
deduced from (\ref{eq:actFix2}) and (\ref{eq:actplus}).
\end{fact} 
\begin{proof} By (\ref{eq:actFix3}), the second projection of 
(\ref{eq:actFix2})
is exactly the
bi-functor described above. For convenience we
can suppose that $*$ represents the identity in
our (repetition free) set of representatives $\Om$
of $\pi_1$, so, without loss of generality, the identity 1-cell 
maps to the identity of $\mathrm{Aut}(R_2)$, and
since
$D(\alpha)=K_3$ we get a map of 2-groups by
construction.
\end{proof}
A useful, and, as it happens, necessary alternative
description of the Postnikov invariant is
\begin{fact}\label{fact:act3}
There is a map of groupoids $R\ra R_0$ 
with $R_2$ the fibre of $R$ over $R_1$ iff $k_3=0$.
\end{fact} 
\begin{proof} Suppose such a $R$ exists. The 
kernel, $P\times\pi_2$, of the stabiliser of
$R_2\ra R_1$ acts, by say left composition, on $R$,
and this action is transitive, {\it i.e.} $R\ra R_0$
must be a $\pi_2$-torsor, whence
we can lift the sections $i_\om:P\ra R_0$ of the
source to sections, $\tilde{i}_\om$ of $R$, 
so we may suppose $F_\om$ is the conjugation
$A\mpo\tilde{i}_\om A\tilde{i}_\om^{-1}$,
and $\alpha_{\tau,\om}$ is given by (\ref{eq:coverFix11}),
albeit now with values in $R_2$ rather than
$R_1$. The vanishing, (\ref{eq:act6}), of $D(\bar{\alpha})$
is, however, a formal consequence of (\ref{eq:coverFix1})
so $k_3=0$. Conversely, suppose that $k_3=0$, and
in the notation of (\ref{eq:coverFix11}) let
$R_{\om}$ be the $\pi_2$ torsor $R_\om:=(c_{1,\om}^{-1})^*t_\om^* R_2$,
and observe that the functors $F_\om:R_1\ra R_1$
of \ref{fact:act2} factor through $R_{1,\om}$
of (\ref{eq:coverFix}) which for convenience
we identify 
to a sub-space of $R_1$ via $t_\om$ of op. cit.,
and write $\tilde{c}_{1,\om}: R_{1,\om}\ts_{R_1} R_2\ra R_\om$
for the natural isomorphism. 
Now for $f,g$ arrows in $R_0$ belonging to
connected components $R_0^\tau$, and $R_0^\om$
respectively we can write $f=Ai_\tau(s(f))$,
$g=Bi_\om(s(g))$ for some unique $A\in R_{1,\tau}$,
respectively unique $B\in R_{1,\om}$. Consequently,
if $f$ and $g$ are compossible
\begin{equation}\label{eq:act10}
fg = A\bar{F}_\tau(B) i_{\tau}i_\om = A\bar{F}_\tau(B) \bar{\alpha}_{\tau,\om}^{-1} i_{(\tau\om)}
\end{equation}
where the arrow $A\bar{F}_\tau(B) \bar{\alpha}_{\tau,\om}^{-1}$
necessarily belongs to $R_{1,(\tau\om)}$. As such,
via $c_{1,(\tau\om)}$ and the definition of $R_{(\tau\om)}$
we get a lifting of the product in $R_0$ via
\begin{equation}\label{eq:act11}
R_{\tau\, s} \ts_t R_{\om}\ra R_{(\tau\om)}:
(y,x)\mpo y\cdot x := \tilde{c}_{1,(\tau\om)} [(\tilde{c}_{1,\tau})^{-1}(y)
F_\tau(\tilde{c}_{1,\om}^{-1}(x)) \alpha_{\tau,\om}^{-1} ]
\end{equation}
Supposing, as we may that $*$ belongs to our set
$\Om$ of representatives of $\pi_1$, conveniently
implies that $R_{1,*}$ is $R_1$, and without
loss of generality, all of  $\alpha_{\tau, *}$,
$\alpha_{*,\om}$, $F_*$ are identities. As such,
$R:=\coprod_{\om\in\Om} R_\om$ is a groupoid,
\'etale fibring over $R_0$ with fibre $\pi_2$,
iff the product (\ref{eq:act11}) is associative.
Unsurprisingly, however,
\begin{equation}\label{eq:act12}
(z\cdot y)\cdot x = z\cdot (y\cdot x) D(\alpha)
\end{equation}
and since group co-homology can be computed
with normalised co-chains we therefore get
an associative product without loosing the
existence of identities in $R$ iff $k_3=0$.
\end{proof}
The only if direction in \ref{fact:act3} is 
wholly non-constructive, and one can do
better, {\it viz}:
\begin{rmk}\label{rmk:act1}
If there 
any groupoid $R\ra R_0$ with fibre $R_2$ over $R_1$,
then up to strict isomorphism 
of categories over $R_0$ it
is necessarily the groupoid $R$
constructed in the course of the proof of \ref{fact:act3}
modulo the possibility of multiplying the
product in (\ref{eq:act11}) by a normalised
group co-cycle. In particular, 
when $k_3=0$
the isomorphism
classes 
(even weakly provided $R_2$ and an isomorphism of the
kernel of $R_2\ra R_1$ with $P\ts\pi_2$ are fixed)
of such $R_0$-groupoids is a principal 
homogeneous space under $\rH^2(\pi_1,\pi_2)$.
\end{rmk}
\begin{proof} We know from the beginning of the
proof of \ref{fact:act3} that each fibre $R^\om$
over $R_0^\om$ must be a $\pi_2$ torsor. 
Now
consider a
based loop, $g:I\ra R_0^\om$. 
To such, there is a unique based
loop $f:I\ra R_{1,\om}$ such that
$g=fi_\om(s(g))$. Now let $\tilde{f}$,
respectively $\tilde{g}$ be any 
liftings to $R_2$, respectively
$R^\om$, then $\tilde{f}^{-1}$,
$\tilde{g}$ are still compossible,
and $\tilde{f}^{-1}\tilde{g}$
lifts $i_\om(s(g))$.
This
latter loop is, however, contractible,
so $\tilde{g}(1)\tilde{g}(0)^{-1}=\tilde{f}(1)\tilde{f}(0)^{-1}$,
and whence $R^\om$ is the torsor $R_\om$
encountered in the proof of \ref{fact:act3}.
The product in $R$ must lift the
product (\ref{eq:act10}) by hypothesis,
so it's necessarily (\ref{eq:act11}) up
to multiplication by a 
normalised
co-chain $S_{\tau,\om}:\pi_1\times \pi_1\ra\pi_2$,
which must in fact be a co-cycle by (\ref{eq:act12}).
Certainly any normalised 1 co-chain defines 
not just a 
functor between such $R_0$ groupoids, but even a 
strict isomorphism, and since the fibre over $R_1$
is fixed any functor is a map of $\pi_2$-torsors,
{\it i.e.} a normalised 1 co-chain, so, a fortiori,
any weak equivalence over $R_0$ fixing $R_2$ is
a strict isomorphism.
\end{proof}  

\subsection{The 2-Galois correspondence}\label{SS:II.5}

The following is a routine exercise
in the definitions.

\begin{fact}\label{fact:cor1}
Let $\cC$ be a 2-category (be it weak or strict) and
$T$ a $0$-cell then,

(a) The sub 2-category, 
$\fAut(T)$,
with $0$-cells $T$, and
$T$-alone, $1$-cells $\Aut_1(T)$, i.e. weakly invertible
morphisms in $\Hom_1(T,T)$, and $2$-cells invertible
2-morphisms in $\cC$ between elements of $\Aut_1(T)$
is a 2-group.

(b) For any $0$-cell $X$ of $\cC$, the category
$\uHom_{\cC} (T,X)$ 
with objects $\Hom_1(T,X)$, and arrows 
2-morphisms between them,
inherits a right $\fAut(T)$ action from $\cC$.

(c) There is a 2-functor,  
$\fHom(T,\_)$,
from $\cC$ to the
2-category of categories with right $\fAut(T)$ action defined
by: a $0$-cell $X$ maps to $\uHom_{\cC}(T,X)$, a $1$-cell $f:X\ra Y$
maps to an $\fAut(T)$-equivariant functor $f_T: \uHom_{\cC}(T,X)\ra 
\uHom_{\cC}(T,Y)$,
and a $2$-cells $\xi:f\Rightarrow g$ maps to an
$\fAut(T)$-equivariant natural transformation
$\xi_T:f_T\Rightarrow g_T$.  
\end{fact}

The particular (strict) 2-category
to which we intend to apply
\ref{fact:cor1} is given by

\begin{defn}\label{defn:cor1} Let $\cX$ be a connected separated
locally 1-connected and semi-locally 2-connected
champ, then the 2-category $\et_2(\cX)$ is the
2-category whose $0$-cells are \'etale fibrations
$q:\cY\ra\cX$, $1$-cells pairs $(f,\eta)$ of a 1-morphism
and a natural transformation $\eta$ forming 
a 2-commutative triangle \`a la 
(\ref{eq:one1}), and $2$-cells natural transformations
$\z$ forming the 2-commutative diagram
\begin{equation}\label{eq:cor1}
 \xy
 (0,0)*+{\cY}="1";
 (36,0)*+{\cY'}="2";
 (18,-20)*+{\cX}="3";
   "1";"2" **\crv{(7,10)} ?(.97)*\dir{>};
   "1";"2" **\crv{(7,-10)} ?(.97)*\dir{>};
     (14, -12)*+{}="D";
     (22,-12)*+{}="C";
     (12 , -9 )*+{}="B";
     (24,-9 )*+{}="E";
     (9,-2 )*+{\scriptstyle f}="F";
     (29, 4 )*+{\scriptstyle g}="G";
     (26,2 )*+{}="H";
       {\ar@{=>}_{\eta} "D";"C"};
       {\ar@{:>}^{\xi} "B";"E"};
       {\ar^{q'} "2";"3"};
       {\ar_{q} "1";"3"};
       {\ar@{=>}^{\zeta} "F";"H"};
\endxy
\end{equation}
\end{defn}
By definition all 2-morphisms in $\et_2(\cX)$ are invertible,
and, in general, we have the simplification that we're interested
in the category of
groupoids,  $\Grpd$, or better $\Pi_2(=\Pi_2(\cX_*))$
right equivariant groupoids, $\Grpd(\Pi_2)$. This said,
we have the 2-Galois correspondence
\begin{prop}\label{prop:cor1}
Let $\cX_*$ be a  
a pointed champ with the prescriptions of
\ref{defn:cor1}, and $\cX_2\ra\cX$
its universal 2-cover, then 
$\fAut_{\et_2(\cX)}(\cX_2)$ is equivalent
to $\Pi_2$ and
the 2-functor,
\begin{equation}\label{eq:cor2}
\fHom(\cX_2,\_):\et_2(\cX)\ra\Grpd(\Pi_2)
\end{equation}
is an equivalence of 2-categories. Better still,
for $q:\cY\ra\cX$ a $0$-cell in $\et_2(\cX)$,
the $\Pi_2$-groupoid
$\uHom_{\et_2(\cX)}(\cX_2,\cY)$ is equivalent
in $\Grpd(\Pi_2)$
to the fibre $q^{-1}(*)$.
\end{prop} 
We establish the proposition in bite size pieces.
To this end, notice that a representative
of the fibre $q^{-1}(*)$ is the category with
objects pairs $(*_\cY, \phi)$ where 
$*_\cY:\rp\ra\cY$,
$\phi:*\Rightarrow q(*_\cY)$, and arrows natural
transformations $\g:*_\cY\Rightarrow *'_\cY$ 
such that $\phi'=q(\g)\phi$. In particular, this
affords a convenient way to prove
\begin{fact}\label{fact:cor2}
Let $q:\cY\ra\cX$ be a $0$-cell in $\et_2(\cX)$,
and $\theta:p(*_2)\Rightarrow *$ a fixed natural
transformation between the projection of a 
base point of $\cX_2$ with that on $\cX$, 
then 
the functor given on objects by
\begin{equation}\label{eq:cor3}
\uHom_{\et_2(\cX)}(\cX_2,\cY)\ra q^{-1}(*):
(r,\eta)\mpo (r(*_2), *{\build\Rightarrow_{}^{\eta_*\theta}} p(*_2))
\end{equation}
and sending a natural transformation 
$\xi$ to its value on $*_2$
is an equivalence of
categories.
\end{fact} 
\begin{proof} That the above is a functor is
automatic from preceding description of $q^{-1}(*)$
and \ref{defn:cor1}. It's plainly full 
and essentially surjective
by
\ref{fact:visit5}; while to see that it's
faithful let $\xi$ be an arrow from
$(r,\eta)$ to $(r',\eta')$ in $\et_2(\cX)$,
and in the notation of \ref{fact:visit5}
take $*_\cY=r(*_2)$, $\rho=\mathbf{1}$,
$\rho'=\xi_{*_2}$ and apply the uniqueness
statement of op. cit.
\end{proof}
One can reason
in the same vein to tidy up a lacuna
\begin{rmk}\label{rmk:cor1}
A connected weakly pointed \'etale fibration $\cX_{2,*}\ra\cX_*$ with the
universal property of \ref{fact:visit5} is unique
up to equivalence.
\end{rmk}
In any case, this establishes that as
categories $\fAut_{\et_2(\cX)}(\cX_2)$
and $\Pi_2$ agree, but, of course
\begin{fact}\label{fact:cor3}
The functor $F$ of \ref{fact:act2} 
together with the natural transformation
$i$ of (\ref{eq:coverFix2})
defines a
(weak) inverse, 
$\om\mpo (F_\om, i_\om)$ on objects,
to the functor of (\ref{eq:cor3})
in the particular case of $\cY=\cX_2=[\rP\cX_*/R_2]$,
which combined with the natural transformation
$\alpha$ of \ref{fact:act2}- 
cf. (\ref{eq:coverplus}), (\ref{eq:coverFix1}), and \ref{cor:one1}-
defines an equivalence of
2-groups
\begin{equation}\label{eq:cor4}
F: \Pi_2\ra \fAut_{\et_2(\cX)}(\cX_2)
\end{equation} 
\end{fact}
\begin{proof}
The fact that (\ref{eq:cor4}) is an equivalence
of categories is a formal consequence of
\ref{fact:cor2} and (\ref{eq:actFix2}), {\it i.e.}
same trick as (\ref{eq:oneFix2}); while the
further fact that it's even an equivalence
of 2-groups is a consequence of
(\ref{eq:actplus}), (\ref{eq:actFix1}), and
the definition of the Postnikov class
in \ref{factdef:act1}.
\end{proof}
Consequently, we define the right action
(or, equivalently a left action by (\ref{eq:LeftRight7}))
of $\Pi_2$ on a fibre by (\ref{eq:cor4}),
and \ref{fact:cor1}.(b), so the part of
the better still in \ref{prop:cor1}
beyond \ref{fact:cor2} is tautologous,
albeit we'll have occasion to explicate
this tautology in the proof of \ref{fact:cor5}. 
Slightly less trivially
\begin{fact}\label{fact:cor4}
The $2$-functor (\ref{eq:cor2}) 
restricted to $1$-cells
is a family of fully faithful functors. More
precisely if $F=(f,\eta)$, $G=(g.\xi)$ are $1$-cells
from $q$ to $q'$ in $\et_2(\cX)$ such that for
some natural transformation $\z_*$ in $\Grpd(\Pi_2)$
\begin{equation}\label{eq:cor5}
 \xy
 (0,0)*+{q^{-1}(*)}="1";
 (36,0)*+{}="2";
 (41,0)*+{{{q'}}^{-1}(*)}="22";
 (18,-20)*+{*}="3";
   "1";"2" **\crv{(7,10)} ?(.97)*\dir{>};
   "1";"2" **\crv{(7,-10)} ?(.97)*\dir{>};
     (14, -12)*+{}="D";
     (22,-12)*+{}="C";
     (12 , -9 )*+{}="B";
     (24,-9 )*+{}="E";
     (9,-2 )*+{\scriptstyle f_*}="F";
     (29, 4 )*+{\scriptstyle g_*}="G";
     (24,2 )*+{}="H";
       {\ar@{=>}_{\eta_*} "D";"C"};
       {\ar@{:>}^{\xi_*} "B";"E"};
       {\ar^{q'} "2";"3"};
       {\ar_{q} "1";"3"};
       {\ar@{=>}^{\zeta_*} "F";"H"};
\endxy
\end{equation}
then there is a unique 
$2$-cell in $\et_2(\cX)$,
$\z: (f,\eta)\Rightarrow(g.\xi)$
such that (\ref{eq:cor1}) 2-commutes,
and restricts to (\ref{eq:cor5})
over the base point.
\end{fact}
\begin{proof}
By  \ref{cor:sep3} there are 
functorial
factorisations
$\cY\xrightarrow{\tilde{q}}\cY_1\xrightarrow{q_1}\cX$,
$\cY\xrightarrow{\tilde{q}'}\cY_1\xrightarrow{q'_1}\cX$
of $q$, respectively $q'$, into a locally
constant gerbe followed by a representable
\'etale cover, so we may suppose that
we have a 2-commutative diagram
\begin{equation}\label{eq:cor6}
 \xy
(-20,0)*+{\cY }="X";
 (-20,-18)*+{\cY'}="Y";
(0,0)*+{\cY_1}="A";
 (20,-9)*+{\cX}="C";
 (0,-18)*+{\cY_1' }="D";
  (10,-5)*+{}="E";
  (4,-14)*+{}="F";
  (-5,-5)*+{}="G";
  (-15,-13)*+{}="H";
    {\ar_{f_1} "A";"D"};
    {\ar^{q_1} "A";"C"};
    {\ar_{q_1'} "D";"C"};
    {\ar^{\tilde{q}} "X";"A"};
    {\ar_{\tilde{q}'} "Y";"D"};
    {\ar^{f} "X";"Y"};
    {\ar@{=>}_{\xi_1} "E";"F"};
    {\ar@{=>}_{\tilde{\xi}} "G";"H"};
 \endxy
\end{equation}
with $\xi=\xi_1q_1'(\tilde{\xi})$ affording $(f,\xi)$,
and similarly for $(g,\eta)$. By \ref{fact:one1}, and
the ubiquitous connectivity argument of
\ref{fact:point1}, our assertion holds if
the $0$-cells are
$q_1$ and $q_1'$, so- \ref{fact:point1} again- it holds
if the $0$-cells are $q$ and $q_1'$, and finally it
holds for $q$ and $q'$ by \ref{lem:sep2} and
path connectedness. 
\end{proof}
Given  \ref{fact:cor4}, 
then again by
\cite[1.5.13]{tom} 
and global choice in NBG- 
cf. end of the
proof of \ref{claim:groupPlus2}-
we will have established \ref{prop:cor1}
if we can prove
\begin{fact}\label{fact:cor5}
The $2$-functor (\ref{eq:cor2}) is essentially
surjective on  
$0$ and $1$-cells.
\end{fact}
\begin{proof}
Consider first the $0$-cells. By \ref{fact:groupPlus1},
or, equivalently \ref{claim:groupPlus1} in the
specific, it will suffice to prove that given
the action of a sub (in the sense of \ref{fact:groupPlus1}) 
2-group $\Pi'_2$ on a group $\rB_\G$ that there
is an \'etale fibration $\cX'\ra\cX$ such that
the restriction to $\Pi'_2$ of the action of $\Pi_2$ via
\ref{fact:cor3} and the 2-functor
 (\ref{eq:cor2}) on $\rB_\G$ is the given one.
To this end we can, therefore, suppose without
loss of generality that $\Pi_2=\Pi'_2$; and we
identify the left action with 
(normalised) co-chains $A:\pi_1\ra \Aut(\G)$,
$\z:\pi_1^2\ra \G$ such that $A_{\tau\om}=\mathrm{Inn}_{\z_{\tau,\om}}
A_\tau A_\om$ and  $D(\z)=(A_2)_* K_3$
for some $\pi_1$-representation $A_2:\pi_2\ra Z$
in the centre of $\G$. Now in the first place,
$R_2\uts\rB_\G$ is certainly a groupoid, and
its stabiliser has a normal sub $P$-group 
\begin{equation}\label{eq:cor7}
P\ts\pi_2\hookrightarrow P\ts\pi_2\ts\G: S\mpo (S, A_2(S))
\end{equation}
so the quotient, $R'_2$, of $R_2\uts\rB_\G$ by
this sub-group is a groupoid fitting into a
diagram
\begin{equation}\label{eq:cor8}
 \xy
 (0,0)*+{R_2}="A";
 (28,0)*+{R'_2}="B";
 (14,-12)*+{R_1}="C";
{\ar_{}^{r'} "A";"B"};
    {\ar^{q'} "B";"C"};
    {\ar_{p'} "A";"C"};
\endxy
\end{equation}
in which $q'$ is naturally, by way
of the action of the stabiliser, $P\ts\G$,  of
$R'_2$, a $\G$ torsor on the left and right.
Similarly, 
for $\om\in\pi_1$,
we dispose of functors $F'_\om:=F_\om\uts A_\om:R_2\uts \rB_\G\ra
R_2\uts \rB_\G$ and natural transformations
$\a'_{\tau,\om}:= \a^{-1}_{\tau,\om}\ts\zeta_{\tau,\om}:
F'_{(\tau\om)}\ra F'_\tau F'_\om$ which
respect the kernel (\ref{eq:cor7}) of $R_2\uts \rB_\G\ra R'_2$,
so these descend to a series (denoted with the
same letters) of functors and
natural transformations on $R'_2$.
In addition 
\eqref{eq:cover6}-\eqref{eq:coverFix} provides a 
series of homotopies between $R_1$, and
the connected components $R^\om_0$ of $R_0$,
so we can 
(notation as per op. cit.) 
find a series of left-right $\G$-torsors
$R'_\om$ and (far from unique) isomorphisms of 
left-right $\G$-torsors
\begin{equation}\label{eq:cor9}
\begin{CD}
R'_{2,\om}:= R'\vert_{R_{1,\om}}@>{\sim}>{c_\om}> R'_\om\\
@VVV @VVV \\
R_{1,\om}@>{\sim}>{c_{1,\om}}> R^\om_0
\end{CD}
\end{equation} 
where, plainly, we take $c_{\mathbf{1}}$ to be
the identity. At which point we define a 
product on $R':=\coprod_\om R'_\om$ by way
of the formula
\begin{equation}\label{eq:cor10}
R'_\tau\ts R'_\om\ra R'_{(\tau\om)}: (g,f)\mpo
c_{(\tau\om)} \bigl( c^{-1}_\tau (g) F'_\tau (c^{-1}_\om(f)) \a'_{\tau,\om}
\bigr)
\end{equation}
the associativity of which is
equivalent to $D(\a')=0$, 
for $D$ as in (\ref{eq:act6}).
By construction, however, $D(\a')$ 
is the
image in $R'_2$ of $K_3\ts\mathrm{obs}_\G$
viewed as an element of the stabiliser of
$R_2\ts\rB_\G$, so this is zero by 
\eqref{eq:cor7} and \ref{fact:group3}.
Identifying $P\ts\G$ with the stabiliser
of $R'_2$, and whence as a subset of $R'$,
one then has that the left product structure
of $P\ts\G$ on $R'_\om$ defined by 
(\ref{eq:cor10}) is exactly that of
the given left $\G$-torsor; while
the right product structure is the 
composition with $A_\om$ of the given
right torsor. As such, the product
(\ref{eq:cor10}) admits left and right
identities, so it's a category fibred
over a groupoid, $R_0$, with fibres
a group, $\G$, whence $R'\rras P$ is
a groupoid. 

To show that this is what we started with,
observe that the sections $i_\om$ of (\ref{eq:cover6})
can be lifted to sections $i'_\om:P\ra R'_\om\sbs R'$,
while exactly as in the proof of 
\ref{fact:act3} a connectedness argument
shows that $F'_\om$ and the conjugation
$x\mpo (i'_\om)x(i'_\om)^{-1}$ agree
on the connected component of the identity,
{\it i.e.} the image of $R_2$ in $R'_2\sbs R'$.
Observe furthermore that in exactly the
same way as we defined $c_{1,\om}$ in
(\ref{eq:coverFix}) we can 
define 
$c_{2,\om}:R'_{2,\om}\ra R'_\om$ via multiplication
by $i'_\om$ on the right, and
that this is a lifting of $c_{1,\om}$ to
an isomorphism of left $\G$-torsors. On the other hand
$c_\om$ of (\ref{eq:cor9}) is an isomorphism
of left/right $\G$-torsors with the same domain
and range, so there is a unique $\gamma_\om\in\G_\om$  
such that $c_{2,\om}(x)=c_\om(x\bullet \gamma_\om)=
c_\om(x)\bullet \gamma_\om$, where here the multiplication
$\bullet$ is to be understood as that of right
$\G$-torsors. Now conjugation by $i'_\om$ restricts
to a constant automorphism $A'_\om$, say, of 
the stabiliser of $R'$ so if $\cdot$,
which coincides with $\bullet$ on $R'_2$, is
the groupoid product (\ref{eq:cor10}) then
for any $h\in\G$ we have
\begin{equation}\label{eq:cor11}
\begin{split}
c_\om(x)\cdot F'_\om (h\gamma_\om) =&
c_\om(x)\bullet (h\gamma_\om)=
c_\om(x\cdot h\cdot \gamma_\om)=
c_{2,\om}(x\cdot h)\\
=&c_{2,\om}(x) \cdot A'_\om (h)
=c_\om(x\cdot \gamma_\om) \cdot A'_\om(h)
=c_\om(x) \cdot (F'_\om(\gamma_\om) A'_\om(h)) 
\end{split}
\end{equation}
from which $F'_\om$ coincides on $P\ts \G$ with
$A'_\om$ up to conjugation by $F'_\om(\gamma_\om)$,
so that changing $i'_\om$ appropriately, we may
suppose that they coincide exactly, and since
we already known that they coincide on the
connected component of the identity, we conclude
that we may suppose that $F'_\om$ is conjugation by $i'_\om$ 
restricted to $R'_2$. This has the convenient
consequence that
\begin{equation}\label{eq:cor12}
\a'_{\tau,\om}= (F'_\om)^*i'_\tau i_\om(i'_{(\tau\om)})^{-1}
\end{equation}
and allows us to calculate the right action of
$\Pi_2$ on the fibres of $q:\cX':=[P/R']\ra \cX$.
To this end we can,
by \ref{fact:cor4}, identify $\rB_\G$ with the
groupoid consisting of the 2-cells
\begin{equation}\label{eq:cor13}
 \xy
 (0,0)*+{R_2}="1";
 (36,0)*+{}="2";
 (41,0)*+{R'}="22";
 (18,-20)*+{R_0}="3";
   "1";"2" **\crv{(7,10)} ?(.97)*\dir{>};
   "1";"2" **\crv{(7,-10)} ?(.97)*\dir{>};
     (14, -12)*+{}="D";
     (22,-12)*+{}="C";
     (12 , -9 )*+{}="B";
     (24,-9 )*+{}="E";
     (9,-2 )*+{\scriptstyle r'}="F";
     (29, 4 )*+{\scriptstyle r'}="G";
     (24,2 )*+{}="H";
       {\ar@{=>}_{\mathrm{id}} "D";"C"};
       {\ar@{:>}^{\mathrm{id}} "B";"E"};
       {\ar^{q'} "2";"3"};
       {\ar_{q} "1";"3"};
       {\ar@{=>}^{\gamma} "F";"H"};
\endxy
\end{equation}
for $r':R_2\ra R$ the natural map constructed above,
and $\gamma\in\G$ any.
The relevant (non-commuting) diagram of 1 and 2-cells in 
$\et_2(\cX)$ is therefore
\begin{equation}\label{eq:cor14}
 \xy
 (0,0)*+{(r'F_\om, i_\om)}="A";
 (60,0)*+{(r'F_\om, i_\om)}="B";
 (60,-35)*+{(r',\mathrm{id})}="C";
(0,-35)*+{(r',\mathrm{id}) }="D";
(50,10)*+{(r'F_\tau F_\om, \om^*i_\tau i_\om)}="E";
 (110,10)*+{(r'F_\tau F_\om, \om^*i_\tau i_\om) }="F";
 (110,-25)*+{(r'F_{(\tau\om)}, i_{\tau\om})}="G";
(50,-25)*+{(r'F_{(\tau\om)}, i_{\tau\om}) }="H";
{\ar^{\gamma} "A";"B"};    
{\ar^{i'_\om } "C";"B"};
    {\ar_{\gamma } "D";"C"};
{\ar_{}_{i'_\om } "D";"A"};
    {\ar^{\gamma} "E";"F"};
    {\ar^{\a_{\tau,\om} } "F";"G"};
{\ar@{-->}_{\gamma} "H";"G"}    
{\ar@{-->}_{\a_{\tau,\om}} "E";"H"}
{\ar@{-->}^{i'_{(\tau\om)}} "D";"H"}
{\ar^{\om^* i'_\tau  } "A";"E"};
    {\ar_{ \om^* i'_\tau } "B";"F"};
{\ar_{i'_{(\tau\om)}} "C";"G"};
\endxy
\end{equation}
where $\gamma\in\G$ should be interpreted as 
an arrow of the stabiliser acting on the
left. These are arrows in a groupoid, so
the right action of $\pi_1$ on $\G$ is
given by moving the foremost top occurrence
of $\g$ to the foremost bottom, whence its
conjugation by $(i'_\om)^{-1}$ which,
as it should be, is
$A_\om^{-1}$ on $\G$. Similarly, the
associator $A_{(\tau\om)}^{-1}\Rightarrow A_\tau^{-1} A_\om^{-1} $
is $(i'_\om)^{-1}((F'_\om)^* i'_\tau)^{-1} \a^{-1}_{\tau,\om} i'_{(\tau\om)}$
which by
(\ref{eq:cor12}) is
\begin{equation}\label{eq:cor15}
(i'_{(\tau\om)})^{-1} \bigl( (\a'_{\tau,\om})^{-1} \a_{\tau,\om}^{-1}\bigr)
i'_{(\tau\om)}
= A_{(\tau\om)}^{-1}(\zeta_{\tau,\om}^{-1})
\end{equation}
which by (\ref{eq:LeftRight8}) is exactly what we
started with expressed as a right rather than
a left action.

Turning to 1-cells, then again by 
\ref{fact:groupPlus1}, or better
\ref{claim:groupPlus2} in the specific,
we may suppose that the 1-cell on
the right of (\ref{eq:cor2}) is a
map between the left action of $\Pi'_2$
on a group $\rB_\G$, and that of
$\Pi''_2$ on another group $\rB_\D$
for $\Pi'_2$ a sub 2-group of $\Pi''_2$,
which itself is a sub 2-group of $\Pi_2$.
As such, we could without loss of
generality suppose that $\Pi''_2=\Pi_2$,
but this doesn't lighten the notation
very much, and the structure is clearer
if we don't. In any case, as above
the action of $\Pi'_2$ will be encoded
in a pair $(A,\z)$, that of $\Pi''_2$
in a pair $(B,\eta)$, and anything not
so far defined for the second action
will be its analogue for the first
action but with a $''$ rather than a $'$. 
The map itself will be written $(x,\xi)$,
in exactly the same notation as the
proof of \ref{claim:groupPlus2} except
that $f(*):\G\ra\D$ is replaced by $x:\G\ra\D$ to avoid
confusion with (\ref{eq:cor10}). This
said, we construct a functor $X:R'\ra R''$
as follows: we have a functor $\mathrm{id}\uts x: 
R_2\uts\G\ra R_2\uts\D$, so by
(\ref{eq:cor7}) we get a functor $X_2:R'_2\ra R''_2$.
Better still, for any $\om\in\pi'_1$ we have
- (\ref{eq:groupPlus29}) and (\ref{eq:cor7}) again-
natural transformations $\xi_\om: F''_\om X_2\Rightarrow X_2 F'_\om$
on identifying $\D$ with the stabiliser in $R''_2$.
Profiting from (\ref{eq:cor11}), and thereabouts,
we may write any arrow in $R'$, respectively $R''$,
as $f'i'_\om$, respectively $f''i''_\om$, for
unique
arrows in the respective fibres over $R_{1,\om}$,
and appropriate liftings of the sections $i_\om$.
As such, the following function is well defined
\begin{equation}\label{eq:cor16}
X:R'\ra R'': f (=f' i'_\om) \mpo X_2(f') \xi_\om i''_\om
\end{equation}
and 
since $\xi_\om: F''_\om X_2\Rightarrow X_2 F'_\om$,
checking that it's a functor amounts to
\begin{equation}\label{eq:cor17}
X_2(\a'_{\tau,\om}) \xi_{\tau,\om}
= \xi_\tau F''_{\tau}(\xi_\om) \a''_{\tau,\om}
\end{equation}
which follows from (\ref{eq:groupPlus30}),
and \ref{fact:visit4}.  By construction,
$X:R'\ra R''$ strictly commutes with the
natural maps of $R'$, and $R''$ to $R_0$, so $(X,\mathbf{1})$
is a 1-cell in $\et_2(\cX)$. To check that
this is what we started with the appropriate
(commuting only on triangles) diagram in $\et_2(\cX)$ is
\begin{equation}\label{eq:cor18}
 \xy
 (25,5)*+{(r''F_\om=Xr'F_\om, i_\om)}="A";
 (85,5)*+{(r''F_\om=Xr'F_\om, i_\om)}="B";
 (60,-35)*+{*=(r'',\mathrm{id})}="C";
(0,-35)*+{(r'',\mathrm{id})=* }="D";
 (110,-25)*+{(r''F_{\om}, i_{\om})}="G";
(50,-25)*+{(r''F_{\om}, i_{\om}) }="H";
{\ar^{x(\gamma)=X(\g)} "A";"B"};    
{\ar^{X(i'_\om) } "C";"B"};
    {\ar_{x(\gamma)=X(\g) } "D";"C"};
{\ar_{}_{X(i'_\om) } "D";"A"};
    {\ar^{\xi_{\om} } "G";"B"};
{\ar@{-->}_{x(\gamma)=X(\g)} "H";"G"}    
{\ar@{-->}_{\xi_{\om}} "H";"A"}
{\ar@{-->}^{i''_{\om}} "D";"H"}
{\ar_{i''_{\om}} "C";"G"};
\endxy
\end{equation}
for any $\g\in\G$, which on the stabiliser of the base point gives,
\begin{equation}\label{eq:cor19}
B_\om^{-1}(\xi_\om^{-1})=(i''_\om)^{-1} \xi_\om^{-1} i''_\om:
B_\om^{-1}x
\Rightarrow
xA_\om^{-1}
\end{equation} 
which is exactly what we started with,
(\ref{eq:groupPlus29}), albeit written
on the right rather than the left by
way of (\ref{eq:LeftRight5}) applied
to $\fAut(\rB_\G)$ and $\fAut(\rB_\D)$,
cf. (\ref{eq:LeftRight8}).
\end{proof}
Let us conclude by 
by way of a couple of remarks.
In the first place
\begin{rmk}\label{rmk:corFunct}
{\bf Functoriality of \ref{prop:cor1}}
.
To the extent that one is prepared to
admit inaccessible cardinals, and 
employ the strict fibre product $\uts$ of
\ref{def:FibreS}, functoriality
of $\cX\mpo\et_2(\cX)$ is trivial.
Indeed for any 1-morphism, $f:\cX\ra\cY$,
we not only have a functorial pull-back,
$f^*:\et_2(\cY)\ra\et_2(\cX):q\mpo q\uts_{\cY}\cX$, but if
$\phi:f\Rightarrow g$ is any 2-morphism,
then $f^*$ and $g^*$ agree up to unique
strict isomorphism on cells. If one doesn't
do this, and arguably one shouldn't, then
it's honest bi-category territory as described
in \cite{sga1}[Expos\'e VI.7-9]. Similarly,
if we work only with (weakly) pointed morphisms
$f:\cX_*\ra\cY_*$, then we have a functor
$f_*:\rP\cX_*\ra\rP\cY_*$ on path spaces, which
by (\ref{eq:pointuni1}) takes natural
transformations $\phi:f\Rightarrow g$ to
strict identities. Consequently, the 
presentation $[\rP\cX_*/R_0(\cX_*)]\xrightarrow{\sim}\cX$
by way of the path fibration is again
functorial with little role for 
2-morphisms $\phi:f\Rightarrow g$ which
for pointed maps are arrows commuting
with the image of the path-groupoid of
their source, so that on this side of the
correspondence \ref{prop:cor1} we
only get into bi-category territory if
we drop the pointing. Indeed, on preserving
the pointing, we can render the construction
of $\Pi_2(\cX_*)$ a little more functorial
by defining the $F_\om$, and $\a_{\tau,\om}$
of \ref{fact:act1} on the whole loop space
$\Om\cX_*$ rather than just a complete repetition
free list of homotopy classes. As such we
get an equivalent 2-group with objects $F_\om$,
arrows natural transformations in $R_2$,
a monoidal product $F_\tau\otimes F_\om:=F_{(\tau\om)}$
via concatenation of paths and associator given
by exactly the
same formula \ref{factdef:act1}. Of course,
$\a^{\cX}_{\tau,\om}$ still involves a choice,
but if $f:\cX_*\ra\cY_*$ is a map, then we
have  $f_2: R_2(\cX_*)\ra R_2(\cY_*)\ts_{R_1(\cY_*)} R_1(\cX_*)$,
and 
\begin{equation}\label{eq:corFunct}
f_2(\a^{\cX}_{\tau,\om})(\a^{\cY}_{f(\tau), f(\om)})^{-1}\in
\pi_2(\cY_*)
\end{equation}
so by \ref{fact:group1}, or
better its origin \cite[Theorem 43]{baez},
and the definition \ref{factdef:act1},
the formula (\ref{eq:corFunct}) actually
defines a map $f_*:\Pi_2(\cX_*)\ra\Pi_2(\cY_*)$,
and this is plainly not just functorial but constant under
homotopies $f\mpo f_t$, $t\in I$ since
$\pi_2(\cY_*)$ is discrete with no role for
2-morphisms, $\phi:f\Rightarrow g$, which
collapse to identities. 
As such we get a strictly commuting-
functorial in $f:\cX_* \ra \cY_*$- diagram of two categories
\begin{equation}\label{eq:corFunct1}
\begin{CD}
\et_2(\cY)@>>{f^*}> \et_2(\cX) \\
@VVV @VVV \\
\mathrm{Grpd}(\Pi_2(\cY_*)) @>> {(f_*)^*}> \mathrm{Grpd}(\Pi_2(\cX_*))
\end{CD}
\end{equation}
on using the strict fibre product, $\uts$,
to define the vertical maps via the fibre
functor (\ref{eq:cor3}). The vertical inverse may,
of course, only be weak, but it's necessarily
unique up to equivalence, so at worst (\ref{eq:corFunct1})
becomes 2-commutative on reversing the vertical
arrows. 
\end{rmk}
Finally some further comment on the
role of pointing is in order
\begin{rmk}\label{rmk:corPoint}
 It's useful to have a more geometric
interpretation of the monoidal variant, \eqref{eq:groupTom},
of equivalence in $\et_2(\Pi_2)$. Specifically,
if we work with pointed 0-cells, $q':\cX'_{*'}\ra\cX_*$
and weakly pointed maps, then  we can
further add to the data by introducing
pairs
\begin{equation}\label{eq:corPoint1}
\text{$(q',\iota')$ where $\iota'$
is an isomorphism of a discrete group
$G'$ with the stabiliser of $*'$.}
\end{equation}
Now to demand that a (weakly pointed) 1-cell
$f':q'\ra q''$ in $\et_2(\Pi_2)$ commutes
with the isomorphisms $\iota'$, $\iota''$
is not a restriction unless one limits
the possibilities for the resulting map
$\G'\ra\G''$. Such limitations can only
ever be natural under restrictive hypothesis, {\it e.g.}
$\G'=\G''$, and, indeed, in this case
we can reasonably talk about
{\it the equivalence class of $0$-cells 
with pointed stabiliser $\G'$}. Unsurprisingly,
such restricted equivalence has the effect
of forcing $f$ to be the identity in
(\ref{eq:groupPlus27}), and corresponds,
as per \ref{sum:group1},  to
the
the equivalence
class of the action $\Pi_2(\cX_*)\ra\fAut(\G')$
qua monoidal functor.
Consequently the
equivalence class of  {\it $0$-cells 
with pointed stabiliser $\G'$} is,
for $Z'$ the centre of $\G'$,
always a principal 
homogeneous space under
$\rH^2(\pi'_1, Z')$. 
\end{rmk}

\subsection{The 2-category 
\texorpdfstring{$\et_2(\Pi_2)$}{\'Et\_2 (Pi\_2)} by examples.}\label{SS:II.6}

The 2-category $\et_2(\rp)$ is, in fact,
a category, {\it i.e.} the category of
groups, while the general 0 dimensional
case of the 2-Galois correspondence 
amounts to some well known, 
at least for $0$ cells,
group theory,
to wit

\begin{example}\label{eg:eg1}
{\bf The 2-category
$\et_2(\rB_{\pi_1})$.} 
The universal 1-cover
of $\rB_{\pi_1}$ is $\rp$, so this is also
the universal 2-cover, and the 2-group
$\Pi_2 ( \rB_{\pi_1})$ can be identified
with the group $\pi_1$, or, better, 
the 2-type $(\pi_1, 0, 0)$. Similarly,
$\rB_{\pi_1}$ is itself
an \'etale fibration over $\rp$, so every
connected 0-cell is just a map of groups $q':\rB_{E'}\ra\rB_{\pi_1}$,
and
the factorisation 
\ref{cor:sep3} is simply the exact sequence
\begin{equation}\label{eq:eg1}
1\ra \G'\ra E'=\pi_1(q') \twoheadrightarrow \pi'_1 \hookrightarrow \pi_1
\end{equation}
for
$\pi'_1$ the image of $E'$ in
$\pi_1$. 
Such an object is synonymous with
a pointed stabiliser (qua 2-group) action of $\pi'_1$
on $\rB_\G$,
and 
by \ref{fact:group3}
it exists iff there is a map of
2-types $\underline{q'}: (\pi'_1, 0, 0)\ra (\Out(\G'), Z', \obs_{\G'})$,
{\it i.e.} a representation 
\begin{equation}\label{eq:eg2}
\text{$\underline{q'}:\pi'_1\ra\Out(\G')$
with $\underline{q'}^*\obs_\G'=0\in \rH^3(\pi'_1, Z')$;
$Z'$ the centre of $\G'$.
}
\end{equation}
Now in speaking of the isomorphism classes of extensions
one normally means up to  diagrams of the form
\begin{equation}\label{eq:egTome4}
\begin{CD}
1@>>> \G'@>>> E'@>>>\pi'_1@>>>1 \\
@. @| @VV{\scriptstyle\mathrm{isomorphism}}V @| @.\\
1@>>> \G'@>>> F'@>>>\pi'_1@>>> 1
\end{CD}
\end{equation}
This is not, however, 
\ref{sum:group1},
the definition of equivalence
of 0-cells in $\et_2(\rB_{\pi_1})$, but
rather, (\ref{eq:groupTom}), the
equivalence class of the representation
$\Pi_2\ra \fAut(\rB_{\G'})$ qua monoidal functor,
or, perhaps better, \ref{rmk:corPoint}, 
of $0$-cells with stabiliser pointed in
$\G'$;
so that on
fixing $\underline{q'}$,
\ref{fact:group1} assures us that
the isomorphism classes of
extensions, (\ref{eq:eg1}),
\begin{equation}\label{eq:eg3}
\text{in the
sense of (\ref{eq:egTome4})
are a principal
homogeneous space under
$\rH^2(\pi'_1, Z')$.}
\end{equation}
To see that (\ref{eq:eg3}) is, in general, 
different from equivalence of 0-cells in
$\et_2(\rB_{\pi_1})$, observe that
a 1-cell $F:q'\ra q''$ is a
pair $(f:E'\rightarrow E'', \om\in\pi_1)$
fitting into a commutative diagram
\begin{equation}\label{eq:eg4}
\begin{CD}
1@>>> \G'@>>> E'@>>>\pi'_1\hookrightarrow\pi_1\\
@. @VV{f'}V @VV{f}V @VV{\mathrm{Inn}_\om}V\\
1@>>> \G''@>>> E''@>>>\pi''_1\hookrightarrow\pi_1
\end{CD}
\end{equation}
where, as indicated, necessarily $f$ restricts to a map
$f':\G'\ra\G''$. 
The intervention of $\mathrm{Inn}_\om$ is something of
a red herring. Indeed by \ref{claim:groupPlus2},
this is just a question of base points, which
can conveniently be suppressed on
 replacing
$\pi'_1$ by its conjugate and $E'$ by its pull-back
(since $\mathrm{Inn}_\om$ is only inner in $\pi_1$)
to the same. Plainly, however,
exactly as per (\ref{eq:groupTom}), there are
more equivalences than  
(\ref{eq:egTome4}) in $\et_2(\rB_{\pi_1})$ since all that's required
is that $f'$ in (\ref{eq:eg4}) be an isomorphism.
Similarly,  
not every map $f'$, be it an isomorphism
or otherwise, can
be realised 
in (\ref{eq:eg4}), {\it e.g.} for $f'$
an isomorphism,
the equivalence class in $\et_2(\rB_{\pi_1})$ of the 
0-cell, $q'$,
is the obstruction.
When the map  $f'$ can be realised,
{\it e.g.} $q'=q''$ and $f'=\mathbf{1}$, by
\ref{sum:group1}, 
the 1-cells modulo equivalence in $\Hom(q',q'')$
\begin{equation}\label{eq:eg5}
\text{with $(f', \om)$ fixed are isomorphic
to $\rH^1(\pi'_1, C_{f'})$; $C_{f'}$
centraliser in $\G''$ of $f'(\G')$ }
\end{equation}
in the $\pi'_1$ module structure (\ref{eq:groupPlus31}). 
Finally, a 2-morphism $\e:(f,\om) \Rightarrow (g,\tau) $ in $\Hom_2(q',q'')$
is just a conjugation by $\e\in E''$ such that
\begin{equation}\label{eq:eg55}
\text{$g=\mathrm{Inn}_\e f$ and $\mathrm{Inn}_\tau=
\mathrm{Inn}_{\bar{\e}} \mathrm{Inn}_\om$;
$\e\mpo\bar{\e}\in\pi''_1$.}
\end{equation}
In the proof of \ref{claim:groupPlus2} we didn't
need to calculate the obstruction, but this lacuna
is easily rectified.
An obvious necessary condition for $\e$ to exist
is that $\mathrm{Inn}_{\tau^{-1}\om}$ is actually
inner in $\pi''_1$. Now as we've said post (\ref{eq:eg4})
we can suppose that $\om$ is the identity, whence
at this point $\tau$ has to be in $\pi''_1$, so
replacing $g$ by a conjugate we can, without loss
of generality, suppose that $\om=\tau=\mathbf{1}$,
and we're looking for some $\e$ such that $\bar{\e}$
is in the centraliser, $\bar{C}$, of $\pi'_1$ in $\pi''_1$.
A similarly obvious necessary condition is that the restrictions
$g'$, $f'$ can be conjugated by an element of $E''$
lying over $\bar{C}$, so supposing this holds we
can conjugate $g$ to obtain $f'=g'$, from which
the remaining obstruction is exactly the 
isomorphism class (\ref{eq:eg5}), so modulo
the preceeding two obvious conditions the
obstruction to $\Hom_2(f,g)\neq \emptyset$ lies in $\rH^1(\pi'_1, C_{f'})$,
and 
\begin{equation}\label{eq:eg6}
\text{
if this vanishes,
$\Hom_2(f,g)$ is a principal
homogeneous space under 
$\rH^0(\pi'_1, C_{f'})$;}
\end{equation}
where, 
\ref{sum:group1},
 the principal homogeneous space structure
is the action of $\Hom_2(f,f)$.
\end{example}
Plainly, it doesn't take much effort to extend
this example to
\begin{example}\label{eg:eg2}
{\bf The 2-category
$\et_2(\cX_*)$ whenever $\pi_2(\cX_*)=0$}, and,
of course, $\cX$ is an \'etale fibration over
a separated champ. Plainly $\Pi_2(\cX_*)$ is
again the group $\pi_1=\pi_1(\cX_*)$, so
the 2-category $\et_2(\cX_*)$ is necessarily
equivalent to $\et_2(\rB_{\pi_1})$ by
\ref{prop:cor1}. As such, it's only a 
question of how the cells manifest themselves.
To fix ideas, suppose in the first instance that
$\cX_*$ is connected and developable. Consequently,
$\cX_1$ is a space, $X$ say, and $\cX$ is the
classifying champ $[X/\pi_1]$ for
some not necessarily faithful representation
$\rho:\pi_1\ra \Aut(X)$. This leads to obvious
candidates for the cells, {\it i.e.} identify
a 0-cell, $q'$, with an extension, $E'$, \`a la
(\ref{eq:eg1}), whence by pull-back we have a representation
$\rho_{E'}:E'\ra \Aut(X)$, so we can form
the classifying champ $X_{q'}:=[X/E']$, then
identify $f$ of (\ref{eq:eg4}) with
a functor $X_{q'}\ra X_{q''}$, and $\om$,
$\e$ with natural transformations in the
obvious ways. In the general case, this may
not have sense, but given $E'$ one can-
notation as per \ref{eq:coverFix}- form
the groupoid
\begin{equation}\label{eq:eg7}
R_{E'}:=\coprod_{e\in E'} R_0^{\bar{e}}\rras \rP\cX_*,\quad
e\mpo \bar{e}\in \pi'_1
\end{equation}
then form the classifying champ 
$[\rP\cX_*/R_{E'}]$, and similarly
for the 1 and 2 cells, which is
not only equivalent to the previous
construction in the developable case,
but proves that in either case this
is the equivalence of 
$\et_2(\cX_*)$ with 
$\et_2(\rB_{\pi_1})$ by the simple
expedient of applying $\pi_0$ at the
groupoid level.
\end{example} 
There is, to some degree, an analogous description in general,
to wit:
\begin{example}\label{eg:eg3}
{\bf The full sub-2-category 
of connected $0$ cells in $\et_2(\cX_*)$
as ``non-associative extensions''}, {\it cf.} \ref{fact:pos1}.
We employ the abbreviation $\Pi_2=\Pi_2(\cX_*)$,
and confuse $\et_2(\Pi_2)$ with $\et_2(\cX_*)$,
as we may legitimately do by \ref{prop:cor1}.
We know that the $0$-cells, $q'$,  
in  $\et_2(\Pi_2)$ 
are determined
by a pointed stabiliser action, \ref{fact:groupPlus1},
so, in particular, we have a map of 2-types
\begin{equation}\label{eq:eg++}
\underline{q}': (\pi'_1,\pi'_2, K'_3) \ra 
(\Out(\G'), Z', \mathrm{obs}_{\G'})
\end{equation}
generalising (\ref{eq:eg2}). The fundamental
group, $\pi_1(q')$, of the cell is 
$\pi_0(R')$ for $R'$ the groupoid
in the proof of \ref{fact:cor5},
so it's
an
extension of $\pi'_1$ by $\G'/\pi'_2$, and
we have an exact sequence
\begin{equation}\label{eq:eg+1}
1\ra \pi'_2 \ra \G'\ra \pi_1(q')\ra \pi'_1\ra 1
\end{equation}
which, in turn,
generalises (\ref{eq:eg1}), and allows
us to identify $\pi'_2$ with a sub-group
of $\G'$.
This sequence is, however, by no means arbitrary
since the condition of a pointed stabiliser
action is  
(for $A_\om$, ${\z}_{\tau,\om}$ as in the proof of
\ref{fact:cor5})
equivalent to 
\begin{equation}\label{eq:eg+5}
\otimes':(y, \tau)\ts (x, \om) \mpo (yA_{\tau}(x) {\z}_{\tau,\om}, \tau\om)
\in E':=\G' \ts \pi'_1\xrightarrow{e\mpo\bar{e}}\pi'_1 
\end{equation}
defining a binary operation 
on $E'$ restricting to the
group operation on $\G'$ such that
\begin{equation}\label{eq:eg+7}
e\otimes'(f\otimes' g) = K'_3( \bar{e},\bar{f},\bar{g}) (e\otimes' f) 
\otimes' g
\end{equation}
which (since all co-chains are normalised) 
has a two sided identity, and inverses on both
the left and right, albeit these can differ
because of the lack of associativity.
Nevertheless (normalisation again) there is
no lack of associativity when one of 
$e,f,g$ of (\ref{eq:eg+7}) belongs to $\G'$,
so we get left and right actions of $\G'$ on
$E'$ and whence an ``exact sequence'',
\begin{equation}\label{eq:eg+6}
1\longrightarrow\G'\longrightarrow E'\xrightarrow{e'\mpo\bar{e'}} 
\pi'_1\longrightarrow 1
\end{equation}
all of which may reasonably be summarised
by saying that (\ref{eq:eg+5})-(\ref{eq:eg+6}) define a
{\it non-associative extension of $\pi'_1$ by $\G'$ with
associator $K'_3$}. Now, plainly, this structure is
just that of a 2-group in disguise; in fact $E'$ is
the set of arrows  
of the groupoid
$\cE':=\G'\ts\pi_1(q')\rras \pi'_1(q')$ -$\G$ acting
on the left by way of its image in
(\ref{eq:eg+1})- and $\otimes'$ of 
(\ref{eq:eg+5}) gives this the structure
of a monoidal category. In particular by
way of the section from $\pi'_1$,
implicit in (\ref{eq:eg+5}), we have
a strict ({\it i.e.} monoidal product goes 
to monoidal product
on the nose) monoidal functor $\Pi'_2\ra\cE'$;
while the 
stabiliser
action of $\Pi'_2$ comes from pulling back 
the conjugation action
\begin{equation}\label{eq:eg++6}
\begin{split} 
 E'\ts \G'\ra\G': (e, \g)  \mpo e\cdot g:  =
e\otimes' & \g \otimes'\{\text{right inverse of $e$}\},\\
\text{or equivalently:}\, \, (x,\om)\ts\g  & \mpo xA_\om(\g)x^{-1}
\end{split}
\end{equation}
which- normalisation 
of co-chains again-
satisfies $(f\otimes' e)\cdot \g= f\cdot(e\cdot \g)$,
As such, 
it's plain from (\ref{eq:eg++6}) that
(\ref{eq:eg+1}) expresses $\G'$ as a
crossed module over the group
$\pi_1(q')=E'/\pi'_2$- 
indeed,
the fundamental group
of the fibre of any fibration 
(here $q'$) 
is a crossed
module over the fundamental group of the total
space- while 
the isomorphism class,
\cite[IV.5.4]{brown},
in $\rH^3(\pi'_1, \pi'_2)$
is, unsurprisingly,
$K'_3$. There is, however, more to a 0 cell up to
equivalence 
than this, which in analogy with
(\ref{eq:eg3}) when understood, 
(\ref{eq:groupTom}),
as an equivalence of monoidal
functors, or, \ref{rmk:corPoint} 
of $0$-cells with stabiliser pointed in
$\G'$,
is again the relation
implied by the diagram (\ref{eq:egTome4}), {\it i.e.}
\begin{equation}\label{eq:eg+8}
\text{strict isomorphism classes of
non-associative extensions with
associator $K'_3$,}
\end{equation}
so that just as in
op. cit. as soon as the
2-type, (\ref{eq:eg++}) exists
(or, equivalently, there is a
crossed module of the form
(\ref{eq:eg+1}) with isomorphism class
$K'_3$)
these form a principal homogeneous
space under $\rH^2(\pi'_1,Z')$.
Again, by \ref{fact:groupPlus1}
the intervention of inner automorphisms
in (\ref{eq:eg4})-(\ref{eq:eg6}) is just
a question of base points, so 
a
 1-cell $F:q'\ra q''$ is a
pair $(E'\xrightarrow{f} E'', \om\in\pi_1)$
fitting into the commutative diagram
(\ref{eq:eg4}) where $f$ takes $\otimes'$
to $\otimes''$ strictly,
and $\pi'_2\ra\pi''_2$ under the
action of the same element of $\pi'_1$ affording
the inner automorphism on the right;
while the obstructions to
existence,  
and isomorphism classes,
are exactly as post
(\ref{eq:eg4}).
Finally, since, as noted, associativity holds for
triples one of which is in
$\G''$, the stabiliser of a
1-cell is again the sub-group (\ref{eq:eg6})
of $\G''$ albeit acting by the conjugation
(\ref{eq:eg++6}).
\end{example}
While retaining exactly the same notation let us
make
\begin{example}\label{eg:eg++1}
{\bf The two type of 0-cells in  $\et_2(\Pi_2)$}
. 
For this calculation we fix:
\begin{equation}\label{eq:eg++101}
\text{$\G'$, the 2-type $\underline{q}'$ (\ref{eq:eg++}),
and $\pi_1(q')=E'/\pi'_2$- wholly described by
(\ref{eq:eg+5})-(\ref{eq:eg+6}).}
\end{equation}
We therefore have, for
$\cX'=[\rP\cX_*/R']$
the source of $q'$, that 
$\pi_2(q'):=\pi_2(\cX')$ is just the
kernel of $\pi_2\ra\pi'_2$, and it remains
to determine the Postnikov class $k_3(q')$
in $\rH^3(\pi_1(q'),\pi_2(q'))$. To this end:
write $\G'\ra \pi_1(q')$ as $x\mpo \bar{x}$, and
choose a (set theoretic) section 
$s:\G'/\pi'_2\ra \G'$, then- 
further notation as per (\ref{eq:cor12})-
we get 
sections, 
\begin{equation}\label{eq:eg+99}
\iota_X:= s(\bar{x})i'_{\om}; \quad X=(\bar{x},\om)\in E'/\pi'_2
\end{equation}
of the source of $R'\rras P$
whose value on $*$ has (in a complete repetition free
way) a sink in each of the connected components
of $R'$ yielding natural transformations in $R'_2$
\begin{equation}\label{eq:eg+999}
\b'_{Y,X}:= ((X^*)\iota_Y)\iota_X(\iota_{(YX)})^{-1}:
\mathrm{Inn}_{\iota_Y}\mathrm{Inn}_{\iota_X}
\Rightarrow \mathrm{Inn}_{\iota_{YX}}
\end{equation} 
which can be calculated 
in essentially group theoretic terms
by observing that
\begin{equation}\label{eq:eg+9}
S'_{Y,X}=S'((\bar{y}, \tau)\ts (\bar{x},\om))
:= s(\bar{y})A_{\tau}(s(\bar{x}))\z_{\tau,\om} 
s(\overline{yx\z_{\tau,\om}})^{-1}\in\pi'_2
\end{equation}
is, a $\pi_1(q')$ 2-co-chain with differential
$K'_3$, and, in the notation of (\ref{eq:cor8}) 
and thereabouts
\begin{equation}\label{eq:eg+10}
\b'_{Y,X}= (\text{image in $R'_2$ of}\, \a^{-1}_{\tau,\om} )S_{Y,X}\in R'_2;
\quad X=(\bar{x},\om),\, \text{{\it etc.}}
\end{equation}
Now \ref{factdef:act1}, assures us that to
calculate $k_3(q')$ we need to lift the
$\b'_{Y,X}$ to $R_2$ and calculate their
$\pi_1(q')$ differential under $D$ of
(\ref{eq:act6}). Once, however, (\ref{eq:eg++101}) is
fixed there is
relatively little possibility for variation.
Already this fixes $R'_2=R_2/\pi_2(q')$;
while \ref{factdef:act1} assures 
us not only (cf. the proof of op. cit.) 
that there is no dependence on the choice of
section $s$, but that the only possibility
to change $k_3(q')$ is to change the 0-cell
itself. The possibilities for doing this
are, \ref{eq:groupTom}, up to automorphisms of $\G'$ 
a principal homogeneous space under
the image of $\rH^2(\pi'_1, \pi'_2)$ in
$\rH^2(\pi'_1, Z')$, {\it i.e.} the former
acting (possibly non-faithfully) by way of
\begin{equation}\label{eq:eg+11}
s'_{\tau, \om}\ts \z_{\tau,\om}\mpo s'_{\tau, \om}\z_{\tau,\om}
\in\Home((\pi'_1)^2, \G')
\end{equation}
which by (\ref{eq:eg+9}) and (\ref{eq:eg+10}) translates
to an action
\begin{equation}\label{eq:eg+12}
\b'_{Y,X}\ts  s'_{\tau, \om}\mpo \b'_{Y,X} s'_{\tau, \om}\in R'_2
\end{equation}
Consequently, once one has a lifting,
$\b'_{Y,X}\curvearrowright \b_{Y,X}\in R_2$,
the lifting of anything else in the
orbit is just given by lifting
$s'_{\tau, \om}\curvearrowright  s_{\tau, \om}\in\pi_2$,
and $\pi_2$ commutes with everything, so:
given 
 (\ref{eq:eg++101})
the possibilities for the 2-type
$\underline{\pi}(q')=(\pi_1(q'), \pi_2(q'), k_3(q'))$, 
or, better $k_3(q')$, which is the only unknown, 
are a principal
homogeneous space under the image of
the composite of:
\begin{equation}\label{eq:eg+13}
\rH^2(\pi'_1,\pi'_2) \ra \rH^3(\pi_1,\pi_2(q'))
\ra \rH^3(\pi'_1(q'),\pi_2(q')) 
\end{equation} 
As such, this part of the structure can often
be controlled by knowledge of $\pi_1$ and $\pi_2$
alone, {\it e.g.} if $\pi_1$ is a finite cyclic
group then, \cite{brown}[VI.9.2], 
the co-homology is periodic of
period 2, whence the first map in
(\ref{eq:eg+13}) is always zero, so a fortiori
we have a principal homogeneous space under $0$.
Even here, however,
$k_3(q')$ can still be arbitrary. Indeed,
suppose that $K_3=0$, then 
the orbit of $q'$ under (\ref{eq:eg+11})
has a distinguished element, $q'_0$ say,
for which
$S'_{Y,X}$ of 
(\ref{eq:eg+9}) is the class of $\G'$, viewed
as an extension of $\pi_1(q')$ by $\pi'_2$-
so for varying $\G'$ this is effectively
arbitrary- and whence 
given
 (\ref{eq:eg++101})
the possible values 
whenever $K_3=0$ of $k_3(q')$ are the 
\begin{equation}\label{eq:eg+14}
\text{orbit under (\ref{eq:eg+13}) of the 
image of $\G'$ via $\rH^2(\pi_1(q'), \pi'_2)
\ra \rH^3(\pi_1(q'), \pi(q'))$.}
\end{equation}
Consequently, 
the only real restriction on
2-types in $\et_2(\Pi_2)$ is $\pi_2$, albeit
that the mechanism for creating non-trivial
Postnikov classes from, say, a trivial $K_3$
is perfectly clear.
\end{example}
Finally, let us conclude this tour by
way of the essential
\begin{example}\label{eg:giraud} 
{\bf The 2-category $\fHom(\Pi_2, \rB\rB_Z)$
}
. Here it's very helpful to take  account of \ref{faq8}
and think of $\Pi_2$ as a 2-groupoid with one 0-cell,
so that $\rB\rB_Z$ means the abelian group $Z$ viewed
as a 2-groupoid with one 0-cell, one 1-cell, and 
2-cells the group $Z$. We could make a similar example
for an arbitrary sheaf of abelian groups on $\cX$,
this, however, slightly obscures the absolute nature
of $\rB\rB_Z$. In any case, let us work through the 
2-categorical notion $\fHom(\Pi_2, \rB\rB_Z)$. In
the first instance, the 0 cells, $q'$, are 2-functors,
so \cite[1.5.9]{tom}, monoidal functors from
$\Pi_2$ to $\rB_Z$. As such, we get a (left) action 
\begin{equation}\label{eq:g1}
\Pi_2\ts \rB_Z \xrightarrow{q'\ts\mathrm{id}} \rB_Z\ts\rB_Z
\xrightarrow{\otimes} \rB_Z
\end{equation}
with an associator deduced from $q'(\otimes)\Rightarrow \otimes(q'\ts q')$
in the natural way. 
At which point,
\cite[1.5.11]{tom}, it's critical 
that $\rB\rB_Z$ has only
one 1-cell and all cells are invertible,
so that again 1-cells $f:q'\Rightarrow q''$
are transformations of monoidal functors.
As such, by (\ref{eq:groupTom}) the natural
2-functor $\fHom(\Pi_2, \rB\rB_Z)\ra\et_2(\Pi_2)$
invariably fails to be faithful. There are
many possibilities for understanding this
in $\et_2(\Pi_2)$, of which the most attractive
is by way of torsors under $\rB_Z$, but 
since we wish to consider $Z$ as a given,
it's perfectly reasonable 
(on identifying $\Pi_2$ with
$\Pi_2(\cX_*)$ of a suitable pointed champ-
in fact even of a 3-dimensional CW-complex,
\cite{WhiteheadMaclane}) 
to appeal to
\ref{rmk:corPoint} and view the 
\begin{equation}\label{eq:g2}
\text{
$0$-cells
in $\fHom(\Pi_2, \rB\rB_Z)$ as 0-cells in
$\et_2(\Pi_2)$ with stabiliser pointed in $Z$.}
\end{equation}
Similarly the 1-cells in $\fHom(\Pi_2, \rB\rB_Z)$
may be described as 1-cells
\begin{equation}\label{eq:g3}
\text{
$q'\xrightarrow{f} q''$ in $\et_2(\Pi_2)$ commuting with
the pointed stabiliser isomorphism (\ref{eq:g2}).}
\end{equation}
Finally 2-cells in $\fHom(\Pi_2, \rB\rB_Z)$ are
modifications, \cite[1.5.12]{tom}, which here
just boils down to elements of $Z$, albeit
if we identify the 1-cells with 1-cells, $f$, $g$
in $\et_2(\Pi_2)$ via (\ref{eq:g3}) then 
\begin{equation}\label{eq:g4}
\text{the $2$-cells $\Hom_2(f,g)$ in
$\fHom(\Pi_2, \rB\rB_Z)$ are precisely
those in $\et_2(\Pi_2)$,}
\end{equation}
since commutativity of these with the pointed 
stabiliser isomorphism (\ref{eq:g2}) is implied
by (\ref{eq:g3}).
Now let us identify the equivalence classes
of (\ref{eq:g2})-(\ref{eq:g4}). 
As representations in $\rB_Z$,
even qua equivalence in the monoidal sense,
(\ref{eq:g2}) certainly doesn't give
all the actions of $\Pi_2$ on $\rB_Z$, since
the two type
$\underline{q'}:\underline{\pi}\ra (0,Z,0)$ is
just a map $q'_2:\pi_2=\pi_2(\cX_*)\ra Z$ of $\pi_1=\pi_1(\cX_*)$ modules,
with $Z$ understood trivially,
but, \ref{fact:group3},
it does give all the actions whose pointed
stabiliser has fundamental group $\pi_1$,
and trivial outer representation, {\it i.e.}
$Z$ is the constant, as opposed to a
locally constant, sheaf. Consequently, the
0 cells in $\fHom(\Pi_2, \rB\rB_Z)$,
up to equivalence,
are described by the exact sequence
\begin{equation}\label{eq:g5}
0\rightarrow \rH^2(\pi_1, Z)\rightarrow
\Hom_0(\Pi_2, \rB\rB_Z)/{\scriptstyle \mathrm{equivalence}}
\rightarrow \Hom_{\pi_1}(\pi_2, Z)
\xrightarrow{q'_2\mpo (q'_2)_*(k_3)} \rH^3(\pi_1, Z)   
\end{equation}
which one recognises as being remarkably
similar to the 3rd-sheet of the Hoschild-Serre
spectral sequence, to wit:
\begin{equation}\label{eq:g6}
0\rightarrow \rH^2(\pi_1, Z)\rightarrow
\rH^2(\cX, Z)
\rightarrow \Hom_{\pi_1}(\pi_2, Z)
\xrightarrow{d_3^{0,2}} \rH^3(\pi_1, Z)   
\end{equation}
and, unsurprisingly, \ref{fact:pos2}, these are the same,
{\it i.e.} the equivalence class of 0-cells (\ref{eq:g2})
is naturally isomorphic to $\rH^2(\cX, Z)$. 
In the \'etale topology, 
albeit unlike \ref{cor:one2} with a 
minor intervention of hypercoverings,
this is,  
\ref{fact:383}-\ref{fact:384},
evident from the local to global C\v{e}ch spectral sequence,
and, indeed, {\it cf. op. cit.}, 
one could reasonably just takes this
as a definition of $\rH^2$. In terms
of the path fibration, however, it merits
postponement till \ref{fact:pos2}. In any case
we already know from \cite[Theorem 43]{baez},
and \ref{sum:group1} that the 1-cells
$\Hom_1(q',q'')$ of (\ref{eq:g3}) are a principal homogeneous
space under $\rH^1(\pi_1, Z)$; while the
2-cells $\Hom_1(q',q'')$ of (\ref{eq:g4})
are a principal homogeneous space under
$\rH^0(\pi_1, Z)$, and whence just as
$\rB_Z$ is not just weakly equivalent to
$\rK(Z,1)$ but carries the extra information
of isomorphisms amongst torsors,  $\rB\rB_Z$ similarly augments
$\rK(Z,2)$, {\it i.e.}
\end{example}
\begin{fact}\label{fact:eg1}
Let $\cX$ be a connected, locally 1-connected, and
semi-locally 2-connected champ \'etale
fibring over a separated champ, then via 
(\ref{eq:g2})-(\ref{eq:g4}) and \ref{prop:cor1}, 
$\fHom(\Pi_2(\cX_*), \rB\rB_Z)$ is equivalent to a (sub)
2-category of locally constant gerbes in $\rB_Z$'s
over $\cX$ with stabilisers pointed in $Z$, and 
for any 0-cell $q$ we have canonical isomorphisms,
\begin{equation}\label{eq:g7}
\pi_p(\fHom(\Pi_2(\cX_*), \rB\rB_Z), q)
\xrightarrow{\sim} \rH^{2-p}(\cX, Z), 
\quad
0\leq p\leq 2
\end{equation}
\end{fact}

\subsection{The Torelli description}\label{SS:II.7} 
The fact that a torus, be it $\bg_m$ in the
algebraic sense, or $\rS^1$ in the topological
sense, is a $\rK(\bz, 1)$ provides an alternative
possibility for describing the universal 2-cover
in certain circumstances. The basic case is:
\begin{example}\label{eg:tor1}
{\bf The universal 2-cover, $\cX_2$,
modulo torsion whenever
$\pi_2(\cX_*)$ is finitely generated modulo torsion}  
($\cX$ 
a separated champ with para-compact moduli). 
Under these hypothesis, 
the universal 1-cover
$\cX_1$ is also a
separated champ with para-compact moduli,
so the sheaf of continuous functions on
$\cX_1$ is flasque, and $\pi_2$ modulo
torsion is free of finite rank.
Consequently, we have a 
sequence of commutative Lie groups,
\begin{equation}\label{eq:eg8}
0\ra \pi'_2:=\pi_2(\cX_*)/\mathrm{Tors} 
\ra L_{\pi_2}:=\bc\otimes\pi_2(\cX_*)
\xrightarrow{l\mpo \bar{l}} \Lambda_{\pi_2}\ra 0 
\end{equation}
where the exactness of the sequence defines
the final group. As such- identifying the
co-homology of a Lie group with the sheaf
of continuous functions to it in the usual
way- we obtain by Huerwicz, \ref{fact:one2},
the habitual isomorphism
\begin{equation}\label{eq:eg9}
\rH^1(\cX_1, \Lambda_{\pi_2}) \xrightarrow{\sim}
\Hom(\pi_2(\cX_{1*}), \pi'_2)= \Hom(\pi'_2, \pi'_2)
\end{equation}
and whence a $\Lambda_{\pi_2}$-torsor $\cT\ra\cX_1$.
Appealing to the exact sequence of a fibration,
\ref{fact:fib1}, we have an exact sequence
\begin{equation}\label{eq:eg10}
0\ra \pi_2(\cT_*)\ra \pi_2(\cX_{1*})\ra \pi'_2 \ra \pi_1(\cT_*) \ra 0
\end{equation}
The map in the middle can be studied
one sphere at a 
time, and for $f:\rS^2\ra\cX_1$ arbitrary, 
$f^*\cT$ is a $\Lambda_{\pi_2}$-torsor, so
if one identifies the 
2nd integral co-homology with values
in $\bz(1)$ (=$\bz(2\pi\sqrt{-1})$ which
is canonically $\pi_1(\rS^1)$)
with line bundles $L$ then
by the functoriality of both the isomorphism
(\ref{eq:eg9}) and the exact sequence \ref{fact:fib1}
(albeit the latter involves a choice of
connecting isomorphism up to $\pm 1$ but
this can be done once and for all independently
of the fibration) $f$,  
viewed as a functional on integral co-homology,
goes 
to $\mathrm{deg}(f^*L)$, from which
the map in the middle is just the quotient
of $\pi_2(\cX_{1*})$ modulo torsion, {\it i.e.}
\begin{equation}\label{eq:eg10plus}
\pi_1(\cT)=0, \quad\text{and,}\quad \pi_2(\cT)=\mathrm{Tors}(\pi_2(\cX))
\end{equation}
Again, to fix ideas, suppose
that $\cT$ is a space $T$, then, plainly, $\cX_1=[T/\Lambda_{\pi_2}]$,
and we can use the natural map (\ref{eq:eg8}) to form the
classifying champ $[T/L_{\pi_2}]$ which by 
\ref{fact:fib1} has the same fundamental groups,
(\ref{eq:eg10plus}), as $T$, and
by construction it's a locally constant
gerbe in $B_{\pi'_2}$'s over $\cX_1$. 
Now, quite generally by,
(\ref{eq:eg++101})
a locally constant
gerbes in $\rB_Z$'s 
over $\cX_1$
is simply
connected 
iff it's given by a surjective map
$q:\pi_2(\cX_{1*})\ra Z$.
Consequently
modulo equivalence, over $\cX_1$
\begin{equation}\label{eq:eg10plus1}
\text{there is a unique simply connected
gerbe in $B_{\pi'_2}$'s, so, it's necessarily $[T/L_{\pi_2}]$}
\end{equation}
In particular, by way of example, the so
called ``bad orbifolds'' on $\rS^2=\bp_\bc^1$ with relatively 
prime signatures $p$, $q$ at $0$ and $\infty$ are
simply connected with $\pi_2\xrightarrow{\sim} \bz(1)$,
and enjoy 
the property that $T=\bc^2\bsh\{0\}$ is a space,
while the classifying champs of the actions
\begin{equation}\label{eq:eg11}
\bg_m\ts T \ni (\lb,x,y)
{\build\rras_{(x,y)}^{(\lb^q x, \lb^p y) }}
T,
\quad
\bg_a\ts T \ni (l,x,y) 
{\build\rras_{ (x,y)}^{(\exp(ql) x, \exp(pl) y) }}
T
\end{equation}
are the 
``bad orbifold'' and its universal 2-cover
respectively.

Now the principle utility of this is 
when the torsor is a space so that descriptions 
such as (\ref{eq:eg11}) become a substitute
for developability. Nevertheless, this 
construction continues to have sense 
even if $\cT$ is not a space, 
{\it i.e.} there is a well defined $L_{\pi_2}$
action on $\cT$,
and continues
to be of some utility since $\cT$ is
almost inevitably more space like than
$\cX_1$. To see this,
observe that the descent spectral sequence
\eqref{eq:one7}, applied to
$\Lambda_{\pi_2}$,  
suggests
\begin{equation}\label{eq:eg10++}
\rH^1(\cX_1, \Lambda_{\pi_2})
\xrightarrow{\sim}
\check{\rH^1}(\Hom_{\mathrm{cts}}(R_1, \Lambda_{\pi_2}))
\end{equation}
{\it i.e.} that the torsor is 
a continuous (normalised) co-cycle, 
$K:R_1\ra \Lambda_{\pi_2}$,
a.k.a. a descent datum for the torsor over the
path fibration. 
There is, however, a subtlety, since, \cite[Theorem 1]{pathNotParacompact},
the path space may not be paracompact. Nevertheless,
the moduli of $\cX_1$ is paracompact, 
so the pull-back of any torsor (under any Lie group)
over $\cX_1$ to any space trivialises over a cover
with a locally finite refinement. On the other hand,
by definition,
\ref{defn:fib}, a champ over $\cX_1$ is a
fibration iff it's pull-back to any space is,
so every torsor over $\cX_1$ is a fibration.
As such, the pull-back to $\rP\cX_*$ of any
torsor over $\cX_1$ is a fibration over a
contractible space, so it has a section,
which justifies  \eqref{eq:eg10++}.
Consequently the torsor, $\cT$, can be 
identified with the diagonal action
\begin{equation}\label{eq:eg12}
R_1\ts \Lambda_{\pi_2} \ni (f,\lb)  
{\build\rras_{(s(f),\lb)}^{(t(f), K(f)\lb) }}
\rP\cX_*\ts \Lambda_{\pi_2}  \ra\cT
\end{equation}
so, for example, a meaningful way to present
$\cT\ts \Lambda_{\pi_2}\rras \cT$ in spaces is
\begin{equation}\label{eq:eg13}
R_1\ts \Lambda_{\pi_2} \ts \Lambda_{\pi_2} \ni (f,\lb)\ts x  
{\build\rras_{(s(f), x)}^{(t(f), K(f)\lb x) }}
\rP\cX_*\ts \Lambda_{\pi_2}  \ra\cT
\end{equation}
with the projection of $\eqref{eq:eg12}\ra\eqref{eq:eg13}$,
{\it i.e.} $\cT\ra[\cT/\Lambda_{\pi_2}]$, given by
the functor $(f,\lb)\ra (f, K(f)^{-1})$;
which remains meaningful for $\cT\ts L_{\pi_2}\rras \cT$, to wit
\begin{equation}\label{eq:eg14}
R_1\ts L_{\pi_2} \ts \Lambda_{\pi_2} \ni (f, l)\ts x  
{\build\rras_{(s(f), x)}^{(t(f), K(f)\bar{l} x) }}
\rP\cX_*\ts \Lambda_{\pi_2}  \ra\cT
\end{equation}
Of course, if one prefers to take the reals
in (\ref{eq:eg8}) so that $\Lambda_{\pi_2}$
is a torus in the topological rather than
algebro-geometric sense, then absolutely
nothing changes, and in either case there
is an elegant variation on the Torelli map 
that yields the co-cycle $K$
as soon as $\cX$ has an \'etale presentation
in $C^1$-maps with an atlas of opens in
some $\br^n$. Under such hypothesis we can
identify $\rH^2(\cX_1, \bc )$ with the complex
De-Rham co-homology $\rH^2_{\mathrm{DR}}(\cX_1)$.
Now, \ref{claim:cover1},
arrows in
$R_1$ are 
the images of the space $M$ of 
continuous maps $I^2\ra\cX_1$
with some extra data, (\ref{eq:cover3}),  at the end points,
and
we can 
(via a Stieltjes integral) on choosing
an orientation of $I^2$
integrate 2-forms over such maps to get
a pairing
\begin{equation}\label{eq:eg15}
M \ts
A^2_{\bc}(\cX_1)\ra \bc:
(f,\om) \mpo \int_{I\ts I} f^*\om 
\end{equation}
This pairing doesn't quite descend to $\rH^2_{\mathrm{DR}}$
even for the restricted class of maps
appearing in 
(\ref{eq:cover3})
since the integral of
the differential, $d\b$ of a 1-form, will be the 
difference of the integrals of $\b$ over the
top and bottom sides of the square (\ref{eq:cover3}).
This is, however, the C\v{e}ch co-boundary
in (\ref{eq:eg10++})- 
or more accurately $d_1^{0,1}$ in
the previous sheet-
of the integral of $\b$
over paths. Consequently, if we choose a
section $h:\rH^2_{\mathrm{DR}}(\cX_1) \ra Z^2_{\bc}(\cX_1)\subseteq A^2_{\bc}(\cX_1)$
with values in closed forms
then we get a map, 
\begin{equation}\label{eq:eg16}
\begin{CD}
\tilde{K}:M @>>>  \rH^2_{\mathrm{DR}}(\cX_1)^{\vee}
@>{\sim}>{\mathrm{Huerwicz}}> L_{\pi_2}
\end{CD}
\end{equation}
and any other section differs by the co-boundary
(defined via integration of the difference of sections)
of some $\Hom_{\mathrm{cts}}(\rP\cX_{1*}, L_{\pi_2})$.
By definition, \ref{factdef:cover2}, $R_2$ is $M$
modulo homotopies, and elements of 
$\rH^2_{\mathrm{DR}}$ are closed, so, by Stokes, 
(\ref{eq:eg16}) descends to a map from $R_2$.
In addition $\tilde{K}$ is an integral, so it's
additive for the groupoid structure in $R_2$,
while by the discussion following
(\ref{eq:eg10}) we know that (for an appropriate
orientation of $I^2$) it's value
on $\pi_2(\cX_{1*})$ (identified with
the stabiliser of $*$ in $R_2$) is 
(independently of the section $h$ since
spheres are closed) its natural image in $\pi'_2$,
and whence we obtain the Torelli type formula
\begin{equation}\label{eq:eg17}
K: R_1 \ra \Lambda_{\pi_2}: f\mpo \int_{I^2} f^*h(\om)  
\end{equation}
\end{example} 
Let us clarify the context in which
this may be employed by introducing
\begin{defn}\label{def:tor1}
Let $\Pi_2$ be a 2-group described as
in (\ref{eq:group1})-(\ref{eq:group2}),
then $\Pi_2/\mathrm{tors}$, or $\Pi'_2$ if there
is no danger of confusion, is the sub, 
\ref{fact:groupPlus1}, 2-group
$\Pi_2/(\mathrm{tors}(\pi_2)$. In particular,
the fundamental group is unchanged, and $K'_3$ is
just $K_3$ pushed forward along $\pi_2\ra\pi_2/\mathrm{tors}$.
\end{defn}
With this in mind we can push 
\ref{eg:tor1} to a description of
\begin{example}\label{eg:tor2}
{\bf $\cX_2/\mathrm{tors}\ra \cX$ via torus
extensions}, with not just the hypothesis
of \ref{eg:tor1}, but also
$\rH^i(\pi_1, L_{\pi_2})=0$, $i\in\{1,2\}$. 
With such suppositions,
(\ref{eq:eg8}) gives 
a canonical (up to a sign)
isomorphism
\begin{equation}\label{eq:tors1}
\rH^2(\pi_1, \Lambda_{\pi_2}) \xrightarrow{\sim}
\rH^3(\pi_1, \pi'_2)
\end{equation}
and whence to the 2-type $\underline{\pi}'=
(\pi_1,\pi'_2,k'_3)$ of $\Pi'_2$ we can
associate a unique extension
\begin{equation}\label{eq:tors2}
0\ra \Lambda_{\pi_2}\ra E \ra \pi_1\ra 0
\end{equation}
Unsurprisingly, therefore, our goal is to express 
the composition
\begin{equation}\label{eq:tors3}
\cT\ra\cX_1\ra\cX
\end{equation}
for $\cT$ as per \ref{eg:tor1} as a torsor under $E$. To
this end observe that the reasoning post \eqref{eq:eg10++}
is wholly general, i.e. for any Deligne-Mumford champ
with paracomapct moduli, $\cX$, and any Lie Group, $G$,
\begin{equation}\label{eq:tors40}
\rH^1(\cX, G) \xrightarrow{\sim} \check{\rH}^1(\Hom_{\mathrm{cts}}(R_0, G))
\end{equation}
so, in a minor change of notation, we require to extend the co-cycle 
$K_1:R_1\ra \Lambda_{\pi_2}$ of \eqref{eq:eg10} to a co-cycle
$K:R_0\ra E$. As such let $K_2: R_2\ra L_{\pi_2}$ be the
lifting of $K_1$ to $R_2$ (or even just $R_2/\mathrm{tors}$)
then by the reasoning post \eqref{eq:eg10}, the restriction
of $K_2$ to the stabiliser may be identified with the
natural map $\pi_2=\pi_2(\cX_*)\ra \pi'_2=\pi_2(\cX_*)/\mathrm{tors}$.
In particular, therefore, the restriction of $K_2$ is equivariant
for the action of $\pi_1$, and whence by
\eqref{eq:eg10++}: for each $\om\in\pi_1$
there is a function $y_\om:P=\rP\cX_*\ra L_{\pi_2}$ such that,
\begin{equation}\label{eq:tors4}
K_2( F_{\om^{-1}}(f))^\om = K_2(f) + y_\om(b) -y_\om(a),
\quad
a\xrightarrow{f} b \in R_2
\end{equation}
Consequently, a direct application of the
definition of $\a_{\bullet,\bullet}$
in \eqref{eq:actplus} implies that
the function
\begin{equation}\label{eq:tors5}
\phi_{\tau, \om}:P\ra\Lambda_{\pi_2}:a\mpo
K_2(\a_{\om^{-1}, \tau^{-1}}(a))^{\tau\om} - y_{\tau\om}(a) + y_\tau(a) + y_\om(\tau^{-1}(a))^{\tau}
\end{equation}
is constant on arrows in $R_2$, thus, it is,
in fact, a function from $\cX_2$, or better
its moduli, so inter alia it's  a function
from $\cX_1$ with (left) $\pi_1$-action 
\begin{equation}\label{eq:tors6}
\pi_1\ts\Hom(\cX_1, L_{\pi_2}):(\om,\phi)\mpo 
\{a\mpo \phi(\om^{-1}a)^\om\}
\end{equation}
Combining this with  the respective definitions,
of the Postnikov class, \ref{factdef:act1},
and the co-opposite isomorphism, \eqref{eq:LeftRight10},
gives the identity
\begin{equation}\label{eq:tors7}
\pa (\phi)_{\s,\tau,\om} (a) = (\mathrm{co-op}(K_3))(\s,\tau,\om), \quad a\in P
\end{equation}
of 3 co-cycles with values in $\Hom(\cX_1, L_{\pi_2})$.
Now the hypothesised vanishing of $\rH^3(\pi_1, L_{\pi_2})$ 
allows us
to
write $K_3$ as the boundary of some $L_{\pi_2}$-valued
co-chain, $\xi$, while the H\"oschild-Serre spectral
sequence yields an exact sequence
\begin{equation}\label{eq:tors8+}
 \rH^0(\pi_1, \rH^1_{\mathrm{cts}}(\cX_1, L_{\pi_2})) \xrightarrow{d_2^{0,1}}
\rH^2(\pi_1, \Hom_{\mathrm{cts}}(\cX_1, L_{\pi_2})) \ra \rH^2(\cX, L_{\pi_2})
\end{equation}
in which the groups on the left and right vanish by
the paracompactness of $\cX_1$ and $\cX$ respectively,
so there are 
functions $z_\om:\cX_1\ra L_{\pi_2}$, $\om\in\pi_1$ such that
for any $a\in P$
\begin{equation}\label{eq:tors8}
((\mathrm{co-op})\xi))_{\tau,\om}= (-\phi + \pa(s^*z))_{\tau,\om}(a), \quad 
(s,t):R_1\rras P
\end{equation} 
Amongst the various choices for applying $\mathrm{co-op}$
to \eqref{eq:tors8} we perform 
the substitution $(\tau,\om) \leftrightarrow (\om^{-1}, \tau^{-1})$
while fixing $a$ to obtain
\begin{equation}\label{eq:tors9}
\xi_{\tau,\om} = -K_2(\a_{\tau,\om}(a))  + x_{\tau\om}(a) -x_\tau (\om a) - x_\om(a)^\tau,
\quad  x_\om(a):= y_{\om^{-1}}(a)^\om +  z_{\om^{-1}}(a)^\om
\end{equation} 
The group law in $E$ can, of course, be written,
\begin{equation}\label{eq:tors9plus9}
(e,\tau)\cdot (h,\om)= (e+h^\tau + \xi_{\tau,\om}, \tau\om), \quad
e,h,\in \Lambda_{\pi_2}, \, \tau,\om\in\pi_1
\end{equation}
from which the fact that $s^*z$ in \eqref{eq:tors8} is equally $t^*z$
combines, in the notation of \ref{factdef:cover1},  with \eqref{eq:tors4} 
to yield that
\begin{equation}\label{eq:tors10}
K_0: R_0\ra E : f (=f_1 i_\om) \mpo (K_1(f_1) + (s^*x_\om) (f), \om ),\quad
f_1\in R_1,\om\in\pi_1
\end{equation}
is an $E$-valued co-cycle over $\cX_0$, which restricted
to $R_1$ is $K_1$, and pushed forward to $\pi_1$ just
sends an arrow to its connected component, so the 
composition \eqref{eq:tors3} has the structure of an
$E$-torsor by way of the descent datum \eqref{eq:tors10}.

All of which is considerably easier in practice, which will
amount to:
an atlas, $U$, of opens in some $\br^n$
with quasi-finite \'etale gluing $R\rras U$ along
$C^1$-maps, and finite $\pi_1$. Under the latter
hypothesis all $\bq[\pi_1]$-modules are acyclic,
so we can suppose that the section $h$ prior
to \eqref{eq:eg16} is $\pi_1$ equivariant, so
$y_\om$ in \eqref{eq:tors4} is identically zero for
all $\om\in\pi_1$, while the Torelli description
\eqref{eq:eg17} of $K$ as an integral 
implies that without loss of generality $a\mpo K(\a_{\tau,\om}(a))$
is independent of $a\in P$, so we can take
$x_\om$ in \eqref{eq:tors10} to be zero too.
In any case, since the conclusions of this section are
preliminary in nature let us summarise them
by way of
\begin{summary}\label{sum:tor}
Under the present hypothesis \ref{eg:tor2}    ,so, in particular those
of \ref{eg:tor1}, the 2-type $\underline{\pi}'$ modulo torsion 
determines, \eqref{eq:tors2},  a unique extension
$E$. Similarly, there is a unique $E$-torsor, $\cT_\cX\ra\cX$,
\eqref{eq:tors10}, such that the classifying champ $[\cT_\cX/L_{\pi_2}]$- so
understood as per \eqref{eq:eg14} if $\cT_\cX$ is not
a space- is the universal 2-cover $\cX_2$.
Equally useful, however, albeit with the same
caveat if $\cT_\cX$ is not a space, $\cX$ itself is
the classifying champ $[\cT_\cX/E]$.  
\end{summary}
\end{example}
Now let us apply this to an example which is
common enough whether in algebraic geometry 
or the study of 4-manifolds, to wit a description of
\begin{example}\label{eg:full}{\bf The full
sub-category, $\et_2^{\mathrm{sep}}(\cX)$,  of separated champ in $\et_2(\cX)$}  whenever
$\cX$ itself is separated with para-compact moduli,
finite fundamental group, and torsion free $\pi_2$.
Plainly the source of every 0-cell, $\cY\ra\cX$, of $\et_2^{\mathrm{sep}}(\cX)$
satisfies the hypothesis of \ref{eg:tor2}, so that
our first task is to show that, qua champ rather
than qua torsor, $\cT_\cY$ does not depend on $\cY$,
where for convenience we'll restrict our attention
to $\cY$ connected.
Thanks to the factorisation \eqref{eq:tors3} we may
suppose that $\cY$ is simply connected, while by
\ref{claim:sep1} $\cY\ra\cX$ factors as $\cY\ra\cX_1\ra\cX$
for $\cY\ra\cX_1$ a finite (since $\cY$ is separated) 
locally constant gerbe in $\rB_{\pi'_2}$'s, 
so we get a short exact sequence
\begin{equation}\label{eq:full1}
0\ra \pi'_2 \ra \Lambda_{\pi_2(\cY)} \ra \Lambda_{\pi_2(\cX)}\ra 0
\end{equation}
by virtue of the long exact sequent of a fibration, \ref{fact:fib1},
and the finiteness of $\pi'_2$. By definition we have
an action $\Lambda_{\pi_2(\cX)}\ts \cT\rras \cT$ (
again via \eqref{eq:eg14} if $\cT$ is not a space), which
by way of
the quotient map in \eqref{eq:full1} yields an action
$\Lambda_{\pi_2(\cY)}\ts \cT\rras \cT$ where by
definition $\cT\ra [\cT/\Lambda_{\pi_2(\cY)}]$ is
a $\Lambda_{\pi_2(\cY)}$-torsor; while $[\cT/\Lambda_{\pi_2(\cY)}]\ra
[\cT/\Lambda_{\pi_2(\cX)}]=\cX_1$ is a locally
constant gerbe in $\rB_{\pi'_2}$'s. Better still, again
by the long exact sequence of a fibration, $[\cT/\Lambda_{\pi_2(\cY)}]$
is simply connected, and the natural map,
$\pi_2([\cT/\Lambda_{\pi_2(\cY)}]) \ra \pi_1(\Lambda_{\pi_2(\cY)})$
is an isomorphism; from which the former gives
that $[\cT/\Lambda_{\pi_2(\cY)}]$ is isomorphic to 
$\cY$
by \ref{eg:eg++1} 
and the 2-Galois correspondence 
\ref{prop:cor1},
while the latter gives that
$\cT\ra [\cT/\Lambda_{\pi_2(\cY)}]=\cY$ is the
torsor $\cT_\cY\ra\cY$ by \eqref{eq:eg9} et sequel.

We can now apply this to give a description of
$\et_2^{\mathrm{sep}}(\cX)$ via extensions and
actions on $\cT$ which closely mimics \ref{eg:eg1},
or perhaps better \ref{eg:eg2}, in which the mantra is
that $\cT$ plays the role of the universal cover.
As such, for technical convenience (albeit that
the general case is covered by \eqref{eq:eg14} let
us suppose that $\cT$ is a space, $T$, which is,
in any case, the principle utility of this discussion.
Consequently, we have a base object consisting of
the pair $(E,\rho)$ where $E$ is the extension
\eqref{eq:tors2} deduced from the 2-type of $\cX$,
and $\rho:E\ra\Aut (T)$ the representation assured
by \ref{sum:tor} such that the resulting classifying
champ $[T/E]$ is isomorphic to $\cX$. 
In particular, the representable \'etale cover $\cY'\ra\cX$
corresponding to a sub-group $\pi'_1$ of $\pi_1=\pi_1(\cX)$
is just $[T/E']$ for $E'=E\ts_{\pi_1}\pi'_1$,
and the factorisation $q':\cY\ra\cY'\ra\cX$ of \ref{claim:sep1}
of the cell into a locally constant gerbe followed
by a representable \'etale cover may,
on writing $\Lambda=\Lambda_{\pi_2(\cX)}$, be read 
-{\it i.e.} for $E(q')$ the extension
of \ref{sum:tor} defined by $\cY$,
and everything else defined therein- 
from
the diagram of exact rows and columns with $\pi'_2$ finite
\begin{equation}\label{eq:fullPlus}
\begin{CD}
@. 0@. 0@. 0 @. \\
@. @VVV @VVV @VVV @.\\
0@>>>\pi'_2@>>>\G' @>>> \G'/\pi'_2 @>>> 0 \\
@. @VVV @VVV @VVV @. \\
0@>>>\Lambda@>>> E(q') @>>> \pi_1(q')@>>> 0 \\
@. @V{q'_{\mathrm{toric}}}VV @VV{q'}V @VV{q'_1}V @. \\
0@>>> \Lambda@>>> E' @>>> \pi'_1 @>>> 0 \\
@. @VVV @VVV @VVV @. \\
@. 0@. 0@. 0 @.
\end{CD}
\end{equation} 
Indeed, the inclusion $E'\hookrightarrow E$ affords not just an
action $E(q')\ts T\rras T$, 
and whence a 0-cell $q':[T/E(q')]\ra [T/E]=\cX$ in
$\et_2^{\mathrm{sep}}(\cX)$
but even the  factorisation \ref{claim:sep1}
by way of $[T/E(q')]\ra [T/E']\ra [T/E]$ in which the
former factor is a locally constant gerbe in $\rB_{\G'}$'s.
 Conversely, 
for $q':\cY\ra\cX$ in $\et_2^{\mathrm{sep}}(\cX)$ we
have a strictly commutative diagram
\begin{equation}\label{eq:full3}  
 \xy  
 (-5,0)*+{T\ts_\cY T\xrightarrow{\sim} T\ts E(q')}="A"; 
 (28,0)*+{T\ts T}="B";  
 (14,-12)*+{T\ts_\cX T\xrightarrow{\sim} T\ts E}="C"; 
{\ar_{}^{} "A";"B"};  
    {\ar_{} "C";"B"};  
    {\ar_{} "A";"C"};  
\endxy
\end{equation}
for $E(q')\rras T$ the action assured by \ref{sum:tor},
and hence, by the functoriality of fibre products, a
functorial assignment $q'\mpo E(q')$ fitting into a
diagram (\ref{eq:fullPlus}) in which the notation is
consistent with the more generally valid descriptions 
\ref{eg:eg3}-\ref{eg:eg++1}, {\it e.g.} the 2-type
$\underline{q}'$ encoded by the exact sequence (\ref{eq:eg+1})
is identically \eqref{eq:fullPlus} read clockwise from
top left to bottom right with the top right omitted,
or, equivalently in the sense of crossed modules,
read anti-clockwise with the bottom left omitted.
In any case, 
under our initial hypothesis, 
the 2-category $\et_2^{\mathrm{sep}}(\cX)$
is equivalent to the 2-category, $\cE$,  
in which:
\begin{enumerate}
\item[0]-cells are maps $q':E(q')\ra E$ of extensions in which
the toric map $q'_{\mathrm{toric}}$ is finite.
\item[1]-cells are pairs $(f,\eta)$:  
$f:E(q')\ra E(q'')$ a map of extensions; 
$\eta\in E$, and
$q''f=\mathrm{Inn}_\eta q'$.
\item[2]-cells $\z:(f,\eta)\Rightarrow (g,\xi)$ are
elements $E(q'')$ such that $\xi=q''(\z) \eta$.
\end{enumerate}
Indeed, we have a 2-functor $\cE\ra\et_2^{\mathrm{sep}}(\cX):
E(q')\mpo [T/E(q')]$ which we've already observed is
essentially surjective on 0-cells. A priori \ref{defn:cor1}
leaves open the possibility that there are more 1-cells
in $\et_2^{\mathrm{sep}}(\cX)$ since these are 2-commutative
diagrams
\begin{equation}\label{eq:full4}  
 \xy  
 (-10,0)*+{[T\ts E(q')\rras T]}="A"; 
 (38,0)*+{[T\ts E(q'')\rras T]}="B";  
 (14,-12)*+{[T\ts E \rras T]}="C"; 
(8, -4)*+{}="D";
     (20,-4)*+{}="E";
{\ar_{}^{f} "A";"B"};  
    {\ar^{q''} "B";"C"};  
    {\ar_{q'} "A";"C"};  
{\ar@{=>}_{\eta} "D";"E"};
\endxy
\end{equation}
so, more succinctly: a 1 co-cycle $f_t:E(q')\ra\Hom(T, E(q''))$
for the $E(q')$-action $f_t\mpo f_{et}$, $e\in E(q')$ affording
the functor of the top row of (\ref{eq:full4}); and a continuous
map $\eta_t:T\ra E$ affording the natural transformation
$\eta$ such that we have the 2-commutativity condition 
\begin{equation}\label{eq:full5}
q''f_t(e)= \eta_{et} q'(e) \eta_t^{-1},\quad t\in T,e\in E(q')
\end{equation}
On the other hand, the image of $\eta$ in 
$E/E''\xrightarrow{\sim} \pi_1/\pi''_1$ must be
constant, so we can write $\eta_t=\bar{\eta}''_t \om$ for some
$E''$-valued function $\bar{\eta}''_t$. Furthermore the space
$T$ is 2-connected, so $\bar{\eta}''_t$ can a fortiori be lifted to an
$E(q'')$-valued function, $\eta''_t$ . Consequently we have a 2-cell
$\eta''_t:(F,\om)\Rightarrow (f_t, \eta_t)$ in 
$\et_2^{\mathrm{sep}}(\cX)$ where $F$ is the functor
afforded by the co-cycle
\begin{equation}\label{eq:full6}
E(q') \ra \Hom(T,E(q'')): e\mpo (\eta''_{et})^{-1} f_t(e) \eta''_t
\end{equation}
which by \eqref{eq:full5} pushes forward to a constant map
in $\Hom(T,E)$, so by (\ref{eq:fullPlus}) and the co-cycle
condition, $F$ is actually a map of extensions $E(q')\ra E(q'')$
and our 2-functor $\cE\ra \et_2^{\mathrm{sep}}(\cX)$ is
essentially surjective on 1-cells. Finally, and plainly, 
2-cells in $\et_2^{\mathrm{sep}}(\cX)$ 
between 1-cells where the underlying functor is 
constant in $t\in T$ are themselves constant in $t$,
so again by \cite[1.5.13]{tom} we deduce that
$\cE\ra \et_2^{\mathrm{sep}}(\cX)$ is an 
equivalence of 2-categories.  
\end{example}

\subsection{The Postnikov class}\label{SS:II.8}

The most expedient way to relate the Postnikov class
as we've defined it in \ref{factdef:act1} to the 
spectral sequence (\ref{eq:one7}) and the Huerwicz
theorem as encountered in (\ref{eq:g5}) and (\ref{eq:g6})
is a more quantitative version of \ref{fact:act3}, {\it i.e.}

\begin{fact}\label{fact:pos1}- 
cf. \ref{eg:eg3}-
Let $R_0\rras P$ be the path
groupoid of a pointed champ $\cX_*$, which itself is an
\'etale fibration over a 
1-connected and locally 2-connected
champ with separated moduli,
then the universal cover $R\ra R_0$ admits a product
$\otimes:R\ts_P R \ra R$ lifting the composition in $R_0$
such that for $\bullet:\pi_2(\cX_*)\ts R\ra R$ the left (or
indeed right since $\pi_2=\pi_1(R_1)$ is commutative)
torsor structure, and $\pi_2$ identified to the stabiliser
$\pi_2\ts P$
of the connected component of the identity, {\it i.e.} $R_2$,
\begin{enumerate}
\item[(a)] $\otimes\vert_{R_2}$ is composition in $R_2$,
and lifts composition in $R_0$.
\item[(b)] $S\otimes f=S\bullet f$; $f\otimes S=f\bullet F_\om(f)$,
$S\in\pi_2$, $f\in R\ts_{R_0} R^\om_0$- notation as per (\ref{eq:coverFix}).
\item[(c)] $h\otimes(g\otimes f)=K_3(\s,\tau,\om) 
\bullet\{(h\otimes g)\otimes f\}$,
where $h,g,f$ are in the fibre of $R\ra R_0$ over the
connected components $R_0^\s$, $R_0^\tau$, and $R_0^{\om}$
respectively; $\s,\tau, \om\in \pi_1(\cX_*)$.
\end{enumerate}
In particular since we can, and do, always work with
normalised co-cycles, associativity holds whenever
any of $h,g$ or $f$ in (c) belong to $R_2$.
\end{fact}
\begin{proof} Just apply the construction (\ref{eq:cor9}) \& (\ref{eq:cor10})
of the proof of \ref{fact:cor5} with $R=R'$ in the 
notation of op. cit.
\end{proof}
The principle 
utility of \ref{fact:pos1} is to allow us to
organise the obvious calculation,
\begin{fact}\label{fact:pos2}
Let $\cX_*$ be a (pointed) champ, which itself is an
\'etale fibration over a  
1-connected and locally 2-connected 
champ with separated moduli, with $\underline{Z}$
a locally constant sheaf of abelian groups
on $\cX$, then for $\pi_p=\pi_p(\cX_*)$ the transgression 
\begin{equation}\label{eq:trans1}
d_3^{0,2}:\Hom_{\pi_1}(\pi_2, \underline{Z})
\ra \rH^3(\pi_1, \underline{Z})
\end{equation}
in the 3rd sheet of the Hoschild-Serre spectral
sequence, or, equivalently, the transgression
$d_2^{1,1}$ in the 2nd sheet of the spectral
sequence (\ref{eq:one7}) is given by
$\Hom_{\pi_1}(\pi_2, \underline{Z})\ni \phi\mpo \phi_* k_3$.
\end{fact}
\begin{proof} We do the $d_2^{1,1}$ term in
the spectral
sequence (\ref{eq:one7}) since it's philosophically
more consistent with our use of the path fibration.
As such, to begin with we have, in the notation
of \ref{fact:pos1}, a commutative diagram
\begin{equation}\label{eq:trans2}
\begin{CD}
R@<{\longleftarrow}<{\longleftarrow}< R\ts_P R\\
@VVV @VVV\\
R_0@<{\longleftarrow}<{\longleftarrow}< R_0\ts_P R_0
\end{CD}
\end{equation}
where the bottom horizontal arrows are the
2-projections and composition while those
on the top are again projections along with
the lifting $\otimes$ of composition. Now, the
group $\pi_2\ts\pi_2$ acts faithfully and 
transitively on $R\ts_P R$ by
way of,
\begin{equation}\label{eq:trans3}
(T,S)\ts (g,f) \mpo (Tg,Sf)
\end{equation}
for $(S,f)\mpo Sf$ the faithful and transitive action
of $\pi_2$ on $R$ itself. These actions afford unique
isomorphisms of the fundamental group of
any connected component of $R_0$, respectively
$R_0\ts_P R_0$ with $\pi_2$, respectively $\pi_2\ts\pi_2$,
and whence maps
\begin{equation}\label{eq:trans4}
\begin{CD}
\pi_2\ts \pi_1 @<{\longleftarrow}<{\longleftarrow}< 
\pi_2\ts\pi_2\ts\pi_1\ts\pi_1
\end{CD}
\end{equation}
corresponding to $\pi_0$ of the 3-maps in (\ref{eq:trans2}),
which for $g$ and $f$ of (\ref{eq:trans3}) belonging
to the components indexed by $\tau$, respectively $\om$,
in $\pi_1$ are
\begin{equation}\label{eq:trans5}
(T,S) \mpo T,\quad (T,S)\mpo S^\tau,\quad (T,S)\mpo T+S^\tau
\end{equation}
Consequently, for some suitable
labelling $p_i$, $1\leq i\leq 3$ of the 
bottom horizontal maps in
\ref{eq:trans2}, $(p_3)_*=(p_1)_* +(p_2)_*$,
and
the next thing to look at is a diagram
\begin{equation}\label{eq:trans6} 
\begin{CD}
U\ts_S U @<{\longleftarrow}<{\longleftarrow}< V\ts_X V\\
\downarrow\downarrow@.\downarrow\downarrow \\
U @<{\longleftarrow}<{\longleftarrow}< V:=(p_1)^*U\ts_X(p_2)^*U\ts_X (p_3)^*U \\
@VVV @VVV \\
S @<{\longleftarrow}<{\longleftarrow}< X
\end{CD}
\end{equation}
where $S$ and $X$ are connected, 
$U/S$ is the universal cover,
and the 
maps $(p_i)_*$ of the
bottom 3-horizontal
arrows, $p_i$, on fundamental groups
satisfy the aforesaid relation implied by (\ref{eq:trans5}).
As such, a convenient identification of $V$ is as the
quotient
\begin{equation}\label{eq:trans7}
Y \ts \pi_1^3(S) /(x,\g_i)\sim (x^\om, (p_i)_*(\om) \g_i),
\quad \om \in \pi_1(X)
\end{equation}
for $Y/X$ the universal cover, while an equally 
convenient representation of $V\ts_X V$ is
\begin{equation}\label{eq:trans8}
Y \ts \pi_1^3(S) \ts\pi_1(S)^3 /(x,\g_i, \g_i')\sim (x^\om, 
(p_i)_*(\om) \g_i, (p_i)_* (\om) \g'_i),
\quad \om \in \pi_1(X) 
\end{equation}
Consequently, if we identify $\pi_0(V)$ with $\pi_1(S)$
by way of $(p_3)_*$ in (\ref{eq:trans7}) then the (C\v{e}ch) 
differential of the identity in $\Hom(U\ts_S U, \pi_1(S))=
\Hom(\pi_1(S),\pi_1(S))$ with respect to the top 3-horizontal
maps in (\ref{eq:trans6}) is 
by (\ref{eq:trans7}) and (\ref{eq:trans8})
the differential of the 
identity 
in $\Hom(V,\pi_1(S))$
with respect to the upper right verticals in (\ref{eq:trans6}).
At which point, it's a little more delicate since
the relevant diagram is,
\begin{equation}\label{eq:trans9}
\begin{CD}
@. V @<{\leftleftarrows}<{\leftarrow}<
\prod^{1\leq i\leq 6}_{R_0\ts_P R_0\ts_P R_0} (q_i)^* R \\
@. @VVV @VVV \\
R @<{\longleftarrow}<{\longleftarrow}<
R\ts_P R  @<{\leftleftarrows}<{\leftarrow}< R\ts_P R\ts_P R \\
@VVV @VVV @VVV \\
R_0 @<{\longleftarrow}<{\longleftarrow}<
R_0\ts_P R_0  @<{\leftleftarrows}<{\leftarrow}< R_0\ts_P R_0\ts_P R_0
\end{CD}
\end{equation}
where the $q_i$ are a repetition free list of the 6
maps from among the 12-maps obtained from 
composing the bottom right 
with the bottom left of the diagram, and we
again use $\otimes$, where appropriate, 
to lift the 4 arrows on the bottom right to
the middle 4, so there are  7 repetition free compositions
from middle right to middle left because of
the non-associativity \ref{fact:pos1}.(c).
In this situation
a better model for $V$ is the fibre product of the
leftmost vertical with the lower middle vertical
followed by composition, which can be
written as 
\begin{equation}\label{eq:trans10} 
V=\{v=(g,f, S(g\otimes f))\} \subseteq R\ts_P R \ts R, \quad S\in\pi_2
\end{equation}
while a similar model for the top right hand corner is
\begin{equation}\label{eq:trans11}
(h,g,f, S_1(h\otimes g), S_2(g\otimes f), S_3(h\otimes(g\otimes f)))
\in R\ts_P R\ts_P R \ts R^3,\quad (S_i)\in\pi_2^3
\end{equation}
As such, the aforesaid intermediary transgression of the
identity in $\Hom(\pi_2,\pi_2)$ is, regardless of the
connected component of $R_0\ts_P R_0$ the map $v\mpo S$
in the notation of (\ref{eq:trans10}) and its differential
with respect to the 4 top most maps in (\ref{eq:trans9})
is exactly the difference between $h\otimes(g\otimes f)$
and $(h\otimes g)\otimes f)$, {\it i.e.} $K_3$ by \ref{fact:pos1}.(c).       
\end{proof}
This spectral sequence interpretation 
can be a rather efficient
way to calculate the invariants
that we've encountered in \ref{eg:eg3},
\ref{eg:eg++1}, and \ref{eg:full}, as
illustrated by
\begin{example}\label{eg:riemman}
{\bf The 2-category $\et_2(\cX)$} whenever $\cX$
is an \'etale fibration over a complex 1-dimensional
orbifold. In the particular case that $\cX$ is an 
\'etale fibration over an orbifold $\cO$ with 
$-\infty\leq\chi(\cO)\leq 0$, orbifold uniformisation
gives that $\cO$ is uniformised by the disc, $\D$, should
the Euler characteristic be strictly negative, and
by the complex plane $\bc$ otherwise. As such $\cO$
is a $\rK(\pi,1)$, and a fortiori so is $\cX$. Consequently
this case is exhaustively described in purely group
theoretic terms by \ref{eg:eg2}, albeit that we can
usefully add that the conformal structure on $\cX$
is inherited from that of the universal cover because
the representations of op. cit. are not arbitrary
topological automorphisms but conformal ones. The
remaining simply connected orbifolds are therefore
parabolic, {\it i.e.} their Euler-characteristic
is positive, and are classifying champs $[T/\bg_m]$
for $T=\bc^2\bsh{0}$ with action as per (\ref{eq:eg11}),
so $\bp^1_\bc$ if $p=q=1$, tear-drops if $p=1, q>1$,
and american footballs otherwise. 
In particular, just as the actions in the
previous case were conformal, the representations
of \ref{eg:full} in the automorphism group of
$T$ are algebraic, {\it i.e.} $\mathrm{GL}_2(\bc)$.
In any case, by
the long-exact sequence of a fibration, \ref{fact:fib1},
any parabolic orbifold has 
(canonically since everything is conformal cf. (\ref{eq:eg10})
et seq.) 
$\pi_2\xrightarrow{\sim}\bz(1)$,
while any \'etale fibration over the same must 
have $\pi_2\xrightarrow{\sim}\bz(1)$ or $0$ by \ref{eg:eg++1}.
At least for separated champs, therefore, we're
in the situation of \ref{eg:full}. Several further
simplifications are, however, possible. In the
first place $\pi_2(\cO)$ is generated by a 
conformal mapping  ($z\mpo z^{pq}$ lifts to a
map $\bp_\bc^1\ra\cO$ of a simply connected
orbifold), while the orbifold fundamental 
group acts by conformal mappings of finite
order, so the action of $\pi_1$ on $\pi_2$ is
trivial. Secondly, the integral co-homology of
an orbifold can be easily computed from the
Leray spectral sequence,
\begin{equation}\label{eq:r1}
\rH^i( \vert\cO\vert, R^j\mu_* \bz(1))\Rightarrow
\rH^{i+j}(\cO, \bz(1))
\end{equation}
of its moduli map $\mu:\cO\ra\vert\cO\vert$. In
the case at hand $\vert\cO\vert\xrightarrow{\sim}\bp^1_\bc$,
and we have some finite (actually at most $3$)
non-scheme like points on $\cO$ with local cyclic
monodromies of order $n_k\in\bn$ at points $p_k:\rp\ra\bp^1_\bc$, so
for $j>0$,
\begin{equation}\label{eq:r2}
R^j\mu_* \bz(1)=\begin{cases}
\coprod_k (p_k)_* \bz/n_k, &\text{$i$ even}, \\
0 &\text{$i$ odd}
\end{cases}
\end{equation}
since, \cite{brown}[VI.9.2], cyclic groups have
2-periodic co-homology, and whence,
\begin{equation}\label{eq:r3}
\rH^n(\cO,\bz(1))=\begin{cases} 
\coprod_k  \bz/n_k, &\text{$n\geq 4$ even}, \\
0 &\text{$n$ odd}
\end{cases}
\end{equation}
Combining this with the Hoschild-Serre spectral
sequence and \ref{fact:pos1} gives an exact sequence
\begin{equation}\label{eq:r4}
0\ra \mathrm{Tors}(\mathrm{Pic}(\cO))
\ra
\mathrm{Pic}(\cO) \ra \mathrm{Pic}(\cO_1)\, (=\bz)
\xrightarrow{k_3}
\rH^3(\pi_1(\cO), \bz(1)) \ra 0
\end{equation}
and, of course, (\ref{eq:r1}), describes the Picard
group by way of
\begin{equation}\label{eq:r5}
0\ra\mathrm{Pic}(\vert\cO\vert)\, (=\bz)\ra 
\mathrm{Pic}(\cO)\ra
\coprod_k  \bz/n_k\ra 0
\end{equation}
from which the middle map in (\ref{eq:r4}) is surjective
if $\pi_1(\cO)$ is cyclic, or a dihedral group $\rD_{2n}$
with $n$ odd, and of index 2 otherwise, {\it i.e.}
$\rD_{2n}$, $n$ even, tetrahedral, $\rT$, octahedral,
$\rO$, or icosahedral $\rI$. Consequently, in the former cases
$\rH^3(\pi_1,\bz(1))$ vanishes, and in the notation of
\ref{eg:eg3} and \ref{eg:eg++1}, by 
(\ref{eq:eg+7}) \& (\ref{eq:eg+13}), associated
to the exact sequence (\ref{eq:eg+1}) there is
a non-unique extension $E'$ of $\pi'_1$ by $\G'$ fitting 
into a diagram
\begin{equation}\label{eq:r6}
\begin{CD}
@. 0@. 0@. 0 @. \\
@. @VVV @VVV @VVV @.\\
0@>>>\pi'_2@>>>\G' @>>> \G'/\pi'_2 @>>> 0 \\
@. @| @VVV @VVV @. \\
0@>>>\pi'_2@>{\mathrm{central}}>> E' @>>> \pi_1(q')@>>> 0 \\
@. @. @VV{q'}V @VV{q'_1}V @. \\
@. @. \pi'_1 @= \pi'_1 @. \\
@. @. @VVV @VVV @. \\
@. @. 0@. 0 @.
\end{CD}
\end{equation} 
If we can compare this with (\ref{eq:fullPlus}), then
the key difference is that 
given the exact sequence (\ref{eq:eg+1})
the extension $E(q')$ of
op. cit. is unique up to isomorphism but the
map $q'$ can vary, whereas in (\ref{eq:r6}) $E'$
is non-unique but once it is fixed $q'$ cannot
vary. Now modulo the distinction (\ref{eq:groupTom})
between equivalence in monoidal as opposed to
2-categories, the respective ways of
varying the respective diagram are  principal
homogeneous spaces under $\rH^2(\pi'_1,\pi'_2)$.
This suggests that if we define a 2-category
$\cE'$ with
\begin{enumerate}
\item[0]-cells diagrams of the form (\ref{eq:r6}), 
equivalently: 
extensions, $E'$, by a cyclic group- $\pi'_2$-
of an extension- $\pi_1(q')$- of $\pi'_1$.
\item[1]-cells pairs $(f,\eta)$ where $f:E'\ra E''$
is a map of extensions of extensions; $\eta\in\pi_1$
$q''f=\mathrm{Inn}_\eta q'$.
\item[2]-cells $\z:(f,\eta)\Rightarrow (g,\xi)$ elements
of $E''$ such that $\xi=q''(\zeta)\eta$.
\end{enumerate}
then we ought, and in fact do, have  
\begin{fact}\label{fact:r1} {\em cf. \cite{noohiuni}[\S 9.1] }
If $\pi_1$ is cyclic (whence $\cO_1$ is arbitrary), or
$\rD_{2n}$, $n$-odd (so $\cO_1\xrightarrow{\sim}\bp^1_\bc$)
then 
we have an an equivalence of 2-categories via
the 2-functor $\cE'\ra \cE$, for $\cE$ of \ref{eg:full}
the 2-category equivalent to $\et_2^{\mathrm{sep}}(\cO)$, given
on $0$-cells by
\begin{equation}\label{eq:r7}
E'\mpo \{\bg_m\coprod_{\pi'_2} E' \xrightarrow{q'} E:= \bg_m\ts\pi_1(\cO) \}
\end{equation}
\end{fact}   
\begin{proof}
The co-homology of finite groups is torsion, and
for any $\ell\in\bn$
\begin{equation}\label{eq:r8}
\rH^2(G, \bz(1)/\ell ) \ra \rH^2(G,\bg_m) \xrightarrow{\ell} \rH^2(G,\bg_m)
\end{equation}
is exact for any group $G$. Consequently the
source of any 0-cell $q'$ in $\cE$ comes from
an extension of the form (\ref{eq:r6}), and
since $k_3(\cO)=0$ its sink is the trivial
extension $E$ of (\ref{eq:r7}). Better still
the 0-cells on either side of (\ref{eq:r7})
are, as we've said, principal homogeneous spaces
under $\rH^2(\pi'_1,\pi'_2)$, and while a simple
diagram chase shows that these homogeneous
space structures are compatible, 
the $\pi'_i$ are finite
so by counting the 
2-functor (\ref{eq:r7}) is essentially 
surjective on 0-cells. As to the 1-cells,
if we examine the proof of the equivalence
of $\cE$ and $\et_2^{\mathrm{sep}}$ in 
\ref{eg:full}, then we can take the $\om$
post (\ref{eq:full5}) to belong to any
section of the map $E\ra\pi_1$ of the
base extension to the base fundamental
group. In general such a section is no
better than set theoretic, but in the
present case $E$ is split, so in the 
definition of the 1-cells, $(f,\eta)$ of $\cE$
we can suppose that $f$ is a map of 
extensions $f:E(q')\ra E(q'')$ and
$\eta\in \pi_1$. As ever, in principle,
there may be more 1-cells in $\cE$ than
$\cE'$ since  the map $f$- 
notation as in \ref{eg:full}-
is synonymous with a
1 co-chain $\psi:\pi_1(q')\ra\Lambda_{\pi_2(q'')}$
whose differential is $(f_{\mathrm{toric}})_* K_3(q')
-(f_1)^* K_3(q')$ upon identifying the respective
Postnikov classes with torus valued 2 co-cycles.
We have, however, the commutativity condition
$\mathrm{Inn}_\om q'=q''f$, so the image of
$\psi$ must take values in the kernel of
$\Lambda_{\pi_2(q'')}\ra \Lambda_{\pi_2}$,
{\it i.e.} $\pi''_2$, and whence  (\ref{eq:r7})
is also essentially surjective on 1-cells.
Finally, the pre-image of $1\times \pi'_1\subset E$
in say $E(q'')$ is $E''$ by construction, so
yet another application of \cite[1.5.13]{tom}
gives an equivalence of 2-categories.
\end{proof}
In the remaining cases 
$\cO_1\xrightarrow{\sim}\bp^1_\bc$,
so $\pi_1=\pi_1(\cO)$ is a sub-group of
$\mathrm{PGL_2}(\bc)$, and we
have 2:1 liftings, $\pi^+_1:=\pi^+_1(\cO) \ra\pi_1(\cO)$, to 
subgroups of $\mathrm{SL_2}(\bc)$,
{\it i.e.} the binary-dihedral group $\rD^+_{2n}$,
albeit $n$ even, the binary tetrahedral group,
$\rT^+$, the binary octahedral group, $\rO^+$,
and the binary icosahedral group $\rI^+$. As such, the tautological
bundle on $\bp^1_\bc$ admits a descent datum
for $\pi^+_1$,  so that by (\ref{eq:r4}) the
Postnikov classes $k_3\in \bz/2$ trivialise
on pull-back and the base of the 2-category
$\cE$, or better $\cE_\cO$ is the extension
\begin{equation}\label{eq:r9}
E_\cO= \bg_m\ts \pi_1^+(\cO)/(-1,-1)
\end{equation}
where, of course, $\cO$ itself is the classifying
champ 
$[T/E_\cO]$ under the action given by sending
$\bg_m$ to a diagonal matrix and $\pi_1^+$ to
the aforesaid representation in $\mathrm{SL_2}(\bc)$.
Now, despite the fact that
all that is at stake is $\pm 1$, this case is
quite different from \ref{fact:r1}. There
is, however, an intermediate sub-case corresponding
to champs whose universal cover is a locally
constant gerbe over $\bp^1_\bc$ of odd degree,
{\it i.e.} the classifying champ of the action
(\ref{eq:eg11}) for $p=q$ odd. In this situation
the 2-type of the pointed stabiliser, (\ref{eq:eg++}),
is $(\pi'_1, \pi'_2, 0)$ for $\pi'_2$ an odd 
cyclic group. Consequently, corresponding to
the exact sequence (\ref{eq:eg+1}) there is
again a diagram exactly as per (\ref{eq:r6}),
and a 2-category $\cE'_{\mathrm{odd}}$ defined
in exactly the same way as $\cE'$ of \ref{fact:r1}
except that we insist that $\pi'_2$ is odd.
There is also a 2-functor $\cE'_{\mathrm{odd}}\ra \cE=\cE_\cO$,
but its definition is a little more complicated,
to wit: start with an extension, $E'$, or more
correctly a diagram (\ref{eq:r6}) and form the
toric extension
\begin{equation}\label{eq:r10}
E^+:= (\bg_m\coprod_{\pi'_2} E')\ts_{\pi_1}\pi^+_1
= \bg_m \coprod_{\pi'_2} (E'\ts_{\pi_1}\pi^+_1)
\end{equation}
which, from $E'\ra\pi'_1$ in (\ref{eq:r6}), also comes equipped with
a map $q^+:E^+\ra \bg_m\ts\pi^+_1$. Further, from
say the first description in (\ref{eq:r10}), the
central element $-1$ in $\pi^+_1$ lifts to a 
central element of order 2 in $E^+$, while the
2-torsion element in the torus (which, by the
way, is not in $\pi'_2$ since this is odd) is also
central. As such, we get an order 2-element $(-1,-1)\in E^+$
which goes to the element $(-1,-1)$ of (\ref{eq:r9})
under $q^+$ since $q^+_{\mathrm{toric}}$ has odd
degree, so,  in total, a ``flip'':
\begin{equation}\label{eq:r11}
\xy  
 (-10,0)*+{\bg_m\coprod_{\pi'_2} E'/\bg_m\ts\pi_1 }="A"; 
 (14,12)*+{E^+/\bg_m\ts\pi^+_1}="B";  
 (38,0)*+{E(q')/E_\cO}="C"; 
{\ar_{}^{} "B";"A"};  
    {\ar^{} "B";"C"};  
    {\ar@{-->}^{} "A";"C"};  
\endxy
\end{equation}
Consequently, and slightly less trivially
than \ref{fact:r1}, we obtain
\begin{fact}\label{fact:r2}
Let $\cO$ be an orbifold other than those
of \ref{fact:r1}, {\it i.e.} $[\bp^1_\bc/G]$
for $G$ any of $\rD_{2n}$, $n$ even, $\rT$,
$\rO$ or $\rI$, and $\et_2^{\mathrm{sep}}(\cO)_{\mathrm{odd}}$
the full sub 2-category of $\et_2^{\mathrm{sep}}(\cO)$
whose universal cover is a locally constant gerbe
over $\bp^1_\bc$ in odd (and necessarily cyclic) groups, then
the 2-functor $\cE'_{\mathrm{odd}}\ra \et_2^{\mathrm{sep}}(\cO)_{\mathrm{odd}}$
envisaged by (\ref{eq:r11}) and \ref{eg:full}
is an equivalence of 2-categories.
\end{fact}
\begin{proof} As per the proof of \ref{fact:r1}
we get essential surjectivity on 0-cells by
the simple expedient that $\rH^2(\pi'_1,\pi'_2)$
is a finite group. As to the 1-cells, we need
to choose a section of $\pi^+_1\ra\pi_1$ even
to define the 2-functor at this level. Such a
choice
in turn defines a section of $E\ra\pi_1$ via
(\ref{eq:r9}), and as we've said in the
proof of \ref{fact:r1} the $\om$ post (\ref{eq:full5}) 
can be taken to belong to any
section of the map $E\ra\pi_1$, so, obviously,
we insist that this is the choice we've
just made. Now plainly a map $f:E(q')\ra E(q'')$
in the definition of a 1-cell of $\cE$ affords
a unique lifting $f^+:E(q')^+\ra E(q'')^+$,
since $E(q')^+=E(q')\ts_E (\bg_m\ts\pi_1^+)$,
and similarly for $E(q'')$. As before such
a map is synonymous with a 1 co-chain, $\psi$, 
with differential $(f_{\mathrm{toric}})_* K_3(q')
-(f_1)^* K_3(q')$,
albeit
here $\psi$ is a priori from $\pi_1(q')^+:=\pi_1(q')\ts_{\pi_1}\pi_1^+$
to $\bg_m$ rather than $\pi_1(q')$. Nevertheless,
the commutativity condition $\mathrm{Inn}_\om q'=q''f$
forces $\psi$ to be $\pi''_2$ valued, and to
be 0 on $-1\in\pi_1^+$, which combined with
the fact that its differential is pulled back
from $\pi_1(q')^2$ implies that $\psi$ itself
is the pull back of a $\pi''_2$-valued $\pi_1(q')$
co-chain, and we have essential surjectivity
on 1-cells. Finally, on 2-cells, we again 
make use of our choice of section to define
the functor, and whence the pre-image in
$E(q'')$ of elements of our section
$\pi_1\curvearrowright E$ are exactly the
image of 2-cells in  $\cE'_{\mathrm{odd}}$,
so we're done by \cite[1.5.13]{tom}. 
\end{proof} 
Arguably, we already knew this from \ref{eg:eg3},
albeit the difference between the 2-functors in
\ref{fact:r1} and \ref{fact:r2} in going from
the diagram (\ref{eq:r6}) to the Torelli description
\ref{eg:full} is instructive.  
Evidently, several remaining cases can be described
via \ref{fact:r1} or \ref{fact:r2} since it suffices
that the 2-type of the pointed stabiliser,
(\ref{eq:eg++}), is $(\pi'_1,\pi'_2, 0)$ which is
usually the case if $\pi'_1\neq \pi_1$. Otherwise,
however, the Postnikov invariant, $k'_3=-1\in \bz/2$,
and such a simple description is not possible. More
precisely, by definition the diagram (\ref{eq:r6})
does not exist. We can, of course, make a diagram
\begin{equation}\label{eq:r12}
\begin{CD}
@. 0@. 0@. 0 @. \\
@. @VVV @VVV @VVV @.\\
0@>>>\pi'_2@>>>\G' @>>> \G'/\pi'_2 @>>> 0 \\
@. @| @VVV @VVV @. \\
0@>>>\pi'_2@>{\mathrm{central}}>> E^+ @>>> \pi_1(q')^+
:=\pi_1(q')\ts_{\pi_1} \pi^+_1
@>>> 0 \\
@. @. @VV{q'}V @VV{q'_1}V @. \\
@. @. (\pi'_1)^+ @= (\pi'_1)^+:=\pi'_1\ts_{\pi_1} \pi^+_1 @. \\
@. @. @VVV @VVV @. \\
@. @. 0@. 0 @.
\end{CD}
\end{equation}
since the Postnikov class of the orbifold pulls
back to trivial over $\pi^+_1$. The number of
such diagrams is, however, insufficient to 
yield all the 0-cells in $\et_2^{\mathrm{sep}}(\cO)$.
Indeed, much as per (\ref{eq:r1})-(\ref{eq:r4}),
Leray for $[\bp^1_\bc/\pi_1^+]\ra\bp^1_\bc$, 
respectively $\rB_{\pi^+_1}\ra\rB_{\pi_1}$,
yields
$\rH^3(\pi_1^+,\bz(1))=0$, respectively
$\rH^2(\pi_1^+,\bz(1))\xrightarrow{\sim} \rH^2(\pi_1,\bz(1))$,
so for $\pi'_1=\pi_1$ and  $\pi'_2$ even
we have an exact sequence
\begin{equation}\label{eq:r13}
0\ra \bz/2 \ra \rH^2(\pi_1,\pi'_2)\rightarrow \rH^2(\pi_1^+,\pi'_2)
\ra 0
\end{equation}
while the isomorphism class of diagrams (\ref{eq:r12})
are only a principal homogeneous space under
$\rH^2(\pi_1^+,\pi'_2)$, {\it i.e.} half as many
as there are 0-cells in $\et_2^{\mathrm{sep}}(\cO)$
with fixed pointed stabiliser. This may be explained
by observing that while $\bz/2\subset (\pi'_1)^+$ can not
by definition be lifted to a normal sub-group of
order 2 in $E^+$ of (\ref{eq:r12}), this does not
exclude the existence of such sub-groups in
$E^+\coprod_{\pi'_2}\bg_m$, and, even
two such since by (\ref{eq:eg+13}) there
are many cases where there are 2-possibilities for the
Postnikov invariant of 0-cells notwithstanding 
that the exact sequence (\ref{eq:eg+1}) is fixed.
As such, the right objects for expressing 0-cells in the
remaining cases in finite group terms are something
like pairs
\begin{equation}\label{eq:r14}
(E^+\coprod_{\pi'_2}\sqrt{\pi'_2}\ra \pi_1^+\ts \bz(1)/2, G)
\end{equation}
for $G$ a normal sub-group of order 2 of the push out
of $E^+$. Plainly, however, apart from illustrating
that this is the wrong way to describe $\et_2^{\mathrm{sep}}(\cO)$,
such a description adds nothing to what we already know
from \ref{eg:eg3} and \ref{eg:full}, so we ignore it.
\end{example}

\newpage
\section{Pro-2-Galois theory}\label{S:III}

\subsection{Flavours and coverings}\label{SS:III.1}
Even if we were to attempt to axiomatise
the situation,
the 1-Galois axiomatisation of
\cite[V.4]{sga1} 
doesn't seem appropriate, since the
more fundamental idea are ``a site 
with enough points'', and, 
``local connectedness'', albeit
this latter condition can, \ref{ProFF}, be dispensed
with.
As such
we'll proceed as follows:
\begin{flav}\label{flav:funct} Our principle
interest is 
the 2-category $\et_2(\cX)$ of champ proper
(so by definition, \ref{mydefnofproperalg}, separated
and quasi-compact) 
and \'etale over a given
connected (but not necessarily separated) 
algebraic (Deligne-Mumford by definition) champ
$\cX$. 
As such the $0$-cells $q:\cY\ra\cX$ are
proper \'etale maps, while the $1$ and $2$-cells
are exactly
as per (\ref{eq:cor1}) in \ref{defn:cor1}.
On the other hand, we essentially only ever
use the fact- \ref{rmk:funct101}- 
that for some atlas $U\ra\cX$,
the fibre of $q$ is equivalent to $U\uts \cG$
for $\cG$ a finite groupoid.  
Consequently, one might be tempted to say
$\cX$ a champ over an arbitrary site, but,
as we've said,
one needs enough points and, in the first
instance in order to considerably
ease the exposition, we'll suppose locally connected
until otherwise stated, {\it i.e.} \ref{SS:III.9}. 
Whence, 
if one takes the definition of $0$-cells
to be that they are locally products with
a finite groupoid 
(which for brevity we'll refer to as 
finite coverings, or just covering if there
is no danger of ambiguity)
then
the 
ensuing ``pro-finite theory''
will be valid 
in some considerable generality, {\it e.g.}
we don't suppose $\cX$ separated, equally
not necessarily separated formal 
(in the obvious sense implied by 
\cite{formalformal})
Deligne-Mumford
champ work too, as do not necessarily separated
topological champs. Consequently, unless specified
otherwise, $\cX$ is  any of the above (finite)
flavours. Subsequently, however, \ref{SS:III.9},
we'll investigate the situation where $\cG$ is
simply a discrete groupoid,
which again although valid in 
some generality, needs connectedness hypothesis
which are topologically 
wholly reasonable, but algebraically
extremely rare. Irrespectively, the $1$-category
with objects $0$-cells of $\et_2(\cX)$, and arrows
$1$-cells 
of the same
modulo equivalence will be denoted
$\et_1(\cX)$, which, of course,
is exactly the original definition,
\ref{factdef:one1}, of the latter.
\end{flav} 
In the algebraic setting it's
appropriate to tweak the definition of
proper a  little so it more closely mimics
that of fibrations, \ref{defn:fib}, {\it i.e.}
\begin{defn}\label{mydefnofproperalg} A 
separated quasi-compact map 
$p:\cE\ra\cX$ of Deligne-Mumford champs will
be called proper (or valued proper if there
is any danger of confusion) if for any map
$f:\mathrm{Spec}(R)\ra\cX$ from a valuation
ring such that there is a lifting 
$\eta:f_K\Rightarrow \tilde{f}_K$ over the
generic point $K$ of $R$ then there is a
map $F:\mathrm{Spec}(R)\ra\cE$ along with
 natural transformations $\xi:\tilde{f}_K\Rightarrow F_K$,
$\a:f\Rightarrow pF$
such that the diagram
\begin{equation}\label{eq:mydefnofproper}
 \xy
 (-12,0)*+{p\tilde{f}_K}="L";
 (12,0)*+{pF_K}="R";
 (0,10)*+{f_K}="T";
    {\ar@{=>}_{\eta} "T";"L"};
    {\ar@{=>}^{\a_K } "T";"R"};
    {\ar@{=>}_{p\xi} "L";"R"}
 \endxy
\end{equation}
commutes. In particular, therefore, the definition
is closed under base change, and a proper map 
in this sense is
universally closed by \cite[7.3.8]{egaII}. 
\end{defn}
As such,
the only difference with the usual definition is
that we don't insist on finite type, whose 
intervention we prefer to note explicitly.
For example,
let us
observe that 
for 
$q:\cY\ra\cX$ a map of not necessarily
separated
 algebraic champ 
\begin{rmk}\label{rmk:funct101}
\'Etale proper 
(sense of \ref{mydefnofproperalg})
and finite type is
equivalent to
having a local fibre 
$U\uts \cG$
for some
\'etale atlas $U\ra\cX$ with
$\cG$ a (finite) groupoid. 
By way of 
\ref{cor:sep1} 
we already have such an equivalence for proper \'etale maps 
in the topological
case, and
plainly, such a local description of
the fibres implies that $q$ is \'etale
and proper, so the issue is the converse.
This follows, however, with exactly the same proof as in
the topological case ({\it i.e.} the easy proper
one \ref{cor:sep1} rather than the more difficult
\ref{cor:sep3}) provided one knows the representable
case, which 
in turn
is an immediate consequence of
\cite[17.4.1.2]{egaIV3} and
\begin{fact}\label{fact:raynaud}
A quasi finite algebraic space $X/S$ of finite
type over an affine $S=\mathrm{Spec}(A)$ is finite
iff it is proper.
\end{fact}
\begin{proof} The if direction is
trivial, and otherwise we need Zariski's lemma that
$X/S$ factors as an open embedding $X\hookrightarrow X'$
of some finite $X'/S$, which is basically
a special case of \cite[16.5]{L-MB}. However,
op. cit. uses \cite[8.12.6]{egaIV2} instead
of \cite[\S 4, Thm. 1]{raynaud}, {\it i.e.}
they only get finite presentation instead of
finite type, cf. \cite[1.8]{milne}. Given
this we can suppose $X$ and $S$ connected,
and since $X'/S$ is proper, $X/S$  proper
in the sense of \ref{mydefnofproperalg}
implies $X\hookrightarrow X'$ proper in the same
sense, so it's closed by \cite[7.3.8]{egaII},
and whence $X'=X$.
\end{proof}
Related to this are a priori competing definitions of
\'etale, {\it i.e.} 
net, flat and finite type, {\it cf.}
\cite[\S 1, Def. 2]{raynaud},
or \cite[\S 3]{milne}, versus net, flat
and finite presentation,  cf. \cite[17.3.1]{egaIV3},
which are a postiori  
equivalent for the same reasons, {\it i.e.} Raynaud's
version
of Zariski's lemma, cf. \cite[3.14-3.16]{milne},
and the fact, \cite[17.4.1.2]{egaIV3},
that net implies locally quasi-finite.
\end{rmk}
In the same vein let us illustrate the
modus operandi, and
clear up a minor point at the same time,
{\it viz:}
\begin{rmk}\label{rmk:flav11}
For any $1$-cell, $(f,\eta)$ in $\et_2(\cX)$,
\begin{equation}\label{eq:flav1}
 \xy 
 (0,0)*+{\cY}="A";
 (28,0)*+{\cY'}="B";
 (14,-12)*+{\cX}="C";
{\ar_{}^{f} "A";"B"};
    {\ar^{q'} "B";"C"};
    {\ar_{q} "A";"C"};
{\ar@{=>}^{\eta} (10,-6);(18,-6)}
\endxy
\end{equation}
$f$ is itself a finite covering. 
In the
algebraic flavour if one were to
have recourse to the \'etale and proper
definition then even 
with \ref{mydefnofproperalg}
one needs a 2-diagram chase to get
proper, net is trivial, while flatness
needs $q'$ faithfully flat, which
is true by \ref{fact:raynaud}. If,
however, one uses the definition of
covering directly then it's trivial.
\end{rmk}

Similarly, a direct appeal to the
definition of covering implies
\begin{fact}\label{fact:flav1}
Let $q:\cY\ra\cX$ be a covering, then there
a factorisation
$\cY\xrightarrow{p}\cY_1\xrightarrow{r}\cX$
into a locally constant gerbe followed by
a representable covering, 
with uniqueness exactly as per \eqref{UniFactor1}-\eqref{UniFactor4},
whence, since
we're supposing $\cX$ connected, the fibres
of $p$ are isomorphic to some fixed group
$\rB_\G$ and the cardinality of the fibres
of $r$ is constant.
\end{fact}
\begin{proof} One argues 
(in fact more easily since
we're taking covering as a given)
as in 
\ref{claim:sep1} using the diagram (\ref{eq:sep3}).
\end{proof}

\subsection{The fibre functor}\label{III.2}
Observe that in any of the above flavours
\begin{claim}\label{claim:funct}
For $q$, $q'$, $q''$ $0$-cells in $\et_2(\cX)$,
and $F'=(f',\xi'):q'\ra q$, $F''=(f'',\xi''):q''\ra q$ $1$-cells the
fibre product (understood in the $2$-category
sense \ref{def:Fibre2}) $q'\ts_q q''$ exists,
or is empty.
Indeed if the sources of $q$, $q'$, $q''$ are
$\cY$, $\cY'$, $\cY''$ respectively, and $\s:f'g'\Rightarrow f'' g''$
form a (2-commutative) fibre square 
\begin{equation}\label{eq:funct101}
\begin{CD}
\cY'\ts_\cY \cY'' @>>{g''}> \cY''\\
@V{g'}VV @VV{f''}V \\
\cY'@>{f'}>> \cY
\end{CD}
\end{equation}
in champ, then in $\et_2(\cX)$
the (albeit ``a'' is much more accurate) 
fibre product is given by: the
$0$-cell $q'g'$, $1$-cells $G'=(g',1)$, 
$G''=(g'', ((g'')^*\xi'')^{-1} q_*(\s) (g')^*\xi')$,
and $2$-cell $\s$;
provided
it is non-empty.
\end{claim}
\begin{proof} Obviously this is just a diagram chase,
so a sufficiently large piece of paper is all that's
required. Nevertheless, there is something to be
checked here, since $q'g'$ and $q''g''$ are not a priori equivalent,
so we're making a choice which we claim doesn't
matter a postiori.
\end{proof} 
Now quite generally, and regardless of
the flavour, for any map of $\cX$-champs
$f:\cY'\ra\cY$ we have a fibre square
\begin{equation}\label{eq:funct102}
\begin{CD}
\cY' @>{\mathbf{1}\ts f}>{:=\G_f}> \cY'\ts_\cX\cY\\
@V{f}VV @VV{f\ts\mathbf{1}}V \\
\cY@>{\D}>> \cY\ts_{\cX}\cY
\end{CD}
\end{equation}
so the graph, $\G_f$, is representable. Usually,
however, it fails to be an embedding. Indeed,
the following sufficient condition is to all
intents and purposes necessary
\begin{fact}\label{fact:funct101}
The graph $\G_f$ in (\ref{eq:funct102}) is an embedding
if $\cY\ra\cX$ is representable and separated.
\end{fact}
\begin{proof} The diagonal in the bottom row of
(\ref{eq:funct102}) is an embedding iff $\cY\ra\cX$
is representable and separated.
\end{proof}
Now just as the definition of fibre products, \ref{def:Fibre2},
in 2-categories is rather different from that in
1-categories so too is the definition of equalisers,
{\it viz:}
\begin{defn}\label{defn:funct101}
Let  $f,g:\cY\ra \cY''$ be $1$-cells in a $2$-category
in which all $2$-cells are invertible,
then an equaliser 
\begin{equation}\label{eq:funct103}
\cY'\xrightarrow{e} \cY {\build\rras_{f}^{g}} \cY''
\end{equation}
of the pair $(f,g)$ is a  $1$-cell
$e:\cY'\ra \cY$ and a $2$-cell $\e:fe\Rightarrow ge$
with the following universal property: if
$h:T\ra \cY$ is a 1-cell such that $\theta:fh\Rightarrow gh$
for some $2$-cell then there is a $1$-cell
$t:T\ra \cY'$ and a $2$-cell $\tau:h\Rightarrow et$
such that we have a commutative diagram
\begin{equation}\label{eq:funct104}
 \xy
 (0,0)*+{fh}="A";
 (20,0)*+{gf}="B";
 (20,-18)*+{get}="C";
 (0,-18)*+{fet}="D";
    {\ar@{=>}_{\theta} "A";"B"};
    {\ar@{=>}_{f_*\tau} "A";"D"};
    {\ar@{=>}^{g_*\tau} "B";"C"};
    {\ar@{=>}^{t^*\e} "D";"C"};
 \endxy
\end{equation}
and for $(s,\s)$ any another such pair there
is a unique $2$-cell $\xi:t\Rightarrow s$ such that
$\s=(e_*\xi)\tau$.
\end{defn}
As ever, the unicity of the 2-cell $\xi$ in 
\ref{defn:funct101} implies that the source
$\cY'$ is unique up to equivalence. Similarly,
one can go from fibre products to equalisers
by way of a fibre square
\begin{equation}\label{eq:funct105}
\begin{CD}
\cY' @>>> \cY\\
@VVV @VV{\G_g}V \\
\cY@>{\G_f}>> \cY'\ts\cY''
\end{CD}
\end{equation}
Consequently, equalisers will fail to be embeddings
in $\et_2(\cX)$ unless \ref{fact:funct101}
applies, {\it i.e.} $q'':\cY''\ra\cX$ is representable,
{\it e.g.} for any 0-cell $q:\cY\ra\cX$, even the equaliser
of the identity with itself fails to be an embedding
if $q$ isn't representable. All of which is the
primary difficulty in establishing the pro-representability
of 
\begin{defn}\label{def:funct102}
the (or more correctly a) fibre functor is the 2-functor $\et_2(\cX)\ra\Grpd$
given on $0$-cells by,
\begin{equation}\label{eq:funct106}
q\mpo q^{-1}(*)
\end{equation}
for some point $*:\rp\ra\cX$, where $\rp$ is understood
to be $\mathrm{Spec}(K)$ for $K$ an algebraically 
(or even just separably) closed
field in the algebraic flavour. 
\end{defn}
Indeed the goal of \cite[V.4]{sga1} is the abstraction
of Galois theory per se, and uses heavily, \cite[expos\'e 195, \S 3]{fga2},
that equalisers are embeddings. In our ad-hoc flavoured
presentation however, the salient point 
is rather mundane, {\it i.e.}
\begin{fact}\label{fact:funct102}
Irrespective of the flavour,
every point, $*_\cY:\rp\ra\cY$ of every 
$0$-cell $q:\cY\ra\cX$ is contained in a unique
connected component. Similarly, if $\cY$
is connected then the equaliser of any pair
of maps to a representable $0$-cell is either
empty or $\cY$ itself.
\end{fact}
\begin{proof} Suppose otherwise, then in (\ref{eq:funct101})
we can take $\cY'\ra \cY$, $\cY''\ra\cY$ to be open
and closed embeddings through which (up to equivalence)
$*_\cY$ factors. Consequently the fibre product of
op. cit. is non-empty, and an open and closed 
embedding in say $\cY'$ by base change, whence $\cY'=\cY''$.
As to the similarly, the equaliser in the representable
case is again an open and closed embedding by (\ref{eq:funct105}) so
if it's non-empty it must be everything.
\end{proof}
Now to apply this to the pro-representability of
the fibre functor in $\et_1(\cX)$ rather than
$\et_2(\cX)$ is straightforward. One first chooses
a set (which is plainly
possible in the pro-finite setting, but,
\ref{fact:395},
is equally true more generally)
of representatives of connected representable
$0$-cells, and similarly
for each 
connected cell $q:\cY\ra\cX$ make a further choice of a set of
representatives of the moduli $q^{-1}(*)$,
{\it i.e.} of pairs $(*_\cY:\rp\ra\cY, \phi)$, $\phi:*\Rightarrow q(*)$,
modulo equivalence- cf. immediately post \ref{prop:cor1}.
Consequently we have a set, $I_{*}$, of such triples
$(q_i, *_i, \phi_i)$ which we partially order
according to the relation $i>j$ if there exists
a 1-cell $F_{ji}:=(f_{ji},\xi_{ji}):q_i\ra q_j$ in $\et_2(\cX)$
and a natural transformation $\zeta_{ji}: *_j \Rightarrow
f_{ji}(*_i)$ such that we have a commutative diagram 
\begin{equation}\label{eq:funct1406}
 \xy 
 (0,0)*+{*}="A";
 (20,0)*+{q_j(*_j)}="B";
 (20,-18)*+{q_jf_{ji}(*_i)}="C";
 (0,-18)*+{q_i(*_i)}="D";
    {\ar@{=>}_{\phi_j} "A";"B"};
    {\ar@{=>}_{\phi_i} "A";"D"};
    {\ar@{=>}^{q_j(\zeta_{ji})} "B";"C"};
    {\ar@{=>}_{\xi_{ji}(*_i)} "D";"C"};
 \endxy 
\end{equation} 
As such 
for any 0-cell $q$,
we get a 2-commutative diagram
\begin{equation}\label{eq:funct1106}
 \xy 
 (0,0)*+{\Hom_{\et_2(\cX)}(q_j,q)}="A";
 (44,0)*+{\Hom_{\et_2(\cX)}(q_i,q)}="B";
 (22,-20)*+{q^{-1}(*)}="C";
{\ar_{}^{F^*_{ji}} "A";"B"};
    {\ar^{} "B";"C"};
    {\ar_{(y_j,\eta_j)\mpo (y_j(*_j), \eta_j(*_j)\phi_j)} "A";"C"};
{\ar@{=>}^{y_j(\zeta_{ji})} (16,-10);(28,-10)}
\endxy
\end{equation}
in which everything is a groupoid, so on
taking moduli (\ref{eq:funct1106})
commutes and on restricting our attention
to representable cells we obtain
\begin{fact}\label{fact:funct103}
The 
moduli of the
fibre functor $\et_1(\cX)\ra\Ens$,
is pro-representable.
More precisely for $I_{*}$ the above partially ordered set
with the further proviso that $q_i$ is representable
we have an isomorphism
\begin{equation}\label{eq:funct1206}
\varinjlim_i\Hom_{\et_1(\cX)} (q_i, q)\ra \vert q^{-1}(*)\vert:
(y_i,\eta_i)\mpo (y_i(*_i), \eta_i(*_i)\phi_i)
\end{equation} 
\end{fact} 
\begin{proof} This is just \cite[expos\'e 195, \S 3]{fga2}
written out ad hoc- surjectivity is the first
part of \ref{fact:funct102}, while, 
the more subtle
injectivity is 
where we use representability via
the similarly
part of {\it op. cit.} since (\ref{eq:funct1406}) and the interchange
rule guarantee a non-empty equaliser.
\end{proof}
Now while we're busy adapting/plagiarising \cite[expose\'e V]{sga1},
it's opportune to recall
\begin{fact}\label{fact:plag1}
Let $q':\cY'\ra\cX$ be a representable $0$-cell with $\cY'$ connected,
then there is a $1$-cell  $q\ra q'$ such that $q$ is 
representable and Galois,
{\it i.e.} for any $1$-cells $F,G:q''\ra q$, there is 
an automorphism $A:q\ra q$ such that we have a
$2$-commutative diagram
\begin{equation}\label{eq:plag1}
 \xy 
 (0,0)*+{q}="A";
 (28,0)*+{q}="B";
 (14,12)*+{q''}="C";
{\ar_{}^{A} "A";"B"};
    {\ar^{F} "C";"B"};
    {\ar_{G} "C";"A"};
{\ar@{=>}^{} (10,6);(18,6)}
\endxy
\end{equation}
\end{fact}  
\begin{proof} We can simultaneously avoid plagiarism, and
do something useful by proving this in a way that avoids
recourse to the finiteness of the fibres. To this end 
following \cite[V.4.g]{sga1}
let
$J$ be a set of representatives of the fibre, 
and $(\cY')^J$ the 
source of the $J$-fold fibre product $(q')^J$
in $\et_2(\cX)$
- so, in general, we need fibre
products up to the cardinality of $J$- and $q:\cY\ra\cX$ the
candidate of {\it op. cit.} for a Galois object,
{\it i.e.} the connected component of 
the moduli $*_{\cY}$ of
$j_{(j\in J)}$ mapping to $\cX$ by some, but not necessarily 
any, {\it cf.} \ref{claim:funct}, projection composed
with $q'$.
Now any other point of the fibre
of $\cY$ over $*$  is
equivalent to some $\s(j)_{(j\in J)}$,
and whose moduli we denote by $*'$,
where $\s:J\ra J$ is injective- otherwise 2-distinct
projections would have a non-empty equaliser to
a representable cell- so certainly,  \cite[V.4]{sga1},  
$s$ is an
isomorphism if $J$ is finite. It is however, a
permutation in general, since if we consider the
fibre product
\begin{equation}\label{eq:plag2}
\begin{CD}
\cJ @>>> \cY'\\
@V{r}VV @VV{q'}V \\
(\cY')^J\hookleftarrow \cY @>q>> \cX
\end{CD}
\end{equation}
then the projections, $p_j$, $j\in J$, of $\cY$ to $\cY'$
furnish sections, $s_j$ of $r$, and by the choice of
$*_{\cY}$ there is exactly one such section for every
point of the fibre, {\it i.e.}
\begin{equation}\label{eq:plag3}
\coprod_{j\in J} \cY \xrightarrow{\sim} \cJ:y\mpo s_j(y)
\end{equation}
is an isomorphism. By construction, however, the 
value of $s_j$ in $*'$ is $\s(j)$, and since by
the definition of fibre products $\cJ_{*'}$ is
naturally isomorphic to $J$, $\s$ is a permutation of $J$.
As such there is a $1$-cell 
$A_{\s}:(q')^J\ra (q')^J$ 
in  $\et_2(\cX)$ (but possibly not in 
$\underline{\mathrm{Cham}}\mathrm{p}\underline{\mathrm{s}}$)
arising from
permutation by $\s$ which sends $j_{(j\in J)}$ to
$\s(j)_{(j\in J)}$ (as opposed to $*_\cY$ to $*'$
which is all that can
be guaranteed in $\underline{\mathrm{Cham}}\mathrm{p}\underline{\mathrm{s}}$).
Now certainly
by \ref{fact:funct102} this restricts to an
automorphism of $\cY$, which by construction
is in fact an automorphism of $q$. Better still,
however, if 1-cells $F,G:q''\ra q$ are given,
then by the simple expedient of choosing a
point, $y$, in the fibre of $q''$ we find
a permutation of $J$ such that $A_\s G(y)$
and $F(y)$ are equivalent points of $q^{-1}(*)$,
whence their equaliser 
in $\et_2(\cX)$ is non-empty while $q$
is manifestly representable, so by \ref{fact:funct102} 
and the definition \ref{defn:funct101} we get
the 2-commutative diagram (\ref{eq:plag1}).
\end{proof}
As such, from \ref{fact:funct103} and \ref{fact:plag1}
or their proofs we obtain
\begin{summary}\label{sum:plag1}
(a) The set $I_{*}$ of \ref{fact:funct103} can be simplified
to a set of representatives of isomorphism classes of
(representable) Galois objects, $q_i$, in $\et_1(\cX)$, albeit 
for each $i$ we 
continue to choose a point $(*_i, \phi_i)$ in the
fibre $q_i^{-1}(*)$ so as to afford the isomorphism
(\ref{eq:funct1206}).

(b) For each Galois object $q_i$, the isomorphism
(\ref{eq:funct1206})
affords an isomorphism 
\begin{equation}\label{eq:plag4}
\Aut_{\et_1(\cX)}(q_i)\xrightarrow{\sim} q_i^{-1}(*)
\end{equation}
which for $i>j$ form a commutative diagram
with left hand side a map of groups
\begin{equation}\label{eq:plag5}
\begin{CD}
\Aut_{\et_1(\cX)}(q_i)@>{\sim}>> q_i^{-1}(*)\\
@VVV @VVV \\
\Aut_{\et_1(\cX)}(q_j)@>{\sim}>> q_j^{-1}(*)
\end{CD}
\end{equation}

(c) In particular both sides of (\ref{eq:funct1206})
are naturally principal homogeneous spaces under
the pro-finite (resp. pro-discrete, 
\ref{fact:395},
in  
the envisaged, 
\ref{flav:funct}, 
extension of
the pro-finite theory) 
group $\pi_1(\cX_*):=\varprojlim_i \Aut_{\et_1(\cX)}(q_i)$
defined via (\ref{eq:plag5}), and (\ref{eq:funct1206})
is actually an isomorphism of functors with values
in $\Ens_{\pi_1(\cX_*)}$, {\it i.e.} sets
with continuous $\pi_1(\cX_*)$-action. 
\end{summary}
\subsection{2-Galois cells}\label{SS:III.3}
Clearly we have done nothing so far except tweak
\cite[V.4-V.5]{sga1} here and there
without even arriving to a conclusion, and certainly
the axiomatisation \cite[V.4]{sga1} applies directly
in the representable case, \cite{noohi1} to
yield that (\ref{eq:funct1206}) in the form
\ref{sum:plag1}.(c) is an equivalence of categories
between $\et_1(\cX)$ and $\Ens_{\pi_1(\cX_*)}$.
This will, however, all come out in the wash
in addressing the analogous theory in $\et_2(\cX)$,
where the principle difficulty is that
\ref{fact:plag1} fails. Indeed it will emerge,
\ref{warn:noplag1}, that the only
generally valid condition for a 0-cell, $q$, in $\et_2(\cX)$
to be Galois in the sense of (\ref{eq:plag1}) is if
$q$ is representable. By way of terminology, 
in the abstract setting of \cite[expos\'e 195]{fga2},
the 0-cells with connected source
are called minimal, and the appropriate (if
rather different) analogue in $\et_2(\cX)$ is
\begin{defn}\label{def:funct103}
For $q:\cY\ra\cX$ a $0$-cell in $\et_2(\cX)$ let
$\cY\xrightarrow{p}\cY_1\xrightarrow{q_1}\cX$
be its factorisation, \ref{fact:flav1}, 
into a locally constant gerbe
followed by a representable \'etale cover 
then we
say that $q$ is quasi-minimal if $\cY$ is connected and
we have an equivalence of 1-categories
\begin{equation}\label{eq:funct107}
\et_1(\cY_1)\ra \et_1(\cY): \cY'/\cY_1\mpo \cY'\ts_{\cY_1}\cY
\end{equation}
\end{defn} 
Unsurprisingly, therefore, we assert
\begin{claim}\label{claim:funct102}
For any connected $0$-cell $q:\cY\ra\cX$
in $\et_2(\cX)$ there is a representable \'etale covering
$\cY_2\xrightarrow{r}\cY$ such that $qr$ is quasi-minimal.
\end{claim}
\begin{proof}
In the notation of \ref{def:funct103} say $p$ is a
locally constant fibration in $\rB_\G$'s for some
finite group $\G$, and for $r:\cY'\ra\cY$ 
with $\cY'$ connected,
consider
the 2-commutative diagram
\begin{equation}\label{eq:funct108}
\begin{CD}
\cY@<{r}<< \cY'\\
@V{p}VV @VV{p'}V \\
\cY_1@<<{r_1}< \cY'_1
\end{CD}
\end{equation}
where $r_1p'$ is the factorisation of $pr$ into
a locally constant gerbe, in, say $\rB_{\G'}$'s
followed by a representable map. Now the fibre
of $p$ is a $B_\G$, and every connected component
of the fibre of such, {\it i.e.} $r^*\rB_\G$,
is a $\rB_{\G'}$, so $\rB_{\G'}\ra \rB_\G$ is
representable. As such we can identify $\G'$
with a sub-group of $\G$, so replacing $\cY$
by a suitable $\cY'$ we can suppose that for 
every diagram of the form (\ref{eq:funct108}),
$\G'=\G$. Now plainly we have a map, 
$r':\cY'\ra \cY'':=\cY\ts_{\cY_1}\cY'_1$,
while the projection, $r'_1$, of the 
fibre product to $\cY$ is representable,
and $r=r'_1r'$, so $r'$ is representable.
On the other hand, $\cY'$ and $\cY''$ are
both locally constant gerbes over $\cY'_1$,
so $\cY''$ is connected, and $r'$ has degree
1, from which $r'$ is an isomorphism. A priori
this only proves that (\ref{eq:funct107}) is
essentially surjective, but the naturality
of \ref{fact:flav1} and fibre products shows
we have an equivalence of categories.
\end{proof}
Now we need to understand maps from quasi-minimal objects, to wit:
\begin{fact}\label{fact:funct104} (cf. \ref{fact:visit4})
Let $q:\cY\ra\cX$ be quasi-minimal and $(f,\eta):q\ra q'$ a
$1$-cell in $\et_2(\cX)$ then for $\cY\xrightarrow{p}
\cY_1\xrightarrow{q_1}\cX$, $\cY'\xrightarrow{p'}
\cY'_1\xrightarrow{q'_1}\cX$ 
the respective
factorisations as compositions
of representable cells and locally constant gerbes in
$\rB_\G$'s, respectively $\rB_{\G'}$'s, the induced map
$f:\G\ra\G'$ has image a central sub-group of $\G'$, so,
in particular, $\G$ is abelian. 
\end{fact}
\begin{proof}
Taking a sufficiently fine \'etale cover
$U\ra\cY'$ we appeal to (\ref{eq:sep3}) to interpret
$\cX$, $\cY$, $\cY'$ as \'etale groupoids $R_0$, $R$,
$R'$ acting on $U$ and $q'$, $f$, {\it etc.} as
functors. Similarly we confuse $\G'$, $\G$
with trivial normal $U$ sub-groups of the stabiliser,
so that we get a co-cycle
\begin{equation}\label{eq:funct111}
R\ra\mathrm{Aut}(\G'): x\mpo \{\g'\mpo f(x)\g' f(x)^{-1}\}
\end{equation}
or, equivalently, a $\mathrm{Aut}(\G')$-torsor over $\cY$.
By definition \ref{def:funct103} this is the pull-back of
a  $\mathrm{Aut}(\G')$-torsor over $\cY_1$, so it's restriction
to any fibre of $\cY\ra\cY_1$ is trivial which is iff
the image of $\G$ in $\G'$ is central.
\end{proof}
Of which a more pertinent variant is
\begin{fact}\label{fact:funct105}
Let $q:\cY\ra\cX$, $q':\cY'\ra\cX$ be $0$-cells
in $\et_2(\cX)$, with $\cY$ connected, then for
any $1$-cell
$F=(f,\eta):q\ra q'$ and any point $*:\rp\ra\cY$,
the restriction map
\begin{equation}\label{eq:funct112}
\Hom_{\et_2(\cX)} (F,F)\ra \mathrm{Ker}\{\Hom_{\cY'}(f(*),f(*))
\ra \Hom_\cX(q(*), q(*))\}
\end{equation}
is always injective. In general it may not be surjective,
but there is a connected representable cover $r:\cY''\ra\cY$
(independent of the point) such that on replacing $f$ by $rf$ it
becomes so.
\end{fact}
\begin{proof} The generally true injectivity statement is
mutatis mutandis the proof of \ref{fact:point1}. Otherwise
we retake the notations of the proof of \ref{fact:funct104},
since  without loss of generality  $\cY$ is quasi-minimal, albeit
this comes out in the wash since equally without loss
of generality the co-cycle (\ref{eq:funct111}) is a 
co-boundary. On the other hand the right hand side of
(\ref{eq:funct112}) is isomorphic to $\G'$. This isomorphism
is, however, dependent on the choice of an isomorphism of
the relative stabiliser of $q'$ with $U\ts\G'$, while the
triviality of (\ref{eq:funct111}) is equivalent to changing
this (uniquely up to automorphisms of $\G'$) in such a 
way that (\ref{eq:funct111}) is not just a co-boundary, but
identically $\mathbf{1}_{\G'}$, so for every $\g'\in\G'$ we
get a 2-cell $\g':F\Rightarrow F$.
\end{proof}
Similarly we can look at maps to a quasi-minimal
cell, {\it viz:}
\begin{fact}\label{fact:noplag1} 
Let $q:\cY\ra\cX$ be quasi-minimal and $(f,\eta):q'\ra q$ a
$1$-cell in $\et_2(\cX)$ then for $\cY\xrightarrow{p}
\cY_1\xrightarrow{q_1}\cX$, $\cY'\xrightarrow{p'}
\cY'_1\xrightarrow{q'_1}\cX$
the respective
factorisations as compositions
of representable cells and locally constant gerbes in
$\rB_\G$'s, respectively $\rB_{\G'}$'s, with
$\cY'$ connected, the induced map (better
maps since there's actually one for every point)
$f:\G'\ra\G$ is surjective.
\end{fact}
\begin{proof} Without loss of generality we can
suppose that $\cY'_1=\cY_1$, and,
up to the obvious reversal of roles, we re-take
the notation of the proof of \ref{fact:funct104}. 
Now observe that the kernel, $K$, say, of 
$f:R'\ra R$ restricted to the stabiliser is
a normal sub $U$-group of $U\ts \G'$. As such,
it is locally constant, and $[R'/K]\rras U$ 
is an element of $\et_2(\cX)$. Consequently,
we can suppose that $\G'\ra \G$ is injective,
but this is equivalent to $f$ representable,
so $\G'=\G$.
\end{proof}
Finally in these variations let's
clear up any ambiguity about the definition
of $\G$ by way of
\begin{fact}\label{fact:noplag2}
Let $\cY\xrightarrow{p}\cY_1\xrightarrow{q_1}\cX$
be the factorisation of a minimal $0$-cell into 
a locally constant gerbe in $\rB_\G$'s followed
by a representable map then there is a 
representable \'etale cover $\cY'_1\ra \cY_1$
such that for any point $*':\rp\ra\cY':=\cY\ts_{\cY_1}\cY'_1$ there
is on identifying $\G$ 
with a sub-group of the stabiliser of $*'$
a unique (dependent on the point) isomorphism
of the stabiliser $\cS_{\cY'/\cX}\ra\cY$ relative to $\cX$ 
with $\cY'\uts\G$.
\end{fact}
\begin{proof}
In principle the stabiliser $\cS_{\cY/\cX}\ra\cY$
is no better than a locally constant $\cY$-group
with fibre $\G$. However amongst such it's
isomorphism class is the co-cycle (\ref{eq:funct111})
whence by the definition of minimality we can kill
this by a (connected) \'etale cover $\cY'_1\ra\cY_1$,
so there is such a trivialisation, and since $\cY'$
is connected it's unique once we fix a fibre.
\end{proof}
We next aim to understand automorphic
1-cells in $\et_2(\cX)$, beginning
with the wholly general
\begin{claim}\label{claim:noplag1}
Let $F=(f,\eta):q\ra q'$ be a $1$-cell
in $\et_2(\cX)$ with
$\cY\xrightarrow{p}
\cY_1\xrightarrow{q_1}\cX$, $\cY'\xrightarrow{p'}
\cY'_1\xrightarrow{q'_1}\cX$ 
the respective
factorisations as compositions
of representable cells and locally constant gerbes,
and $F_1:q_1\ra q'_1$ the resulting representable
$1$-cell, the space 
(modulo equivalence)
of $1$-cells between $q$ and $q'$ which 
(modulo equivalence) induce $F_1:q_1\ra q'_1$
is a principal homogeneous space
under $\rH^1(\cY,f^*\cS_{\cY'/\cX})$.
\end{claim}
\begin{proof} In the notation of the proof of
\ref{fact:funct104} we interpret everything in
terms of groupoids and functors, so, in particular
the locally constant sheaf structure on $f^*\cS_{\cY'/\cX}$
is given by (\ref{eq:funct111}). Now if $g:R\ra R'$
is a functor such that $\xi:f_1\Rightarrow g_1$ for
some natural transformation, then since 
$R'\ra R'_1$ is a principal homogeneous space
under $\G'$ for a sufficiently fine \'etale atlas
$U\ra\cY$, $\xi$ can be lifted 
to a map from
$U\ra R'$- which, by the way, is irrespective of
whether $R'\ra R'_1$ is trivial over connected
components which isn't a priori true in an
arbitrary Grothendieck topology, or even in
a classical topology if $\cX$ doesn't have a
good cover. In any case, we can therefore suppose
that $g_1=f_1$, so we must have $g(x)=\g_x f(x)$
for all arrows $x$, for some stabiliser $\g_x$ of
the sink of $f(x)$, and for $g$ to be a 
functor is equivalent to $x\mpo \g_x \in f^*\cS_{\cY'/\cX}$
being a co-cycle, while $g$ being equivalent to $f$
is iff this co-cycle is a co-boundary.
\end{proof}
Which in the context of automorphisms of
quasi-minimal cells yields
\begin{claim}\label{claim:noplag2}
Every endomorphism of every quasi minimal
$1$-cell, $q:\cY\ra\cX$, is an automorphism, and if,
moreover \ref{fact:noplag2} holds, 
{\it i.e.} $\cS_{\cY/\cX}$ is trivial, with
$\cY\xrightarrow{p}\cY_1\xrightarrow{q_1}\cX$ a 
factorisation of 
$q$ into a locally constant gerbe in $\rB_\G$'s
followed by a representable map, then the $2$-group
$\fAut_{\et_2(\cX)}(q)$ has
$\pi_2=\G$ while its fundamental group 
fits into an exact sequence
\begin{equation}\label{eq:noplag2}
1\ra \rH^1(\cY_1,\G) \ra \pi_1 \fAut_{\et_2(\cX)}(q)
\ra \pi_1 \fAut_{\et_2(\cX)}(q_1)
\end{equation}
\end{claim}
\begin{proof} By \ref{rmk:flav11} every endomorphism
is a covering and quasi-minimal implies connected,
so it's surjective on points, and by \ref{fact:funct102}
applied to the representable quotient it's also
injective on points, while by \ref{fact:funct105}
and \ref{fact:noplag1} it's an isomorphism on
stabilisers. Consequently, {\it cf.} \ref{rmk:flav11},
it's an isomorphism on fibres over any \'etale atlas
$U\ra\cX$, and whence an equivalence of categories.
As to the rest: by  definition $\pi_2$ of a 2-group
is the stabiliser of the identity, and this is $\G$
by \ref{fact:noplag2}, while its $\pi_1$ is 1-cells
modulo equivalence, so we're done by by \ref{def:funct103},
\ref{fact:noplag2} and \ref{claim:noplag1}.
\end{proof} 
Now currently, we're quite far from being able
to assert surjectivity on the right of
(\ref{eq:noplag2}) since in the notation
of {\it op. cit.}, albeit independent of
any quasi-minimality hypothesis, lifting
an automorphism $F_1=(f_1,\eta_1):q_1\ra q_1$
is (since lifting natural transformations
is a non-issue, {\it cf.} the proof 
of  \ref{claim:noplag1})  
equivalent to the existence of a 2-commutative
fibre square
\begin{equation}\label{eq:noplag3}
\begin{CD}
\cY@<<< \cY \\
@V{p}VV @VV{p}V \\
\cY_1@<{f_1}<< \cY_1
\end{CD}
\end{equation}
Which for ease of reference we encode in
\begin{defn}\label{def:noplag1}
If $q=q_1p$ is the factorisation of a $0$-cell
in $\et_2(\cX)$ into a representable map and
a locally constant gerbe, then we say that
$q$ is quasi-Galois if 
$q_1$ is Galois and
for every automorphism, $F_1$
of $q_1$ in $\et_2(\cX)$ we have the $2$-commutative
fibre square (\ref{eq:noplag3}).
\end{defn}
From which a couple of easy consequences are
\begin{fact}\label{fact:noplag3}
\ref{claim:funct102} holds
with ``quasi-minimal'' replaced by ``quasi-Galois''.
\end{fact}
\begin{proof} By \ref{fact:plag1}, we can, without
loss of generality suppose that $q_1$ is Galois, so
if we choose a complete repetition free list 
$\{F_\om=(f_\om, \eta_\om):\om\in\pi_1(\fAut(q_1))\}$ of 1-cells
of $q_1$ lifting the equivalence classes $\pi_1(\fAut(q_1)) $,
then the fibre product relative to $q_1$ of all
of the $F_\om^* q$ does the job.
\end{proof} 
\begin{fact}\label{fact:noplag4}
Let $q=q_1p$ be the factorisation of a 
quasi-Galois $0$-cell
in $\et_2(\cX)$ into a representable map and
a locally constant gerbe, then for any
$1$-cell $G:r_1\ra q_1$ with $r_1$
Galois and representable, $G^*q$ is quasi-Galois.
\end{fact}
\begin{proof} For any automorphism $H_1:r_1\ra r_1$
we have by \ref{sum:plag1}.(b) a 2-commutative 
fibre square
\begin{equation}\label{eq:noplag4}
\begin{CD}
r_1 @>>{H_1}> r_1 \\
@V{G}VV        @VV{G}V \\
q_1 @>{F_1}>> q_1
\end{CD}
\end{equation}
for some automorphism $F_1$ of $q_1$, so
we conclude
from the obvious 2-commutative cube.
\end{proof}
If, furthermore, we combine the notions
of quasi-Galois and quasi-minimal then
\begin{fact}\label{fact:noplag5}
\ref{claim:funct102} holds 
with ``quasi-minimal'' replaced by ``quasi-minimal and quasi-Galois''.
\end{fact}
\begin{proof} By \ref{claim:funct102} (or better
its proof), \ref{fact:noplag3} and \ref{fact:noplag4}
we can suppose that $q$ is quasi-Galois and admits
a representable map $q'\ra q$ from a quasi-minimal cell where
both $q$, $q'$ are locally constant gerbes over 
(some necessarily Galois) $q_1$. Now let $F_1:q_1\ra q_1$
be an automorphism, then since $q$ is quasi-Galois, 
pulling back along $F_1$ yields a representable
map $(F_1)^*q'\ra q$, and we can form a 2-commutative
fibre square
\begin{equation}\label{eq:noplag5}
\begin{CD}
(F_1)^*q'@<<< q'' \\
@VVV @VVV \\
q @<<< q'
\end{CD}
\end{equation}
in which every map is representable. Consequently
by \ref{fact:plag1} and the definition of minimal,
\ref{def:funct103}, there is a representable 
Galois covering $G:r_1\ra q_1$ such that if we form
the 2-commutative diagram of fibred squares
\begin{equation}\label{eq:noplag6}
\begin{CD} 
r_1 @<<< r' @<<< r'' \\
@VVV @VVV @VVV \\
q_1@<<< q' @<<< q''\xrightarrow{\sim} r'
\end{CD}
\end{equation}
then every connected component, $r''_k$,
$k\in K$, of $r''$ is,
via the right top horizontal arrow, isomorphic
to $r'$. Consequently if we look at the base
change of (\ref{eq:noplag5}), {\it i.e.}
the fibre square
\begin{equation}\label{eq:noplag7}
\begin{CD} 
G^*(F_1)^*q'@<<< r'' \xrightarrow{\sim} \coprod_{k\in K} r' \\
@VVV @VVV \\
r @<<< r' 
\end{CD} 
\end{equation}
then we get many 
representable
maps from $r'$ to $G^*(F_1)^*q'$,
and since both of these are quasi-minimal
locally constant gerbes over $r_1$ 
under the same $\rB_{\G'}$ they're
all isomorphisms by \ref{fact:noplag1}.
By the simple expedient of taking 
$G:r_1\ra q_1$ sufficiently large
we can suppose that this holds for
a complete repetition free list of
automorphisms $F_1$ representing
the elements of $\pi_1(\fAut(q_1))$,
while for every $H_1$ in $\pi_1(\fAut(r_1))$
there is a $F_1$ such that (\ref{eq:noplag4})
is a 2-commutative fibre square,
and whence any $(H_1)^*r'$ is isomorphic
to some $G^*(F_1)^*q'$ for general
nonsense reasons.
\end{proof}
This affords some strengthening of
\ref{claim:noplag2} by way of an
exercise in the definitions, {\it i.e.}
\begin{fact}\label{fact:noplagAct}
(cf. (\ref{eq:act1}) et seq.) 
Let everything 
(so in particular trivial stabiliser)
be as in \ref{claim:noplag2}
then if in addition $q$ is quasi-Galois
then (\ref{eq:noplag2}) extends on the
right to a short sequence. Furthermore,
the wholly general action of $\pi_1(\fAut_{\et_2(\cX)}(q))$
on $\pi_2(\fAut_{\et_2(\cX)}(q))$ is, in fact,
afforded by an action of $\pi_1(\fAut_{\et_2(\cX)}(q_1))$
on the same which on choosing a point $*':\rp\ra\cY$
can be written via  the isomorphism of \ref{fact:noplag2}
as
\begin{equation}\label{eq:noplagAct1}
\G=\mathrm{stab}_{\cY/\cX}(*')\ra \mathrm{stab}_{\cY/\cX}(f_\om(*'))
\xrightarrow{\text{\ref{fact:noplag2}}} 
\G: \g\ra f_\om(\g)
\end{equation}
for any $1$-cell $(f_\om, \eta_\om)$ lifting
$\om\in\pi_1(\fAut_{\et_2(\cX)}(q_1))$.
\end{fact}
Nevertheless this is still short of
our final definition, {\it i.e.}
\begin{defn}\label{def:noplag2} A 
quasi-minimal and quasi-Galois $0$-cell
$q$ with trivialisable (relative to $\cX$)  stabiliser is said
to be 2-Galois if the sequence
(\ref{eq:noplag2}) extended to a
short exact sequence is split exact,
or, equivalently 
if $q_1p$ is its factorisation as
a locally constant gerbe over a 
representable Galois cover then
for some (whence any)
complete repetition free list
$\{F_\om:q_1\ra q_1\}$ of $1$-cells
representing the elements $\om\in\pi_1(\fAut(q_1))$
there exists, 
for every pair $\tau,\om\in\pi_1(\fAut(q_1))$,
a $2$-commutative diagram
in which every square is fibred
\begin{equation}\label{eq:noplag8}
 \xy
 (38,0)*+{q}="B";
 (0,-6)*+{q}="C";
(15,0)*+{q}="D";
 (38,-12)*+{q_1}="F";
 (0,-18)*+{q_1}="G";
(15,-12)*+{q_1}="H";
{\ar_{}_{p} "B";"F"};
    {\ar_{p} "C";"G"};
    {\ar^{p} "D";"H"};
{\ar_{}^{F_{\tau} } "F";"G"};
    {\ar@{-->}_{F_{(\tau\om)} } "H";"G"};
    {\ar^{F_{\om} } "H";"F"};
{\ar_{}^{ } "B";"C"};
    {\ar_{ } "D";"C"};
    {\ar^{ } "D";"B"};
\endxy
\end{equation}
\end{defn}
The advantage of the formulation (\ref{eq:noplag8})
is that it  evidences
\begin{fact}\label{fact:noplag5noplag5}
\ref{claim:funct102} holds
with ``quasi-minimal'' replaced by ``2-Galois''.
\end{fact}
\begin{proof} 
By \ref{fact:noplag2} and \ref{fact:noplag5}
we can suppose that the cell is quasi-Galois and
quasi-minimal with trivialisable relative stabiliser. 
As such there are liftings, $\tf_{\bullet}$, of
the 1-cells in the bottom face of \eqref{eq:noplag8}.
Thus, by definition, there are 2-cells,
$\phi_{\bullet}: p\tf_{\bullet}\Rightarrow F_{\bullet} p$,
but a 2 cell filling the top face of {\it op. cit.}
need not exist. Nevertheless, 
if $\a^1_{\tau\om}$ is the 2-cell on the bottom
face, 
there is, \ref{claim:noplag1}, a 
1 cell $Z_{\tau,\om}:q\ra q$ such that $pZ_{\tau,\om}=p$, along 
with a 2-cell $\a_{\tau,\om}:\tf_\tau\tf_\om\Rightarrow
\tf_{(\tau\om)}Z_{\tau,\om}$; and since $q_1$ has no
relative stabiliser, {\it cf.} \eqref{eq:plagAct15},
we have a commutative diagram
\begin{equation}\label{moredetail1}
\xy 
 (0,0)*+{p\tf_{\tau}\tf_{\om}}="A";
 (30,0)*+{p\tf_{(\tau\om)}Z_{\tau,\om}}="B";
 (30,-18)*+{F_{(\tau\om)}p}="C";
 (0,-18)*+{F_{\tau}F_{\om}p}="D";
    {\ar@{=>}_{p_*\a_{\tau\om}} "A";"B"};
    {\ar@{=>}_{(F_{\tau})_*(\phi_\om)(\tf_\om)^*(\phi_\tau)} "A";"D"};
    {\ar@{=>}^{Z_{\tau,\om}^*\phi^{}_{\tau\om}} "B";"C"};
    {\ar@{=>}_{ (p)^*\a^1_{\tau\om}} "D";"C"};
 \endxy 
\end{equation}
Now let $G_1:r_1\ra q_1$ be a Galois cover, so
identifying $\tau,\om$ with elements of $\pi_1(\cX_*)$
we 
have 1-cells $F'_{\tau}$, $F'_\om:r_1\ra r_1$;
2-cells 
$\phi'_\om: G_1 F'_\om\Rightarrow F_\om G_1$;
and we 
can form 2 commutative diagrams
\begin{equation}\label{moredetail2}
 \xy 
 (0,0)*+{r}="A";
 (30,0)*+{q}="B";
 (30,-18)*+{q_1}="C";
 (0,-18)*+{r_1}="D";
 (22, -6)*+{}="E";
 (8,-12)*+{}="F";
 (-20,15)*+{r}="H";
 (-2,4)*+{}="a";
 (2,8)*+{}="b";
 (-2,-2)*+{}="c";
 (-6,-6)*+{}="d";
    {\ar_{G} "A";"B"};
    {\ar^{p'} "A";"D"};
    {\ar^{p} "B";"C"};
    {\ar^{G_1} "D";"C"};
 {\ar@{=>}_{ \b} "E";"F"}; 
 {\ar_{F'_{\om}p'} "H";"D"};
 {\ar^{\tf_{\om}G} "H";"B"};
  {\ar_{\tf'_\om} "H";"A"};
 {\ar@{=>}_{ \psi_\om} "a";"b"}; 
 {\ar@{=>}_{ \theta_\om} "c";"d"}; 
\endxy 
\end{equation}
for any $\om\in\pi_1(\cX_*)$, by applying the universal property, 
\eqref{FibreA3}-\eqref{FibreA5},
to the 2-cell
\begin{equation}\label{moredetail3}
p\tf_\om G 
{\build\Longrightarrow_{}^{G^*\phi_\om}}
F_\om pG
{\build\Longrightarrow_{}^{(F_\om)_*\b}}
F_\om G_1 p'
{\build\Longleftarrow_{}^{(p')^*\phi'_\om}}
G_1 F'_\om p'
\end{equation}
On the other hand, $Z_{\tau,\om}$ is given by
\ref{claim:noplag1}, so for any $G$
sufficiently large independently of $\tau$ or $\om$,
there are 2-cells $\z_{\tau,\om}:G\Rightarrow Z_{\tau,\om}G$,
while $r_1$ is Galois, so there are also 2-cells
$\a'_{\tau,\om}:F'_\tau F'_\om\Rightarrow F'_{\tau\om}$,
and a small diagram chase reveals that if we
take $(\tau\om)$ in \eqref{moredetail2} then
$\tf'_\tau\tf'_\om$ satisfies the universal property
with
\begin{equation}\label{moredetail4}
\begin{split}
\psi_{(\tau\om)}: & \quad G\tf'_\tau\tf'_\om 
{\build\Longrightarrow_{}^{(\tf'_\om)^*\psi_\tau}}
\tf_\tau G \tf'_\om
{\build\Longrightarrow_{}^{(\tf_\tau)_*\psi_\om}}
\tf_\tau\tf_\om G
{\build\Longrightarrow_{}^{G^*\a_{\tau,\om}}}
\tf_{\tau\om} Z_{\tau,\om}G
{\build\Longleftarrow_{}^{(\tf_{\tau\om})_*\z_{\tau,\om}}}
\tf_{\tau\om} G\\
\theta_{(\tau\om)}: & \quad p'\tf'_\tau\tf'_\om 
{\build\Longrightarrow_{}^{(\tf'_\om)^*\theta_\tau}}
F'_\tau p' \tf'_\om
{\build\Longrightarrow_{}^{(F'_\tau)_*\theta_\om}}
F'_\tau F'_\om p'
{\build\Longrightarrow_{}^{(p')^*\a'_{\tau,\om}}}
F'_{\tau\om} p'
\end{split}
\end{equation}
The elements $\tau,\om\in\pi_1(\cX_*)$ were, however,
arbitrary, so $r$ is 2-Galois.
\end{proof}
Similarly, for much the same, and in fact easier 
{\it cf.} \ref{5.4}, reason
\begin{fact}\label{fact:noplag6}
\ref{fact:noplag4} holds with ``quasi-Galois''
replaced by ``2-Galois''.
\end{fact}
As such, let us close this section with the
first of many
\begin{warning}\label{warn:noplag1}
While 2-Galois will prove (and in any case it's
manifestly the best that one can do) to be the right
generalisation of Galois, the 2-Galois cells
are not Galois in the sense of (\ref{eq:plag1}).
Indeed for $q$ 2-Galois, expressed as ever
via a factorisation $q_1p$ of a locally 
constant gerbe in $\rB_\G$'s and a representable map, we
can take any Galois cover $G_1:q''_1\ra q_1$,
look at the pull-back $G:q'':=(G_1)^*q\ra G$,
and take $F:q''\ra q$ in (\ref{eq:plag1}) to
be the composition of $G$ with any automorphism
of $q''$ in the kernel (modulo notation, {\it i.e.}  $q$ replaced by
$q''$) of (\ref{eq:noplag2}). There is, however,
no way that this can be completed to the 2-commutative
triangle (\ref{eq:plag1}) unless (in the obvious notation)
\begin{equation}\label{eq:noplag88}
(G_1)^*:\rH^1(\cY_1, \G) \ra \rH^1(\cY''_1, \G)
\end{equation}
is surjective, which, in turn, is very often false.
Indeed, the only generally true way to guarantee this
is if $\G$ is trivial, but then $q=q_1$ is representable
and Galois.
\end{warning}

\subsection{\texorpdfstring{$\pi_2$}{pi\_2} and its Galois 
module structure}\label{SS:III.4} 

If $q_1$ is a representable $1$-cell, then the
definition of $\et_2(\cX)$ is such that any
lifting of the natural surjection
\begin{equation}\label{eq:plagAct1}
\fAut_{\et_2(\cX)}\ra \pi_1(\fAut_{\et_2(\cX)}(q_1))\ra 1
\end{equation}
defines an action of 
$ \pi_1(\fAut_{\et_2(\cX)}(q_1))$ on $q_1$.
Indeed, if we adapt
the notation of \ref{factdef:cover1}, and write
$F_\om=(f_\om, i_\om):q_1\ra q_1$ for 1-cells lifting
$\om$ on the right of (\ref{eq:plagAct1})
then if $\a_{\tau,\om}:F_\tau F_\om\Rightarrow F_{(\tau\om)}$
is any 2-cell, (\ref{eq:coverFix1}) must hold, and
since $q_1$ is representable this both
determines $\a_{\tau,\om}$ and
forces it to be a  $ \pi_1(\fAut_{\et_2(\cX)}(q_1))$ 
2-co-cycle, which is exactly,  {\it cf. op. cit.},
the condition for an action. In the particular case
that $q_1$ is Galois therefore
\begin{fact}\label{fact:plagAct1} Let $q_1:\cY_1\ra\cX$ 
be a representable Galois 1-cell in $\et_2(\cX)$ then,
in the above notation, the action yields a
fibre square in
$\underline{\mathrm{Cham}}\mathrm{p}\underline{\mathrm{s}}$
\begin{equation}\label{eq:plagAct2}
 \xy
 (0,0)*+{\pi_1(\fAut_{\et_2(\cX)}(q_1))\uts \cY_1}="A";
 (38,0)*+{\cY_1 }="B";
 (38,-18)*+{\cX}="C";
 (0,-18)*+{\cY_1 }="D";
  (23,-5)*+{}="E";
  (15,-13)*+{}="F";
    {\ar_{\mathrm{trivial}} "A";"D"};
    {\ar^{\quad\quad f_\om} "A";"B"};
    {\ar^{q_1} "B";"C"};
    {\ar_{q_1} "D";"C"};
    {\ar@{=>}_{i_\om} "F";"E"};
 \endxy
\end{equation}
which is equally, \ref{claim:funct}, interpretable as a fibre
square in $\et_2(\cX)$, and (modulo replacing $\cX_1$ by
$\cY_1$) all of \ref{factdef:cover1} holds, so, in particular,
\ref{sum:plag1}.(b), we have a faithful action in $\et_2(\cX)$.
\end{fact}
\begin{proof} As observed, we have an action,
and since $\pi_1(\fAut_{\et_2(\cX)}(q_1))=\fAut_{\et_1(\cX)}(q_1)$
it's necessarily faithful by \ref{sum:plag1}.(b). Consequently,
the only thing to note is that the square \eqref{eq:plagAct2}
is fibred, but since this maps to an actual fibre square
in which everything is a representable \'etale cover, we're
done by the simple expedient of observing that the fibres
of say the leftmost vertical have the same cardinality as
the rightmost vertical.
\end{proof} 
The situation for non-representable maps is,
however, far more subtle. To fix ideas suppose
that $q:\cY\ra\cX$ is quasi-Galois
with a factorisation $q_1p$ into a locally
constant gerbe in $\rB_\G$'s followed
by a representable map, then the fibre $q^{-1}(*)$
can be identified with the groupoid 
$\pi_1(\fAut_{\et_2(\cX)}(q_1))\uts \rB_\G$, and
we have (essentially the same argument as the proof of
\ref{fact:plagAct1}) an isomorphism
\begin{equation}\label{eq:plagAct3}
\cY\ts_{\cX}\cY \xrightarrow{\sim} \cY\uts q^{-1}(*)
\end{equation}
whose existence is basically the quasi-Galois
condition. Nevertheless, this does not, in
general, translate into the action of a 2-group.
In fact, quasi-Galois doesn't even imply $\G$
abelian, so there may be no 2-group structures on
$q^{-1}(*)$, and indeed
\begin{fact}\label{fact:plagAct2}
Notation as above for $q:\cY\ra\cX$ quasi-Galois
and quasi-minimal with trivial stabiliser then
there is a 2-commutative fibre square
\begin{equation}\label{eq:plagAct4}
 \xy
 (0,0)*+{q^{-1}(*)\uts \cY}="A";
 (38,0)*+{\cY}="B";
 (38,-18)*+{\cX}="C";
 (0,-18)*+{\cY }="D";
  (23,-5)*+{}="E";
  (15,-13)*+{}="F";
    {\ar_{\mathrm{trivial}} "A";"D"};
    {\ar^{\quad\quad f_\om} "A";"B"};
    {\ar^{q_1} "B";"C"};
    {\ar_{q_1} "D";"C"};
    {\ar@{=>}_{i_\om} "F";"E"};
 \endxy
\end{equation}
lifting (\ref{eq:plagAct3}) to the action of a
$2$-group $\gP\uts q\ra q$, cf. \ref{fact:act1}, 
iff $q$ is $2$-Galois.
\end{fact}
\begin{proof} If there is such a 2-group, then we
must have maps
\begin{equation}\label{eq:plagAct5}
\pi_1(\gP) \ra \pi_1(\fAut_{\et_2(\cX)}(q)) 
\ra  \pi_1(\fAut_{\et_2(\cX)}(q_1))
\end{equation}
whose composition is an isomorphism, {\it i.e.}
by definition,
\ref{def:noplag2}, $q$ is 2-Galois. 
Conversely,
by hypothesis there exist $1$-cells
$F_\om:q\ra q$ in $\et_2(\cX)$ for each $\om\in\pi_1(\fAut(q_1))$
along with $2$-cells $\a_{\tau,\om}:F_{\tau}F_\om\Rightarrow F_{(\tau\om)}$.
On the other hand $\fAut(q)$ is a strictly associative
2-category so $D(\a)$ as defined in (\ref{eq:act6})
is an element of the stabiliser in $\fAut(q)$ of
$F_{(\s\tau\om)}$ for any $\s,\tau,\om\in \pi_1(\fAut(q_1))$.
The stabiliser of any (weak) automorphism is naturally
isomorphic to $\pi_2(\fAut(q))$ so (following the
convention of \ref{factdef:act1}) $K_3(\a):=-D(\a)$
defines a 3 co-cycle for the action, 
\ref{fact:noplagAct},
of $\pi_1(\fAut(q_1))$
on $\pi_2(\fAut(q))$. By definition, (\ref{eq:group2}) \&
\ref{def:groupNick}, if, therefore, we give $q^{-1}(*)$
the $2$-group structure defined by the above
$K_3(\a)$ then, 
{\it cf.} \ref{factdef:act1},
(\ref{eq:plagAct4})  together with the $2$-cells
$\a_{\tau,\om}$ is a 2-group action.
\end{proof}
At which point its time for another
\begin{warning}\label{warn:plagAct1}
There are many non-equivalent 2-types,
let alone 2-group actions 
satisfying \ref{fact:plagAct2}, {\it i.e.} a
2-Galois cell, $q$ 
need not have a well defined Postnikov class in
\begin{equation}\label{eq:plagAct55}
\rH^3(\pi_1(\fAut(q_1)), \pi_2(\fAut(q)))
\end{equation}
The effect of equivalences in $\et_2(\cX)$ is
covered by \ref{factdef:act1}, but the argument
for dealing with a different $\a'_{\tau,\om}$
to the 2-cells $\a_{\tau,\om}$ encountered in the 
proof of \ref{fact:plagAct2} is a little different,
to wit: $\a'_{\tau,\om}\a^{-1}_{\tau,\om}$ stabilises
the automorphic 1-cell $F_{(\tau\om)}$, whence it
defines an element of $\pi_2(\fAut(q))$, so that
the difference between $K_3(\a')$ and $K_3(\a)$
is a $\pi_1(\fAut(q_1))$ co-boundary with values
in $\pi_2(\fAut(q))$. 
In neither case therefore can such effects change
the class in \eqref{eq:plagAct55} afforded by
\ref{fact:plagAct2} via
a choice of section of the rightmost arrow in
(\ref{eq:noplag2}), but rather to the
fact
that without
a universal cover there is a kernel in 
(\ref{eq:noplag2}), 
so that, even modulo equivalnce, the functors
$F_{\om}$ appearing in \ref{fact:plagAct2}
are very far from unique. As such,
since equivalence is
not the cause of the problem, it can be conveniently
summarised modulo the same, {\it i.e.} in co-homology.
To this end, pro tempore, denote the groups in
(\ref{eq:plagAct55}) by $\pi_1$, respectively $\pi_2$,
and observe that the Leray spectral sequence for
$p:\cY\ra\cY_1$ and the definition of minimal
yield a (canonical) element
\begin{equation}\label{eq:plagAct6}
\Hom_{\mathrm{grps}}(\pi_2,\pi_2)\ni \mathbf{1}\mpo  
\mathrm{Leray}_2^{0,1}(\mathbf{1})
\in \rH^2(\cY_1, \pi_2)
\end{equation} 
Now if in addition $q$ is quasi-Galois then
we have the action \ref{fact:noplagAct} and
$\pi_2$ defines a locally constant sheaf on
$\cX$ (albeit, strictly speaking, at this 
stage one has to either prove this or suppose
that $\et_1(\cX)$ with the fibre functor
satisfies the axioms \cite[V.4]{sga1}) and we
have the H\"oschild-Serre spectral sequence
\begin{equation}\label{eq:plagAct7}
\mathrm{hs}:\quad
\rH^i(\pi_1, \rH^j(\cY_1,\pi_2))\Rightarrow \rH^{i+j}(\cX,\pi_2)
\end{equation}
In particular, the obstruction to the 
vanishing of the composite
\begin{equation}\label{eq:plagAct8}
\mathrm{hs}_2^{0,2}(\mathrm{Leray}_2^{0,1}(\mathbf{1}))
\in \rE^{2,1}_{2, \mathrm{hs}}
\end{equation}
is exactly the obstruction to being 2-Galois.
Consequently, 
if $q$ is 2-Galois
we get something that looks like the Postnikov
class, {\it cf.} \ref{fact:pos2},  
\begin{equation}\label{eq:plagAct9}
\mathrm{hs}_3^{0,2}(\mathrm{Leray}_2^{0,1}(\mathbf{1}))
\in \rE^{3,0}_{3, \mathrm{hs}}
\end{equation} 
Nevertheless, since $\cY_1$ is very likely not
simply connected the $\rE^{3,0}$ will, in general,
change between the 2nd and 3rd sheets, {\it i.e.}
we have an exact sequence
\begin{equation}\label{eq:plagAct10}
\rH^1(\pi_1, \rH^1(\cY_1, \pi_2))
\xrightarrow{ \mathrm{hs}_2^{1,1}}
\rH^3(\pi_1,\pi_2) \ra 
\rE^{3,0}_{3, \mathrm{hs}}\ra 0
\end{equation}
in which one recognises the left hand group as
the isomorphism classes of splittings of the
(short exact because $q$ is 2-Galois)
sequence (\ref{eq:noplag2}), and the translates
under the image of this group in 
(\ref{eq:plagAct55}) of the class $K_3(\a)$
occurring in the proof of 
\ref{fact:plagAct2} describes the possible
2-types of 2-groups which can (and do)
satisfy {\it op. cit.}
\end{warning}
On the other hand the above difference, (\ref{eq:plagAct10})
between any two possible Postnikov classes dies
after a representable \'etale cover, so to get
something well defined is only a question of being
careful about limits. To this end observe
\begin{rmk}\label{rmk:plagAct1}
For any representable 0-cell, $q_1$, in $\et_2(\cX)$,
\begin{equation}\label{eq:plagAct11}
\pi_1 (\fAut_{\et_2(\cX)}(q_1)) \xrightarrow{\sim} 
\Aut_{\et_1(\cX)}(q_1),\,\,\mathrm{and}\,\,\,
\pi_2 (\fAut_{\et_2(\cX)}(q_1))=0
\end{equation}
while for any map $G_1:r_1\ra q_1$ of representable
Galois 0-cells, given an automorphism $H_1:r_1\ra r_1$
the $F_1$ affording the fibre square (\ref{eq:noplag4})
is unique up to equivalence and defines
\begin{equation}\label{eq:plagAct12}
(G_1)_*: \pi_1 (\fAut_{\et_2(\cX)}(r_1))
\ra 
\pi_1 (\fAut_{\et_2(\cX)}(q_1))
\end{equation}
so, slightly more intrinsically, for
$I_{*}$ the poset of \ref{sum:plag1}.(a)
\begin{equation}\label{eq:plagAct13}
\pi_1(\cX_*) := \varprojlim_{i} \pi_1 (\fAut_{\et_2(\cX)}(q_i)) 
\end{equation}
or, perhaps better, idem but for $I_{*}'$ 
rather than $I_{*}$ where the former is as per
\ref{fact:funct103} but with the 1-cells
$q_i$ Galois, {\it i.e.} avoid the 
specific choice of a point in the fibre
of \ref{sum:plag1}.(a). Nevertheless,
although $I_{*}'$ may appear more functorial
than $I_{*}$ 
of \ref{sum:plag1}.(a) they
are actually mutually co-final, since,
as the notation suggests,
by (\ref{eq:plagAct12}) the dependence
on the point whether in 
(\ref{eq:plagAct13}) or
\ref{sum:plag1}.(c) comes from the 
existence of the partially ordered
set $I_{*}$ itself, and not from the 
transition maps. Of course, one could
aim to take a limit over something
more general than a po-set to avoid such
dependence, but that would create a
swings and round-abouts issue since
a more ``general limit''
would fail to be uniquely unique.
\end{rmk}
Next let us consider a 2-Galois cell $q$ factorised
as $q_1p$ into a representable Galois cell and a
locally constant gerbe in $\rB_\G$'s. Independently
of any choice of points, we have, 
\ref{fact:noplagAct}, a canonical action of
$\pi_1(\fAut_{\et_2(\cX)}(q_1))$ on $\pi_2(\fAut_{\et_2(\cX)}(q))$,
which, inter alia, does not depend on the equivalence
class of the 1-cell $p$.
As such, consider a, not necessarily fibred, 
2-commutative square
\begin{equation}\label{eq:plagAct14}
\begin{CD}
q@<<{G}< r \\
@V{p}VV @VV{p'}V\\
q_1@<<{G_1}< r_1
\end{CD}
\end{equation}
where the 1-cells in the top row are 2-Galois,
the verticals locally constant gerbes in $\rB_\G$'s
respectively $\rB_{\G'}$'s, and, whence, the
lower cells are representable Galois. Now quite
generally- {\it i.e.} any 2-category- we have maps
\begin{equation}\label{eq:plagAct15} 
\begin{CD}
\pi_2(\fAut_{\et_2(\cX)}(r))@>{G_*}>> \mathrm{stab}(G)
@<{G^*}<< \pi_2(\fAut_{\et_2(\cX)}(q)) 
\end{CD}
\end{equation}
where $G^*$ is injective if $G$ is essentially
surjective on objects, and, if moreover $G$ is
a surjection on $\Hom$-sets then $G^*$ is surjective,
so that by \ref{fact:noplag1}
\begin{fact}\label{fact:plagAct3}
Under the hypothesis of (\ref{eq:plagAct14}),
$G^*$ of
(\ref{eq:plagAct15}) is an isomorphism affording,
\begin{equation}\label{eq:plagAct16}
G_2=(G^*)^{-1}G_*: \pi_2(\fAut_{\et_2(\cX)}(r))
\ra \pi_2(\fAut_{\et_2(\cX)}(q))
\end{equation}
\end{fact}
Better still,
\begin{fact}\label{fact:plagAct4}
$G_2$ of (\ref{eq:plagAct16})  
does not depend on the equivalence 
class of the 1-cell $G$ in
$\et_2(\cX)$ and yields 
a map of Galois modules,
{\it i.e.} for $G_1=(G_1)_*$ of (\ref{eq:plagAct12})
\begin{equation}\label{eq:plagAct17}
G_2(S^\om) = G_2(S)^{G_1(\om)},
\quad S\in \pi_2(\fAut_{\et_2(\cX)}(r)),\,
\om\in \pi_1(\fAut_{\et_2(\cX)}(r))
\end{equation}
\end{fact}
\begin{proof} A diagram chase (using \eqref{eq:plagAct15} b.t.w.)
reveals that $G_*(S^\om)$ is $G^*(S)^{G_1(\om)}$ up to
conjugation by an element of $\mathrm{stab}(G)$, which,
as we've observed is abelian because $G^*$ is an isomorphism,
whence the identity (\ref{eq:plagAct15}), and 
(same argument with $\om=\mathbf{1}$)
independence of $G_2$ from equivalences in $\et_2(\cX)$. 
\end{proof}
Putting all of which together, we obtain
\begin{factdef}\label{factdef:plagAct1}
For each isomorphism class of $2$-Galois cells,
choose one, and a point in the fibre of the 
same, so as to get, {\it cf.} \ref{sum:plag1}.(a),
a set, $I_*$, of triples $(q_i, *_i, \phi_i)$
(alternatively, \ref{rmk:plagAct1}, take the
moduli of the fibre $q_i^{-1}(*)$ so as to make
less choices) which we partially order according
to (\ref{eq:funct1406}). The resulting factorisation
of the $q_i$ into locally constant gerbes and representable
maps implies that the above directed 
({\it i.e.} partially ordered and right co-filtering)
set, $I_*$,
is mutually co-final with that of \ref{sum:plag1}.(a),
while by \ref{fact:plagAct4},
\begin{equation}\label{eq:plagAct18}
\pi_2(\cX_*):=\varprojlim_{i} \pi_2 (\fAut_{\et_2(\cX)}(q_i)) 
\end{equation}
is a continuous $\pi_1(\cX_*)$-module. Equivalently,
if we (uniquely) identify $ \pi_2 (\fAut_{\et_2(\cX)}(q_i))$
with a sub-group $\G_i$ of the stabiliser of
$*_i$ via the isomorphism \ref{fact:noplag2}, then
whenever $i>j$ there is a unique map $\G_i\ra \G_j$
afforded by (\ref{eq:funct1406}),
but independent of $\z_{ji}$ in op. cit.,
 and
\begin{equation}\label{eq:plagAct18plagAct18}
\pi_2(\cX_*):=\varprojlim_{i} \G_i
\end{equation}
with continuous-$\pi_1(\cX_*)$ action given by (\ref{eq:noplagAct1}).
\end{factdef}
\subsection{The Postnikov class}\label{SS:III.5}

We turn to defining the Postnikov
class, and, in particular,  the problem, \ref{warn:plagAct1},
posed by the lack of the uniqueness in
\ref{fact:plagAct2}. To this end let us make
\begin{defn}\label{def:PosSeq}
Let $I$ be a 
(not necessarily co-final nor right co-filtering nor full)
sub-category of 
the directed set, $I_*$, of
\ref{factdef:plagAct1}, viewed as a category, 
then by a {\it Postnikov sequence},
or $I$-Postnikov sequence if there is a danger
of confusion, 
is to be understood the following data

(a) For each $\om\in\pi_1(\cX_*)$ and each $i\in I$ a
weak automorphism $\om^i:q_i\ra q_i$ in $\et_2(\cX)$ of
the 2-Galois cell $q_i$ whose action on the moduli of
the fibre, $\vert q_i^{-1}(*)\vert$, is that of $\om$
implied by \eqref{eq:plag5} and \ref{def:funct103},
with the further proviso that $\om^i$ is trivial
should this action be trivial.

(b) For each $i\geq j\in I$, a (unique) $1$-cell $F_{ji}:=(f_{ji},\xi_{ji}):q_i\ra q_j$
in $\et_2(\cX)$ 
satisfying the base point condition \eqref{eq:funct1406},
along with 2-cells 
$\phi^{ji}_\om: F_{ji} \om^i\Rightarrow \om^j F_{ji}$
in $\et_2(\cX)$, $\om\in\pi_1(\cX_*)$, with all of
$F_{ji}$, $\phi^{ji}_\om$ identities if $i=j$.

(c) For all $i\in I$ and $\tau,\om\in\pi_1(\cX_*)$ 
2-cells $\a^{i}_{\tau,\om}: \tau^i\om^i\Rightarrow (\tau\om)^i$
such that whenever $i\geq j\in I$ the following diagram commutes
\begin{equation}\label{349}
 \xy 
 (0,0)*+{F_{ji}\tau^i\om^i}="A";
 (30,0)*+{F_{ji}(\tau\om)^i}="B";
 (30,-18)*+{(\tau\om)^jF_{ji}}="C";
 (0,-18)*+{\tau^j\om^jF_{ji}}="D";
    {\ar@{=>}_{(F_{ji})_*\a^i_{\tau\om}} "A";"B"};
    {\ar@{=>}_{(\tau^j)_*(\phi^{ji}_\om)(\om^i)^*(\phi^{ji}_\tau)} "A";"D"};
    {\ar@{=>}^{\phi^{ji}_{\tau\om}} "B";"C"};
    {\ar@{=>}_{ (F_{ji})^*\a^j_{\tau\om}} "D";"C"};
 \endxy 
\end{equation} 

(d) For every $i\geq j\geq k\in I$ there is a 2-cell
$\g_{kji}:F_{ki}\Rightarrow F_{kj}F_{ji}$ such that
for all $\om\in\pi_1(\cX_*)$ the following diagram 2-commutes
\begin{equation}\label{PosSeq1}
 \xy
 (38,0)*+{q_j}="B";
 (0,-6)*+{q_k}="C";
(15,0)*+{q_i}="D";
 (38,-12)*+{q_j}="F";
 (0,-18)*+{q_k}="G";
(15,-12)*+{q_i}="H";
{\ar_{}_{\om^j} "B";"F"};
    {\ar_{\om^k} "C";"G"};
    {\ar_{\om^i} "D";"H"};
{\ar_{}^{F_{kj} } "F";"G"};
    {\ar@{-->}_{F_{ki} } "H";"G"};
    {\ar^{F_{ji} } "H";"F"};
{\ar_{}^{F_{kj}} "B";"C"};
    {\ar_{F_{ki}} "D";"C"};
    {\ar^{F_{ji}} "D";"B"};
\endxy
\end{equation}
wherein the implied rectangular cells are given
by (c) so, equivalently, the following commutes
\begin{equation}\label{PosSeq2}
 \xy 
 (0,0)*+{F_{ki}\om^i}="A";
 (20,0)*+{\om^k F_{ki}}="E";
 (40,0)*+{\om^k F_{kj}F_{ji} }="B";
 (40,-18)*+{F_{kj} \om^j F_{ji}}="C";
 (0,-18)*+{F_{kj}F_{ji}\om^i }="D";
    {\ar@{=>}_{\phi^{ki}_{\om}} "A";"E"};
    {\ar@{=>}_{\om^k_*\g_{kji}} "E";"B"};
    {\ar@{=>}^{F^*_{ji}\phi^{kj}_{\om}} "B";"C"};
    {\ar@{=>}_{ (F_{kj})_*\phi^{ji}_{\om}} "D";"C"};
    {\ar@{=>}_{ (\om_{i})^*\g_{kji}} "A";"D"};
 \endxy 
\end{equation} 

(e) For $i\geq j\geq k\geq l\in I$ the following 
diagram 2 commutes
\begin{equation}\label{PosSeq3}
 \xy
 (20,-24)*+{q_i}="A";
 (38,6)*+{q_j}="B";
 (0,0)*+{q_k}="C";
(15,12)*+{q_l}="D";
{\ar_{}_{F_{ji} } "A";"B"};
    {\ar^{F_{kj} } "B";"C"};
    {\ar^{F_{ki} } "A";"C"};
    {\ar^{F_{lk} } "C";"D"};
    {\ar_{F_{lj} } "B";"D"};
{\ar@{-->}_{F_{li}} "A";"D"}
\endxy
\end{equation}
or, again, since the triangles are given by (d)
the following commutes
\begin{equation}\label{PosSeq4}
 \xy 
 (0,0)*+{F_{lj}F_{ji}}="A";
 (30,0)*+{F_{li} }="B";
 (30,-18)*+{F_{lk}F_{ki}}="C";
 (0,-18)*+{F_{lk}F_{kj}F_{ji} }="D";
    {\ar@{=>}^{\g_{lji}} "B";"A"};
    {\ar@{=>}^{\g_{lki}} "B";"C"};
    {\ar@{=>}_{ (F_{lk})_*\g_{kji}} "C";"D"};
    {\ar@{=>}_{ (F_{ji})^*\g_{lkj}} "A";"D"};
 \endxy 
\end{equation} 
\end{defn}
In order to see where this leads, some commentary is in order 
by way of
\begin{rmk}\label{5.2}
Exactly as in \ref{factdef:act1}, 
should we have a Postnikov sequence in the
sense of \ref{def:PosSeq}
with $I$ co-final, then
for each $i\in I_*$ there is a 3 co-cycle,
\begin{equation}\label{350}
K_3^i= -D(\a^i):\pi_1((\fAut(q_i))^3\ra \pi_2(\fAut(q_i))
\end{equation}
or, slightly more correctly, if 
$D$ is defined as in \eqref{eq:act6} then we
should really 
employ the isomorphism \eqref{eq:plagAct15} in order to
write $-D(\a^i)=(\s\tau\om)^i_* (K_3^i)$. However,
in the presence of the compatability condition \eqref{349},
$(F_{ji})_* D(\a^i)
\phi_{(\s\tau\om)}^{ji}
=\phi_{(\s\tau\om)}^{ji} (F_{ji})^* D(\a^j)$-
{\it cf.}, modulo the 
change of notation, the proof of \ref{factdef:act1} post
\eqref{eq:act9}-whence,
under the unique map $\G_i\ra\G_j$
of \ref{factdef:plagAct1}, $K_3^i\mpo K_3^j$. Consequently,
a Postnikov sequence affords a continuous 3 co-cycle,
\begin{equation}\label{351}
K_3: \pi_1(\cX_*)^3\ra\pi_2(\cX_*)
\end{equation}
which, as the notation suggests, will lead to the definition of
the pro-finite Postnikov class, 
so what we need to do is prove 
(appropriate) existence and 
uniqueness of Postnikov sequences.
\end{rmk}
This much only requires \ref{def:PosSeq}(a)-(c), and
where we need (d)-(e) is to get good
behaviour of Postnikov sequences under base change,
which, if nothing else, is a strain on notation,
so we spell out in detail our
\begin{setup}\label{5.3} Let $i,j,k,l\in I_*$, not
just $I$, be given,
{\it i.e.} our cells are fairly arbitrary,
with $i>j$, $k>j$, and
suppose

(a) For each of $i,j,k$ but not $l$, the condition
\ref{def:PosSeq}.(a) is given.

(b) For the pairs, 
$ji$, and $jk$ the conditions \ref{def:PosSeq}.(b)-(c)
are given.

(c) There is given in $\et_2(\cX)$ a 2-commutative fibre
square in which the compatability condition,
\eqref{eq:funct1406}, on base points is supposed
\begin{equation}\label{352}
 \xy 
 (0,0)*+{q_l}="A";
 (30,0)*+{q_i}="B";
 (30,-18)*+{q_j}="C";
 (0,-18)*+{q_k}="D";
 (22, -6)*+{}="E";
 (8,-12)*+{}="F";
    {\ar_{F_{il}} "A";"B"};
    {\ar_{F_{kl}} "A";"D"};
    {\ar^{F_{ji}} "B";"C"};
    {\ar_{ F_{jk}} "D";"C"};
 {\ar@{=>}_{ \b} "E";"F"}; 
 \endxy 
\end{equation}
\end{setup}
With such a set up, 
we have, for each $\om\in\pi_1(\cX_*)$, 
a diagram of 2-cells,
\begin{equation}\label{353}
F_{ji}\om^i F_{il} 
{\build\Longrightarrow_{}^{F^*_{il}\phi^{ji}_\om}}
\om^j F_{ji}F_{il}
{\build\Longrightarrow_{}^{(\om^j)_*\b}}
\om^j F_{jk} F_{kl}
{\build\Longleftarrow_{}^{F^*_{kl}\phi^{jk}_\om}}
F_{jk}\om^k F_{kl}
\end{equation}
so by the universal property, \eqref{FibreA3},
of fibre squares: \eqref{353}
affords a 2-commutative diagram
\begin{equation}\label{354}
 \xy 
 (0,0)*+{q_l}="A";
 (30,0)*+{q_i}="B";
 (30,-18)*+{q_j}="C";
 (0,-18)*+{q_k}="D";
 (22, -6)*+{}="E";
 (8,-12)*+{}="F";
 (-20,15)*+{q_l}="H";
 (-2,4)*+{}="a";
 (2,8)*+{}="b";
 (-2,-2)*+{}="c";
 (-6,-6)*+{}="d";
    {\ar_{F_{il}} "A";"B"};
    {\ar^{F_{kl}} "A";"D"};
    {\ar^{F_{ji}} "B";"C"};
    {\ar_{ F_{jk}} "D";"C"};
 {\ar@{=>}_{ \b} "E";"F"}; 
 {\ar^{\om^i F_{il}} "H";"B"};
 {\ar_{\om^k F_{kl}} "H";"D"};
  {\ar^{\om^l} "H";"A"};
 {\ar@{=>}_{ \phi^{il}_\om} "a";"b"}; 
 {\ar@{=>}_{ \phi^{kl}_\om} "c";"d"}; 
\endxy 
\end{equation}
Now, in principle, there are many such diagrams,
equivalently, many ways to choose $\om^l$,
$\phi^{il}_\om$, and $\phi^{kl}_\om$. Nevertheless,
on choosing one
\begin{fact}\label{5.4}
For every $\tau,\om\in\pi_1(\cX_*)$ there are
unique 2-cells, $\a^{l}_{\tau,\om}: \tau^l\om^l\Rightarrow (\tau\om)^l$ 
such that both the diagrams obtained on replacing
$(i,j)$ by $(l,i)$, respectively $(l,k)$ in \eqref{349}
commute.
\end{fact}
\begin{proof} Get a large piece of paper, and
apply the uniqueness of 2 cells in
\eqref{FibreA3}-\eqref{FibreA5}.
\end{proof}

We also need a related unicity statement for the
diagrams \eqref{353} and \eqref{354} themselves, 
which we
\begin{setup}\label{5.5} as follows: five 2-Galois
cells $q_i$, $q_j$, $q_k$, $q_0$, $q'$ are given in $I_*$, with all
subsequent 1-cells respecting the base point condition
\eqref{eq:funct1406}. Moreover,

(a) The condition \ref{def:PosSeq}.(a) is given for all
five cells $q_i$, $q_j$, $q_k$, $q_0$, $q'$.

(b) The condition \ref{def:PosSeq}.(b) is given only for the
pairs $(i,j)$, $(i,0)$, $(j,0)$, $(k,0)$, and $(',0)$.

(c) For $\iota\in\{i,j,k \}$, there is given  a 2-commutative
fibre square
\begin{equation}\label{355}
 \xy 
 (0,0)*+{q'_{\iota}}="A";
 (30,0)*+{q'}="B";
 (30,-18)*+{q_0}="C";
 (0,-18)*+{q_i}="D";
 (22, -6)*+{}="E";
 (8,-12)*+{}="F";
    {\ar_{F'_{0\iota}} "A";"B"};
    {\ar_{p_\iota} "A";"D"};
    {\ar^{p} "B";"C"};
    {\ar^{F_{0\iota}} "D";"C"};
 {\ar@{=>}_{ \b_{\iota}} "E";"F"}; 
 \endxy 
\end{equation}

(d) The condition \ref{def:PosSeq}.(d) is given for
all possible triples in $\{i,j,k,0\}$.

(e) The condition \ref{def:PosSeq}.(e) is given,
for the quadruple $i\geq j\geq k\geq 0$.
\end{setup}
Now, in the presence of \ref{5.5} we can,
modulo notation ({\it i.e.} $q_k=q_i$, $q'=q_i$, 
$\b=\b_{ji}:p_jF'_{ji}\Rightarrow F_{ji}p_i$) suppose
that the square \eqref{352} fits 
with the squares \eqref{355}
into a 2-commutative
diagram
\begin{equation}\label{356}
 \xy
 (38,0)*+{q'_j}="B";
 (0,-6)*+{q'}="C";
(15,0)*+{q'_i}="D";
 (38,-12)*+{q_j}="F";
 (0,-18)*+{q_0}="G";
(15,-12)*+{q_i}="H";
{\ar_{}_{p_j} "B";"F"};
    {\ar_{p} "C";"G"};
    {\ar_{p_i} "D";"H"};
{\ar_{}^{F_{0j} } "F";"G"};
    {\ar@{-->}_{F_{0i} } "H";"G"};
    {\ar^{F_{ji} } "H";"F"};
{\ar_{}^{F'_{0j}} "B";"C"};
    {\ar_{F'_{0i}} "D";"C"};
    {\ar^{F'_{ji}} "D";"B"};
\endxy
\end{equation}
where in the bottom face
is $\g_{0ji}$,
so naturally we similarly denote the top most 2-cell in \eqref{356}
by $\g'_{0ji}:r_i\Rightarrow r_j F_{ji}$, and
\begin{fact}\label{fact:PosSeq1}
There is a unique natural transformation $\g'_{kji}:F'_{ki}
\Rightarrow F'_{kj}F'_{ji}$ such that both of the following
diagrams commute
\begin{equation}\label{PosSeq5}
 \xy 
 (0,0)*+{F'_{0k}F'_{ki}}="A";
 (30,0)*+{F'_{0i}}="B";
 (30,-18)*+{F'_{0j}F'_{ji}}="C";
 (0,-18)*+{F'_{0k}F'_{kj} F'_{ji}}="D";
    {\ar@{=>}^{\g'_{0ki}} "B";"A"};
    {\ar@{=>}^{\g'_{0ji}} "B";"C"};
    {\ar@{=>}_{(F'_{ji})^* \g'_{0kj}} "C";"D"};
    {\ar@{=>}_{(F'_{0k})_*\g'_{kji} } "A";"D"};
 \endxy 
\end{equation} 
\begin{equation}\label{PosSeq6}
 \xy 
 (0,0)*+{p_k F'_{ki}}="A";
 (40,0)*+{F_{kj}F_{ji} p_i}="B";
 (40,-18)*+{F_{kj} p_j  F'_{ji}}="C";
 (0,-18)*+{p_k F'_{kj} F'_{ji}}="D";
 (20,0)*+{F_{ki} p_i}="E";
    {\ar@{=>}_{(p_k)_*\g'_{kji} } "A";"D"};
    {\ar@{=>}_{(F_{kj})_* \b_{ji}} "C";"B"};
    {\ar@{=>}_{\b_{ki}} "A";"E"};
    {\ar@{=>}_{p_{i}^*\g_{kji}} "E";"B"};
    {\ar@{=>}^{(F'_{ji})^* \b_{kj}} "D";"C"};
 \endxy 
\end{equation} 
\end{fact}
\begin{proof} Apply the definition of fibre
products, \ref{def:Fibre2}, and \ref{5.5}.(e).
\end{proof}
Now for $\iota\in\{i,j, k\}$, form the diagram \eqref{354},
but with the square of {\it op. cit.} replaced by the
square(s) \eqref{355}, and change the notation as follows
\begin{newnot}\label{5.7} 
Denote the lower, and upper, triangular 
cells by 
$\theta^{\iota}_\om$, respectively $\psi^{\iota}_\om$,
and $G^{\iota}_\om$ the lifting of $\om\in\pi_1(\cX_*)$ to
a weak automorphism of $q'_{\iota}$, with $\om':q'\ra q'$.
\end{newnot}
As such another exercise in \ref{def:Fibre2},
where now the salient condition is \ref{5.5}.(d),
implies
\begin{fact}\label{5.8} Notation as in \ref{5.7}, then 
for each $\om\in\pi_1(\cX_*)$ there is a unique 2 cell
in $\et_2(\cX)$,
$\Phi^{ji}_\om: F'_{ji} G^i_\om\Rightarrow G^j_\om F'_{ji}$,
such that both of the following diagrams commute
\begin{equation}\label{357}
 \xy 
 (0,0)*+{F'_{0j}F'_{ji}G_{\om}^i}="A";
 (60,0)*+{\om' F'_{0i}}="B";
 (60,-18)*+{ \om' F'_{0i}}="C";
 (0,-18)*+{F'_{0j}G^j_\om F'_{ji}}="D";
 (30,0)*+{F'_{0i}G^i_\om}="E";
 (30,-18)*+{\om'F'_{0j}F'_{ji}}="F";
    {\ar@{=>}_{(F'_{0j})_*\Phi^{ji}_\om } "A";"D"};
    {\ar@{=}_{} "B";"C"};
    {\ar@{=>}^{(G^i_\om)^*\g'_{0ji}} "E";"A"};
    {\ar@{=>}_{\psi^i_\om} "E";"B"};
    {\ar@{=>}_{(\om')_* \g'_{0ji}} "C";"F"};
    {\ar@{=>}^{ (F'_{ji})^*\psi^j_\om} "D";"F"};
 \endxy 
\end{equation} 
\begin{equation}\label{358}
 \xy 
 (0,0)*+{p_jF'_{ji}G_{\om}^i}="A";
 (60,0)*+{F_{ji}\om^i p_i}="B";
 (60,-18)*+{\om^j F_{ji} p_i}="C";
 (0,-18)*+{p_jG^j_\om F'_{ji}}="D";
 (30,0)*+{F_{ji} p_i G^i_\om}="E";
 (30,-18)*+{\om^j p_j F'_{ji}}="F";
    {\ar@{=>}_{(p_j)_*\Phi^{ji}_\om } "A";"D"};
    {\ar@{=>}^{(p_i)^*\phi^{ji}_\om} "B";"C"};
    {\ar@{=>}_{(G^i_\om)^*\b_{ji}} "A";"E"};
    {\ar@{=>}_{(F_{ji})_*\theta^i_\om} "E";"B"};
    {\ar@{=>}^{(\om^j)_* \b_{ji}} "F";"C"};
    {\ar@{=>}^{ (F'_{ji})^*\theta^j_\om} "D";"F"};
 \endxy 
\end{equation} 
\end{fact}
To which it's worth adding a clarifying
\begin{rmk}\label{rmk5.8}
In the particular case that $j=0$ we can, and do, suppose
that $q'_0=q'$ exactly, so that everything collapses, and
whence $\psi^i_\om= \Phi^{0i}_\om$ by \eqref{357}.
\end{rmk} 
Which is not without relevance to the statement of
\begin{fact}\label{5.9} 
Let a Postnikov sequence, $I$, with final object $0$ be given
along with a $1$-cell $p:q'\ra q_0$ satisfying \ref{def:PosSeq}.(a)-(c);
form the fibre squares \eqref{355}, and suppose that
each $q'_i$ therein is quasi-minimal
with trivialisable stabiliser, then the following
base change data defines an $I'$ Postnikov sequence for
some appropriate subset of $I_*$ in $1$-to-$1$ correspondence
with $I$,

(a) The weak automorphisms of \ref{def:PosSeq}.(a) for the
new $q'_i$ 
are defined as in \ref{5.7}.

(b) For $i\geq j\in I$,
we take $F'_{ji}:q'_i\ra q'_j$
the 1-cell in $\et_2(\cX)$ defined by \ref{5.8},
with transition 2 cell the $\Phi^{ji}_\om$.

(c) The two cells $(\a')^i_{\tau,\om}$ of \ref{def:PosSeq}.(b)
are defined using
\ref{5.4} albeit with the square \eqref{352} replaced 
by \eqref{355}.

(d) The two cells $\g'_{kji}$ of \ref{def:PosSeq}.(d)
are defined in \ref{fact:PosSeq1}. 
\end{fact}
\begin{proof} By hypothesis each $q'_i$ is quasi-minimal
with trivialisable stabiliser,
so by (a) and (b) above, they're equally 2-Galois. 
It thus remains to check the various commutativity
conditions of \ref{def:PosSeq}. Each of which
is an exercise in \ref{def:Fibre2}, and 
the salient points are

(c)$'$ The resulting \eqref{349} commutes by \ref{5.4}
and \ref{5.8}.

(d)$'$ The commutativity of the new \eqref{PosSeq2} 
follows from \ref{fact:PosSeq1} and \ref{5.8}.

(e)$'$ The tetrahedron condition \eqref{PosSeq1} 
follows from \ref{fact:PosSeq1}.
\end{proof}

Before proceeding to a conclusion an observation is
in order, {\it i.e.}

\begin{rmk}\label{5.11} If we take $I_1$ to be the
sub-category of $I_*$ to be the representable Galois
cells, then it forms a Postnikov sequence. Indeed,
by \ref{sum:plag1} if the $F_{ji}$ are any 1-cells
respecting the base point condition, \eqref{eq:funct1406},
then  by \ref{fact:funct102}, there certainly exist
$\g_{kji}:F_{ki}\Rightarrow F_{kj}F_{ji}$, and everything
else, {\it i.e.} the commutativity of \eqref{349},
\eqref{PosSeq2}, and \eqref{PosSeq4}, is trivial by
\eqref{eq:plagAct15} since representable cells have
no stabiliser in $\et_2(\cX)$. Indeed, had we wished
to, we could even suppose in $I_1$ that the diagrams
\eqref{PosSeq1} and \eqref{PosSeq3} were strictly,
and not just 2-commutative. Similarly, the $\om^l$,
$\phi^{il}_\om$, and $\phi^{kl}_\om$ in \eqref{354}
are to all intents and purposes unambiguous if $F_{jk}$
is in $I_1$. More precisely, as soon as one has the
2 triangular cells in {\it op. cit.} then the
universal property determined by \eqref{353} is
automatically satisfied since $q_j$ has no stabiliser.
\end{rmk}
These consideration lead us to make one final
\begin{defn}\label{5.12} A Postnikov sequence $I$ is
stable under representable base change if,

(a) For $q_i$ a 2-Galois cell,  $i\in I$,
there is an arrow in $I$ corresponding to the map 
$q_i\ra\bar{q}_i$ of the factorisation of $q_i$ as
a locally constant gerbe over a representable Galois cell.

(b) For every (representable by definition) Galois cover 
$G_1:r_1\ra \bar{q}_i$, there is an arrow in $I$
between the class of $(G_1)^*q_i$ and $q_i$.

(c) The natural transformation between the sides
of a fibre square implied by the above (a) and (b),
and \ref{def:PosSeq}.(d) has the universal property,
\eqref{FibreA3}-\eqref{FibreA5}, of a fibre square. 
\end{defn}
To apply all of this put an ordering, $I<I'$, on Postnikov
sequences by insisting in the first instance that we
have an inclusion of sub-categories of $I_*$, and in the
second instance that all of the extra data \ref{def:PosSeq}(a)-(d)
for $I'$ restricts identically ({\it i.e.} no extra equivalences
whatsoever) to that for $I$. As such, by choice and \ref{5.11}
a Postnikov sequence which is maximal amongst those
which are right co-filtering, stable under representable
base change, and contain
$I_1$ certainly exists. Better still,
\begin{fact}\label{5.13} If $I$ is such a Postnikov sequence,
then for every 2-Galois cell $q'$ in $\et_2(\cX)$, there is
an $i\in I$ such that the corresponding 2-Galois cell
$q_i$ admits a map $q_i\ra q'$ in $\et_2(\cX)$ respecting
the base point condition \eqref{eq:funct1406}. In particular,
as a partially ordered subset of $I_*$, $I$ is co-final,
and right co-filtering.
\end{fact}
\begin{proof} Suppose there is no such map. There are, however,
cells in $I$ under $q'$, {\it e.g.} $\bar{q}'$, for $q'\ra\bar{q}'$
it's factorisation as a locally constant gerbe over a representable
Galois cell, and
we first reduce to the case that there is a cell
under $q'$ with
some further pleasing properties
\begin{claim}\label{5.14} Without loss of generality 
$\exists 0\in I$ whose class in $I_*$ is less
than $q'$ such that,

(a) For every $i>0\in I$, the base change $F_{i0}^*q'$ is
quasi-minimal with trivialisable stabiliser.

(b) There is a map, as ever satisfying the base point
condition \eqref{eq:funct1406}, $p:q'\ra q_0$, together with a
choice of weak automorphisms, $\om':q'\ra q'$, and natural
transformations $\theta_\om: p\om'\Rightarrow \om^0 p$;
$\a'_{\tau,\om}:\tau'\om'\Rightarrow (\tau\om)'$ such 
that (for $F_{ji}=p$ {\it etc.}) \eqref{349} commutes.

\end{claim} 
\begin{proof}[Sub-proof] Consider, quite generally, maps of
2-Galois cells
\begin{equation}\label{369}
q_a\leftarrow q \rightarrow q_b
\end{equation}
satisfying the base point condition \eqref{eq:funct1406},
with $\G_a$, $\G_b$ the stabilisers of the base points of
$q_a$, $q_b$ so the stabiliser of the
base point, $*$, of their fibre product is $\G_a\ts\G_b$, 
and- proof of \ref{claim:funct102}- there's
a quasi minimal cell $r$ which is Galois over 
the connected component of $*$ such that  
\begin{equation}\label{370}
A:=\pi_2(\fAut_{\et_2(\cX)}(r))\hookrightarrow \G_a\ts \G_b
\end{equation}
On the other hand a representable Galois cover
of $q$ maps to $r$, so we must have a factorisation
\begin{equation}\label{371}
 \xy 
 (0,0)*+{\pi_2(\fAut_{\et_2(\cX)}(q))}="A";
 (22,0)*+{A}="B";
 (11,-10)*+{\G_a\ts \G_b}="C";
{\ar_{}^{} "A";"B"};
    {\ar^{} "B";"C"};
    {\ar_{} "A";"C"};
\endxy
\end{equation}
so, as far as the stabilisers of the base points
of 2-Galois cells are concerned, $A$, is the 
fibre product of $\G_a$ and $\G_b$. To identify
it more precisely, observe that by \ref{fact:noplag1},
we have a commutative diagram 
($K_a$, $K_b$, $\G''$ defined therein)
with exact rows and
columns,
\begin{equation}\label{372}
\begin{CD}
@. @. 0 @. 0 @. \\
@. @. @VVV @VVV @.\\
@. @.  K_a @= K_a @.\\
@. @. @VVV @VVV @. \\
0@>>> K_b @>>> A@>>> \G_b@>>> 0\\
@. @| @VVV @VVV @.\\
0@>>> K_b @>>> \G_a @>>> \G''@>>> 0 \\
@. @. @VVV @VVV @. \\
@. @. 0 @. 0 @.
\end{CD}
\end{equation}
so $A\xrightarrow{\sim} \G_a\ts_{\G''} \G_b$. To apply
this to the case in point let $\G_i$ be the stabilisers
of the base points of the $q_i$, $i\in I$, $\G'$ that
of $q'$, and forms much the same sort of diagram, {\it i.e.}
\begin{equation}\label{373}
\begin{CD}
@. 0 @. 0 @. 0 @. \\
@. @VVV @VVV @VVV @.\\
0@>>> K'\cap K @>>> K'@>>>   K'' @>>> 0\\
@. @VVV @VVV @VVV @. \\
0@>>> K @>>> \pi_2(\cX_*) @>>>\varprojlim_{i\in I} \G_i@>>> 0\\
@. @VVV @VVV @VVV @.\\
0@>>> \bar{K} @>>> \G' @>>> \G''@>>> 0 \\
@. @VVV @VVV @VVV @. \\
@. 0 @. 0 @. 0 @.
\end{CD}
\end{equation}
In which everything is continuous, so there is a $0\in I$
such that $\G_0 \twoheadrightarrow \G''$ under the natural
vertical in the rightmost column of \eqref{373},  
so by the universal property of $A$ in \eqref{371} the
definitions of $\G''$ in \eqref{372} and \eqref{373}
coincide for $q_a=q'$, $q_b=q_i$, 
$i\geq 0\in I$ in \eqref{369}. In particular, therefore,
by the stability of $I$ under representable base change,
we can suppose that there is a 2-Galois cell, $q$, mapping to
$q'$ and $q_0$ such that $q$ and $q_0$ are both locally
constant gerbes over the same Galois cell, $\bar{q}_0$, 
and the stabiliser
of the base point of $q$ is (canonically) isomorphic to
$\G_0\ts_{\G''} \G'$. Now for $i>0\in I$, the stabiliser
of the base point of $q_i\ts_{q_0} q$ is certainly trivialisable,
and canonically isomorphic to
\begin{equation}\label{374}
\G_i\ts_{\G_0} (\G_0\ts_{\G''} \G') \xrightarrow{\sim}
\G_i\ts_{\G''} \G'
\end{equation}
which is equally the stabiliser of the base point of the
smallest (pointed) 2-Galois cell mapping to $q_i\ts_{q_0} q$,
and since $q$, $q_0$ are both locally constant gerbes over
the same Galois cover, the said fibre product is connected.
Consequently, on replacing $q'$ by $q$,
\ref{5.14}.(a) holds. 

As to item (b), with $q'=q$, $q_0$ as above, we  certainly
have that the 1-cells $\om'$, and 2-cells $\a'_{\tau,\om}$
exist since $q'$ is 2-Galois. The $\theta_\om$ however may
not, since a priori it's only the case that $p\om'$ is
equivalent to some $\om p Z_\om$ for some 1-cell, $Z_\om$,
described by \ref{claim:noplag1}. On the other hand all the
$Z_\om$ become trivial after any sufficiently large Galois
covering $G:\bar{r}\ra \bar{q}_0=\bar{q'}$, so one forms
the 2-commutative diagram \eqref{356} albeit with $\bar{q}_0$,
$q_0$, $q'$ in the bottom triangle, and say $\bar{r}$, $r_0$,
$r'$ on the top. As such by \eqref{5.4} applied to the
unique fibre squares in such a diagram with $\bar{r}$ as
a vertex, there are liftings $G'_\om$, $G^0_\om$ of $\om'$,
respectively $\om^0$, along with the $\a'_{\tau,\om}$,
respectively $\a^0_{\tau,\om}$. Now the definition,
\eqref{354}, of these liftings is such that one
gets rather cheaply, {\it cf.} \ref{5.11}, that the
push forward of the $G'_\om$'s are 
equivalent to the pull-backs of
the $G^0_\om$. Similarly, again
by \ref{5.11},  the $G^0_\om$ may harmlessly be
confused with any existing liftings to $r_0$ defined
by the Postnikov sequence structure on $I$, 
so the $\theta_\om$ of \ref{5.14}.(b) exist. Now while
such a confusion with the pre-defined $\a_{\tau,\om}$
on $r_0$ isn't legitimate, it's nevertheless true
that the only remaining problem is that \eqref{349}
understood for $r'\ra r_0$ may fail to commute.
Should this occur, however, one just changes the
definition of the $\a'_{\tau,\om}$ on $r'$ by way
of \eqref{eq:plagAct15} and the surjectivity of
$\G'\ra \G_0$. As such, replacing $q_0$ by $r_0$,
and $q'$ by $r'$ we get item (b).
\end{proof}
Now we apply this to the proof of \ref{5.13} in the
obvious way, {\it i.e.} $I_0:=\{i\geq 0\}$ is naturally
a sub Postnikov sequence of $I$, and we can form the
base change $I'_0$ of \ref{5.9} along $q'\ra q_0$. As such
we require to splice $I'_0$ and $I$ together to a new
Postnikov sequence $I_\mathrm{big}$. Plainly the objects
of the latter are $I\cup I'_{0}\subseteq I_*$, which is,
in fact, a disjoint union by hypothesis, and we dispose
of the maps $p_i:q'_i\ra q_i$ of \eqref{355}. Confusing
cells with their classes in $I_\mathrm{big}$ we add 
new arrows beyond those already existing in $I$ or $I'_0$
according to the rule
\begin{equation}\label{375}
q'_i \geq q_j \, \text{iff}\, i \geq j \in I
\end{equation}
which not only defines a right co-filtered structure on $I_\mathrm{big}$
as a sub-category of $I_*$, but also satisfies conditions
\ref{5.12}.(a)-(b) for stability under representable base change.
The additional items implied by \eqref{375} are

(b) The 1-cell in question is $F_{ji}p_i$ and for $\om\in\pi_1(\cX_*)$
the transition 2-cells are, {\it cf.} \ref{5.7}, are the composition of
\begin{equation}\label{376}
F_{ji} p_i G^i_\om {\build\Longrightarrow_{}^{(F_{ji})_*\theta^i_\om}}
F_{ji} \om^i p_i  {\build\Longrightarrow_{}^{p_i^*\phi^{ji}_\om}}
\om^j F_{ji} p_i
\end{equation}
(d) The two cells, $\g_{abc}$, say, of \ref{def:PosSeq}.(d) are
according to the cases,
\begin{equation}\label{377}
\begin{split}
a=q'_i,\, b=q_j,\, c=q_k: & \quad
F_{ki} p_i  {\build\Longrightarrow_{}^{p_{i}^*\g_{kji}}}
F_{kj} F_{ji} p_i  \\
a=q'_i,\, b=q'_j,\, c=q_k: & \quad
F_{ki} p_i  {\build\Longrightarrow_{}^{p_{i}^*\g_{kji}}}
F_{kj} F_{ji} p_i {\build\Longleftarrow_{}^{(F_{kj})_*\b_{ji}}}
F_{kj} p_j F'_{ji}  
\end{split}
\end{equation}
So, by construction, 
we certainly get the outstanding item, \ref{5.12}.(c), of
stability under representable base change.
As such, it only remains to check the various commutativity
conditions of \ref{def:PosSeq} for these additional cells.
Amongst which the commutativity of \eqref{349} for \eqref{376}
just comes out in the wash in doing (c)' of the proof of
\ref{5.9}, and while there are several cases of \ref{PosSeq2}
and \ref{PosSeq4} to do, they're all immediate from either
the definitions; \ref{fact:PosSeq1}, or \ref{5.8} in combination
with the zig-zag axiom for 2-categories.
\end{proof}
This establishes the existence part of
\ref{5.2}, so there remains uniqueness, {\it i.e.}
\begin{fact}\label{5.15} Let $I$, $I'$ be Postnikov sequences
which are maximal amongst those which are right co-filtering and
stable under base change, with $K$, respectively $K'$
the continuous $\pi_2(\cX_*)$ valued 3 co-cycles they
define according to the prescription of \eqref{350}, then
there is a $\pi_2(\cX_*)$ valued continuous 2 co-chain, $z$,
whose boundary is $K-K'$.
\end{fact}
\begin{proof} By choice and \ref{5.13} there is an increasing
map $a:I\ra J$ such that $a(i)\geq i\in I_*$, for all $i\in I$.
In particular, therefore, there are projections,
$p_i: q'_i:=q_{a(i)}\ra q_i$ satisfying the base point condition
\eqref{eq:funct1406}. In order to lighten the notation we'll
confuse $I$, and $J$ with appropriate co-final systems defining
$\pi_2(\cX_*)$, equivalently take the consequences of stability
under base change for granted. Thus, just as in the proof of
item (b) in \ref{5.14} we can assert that for all $i,j\in I$
there are 2-commutative diagrams
\begin{equation}\label{378}
 \xy 
 (0,0)*+{q'_{i}}="A";
 (30,0)*+{q'_j}="B";
 (30,-18)*+{q_j}="C";
 (0,-18)*+{q_i}="D";
 (22, -6)*+{}="E";
 (8,-12)*+{}="F";
    {\ar_{F'_{ji}} "A";"B"};
    {\ar_{p_i} "A";"D"};
    {\ar^{p_j} "B";"C"};
    {\ar^{F_{ji}} "D";"C"};
 {\ar@{=>}_{ \xi_{ji}} "E";"F"}; 
 \endxy 
\end{equation}
and similarly for every $\om\in\pi_1(\cX_*)$, $i\in I$
another 2-commutative square
\begin{equation}\label{379}
 \xy 
 (0,0)*+{q'_{i}}="A";
 (30,0)*+{q'_i}="B";
 (30,-18)*+{q_i}="C";
 (0,-18)*+{q_i}="D";
 (22, -6)*+{}="E";
 (8,-12)*+{}="F";
    {\ar_{\om'_{i}} "A";"B"};
    {\ar_{p_i} "A";"D"};
    {\ar^{p_i} "B";"C"};
    {\ar^{\om_{i}} "D";"C"};
 {\ar@{=>}_{ \theta^\om_{i}} "E";"F"}; 
 \endxy 
\end{equation}
where, given the obvious strain on notation, we try
to keep the notation of \ref{def:PosSeq} for $I$,
while indicating the same for $I'$ by $'$. In any case,
the problem is that \eqref{378} and \eqref{379} may
fail to be compatible with the commutativity conditions 
\eqref{PosSeq2} and \eqref{PosSeq4} in the definition
of Postnikov sequences. In the first place, therefore,
for $i\geq j\geq k\in I$, define,
using \eqref{eq:plagAct15}, an element $z_{kji}$
of the stabiliser, $\G_k$, of the base point of $q_k$ by the
diagram
\begin{equation}\label{380}
 \xy 
 (0,0)*+{F_{ki}p_i}="A";
 (60,0)*+{F_{kj}p_j F'_{ji} }="B";
 (60,-18)*+{p_k F'_{kj}F'_{ji}}="C";
 (0,-18)*+{F_{ki}p_i}="D";
 (30,0)*+{F_{kj}F_{ji}p_i}="E";
 (30,-18)*+{p_k F'_{ki}}="F";
    {\ar@{=>}^{(F_{ki}p_i)^*z_{kji}} "D";"A"};
    {\ar@{=>}_{(F'_{ji})^* \xi_{kj}} "C";"B"};
    {\ar@{=>}_{p_i^*\g_{kji}} "A";"E"};
    {\ar@{=>}^{(F_{kj})_* \xi_{ji}} "B";"E"};
    {\ar@{=>}_{(p_k)_* \g'_{kji}} "C";"F"};
    {\ar@{=>}_{ \xi_{ki}} "F";"D"};
 \endxy 
\end{equation} 
The function $z_{kji}$ is, therefore, a 2 co-cycle in
the C\v{e}ch description, \cite[pg. 4]{jensen}, of the 
derived inverse limits of the $\G_i$ of the directed
system $I$. The groups $\G_i$ in question are, however,
finite, so \cite[Th\'eor\`eme 7.1]{jensen},  all the
higher direct limits vanish, and $z_{kji}$ is a 
co-boundary. As such by adjusting the $\xi_{\bullet}$'s
appropriately, we can suppose that $z_{kji}$ in \eqref{380}
is the identity. Similarly if $\Phi^{ji}_\om$ are 
the transition 2-cells of \ref{def:PosSeq}.(b) for
$I'$, we can define $z^{ji}_\om \in \G_j$ by
\begin{equation}\label{381}
 \xy 
 (0,0)*+{\om^j p_j F'_{ji}}="A";
 (25,0)*+{p_j\om'_j F'_{ji} }="B";
 (50,0)*+{p_j F'_{ji}\om'_i}="C";
 (75,0)*+{F_{ji}p_i\om'_i}="D";
 (75,-18)*+{F_{ji}\om^i p_i}="E";
 (37.5,-18)*+{\om^j F_{ji} p_i}="F";
 (0,-18)*+{ \om^j p_i F'_{ji} }="G";
    {\ar@{=>}^{( \om^j p_j F'_{ji})^*z^{ji}_\om} "G";"A"};
    {\ar@{=>}_{(F'_{ji})^* \theta_\om ^{j}} "B";"A"};
    {\ar@{=>}_{(p_j)_*\Phi^{ji}_\om} "C";"B"};
    {\ar@{=>}^{(\om'_{i})^* \xi_{ji}} "C";"D"};
    {\ar@{=>}^{(F_{ji})_* \theta^{i}_\om} "D";"E"};
    {\ar@{=>}_{ p_i^*\phi^{ji}_\om } "E";"F"};
    {\ar@{=>}^{ (\om^j)_* \xi_{ji}} "G";"F"};
 \endxy 
\end{equation} 
which by \ref{def:PosSeq}.(d) is a 1 co-cycle for the
directed system of the $\G_i$'s as soon as $z_{kji}$ in \eqref{380}
is the identity, so, it's a co-boundary, and adjusting
the $\theta^{\bullet}_\om$'s appropriately, we can
suppose that $z^{ji}_\om$ is the identity. Finally, therefore,
if we define $z^i_{\tau,\om}\in\G_i$ by the commutativity of
\begin{equation}\label{382}
 \xy 
 (0,0)*+{p_i (\tau\om)'_{i}}="A";
 (60,0)*+{p_i \tau'_{i}\om'_i}="B";
 (60,-18)*+{\tau^{i}\om^i p_i}="C";
 (30,-18)*+{(\tau\om)^i p_i}="D";
 (0,-18)*+{  (\tau\om)^i p_i}="E";
    {\ar@{=>}^{  (p_i)_* (\a')^{i}_{\tau,\om}} "B";"A"};
    {\ar@{=>}^{\tau^{i}_* \theta_\om ^{i}(\om'_i)^*\theta^i_\tau} "B";"C"};
    {\ar@{=>}_{p_{i}^* \a^i_{\tau,\om}} "C";"D"};
    {\ar@{=>}^{    ((\tau\om)^i p_i)^* z^i_{\tau,\om}} "E";"D"};
    {\ar@{=>}_{ \theta^i_{\tau\om} } "A";"E"};
 \endxy 
\end{equation} 
then since $z^{ji}_\om$ is the identity in \eqref{381},
the Postnikov condition \eqref{349} for $I$ and $I'$
implies that $z^i_{\tau,\om}\mpo z^j_{\tau,\om}$ under
the natural map $\G_i\ra \G_j$ of \ref{factdef:plagAct1}.
Better still, this holds for all $i\geq j$ such that the
2-commutative squares \eqref{378} and \eqref{379} exist,
so by stability under base change, {\it cf.} \ref{5.11},
$z^i:\pi_1(\cX_*)^2\ra \pi_2(\cX_*)$ is continuous, while
by \eqref{350} and \eqref{382} its differential is $K-K'$.
\end{proof}
Finally therefore we can make,
\begin{defn}\label{def:Pro2type}
The pro-2-type of $\et_2(\cX)$ is the triple 
$(\pi_1(\cX_*), \pi_2(\cX_*), k_3)$ consisting of
the pro-finite fundamental group, \ref{sum:plag1},
the pro-finite continuous $\pi_1(\cX_*)$-module,
$\pi_2(\cX_*)$, \ref{factdef:plagAct1}, and the
unique continuous co-cycle,
\begin{equation}\label{383}
k_3 \, \in\,  \rH^3_{\mathrm{cts}}\bigl( \pi_1(\cX_*), \pi_2(\cX_*)\bigr)
\end{equation}
defined as in \ref{5.2} by way of \ref{5.13} (existence)
and \ref{5.15} (uniqueness). Similarly the pro-finite
fundamental 2-group $\Pi_2(\cX_*)$ is (up to equivalence
of continuous 2-groups) defined 
exactly as in \eqref{eq:group1}-\eqref{eq:group3},
up to the simple further proviso of insisting that the
associator, \eqref{eq:group2}, is a continuous co-cycle
lifting \eqref{383}- by construction/definition \ref{factdef:plagAct1},
the group law \eqref{eq:group1} is already continuous. 
In particular, $\Pi_2(\cX_*)$ is a continuous 2-group
in the obvious internalisation of the term. 
\end{defn}
Notice that as an immediate consequence of the definitions
we have
\begin{cor}\label{Whitehead} (``Whitehead Theorem'')
Let $F:\cX\ra\cX'$ be a (pointed) map of connected champs then
$F_*:\Pi_2(\cX_*)\ra \Pi_2(\cX'_*)$ is an equivalence of
$2$-groups iff $F_*$ is an isomorphism on homotopy
groups $\pi_q$, $q=1$ or $2$.
\end{cor}
\begin{proof} Necessity is obvious, and sufficiency almost
as obvious since it implies that the pull back under $F$
of a 
maximal right co-filtering stable under base change
Postnikov sequence on $\cX'$ is 
again such a Postnikov sequence
over $\cX$, but, \ref{5.15}, these are unique.
\end{proof}  
\subsection{The fibre functor revisited}\label{SS:III.6}
The functors $F^*_{ji}$ of \eqref{eq:funct1106} and the
natural transformation enclosed by the triangle of
{\it op. cit.} do not a priori patch to a map from
$\varinjlim_i\Hom_{\et_2(\cX)} (q_i, q)$ to $q^{-1}(*)$.
By the universal property, \eqref{LimitA3}-\eqref{LimitA4}, of direct
limits of categories, the required condition for 
achieving this is that 
for $i\geq j\geq k$ we have a commutative diagram
of isomorphisms of base points, which at this juncture
are most conveniently understood as functors from
the trivial category $\rp$, {\it i.e.}
\begin{equation}\label{384}
 \xy 
 (0,0)*+{*_k}="A";
 (30,0)*+{F_{kj}(*_j)}="B";
 (30,-18)*+{F_{kj} F_{ji}(*_i)}="C";
 (0,-18)*+{F_{ki}(*_i)}="D";
    {\ar@{=>}_{\z_{kj}} "A";"B"};
    {\ar@{=>}_{\z_{ki}} "A";"D"};
    {\ar@{=>}^{F_{kj}(\z_{ji})} "B";"C"};
    {\ar@{=>}^{\g_{kji}(*_i)} "D";"C"};
 \endxy 
\end{equation} 
where $\g_{kji}$ is as per item (d) of the definition,
\ref{def:PosSeq}, of a Postnikov sequence. On the other
hand, if this diagram doesn't commute then following
it round, say clockwise, defines not only an element
$z_{kji}$ in the stabiliser, $\G_k$, of $*_k$, but as
immediately post \eqref{380} a  2 co-cycle in
the C\v{e}ch description, of the 
derived inverse limits of the $\G_k$, so
again by \cite[Th\'eor\`eme 7.1]{jensen} this
vanishes. Thus, without loss of generality we
may suppose that \eqref{384} commutes, and whence
we obtain a functor
\begin{equation}\label{385}
\varinjlim_i\Hom_{\et_2(\cX)} (q_i, q)\ra q^{-1}(*):
(y_i,\eta_i)\mpo (y_i(*_i), \eta_i(*_i)\phi_i)
\end{equation}
Consequently on introducing only the
the $\g_{kji}$ satisfying \ref{def:PosSeq}.(e) 
and the condition, \eqref{PosSeq3}, therein we have
\begin{fact}\label{fact:funct106}
The fibre 2-functor, \ref{def:funct102},
is pro-representable, {\it i.e.}
for $I_{*}$ the partially ordered set
of \ref{fact:funct103}, but with the $q_i$
not just connected but 2-Galois the functor
\eqref{385} is 
an equivalence of categories.
\end{fact}
\begin{proof} Essential surjectivity is 
\ref{fact:funct103} and \ref{fact:noplag5noplag5}, 
while by \ref{fact:funct105}
it's fully faithful.
\end{proof}
This equivalence of categories notwithstanding, the
fibre
functor in $\et_2(\cX)$- unlike that, {\it cf.} \ref{sum:plag1}, 
in $\et_1(\cX)$- when restricted
to a 2-Galois cell, $q_i$, yields a map
\begin{equation}\label{386}
\Hom_{\et_2(\cX)} (q_i, q_i)\ra q_i^{-1}(*) 
\end{equation}
which by \eqref{eq:noplag2} is potentially
very far from an equivalence of categories. In particular therefore,
on the one hand, {\it cf.} \ref{warn:plagAct1}, 
the right hand side of \eqref{386} does not trivially
inherit a 2-group structure from the left,
while, on the other: the natural restrictions  
$q_i^{-1}(*)\ra q_j^{-1}(*)$ on the right do not afford
a series of similarly compatible functors on the left.
Similarly, even given a Postnikov sequence, it's not
even true that the right hand side of \eqref{385}
has a unique $\Pi_2$ action. More precisely 
according to the definition of the same,
\ref{def:LimitB}, directed limits of categories
only make sense- 
{\it cf.} \ref{def:FibreS}-
up to equivalence, so there's
no reason for a $\Pi_2$ action whether on the
left or right of \eqref{385} to exhibit any further
uniqueness. On the other hand if the functors
and natural transformations
\begin{equation}\label{387}
U_i:\cH_i:=\Hom_{\et_2(\cX)} (q_i, q)\ra
\cH:=\varinjlim_i\Hom_{\et_2(\cX)} (q_i, q),\quad
u_{ij}:U_j\Rightarrow U_i F^*_{ji},\, i\geq j
\end{equation}
define the direct limit
in the sense of \ref{def:LimitB},
then for each $\om\in\pi_1(\cX_*)$
we have, for a Postnikov sequence, a further sequence 
of functors and natural transformations 
\begin{equation}\label{388} 
U_i\om_i^*:\cH_i\ra\cH,\quad \psi_\om^{ij}:=\bigl(U_j \om_j^* 
{\build\Longrightarrow_{}^{\om_j^*u_{ij}}} U_i F^*_{ji} \om_j^*
 {\build\Longleftarrow_{}^{(U_i)_*\phi_\om^{ji}} } U_i\om_i^* F^*_{ji} 
\bigr)
\end{equation}
which by \ref{def:PosSeq}.(b) and (d)
satisfy the requisite conditions of
\ref{def:LimitB} to yield a functor,
\begin{equation}\label{389}
\om:\cH\ra\cH,\,\text{and natural transformations,}\,\,
\theta^i_\om: U_i\om_i^*\Rightarrow \om U_i
\end{equation}
according to the prescriptions of the universal
property of {\it op. cit.}. Finally by \ref{def:PosSeq}.(c),
and a moderate diagram chase, the endomorphisms
$\om\tau$ and $(\tau\om)$ of $\cH$ solve the same universal
problem
so by \eqref{LimitA4} there is a unique natural transformation
$\lb_{\om,\tau}$ between them such that the
diagram
\begin{equation}\label{390}
 \xy 
 (0,0)*+{U_i(\tau_i\om_i)^*}="A";
 (30,0)*+{\om U_i \tau_i^*}="B";
 (60,0)*+{\om \tau U_i}="E";
 (60,-20)*+{(\tau\om)U_i}="C";
 (0,-20)*+{U_i(\tau\om)_i^* }="D";
    {\ar@{=>}_{(\tau_i^*)^*\theta^i_\om } "A";"B"};
    {\ar@{=>}_{(U_i)_* \a^i_{\tau,\om}} "A";"D"};
    {\ar@{=>}_{\om_* \theta^i_\tau } "B";"E"};
    {\ar@{=>}^{\theta^i_{\tau\om}} "D";"C"};
    {\ar@{=>}^{\lb_{\om,\tau}} "E";"C"};
 \endxy 
\end{equation}
commutes for all $i$, where $\a^i_{\tau,\om}$ is as per
item (c) of \ref{def:PosSeq}. Consequently the definition,
\eqref{350}, of the pro-finite Postnikov class combines
with the uniqueness of $\lb_{\om,\tau}$ to yield
\begin{fact}\label{fact:362} 
The solution of the universal problem \eqref{387} implied
by \eqref{385} affords (tautologically) on either side of
\eqref{385} a $\Pi_2(\cX_*)$-action (unique up to unique
equivalence given a Postnikov sequence) such that the equivalence
\ref{fact:funct106} is $\Pi_2$ equivariant.
\end{fact}
The discussion becomes more satisfactory still if we introduce
\begin{defn}\label{def:363} Let $\het_2(\cX)$ be the 
(big unless one supposes the existence
of a universe, or places some limiting cardinality on the
indexing sets, $I$, 
which follow,
such as the cardinality of 
something small and weak
equivalent to $\et_2(\cX)$ ) pro-2-category 
associated to $\et_2(\cX)$, {\it i.e.}
0-cells are 
2-functors from  some co-filtered
partially ordered set $I$ 
(albeit that 
directed, {\it i.e.}
right co-filtered
and partially ordered
are all one needs, and
arguably should, suppose
in the 
current
pro-finite, rather
than pro-discrete \ref{SS:III.9},
context)
viewed
as a category in the usual way, so,
specifically, \eqref{LimitA1},
sequences $\hat{q}=(q_i, F_{ji}, \g_{kji})$ 
for $i\geq j\geq k\in I$
belonging to a co-filtering partially ordered set (depending on
$\hat{q}$)
satisfying the prescriptions (a)-(c) of {\it op. cit.},
and whose 1/2-cells are given by understanding the bi-limit
\begin{equation}\label{391}
\Hom_{\hat{\cC}}(\hat{q}, \hat{q}'):=
\varprojlim_\a\varinjlim_i \Hom_{\cC} (q_i, q'_\a)
\end{equation}
in the (strict) 2-category of groupoids. In particular, there is something
here to check (a.k.a. a diagram chase for the reader), {\it i.e.}
that for $\a$ fixed the obvious sequences of categories, functors,
and natural transformations satisfy the requisite condition,
\ref{def:LimitA}, for the direct limit in $i$ to exist, and
that it's existence (up to unique equivalence) is sufficient
for the resulting sequence in $\a$ to satisfy the dual condition
for the existence of the inverse limit.
\end{defn}
By way of application 
one observes
\begin{fact}\label{fact:364} Let $\hat{q}$ be the $0$-cell in
$\het_2(\cX)$ associated to a Postnikov sequence
via \ref{def:PosSeq}.(d)-(e), then 
the fibre functor \eqref{385} extends to an
equivalence of categories
\begin{equation}\label{392}
\Hom_{\het_2(\cX)}(\hat{q}, \hat{q}')\ra \varprojlim_\a (q'_a)^{-1}(*)
\end{equation}
while the $\Pi_2(\cX_*)$-action of \ref{fact:362} extends
to a continuous action on (say the right) understood as
a continuous (small) groupoid.
\end{fact}
\begin{proof} That this is an equivalence of categories 
is automatic from \ref{fact:funct106} and the definition
\ref{def:363}, and, similarly, if the $\Pi_2$ action
exists then it's certainly continuous because \ref{fact:362}
is. Otherwise, the existence of an action is an exercise
in the universal properties of direct and inverse limits
of categories. More precisely if $\hat{q}'=(q'_a, P_{\a\b}, \rho_{\a\b\g})$,
for $\g\geq\b\geq\a$ in some directed set $A$, then denoting
the defining Hom-categories in \eqref{391} by $\cH_i^\a$,
and their limits in $i$ by $\cH^\a$, we have on adding
the superscript $\a$ in \eqref{387}-\eqref{390}, 
a set of solutions to a set of universal problems indexed
by $\a$. Similarly a small (but non-empty, since one must
respect the conditions in \ref{def:LimitA}) diagram chase
shows that the universal property of direct limits affords
a 2-commutative square
\begin{equation}\label{393}
 \xy 
 (0,0)*+{\cH_i^\b}="A";
 (30,0)*+{\cH_i^\a}="B";
 (30,-18)*+{\cH^\a}="C";
 (0,-18)*+{\cH^\b}="D";
 (22, -6)*+{}="E";
 (8,-12)*+{}="F";
    {\ar_{P_{\a\b}} "A";"B"};
    {\ar_{U_{i}^\b} "A";"D"};
    {\ar^{U_{i}^\a} "B";"C"};
    {\ar_{ P_{\a\b}} "D";"C"};
 {\ar@{=>}_{ u_i^{\a\b}} "F";"E"}; 
 \endxy 
\end{equation}
where $P_{\a\b}$ is understood in the obvious way. To all
of which one can apply the unicity, \eqref{LimitA4}, 
of natural transformations
of functors from $\cH_\b$ to obtain 
for each $\om\in\pi_1(\cX_*)$ unique natural
transformations $\psi_\om^{\a\b}$ such that the
following diagram commutes
\begin{equation}\label{394}
 \xy 
 (0,0)*+{P_{\a\b}U_i^\b\om_i^\b}="A";
 (50,0)*+{U_i^\a P_{\a\b}\om_i^\b=U_i^\a\om_i^\a P_{\a\b}}="B";
 (100,0)*+{\om_i^\a U_i^\a P_{\a\b}}="E";
 (100,-20)*+{\om^\a P_{\a\b}U^\b_i }="C";
 (0,-20)*+{P_{\a\b} \om^b U^\b_i}="D";
    {\ar@{=>}_{(\om_i^\b)^* u_i^{\a\b} } "A";"B"};
    {\ar@{=>}_{(P_{\a\b})_* \theta^i_{\om^\b}} "A";"D"};
    {\ar@{=>}_{P_{\a\b}^* \theta^i_{\om^\a} } "B";"E"};
    {\ar@{=>}_{\om^\a_* u^{\a\b}_i} "C";"E"};
    {\ar@{=>}^{(U_i^\b)^*\psi^{\a,\b}_\om} "D";"C"};
 \endxy 
\end{equation}
for all $i$. Now the unicity (and some diagram chasing)
does the rest, {\it i.e.} 
the categories $\cH^\a$;the functors $P_{\a\b}$, $\om_\a$; and
natural transformations $\psi_\om^{\a\b}$, $\lb^\a_{\om,\tau}$,
$\rho_{\a\b\g}$ satisfy (modulo the change from right to left)
the same diagrams as found in the definition \ref{def:PosSeq}.(a)-(e)
of a Postnikov sequence, which by the same/dual argument to
\eqref{387}-\eqref{390} are exactly the required conditions
for defining a $\Pi_2$-action on $\varprojlim_\a\cH^\a$.
\end{proof}
From which we deduce the conceptually satisfactory
\begin{cor}\label{fact:365} Again let $\hat{q}$ be the $0$-cell in
$\het_2(\cX)$ defined by a Postnikov sequence, then the 
fibre functor affords an equivalence of 2-groups
\begin{equation}\label{395}
\Pi_2(\cX_*)\xleftarrow{\sim}
\fAut_{\het_2(\cX)} (\hat{q}) \, 
\bigl(=\Hom_{\het_2(\cX)}(\hat{q}, \hat{q})\bigr)\,
\end{equation}
\end{cor}
\begin{proof}
Let $\gP$ be the 2-group on the right, then, 
in the notation of the proof of \ref{fact:364}, its
identity can be identified with the functors
$F_{\a i}\in\cH_i^\a$ 
of \ref{def:PosSeq}.(b)
and (when taking the
inverse limit) the natural transformations $\g_{\g\b\a}$
of {\it op. cit.} (d). In particular it
maps under \eqref{395} to an object equivalent to the object $\hat{*}$
in $\varprojlim_\a q_\a^{-1}(*)$ defined by the square
\eqref{384}- {\it cf.} \eqref{LimitA3}-\eqref{LimitA4}. Similarly the
objects $\om_\a F_{\a i}$ (which are equivalent to
the sequence of objects $F_{\a i}\om_i$ by \ref{def:PosSeq}.(b) \& (d) )
map to an object isomorphic to the object $\om(\hat{*})$
obtained by applying $\om$ on the left to the said
diagram \eqref{384}. By \ref{fact:364} this is all
the objects in $\Hom_{\het_2(\cX)}(\hat{q}, \hat{q})$
up to equivalence, so every endomorphism of $\hat{q}$
is a (weak)  automorphism, and whence by \ref{sum:plag1}
\begin{equation}\label{396}
\pi_1(\gP)\xrightarrow{\sim} \pi_1(\cX_*)
\end{equation}
Similarly, \ref{fact:364} and \eqref{eq:plagAct18plagAct18} yields
a group automorphism 
\begin{equation}\label{397} 
\pi_2(\gP)\xrightarrow{\sim} \pi_2(\cX_*) 
\end{equation}
in which both sides are $\pi_1(\gP)$ (respectively
$\pi_1(\cX_*)$) modules in a compatible way thanks
to the equivariant part of \ref{fact:364}. Finally
the definition of their respective Postnikov classes
is the same by construction, {\it i.e.} \ref{def:PosSeq}.(c) and
\eqref{350}.
\end{proof}
Finally, let's tie things up with 
the actions encountered in
\ref{fact:plagAct2}-\ref{warn:plagAct1} by way of
\begin{rmk}\label{rmk:366}
We have of course the possibility to identify the fibre
over the base point of any 0-cell with a unique (if not
quite uniquely unique) skeletal groupoid, {\it cf.} \ref{SS:II.1}.
In the particular case that the cell is 2-Galois, or
even just connected with abelian stabiliser, this is
sufficient to get a unique fibre functor in \eqref{385}
given \eqref{384}, and similarly on taking inverse limits
of the same. There is therefore little, or no, ambiguity
for the vertices in the action diagram \eqref{eq:plagAct4}.
The top horizontal arrow of {\it op. cit.} is, however,
another matter, \ref{claim:noplag1}. Nevertheless, all of
the problems that it threatens to cause, \ref{warn:plagAct1},
result from the lack of an honest universal cover, which
not only, as one might expect, disappear in the limit, \ref{SS:III.5},
by virtue of the existence of a Postnikov sequence, but
even the said existence itself admits the conceptually
pleasing description \eqref{395}- {\it cf.} \cite[V.4.g]{sga1}
for $\pi_1$- which parallels the topological case,
\ref{fact:cor3} as closely as is possible.
\end{rmk}

\subsection{The pro-2-Galois correspondence}\label{SS:III.7}
The pro-finite structure of $\Pi_2:=\Pi_2(\cX_*)$ naturally
affords it the structure of a continuous 2-group, and, similarly,
one defines a continuous
groupoid with continuous $\Pi_2$-action by the simple 
expedient of supposing that arrows in \ref{SS:II.1} are,
where appropriate, continuous, and 
amongst all such we may distinguish
\begin{defn}\label{def:371}
The 2-category $\Grpd_{\cts}(\Pi_2)$ is the full 
({\it i.e.} the same 1 and 2 cells) sub 2-category
of continuous $\Pi_2$ equivariant groupoids in
which the zero cells $\cF$ are 
groupoids satisfying
\begin{itemize}
\item[(a)] The topology of the space of
arrows, $\mathrm{Fl}(\cF)$ is 
compact separated and totally disconnected,
or, equivalently, given compact separated,
totally separated, {\it cf.} \ref{def:391} and \eqref{ProF1}.
\item[(b)] The identity map, 
$\mathrm{Ob}(\cF)\xrightarrow{\mathbf{1}} \mathrm{Fl}(\cF)\ts \mathrm{Fl}(\cF)$
is a closed embedding.
\end{itemize}
\end{defn}
An alternative description of the above conditions is therefore
\begin{fact}\label{fact:372} 
The 2-category of
topological groupoids $\cF$ satisfying \ref{def:371}.(a)-(b) 
is equivalent to the (pro) 2-category of pro-finite
groupoids, {\it i.e.} $\het_2(\rp)$ in the sense
of \ref{def:363}. 
\end{fact}
\begin{proof} There is, 
\eqref{ProF1}, a functor
\begin{equation}\label{1372}
X\mpo \pi_0^\fin (X)
\end{equation}
from topological spaces
to pro-finite sets, which is an equivalence
when restricted to the full subcategory of
topological spaces satisfying \ref{def:371}.(a). As such
the space of arrows of
a groupoid $(s,t):R\rras X$ satisfying (a) is a pro-finite
set, while if it also satisfies (b) then $X$ is a pro-finite
set too. Consequently $R\rras X$ is equivalent to a
pro-finite groupoid object, so by \cite[A.3.3]{artin-mazur}
it's a 0-cell in pro-finite groupoids of a 
particularly simple form, {\it i.e.}
\begin{equation}\label{2372}
\begin{split}
& \text{(a) There are strict surjections $R\ra R_i$ 
to finite quotient groupoids $R_i$, $i\in I$.}\\
& \text{(b) The induced inverse system $F_{ji}:R_i\ra R_j$
of functors satisfies $F_{ki}=F_{kj}F_{ji}$.}\\
& \text{(c) The set map $R\ra \varprojlim_i R_i$ implied
by (b) is a strict isomorphism of groupoids.}
\end{split}
\end{equation}
Such simplicity notwithstanding, we therefore have
a 2-functor
\begin{equation}\label{3372}
\text{Groupoids satisfying \ref{def:371}.(a)-(b)}
\ra \het_2(\rp): [R\rras X]\mpo [R_i\rras X_i]_{i\in I}
\end{equation}
for $I$ (depending on $R$) the co-filtered directed set
of \eqref{2372}. In principle there could be many more
0-cells in $\het_2(\rp)$, but \eqref{3372} is actually
essentially surjective. Indeed a 0-cell, 
$(R_i\rras X_i, F_{ji}, \g_{kji})$, $i,j,k\in I$
of $\het_2(\rp)$
as described by \eqref{LimitA1}, may by the fact that
\eqref{1372} is an equivalence of categories under
the hypothesis \ref{def:371}.(a), be identified with 
the inverse limit, 
$\hat{R}\rras \hat{X}$ of
\ref{exInv} albeit with the
proviso that arrows and objects are understood to have
the obvious topology. Supposing, as we may, 
that $R_i\rras X_i$ is skeletal, we may, without
loss of generality, further suppose, \eqref{LimitA10},
that if
$x=(x_i) \in\varprojlim_i X_i$ then
there are automorphisms 
\begin{equation}\label{4372}
\xi_{ji}(x) \in \Aut (x_j), \,\, x_j=F_{ji}(x_i),\,\,
\text{such that}\,\, \xi_{ki}(x)=\xi_{kj}F_{kj}(\xi_{ji})\g_{kji}(x)
\end{equation}
Consequently if we define functors $F'_{ji}$ by the
rule,
\begin{equation}\label{5372}
F_{ji} {\build\Rightarrow_{}^{\xi_{ji}(x)}} F'_{ji}
\end{equation}
then the $F'_{ji}$ satisfy the strict commutativity
condition \eqref{2372}.(b), while for $a\in I$ fixed,
with $i\geq a$,
the functors $F'_{ai}:R_i\ra R_a$ and natural transformations
\begin{equation}\label{6372}
F_{ji}^*(\xi_{aj})\g_{aji} \xi_{ai}^{-1}: F'_{ai}\Rightarrow F'_{aj}F_{ji}
\end{equation}
define an isomorphism in $\het_2(\rp)$ between our
given $0$-cell and that defined by
$(R_i\rras X_i, F'_{ji}, \mathbf{1})$.
In particular, therefore, to each object $x\in\varprojlim_i X_i$
there is associated a pro-finite group,
\begin{equation}\label{7372}
S_x := \varprojlim_i [\Hom_{R_i}(x_i, x_i)\xrightarrow{F'_{ji}(x_i)}
\Hom_{R_j}(x_j, x_j)]
\end{equation}
which (by definition) glue together to a continuous 
$\varprojlim_i X_i$-group $S:=\cup_x S_x$. Nevertheless,
there could still be, \eqref{LimitA10},  
many more objects in $\hat{X}$ than
$\varprojlim_i X_i$, but in the presence 
of the strict commutativity
condition \eqref{2372}.(b)   
and \eqref{4372} 
(with $\g=\mathbf{1}$)
every object as
described by \eqref{LimitA10} is isomorphic
to one in $\varprojlim_i X_i$ by
\cite[Th\'eor\`eme 7.1]{jensen} (or more
accurately the non-commutative variant
for $\rH^1$). As such, if we further replace
the object sets $X_i$ by the image of $\varprojlim_i X_i$,
and the arrows by the image of $S$ defined 
by \eqref{7372} then we not only deduce
that the 2-functor \eqref{3372} is essentially
surjective on 0-cells, but that we can add
to \eqref{2372} the further property
\begin{equation}\label{8372}
\text{(d). Whether $R\rras X$, or the
$R_i\rras X_i$ are skeletal groupoids.}
\end{equation}
This said, in order to prove that \eqref{3372}
is essentially surjective on $1$-cells, {\it i.e.}
the general prescription \eqref{391} may
be reduced to that of a 
net of functors $A_i:R'_i\ra R''_i$
forming strictly commutating diagrams
\begin{equation}\label{9372}
\begin{CD}
R'_i @>>{F'_{ji}} > R'_j \\
@V{A_i}VV @VV{A_j}V \\
R''_j @>{F''_{ji}} >> R''_j
\end{CD}
\end{equation}
wherein a priori the implied 0-cells
only satisfy \eqref{2372}.(a)-(c), but
a postiori also satisfy \eqref{8372}.(d),
so essential surjectivity again reduces
to the non-commutative (for $\rH^1$)
variant of \cite[Th\'eor\`eme 7.1]{jensen}.
Similarly the last thing we have to check,
\cite[1.5.13]{tom}, {\it i.e.} maps between
1-cells are a family of faithful functors
is again a non-commutative 
variant of \cite[Th\'eor\`eme 7.1]{jensen},
albeit, on this occasion, for $\rH^0$.
\end{proof}
To characterise the action of $\Pi_2$ in pro-finite
terms is similar, to wit:
\begin{fact}\label{fact:373} 
Let $\Grpd_{\mathrm{fin}}(\Pi_2)\subset \Grpd_{\cts}(\Pi_2)$
be the sub 2-category of finite groupoids with continuous $\Pi_2$-action
then there is an equivalence of 2-categories,
\begin{equation}\label{1399}
\Grpd_{\cts}(\Pi_2)\xrightarrow{\sim}
\mathrm{pro}-\Grpd_{\mathrm{fin}}(\Pi_2)
\end{equation}
where, following \ref{def:363} the
$0$-cells, $\cF$, of the latter are equivalent to sequences
(indexed by a co-filtered directed set $I$ in the sense of
\ref{def:LimitA} in an arbitrary 2-category)
$(\cF_i, f_{ji}, \xi_{kji})$ of $0$, $1$, and $2$-cells in 
$\Grpd_{\mathrm{fin}}(\Pi_2)$, with $1$ and $2$-cells given
by the bi-limit of categories 
\begin{equation}\label{2399}
\varprojlim_\a\varinjlim_i \Hom_{\Pi_2} (\cF_i, \cF'_\a)
\end{equation}
understood in the sense of \ref{def:LimitB}.
%
\end{fact}
\begin{proof} 
By \ref{fact:372} we need only check that the
mutually inverse 2-functors implicit in 
\eqref{3372} may be taken to be $\Pi_2$-equivariant.
On the one hand, if we start with a $\Pi_2$ equivariant
groupoid on the left hand side of \eqref{1399}, then
since 
(by definition \ref{def:Pro2type})
$\Pi_2$ itself satisfies \ref{def:371}.(a) we
can again appeal to the
fact that \eqref{1372} is an equivalence between
pro-finite sets and topological spaces satisfying
\ref{def:371}.(a) so \cite[A.3.3]{artin-mazur} applies
once more, and whence,
in the equivariant case, we can adjoin to
\eqref{2372}.(a)-(c) and \eqref{8372}.(d),
\begin{equation}\label{3399}
\text{(e). The action may be expressed as the
limit of finite actions $\Pi^i_2 \uts R_i\ra R_i$.}
\end{equation}
for a possibly different co-filtering directed set $I$.
Conversely, if one writes out what the sequence condition,
\ref{def:LimitA},
means in $\Grpd_{\mathrm{fin}}(\Pi_2)$ by way of 
\eqref{eq:groupPlus27}-\eqref{eq:groupPlus28} then (modulo
notation) one finds the conditions (a)-(e) of 
the definition, \ref{def:PosSeq}, of a Postnikov sequence
together with the proviso that the differential of the
natural transformations 
$\a_{\bullet,\bullet}$
of item (c) of {\it op. cit}
should be the image of the Postnikov class of $\Pi_2$.
As such, exactly the same arguments via the 
universal property of limits employed in \ref{fact:364},
imply that the inverse limit,  
$\varprojlim_i \cF_i$, in the sense of \ref{def:LimitB},
of a 0-cell on the right of
\eqref{1399} admits a
unique, up to equivalence, $\Pi_2$-action. 
\end{proof}
Given this is the sum total of what we'll need to
know about continuous actions it's appropriate to
emphasise a couple of points, {\it i.e.}
\begin{rmk}\label{rmk:0374}
Since all of \eqref{2372}.(a)-(c), \eqref{8372}.(d),
\eqref{3399}.(e) may be supposed to hold simultaneously,
the structure of pro-$\Grpd_{\mathrm{fin}}(\Pi_2)$ is
a lot simpler than that of an arbitrary 
(defined in the the obvious way suggested by
\ref{def:363}, \eqref{391}  \& \eqref{2399})
pro-2-category. We will, however, always eschew
(unless it's logically impossible) any appeal to
such simplifications which are not available
in the general pro-discrete case of \ref{SS:III.9}
\end{rmk}
\begin{rmk}\label{rmk:374}
Evidently 
most 
of \ref{SS:II.1} just extends to the continuous
case by the simple expedient of insisting that all
maps occurring therein are continuous.
The one point where one should exercise caution
is that it's not always possible to reduce to
the transitive case since $\pi_0$ of a pro-finite
set/equivalently a space satisfying \ref{def:371}.(a)
is naturally topologised, \ref{def:391}, \ref{def:392},
\ref{ProFF}. If, however, one does have a 
transitive action then \ref{fact:groupPlus1}
and the resulting description
\ref{sum:group1} 
of $\Grpd(\Pi_2)$ in terms of group co-homology 
are equally valid for $\Grpd_{\cts}(\Pi_2)$ with
exactly the same statements albeit understood 
in continuous group co-homology.
\end{rmk}
With this in mind we have
\begin{prop}\label{prop:376} 
Let $\hat{q}$ be the $0$-cell in
$\het_2(\cX)$ defined by a Postnikov sequence, then the 
2-functor ({\it cf.} \ref{fact:cor1}) 
\begin{equation}\label{3103}
\fHom_{\het_2(\cX)}(\hat{q}, \_):\het_2(\cX)\ra \Grpd_{\cts}(\Pi_2)
\end{equation}
affords an equivalence of 2-categories, wherein the $\Pi_2$-groupoid
$\Hom_{\het_2(\cX)}(\hat{q}, \hat{q}')$ is, \ref{fact:364},
equally described by the fibre functor \eqref{385}.
\end{prop}
An immediate simplification can be obtained by observing
\begin{lem}\label{lem:377} Let 
everything be as above, then \ref{prop:376} holds iff 
\begin{equation}\label{3104}
\fHom_{\et_2(\cX)}(\hat{q}, \_):\et_2(\cX)\ra \Grpd_{\mathrm{fin}}(\Pi_2)
\end{equation}
affords an equivalence of 2-categories, wherein the $\Pi_2$-groupoid
$\Hom_{\et_2(\cX)}(\hat{q}, q')$ is described 
by the fibre functor \eqref{385} and \ref{fact:funct106}-\ref{fact:362}.
\end{lem}
\begin{proof} Necessity is obvious, and otherwise, 
if \eqref{3104} is an equivalence, then it affords
an equivalence between pro-$\et_2(\cX)$ and
pro-$\Grpd_{\mathrm{fin}}(\Pi_2)$ which are
the left hand side of \eqref{3103} by definition, \ref{def:363},
respectively the right hand side by 
\ref{fact:373}.
\end{proof}

Following the schema of the topological proof, \ref{SS:II.5},
the step corresponding to \ref{fact:cor2} has already, \ref{fact:funct106}, 
been done while the appropriate variant of \ref{fact:cor4} is

\begin{fact}\label{fact:378}
The $2$-functor \eqref{3104} restricted to $1$-cells is a
family of fully faithful functors.
\end{fact}
\begin{proof} Let $F$, $G$ be a pair of 2-cells as in
\ref{fact:cor4} but with source a
2-Galois $q$. By \ref{fact:flav1} the diagram \eqref{eq:cor6}
is valid and, \ref{rmk:funct101}, all the maps are
proper even if $\cX$ isn't separated. Consequently,
in the notation of \eqref{eq:cor6} the support (on $\cY_1$) of
the sheaf of sets
\begin{equation}\label{3105}
\z_1: f_1\Rightarrow g_1,\quad q'_1(\z_1)\xi_1=\eta_1
\end{equation}
is closed. It is manifestly also open, and by construction
every stalk has at most one element, so it's the trivial
sheaf since it's non-empty at the base point. The related sheaf
\begin{equation}\label{3106}
\z: f\Rightarrow g,\quad q'(\z)\xi=\eta
\end{equation}
is therefore a locally constant sheaf on $\cY$ with 
fibres isomorphic to the stabiliser of $f(*)$. In 
principle, on a given 2-Galois cell, it could be
non-trivial, but on replacing $q$ by a Galois cover $q_j$
we can suppose that \eqref{3105} is generated by
global sections, which are isomorphic in a unique
way to the stabiliser of $f(*_j)$, so on passing to
the limit $\hat{q}$ all such isomorphisms are compatible.
\end{proof}
It therefore remains to prove
\begin{fact}\label{fact:379} The 2-functor \eqref{3104} is essentially
surjective on 0 and 1-cells.
\end{fact}
Now, as it happens, we still haven't proved this at the
level of $\et_1(\cX)$, and while we could have recourse
to \cite{noohi1}, it's not only appropriate to give
a proof in order to set up our notation and exhibit the
similarity with \ref{SS:I.5}, but also to clear up a 
certain amount of misconception about the axioms in 
\cite[Expos\'e V.4]{sga1} a number of which become (even
in an arbitrary site with enough points) superfluous on
working with champs. In particular, as the proof of the
following will demonstrate the conditions on the existence
of quotients by finite groups of automorphisms can be dropped.
\begin{fact}\label{fact:380} At the level of $\et_1(\cX)$ the
2-functor \eqref{3104} defines an essentially surjective 
functor to finite 
$\pi_1(\cX_*)$-sets.
\end{fact}
\begin{proof} We immediately reduce to the case where the
action is transitive, so what has to be shown is that if $q_1:\cY_1\ra\cX$ 
is a representable \'etale Galois cover with group $G$, then
for any finite sub-group $K$ there is an \'etale cover
$q'_1:\cY'\ra \cX$ with fibre $G/K$, which, of course, one
exhibits as a quotient of $q$ in $\et_1(\cX)$. 

To this end, 
profiting from \ref{fact:funct101} (with $f$ the identity)
we can, as in \eqref{eq:sep3}, represent $q_1$ by an embedding
of groupoids
\begin{equation}\label{3107}
q_1: R_1 \hookrightarrow R_0
\end{equation}
with \'etale actions on the same atlas $U$, while, similarly
an element, $\om$, of $G$ can be identified with a 
2 commutative diagram
\begin{equation}\label{3108}
 \xy 
 (0,0)*+{R_1}="A";
 (36,0)*+{R_1}="B";
 (18,-12)*+{R_0}="C";
{\ar_{}^{F_{\om}} "A";"B"};
    {\ar^{q_1} "B";"C"};
    {\ar_{q_1} "A";"C"};
{\ar@{=>}^{\iota_\om} (14,-6);(22,-6)}
\endxy
\end{equation}
Furthermore the definition of $\et_1(\cX)$ implies, 
just as in \eqref{eq:coverFix1}, that the natural
transformations $\a_{\tau,\om}:F_{\tau}F_\om\ra F_{\tau\om}$
of the $G$ action
are, in fact, given by,
\begin{equation}\label{3109}
\a_{\tau,\om} = \iota_{\tau\om}\iota_\om^{-1} (F_\om^* \iota_\tau)^{-1}
\end{equation}
Now, {\it cf.} \eqref{eq:coverFix}, consider the fibre square
\begin{equation}\label{3110}
\begin{CD}
R_{1,\om} @>{t_\om} >> R_1\\
@VVV @VV{s=\mathrm{source}}V \\
U@>{\om}>> U
\end{CD}
\end{equation}
where, unlike \eqref{eq:coverFix}, $t_\om$ 
need not be an embedding. Nevertheless we still have an isomorphism
\begin{equation}\label{3111}
c_{1,\om}:
R_{1,\om} \xrightarrow{\sim} R_0^\om \subseteq R_0: (f, x) \ra f\iota_\om(x)
\end{equation}
which also serves as a definition of the image $R_0^\om$. 
On the other hand, the fibre of $q_1$ over a point $x\in U$
can be identified with arrows in $R_0$ with source $x$ 
modulo the relation that there difference belongs to $R_1$,
and
by \eqref{eq:plag4} the $\iota_\om (x)$ are a complete
repetition free list of representatives of such equivalence
classes. Consequently, {\it cf.} \ref{lem:cover2}, as a set
\begin{equation}\label{3112}
R_0 = \coprod_{\om\in G} R_0^\om
\end{equation}
or, what is slightly better that the groupoid $\coprod_\om R_{1,\om}$
of which an element $(f,x)_\om$ has source $x$ and sink that of $f$,
with multiplication defined by
\begin{equation}\label{3113}
R_{1,\tau\, t}\ts_s R_{1,\om} \ra R_{1,\tau\om}: (g,z)_\tau \ts (f,x)_\om
\mpo (gF_\tau(f) \a_{\tau,\om}^{-1}, x)_{\tau\om}
\end{equation}
is strictly isomorphic to $R_0$.
\end{proof}
All of which can usefully considered as a lemma in
the main business, {\it i.e.}
\begin{proof}[proof of \ref{fact:379}]
We follow the proof of the topological case, \ref{fact:cor5},
with minor modifications. In the case of $0$-cells by
\ref{fact:380} we can, again, without loss of generality,
suppose that the representation has values in $\rB_\G$ for
some finite group $\G$. Consequently there is a 2-Galois
cell $q$, factoring as a locally constant gerbe $q\ra q_1$
in $\rB_{\pi'_2}$'s and a representable Galois cover $q_1\ra \cX$
with group $\pi'_1$ such that the action may be identified
with normalised co-chains $A:\pi'_1\ra\Aut(\G)$, $\z:(\pi'_1)^2\ra \G$,
and a map of $\pi'_1$ modules, $A_2:\pi'_2\ra Z$. with values
in the centre of $\G$ satisfying exactly the relations detailed
in \ref{fact:cor5}. In particular, we can identify $q_1$ with
the embedding  of groupoids \eqref{3107} and $q$ with a 
groupoid $R_2\ra R_1$ which is equally an \'etale $\pi'_2$-torsor.
By \ref{fact:380} (applied to $R_1$, and without, by the way
way, having to suppose that quotients by finite groups exist
in the ambient big site, {\it cf.} \cite[Expos\'e XI.5]{sga1}) 
there is, {\it cf.} \eqref{eq:cor7}-\eqref{eq:cor8}, a left right $\G$ 
torsor
\begin{equation}\label{3114}
R'_2 := R_2 \ts \G/\bigl( (f, \g) \sim (f S, \g A_2(S) ),
\quad S\in \pi'_2
\end{equation}
which equally defines a groupoid with action 
factoring as $R'_2 \ra R_1 \rras U$, along with
functors $F'_\om$, $\om\in G$,
and natural transformations $\a'_{\tau,\om}: F'_{\tau\om}
\Rightarrow F'_\tau F'_\om$ exactly as in  {\it op. cit.}.

Now consider adding a further fibre square in \eqref{3110}, to wit 
\begin{equation}\label{3115}
\begin{CD}
R'_{2,\om} @>>> R_{1,\om} @>>> U \\
@VVV @VVV @VV{\om}V \\
R'_2 @>>> R_1 @>{s=\mathrm{source}}>> U 
\end{CD}
\end{equation}
and which equally serves as the definition of the term
in the top left. In particular therefore we can write
elements of $R':=\coprod_\om R'_{2,\om}$ as pairs
$(f,x)_\om$ indexed by $\om\in\pi'_1$ with $f\in R'_2$
enjoying source $\om x$, so that exactly as in \eqref{eq:cor10}
\begin{equation}\label{3116}
R'_{2,\tau\, t}\ts_s R'_{2,\om} \ra R'_{2,\tau\om}: (g,z)_\tau \ts (f,x)_\om
\mpo (gF_\tau(f) \a'_{\tau,\om}, x)_{\tau\om}
\end{equation}
defines a groupoid structure which factors as 
$R'\ra R_0\rras U$ where the $\G$-torsor structure
of $R'\ra R_0$ is equally the quotient by the 
relative stabiliser thanks to \eqref{3113}. Furthermore,
since all co-chains are normalised, identities in $R'_2$
are the identities in $R'$, while the 
inverse of $(f,x)_\om$ is given by
\begin{equation}\label{3117}
\bigl((\a'_{\om^{-1},\om})^{-1} F'_{{\om}^{-1}} (f)^{-1}, t(f)\bigr)
\end{equation}
As such an explicit lifting of the natural transformation
\eqref{3108} is 
\begin{equation}\label{3118}
\iota'_\om: U\ra R': x\ra (\mathbf{1}_{\om x}, x)_{\om}
\end{equation}
which by \eqref{3117} satisfies \eqref{eq:cor12}, so by
exactly the same diagram chase, \eqref{eq:cor13}-\eqref{eq:cor15},
the 0-cell $[U/R']$ in $\et_2(\cX)$ affords the $\Pi_2$-action
on $\rB_{\G}$ that we started with.

As to the essential surjectivity on 1-cells, modulo replacing
the universal widget $R_2\ra R_1$ in the topological
proof, by groupoids $R_2\ra R_1 \rras U$ corresponding
to the factorisation of a (sufficiently large) 2-Galois cell
into a locally constant gerbe followed by a representable
map as above, everything is, otherwise, formally the
same as the topological case, \eqref{eq:cor16}-\eqref{eq:cor19}.
\end{proof}

\subsection{Relation with abelian and non-abelian co-homology}\label{SS:III.8}

The lack of an universal cover in the pro-finite theory 
notwithstanding, the 
relation of $\et_2(\cX)$ with co-homology remains 
essentially the same as the topological class, {\it i.e.}
\ref{fact:pos2}. To the end of describing this let us
introduce
\begin{defn}\label{def:381}
Let $\underline{Z}$ be a locally constant sheaf of abelian
groups on $\cX$ then the set of locally constant gerbes
$q:\cY\ra\cX$ with relative stabiliser (which always
descends to a sheaf on $\cX$ whenever it's abelian) 
$S_{\cY/\cX}\xrightarrow{\sim} q^*\underline{Z}$ 
(implicitly with an isomorphism of the stabiliser
of a base point with a fixed abelian group $Z$
to avoid the problem posed by \eqref{eq:groupTom})
modulo
equivalences in $\et_2(\cX)$ will be written $\htg(\cX, \underline{Z})$.
On representing $q$, respectively $q'$ by
groupoids $R$, $R'$ over $R_0$ \`a la \eqref{eq:sep3}
which are torsors under $\uZ$;
this set admits a natural abelian group structure,
\cite[IV.3.3.2]{giraud},
by way of Giraud's contracted product
\begin{equation}\label{3119}
q\wedge^{\uZ} q' \, := \, R\ts_{R_0} R'/\uZ
\end{equation}
\end{defn}
Now while we've only defined this for 
what is largely our context, {\it i.e.} locally constant sheaves,
it continues to have sense for arbitrary sheaves, and
to any short exact sequence of sheaves,
\begin{equation}\label{3120}
0\ra \uZ'\ra \uZ\ra \uZ''\ra 0
\end{equation}
one has, \cite[IV.3.4.1]{giraud}, 
\ref{GiraudWarn},
a tautological connecting homomorphism
\begin{equation}\label{3121}
\rH^1(\cX, \uZ') \xrightarrow{\d} \htg (\cX, \uZ'')
\end{equation} 
which sends a $\uZ'$-torsor to the gerbe whose objects
are its possible liftings as a $\uZ$ torsor. Unsurprisingly,
{\it cf. op. cit.}, this affords a $\d$-functor in degrees
$\leq 2$, and whence a natural map, or better a natural
transformation of functors
\begin{equation}\label{3122}
\eta_{\uZ}:\rH^2 (\cX, \uZ) \ra \htg (\cX, \uZ)
\end{equation}
and one can certainly show, \cite[IV.3.4.3]{giraud}, that
$\htg$ vanishes on injectives, and, whence, that \eqref{3122}
is an isomorphism. In our current considerations, however,
it better illustrates the connection with \cite{artin-mazur}
to proceed via C\v{e}ch co-homology, {\it cf.} \cite[IV.3.5]{giraud}.
The principle complication here comes from the 
fact that in degrees $\geq 2$ the C\v{e}ch co-homology 
need not compute 
co-homology in an arbitrary site- indeed this can
already occur for the \'etale site of schemes
over an algebraically closed field
which aren't quasi-projective. Nevertheless, we 
always have
\begin{notation}\label{not:382} Let $U\ra\cX$ be an
\'etale atlas, and $\underline{Z}$ a sheaf of abelian
groups then $\rE_r^{p,q}(U,\underline{Z})$ ($\Rightarrow
\rH^{p+q}(\cX, \underline{Z})$) will be the C\v{e}ch
spectral sequence for the cover $U$, and the directed
limit (spectral sequence) of such over all open covers 
$U$ will be written
$\check{\rE}_r^{p,q} (\cX, \underline{Z})$ or, possibly,
in the special case that $q=0$, $r=2$, 
$\check{\rH}^p (\cX, \underline{Z})$.  
\end{notation}
In particular therefore we have an exact sequence
\begin{equation}\label{3123}
0\ra \check{\rH}^2 (\cX, \underline{Z})\ra
\rH^2(\cX, \uZ) \ra \check{\rE}_2^{1,1} (\cX, \underline{Z})
\xrightarrow{d^{1,1}_2}
\check{\rH}^3 (\cX,\uZ)
\end{equation}
whose manifestation in $\htg$ can be constructed rather
explicitly. As ever we choose a sufficiently fine atlas
$U\ra \cX$ so that both $\cX$ and the class of $q:\cY\ra\cX$
can be represented by groupoids $R\xrightarrow{q} R_0 \rras U$,
with $q$ a $\uZ$-torsor.
Consequently if we have a normalised 
({\it i.e.} vanishing on identities) 
C\v{e}ch co-cycle, $\z:R_0\ra s^* \uZ$,
(which is, in fact, isomorphic to the constant sheaf $Z$ for $U$
sufficiently fine)
then we can define another
groupoid $R'$ on changing the multiplication rule in $R$
by way of,
\begin{equation}\label{3124}
(g\cdot f)_{\mathrm{new}} \, =\, (g\cdot f\cdot \z(g,f))_{\mathrm{old}}
\end{equation}
which is associative because $\z$ is a co-cycle. On the
other hand by \eqref{eq:sep3} to give a map between $q=[U/R]$
and $q'=[U/R']$ in $\et_2(\cX)$ is, on refining $U$ sufficiently,
equivalent to giving a functor $F:R\ra R'$ along with a
natural transformation $q\Rightarrow q'F$. As a map
from $U$ to $R_0$, however, any natural transformation
can, on refining $U$ as necessary, be lifted to a map
from $U$ to $R'$, so without loss of generality $q=q'F$,
and the functor has the form $f\mpo f z(f)$ for $z(f)$
in the stalk $\uZ_{s(f)}$. In particular, therefore,
$\z$ must be the differential of $z$, and whence
\begin{fact}\label{fact:383}
The change of multiplication rule \eqref{3124} makes 
$\htg(\cX, \uZ)$ a  
$\check{\rH}^2(\cX, \uZ)$-set 
with faithful action, which applied to the trivial
class yields a commutative square
\begin{equation}\label{3125}
\begin{CD}
\check{\rH}^2(\cX, \uZ) @>{\eqref{3124}}>> \htg(\cX, \uZ) \\
@| @AA{\eqref{3122}}A \\
\check{\rH}^2(\cX, \uZ) @>{\eqref{3123}}>> \rH^2(\cX, \uZ)
\end{CD}
\end{equation}
\end{fact}
\begin{proof} We've done everything except the commutativity,
which by the naturality of \eqref{3122}-\eqref{3123} amounts
to checking that the composite of the bottom row with the 
right vertical is what one gets from the tautological
$\d$-functor \eqref{3121} applied to the C\v{e}ch complex,
which, by the way, is a slightly different line of
reasoning to the same assertion in \cite[IV.3.5.1.6]{giraud}.
\end{proof}
Maintaining the same notation, and 
bearing in mind that $s^* \uZ$ is
the constant sheaf $Z$, we equally have 
the intervention of hyper-coverings via the
composition:
\begin{equation}\label{3126}
\htg(\cX, \uZ) \ra \rH^1 (R_0, Z) \ra 
\check{\rE}^{1,1}_2:  q=[U/R] \ra  \bigl[ \text{class of $R/R_0$
as a $Z$-torsor}\bigr]
\end{equation}
whose kernel is the orbit of the trivial class
under \eqref{3124}. Better still, as in the proof
of \ref{fact:383} since $\rH_{\mathrm{Giraud}}$
is a $\d$-functor in degrees $\leq 2$, 
one can apply it directly to the C\v{e}ch complex
to get an exact sequence
\begin{equation}\label{3127}
0\ra \check{\rH}^2 (\cX, \underline{Z})\xrightarrow{\eqref{3124}}
\htg(\cX, \uZ) \xrightarrow{\eqref{3126}} 
\check{\rE}_2^{1,1} (\cX, \underline{Z})
\end{equation}
which commutes with \eqref{3123}- in the obvious way akin to \eqref{3125}-
by the naturality of \eqref{3122}. On the other hand \eqref{3122}
is always injective for general nonsense reasons, so it's an
isomorphism iff $d_2^{1,1}$ of \eqref{3123} vanishes on the
image of \eqref{3126}. On replacing $R$ in \ref{fact:pos1}
by that of \eqref{3126} this is formally the same as the proof
of \ref{fact:pos2}- 
replace $U/S$ in \eqref{eq:trans6} by a $Z$ torsor; $Y/X$ 
by a $Z\ts Z$ torsor and $\pi_1(X)$ by $Z\ts Z$ in 
\eqref{eq:trans7}-\eqref{eq:trans8}; and, of course, $P$ by a sufficiently
fine atlas $U$. In particular,
since one no longer has the lack 
of associativity in item.(c) of \ref{fact:pos1}
there is no transgression, and whence
\begin{fact}\label{fact:384}
The natural map \eqref{3122} is an isomorphism, and, better still,
the exact sequence \eqref{3123} (which is not a trivial isomorphism
between $\check{\rH}^2$ and $\rH^2$ iff hyper-coverings intervene
in a non-trivial way) is simply the extension of
\eqref{3127} on the right by $d^{1,1}_2$.
\end{fact}
This established we change notation accordingly, {\it i.e.}
\begin{newnot}\label{newnot:385}
We no longer distinguish $\htg$ from $\rH^2$ regardless
of whether it be in the
abelian or non-abelian case. The latter is defined
as follows: by $\uuG$ is to be understood a 
{\it locally constant (finite) link}, {\it i.e.} a
finite group
$\G$ together with a representation of $\pi_1(\cX_*)$ 
in $\Out(\G)$, 
\cite[IV.1.1.7.3]{giraud},
if this lifts to a representation in
$\Aut(\G)$ then we have a locally constant sheaf which 
will be written $\uG$, and in either case the centre
always defines a locally constant sheaf, $\uZ$, of
abelian groups. As ever, $\rH^1(\cX, \uG)$ is the set
of $\G$-torsors modulo isomorphism in $\et_1(\cX)$, 
while, following Giraud, $\rH^2(\cX, \uuG)$ is the
set of equivalence classes of $0$ cells $q:\cY\ra\cX$ in
$\et_2(\cX)$ of locally constant gerbes in $\rB_\G$'s
such that the locally constant sheaf $\cS_{\cY/\cX}$ 
(again, implicitly with an isomorphism of the stabiliser
of a base point with a fixed abelian group $\G$
to avoid the problem posed by \eqref{eq:groupTom})
affords the
outer representation $\uuG$.
\end{newnot}
These groups are easily described by the 2-Galois correspondence,
{\it i.e.}
\begin{fact}\label{fact:386} The $1$-category of (finite) locally constant 
sheaves is equivalent to the category of finite groups with
continuous
$\pi_1=\pi_1(\cX_*)$-action, and we have the Huerwicz isomorphism
\begin{equation}\label{3128}
\rH^1(\cX, \uG) \xrightarrow{\sim} \rH^1_{\mathrm{cts}}(\pi_1, \uG)
\end{equation}
Similarly the $2$-Category of locally constant gerbes 
(in finite groups) over $\cX$ is
equivalent to continuous actions of $\Pi_2(\cX_*)$ on groupoids with
one (up to isomorphism) object while $\rH^2(\cX, \uuG)$ is
formally ({\it i.e.} may be empty) a faithful
$\rH^2_{\mathrm{cts}}(\pi_1, \uZ)$-set and we have an exact sequence
of sets
\begin{equation}\label{3129}
0\ra \rH^2_{\mathrm{cts}}(\pi_1, \uZ)
\ra \rH^2(\cX,\uuG) \ra \Hom^{\pi_1}_{\mathrm{cts}}(\pi_2, \uZ)
\xrightarrow{A\mpo A_* K_3} 
\bigl( \rH^3_{\mathrm{cts}}(\pi_1, \uZ), \mathrm{obs} \bigr)
\end{equation}
where $\mathrm{obs}$ is the pull-back to $\pi_1$ of the
obstruction class of \ref{fact:group3}, so exactness
on the right is to be understood as image equals pre-image
of $\mathrm{obs}$ and elsewhere as the fibring of the
action over the orbits, so, in particular, 
$\rH^2(\cX, \uuG)\neq\emptyset$ iff the pre-image of the
obstruction class is non-empty. Better still, if $\G$ is
abelian then \eqref{3129} is an exact sequence of groups,
and $\rH^2(\cX, \uuG)$ is formally a principal homogeneous
space 
under $\rH^2(\cX, \uZ)$
via Giraud's concatenated product,
\begin{equation}\label{3130}
\rH^2(\cX, \uZ)\ts \rH^2 (\cX, \uuG) \ra \rH^2 (\cX, \uuG):
q\ts q'\mpo q\wedge^{\uZ} q'
\end{equation}
\end{fact}
\begin{proof} \eqref{3128} is an immediate tautology 
of 1-Galois theory, \cite[Expos\'e XI.5]{sga1}, while
\eqref{3129} is a similar tautology for 2-Galois
theory profiting from the explicit description, \ref{sum:group1}, of
actions of 2-groups on groupoids with 1 object. It
then follows from \eqref{3129} that if $\rH^2(\cX, \uuG)\neq\emptyset$
then it is a principal homogeneous space under $\rH^2(\cX, \uZ)$
as soon as the action \eqref{3130} is faithful. 
As such suppose there are  $q\in \rH^2(\cX, \uZ)$,
and $q'\in\rH^2(\cX, \uuG)$ such that
$q\wedge^{\uZ} q' =  q'$.
Now 
the restriction of the stabiliser to any sufficiently
fine open cover is always the constant sheaf, so
identifying $q$, $q'$, {\it etc.} with groupoids
$R$, $R'$ we trivially have the non-abelian 
variant of \eqref{3126}, {\it i.e.}
\begin{equation}\label{3131}
\rH^2(\cX, \uuG) \ra \rH^1 (R_0, \G): q'\ra 
\bigl[\text{class of $R'/R_0$ as a $\G$-torsor}\bigr]
\end{equation}
which in turn factors through the abelianisation of
the fundamental group of every component of $R_0$
by \eqref{3114}. Consequently by \eqref{3128} applied
to $R_0$, the image of $q$ under \eqref{3126} is
trivial, so it's a C\v{e}ch class. Modulo notation,
the action of C\v{e}ch classes is given by \eqref{3124},
and it's always non-trivial if $q$ is non-trivial for the same reason as
discussed post {\it op. cit.}.
\end{proof}
Plainly \eqref{3129} is crying out to be dualised
and/or take values in linearly topologised sheaves,
so let us make a brief
\begin{scholion}\label{schol:386}{\it Continuous \'etale homology
and continuous Giraud co-homology}
Suppose for the sake of argument we simply define
\begin{equation}\label{3132}
\begin{split}
& \rH_q^\cts (\cX, \hat{\bz}):= \Hom_\cts (\rH^q(\cX, \bq/\bz), \bq/\bz),
\\
& 
\rH_q^\cts (\pi_1, \hat{\bz}):= \Hom_\cts (\rH^q_\cts(\pi_1, \bq/\bz), \bq/\bz)
\end{split}
\end{equation}
\end{scholion}
where even though the co-homology groups being dualised
on the right are discrete, we write $\Hom_\cts$ to emphasise
that the dual is to be understood as a pro-finite group,
{\it i.e.}  the topological dual of the direct limit over
finite sub-groups.
As such, \eqref{3129}, affords an exact sequence of
topological groups 
\begin{equation}\label{3133}
0\leftarrow \rH^\cts_2 (\pi_1, \hbz)
\leftarrow \rH^\cts_2 (\cX,\hbz)
\leftarrow \pi^\cts_2(\cX)_{\pi_1}
\leftarrow \rH^\cts_3 (\pi_1, \hbz)
\end{equation}
Of which the easiest term to understand is the 
topological co-invariants
of $\pi_2$, 
{\it i.e.} the maximal invariant topological quotient, or
equivalently the quotient by the smallest closed sub-module
containing (the closed set) $S^\om-S$, for $S\in \pi_2$, and
$\om\in\pi_1$. As such the topological co-invariants are
potentially smaller than the co-invariants, albeit everything
is profinite so by \cite[Th\'eor\`eme 7.1]{jensen}
topological co-invariants have the right behaviour, {\it i.e.}
\begin{equation}\label{TopCo}
\pi^\cts_2(\cX)_{\pi_1} =\varprojlim_i \pi_2^i(\cX)_{\pi_1}  
\end{equation}
for $i$ ranging over the finite quotients of $\pi_2(\cX)$.
Similarly the {\it continuous group homology} term
will invariably differ from the usual group homology.
Specifically,
one computes continuous group co-homology 
with values in a pro-finite
$\pi_1$-module $A$
by applying the functor $\Hom^\cts_\Ens(-,A)$
to the
continuous complex of sets
\begin{equation}\label{3134}
\rp  \leftleftarrows \pi_1 
{\build\leftarrow_{\leftarrow}^{\leftarrow}} \pi_1\ts\pi_1\cdots
\end{equation}
On the other hand, the notion of group ring,
a.k.a. adjoint to the forgetful functor to $\Ens$,
continues to have perfect sense, {\it i.e.}
\begin{claim}\label{claim:388} Let $F=\varprojlim_{i\in I} F_i$ be a 
compact separated totally disconnected space in $\Top$, {\it i.e.}
$F_i$ finite with the limit taken in $\Top$, or, 
respectively and 
more generally 
a prodiscrete set $F_{i\in I}$ in $\pro-\Ens$,
then there is a (uniquely unique) pair
consisting of a profinite
group $C(F)$, respectively a protorsion
group $C_{j\in J}$,
and a continuous map $c:F\ra C(F)$, respectively a 
map in $\pro-\Ens$,
such  that if 
$f:F\ra A=\varprojlim_\a A_\a$ is any other continuous 
map to a pro-torsion abelian group,
respectively a map to a pro-group $A_{\a\in\aleph}$,
 then there is a unique
continuous group homomorphism $\phi_f:C(F)\ra A$, 
respectively pro-homomorphism,
such
that $f=\phi_f c$.
\end{claim}
\begin{proof} Partially order the set $I\ts\bn$ by
$(i,m) > (j,n)$ iff $i>j$ and $m>n$,
which is right co-filtering, respectively co-filtering, if $I$ is, 
then
\begin{equation}\label{3135}
C(F) := \varprojlim_{(i,n)} \bz/n\otimes_\bz \bz[F_i] 
\end{equation}
for $\bz[-]$ the free abelian group on a set, 
respectively the implied pro-object,
does
the job.
\end{proof}
Applying this to the complex \eqref{3134} yields 
a complex of pro-finite abelian groups
 \begin{equation}\label{3136}
\hbz=C_0:=C(\rp) \leftleftarrows 
C_1:=C(\pi_1)
{\build\leftarrow_{\leftarrow}^{\leftarrow}} 
C_2:=C(\pi_1\ts\pi_1)\cdots
\end{equation}
which for a finite group is just $\hbz\otimes_\bz-$
applied to the standard homology complex, but,
otherwise may be a lot more complicated. Irrespectively,
the complex \eqref{3136} shares the following pleasing
features with the standard homology complex

(a) If $A$ is a pro-torsion module,
there is an identity of complexes
\begin{equation}\label{3137}
\Hom^\cts_\Ens(\pi^\bullet, A) = \Hom^\cts_{\mathrm{Ab}}(C_\bullet, A)
\end{equation}
which is the right widget for 
continuous group co-homology if $A$
is discrete or pro-finite.

(b) In the particular case of $A=\bq/\bz$ Pontryagin duality implies
that $\rH^\cts_q(\pi_1, \hbz)$ is the homology of the
complex \eqref{3136}.

Now there are 2, and to some extent 3,  ways
in which this rather
satisfactory discussion starts to break down.
In the first place, one cannot find a compact
complex of compact totally disconnected spaces
such as \eqref{3134}, {\it e.g.} already 
for $\pi_1$ a discrete group this is an issue,
and pro-discrete is much worse, albeit this
isn't really an algebraic problem since under
very reasonable hypothesis 
({\it e.g.} Noetherian in characteristic zero, champs whose moduli
are quasi-projective over affine- use the descent
spectral sequence
of \cite[5.3.5]{deligne3}
as in \ref{cor:l12} or \ref{cor:SmoothBaseChange}
to reduce to corresponding statement for schemes)
the inverse limit
\begin{equation}\label{3138}
F_0(\cX)=\varprojlim_U \pi_0(U)  \leftleftarrows 
F_1(\cX)= \varprojlim_U \pi_0(U\ts_\cX U)
{\build\leftarrow_{\leftarrow}^{\leftarrow}} 
F_2(\cX)= \varprojlim_U\pi_0(U\ts_\cX U\ts_\cX U)
\end{equation}
over finite coverings does the job. Nevertheless it can
fail, and it isn't so easily remedied in pro-Top rather
than the more complicated pro-homotopy category, and in
either case, one looses the very useful freedom of
working with right co-filtered rather than co-filtered partially
ordered sets
that one has in the 
pro-finite case, equivalently compactness allows one to pass
from right co-filtered to co-filtered. In a similar vein
is the issue of hypercoverings, but since we're only
concerned with gerbes and torsors it's intervention,
\eqref{3126}, 
is very limited, 
and we could just change \eqref{3138}
to the complex of pro-finite sets associated
to the obvious double complex of such suggested by
\eqref{3127} {\it et seq.} for
the general Noetherian algebraic case. Consequently, the 
more serious way in which the discussion begins 
to break down is when one tries to generalise
from pro-finite groups to pro-group. In order to see
what's involved let's put ourselves in a
good situation such as 
pro-finite group cohomology or
algebraic varieties where \eqref{3138}
is the right widget and use \ref{claim:388}
to define
the ``continuous C\v{e}ch homology complex''
\begin{equation}\label{3139}
C_0(\cX):=C(F_0(\cX))\leftleftarrows 
C_1(\cX):=C(F_1(\cX))
{\build\leftarrow_{\leftarrow}^{\leftarrow}} 
C_2(\cX):=C(F_2(\cX))
\cdots
\end{equation}
From the simple expedient of applying the definitions,
we have for $A=\varprojlim_\a A_\a$
a locally constant sheaf of pro-finite abelian groups,
\begin{equation}\label{3140}
\check{C}^q_\cts (\cX, A):=
\Hom^\cts_{\mathrm{Ab}} (C_q (\cX), A)=
\Hom^\cts_\Ens(F_q, A)
=\varprojlim_\a \varinjlim_U
\rH^0(U^{(q+1)}, A_\a)
\end{equation}
and this sort of behaviour can be replicated in
pro-Sheaf, but what can't be replicated is
\begin{claim}\label{claim:ProExtra1} If
\eqref{3138} is a complex of compact totally
disconnected spaces, and
$A=\varprojlim_\a A_\a$
is a sheaf of locally constant profinite groups then
\begin{equation}\label{3144}
{\varprojlim_\a}^{(p)} \varinjlim_U
\rH^0(U^{(q+1)}, A_\a) =0, \forall p>0
\end{equation}
\end{claim}
\begin{proof}
The good description is the
final term in \eqref{3140}, since as far as the 
right co-filtered and partially ordered
set ${\a}$ is concerned the sheaves
$\Hom^\cts_\Ens(F_\bullet, A)$ are
weakly flasque in the sense of \cite[Th\'eor\`eme 1.9]{jensen}
because a topological surjection of pro-finite
groups has a section- \cite[I.Proposition 1]{Serre}.
\end{proof}
Consequently,
under such hypothesis the co-homology,
$\check{\rH}^\bullet_\cts$,
of the complex \eqref{3140} exhibits exactly the
same good features of continuous group
co-homology, {\it e.g.} turns short
exact sequences of locally constant
pro-finite sheaves into long exact sequences,
and is the abutment of a spectral sequence
\begin{equation}\label{31410} 
{\varprojlim_\a} ^{(p)} \check{\rH}^{q} (\cX, A_\a)
\Rightarrow \check{\rH}_\cts^{p+q}(\cX, A)
\end{equation} 
and there are two distinct points to note, {\it i.e.}
\begin{rmk}\label{rmk:Homology} The relation between homology
and co-homology is a tautology, which is best evidenced by
a systematic use of Pontryagin duality, rather than the 
so called `universal coefficient theorem', which in the
above situation takes the form of a spectral sequence
\begin{equation}\label{3142}
\mathrm{Ext}^p_\cts (\check{\rH}^\cts_{-q}(\cX, \hbz), A)
\Rightarrow \check{\rH}_\cts^{p-q}(\cX, A)
\end{equation}
for 
$\check{\rH}^\cts_{\bullet}$ the homology of the complex
\eqref{3139} and
$A$ either a pro-finite group or indeed a discrete
discrete group, where in the latter case $\rH_\cts^\bullet$
just reduces to ordinary co-homology. 
In particular, therefore  $\check{\rH}^\cts_q(\cX, \hbz)$
is the topological dual of
the discrete group $\check{\rH}^q(\cX, \bq/\bz)$, as
proposed in the definition
\eqref{3132}. As such, while it's true that in absolute 
generality, \cite[\S 2]{artin-mazur}, one cannot make such
a simplistic definition, and the `homology complex' will
only ever be unique up to homotopy, all sensible definitions
yield \eqref{3132} tautologically, {\it i.e.} duality between
homology and co-homology should not be confused with an
actual theorem, {\it viz:} duality between co-homology and
co-homology with supports.
\end{rmk} 
\begin{rmkdef}\label{defn:3810} Consequently homology is
nothing more than notation, whereas it's preferable to have
a definition of $\rH^\bullet_\cts$ which is beyond ambiguity.
This can be done as follows: in Pro-sheaf (or perhaps better
$\pro^{\leq\aleph}$, {\it i.e.} indexing sets of cardinality
as most some fixed accessible ordinal) a short exact sequence
\begin{equation}\label{MorePro1}
0\ra \hat{\cF}'\ra \hat{\cF} \ra \hat{\cF}''\ra 0
\end{equation}
can be identified with a net of exact sequences of sheaves
of abelian groups
\begin{equation}\label{MorePro2} 
0\ra \cF'_i\ra \cF_i \ra \cF''_i\ra 0 
\end{equation}
with implied transition functions
for the same filtering $I$,
\cite[A.3.2-4.6]{artin-mazur}. Thus in the first place
there are enough injectives- same reasoning as \cite[\S 1]{jensen}- 
and in the second place the derived functors of
\begin{equation}\label{ProExtra1}
\rH^0: \mathrm{proSheaf}\ra \mathrm{proGroup}:
\hat{\cF}:=\cF_{i\in I}\ra \rH^0(\cX,\cF_i)_{i\in I}
\end{equation}
are just the pro-groups $\rH^q(\cX,\cF_i)_{i\in I}$.
while those of 
\begin{equation}\label{ProExtra2}
\mathrm{proGroup}\ra\mathrm{Group}:A_{i\in I}\ra\varprojlim_i A_i
\end{equation}
are ${\varprojlim_i}^{(p)}A_i$. From which, if we define
$\rH^n_\cts(\cX, \hat{\cF})$
to be the derived functors
of $\rH^0_\cts(\cX, \hat{\cF}):=\varprojlim_i\rH^0(\cX,\cF_i)$ then we get
a spectral sequence
\begin{equation}\label{3141}
{\varprojlim_\a} ^{(p)} {\rH}^{q} (\cX, A_\a)
\Rightarrow {\rH}_\cts^{p+q}(\cX, A)
\end{equation} 
So by \eqref{31410} this more generally valid
definition coincides with the previously
envisioned good cases.
\end{rmkdef}
Concentrating, therefore, 
on continuous co-homology
we can, and should,
extend the C\v{e}ch 
discussion
to a non-abelian locally 
constant sheaf $\uG=\varprojlim_i \uG_i$ of 
pro-finite groups. Indeed, $\rH^1$ is always
C\v{e}ch, and we just define 
$\rH^1_\cts (\cX, \uG)$,
or for that matter $\rH^1_\cts(\pi_1, \uG)$,
which is the same,  
to be continuous C\v{e}ch 1 co-cycles
modulo continuous C\v{e}ch boundaries
exactly as in \eqref{3140}, albeit via
the final identity over all possible
covers should there be need for covers
with infinitely many connected components.
This said, supposing further that the hypothesis
of \ref{claim:ProExtra1} are valid, it's equally
true in the non-abelian case, whence
the spectral sequence \eqref{31410} holds
as stated, up to an appropriate change in the
definition of the $d^{0,1}_2$ 
\`a la \eqref{l14}. In particular, therefore, since
the spectral sequences are the same it's
a formal consequence 
(in the presence of a sane definition of
$\pi_0(\cX)$, {\it cf.} \ref{fact:392})
that the Huerwicz
isomorphism, \eqref{3128}, self generalises to
\begin{equation}\label{3145}
\rH^1_\cts(\cX, \uG) \xrightarrow{\sim} 
\varprojlim_i \rH^1(\cX, \uG_i) \xrightarrow{\sim}
\varprojlim_i \rH^1_{\mathrm{cts}}(\pi_1, \uG_i)
\xleftarrow{\sim}\rH^1_\cts (\pi_1, \uG)
\end{equation}
for arbitrary pro-finite $\G$ thanks to the
ubiquitous \cite[Th\'eor\`eme 7.1]{jensen}.
Consequently, the continuous co-homology classes
in degree 1
are 
represented by
nets $q_i$ where $q_i:\cE_i\ra\cX$ is a $\G_i$-torsor
along with 1-cells 
$F_{ji}:q_i\ra q_j$, $i\geq j$, which are equivalences
on quotienting $\cE_i$ by the kernel, 
$\uG_i^j$,
of $\uG_i\ra\uG_j$
such that there is a 2 commutative diagram in $\et_2(\cX)$ 
\begin{equation}\label{3146}
 \xy
 (-18,0)*+{q_i}="L";
 (18,0)*+{q_k}="R";
 (0,16)*+{q_j}="T";
    {\ar^{F_{ji}} "L";"T"};
    {\ar^{F_{kj}} "T";"R"};
    {\ar_{F_{ki}} "L";"R"};
    {\ar@{=>}^{\xi_{kji}} (0,2);(0,12)}
 \endxy
\end{equation}
whenever $i\geq j\geq k$, or, equivalently a strictly
commutative diagram in $\et_1(\cX)$, a.k.a. drop the
$\xi_{kji}$. Now while we supposed $\G$ pro-finite to
arrive at this conclusion, it is in fact always valid,
{\it i.e.} 
\begin{rmk}\label{fact:ProExtra1} 
Even for a general pro-sheaf of groups, $\hat{\uG}=\uG_{i\in I}$-
say locally constant to fix ideas, but this isn't important- 
commutative 
diagrams in $\et_1(\cX)$ of torsors of the form \eqref{3146}
modulo equivalence, {\it i.e.} commuting squares
\begin{equation}\label{Pro3150}
\begin{CD}
q'_i  @>>{B_i}> q_i \\
@V{F'_{ji}}VV @VV{F_{ji}}V \\
q'_j @>>{B_j}> q_j
\end{CD}
\end{equation}
where the $A_i$ are isomorphisms of torsors represent 
exactly
co-homology classes in $\rH^1_\cts(\cX, \hat{\uG})$
defined, for example, via \eqref{ProExtra1}-\eqref{3141} in
the abelian case, and otherwise by anything satisfying
the non-abelian spectral sequence, {\it cf.} \eqref{l14}.
\end{rmk}
\begin{proof} By the spectral sequence \eqref{3141},
or it's non-abelian variant, it suffices to show that
if all the $q_i$, $q'_i$ are trivial in \eqref{Pro3150},
then the resulting diagrams modulo isomorphisms of the
trivial torsor ({\it i.e.} $B_i$ multiplication by
$\g_i\in\rH^0(\cX, \G_i)$ and $F_{ji}$ the projections, a.k.a.
$\varprojlim_i \rH^0(\cX, \G_i)$ ) are 
${\varprojlim}^{(1)} \rH^0(\cX, \G_i)$. On the one hand,
given the diagram there is a unique global section $\g_{ji}$
of $\G_j$ such that $F_{ji}B_i$ is $B_jF'_{ji}$ followed
by multiplication by $\g_{ji}$, {\it i.e.} a class in
the C\v{e}ch description of $\lim^{(1)}$, while, conversely
given such a class we can take the right   hand
vertical in \eqref{Pro3150} to be the trivial projection
followed by multiplication by it, and everything else
trivial.
\end{proof}
Similarly
restricting ourselves momentarily to the profinite abelian
case, the spectral sequence \eqref{3141} implies-
\ref{claim:noplag1} for the $F^1\rH^2_\cts$ term-
that classes in $\rH^2_\cts$ 
are themselves
2 commutative triangles of the form \eqref{3146},
but now with
$q_i:\cE_i\ra\cX$ a locally constant gerbe
for which $S_{\cE_i/\cX}\xrightarrow{\sim} \uG_i$,
where
again
$F_{ji}$ are equivalences modulo $\G_i^j$, and
these triangles themselves 
form a 2-commutative tetrahedron
\begin{equation}\label{3148}
 \xy
 (20,-24)*+{q_i}="A";
 (38,6)*+{q_j}="B";
 (0,0)*+{q_k}="C";
(15,12)*+{q_l}="D";
{\ar_{}_{F_{ji} } "A";"B"};
    {\ar^{F_{kj} } "B";"C"};
    {\ar^{F_{ki} } "A";"C"};
    {\ar^{F_{lk} } "C";"D"};
    {\ar_{F_{lj} } "B";"D"};
{\ar@{-->}_{F_{li}} "A";"D"}
\endxy
\end{equation}
as already encountered in the condition
\ref{def:PosSeq}.(e) of Postnikov sequences, which,
{\it cf.} \ref{5.11}, is always true if $q_i$ is a torsor rather
than a gerbe. Again, however, the situation is
more general than pro-finite, {\it i.e.}
\begin{rmk}\label{fact:ProExtra2}
Even for a general pro-sheaf of groups, $\hat{\uG}=\uG_{i\in I}$,
say abelian for the moment triangles, \eqref{3146}, 
of $\uG_i$ gerbes forming a 2-commutative diagram,
\eqref{3148} modulo equivalences, {\it i.e.}
2-commutative diagrams \eqref{Pro3150} with
$\b_{ji}: B_jF'_{ji}\Rightarrow F_{ji}B_i$ such
that
\begin{equation}\label{Pro3152}
 \xy
 (38,0)*+{q'_j}="B";
 (0,-6)*+{q'_k}="C";
(15,0)*+{q'_i}="D";
 (38,-12)*+{q_j}="F";
 (0,-18)*+{q_k}="G";
(15,-12)*+{q_i}="H";
{\ar_{}_{B_j} "B";"F"};
    {\ar_{B_k} "C";"G"};
    {\ar_{B_i} "D";"H"};
{\ar_{}^{F_{kj} } "F";"G"};
    {\ar@{-->}_{F_{ki} } "H";"G"};
    {\ar^{F_{ji} } "H";"F"};
{\ar_{}^{F'_{kj}} "B";"C"};
    {\ar_{F'_{ki}} "D";"C"};
    {\ar^{F'_{ji}} "D";"B"};
\endxy
\end{equation}
2-commutes represent exactly the classes in $\rH^2_\cts(\cX,\hat{\uG})$.
\end{rmk}
\begin{proof} Again, from what we've already said it suffices
to look at triangles $(\cE_i, F_{ji}, \xi_{kji})$ with
the first 2 entries trivial, so $\xi_{kji}$ is the
C\v{e}ch co-cycle in 
${\varprojlim_i}^{2}\rH^0(\cX, \G_i)$,
and conversely, of which the remaining term, 
$F^2\rH^2_\cts(\cX,\hat{\uG})$, in 
the spectral sequence, \eqref{3141}, is a 
quotient by ${\varprojlim_i}\rH^1(\cX, \G_i)$.
As such the right interpretation of
this latter term is \ref{claim:noplag1},
{\it i.e.} its exactly the effect of changing the $B_i$'s. 
\end{proof}

Manifestly this description of $\rH^2_\cts$
in terms of gerbes is independent of whether the $\G_i$
are abelian or not and whence \ref{newnot:385} extends
continuously to
\begin{rmkdef}\label{def:3811} Denote by $\uuG=\varprojlim_i \uuG_i$
a 
{\it locally constant profinite link}, {\it i.e.} a
profinite group, $\G$, together with a continuous 
representation of $\pi_1(\cX_*)$ in $\Out(\G)$ using
the equivalence of categories \eqref{1372}
in the presence of
\ref{def:371}.(a) 
and \cite[A.3.3]{artin-mazur}
if necessary 
to replace $I$ by a co-final subset, 
or, more generally, a {\it pro-link}, {\it i.e.}
a pro-object in the category of links, \cite[IV.1.1.6]{giraud},
over $\cX$ albeit the most complicated example we
have in mind is a representation of a discrete
group in the pro-outer automorphisms of a pro-group, 
we define
$\rH^2_\cts(\uuG)$ 
exactly as in \ref{fact:ProExtra2}.
In particular
for $\uZ$ 
(locally constant sheaf of pro-finite abelian groups)
the centre of a profinite link $\uuG$, 
\eqref{3129} and \eqref{3130} hold (same proof) 
on replacing the co-homology groups $\rH^2_\cts(\cX, \uuG)$ 
and $\rH^2_\cts(\cX, \uZ)$ with values in finite groups
by their pro-finite analogues
$\rH^2_\cts(\cX, \uuG)$ and $\rH^2_\cts(\cX, \uZ)$.
\end{rmkdef}
The highly conceptual nature of Giraud co-homology
equally implies that the essentially tautological
connecting morphisms extend to the pro setting.
Specifically say,
\begin{equation}\label{3149}
0\ra\uG'\xrightarrow{a}\uG\xrightarrow{b}\uG''\ra 0
\end{equation}
a topologically exact sequence of sheaves for
pro-finite groups, or pro-exact otherwise,
so in either case we can describe the 
sequence 
as a net of short exact sequences
\begin{equation}\label{3150} 
0\ra\uG'_i\xrightarrow{a_i}\uG_i\xrightarrow{b_i}\uG''_i\ra 0
\end{equation}
over the 
same right co-filtered
set, respectively co-filtered partially ordered set, $I$-
in the former case because
$a$ is closed, and
$b$ has the quotient topology, 
and 
\eqref{MorePro2} for the latter.
Consequently if $\hat{\g}''=(\g''_i)$ is a global section
of $\uG''$ then,
\begin{equation}\label{3151}
U\mpo E_i(U):=\{ \g_i\in \G_i(U)\,\,\vert\,\, \g_i\mpo \g''_i\mid_U\}
\end{equation}
defines a net of $\G'_i$-torsors satisfying \eqref{3146},
whose triviality in the sense of \ref{fact:ProExtra1} is
necessary and sufficient to lift $\g''$ to a global section
of $\uG$.
Similarly, if $\hat{q}''=(q''_i, F''_{ji},
\xi''_{kji})$ is
a net of $\G''_i$-torsors satisfying the continuity 
condition \eqref{3146}, then for each $i\in I$ we can
define a gerbe $\cE(q''_i)$ whose objects over an open
set $U$ are pairs $(q_i, Q_i)$ consisting of a
$\uG_i\vert_U$ torsor $q_i:E_i\ra\cX$ and a 1-cell
$Q_i:q_i\ra q''_i$ which is an equivalence modulo $\G'_i$,
while arrows are maps of $\uG_i\vert_U$-torsors satisfying
the obvious commutativity conditions. As such there are
functors and natural transformation
\begin{equation}\label{3152}
F_{ji}: (q_i, Q_i) \mpo (q_i\, \mathrm{mod}\, \G_i^j, F''_{ji} Q_i),
\quad 
\xi_{kji}:F_{ki}\Rightarrow F_{kj}F_{ji}: (q_i, Q_i)\mpo Q_i^*\xi''_{kji}
\end{equation} 
satisfying the tetrahedron condition \eqref{3148} 
because, as we've said, it's always true for torsors
and each $\cE(q''_i)$ is a (fibred) category of torsors.  
By construction, therefore, the pro-torsor
$\hat{q}''$ can be lifted to a pro-torsor
$\hat{q}=(q_i, G_{ji}, \eta_{kji})$ iff the 
pro-gerbe $\d_2(\hat{q}''):=(\cE(q''_i), F_{ji}, \xi_{kji})$ 
is trivial in the sense of
\ref{fact:ProExtra2}. Nevertheless, the usual non-abelian
subtlety merits
\begin{warning}\label{GiraudWarn} Already for $i$ fixed, it is by no means 
the case that the
stabilisers $S_{\cE(q''_i)/\cX}$ are isomorphic to
$\uG'_i$, but only to some locally constant sheaf
whose stalks are isomorphic to the same finite group
$\G'_i$. Nor need it even be true that $S_{\cE(q''_i)/\cX}$
defines a sheaf rather than a link $\uuG'_i$ on $\cX$.
Consequently the conditions for the annihilation
of the connecting homomorphism $\d_2(\hat{q}'')$
are twofold: the continuous link $\uuG=\varprojlim_i\uuG'_i$ lifts to a 
(necessarily locally isomorphic to $\uG'$) continuous
sheaf, whence $\rH^2_\cts (\cX, \uuG)$ has a distinguished
trivial class and $\d_2(\hat{q}'')$ must be equal to it.
As such, 
irrespective of continuity,
the only general hypothesis
where one can write 
\begin{equation}\label{3153}
\rH^1_\cts (\cX, \uG') \xrightarrow{\d_2} \rH^2_\cts (\cX, \uG'')
\end{equation}
is if $\uG'$ is central in $\uG$.
\end{warning} 
\subsection{Pro-discrete homotopy}\label{SS:III.9}
To begin with let $\cX$ be a topological champ-
so, no hypothesis beyond $\cX$ being the classifying
champ $[U/R]$ for $(s,t):R\rras U$ an \'etale groupoid
in separated spaces. Now the set of sub-champs
$\cV\hookrightarrow\cX$ which are open and closed,
define subsets $A$, respectively
$B$, of the set of open and closed subsets of
$U$, respectively $R$.  On the other hand if $a\in A$,
it's complement is in $A$, so $x\sim y$ iff
$\forall a,\, x\in a\Rightarrow y\in a$, is an
equivalence relation on $U$, and similarly on $R$.
Better, for all $a\in A$,  
$s^{-1}a=t^{-1}a\in B$, and every element of $B$
is of this form, so $R\xrightarrow{s\ts t} U\ts U$ factors
through $\sim$ and we make
\begin{rmkdef}\label{def:391} The topological quotient,
$U/\sim$ is the maximal totally separated ({\it i.e.}
open and closed sets separate points) of $\cX$ and will
be written $\pi^\cts_0 (\cX)$, so, for example if 
$F=\varprojlim_i F_i$ is a pro-finite set then $\pi_0^\cts(F)$
is the space $F$ itself, while for a 
non-compact
(separated) totally disconnected space 
$F\ra\pi^\cts_0 (F)$ need not even be a set isomorphism.
\end{rmkdef}

Alternatively, one can do a similar thing in
pro-topological spaces, {\it i.e.}

(a) Any (continuous) surjective map $\cX\xrightarrow{a} A$ to a discrete
space is a set quotient, so it's also a topological quotient
since no topology is finer than the discrete topology.

(b) Any  map to a discrete space factors through a unique
surjective map, a.k.a. the image which is also the image
of $U$.

(c) Such quotients form a co-filtered directed set. Indeed if
$\cX\xrightarrow{a} A$, $\cX\xrightarrow{b} B$,
then the image $\cX\xrightarrow{c} C$ of $a\ts b$
does the job.

(d) Any such quotient is a fortiori a quotient of $U$, so
up to isomorphism there are only a sets worth, $\aleph$, of
such quotients.

So as a result we we can make

\begin{rmkdef}\label{def:392} Define $\pi^\pro_0(\cX)$ to be 
the 
universal map from $\cX$ to a pro-discrete space, {\it i.e.}
the pro-topological space $\hat{A}:=A_{{a\in\aleph}}$,
equivalently the 
opposite of the
functor $\varinjlim_{a\in\aleph}\Hom(A,-):
\mathrm{Top}\ra\Ens$. 
Thus, $\pi_0^\cts\ra\pi_0^\pro$, while, for example,
if 
$\hat{F}={F}_{i\in I}$ is a pro-object in discrete sets;
$F=\varprojlim_i F_i$ the limit in $\mathrm{Top}$; and
$\Lambda$ a discrete space
\begin{equation}\label{3154}
\Hom_\pro(\hat{F}, \Lambda)\ra \Hom_\cts (F, \Lambda)
\end{equation}
with equality if $F$ is pro-finite, or $\Lambda$ is finite.
\end{rmkdef} 
In a similar vein \ref{def:391} and \ref{def:392} 
do the same thing in slightly different ways
\begin{fact}\label{fact:392}
(Huerwicz-0) 
Let $\Lambda$ be
a discrete space then
\begin{equation}\label{31154}
\rH^0(\cX,\Lambda):=
\Hom_\cts (\cX, \Lambda)=
\Hom_\cts (\pi^\cts_0(\cX), \Lambda)=
\Hom_\pro (\pi^\pro_0(\cX), \Lambda)
\end{equation}
and 
$\cX$ is connected iff both $\pi_0^\cts(\cX)$,
and $\pi_0^{\pro}(\cX)$ are reduced to points.
Furthermore,
the following are equivalent

(a) $\cX$ admits a universal map $p_0:\cX\ra X_0$ to the
category of discrete spaces.

(b) $\pi_0^\cts (\cX)$, respectively $\pi_0^\pro (\cX)$,
is, respectively is representable by, a discrete space.

(c) $\cX$ is continuously, resp. pro, 
semi-locally $0$-connected, {\it i.e.} every
point has an \'etale neighbourhood $V$ such
that the image of $\pi_0^\cts (V) \ra \pi_0^\cts (\cX)$,
resp. $\pi_0^\pro (V) \ra \pi_0^\pro (\cX)$,
is a point.
\end{fact}
\begin{proof} Everything prior to the 
equivalence of (a)-(c) is just the
definitions, as is (b)$\Rightarrow$(a). For (a)$\Rightarrow$(c),
and $x$ a point take $\cV$ to be $p_0^{-1}(p_0(x))$, so
any \'etale neighbourhood $V\ra\cV$ of $x$ does the
job. Finally (c) for $\pi^\cts_0$ implies that the
equivalence classes in \ref{def:391} are open, {\it i.e.}
$\pi^\cts_0$ is discrete, while in the respective case
the image of $\pi^\pro_0$ trivial implies the same for
$\pi^\cts_0$.
\end{proof}
Now while \ref{def:391} 
is wholly in $\Top$, and whence
appears
more satisfactory than \ref{def:392},
a priori it is the
latter rather than the former which
has the good Galois analogue, {\it i.e.}
\begin{rmkdef}\label{def:394} 
If $\pi_1^{\pro}=\pi^i_{1,i\in I}$ is a pro-discrete
group, then we'll write $\pi^\cts_1$ for its limit in 
$\Top$. Plainly the topology of the latter has a basis which are
translations of open and closed normal sub-groups, nevertheless
the map $\pi^\cts_1\ra \pi_1^\pro$ has much worse behaviour
than the map $\pi^\cts_0(\cX) \ra \pi_0^\pro(\cX)$, indeed
all transitions $\pi^i_1\ra \pi^j_1$ can be non-trivial surjections,
yet $\pi^\cts_1=\mathbf{1}$.
\end{rmkdef}
This said we can extend \ref{fact:392} to the
next level via
\begin{fact}\label{fact:395} (Huerwicz-1) Suppose $\cX$ is 
connected and
everywhere
locally connected ({\it i.e.} every point has a co-final
system of connected neighbourhoods) then for $*:\rp\ra \cX$
any point and $I_*$ the co-filtered directed class of
\eqref{eq:funct1406}-\eqref{eq:funct1106}
restricted to Galois objects, the category
$I_*$ is equivalent to a small category, and
there is a pro-discrete group $\pi_1^{\pro}(\cX_*)
=\pi^i_{1,i\in I_*} $ such
that
the fibre functor \eqref{eq:funct1206} affords
an equivalence of categories between
$\et_1(\cX)$ and
pro-discrete sets with $\pi_1^{\pro}(\cX_*)$ action
(which may seem a convoluted statement but it implies
a non-trivial relation between base points),
so, in particular a locally constant sheaf $\uG$
of 
discrete groups is
equivalent to a 
discrete group $\G$ with $\pi_1^{\pro}(\cX)$ action,
while
\begin{equation}\label{31156}
\rH^1 (\cX, \uG)\xrightarrow{\sim}
\varinjlim_i \rH^1 (\pi^i_1, \G)
\end{equation}
and the following are equivalent

(a) $\cX$ admits a universal 
representable \'etale covering $p_1:\cX_1\ra\cX$.

(b) $\pi_1^\pro(\cX_*)$ is representable by a discrete group.

(c) $\cX$ is semi-locally 1-connected, {\it i.e.} every
point, $x:\rp\ra\cX$, 
has an \'etale neighbourhood such that the image
of $\pi^\pro_1(V_x)\ra\pi^\pro_1(\cX_x)$ is trivial.

(d) Every point, $x:\rp\ra\cX$,
has an \'etale neighbourhood, $V$, such that any
\'etale cover of $\cX$ can be trivialised over $V$.
\end{fact}
\begin{proof}
By \ref{fact:plag1} every representable $0$-cell
$q':\cY'\ra\cX$ in $\et_2(\cX)$ is covered by a
representable Galois cell, {\it i.e.} a connected
torsor $q:\cY\ra\cX$ under a discrete group, say,
$\G$ at the risk of minor notational confusion. As
such everything prior to the equivalence of (a)-(d)
follows exactly as in the pro-finite case provided
we know that the Galois coverings form a set. Now
to give a $\G$ torsor is subordinate to
giving an atlas $U\ra\cX$ over which it trivialises,
and a 1 co-cycle,  which for $(s,t):R=U\ts_\cX U\rras U$ is,
in turn, subordinate to giving a map $\pi_0(R)\ra \G$.
If, however, the torsor is connected then $\pi_0(R)$
must generate $\G$, so, up to isomorphism, Galois
coverings are a set. As to the second part, (a) iff
(b) is the definitions; (a)$\Rightarrow$(d)
by taking $V\ra \cX$ a neighbourhood over which
$p_1$ trivialises; 
and (d) iff (c) is again the definitions modulo
the pro-Galois correspondence of the pre-amble.
As such suppose (d), choose a
base point $*:\rp\ra\cX$ 
then we have a net $q_i:\cE_i\ra\cX$, $i\in I_*$, of connected 
$\pi^i_1$- identified with $q^{-1}_i(*)$- 
torsors,
along with transition functions $F_{ji}$,
commuting
as in \eqref{3146}. 
Without loss of generality we can 
identify the
$\cE_i$ with locally constant sheaves, equivalently
make a legitimate but simplifying a priori choices,
of the pull-back of $\cE_i$ to any $V\ra\cX$ in
the small \'etale site of $\cX$,
and, similarly, not only identify any point, $x$, with a point on
an \'etale neighbourhood of the same, 
but also $I_x$ with $I_*$
as co-filtered directed sets. There are, however, by 
the pro-Galois correspondence, a (not necessarily
unique, albeit it's implicitly rigidified by
the points $*_i$, and $x_i$, $i\in I$, but
we don't need this) equivalences of fibre functors
\begin{equation}\label{31570}
\a_i(x): I_* \ra \Hom_{\pi^i_1} (q_i^{-1}(*), q_i^{-1}(x)),
\quad \a_i(*)=\mathbf{1}
\end{equation}  
as $x$ varies over a sets worth of points 
covering $\cX$. On the other hand by the
hypothesis of (d) each $x$ has a
connected
\'etale neighbourhood $U_x\ni x$ 
over which $\cE_i$ trivialises,
so there is a unique net of
trivialisations of torsors such that 
\begin{equation}\label{3157}
A_i: U
\bigl(:=\coprod_x U_x\bigr)
\ts \pi^i_1\xrightarrow{\sim} \cE_i\vert_U,\quad A_i(x)=\a_i(x)
\end{equation}
Thus, for $(s,t):R:=U\ts_\cX U\rras U$, 
we not only get
a $\pi^i_1$-valued co-cycle
\begin{equation}\label{3158}
c_i:=(t^*A_i)^{-1}s^*A_i: \pi_0(R) \ra \pi^i_1
\end{equation}
but $c_j=F_{ji}c_i$ for all $i\geq j$.
As such, $c:=\varprojlim_i c_i$,
defines a map from $\pi_0(R)$ to $\pi^\cts_1$, whose image
generates a discrete group $\pi_1$, while the induced
map $c:\pi_0(R)\ra \pi_1$ is a co-cycle. Consequently,
$c$ defines a connected $\pi_1$ torsor, $p_1:\cX_1\ra\cX$,
and since each $\cE_i$ is connected, the map induced
from $\pi_1\ra\pi^i_1$ by $c_i$ is surjective, so 
$\pi_1=\pi^\cts_1=\pi^\pro_1$, with
$p_1$
the desired universal widget. 
\end{proof}
Unsurprisingly, albeit bearing in mind that in
this generality fibration may be an empty notion
so the definition 
of $\et_2$ is, \ref{flav:funct}, 
via local trivialisation
and, of
course, continuing to employ the notation
\ref{def:394}, 
we move to 
\begin{fact}\label{fact:396} (Huerwicz-2) Suppose $\cX$ is 
connected and
everywhere
locally 
simply
connected ({\it i.e.} every point has a co-final
system of simply connected neighbourhoods) then for $*:\rp\ra \cX$
any point and $I_*$ the co-filtered directed class of
\eqref{eq:funct1406}-\eqref{eq:funct1106}
restricted to the 2-Galois objects
over the universal, \ref{fact:395}.(a),  cover, 
$\cX_1$, the category
$I_*$ is equivalent to a small category, and
there is a pro-discrete 2-group 
$\Pi_2^{\pro}(\cX_*)$ such
that
the fibre functor, 
\eqref{385} \& \eqref{392}, affords
an equivalence of categories between
$\et_2(\cX)$ and
pro-discrete groupoids with $\Pi_2^{\pro}(\cX_*)$ action.
Better still: the 2-type is
the discrete group $\pi_1(\cX_*)$,
a pro-abelian second homotopy group 
$\pi^\pro_2(\cX_*)={\pi_2^i}_{i\in I}$, and a Postnikov
class, $k_3\in\rH^3_\cts(\pi_1, \pi^\pro_2)$,
\ref{defn:3810}, 
representable by a 3 co-cycle, $K_3$, in 
$\Hom^\cts_\Ens(\pi_1, \pi^\cts_2)$; the Huerwicz theorem
holds, {\it i.e.} for a locally constant sheaf $\uZ$
of abelian groups on $\cX_1$,
\begin{equation}\label{31166}
\rH^2 (\cX_1, \uZ)\xrightarrow{\sim}
\varinjlim_i \Hom (\pi^i_2, Z)
\end{equation}
so, the pro-discrete analogue of \ref{fact:386}, can
be read directly from the Hoschild-Serre spectral
sequence, \ref{fact:pos2};
and the following are equivalent

(a) $\cX$ admits a universal \'etale 2-covering $p_2:\cX_2\ra\cX$.

(b) $\pi_2^\pro(\cX_*)$ is representable by a discrete group.

(c) $\cX$ is semi-locally 2-connected, {\it i.e.} every
point, $x:\rp\ra\cX$, 
has an \'etale neighbourhood such that the image
of $\pi^\pro_2(V_x)\ra\pi^\pro_2(\cX_x)$ is trivial.

(d) Every point, $x:\rp\ra\cX_1$,
has an \'etale neighbourhood, $V$, such that any
locally constant gerbe on
$\cX_1$ can be trivialised over $V$.
\end{fact}
As a synopsis and generalisation of the
entire chapter, it seems appropriate to
break into pieces the pro-2-Galois pre-amble,
of \ref{fact:395}
beginning with
\begin{DefPi2}
By \ref{fact:395} everything in $\et_2(\cX)$ has a
universal cover, so in the construction, \ref{SS:III.3}-\ref{SS:III.4},
of $\pi_2$ we can replace the starting point, 
\ref{def:funct103}, of quasi-minimal cells by 
the easier condition of simply
connected cells. The definitions \ref{def:noplag1},
\ref{def:noplag2} and proofs of existence,
\ref{fact:noplag5}, \ref{fact:noplag5noplag5}, 
stand as given, but, with the simplification
that simply connected 2-Galois cells can be
expressed as locally constant gerbes, $q:\cY\ra\cX_1$,
in $\rB_Z$'s for some discrete abelian group, $Z$,
while the relative stabiliser, $S_{\cY/\cX_1}$, is
exactly the constant sheaf $\uZ$ rather than something
locally isomorphic to it. Now, let's examine
these for a given $Z$: up to isomorphism \'etale atlases
of $\cX$ are a set, and for a sufficiently fine one,
$U\ra\cX$, $q|_U$ has sections everywhere, so this
cover factors as $U\ra\cY\ra\cX_1$, and whence,
\eqref{eq:sep3},
$q$ is equivalent to a map $R\ra R_1\rras U$ of groupoids,
which exhibits the former  as a $Z$-torsor over
the latter. Better still, by \ref{fact:noplag2},
any sub-group of $Z$ is (up to the choice of a
point) naturally a normal sub-groupoid of $R$, so
there is a maximal intermediary $R\ra R'\ra R_1$
with the property that the latter map has a
section, and the possibilities for $R'$ are
subordinate to the universal covers of
the components of $R_1$, which are independent 
of $Z$. Consequently,
to prove that the 2-Galois cells form a set,
we can suppose that $R'=R_1$, and $q$ is
subordinate to a 2 co-cycle $\pi_0(R_{1\,t}\ts_{s}R_1)\ra Z$.
Such a co-cycle must, however, generate $Z$- otherwise
$\cY$ wouldn't be simply connected- so, in totality,
the equivalence classes of 2-Galois cells are a set.
Consequently, the definition, \ref{factdef:plagAct1}, of
$\pi_2$ holds up to limiting $I$ of {\it op. cit.}
to simply connected 2-Galois cells, and checking that
this is co-filtering. As such suppose a pair, 
$q':\cY'\ra\cX_1$, $q'':\cY''\ra\cX_1$ of simply
connected 2-Galois
cells are given, form (over $\cX_1$) their fibre product, $q_1$,
then any other connected 2-Galois cell must factor through
the universal cover of the quasi-Galois cell constructed in
the proof of \ref{fact:noplag3}, which by the functoriality
of universal covering is itself 2-Galois, so $I$ is indeed co-filtering.
\end{DefPi2}
\begin{DefPos} Most of the conditions, \ref{def:PosSeq}. (a)-(d),
in the definition of a Postnikov sequence have no sense
for pointed maps. The condition \ref{def:PosSeq}.(e) not only
has a pointed sense, but is even trivial when everything
is pointed and connected. In general a `point' is really
the pair, $(*_i, \phi_i)$ of \eqref{eq:funct1406}, {\it i.e.}
a point of the fibre
{\it cf.} \ref{fact:cor2}, which in turn is best thought
of as being a pointed 0-cell in the obvious variation of
$\underline{\mathrm{Cham}}\mathrm{p}\underline{\mathrm{s}}/\cX$
wherein one takes diagrams exactly as in 
\eqref{eq:cor1}. Irrespectively, a pointed map, is a
pair $\uF_{ji}:=(F_{ji}, \z_{ji})$ as in 
\eqref{eq:funct1406}, equivalently the triple
$(f_{ji}, \xi_{ji}, \z_{ji})$ of {\it op. cit.}, with
the commutativity condition therein. Consequently,
given $i\geq j\geq k$ in the implied co-filtered directed set,
along with pointed maps $\uF_{ji}$ {\it etc.}, 
as soon as $q_k$  is 2-Galois there is 
(up to representable base change should the 
universal cover not exist) a unique 2-cell
$\g_{kji}:F_{ki}\Rightarrow F_{kj}F_{ji}$ compatible
with the pointing, 
{\it i.e.} \eqref{384} holds,
and the tetrahedron condition,
\ref{def:PosSeq}.(e) trivially follows. Therefore, an
alternative to the definition of Postnikov sequence
is to replace \ref{def:PosSeq}.(e) by the a priori
condition that all maps are pointed, say $\uF_{ji}$,
but the conditions \ref{def:PosSeq}.(a)-(d) should
be understood in an un-pointed way modulo the
2-cell $\g_{kji}$
being afforded by the pointing. In the pro-finite
case, \eqref{384}, this is equivalent to the 
non-pointed definition, but otherwise there is
a difference as soon as $I$ is large enough
that ${\varprojlim_i}^{(2)}$ may be non-zero. Irrespectively,
{\it i.e.} without changing the definition,
the proof, 
\ref{5.3}-\ref{5.8}, of the existence
of Postnikov sequences is valid as stated, but since
the pointed 2-category also has fibre products,
\ref{fact:PosSeq1} holds in a pointed way, so
one deduces the slightly stronger statement that
there exist co-final pointed Postnikov sequences,
equivalently \ref{def:PosSeq}.(a)-(e) and \eqref{384}
hold. 
\end{DefPos}
\begin{UniPos} A co-final Postnikov sequence defines,
\eqref{350}, a 3 cocycle, 
$K_3$, in $\Hom_{\Ens}(\pi_1^3, \pi_2^\cts)$,
but it needn't be unique in the simple way, 
\ref{5.15}, which one encounters in the pro-finite case.
Let us first consider the image of 2 such co-cycles,
$K$, $K'$ in $\rH^3_\cts (\pi_1, \pi_2^\pro)$ by
way of the filtration, $F^p$, afforded by the spectral 
sequence \eqref{3141}. To this end let $\d$ be 
the differential in the C\v{e}ch description of
the higher lims, $D$ the group differential, and
otherwise notation as in the proof of \ref{5.15}.
Plainly $K-K'$ is zero in $F^0$, and it's image
in $F^1$ is $\d (z^i_{\tau,\om})$- \eqref{382}. 
This vanishes, however, because 
$\d (z^i_{\tau,\om})=D (z^{ji}_\om)$- \eqref{381}-
and the image in $F^2$ is $\d(z^{ji}_\om)$, which
in turn is $D(z_{kji})$- \eqref{380}. Now if we're
in a pointed situation, then we can choose the
$\xi_{ji}$ of \eqref{378} in a unique way compatible
with the pointing, $z_{kji}=0$, and we stop. Otherwise,
the image in $F^3$ is $\d(z_{kji})$, and this is zero
by \ref{def:PosSeq}.(e). Identifying the co-cycles
$K$, $K'$ with nets of 2-Groups,
$\gP_i$, $\gQ_i$, and functors, this
can, \ref{fact:group1}, be usefully visualised as 
\begin{equation}\label{3167}
 \xy
 (38,0)*+{\gP_j}="B";
 (0,-6)*+{\gP_k}="C";
(15,0)*+{\gP_i}="D";
 (38,-12)*+{\gQ_j}="F";
 (0,-18)*+{\gQ_k}="G";
(15,-12)*+{\gQ_i}="H";
{\ar_{}_{B_j} "B";"F"};
    {\ar_{B_k} "C";"G"};
    {\ar_{B_i} "D";"H"};
{\ar_{}^{F_{kj} } "F";"G"};
    {\ar@{-->}_{F_{ki} } "H";"G"};
    {\ar^{F_{ji} } "H";"F"};
{\ar_{}^{F'_{kj}} "B";"C"};
    {\ar_{F'_{ki}} "D";"C"};
    {\ar^{F'_{ji}} "D";"B"};
\endxy
\end{equation}
The vanishing in $F^0$ says that the arrows $B_i$ exist;
the vanishing in $F^1$ says that the 
horizontal faces are transforms;
the vanishing in $F^2$ that the interior of the
diagram is a modification, \cite[1.5.12]{tom}; and the
vanishing in $F^3$ that this is all coherent in 
the way a pro-bicategory
should be. A priori this is perfectly sensible,
but it 
means that \eqref{3167} is
a 3-commutative diagram and we're
in pro-2-Cat, whereas if we point
the modification is trivial, \eqref{3167} becomes
a 2-commutative diagram, 
and we can shoehorn the thing into pro-Cat,
{\it cf.} \eqref{eq:groupTom} in \ref{sum:group1}.
In either case the Whitehead
theorem, \ref{Whitehead}, 
which is just the uniqueness of Postnikov
sequences in other clothes,
holds as stated.
\end{UniPos}
\begin{Pro2G}\label{Pro22G}
Since the pointing allows us to work in pro-Cat,
rather than pro-2-Cat, \ref{SS:III.6}-\ref{SS:III.7},
hold, actually with some
simplification because there
is a universal cover, up to the expedient of replacing $\cts$ by
$\pro$ where appropriate.
\end{Pro2G}
This dealt with we can move to  
\begin{proof}[proof of equivalence of (a)-(d) in \ref{fact:396}]
(a) iff (b) and (c) iff (d) are tautologies modulo 
the 2-Galois correspondence \ref{Pro22G}, while 
(a)$\Rightarrow$(d) by taking $V$ to be an
\'etale neighbourhood over which $p_2$-trivialises.
As such suppose (d), then 
for $I_*$ the co-filtered directed set of \ref{fact:396},
we have simply connected
gerbes, $q_i:\cE_i\ra \cX_1$ with fibres $\rB_{\pi_2^i}$'s;
1-cells $F_{ji}$; 2-cells $\g_{kji}$ satisfying the
commutativity conditions encountered in \ref{def:3811};
and an \'etale neighbourhood $U_x\ra \cX_1$ of a
covering set of points, $x$, such that each $q_i$
is trivialisable over $U_x$. Consequently, 
each $q_i$ has a section over 
$U:=\coprod_x U_x\ra \cX_1$, so for $(s,t):R_1:=U\ts_{\cX_1} U\ra U$,
each $q_i$ is, \eqref{eq:sep3}, representable by
an \'etale functor, $q_i:R^i\ra R_1$ of groupoids,
respectively the $F_{ji}$'s,
which, rather canonically, 
\ref{fact:noplag2} \&
\ref{fact:plagAct3}, is a torsor under $\pi_2^i$,
respectively the kernel $\pi_2^{ij}$ of
$\pi_2^i\ra\pi^j_2$.
By hypothesis each connected component of $R_1$
has a universal cover by \ref{fact:395}, and $I_*$
is co-filtering, so there is a unique \'etale groupoid
$R'\ra R_1$ through which each $q_i$ factors such
that every $q'_i:R^i\ra R'$ has a section. Consequently
on replacing $\cX_1\xrightarrow{\sim} [U/R_1]$ by
$[U/R']$, we can suppose that every $q_i$ has 
a section, or, equivalently, 
\eqref{3127},
that it's class in $\rH^2(\cX_1, \pi^i_2)$
actually belongs to $\check{\rH}^2(\cX_1, \pi^i_2)$. As
such consider pairs $(J,\z)$ where
$J\subseteq I$ is a 
right co-filtered
sub-partially ordered set; 
\begin{equation}\label{FixPro2}
\z^j_{g,f}:\pi_0(R_{1\, t}\ts_s R_1)\ra\pi_2^j
\end{equation}
is a co-cycle describing the gerbe $q_j$; and
for every pair $i\geq j$ in $J$
\begin{equation}\label{FixPro22}
F_{ji}\z^i=\z^j
\end{equation}
Such pairs may be partially ordered by
$\underline{<}$ where
\begin{equation}\label{FixPro1}
(J,\z) \underline{<} (K, \xi )\,\,\text{iff}\quad\,\, 
J\subseteq K,\, \z^j=\xi^j,\,\, \forall j\in J
\end{equation}
and, 
by choice, there is at least one maximal pair, say $(J,\z)$.
Consequently we get a co-cycle
\begin{equation}\label{FixPro3}
\z:=\varprojlim_{j\in J} \z^j: \pi_0(R_{1\, t}\ts_s R)\ra
\varprojlim_{j\in J}\pi_2^j
\end{equation}
whose image generates a discrete group, $\pi^\ell_2$, so
the resulting gerbe, $p_2:\cX_2\ra\cX_1$,
in $\rB_{\pi^\ell_2}$'s is simply connected. Consequently,
as the notation suggests, $\ell\in I$, and, 
since every
$\cE_j$ is simply connected,
the image of each $\z^j$ generates $\pi^j_2$
so $\ell\geq j$
for all $j\in J$. As such, the pair $(J,\z)$ is simply 
the set 
$J=\{i\mid i\leq \ell\}$, while
$\z^j=F_{\ell j} \z^\ell$. Now, as ever, there is some
section $T^\ell$ of $q_\ell$ such that
\begin{equation}\label{FixPro4}
T^\ell_{gf}= T^\ell_g T^\ell_f \z_{g,f}
\end{equation}
for all compossible arrows $f,g\in R_1$. If, however, $\ell$
were not the max of $I$ there would be an $i>\ell$ such
that $F_{\ell i}:R_i\ra R_\ell$ is a trivial $\pi_2^{i\ell}$
torsor, so there would be a section $T^i$ of $q_i$ such
that $F_{\ell i}T^i=T^\ell$, and whence a resulting co-cycle
$\z^i$ contradicting the maximality of $(J,\z)$.
\end{proof}
Now while the existence of the universal $q$-cover
doesn't imply everywhere locally $q-1$-connected
there is a certain naturality in such hypothesis,
and, as far as the proof is concerned, even a moral
necessity. There is, however, no such hypothesis 
in the passage from 1-Galois to 2-Galois in the
pro-finite case, wherein both were developed under
the hypothesis of local connectedness. Nevertheless
this is a 0-connectedness hypothesis, which is
in fact logically un-necessary since it is possible
to make
\begin{ProF}\label{ProFF}
Since locally Noetherian implies locally connected
in the algebraic flavour, let us
stick to champs, $\cX$, in 
separated topological spaces for the purposes
of the exposition, even though the discussion is
manifestly valid on fairly arbitrary sites. As such
the starting point is to realise that there is 
always a universal map, 
\begin{equation}\label{ProF1}
p_0^\fin:\cX\twoheadrightarrow \pi_0^\fin(\cX)
\end{equation}  
to compact seprated totally dis-connected spaces, 
{\it i.e.} exactly as in \ref{def:392} but use
finite sets instead of discrete ones, and profit
from the fact that pro-finite widgets are representable
in Top. Or, perhaps logically better,
as in \ref{def:391} since the key point in
the representability theorem is
that a compact separated space is
totally dis-connected iff it's totally separated.
Now, in an intended notational
confusion, a fibre, $\cX_*$, of $p_0$ need be
neither connected nor locally connected, but
it has all the necessary properties for doing
pro-finite Galois theory. Specifically, let
$q:\cY\ra\cX$ be a locally finite \'etale cover,
{\it i.e.} locally a product with a 
(possibly empty) finite
discrete groupoid,
then for every $m\in\bn$ the function,
\begin{equation}\label{ProF2}
\cX\ra \{0,1,\cdots,n\}:x\mpo \min\{m, \vert q^{-1}(x)\vert\}
\end{equation}
is continuous, so, taking $m\gg 0$ the function,
\begin{equation}\label{ProF3}
n(q,*): \cX_*\ra\bn: x\mpo \vert q^{-1}(x)\vert
\end{equation}
must be a constant. Quite possibly $n(q,*)=0$,
equivalently the fibre $\cY_*$ is empty, but 
this is the only difference with the connected
case, {\it i.e.}
\begin{fact}\label{fact:proF1}
For every locally finite \'etale cover $q:\cY\ra\cX$
with $n(q,*)\neq 0$ there exists a unique decomposition
of the fibre
\begin{equation}\label{ProF4}
\cY_*:= \cY\ts_\cX \cX_* = \coprod_i \cY_{i*} 
\end{equation}
where $\cY_{i*}\neq\emptyset$ is the fibre of an
open and closed subset of $\cY_i$ which itself is a 
$\underline{*\,\,\text{minimal coverin}}\text{g}$,
$q_i:\cY_i\ra\cX$, {\it i.e.} if $\cU\subseteq \cY_i$
is open and closed then either $\cU_*=\cY_{i*}$ or
$\cU_*$ is empty.
\end{fact}
\begin{proof} Let $\cU\hookrightarrow\cY$, then for 
$m\in\bn$ if $\cU$ is open, 
respectively closed,
points where the moduli
of the fibre $\cU_x$, $x:\rp\ra\cX$, has cardinality
at least $m$ are open, respectively closed.
Consequently if $\cU$ is both open and closed, 
so are the points where the fibres have 
cardinality at least $m$, and whence the function
$x\mpo \vert\cU_x\vert$ is continuous. Arguing as
in \eqref{ProF2}, the restriction of this function
to $\cX_*$ is, therefore, a constant, say $n(\cU, *)$.
Now proceed in the obvious way: if $\cY\ra\cX$ is
a $*$-minimal covering, stop, otherwise $\cY=\cY_1\coprod \cY_2$
where each $\cY_i\ra\cY$ is open and closed with
a non-empty fibre so 
each $\cY_i\ra\cX$ is a locally finite \'etale 
covering, and
$n(-, *)$ decreases strictly.
This proves the existence of the decomposition
\eqref{ProF4}, which in turn is unique since 
an intersection of finitely many open and closed
sets is again open and closed.
\end{proof}
Better still it follows from the proof that
\begin{cor}\label{cor:proF1} Every $\cY_{i*}\ra\cX_*$
occurring in the \eqref{ProF4} is an \'etale covering
with exactly the constant $n(\cY_i, *)$ points in
the moduli of every fibre.
\end{cor}
This is plainly enough to do 1-Galois theory, {\it i.e.}
there is a pro-finite group $\pi_1(\cX_*)$ such that
the image of the functor
\begin{equation}\label{ProF5}
\het_1(\cX)\ra \het_1(\cX_*)
\end{equation}
is, via the fibre functor on choosing a base
point $*:\rp\ra\cX_*$, equivalent to the category
of pro-finite sets with continuous $\pi_1(\cX_*)$-action.
It's not quite enough to do 2-Galois theory, but
going through the same rigmoral with the cardinality
of the relative stabiliser rather than the cardinality
of the fibre, 
then the isomorphism class of the stabiliser amongst
(the finitely many) groups of the same finite cardinality,
it's also true that
\begin{cor}\label{cor:proF2} Every $\cY_{i*}\ra\cX_*$
occurring in the \eqref{ProF4} 
admits a factorisation $\cY_{i*}\ra\cY'_{i*}\ra\cX_*$
where the latter 
is a representable \'etale covering
with exactly the constant $n(\cY_i, *)$ points in
the moduli of every fibre, and the former a locally
constant gerbe in $\rB_{\G_i}$'s for some finite 
group $\G_i$ depending only on $i$.
\end{cor} 
So that we get the un-conditional statement: 
there is a pro-finite 2-group $\Pi_2(\cX_*)$ such that
the image of the functor 
\begin{equation}\label{ProF6} 
\het_2(\cX)\ra \het_2(\cX_*) 
\end{equation} 
is, via the fibre functor on choosing a base 
point $*:\rp\ra\cX_*$, equivalent to the category
of pro-finite groupoids with continuous $\Pi_2(\cX_*)$-action.
\end{ProF}

\newpage
\section{Algebraic champs}\label{S:IV}

\subsection{GAGA and specialisation}\label{SS:IV.1}
In the first place we wish to consider
the 
relation between the
topological theory of \S.\ref{S:II} or,
more generally,
\ref{SS:III.9} and
the pro-finite 
(in the algebraic category)
theory of \S.\ref{S:III}.
Our concern is, therefore,
algebraic Deligne-Mumford champs
of finite type over $\bc$, 
which, in principle need
not even be separated if 
we employ the theory as developed
in \S.\ref{S:III}, albeit if one
wants to do calculations, equivalently
a description in terms of loops and spheres,
then this requires the separation
hypothesis of \S.\ref{S:II}.
A further a priori issue  with the sphere theory
of \S.\ref{S:II} 
is that it doesn't allow for nilpotent structure,
{\it i.e.} it's only been defined for reduced
holomorphic champs,
whereas the pro-finite, or
indeed pro-dicrete, \ref{SS:III.9}, 
theory does.
In reality, however, this is a non-issue since 
\begin{fact}\label{fact:411}
Let $X$ be a holomorphic, respectively algebraic, space
and $\et'(X)$ the category of holomorphic, respectively
algebraic, spaces which are \'etale (but, not fibrations
in the holomorphic case, respectively not proper in
the algebraic case) with morphisms understood to be
$X$ morphisms, then there  is an
equivalence of categories
\begin{equation}\label{411}
\et'(X)\xrightarrow{\sim}\et'(X_{\red}):\,
\bigl( Y\xrightarrow{q}X\bigr) \, \mpo\, 
\bigl( Y_{\mathrm{red}}\xrightarrow{q_{\red}}X_\red \bigr) 
\end{equation}
\end{fact}
\begin{proof}
The case of schemes is \cite[Expos\'e I.8.3]{sga1}, and the
proof goes through verbatim whether for holomorphic or
algebraic spaces.
\end{proof}
From which we have the following 
\begin{cor}\label{cor:412}
Let everything be as above, but now with $\et_2'(\cX)$
the 2-category of \'etale (again neither fibrations nor
proper) champs over a holomorphic, respectively algebraic
champs, $\cX$ with cells diagrams of the form \eqref{eq:cor1},
then there is an equivalence of 2-Categories
\begin{equation}\label{412}
\et'_2(\cX)\xrightarrow{\sim}\et'_2(\cX_{\red}):\,
\bigl( \cY\xrightarrow{q}\cX\bigr) \, \mpo\, 
\bigl( \cY_{\mathrm{red}}\xrightarrow{q_{\red}}\cX_\red \bigr) 
\end{equation}
so that in particular the full sub-categories $\et_2(\cX)$
(which in the case that hasn't been defined, {\it i.e.}
non-reduced holomorphic spaces, one takes the $0$-cells
to be locally of the form $U\uts \cG$ for $\cG$ a 
discrete groupoid) and $\et_2(\cX_\red)$ are equivalent
under the restriction of \eqref{412}.
\end{cor}
\begin{proof} By \eqref{eq:sep3} every diagram of the 
form \eqref{eq:cor1} on the left hand side of
\eqref{412} can be represented by functors
and natural transformations of groupoids acting on 
the same \'etale cover $U\ra \cX$, so the fact that
\eqref{412} on 1 and 2 cells is a family of fully faithful
functors follows on applying \ref{fact:411} with $X=U$
over all such $U$. To show essential surjectivity on
0-cells observe that if an atlas $U\ra\cX$ and a
0-cell $q_\red:\cY_\red\ra\cX_\red$ are given, then
there is an \'etale atlas $V_\red\ra\cY_\red$ 
factoring through $U_\red$ which by \ref{fact:411}
defines a unique \'etale $V/U$ with reduced structure
$V_\red/U_\red$. In particular, therefore, $V\ra\cX$
is a cover, and $\cX$, respectively $\cY_\red$, $\cX_\red$
can be identified with groupoid(s) acting on $V$,
respectively $V_\red$, so 
we conclude by
applying \ref{fact:411}, again, but now with
$X=V$. 
\end{proof}
Now there is an evident 2-functor $\cX\mpo\cX^{\hol}$
which associates to any algebraic (of
finite type over $\bc$) champ the corresponding
champ in holomorphic spaces. By definition,
\ref{defn:cor1}, all \'etale fibrations over
$\cX^{\hol}$ are 0-cells in $\et_2(\cX^{\hol})$,
and so we introduce
\begin{defn}\label{defn:413}
The 2-category $\et_2(\cX^{\hol})_{\mathrm{prop}}$ is
the full sub 2-category of $\et_2(\cX_{\hol})$ in
which the $0$-cells are proper over $\cX$. By
\ref{cor:412} and \ref{prop:cor1} it is equivalent
to the 2-category of finite groupoids with 
$\Pi_2^{\hol}:=\Pi_2(\cX^{\hol}_\red)$
action, which is equally the 2-category of finite
groupoids with continuous action under the
pro-finite completion $\hat{\Pi}_2^{\hol}$, which
could either be defined as the fundamental pro-2-group
of $\et_2(\cX_{\hol})_{\mathrm{prop}}$ in the sense of
\ref{def:Pro2type}, or, equivalently, but more usefully,
by \ref{prop:cor1},
as the unique (up to equivalence) pro-2-group obtained
on rendering the topological 2-type, $(\pi_1, \pi_2, k_3)$,
profinite, {\it i.e.} replace each $\pi_i$ by it's
completion, $\hat{\pi}_i$, in finite index (normal) sub-groups, and
$k_3$ by the resulting class $\hat{k}_3\in\rH^3_{\mathrm{cts}}
(\hat{\pi}_1, \hat{\pi}_2)$.
\end{defn}
With this notation we assert 
\begin{prop}\label{prop:415}
Let $\cX$ be a separated algebraic Deligne-Mumford champs
of finite type over $\bc$ then the natural 2-functor
\begin{equation}\label{413}
\et_2(\cX)\ra \et_2(\cX^{\hol})_{\mathrm{prop}}
\end{equation}
 is
an equivalence of 2-categories affording a continuous equivalence
of (topological) 2-groups
\begin{equation}\label{414}
\Pi_2(\cX) \xrightarrow{\sim} \hat{\Pi}_2^{\hol}
\end{equation}
\end{prop}
To this end observe the following reductions
\begin{lem}\label{lem:416}
It's sufficient to prove that \eqref{413} is essentially
surjective on 2-Galois cells.
\end{lem}
\begin{proof} By \ref{sum:plag1} and \ref{factdef:plagAct1}
this certainly implies that the $\pi_i(\cX)$ and the
$\hat{\pi}_i$, $i=1$ or $2$, are isomorphic. In addition,
as in the 
Whitehead theorem \ref{Whitehead}, Postnikov sequences are unique (and as it
happens exist trivially in the holomorphic case by
\ref{prop:cor1} ) so $k_3(\cX)$ and $\hat{k}_3$ would coincide too.
\end{proof}
\begin{lem}\label{lem:417}
A $0$-cell in 
$\et_2(\cX^{\hol})_{\mathrm{prop}}$
is algebraic iff it's covered by
an element of $\et'_2(\cX)$.
\end{lem}
\begin{proof} Necessity is obvious, and suppose conversely
that a holomorphic 0-cell $q:\cY\ra \cX^{\hol}$ is covered
by some algebraic space $V$ in $\et'_2(\cX)$. Consequently
by \eqref{eq:sep3} we can represent $q$ as a functor between
groupoids
\begin{equation}\label{415}
R\xrightarrow{q} R_0 \rras V
\end{equation}
where $R_0$ is an algebraic space. Now consider first the
case that $q$ is representable, then, by \eqref{3112},
$R$ is a sum of connected components of $R_0$ so it's algebraic.
As such by the factorisation, \ref{fact:flav1}, of $q$
into a locally constant gerbe and a representable map,
we can, without loss of generality, suppose that $q$
in \eqref{415} is a proper \'etale map, so we're done
by the Riemann existence theorem.
\end{proof}
At which point it's easy to give
\begin{proof}[proof of \ref{prop:415}]
By \ref{lem:417}, showing that the condition  
of \ref{lem:416} holds is local on $\cX$, so
it will suffice to prove that a proper \'etale
map $q:\cY\ra V$ from a holomorphic champ with
finite abelian stabiliser $\uZ$ to
a (reduced) algebraic variety is algebraic.
Factoring $q$ into a locally constant gerbe
and a representable map, the Riemann existence
theorem reduces us to the case that 
$\uZ$ is a constant sheaf, and
$q$ is a
locally constant gerbe in $\rB_Z$'s. The set of
all such gerbes modulo equivalence is, by \ref{fact:384},
isomorphic to $\rH^2(V,Z)$ in the classical
topology, which is equally $\rH^2_{\acute{\mathrm{e}}\mathrm{t}}(V,Z)$
by \cite[Expos\'e XVI.4.1]{sga4}.
\end{proof}  
It plainly, therefore, follows that whatever we
know about homotopy groups of algebraic varieties
over $\bc$ can simply be imported over any 
algebraically closed field of characteristic zero,
{\it e.g.} example \ref{eg:riemman}, immediately
implies by \ref{prop:415}, an exhaustive description
of all smooth 1-dimensional algebraic champs of
finite type over a field of characteristic zero,
and one might, 
irrespectively of the dimension,
be tempted to think following 
\cite[Expos\'e X]{sga1} that in characteristic
$p$ the prime to $p$ part of $\Pi_2$ will be
exactly as in characteristic zero. Certainly such
a supposition is reinforced by the fact that
this holds in co-homology, \cite[Expos\'e XVI.2.2]{sga4}.
Nevertheless, it's well wide of the mark in
homotopy, which is surprising enough to merit
\begin{scholion}\label{schol:GAGA} {\it Specialisation
of $\Pi_2$ in mixed characteristic.}
The set up is as follows: $\cX/S$ is a 
proper smooth geometrically
connected champs over a complete DVR, $S=\mathrm{Spec}(R)$,
of mixed characteristic, so, say $i:s\hookrightarrow S$
the closed point of characteristic $p$, and 
$j:U=\mathrm{Spec}(K)\hookrightarrow S$
the generic point of characteristic $0$. As such 
$\Pi_2(\cX_U)$ is wholly described by \ref{prop:415},
or some minor variant thereof if $K$ isn't algebraically
closed,
and
we
have maps
\begin{equation}\label{scholG1}
\Pi_2(\cX_U)\xrightarrow{j_*} \Pi_2(\cX) \xleftarrow{i_*}
\Pi_2(\cX_s)
\end{equation}
amongst which $i_*$ is an isomorphism by \ref{fact:Lef13},
so that having supposed geometrically connected fibres
to lighten the notation and remove any problems with 
base points, we get the specialisation map
\begin{equation}\label{scholG2}
\s: \Pi_2(\cX_U) \xrightarrow{(i_*)^{-1}j_*} \Pi_2(\cX_s)
\end{equation}
Of course, already for curves specialisation on $\pi_1$
fails to be an isomorphism, but what is true-
either by the non-abelian case of smooth base
change, \ref{cor:SmoothBaseChange}, or 
\cite[Expos\'e X.3.8]{sga1} which is a priori only
for spaces but is equally valid for champs with
the same proof, {\it op. cit.} 3.7, given \ref{cor:l13}- 
is that we have an isomorphism
\begin{equation}\label{scholG3}
\s': \pi'_1(\cX_U) \xrightarrow{\sim} \pi'_1(\cX_s)
\end{equation}
where $'$ denotes the maximal prime to $p$ quotient.
Now irrespective of how we may wish to define the
`maximal prime to $p$ quotient' of a 2-group it's
beyond doubt that it coincides with $\pi'_1$ for
$\rK(\pi_1, 1)$'s, while $\pi_2$ is, \ref{lem:l21},
unchanged by separably closed extensions of the
field of definition. Consequently if $X_{\bc}$ 
is a $\rK(\pi_1, 1)$, then the only possible sense
of $\Pi'_2(\cX_U)$ is $\pi'_1(\cX_U)$, and, in particular
$\pi'_2(\cX_U)=0$. On the other hand
\begin{fact}\label{fact:sp1}
Let $X/S$ be a $p$-adic model (with good reduction)
of a smooth bi-disc quotient which isn't a product
of curves then for $p$ belonging to a set of primes
of density $1/2$,
\begin{equation}\label{scholG4}
\pi'_2(\cX_s)\supseteq \prod_{\ell\neq p} \bz_\ell (1)
\end{equation}
\end{fact}
\begin{proof} By \eqref{3129} the right hand side
of \eqref{scholG4} is $\pi'_2(\bp^1_k)$ for $k$ any
field of characteristic $p$. If, however, $f:\bp^1_k\ra X_k$
is a rational curve on a projective variety over
$k$ then- again by \eqref{3129}- $f_*$ is injective
on $\pi'_2$, so it suffices to exhibit rational
curves on the reduction of $X$ modulo $p$, which is
the content of \cite[III.5]{charp}.
\end{proof}
As such, the good hypothesis appear to be of
Huerwicz type, to wit
\begin{fact}\label{fact:sp2}
Suppose the specialisation map, $\s_1$,  
\eqref{scholG2} is an isomorphism on all of $\pi_1$
(and not just the prime to $p$ part $\pi'_1$)
then there is an induced isomorphism
\begin{equation}\label{scholG5}
\s'_2:\pi'_2(\cX_U)\xrightarrow{\sim} \pi'_2(\cX_s)
\end{equation}
on the maximal prime to $p$ quotient of $\pi_2$.
\end{fact}
\begin{proof} Quite generally let $\cE\xrightarrow{g}\cX_1\xrightarrow{r} \cX$
be the factorisation of a champs proper and \'etale over $\cX$
into a representable cover and a locally constant gerbe
in $\rB_\G$'s for $\uG$ a sheaf of groups on $\cX_1$ of
some 2-cell $q$ in $\et_2(\cX)$. The construction of $r$,
\ref{claim:sep1}, is by the simple expedient of moding
out by the relative stabiliser $S_{\cE/\cX}$ after 
expressing $\cE$ and $\cX$ as \'etale groupoids acting
on the same atlas. As such the same construction works
for any other normal sub-sheaf of the locally constant
sheaf $S_{\cE/\cX}$. The stalks, however, of this sheaf
are isomorphic to $\G$ so the kernel $\G''$ of the
maximal prime to $p$ quotient $\G'$ sheafifies
to define such a sheaf, $S''_{\cE/\cX}$, and whence a
further factorisation
\begin{equation}\label{scholG6}
\cE\xrightarrow{g''}\cE'\xrightarrow{g'}\cX_1
\end{equation}
into locally a constant gerbe, $g''$, in $\rB_{\G''}$'s
followed by the maximal prime to $p$ gerbe, $g'$, in 
$\rB_{\G'}$'s. Now since the pull-back of connected
representable covers along gerbes are connected
representable covers, we always have surjections
\begin{equation}\label{scholG7} 
\pi_1(\cE_*) {\build\twoheadrightarrow_{}^{g''_*}}\pi_1(\cE'_*)
{\build\twoheadrightarrow_{}^{g'_*}}\pi_1(\cX_{1*}) 
\end{equation}
and the quasi-minimality, \ref{def:funct103}, condition for the 0-cell,
$q$ is that the composite in \eqref{scholG7} is injective, so
a fortiori we've proved
\begin{lem}\label{lem:scholG1}
Any quasi-minimal (respectively quasi-Galois, respectively
2-Galois) cell $q:\cE\ra\cX$ in $\et_2(\cX)$ admits a
factorisation 
\begin{equation}\label{scholG8}
\cE\xrightarrow{g''}\cE'\xrightarrow{q'}\cX
\end{equation}
such that the relative stabiliser $S_{\cE'/\cX}$ is not only prime
to $p$ but it's restriction to the base point $*$ is the
maximal prime to $p$-quotient of $\G=S_{\cE/\cX}(*)$, and
$q'$ is itself quasi-minimal (respectively quasi-Galois,
respectively 2-Galois).
\end{lem}
\begin{proof} We've just done the quasi-minimal case, 
and the procedure, \ref{def:noplag1}-\ref{fact:noplag5noplag5},
of going from this to quasi-Galois, respectively 2 Galois,
has no effect on the described property of the stabiliser
of $q'$.
\end{proof}
Consequently 
by \ref{fact:Lef13}- {\it cf.} the proof of \ref{cor:l21}-
to prove \eqref{scholG5} is surjective, it will suffice
to show that a 2-Galois cells $q'$, which factors
$\cE'\xrightarrow{g'}\cX_1\xrightarrow{r}\cX/S$ as
a prime to $p$-gerbe followed by a representable map
is again 2-Galois on restricting to the generic fibre.
We have, however, a commutative square of surjective maps
\begin{equation}\label{scholG9}
\begin{CD}
\pi_1(\cE')@<<{j_*}< \pi_1(\cE'_U)\\
@V{\wr}VV @VVV \\
\pi_1(\cX_1)@<{\sim}<{j_*}< \pi_1(\cX_{1,U})
\end{CD}
\end{equation}
Now suppose that
the right vertical isn't an isomorphism 
then there is a Galois
cover $\cF_U\ra\cE'_U$ which isn't
the pull back of a cover of $\cX_{1,U}$.
It has, however, a factorisation $\cF_U\ra\cY_U\ra\cX_{1,U}$
into a locally constant gerbe followed by
a representable cover, while, by hypothesis
$\et_1(\cX)\ra\et_1(\cX_U)$ is an equivalence
of categories, so $\cY_U$ extends to a unique
cover $\cY\ra\cX_1$. Better still,
by the functoriality of factorisation into
representable maps and locally constant gerbes,
$\cY\ra\cX_1$ is Galois,
so taking fibre products
with $\cY$ we can, without loss of generality,
suppose that $\cY=\cX_{1}$. As such, $\cF_U\ra\cE'_U$
is a Galois cover whose fibre over a geometric point
$*$ may be canonically identified- {\it cf.} proof
of \ref{claim:funct102}- with a quotient group
of $S_{\cE'/\cX}(*)$. In particular therefore it
is prime to $p$, so  \eqref{scholG3} applies to
afford an extension $\cF\ra\cE'$ of $\cF_U$ over
all of $S$.
Given which- {\it cf.} proofs of \ref{fact:Lef13} \&
\ref{claim:Lef21}- we certainly have the surjectivity
of \ref{scholG5} and
injectivity follows if for $\uZ$ the sheaf of
functions of a prime to $p$ finite abelian
group $Z$,
\begin{equation}\label{scholG10}
j^*:\rH^2(\cX_{1,S},\uZ)\xrightarrow{\sim} \rH^2(\cX_{1,U},\uZ)
\end{equation}
for any representable cover $\cX_1\ra\cX$. In the
particular case that $\cX$ is a scheme this is
\cite[Expos\'e XVI.2.2]{sga4}, and while, {\it op. cit.}
XVI.2.3, there are more intelligent ways of going
about it, it's equally, \cite[VI.4.2]{milne}, a 
formal consequence of proper base change,
\ref{cor:l12}, and smooth
base change,
\ref{cor:SmoothBaseChange}.
\end{proof}
Notice that en passant we've proved
\begin{cor}\label{cor:prime}
Quite generally,
let $\et'_2$ be the sub-2-category of $\et_2$
whose $0$-cells are the champs of \ref{lem:scholG1}
with prime to $p$ 
-for any given prime $p$, or even infinite set thereof-
relative stabiliser then $\et'_2$
is equivalent to the 2-category of groupoids with
prime to $p$-stabilisers on which a 2-group $\Pi'_2$-acts,
and the natural inclusion $\et'_2\hookrightarrow\et_2$
induces maps
\begin{equation}\label{prime1}
\pi_1(\Pi_2)\xrightarrow{\sim} \pi_1(\Pi'_2);\,\,
\pi'_2(\Pi_2)\xrightarrow{\sim} \pi_2(\Pi'_2);\,\,
k'_3(\Pi_2) \mpo k_3(\Pi_2)
\end{equation}
where $\pi'_2$ is the maximal prime to $p$ quotient
of $\pi_2$, and $k'_3$ the image of the Postnikov
class in $\pi'_2$.
In particular, under the (extremely strong) hypothesis 
on $\s_1$ of \ref{fact:sp2},
specialisation affords an equivalence of 2-Categories
\begin{equation}\label{prime2}
j^*(i^*)^{-1}: \et'_2(\cX_k)\xrightarrow{\sim} \et'_2(\cX_U)
\end{equation}
\end{cor} 
\begin{proof} Just apply the 2-Galois correspondence,
\ref{prop:376},
and the Whitehead theorem, a.k.a. unicity of Postnikov sequences,
\ref{Whitehead}.
\end{proof} 
It therefore follows that we can't just jump
from the holomorphic \ref{eg:riemman}
to a description of
\begin{ex}\label{ex:RiemannPos}
{\it The 2-category $\et_2(\cX)$ whenever $\cX$
is connected \'etale and proper over 
a positive characteristic orbifold $\cO$, i.e. a
smooth 1-dimensional (separated)
tame champs $\cO$ over a separably closed field
$k$ of characteristic $p>0$.}
A tame orbifold can always be lifted to the
spectrum, $S$, of the Witt vectors of $k$-
\cite[Expos\'e III.7.4]{sga1}, \cite[1.1]{km},
and the fact that there are no obstructions
to lifting points. As such by \ref{fact:Lef13},
we can identify $\et_2$ of $\cO$ with that of
its lifting, and similarly for the sub-2-category
$\et_2(\cX)$. Consequently, rather than supposing
things are defined over $k$, we may as well a priori
suppose everything defined over $S$ and proceed
with our initial/immediately preceeding notation,
\ref{schol:GAGA}, in the obvious way. 

This said, a simple inspection of cases shows
that if $\chi(\cO)>0$ then $\pi_1(\cO_U)$ is
prime to $p$, so by \eqref{scholG3}, such an $\cO$
satisfies the hypothesis of \ref{fact:sp2},
and, otherwise the said hypothesis are hopelessly false.
On the bright side however,
\begin{lem}\label{lem:CurveP1} 
$\pi_2(\cX)\xrightarrow{\sim}\pi'_2(\cX)$
\end{lem}
\begin{proof} By \ref{prop:376}
$\pi_2(\cX)$ is always a sub-group
of $\pi_2(\cO)$- in fact, \ref{eg:eg++1}, the kernel
of the corresponding pointed stabiliser
representation- so it'll suffice to prove
the lemma for $\cX=\cO$. To this end
let 
$\cY/k$ be a tame champ,
$\uV$ a locally constant 
sheaf of finite dimensional $\bF_p$-vector spaces
on the same; with
$V:=\uV\otimes_{\bF_p}\cO_{\cY}$ the
corresponding locally free sheaf; and
\begin{equation}\label{CurveArtinSch}
0\ra \uV \ra V \xrightarrow{AS} V\ra 0
\end{equation}
the resulting Artin-Schreier sequence. Now
since $\cY$ is tame, the cohomology, $\rH^q$, of any
coherent sheaf vanishes for $q\geq 2$, so 
we get an exact sequence
\begin{equation}\label{CurveArtinSch1}
 \rH^1(\cY, V) \xrightarrow{AS^1} \rH^1(\cY, V)\ra 
\rH^2(\cY,\uV)\ra 0
\end{equation}
where $AS^1$ is surjective by
\cite[III.4.13]{milne}, and whence
$\rH^2(\cY, \uV)=0$. On the other hand 
for any finite locally constant sheaf, $\uZ$
with torsion a power of $p$, the $p$-torsion
elements form a constant sub-sheaf of $\bF_p$
vector spaces, so by induction on the maximal
order of torsion
$\rH^2(\cY, \uZ)$ vanishes. If, however,
$\cE\xrightarrow{g''}\cE'\xrightarrow{q'}\cO$
is the prime to $p$-factorisation,
\ref{lem:scholG1}, of a 2-Galois cell
then $\cE'$ is tame while $g''$ is a locally
constant gerbe in $\rB_Z$'s for some abelian
$p$-group $Z$. Consequently $g''$ is trivial,
and $Z=\mathbf{1}$.
\end{proof}
Which in turn implies that even though $\cX$
may not be tame
\begin{fact}\label{fact:curveP2}
The conditions of \ref{fact:sp2} are satisfied
if $\chi(\cX)>0$, and better still there are
equivalences of 2-categories
\begin{equation}\label{prime22}
\begin{CD}
\et_2(\cX) 
@<{\sim}<< \et'_2(\cX_k)@>{\sim}>{j^*(i^*)^{-1}}> \et'_2(\cX_U)
\end{CD}
\end{equation}
In particular, therefore, $\et_2(\cX)$ is wholly described
by the holomorphic result, \ref{eg:riemman}.
\end{fact}
\begin{proof} Let $\cX\ra\cX_1\ra\cO$
be the factorisation into a locally constant gerbe
in
$\rB_G$'s followed by a representable cover.  
By \ref{prop:376} and
\ref{eg:eg++1} (which although a priori topological is just some general
nonsense about how to
read $\Pi_2$ of the sub-2-category
$\et_2(\cX)$ associated to a 0-cell from
the representation), $\pi_1(\cX)$, respectively
$\pi_1(\cX_U)$ is an extension of $\pi_1(\cX_1)=\pi_1(\cX_{1U})$
by quotients $\G'$, respectively $\G'_U$ of $\G$
by some central sub-groups $Z$, respectively $Z_U$.
These groups themselves are quotients of $\pi_2(\cX_1)$,
respectively $\pi_2(\cX_{1U})$ so by \ref{lem:CurveP1}
$Z_U\ra Z$ is surjective, whence it's an isomorphism.
Consequently the conditions of \ref{fact:sp2} are
satisfied and \ref{cor:prime} applies to give
the second isomorphism in \eqref{prime22}, while
the first isomorphism follows 
from \ref{lem:CurveP1} and the 2-Galois correspondence
\ref{prop:376}.
\end{proof} 
Otherwise, by \ref{fact:sp1} we have to calculate
$\pi_2$ by hand, beginning with
\begin{fact}\label{fact:IV115} Suppose, 
$\chi(\cX)\leq 0$ and
the orbifold
$\cO$ is actually a curve then $\pi_2(\cX)=0$.
\end{fact}
\begin{proof} Again $\pi_2(\cX)$ is a sub-group of
$\pi_2(\cO)$, so by \ref{lem:CurveP1} we just need
to do $\pi'_2(C)$ for $C/k$ a smooth curve of non-positive
Euler characteristic over a separably closed field.
As such let $\cE\xrightarrow{g} C_1\xrightarrow{r} C$ be
the factorisation of a 2-Galois cell into a representable
cover and a locally constant gerbe in $\rB_Z$'s for
some prime to $p$ abelian group $Z$. Base changing 
by representable covers as
necessary, we can suppose that $C=C_1$, and the stabiliser 
$S_{\cE/C}$ is the trivial sheaf $\uZ$. If $Z=\mathbf{1}$
we're done, so otherwise choose a 
(non-trivial) cyclic
quotient $Z\twoheadrightarrow \bz/\ell$, and argue as
in \ref{lem:scholG1} in order to factor $g$ as 
$\cE\xrightarrow{h}\cE_1\xrightarrow{g_1} C$
where $g_1$ is a locally constant 2-Galois gerbe in $\rB_{\bz/\ell}$'s.
Consequently, without loss of generality, $g=g_1$, and
the isomorphism class of the gerbe $g$ belongs to
$\rH^2(C,\bz/\ell)$. As such, there is a trivial sub-case
where $C$ is affine, according to which the said group
vanishes, {\it i.e.} $g$ has a section, so it can never
be 2-Galois. Slightly (but not much) less trivially, $C/k$ is proper,
and up to a non-canonical identification of $\bz/\ell$
with $\bz/\ell(1)$, $\rH^2(C,\bz/\ell)$ is generated
by the class, $\d_c$, of a $k$-point $c\in C$. If, however,
$f:\tilde{C}\ra C$ is an \'etale cover of order $\ell$-
{\it e.g.} as afforded by the Jacobian- then $f^*\d_c$
vanishes in $\rH^2(\tilde{C},\bz/\ell)$, {\it i.e.} 
$\cE\ts_C\tilde{C}$ has a section, and, again, $\cE$ isn't 2-Galois.
\end{proof}
Plainly this can be extended in the obvious way, {\it i.e.}
\begin{fact}\label{fact:IV116} Suppose more generally that
$\chi(\cX)\leq 0$ and
the orbifold $\cO$ admits an \'etale cover by a curve, 
for example the moduli 
is affine or
satisfies $\chi(\vert\cO\vert)\leq 0$,
then
$\pi_2(\cX)=0$. 
\end{fact}
\begin{proof} Again, $\pi_2(\cX)$ is a sub-group of $\pi_2(\cO)$,
which in turn is unchanged by \'etale covers, so 
\ref{fact:IV115} or its
proof reduces us to proving the non-affine part of the ``for example''. 
As such denote the moduli 
(over a separably closed field $k$) by $C$, with $c_i$, $1\leq i\leq n$
the finitely many non-scheme like $k$-points of $\cO$.
By the tameness hypothesis, \ref{ex:RiemannPos}, the
fibres of $\cO$ over $c_i$ are $\rB_{\mu_{\ell_i,k}}$'s
for some integers $\ell_i$ prime to $p$. In particular there
are maps 
\begin{equation}\label{IV20}
\mu_{\ell_i}(k)\xrightarrow{\sim} \pi_1(\rB_{\mu_{\ell_i,k}})
\ra \pi_1(\cO)
\end{equation}
which by \eqref{scholG3} (or the theory of generalised
Jacobians if one wants to proceed purely algebraically)
generate, for $\ell$ the l.c.m. of the $\ell_i$, a quotient group
\begin{equation}\label{IV21}
0\ra \mu_{\ell}(k)\xrightarrow{x\mpo \sum x^{\ell/{\ell_i}}}
\coprod_i \mu_{\ell_i}(k)\rightarrow Q\ra 0
\end{equation}
of the maximal abelian quotient of $\pi'_1(\cO)$.
Notation established, we proceed by
induction on $\ell_{\max}:=\max_i \ell_i$ to show that $\cO$
is covered by a curve. Plainly $\ell_{\max}=1$ is trivial,
and otherwise let $J$ be the set of $i$'s such that
$\ell_i=\ell_{\max}$, and $1\leq m\leq n$ its cardinality. In the event that
$m\geq 2$, the quotient $Q'$ of $Q$ by
the image of $\coprod_{i\notin J} \mu_{\ell_i}(k)$
has the property that for every $j\in J$ the 
induced map $\mu_{\ell_j}(k)\ra Q'$ is injective.
Consequently if we take the Galois cover $\tilde{\cO}\ra \cO$
defined by $Q'$, the moduli of $\tilde{\cO}$ is everywhere
scheme like over the pre-image of every $c_j$, $j\in J$,
and $\ell_{\max}$ decreases strictly. Otherwise $m=1$, and
we use the hypothesis that $\chi(C)\leq 0$ to 
find a non-trivial (again abelian is fine) 
\'etale cover $\tilde{C}\ra C$,
so the pre-image of the unique point $c_j$ in
$\tilde{C}\ts_C\cO$ has cardinality $> 1$, and
we reduce to the previous case.
\end{proof}
This leaves us with the notoriously thorny case
that the moduli of $\cO$ is $\bp^1$, where
\begin{lem}\label{lem:IV117} The 
following are equivalent for tame orbifolds, $\cO$,
with $\chi(\cO)\leq 0$ and moduli $\bp^1_k$
over a separably closed field $k$ of characteristic $p$,

(a) Every such $\cO$ admits an \'etale cover
by a curve.

(b) Every p\underline{rime h}yp\underline{bolic trian}g\underline{le},
{\it i.e.} the orbifold, $\rT_\ul$, with
signature $\ul=(\ell_1, \ell_2, \ell_3)$ at $0,1,\infty\in\bp^1_k$
for $p, \ell_i$ 4 distinct primes and $\sum_i 1/\ell_i <1$,
admits an \'etale cover by a curve.

(c) For every prime hyperbolic triangle, $\rT_\ul$,
$\pi_1(\rT_\ul)\neq 1$.

(d) For every prime hyperbolic triangle, $\rT_\ul$,
$\pi_2(\rT_\ul)= 0$.
\end{lem}
\begin{proof} Trivially (a)$\Rightarrow$(b)$\Rightarrow$(c),
while  (b)$\Rightarrow$(d) by \ref{fact:IV115}, and
(d)$\Rightarrow$(c) by \eqref{3129}.
Otherwise observe that in the proof of \ref{fact:IV116} the only
placed that we used that the moduli wasn't $\bp^1$ was
in order to start the induction at $\ell_{\max}=1$, and at 
the last stage in order to ensure that $m$ of {\it op. cit.}
was greater than $1$. In respect of the start of the 
induction this either has to start at $\ell_{\max}=2$,
for $p\neq 2$, or $\ell_{\max}=3$ in characteristic $2$.
In either case, all the weights at all the non-scheme
like points are the same, so the cover defined by $Q$ of
\eqref{IV21} does the job.
As to the end of the induction,
if some not necessarily maximal weight, $\ell_k$,
occurs at least twice, then the cover defined by
the quotient
\begin{equation}\label{IV23}
\coprod_{\ell_i\neq \ell_k} \mu_{\ell_i} \ra Q \ra Q_k\ra 0
\end{equation}
not only increases $m$, but ensures that on the resulting
cover the number of points with a given weight are at
least 2. This latter property is preserved under
coverings, so we've already proved
\begin{equation}\label{IV24}
\begin{split}
&\text{$\cO$ is \'etale covered by a curve if the same weight, $\ell_k$,
occurs more than}\\ 
&\text{once, or, more generally a pair of weights
have a common factor,}
\end{split}
\end{equation}
which by inspection covers all cases with $\chi(\cO)=0$.
Similarly if for some prime hyperbolic triangle,
we have a non-trivial
Galois covering $\cT\ra \rT$, its moduli is either
$\bp^1_k$ and \eqref{IV24} is satisfied, or it's
non-rational and \ref{fact:IV116} applies. As such
not only does
(c)$\Rightarrow$(b), but it implies that 
for every $N\in\bn$, every
hyperbolic triangle has a Galois cover such that
the fibre of every geometric point has cardinality
at least $N$. Consequently in the final stage of the
proof of \ref{fact:IV116} we can raise $m$ of {\it op. cit.} 
by way 
an \'etale  cover of the form $\cT\ts_{\rT_\ul} \cO\ra \cO$
whenever $\cO$ admits a map to a hyperbolic triangle
$\rT_\ul$. By \eqref{IV24}, however, $\cO$ has at least
$3$ pairwise relatively prime weights, so the only scenario
in which this can't be done is if there are exactly 3
non-scheme like points with weights powers of 2,3, and 5
respectively. As we've already observed, however, the
prime but parabolic triangle $\rT_{2,3,5}$ has fundamental
group $\mathrm{A}_5$ of order 60 in characteristic 0, whence, by hypothesis,
prime to $p$, so we can use \eqref{scholG3} instead.
\end{proof}
Now,  \ref{lem:IV117}.(c) may look simple enough
but it's already a pain in characteristic 0, and in
characteristic $p$ the situation is much worse since
one is a priori limited to \eqref{scholG3}, while
\begin{fact}\label{fact:IV118}
If $\rT_\ul$ is a prime hyperbolic triangle in
characteristics distinct from $2$ or $3$, 
its prime to $p$ fundamental group $\pi'_1(\rT_\ul)\neq 1$.
In particular, therefore, if $p\neq 2$ or $3$, 
and $\chi(\cX)\leq 0$, $\Pi_2(\cX)$
can be identified with its (
certainly more complicated than in
characteristic zero, but still topologically finitely generated)
fundamental group $\pi_1(\cX)$, and $\et_2(\cX)$ is
equivalent to $\et_2(\rB_{\pi_1})$ 
as described in \ref{eg:eg1}-\ref{eg:eg2}
up to the evident
addition of residual finiteness in the latter's definition.
\end{fact}
\begin{proof} Let $\lb$ be a large prime to be
specified, then 
by a theorem of Macbeath, \cite{macbeath}, the finite simple
group $\PSL(2,\lb)$ is a quotient of the
characteristic zero fundamental group $\pi_1(\rT_{\ul, U})$
whenever it contains elements of order $\ell_i$,
$1\leq i\leq 3$. The said group had order
$(\lb-1)\lb(\lb+1)/2$, so we're done a fortiori if we
can find $\lb$ such that 
\begin{equation}\label{IV25}
\lb = 1,\,\, \mathrm{mod} \, \ell_i, \, 1\leq i\leq 3,\,\,
\text{but}\,\,\, \lb \neq \pm 1 \mathrm{mod} \, p
\end{equation}
Now for $\ell=\ell_1\ell_2\ell_3$, by Dirichlet's
theorems primes are equidistributed over the
congruence classes $(\bz/p\ell)^\ts$, which in
turn is canonically $(\bz/p)^\ts\ts(\bz/\ell)^\ts$,
so there are more primes
than there are conditions in \eqref{IV25} provided
$\phi(p) > 2$. 
\end{proof}
Given the rather unsatisfactory restriction on
$p$ in \ref{fact:IV118} let us make
\begin{rmk}\label{rmk:FeitThomson}
The Feit-Thomson theorem trivially implies that 
for all prime hyperbolic triangles in characteristic 2,
$\pi'_1(\rT_\ul)=1$, so some intelligence 
(or, alternatively, tedious grind to force
an ad hoc version of \cite[Expos\'e X.3.8]{sga1}
according to some specific facts about
$\PSL(2,\lb)$ or, more likely, 
some other relevant group since given its relation
to the modular group the failure of the proof
of \ref{fact:IV118} for $p=2$ or $3$ doesn't look like a
coincidence)
is
required to show that \ref{lem:IV117}.(c) holds.
Similarly, the classification of finite simple
groups doesn't imply $\pi'_1(\rT_\ul)=1$ for all
prime hyperbolic triangles in characteristic 3,
but a cursory glance at the tables does suggest
that it happens very very often.
\end{rmk}
\end{ex}
\end{scholion}

\subsection{Lefschetz for \texorpdfstring{$\pi_0$}{pi\_0}}\label{SS:IV.2}

By way of notation let us spell out our
\begin{setuprev}\label{setup:l01}
For $\cX/S$  a separated algebraic Deligne-Mumford champ of finite
type over a locally Noetherian algebraic space $S$, 
there is \cite[1.1]{km}, 
a proper moduli map $\mu_\cX:\cX\ra \vert \cX\vert$, 
or just $\mu:\cX\ra X$, or similar, if there is
no danger of confusion, which is universal for maps
to algebraic spaces. By an ample (relative to $S$) bundle,
$H$, is to be understood a line bundle on $\cX$ such 
that for some $n\in\bn$, $nH$ is the pull-back of an
ample divisor on $X$. Needless to say, if 
$S$ is quasi-compact, there is a $n\in\bn$ such for any line
bundle, $L$, on $\cX$, $nL$ is the pull-back of a 
bundle on $X$, but the existence of such an integer
even for a single bundle is a non-void condition
without quasi-compactness. In any case, by a {\it hyperplane
section}, is to be understood the zero locus $\cH\hookrightarrow\cX$
of a not
necessarily
regular section, $h$, of an ample line bundle, {\it i.e.}
the local defining equation might be a divisor of zero,
so the habitual exact sequence
\begin{equation}\label{l01+}
L^\vee\xrightarrow{h} \cO_\cX \ra \cO_\cH\ra 0
\end{equation}
may fail to be exact on the left.
\end{setuprev}
Before proceeding let us make a remark in the form of a
\begin{warning}\label{warn:l01}
In so called wild situations where the order of the
monodromy group of a geometric point $x:\rp\ra\cX$
is not invertible in $\cO_{X,x}$, the moduli $Y\ra X$
of an embedding $\cY\hookrightarrow\cX$ need not be
an embedding. We are, however, only interested in
the topological properties of hyperplanes $\cH\ra\cX$,
so by \ref{fact:411} there is (and we may well do
so without comment)
no loss of generality
in replacing $\cH$ by as large a multiple as we wish.
In particular therefore, we can suppose that $H$ is
an ample bundle on $X$, so by
\cite[1.8.F]{km}, $\cH$ is locally defined by 
the pull-back of a 
function on $X$. As such,  
for $H\hookrightarrow X$ a hyperplane section
there is a fibre square
\begin{equation}\label{l01}
\begin{CD}
\cX@<<<\cH \\
@V{\mu_X}VV @VV{\nu}V\\
X@<<< H
\end{CD}
\end{equation}
where although $\nu$ still may not be the moduli
map $\mu_{\cH}:\cH\ra\vert\cH\vert$, it is, nevertheless, the case,
{\it cf.} \cite[1.9]{km},
that the induced map $\vert\cH\vert\ra H$ is a universal
homeomorphism. 
\end{warning}
Now the plan is to reduce the Lefschetz theorem to
\begin{fact}\label{fact:l01}
Let 
everything be as in \ref{setup:l01}-
so, in particular $\cX$ is locally
connected and whence globally
a disjoint union of  connected components-
with $\cX/S$ proper enjoying
for all $s\in S$
a fibre $\cX_s$ which is everywhere of 
dimension
at least 1
then the inclusion of a hyperplane section yields 
a surjection
\begin{equation}\label{l02}
\pi_0 (\cH) \twoheadrightarrow \pi_0(\cX)
\end{equation}
\end{fact}
\begin{proof}
What the statement means is that if $\cX$ is
connected then $\cH$ is non-empty.
\end{proof}
Concretely, therefore, the first stage of the plan is
to bump this up to
\begin{fact}\label{fact:Lef02}
Let everything be as in \eqref{fact:l01} but suppose
further that the 
dimension 
of every fibre
is everywhere at least $2$ then there is a
Zariski open neighbourhood, $\cU\hookrightarrow \cX$,
of $\cH$ such that the inclusion $\cH\hookrightarrow\cU$
affords 
an
isomorphism
\begin{equation}\label{Lef01}
\pi_0 (\cH) \xrightarrow{\sim} \pi_0(\cU)
\end{equation}
Better still, there is a unique maximal Zariski open
neighbourhood $\cU_{\mathrm{max}}\supseteq \cH$ satisfying \eqref{Lef01},
and for any intermediate Zariski open neighbourhood, 
$\cH\subseteq \cU\subseteq \cU_{\mathrm{max}}$, \eqref{Lef01}
still holds. 
\end{fact}
As we've observed immediately post \eqref{l01} this
reduces to the purely schematic question for 
the hyperplane section $H\hookrightarrow X$ of {\it op. cit.}.
As such,
it's basically the easy part of \cite[Expos\'e XII.3.5]{sga2}.
We have, however, 
substantially weakened the hypothesis, {\it i.e.}
\begin{itemize}
\item The hyperplane section is no longer regular, \eqref{l01+}.
\item There is no hypothesis of flatness over $S$.
\end{itemize}
This said let us proceed to
\begin{proof}[proof of \ref{fact:Lef02}]
Let $\uZ_{\bullet}$ be the constant sheaf associated
to the torsion abelian group $Z$; $\s$
the structure map to $S$ of whatever; $s\in S$
fixed; and, for the moment, $S$ Noetherian. 
By hypothesis $H$ is a 
 a not necessarily regular
({\it i.e.} $H$ vanishing on 
components of $X$ is allowed) section
of an ample bundle $L$. As such, for some 
large $n$, to be decided, consider the space
$P$ of hyperplanes in $X$- {\it i.e.}
$\bp(\underline{\Hom}_S (\s_* L^{\otimes n}, \cO_S))$
in EGA notation-with
\begin{equation}\label{l0p2}
\begin{CD}
F@>>{p}>  P\\
@V{q}VV @.\\
X
\end{CD}
\end{equation}
the universal family, so for $n\gg 0$, $q$ is smooth.
Now identify $nH_s$ with a $k(s)$-point, $0$, of $P$,
then 
for $D=q^* H$
by 
(a double application of)
proper \'etale base change there is an \'etale 
neighbourhood $V\ra P$ of $0$ such that
\begin{equation}\label{l0p3}
\pi_0( D_V:= D\ts_P V)
\xleftarrow{\sim} \pi_0(H)\xrightarrow{\sim} \pi_0( F_V:= F\ts_P V)
\end{equation}
Consequently if $U$ is the
(Zariski open) image of $V$ in $P$,
and
$(s,t):V\ts_U V\rras V$,
then we have a  commutative diagram
\begin{equation}\label{l0p4}
\begin{CD}
0@>>> \rH^0 (D_U, \uZ) @>>> \rH^0 (D_V, \uZ) @>>{s^*-t^*}> 
\rH^0(D_{V\ts_U V},\uZ)\\
@. @AAA @AAA @AAA \\
0@>>> \rH^0 (F_U, \uZ) @>>> 
\rH^0 (F_V, \uZ) @>{s^*-t^*}>> \rH^0(F_{V\ts_U V},\uZ)
\end{CD}
\end{equation}
where the various subscripts, $_V$, {\it etc.}, denote the fibre over
the same. By construction $D$ is a hyperplane section
of $F/P$ all fibres of
which have dimension at
least 1, so \ref{fact:l01} applies to
yield that all the vertical maps are injective.
On the other hand, by  
\eqref{l0p3} the middle
vertical is an isomorphism, so
the leftmost vertical 
is too. 

At which juncture, we switch attention to 
the (smooth) projection $q:F_U\ra X$ whose
image $U'$ is a Zariski open neighbourhood
of $H$, so that for $(s,t):F_U\ts_X F_U\rras X$
we have a commutative diagram
\begin{equation}\label{l0p44}
\begin{CD}
0@>>> \rH^0 (H, \uZ) @>>> \rH^0 (D_U, \uZ) @>>{s^*-t^*}> 
\rH^0(D_U\ts_H D_U)\\
@. @AAA @AAA @AAA \\
0@>>> \rH^0 (U', \uZ) @>>> 
\rH^0 (F_U, \uZ) @>{s^*-t^*}>> \rH^0(F_U\ts_X F_U,\uZ)
\end{CD}
\end{equation}
in which the middle vertical is an isomorphism.
The necessary dimension condition to apply \ref{fact:l01}
is no longer  valid for the rightmost vertical,
but observe that if $z$ is a locally constant
function on $H$ viewed as such on $F_U$ then
$s^*z-t^*z$ vanishes a fortiori on the fibre
over $0\ts 0$. As such, by 
proper base change for $F\ts_X F/P\ts_S P$, 
$s^*z-t^*z$ vanishes on all the fibres of the
same over a Zariski open neighbourhood $W\ni 0\ts 0$.
This is, however, a classical topology so products
of neighbourhoods are co-final, {\it i.e.} without
loss of generality $W=U\ts_S U$, and whence the
leftmost vertical is an isomorphism. 

This proves \ref{fact:Lef02} locally on $S$, {\it i.e.}
on replacing $S$ by a Zariski open neighbourhood $S'\ni s$
in {\it op. cit.}. If, however, $U'\supseteq H\ts_S S'$,
respectively $U''\supseteq H\ts_S S''$ are Zariski
neighbourhoods satisfying \eqref{Lef01},  for
some Zariski opens $S',S''\subset S$ then 
$U'\cup U''$ also satisfies \eqref{Lef01}- 
argue exactly as in \eqref{l0p4} but with
$U'\coprod U''$ instead of $V$ while observing
that 
one can can replace $\cX$ by a Zariski open, $\cU$, in
\ref{fact:l01} provided that $\cU\supseteq \cH$. Consequently,
and irrespective of any Noetherian hypothesis, there is
not only a Zariski open satisfying 
\eqref{Lef01} but even a unique maximal one, $\cU_{\mathrm{max}}$.
Finally: if $\cH\subseteq \cU\subseteq \cU_{\mathrm{max}}$
is intermediate, then the map on $\pi_0$ is necessarily
injective, but, as above, \ref{fact:l01} is valid
for every Zariski neighbourhood of $\cH$.
\end{proof}

To extract a more standard looking corollary requires
\begin{defn}\label{rmk:l01} A champ $\cX/k$ of finite
type over a field is said, {\it cf.}
\cite[Expos\'e XIII.4.3]{sga2},
to have 
homotopy depth at least $2$ 
in the Zariski sense if for $X$ its moduli, 
at every 
closed (so quite possibly not geometric)
point $x\in X$ and
every 
sufficiently small (connected)
Zariski open neighbourhood $U\ni x$
\begin{equation}\label{hd1}
\pi_0(U\bsh x) \xrightarrow{\sim} \pi_0(U),\,\,
\end{equation}
which, \ref{schol:l01}, is the weakest local
connectivity hypothesis possible.
\end{defn}
According to which we have
\begin{cor}\label{fact:l02} Let everything be as in
\ref{fact:l01} and suppose further that for every
$s\in S$ the fibre $\cX_s$ is everywhere of dimension,
and homotopy depth in the Zariski sense, at least $2$,
then the inclusion of a hyperplane $\cH\hookrightarrow\cX$
affords an isomorphism
\begin{equation}\label{l03}
\pi_0(\cH)\xrightarrow{\sim} \pi_0(\cX)
\end{equation}
\end{cor}
\begin{proof}
Continuing in the notation of the proof of \ref{fact:Lef02},
we switch to the Nistnevich topology, $\mathrm{Nis}$,
where
by \ref{fact:l01} and
proper base change for the same
({\it i.e.} idem
for the \'etale topology plus 
Leray for \'etale $\ra$ Nistnevich)
we already know that
\begin{equation}\label{l0p1}
\s_*^{\mathrm{Nis}} \uZ_\cX \ra \s_*^{\mathrm{Nis}}  \uZ_\cH
\end{equation}
is injective. 
On the other hand by \ref{fact:Lef02}, applied
at $k(s)$, 
it's a fortiori
surjective at $s\in S$ whenever every Zariski
open neighbourhood $U_s$ of $H_s$, affords an
isomorphism
\begin{equation}\label{l0Cor1}
\pi_0(U_s)\ra \pi_0(X_s)
\end{equation}
The complement $X_s\bsh U_s$
has, however, dimension $0$, so this is immediate
by \eqref{hd1}.
\end{proof}
Before progressing let us observe
\begin{rmk}\label{rmk:l02}
By the simple expedient
of replacing $P$ and $F$ in \eqref{l0p2} by $P^m$, $m\in\bn$,
and the (not necessarily flat) universal intersection
of $m$ hyperplanes, one gets at no extra cost
\begin{fact}\label{GrothCon} cf. \cite[Expos\'e XIII.2.3]{sga2},
Let everything be as in \ref{setup:l01}, with
$\cH_1,\hdots, \cH_m$ not necessarily regular
hyperplane sections, then there is a Zariski
open neighbourhood $\cU\supseteq \cH_1\cap\hdots\cH_m$
such that the inclusion affords an isomorphism
\begin{equation}\label{GrothCon1}
\pi_0(\cH_1\cap\hdots\cap\cH_m)\xrightarrow{\sim}
\pi_0(\cU)
\end{equation}
In particular if for each $s\in S$ the moduli
of the fibre $\cX_s$ is everywhere 
of dimension at least $m$ and,
in
the Zariski sense, everywhere locally
connected in dimension $m$- \ref{rmk:l11}-
then
\begin{equation}\label{GrothCon2} 
\pi_0(\cH_1\cap\hdots\cap\cH_m)\xrightarrow{\sim} 
\pi_0(\cX) 
\end{equation}
\end{fact}
Of which, \eqref{GrothCon2}, equally follows immediately
from \ref{fact:l02}, {\it i.e.} the $m=1$ case and
\eqref{Csch2}. If, however, one is prepared to assume
\eqref{Csch2}, then there are, at least over a field,
other, {\it e.g.} \cite[\S 3.1]{joinsAndIntersections}, ways to proceed.  
\end{rmk}
Finally let's aim to clarify the ``homotopy
depth at least 2 hypothesis'' by way of
\begin{scholion}\label{schol:l01}
Plainly for any scheme, $X$, there are various senses in which
\eqref{hd1} can be understood at a point $x\in X$, according
to whether $U\ni x$ is Zariski; Nistnevich or \'etale
({\it i.e.} in either of the latter cases $U$ is \'etale,
but in the \'etale case we replace $x$ by its separable
closure $\bar{x}$). Equally plainly if $E\ra E'$ is a
surjective (in any reasonable sense) map of sites
then $E'$ disconnected implies $E$ disconnected, so
that we trivially have
\begin{equation}\label{hd2}
\text{\ref{hd1} for \'etale $\Rightarrow$
\ref{hd1} for Nistnevich $\Rightarrow$
\ref{hd1} for Zariski }
\end{equation}
while all the reverse implications are trivially false,
{\it e.g.} the latter fails for a nodal cubic
over $\bc$, and the former for $x^2-dy^2=0$
over $\bq$ whenever $d\in\bq $ is not a square. 
Similarly in all cases they are equivalent to
the vanishing of
\begin{equation}\label{l0exc1}
\rH^q_x(U, \uZ)=0,\quad q=0\,\, \text{or}\,\, 1
\end{equation}
in the appropriate site for
$\uZ$ a locally constant sheaf of abelian
groups, wherein
\eqref{l0exc1} always holds for $q=0$ if
the dimension at $x$ is at least 1. 
Of course, we should also occupy ourselves
about the difference between a champ and
its moduli. Here, a priori, only the \'etale
definition has sense, and locally 
around  $\bar{x}\in\cX$
the moduli is given by a quotient of a finite
group $G$ acting on some strictly Henselian 
neighbourhood $W\ni\bar{x}$ of $\cX$ while
fixing $\bar{x}$, so
$W\bsh \bar{x}$ connected certainly implies $W^G\bsh\bar{x}$
connected, {\it i.e.}  
\begin{equation}\label{hd3}
\text{\ref{hd1} \'etale locally on $\cX$ $\Rightarrow$
\ref{hd1} \'etale locally on $\vert \cX\vert$}
\end{equation}
and, in any case, the nature of 
0-connectivity in champs is such, {\it cf.}
\eqref{eq:point8}, that it's equivalent
to that of the moduli, so \eqref{hd1} in
the Zariski sense is the weakest possible
when the moduli is a scheme.

There is, however a somewhat more interesting
ambiguity between \eqref{hd1} for Henselian
local rings, $A$, and their completion, $\hat{A}$,
in the maximal ideal. This in turn may be 
subordinated to a lemma of independent utility,
{\it viz:}
\begin{lem}\label{lem:hd1}
Let $A$ be an excellent reduced Henselian local
ring, then $A$ is the set of $x\in \hat{A}$ such
that for some $n\in\bn$, $a_i\in A$, $0\leq i\leq n$
\begin{equation}\label{hde}
a_0 x^n + a_1 x^{n-1} +\hdots a_n=0,\,\, a_0\,\,
\text{not a zero divisor}
\end{equation}
In particular $A$ is integrally closed in $\hat{A}$,
and if $A$ is a domain 
it's algebraically closed.
\end{lem} 
\begin{proof}
For $A\sbs B\sbs \hat{A}$ a finite  $A$ algebra to
be determined consider the base change diagram
\begin{equation}\label{hd4}
\begin{CD}
\hat{A} @>>> \hat{B}:= B\otimes_A \hat{A} @>>> \hat{A}\otimes_A \hat{A}\\
@AAA @AAA @AAA \\
A @>>> B @>>> \hat{A}
\end{CD}
\end{equation}
wherein the diagonal $\hat{A}\otimes_A \hat{A}\ra \hat{A}$ yields
a retraction of the top left horizontal. This said we
proceed by cases- throughout which $K$ will
be the total ring of quotients of $A$- starting with
\begin{claim}\label{claim:hd1} The lemma holds if $A$ is normal,
{\it i.e.} $R_1 + S_2$.
\end{claim}
\begin{proof}[sub-proof] 
Observe that a normal local ring is a domain, while since
$A$ is excellent
$\hat{A}$ is also normal, so, in particular a domain. 
Now let $x\in\hat{A}$ satisfy \eqref{hde},
then the field $K(x)$ is a sub-field of the quotient
field of $\hat{A}$, whence the integral closure, $B$, of $A$
in $K(x)$ is naturally a sub-ring of $\hat{A}$. Excellent
implies universally Japanese, so $B$ is a finite $A$-algebra.
The latter is, however, Henselian so $B$ is a product of
finite Henselian local $A$-algebras, and since $B$ is a
domain, it must, therefore, be a local ring. Furthermore
it's of finite type over $A$ so it's also excellent, and
whence $\hat{B}$ is 
normal, so, in particular it's a domain. On the other hand
\begin{sublem}\label{sublem:hd1}
If $B/A$ is a finite algebra with a retraction, and
$B$ is a domain, then $A=B$.
\end{sublem}
\begin{proof}[sub-sub-proof]
Let $I$ be the kernel of the retraction $\b$, and for
$b\in I$ let $n\geq 1$ be the degree of a minimal monic polynomial,
\begin{equation}\label{hd5}
f(X)=X^n + a_1 X^{n-1} +\hdots + a_n \in A[X]
\end{equation}
such that $f(b)=0$, then applying $\b$ yields $a_n=0$,
so $f(b)=bg(b)$ where $g$ is monic of smaller degree,
and since $B$ is a domain $b=0$.\end{proof}
Applying the sub-lemma to $\hat{A}\ra \hat{B}$, shows
that they're equal, and since completion is faithfully
flat $A=B$.
By \eqref{hde}, however,  $a_0 x\in B$,
so $a_0 x=a\in A$, and whence we have an exact sequence
\begin{equation}\label{hd6}
0\ra (a_0, a)\ra (a_0) \ra Q\ra 0
\end{equation}
of $A$-modules
such that $Q\otimes_A \hat{A}=0$, so $Q=0$, 
{\it i.e.} $a_0\mid a$ in $A$, and $x\in A$,
since $a_0$ is not a divisor of zero.
\end{proof}

Unsurprisingly the next case is
\begin{claim}\label{claim:hd2} The lemma holds if $A$ is 
a domain.
\end{claim}
\begin{proof}[sub-proof] As before since $A$ is
Japanese: the 
normalisation $A'$ in $K$ is finite
over $A$, so it's a local Henselian domain.  
As such by \cite[7.6.2]{EGAIV}, the completion
$\hat{A}$ has 1-minimal prime, so it's a domain.
Similarly if $x\in\hat{A}$ is algebraic, then
for some $a_1\in A$, $B:=A[a_0 x]$ is finite over
$A$ (whence excellent) with quotient field $K(x)$. Now, again,
the normalisation $B'$ in $K(x)$ is finite,
thus $B'$ is a Henselian local ring, so $\hat{B}$
is a domain and \ref{sublem:hd1} applies as
before.
\end{proof}
 
This leaves us with the general case, so $K$ is
a (finite) product of fields $K_i$ corresponding
to minimal prime ideals/components $\gp_i$ of $A$.
In particular the normalisation, $A'$, is a product of
normal domains $A_i$, whence these are all Henselian
local rings, so- \cite[7.6.2]{EGAIV} again,
and, as above each $\hat{A}_i$ is a domain- the
minimal primes of $\hat{A}$ 
(which is reduced by \cite[7.6.1]{EGAIV}) 
are exactly of
the form $\hat{\gp}_i:=\hat{A}\gp_i$. 

Now exactly as in \eqref{hd6} it will suffice to
prove that the (necessarily reduced and finite) sub-ring 
$B:=A[a_0x]\sbs \hat{A}$ is equal to $A$. 
Each $\gq_i:=\hat{\gp}_i\cap B\supseteq \gp_i$ is
prime, whence since $B/A$ is finite, it is a minimal
prime over $0$ while
\begin{equation}\label{hd7}
\cap_i \gq_i= \cap_i \hat{\gp}_i\cap B =0
\end{equation}
so these are all the minimal primes. Exactly as above,
therefore, $\hat{B}$ is reduced with minimal primes
$\hat{\gq}_i$. 
To conclude: 
observe that since the retraction
$\b:\hat{B}\ra \hat{A}$ comes from the diagonal
in \eqref{hd4}, it sends $\hat{\gq}_i$ to $\hat{\gp}_i$,
so for $I$ the kernel of $\b$, \ref{sublem:hd1}
implies $I\sbs \hat{\gq}_i$, for all $i$, so $I=0$.
\end{proof} 
This may be applied to the question of ``homotopy
depth 2'' as follows: 
for $A$ as in \ref{lem:hd1},
let $U=\mathrm{Spec}(A)$,
and $\hat{U}=\mathrm{Spec}(\hat{A})$, with $x$
the closed point then
\begin{cor}\label{cor:hd1}
There is a  natural isomorphism
\begin{equation}\label{hd8}
\pi_0(\hat{U}\bsh x) \xrightarrow{\sim} \pi_0(U\bsh x)
\end{equation}
\end{cor}
\begin{proof} 
We may, \ref{fact:411}, suppose that everything is reduced,
while by hypothesis
everything is Noetherian,  so either side
is a finite direct sum of connected components. 
The trivial direction is: any open
and closed subset non-empty subset on the right pulls back
to the same on the left since completion is faithfully flat,
so \eqref{hd8} is always surjective. Conversely the connected
components can be identified with finitely many
indecomposable idempotents
$e_i\in\G(\hat{U}\bsh x)$, which, of course satisfy $e_i^2=1$,
whence they're integral over $U\bsh x$, so we conclude by
\ref{lem:hd1}.
\end{proof}
The intervention of idempotents suggests making 
\begin{rmk}\label{rmk:hd1}
Independent of any hypothesis of excellence:
just as in the trivial implications  \eqref{hd2}, 
\eqref{hd8} is, as noted, always surjective, so
the proof of \ref{cor:hd1} equally shows that \eqref{hd1}
holds as soon as $\rH^1_x(\cO_{\hat{U}})=0$, {\it i.e.}
the (formal) algebraic depth is at least 2. While sufficient
for homotopy depth $2$ it is far from necessary, and, in many
ways, {\it cf.} \cite[Expos\'e XIII.2]{sga2}, not particularly
desirable due to its non-topological nature.  
\end{rmk}
\end{scholion}
\subsection{Lefschetz for \texorpdfstring{$\pi_1$}{pi\_1}}\label{SS:L1}

We continue with the inductive strategy for bumping
up by way of
\begin{cor}\label{cor:l11}
Let everything be as in \ref{fact:l02} 
but 
with homotopy depth at least $2$
in the \'etale sense, and
suppose further that $\cX$ is 
connected (whence a hyperplane
section is connected by {\it op. cit.})
then the 
inclusion $\cH\ra\cX$ of a hyperplane section
affords a fully faithful functor
\begin{equation}\label{l11}
\et_1 (\cX) \ra \et_1(\cH)
\end{equation}
\end{cor}
\begin{proof} 
Since either sides is equivalent to a category
of finite sets on which $\pi_1(\cX)$, respectively
$\pi_1(\cH)$ acts, it's sufficient (and in fact
equivalent from the description \eqref{eq:plag1}
of Galois objects) to prove
\begin{equation}\label{l12}
\pi_1(\cH) \twoheadrightarrow \pi_1 (\cX)
\end{equation}
If, however, $\cY\ra\cX$ is a representable Galois
cover under a finite group $G$, then $G$ certainly
acts transitively on $\cY':=\cY\ts_\cX \cH$, so,
\ref{sum:plag1}.(b), $\cY'$ is Galois and \eqref{l12}
holds iff $\cY'$ is connected for all Galois $\cY/\cX$,
which, \ref{fact:l02}, is indeed the case since $\cY'$
is a hyperplane section of $\cY$, which certainly has
homotopy depth $2$ in the Zariski sense since
$\cX$ does in the \'etale sense.
\end{proof}
Of which a useful variant in the spirit of
Grothendieck's Lefschetz condition, \ref{fact:Lef02},
is
\begin{cor}\label{cor:l1new}
Suppose the simpler hypothesis that everything
is as the set up \ref{setup:l01}, with $\cX/S$ proper
enjoying  
for all $s\in S$
fibres, $\cX_s$, which are everywhere of dimension
at least 2, and define a category $\et_{1,\cH}^{\mathrm{Zar}}$ whose objects
are representable covers of a Zariski neighbourhood 
of $\cH$ modulo isomorphism over such neighbourhoods
then restriction affords a fully faithful functor
\begin{equation}\label{l1new}
\et_{1,\cH}^{\mathrm{Zar}}\ra \et_1(\cH)
\end{equation}
\end{cor}
\begin{proof} 
Let $*:\rp\ra\cH$ be given, and observe that if $\cE_*\ra\cU_*$ is a finite
pointed cover of a Zariski neighbourhood $\cU\supseteq \cH\ni *$,
then as $\cU$ decreases the cardinality of the fibre of
the connected component of $\cE_*$ eventually stabilises.
As such, 
and irrespective of the hypothesis on fibre dimension,
the proof, \ref{fact:plag1}, of the existence of
Galois objects goes through verbatim to show that
the full sub-category
\begin{equation}\label{LefCor1} 
\et_{1,\cH}^{\mathrm{Zar}}(*)\subseteq \et_{1,\cH}^{\mathrm{Zar}}
\end{equation}
whose objects are those every connected component of
which has a non-empty fibre over $*$
is equivalent to the category of finite sets on which
some pro-finite group $\pi_{1,\cH}^{\mathrm{Zar}}(*)$ acts.
Now by \ref{fact:Lef02}, 
$\et_{1,\cH}^{\mathrm{Zar}}$ is equally the direct sum
of the categories
$\et_{1,\cH}^{\mathrm{Zar}}(*)$ as $*$ runs through 
a set of base points in 1-1 correspondence with the
connected components, say $\cH_*$ in
a minor abuse of notation, of $\cH$, so, again
what has to be proved is that if $\cE\ra\cU$ 
is connected for all
sufficiently small $\cU$ 
Zariski neighbourhoods of $\cH$
then the fibre, $\cE_\cH$, over $\cH$ is
connected. Now everything is representable so,
\cite[16.5]{L-MB}, 
we can find some finite $\bar{\cE}\ra \cX$ containing $\cE$
as a dense Zariski open, and indeed independently of $\cU$
since $\bar{\cE}$ is basically just the integral closure
of $\cO_\cX$ in $\cO_\cE$. As such by \ref{fact:Lef02} there
is a Zariski open neighbourhood $\cV_{\mathrm{max}}$ of $\cE_{\cH}$ such
that 
\begin{equation}\label{l12new}
\pi_0(\cE_\cH)\xrightarrow{\sim}\pi_0(\cV)
\end{equation}
for all Zariski opens $\cV\subset \cV_{\mathrm{max}}$ 
containing $\cE_{\cH}$. In particular, therefore, 
it holds for 
$\cV=\cE\ts_{\cX} \cU$ and $\cU$ sufficiently small.
\end{proof} 

In so much as we're now concerned with $\pi_1$
rather than $\pi_0$ we'll require to pay more attention
to the difference between a champ and its
moduli, and will have need of
\begin{lem}\label{lem:l11}
Let everything be set up as in \ref{setup:l01};
$Y\ra\cX$ an effective
descent morphism for
the \'etale topology 
(e.g. a finite map \cite[Expos\'e IX.4.7]{sga1})
from an algebraic space;
\begin{equation}\label{l13}
\cdots Y^{(2)}:= Y\times_\cX Y\times_\cX Y {\build\ra_{\ra}^{\ra}}
Y^{(1)}:= Y\times_\cX Y\, (=R_0){\build\rras_{s}^{t}} Y^{(0)}:= Y
\end{equation}
the resulting simplicial space; and $\uG$ a locally
constant group 
(or more generally just a constructible sheaf of groups)
over $\cX$, then there is an exact
sequence of sets
\begin{equation}\label{l14}
1\ra \check{\rH}^1(\rH^0(Y^{(1)}, \uG)) \ra \rH^1(\cX,\uG)
\ra \check{\rH}^0(\rH^1(Y^{(0)}, \uG))  
\ra \check{\rH}^2(\rH^0(Y^{(2)}, \uZ))  
\end{equation}
where $\uZ$ is the centre of $\uG$.
\end{lem}
\begin{proof}
Plainly this is just the non-abelian case of
Deligne's descent spectral sequence, \cite[5.3.3]{deligne3},
as already encountered in, say, \eqref{eq:one7}. 
Furthermore, since $Y\ra\cX$ 
is an effective descent
morphism for the \'etale topology
we have a number of simplifications. Indeed,
on identifying
constructible sheaves with their espace \'etal\'e,
there is an identity between  constructible
sheaves on $\cX$, and constructible sheaves on
$Y$ with a descent datum. In particular, the
first term in \eqref{l14} is just a descent datum
for the trivial $\uG$ torsor over $Y$, which then
maps to the resulting $\uG$ torsor over $\cX$, which
in turn goes to the class of the same over $Y$.  
The remaining terms are a bit more interesting.
Specifically given a $\uG$-torsor, $E$, over $Y$
one constructs a groupoid, $e:R\ra R_0$, 
with fibre $\uG$,
whose 
objects are those of $Y$, but whose arrows are
given by
\begin{equation}\label{l15}
\Hom_R (x,y) \,:=\, \Hom_{\uG} (E_x, E_y)
\end{equation}
The further condition that $e$ has a section, $r$,
in spaces, but not necessarily a section of
groupoids, is that $E$ belongs to the pen-ultimate
group in  \eqref{l14}.
Now the stabiliser of $R/R_0$ is canonically
the automorphisms of $E$, which may well be
different from $\uG$ but it has the same centre, so
the fact that one has a section implies
(for more or less the same reason encountered
in group extensions, {\it i.e.} $Y=\rp$, otherwise
there would be a hypercovering issue)
that the equivalence classes of such extensions
of $R$ by $R_0$ are classified by the final
group in \eqref{l14}- {\it cf.} 
\eqref{3124} \& \cite[IV.3.5]{giraud}.
The class of $E$ in this final group is, therefore, trivial in
this final group iff $r$ can be taken to be 
a section of groupoids, a.k.a. defines a descent
datum for $E$. 
\end{proof} 
In consequence we have
\begin{cor}\label{cor:l12} 
Again, let everything be as in \ref{setup:l01}, with
$\uZ$  a constructible sheaf of sets (respectively groups,
respectively abelian groups) then for any 
fibre square of $S$-champs
\begin{equation}\label{l16}
\begin{CD}
\cX@<< g< \cX'\\
@VfVV @VV{f}V\\
\cS @<g<< \cS'
\end{CD}
\end{equation}
with the left vertical proper there are isomorphisms
of constructible sheaves
\begin{equation}\label{l17}
\text{$g^*(R^qf_*\uZ)\xrightarrow{\sim} R^qf_*(g^*\uZ)$, for
$q=0$, respectively $q=1$, respectively $q\geq 2$.}
\end{equation}
\end{cor}
\begin{proof} We've already used the $q=0$ case
several times by reduction to the moduli $\vert\cX\vert$.
The abelian case is \cite[18.5.1]{L-MB}, and in
our situation, {\it i.e.} $\cX$ Deligne-Mumford,
is independent of the error, \cite{olsson}, 
in \cite[12.2]{L-MB}. Indeed, whether in the abelian
or non-abelian case, the demonstration is the same:
{\it i.e.} without loss of generality $\cS=S$, $\cS'=S'$ are
Noetherian and affine, take $Y\ra\cX$ as in
\cite[16.6]{L-MB}, and reduce to usual proper
base change theorem whether by the descent spectral
sequence, \cite[5.3.5]{deligne3} in the abelian case,
or \eqref{l14} in the non-abelian case.
\end{proof}
Similarly, there is a smooth base change theorem
\begin{cor}\label{cor:SmoothBaseChange}
Again, let everything be as in \ref{setup:l01}, with
$\uZ$  a constructible sheaf of sets (respectively groups,
respectively abelian groups) then for any 
fibre square of $S$-champs
\begin{equation}\label{lll16}
\begin{CD}
\cX@<< g< \cX'\\
@VfVV @VV{f}V\\
\cS @<g<< \cS'
\end{CD}
\end{equation}
with the bottom right horizontal smooth there are isomorphisms
of constructible sheaves
\begin{equation}\label{lll17}
\text{$g^*(R^qf_*\uZ)\xrightarrow{\sim} R^qf_*(g^*\uZ)$, for
$q=0$, respectively $q=1$, respectively $q\geq 2$.}
\end{equation}
provided that for $q\geq 1$ every element of every
stalk of $\uZ$ has order prime to the residue 
characteristics of the points of $\cS$.
\end{cor}
\begin{proof} Exactly as above this reduces to
smooth base change of schemes.
\end{proof}
Notice also the pertinent corollary to the corollary \ref{cor:l12} is
\begin{cor}\label{cor:l13}
Let $\cX/V$ be a proper connected Deligne-Mumford champ over
the spectrum of a Noetherian local ring, with $0\hookrightarrow V$
the closed point, then the functor
\begin{equation}\label{l18}
\et_1(\cX) \ra \et_1 (\cX_0)
\end{equation}
is fully faithful, and an equivalence of categories
if $V$ is (not necessarily strictly) Henselian.
\end{cor}
\begin{proof} Just as in \ref{cor:l11}, fully faithfulness
follows from knowing that under the said hypothesis,
$\cX_0$ is connected, which reduces to the same for
the moduli in light of \eqref{l01}, which in turn is
\cite[XII.5.8]{sga4}. Similarly, since both sides are
equivalent to finite sets on which the respective
fundamental groups act, essential surjectivity follows if
every Galois cover of $\cX_0$ lifts to $\cX$. Necessarily,
however, a Galois cover is a torsor under some finite
group $\G$, so this is immediate from \eqref{l17} if
$V$ is strictly Henselian. To get the Henselian case
from the strict one, $(\bar{V}, \bar{0})$, say, observe
that the proof of \ref{lem:l11} equally gives an
exact sequence,
\begin{equation}\label{l19}
 1\ra {\rH}^1(G, \rH^0(\cX_{\bar{V}}, \G)) \ra \rH^1(\cX,\G)
\ra {\rH}^0(G, \rH^1(\cX_{\bar{V}}, \G))  
\ra {\rH}^2(G, \rH^0(\cX_{\bar{V}}, Z))  
\end{equation}
and similarly for $\cX_0$, where
$G$ is the Galois group of $\bar{0}/0$, and $Z$ the centre
of $\G$. 
\end{proof}
The final preliminary we need is
\begin{defn}\label{def:l11}
Let 
$(U,x)$ be the spectrum of the (Zariski) local ring
of a variety over a field $k$ (respectively the
Henselian local ring, respectively the strictly
Henselian local ring) then we say that $U$ has
homotopy depth $d$ 
(in practice at most 4)
at $x$ in the Zariski 
(respectively Nistnevich, respectively \'etale)
sense if
\begin{equation}\label{l110}
\pi_q(U\bsh x) \xrightarrow{\sim} \pi_q(U), \quad
\forall\,\, q + \mathrm{Trdeg}_k k(x) < d-1
\end{equation}
A $k$-scheme of finite type is said to have homotopy
depth $d$ in the respective sense according as all
of its  
local rings do. A $k$-champ, $\cX$, whose moduli is
a scheme, $X$, is said to have Zariski, respectively
Nistnevich, respectively {\it weakly} \'etale, 
homotopy depth $d$ at $x\in X$ if \eqref{l110}
holds on replacing whether $U$, 
or $U\bsh x$ by their pre-image in $\cX$, while
the depth is said to be $d$ in the \'etale sense
if \eqref{l110} holds in all the (by definition
strictly Henselian) local rings of $\cX$. 
\end{defn}
Even without the plethora of respectives there's
plenty of possibility for confusion
about the terminology, so it's
opportune to make
\begin{rmk}\label{rmk:l11}
The definitions \cite[Expos\'e XIII.4.3]{sga2},
or {\it op. cit.} Expos\'e XIV.1.2 are only designed
to work in the \'etale, 
respectively Nistnevich, case. As we've already seen,
however, \ref{fact:l02} only requires homotopy
depth 2 in the Zariski sense, and \ref{schol:l01},
this is weaker than the other senses of the term.
Worse, even in the \'etale sense what we've called
homotopy depth is called rectified homotopy depth
in \cite[Expos\'e XIII.4.3.D\'efinition 2]{sga2}-
{\it cf.} \ref{cschol:Lef11} for some examples/motivation
for dropping the word `rectified'-
while the condition \eqref{l110} for $q=0$ is, 
{\it op. cit.} Expos\'e XIII.2.1,
called {\it connected in dimension $d-1$}, where the role
of the flavours has been covered 
in \ref{schol:l01}.  A useful example for avoiding confusion is
\begin{ex}\label{ex:lci}
Let $X/k$ be everywhere a local complete intersection
of dimension at least $d$, then 
(modulo the practical difficulty that $q\leq 2$ in
\eqref{l110})
the homotopy depth is
everywhere \'etale locally at least $d$.
\end{ex}
\begin{proof}
The question is local, so, say $(U,x)$ the strict
Henselisation of some scheme point $x$ with
$\d(x)=\mathrm{Trdeg}_k k(x)$. 
Consequently if $\d(x)\leq d-2$, then the local
ring is $S_2$, so $U\bsh x$ is connected- albeit
{\it cf.} \ref{cschol:l11} for a much better result.
Similarly, if $\d(x)\leq d-3$, then the local
ring is a complete intersection of dimension at
least $3$, so this is \cite[Expos\'e X.3.4]{sga2}.
Finally consider the case $\d(x)\leq d-4$, 
{\it i.e.} a l.c.i. local ring of dimension at least $4$.
Given the previous cases, it suffices, {\it cf.}
\ref{cor:l21}-\ref{fact:Lef13} and \ref{claim:Lef21},
to prove that the local \'etale co-homology
$\rH^3_x(U,\uZ)$ vanishes for $\uZ$ a finite constant
sheaf of abelian groups. Plainly we can divide $\uZ$
up into its $p=\mathrm{char}(k)$-part and its prime
to $p$-part. To do the former, observe that the
local ring is $S_4$, so this
follows from 
the Artin-Schreier exact sequence, while the latter
is \cite[1.3]{IllusiePerversiteEtVariation}, and
one can usefully note that for the same reasons
it holds quite generally that
\begin{equation}\label{lci1}
\rH^q_x (U, \uZ)=0, \quad q<d-\d(x)=\mathrm{dim}\cO_{U,x}
\end{equation}
for any torsion sheaf of abelian groups.
\end{proof}
Rather amusingly, therefore, \ref{ex:lci} is still
true for $d=1$, provided that one interprets {\it
homotopy depth 1} to mean that \eqref{l110} is
surjective for $q=\d(x)=0$, so homotopy depth 1,
and dimension 1 coincide. In general, however, 
\ref{def:l11} is a dimension free definition, but
at the same time it only behaves well if the
ambient dimension is at least $d$ as the case
of $q=0$ illustrates, {\it i.e.} 
\begin{Cscholion}\label{cschol:l11}\ref{schol:l01} The point
in question is the behaviour of the condition 
{\it connected in dimension $d-1$} for a local, 
according to all possible flavours,
ring $(U,x)$ and a hypersurface $V:f=0\ni x$. The
best that one could hope for is
\begin{equation}\label{Csch1}
\text{$\mathrm{dim}_x U> d-1$ implies $V$
connected in dimension $d-2$}
\end{equation}
since a smooth (indeed irreducible) variety of dimension $d-1$ is
equally connected in dimension $d-1$, while a node can be
disconnected by a subvariety of dimension $d-3$. Similarly,
for stupid reasons, {\it i.e.} $2-3=-1$
a.k.a. \eqref{l110} becomes empty,
nodal curves in surfaces are connected in dimension $0$, so
to use the hypothesis in \eqref{Csch1} in order to establish
the same by induction, one has to start with the case of
$U$ of dimension $3$ rather than the empty $2$-dimensional case.
Incredibly, however, the ``proof'', \cite[Expos\'e XIII.2.1]{sga2},
of \eqref{Csch1} 
for complete local rings
makes exactly this mistake and is completely
wrong- the induction step is correct, but the initial case,
which (because of a slightly different set up) 
is even false in {\it op. cit.} rather than just empty for
surfaces, has to be done in dimension 3 rather than dimension
$2$. Fortunately this has been corrected (and it requires
a non-trivial trick) in \cite[3.1.7]{joinsAndIntersections},
where irrespectively of whether we're in the Zariski,
Nistnevich, Henselian or complete flavours we have 
\begin{equation}\label{Csch2}
\text{\eqref{Csch1} holds if $U$ is embeddable in a 
Gorenstein, Noetherian scheme.}
\end{equation}
even though the statement, rather then the proof, of
{\it op. cit.} is only given in the complete case. Indeed
even the Gorenstein embedding condition could be dropped
if one knew 
\begin{equation}\label{Csch3}
\text{$(U,x)$ normal, dimension $\geq3\Rightarrow$ 
$\rH^2_x(U,\cO_U)$ an $\cO_{U,x}$ module of finite
length.}
\end{equation}
Irrespectively, in all cases where we'll need it, we have
from \eqref{Csch2} that
\begin{equation}\label{Csch4}
\text{ $\mathrm{dim}_x (U)$ 
and homotopy depth $\geq 3\Rightarrow$
$\mathrm{dim}_x (V)$
and it's homotopy depth $\geq 2$,}
\end{equation}
and, as it happens, we only need the much weaker condition
that \eqref{Csch4} holds for a generic hypersurface $f=0$-
{\it cf.} \ref{cor:Lef11}-
with everything of finite type over a field, so, ironically,
 the correct part of \cite[Expos\'e XIII.2.1]{sga2}
would be good enough.
\end{Cscholion}
\end{rmk}
Having thus gone through connected in dimension 2, the
remaining condition in
our immediate interest, {\it i.e.} homotopy depth 3,
is
\begin{equation}\label{l111}
\pi_1(U\bsh x)\xrightarrow{\sim} \pi_1(U)
\end{equation}
for every closed point $x$, which simplifies
according to,
\begin{fact}\label{fact:l11}
Suppose that the homotopy depth is at least 2
in the \'etale sense then
\begin{equation}\label{l112}
\text{\eqref{l111} for \'etale $\Rightarrow$
\eqref{l111} for Nistnevich $\Rightarrow$
\eqref{l111} for Zariski }
\end{equation}
and similarly if $\cU$ is a champ with moduli
the Zariski local ring $U$
\begin{equation*}
\text{\eqref{l111} for \'etale $\Rightarrow$
\eqref{l111} for weakly \'etale $\Rightarrow$
\eqref{l111} for Nistnevich $\Rightarrow$
\eqref{l111} for Zariski }
\end{equation*}
wherein for the last 3 implications one should
understand \eqref{l111} with $\cU$, and $\cU':=\cU\ts_U U\bsh x$.
In particular, \'etale homotopy depth at least 3
implies 
the same in
all other possible senses whether for schemes 
or champs.
\end{fact}
\begin{proof} All the implications are essentially the
same. We do the \'etale $\Rightarrow$ Zariski case for
champs since it's both the most difficult and most 
relevant. To this end let $(V,\bar{x})$ be the strictly
Henselian neighbourhood of $\cU$;  $\cE'\ra \cU$
a representable Galois cover; and consider the fibre
square
\begin{equation}\label{l113}
\begin{CD}
\cE'@<<< E'\\
@VVV @VVV\\
\cU@<<< V\bsh\bar{x}
\end{CD}
\end{equation}
where, by hypothesis, $E'\xrightarrow{\sim}\pi_0(E')\ts V\bsh\bar{x}$.
As such $\cE'$ extends locally across the puncture, and
there is a myriad of ways to extend globally to an \'etale
cover of $\cU$, {\it e.g.} 
use \ref{cor:412} to suppose everything reduced,  
take the integral closure in $\cE'$, and use \cite[6.14.4]{EGAIV},
or, better use that the above local extension is uniquely
unique because the \'etale homotopy depth is at least 2
to get a descent datum for the extension over some \'etale
neighbourhood $W\ra\cU$. The latter strategy is arguably better
since it also gives the unicity of the extension.
\end{proof}
This pretty much covers everything one might
want to know about homotopy depth 3 beyond
\begin{rmk}\label{rmk:l12}
The relation between \eqref{l111} for $U$ 
(not necessarily strictly) Henselian and
its completion, $\hat{U}$, is plausibly 
more subtle than that encountered for 
$\pi_0$ in \ref{cor:hd1}. For exactly the
same reason as \ref{fact:l11}, the implication
\begin{equation}\label{l114}
\text{\eqref{l111} for $\hat{U}$ $\Rightarrow$
\eqref{l111} for $U$ 
}
\end{equation}
always holds, but the other way round is a priori
open in full generality, albeit \cite{ArtinLocalMonodromy}
covers a lot of cases.
\end{rmk} 
All of this said we have,
\begin{prop}\label{prop:Lef11}
Let everything be as in the set up \ref{setup:l01},
with 
$\cX/S$ proper enjoying 
for all $s\in S$
fibres, $\cX_s$,
everywhere of dimension 3 and connected
in dimension 2 in the \'etale sense, then for
$i:\cH\hookrightarrow \cX$  a hyperplane section 
there is an equivalence of
categories
\begin{equation}\label{Lef115}
i^*:\et_{1,\cH}^{\mathrm{Zar}} \xrightarrow{\sim} \et_1(\cH)
\end{equation}
\end{prop}
\begin{proof} 
We already know, \ref{cor:l1new}, that the
functor \eqref{Lef115} is fully faithful, and
we require to prove essential surjectivity,
{\it i.e.} lift a $\G$-torsor, $\cE\ra\cH$,
to a $\G$-torsor over a Zariski neighbourhood for any finite group $\G$.  
To this end, we
retake the notations of the proof of \ref{fact:Lef02};
profiting from \eqref{l01} on replacing
$\cH$ by a sufficiently large multiple, 
we add a further fibre square to the diagram \eqref{l0p2}
to obtain
\begin{equation}\label{l118}
\begin{CD}
p: \cF@>>> F@>>{p_X}>  P\\
@V{q}VV @VV{q_X}V  @.\\
\cX@>{\mu_\cX}>> X
\end{CD}
\end{equation}
As such if we again 
fix $s\in S$ and
identify $\cH_s$ with a $k(s)$-point 
$0$ of $P$; put $\cD=q^* \cH$; 
write fibres of $p$ as sub-scripts;
and denote by $N$
the strictly Henselian local neighbourhood of $0\in P$, then
a couple of applications of \eqref{cor:l13}
gives us isomorphisms
\begin{equation}\label{l119}
\pi_1(\cD_N)\xleftarrow{\sim} \pi_1(\cH) \xrightarrow{\sim} 
\pi_1(\cF_N)
\end{equation}
Consequently if $\cE\ra\cH$ is a $\G$-torsor for
some finite group $\G$, there is an \'etale
neighbourhood $V\ni 0$, and a $\G$-torsor $\cG$
over $\cF_V$ together with an isomorphism
\begin{equation}\label{l120} 
G:\cG\mid_{\cD_V}\xrightarrow{\sim} q^*\cE_V
\end{equation}
in $\et_1(\cF_V)$;
where despite the ambiguous
nature of fibre products in 2-categories 
every $S$-champ has a clivage by definition,
so $\cF_U$, $q^*\cE$ {\it etc.} are unambiguously defined,
and even $\cG\mid_{\cD_V}$ since we've replaced $\cH$
by a large multiple so we can use \eqref{l01}.
Now 
by 
\eqref{Csch4}
the universal family
$\cF/P$ satisfies (universally) the hypothesis of \ref{cor:l11},
so if
$U$ is the Zariski
image of $V$ in $P$; $(s,t):R:=V\ts_U V\rras V$ the resulting
groupoid then there is a unique arrow, $A$,
in $\et_1(\cF_{V\ts_U V})$ such that
\begin{equation}\label{l121}
\begin{CD}
s^* \cG\mid_{\cD_V} @>>{A\mid_{\cD_{V\ts_U V}}}> t^*\cG\mid_{\cD_V} \\
@V{s^*G}VV @VV{t^*G}V \\
s^*q^*\cE_V @=t^*q^*\cE_V
\end{CD}
\end{equation}
commutes. In particular, therefore,
since $\g A\g^{-1}$, $\g\in\G$ also satisfies \eqref{l121},
$A$ is a map of $\G$-torsors,
and since \ref{cor:l11}
is equally valid
over $V\ts_P V\ts_P V$, we may argue similarly to
conclude that $A$ satisfies the descent/1 co-cycle condition,
\begin{equation}\label{l117}
p^*_{31} A = p^*_{32} A p^*_{21} A
\end{equation}
in $\et_1(\cF_{V\ts_P V\ts_P V})$, wherein the indices
denote the various projections to $\cF_R$.
If $\cF_U$ were a space this concludes the discussion,
{\it i.e.} there is a
$\G$-torsor $\cG'\ra\cF_U$ whose restriction
to $\cF_V$ is $\cG$. In general, however, one
needs to be more careful about the lack of
unicity in fibre products which can result 
in a non-trivial 2 co-cycle condition, {\it cf.}
\eqref{PosSeq3} or \cite[IV.3.5.1]{giraud}.
Fortunately, however, 
\eqref{eq:cor1},
a natural transformation
in $\underline{\mathrm{Cham}}\mathrm{p}\underline{\mathrm{s}}$
which is a 2-cell in $\et_2$ between 1-cells
which themselves are maps between representable
0-cells is unique, so the said 2 co-cycle
condition is implied by \eqref{l117}, 
{\it cf.} \ref{5.11} and immediately post \eqref{3148}.
Consequently, we can
replace $V$ by the Zariski neighbourhood $U$
in \eqref{l120}.

At which point, as in \eqref{l0p44}, we switch  tact and notation
with $(s,t)$ now from $\cF_U\ts_\cX \cF_U\rras \cF_U$, to find an
isomorphism
\begin{equation}\label{Lef117}
B:s^*\cG \xrightarrow{\sim} t^*\cG
\end{equation}
restricting to the identity on all fibres over an \'etale
neighbourhood $V\ra U\ts_S U$ of $0\ts 0$. The restriction
of $B$ to the fibres of $\cD_U\ts_\cH \cD_U$ via
\eqref{l120} can equally be seen as the descent data
for $s^*q^*\cE_U=t^*q^*\cE_U$ for the covering of
a Zariski neighbourhood of $\cD_U\ts_\cH \cD_U$ induced
by $V$. On the other hand the hyperplane 
$\cD_U\ts_\cH \cD_U\hookrightarrow \cF_U\ts_\cX \cF_U$
satisfies exactly the conditions of \ref{fact:l01},
so arguing exactly as in \eqref{l0p44} but for the
locally constant sheaf $\underline{\Hom}(s^*\cG, t^*\cG)$
instead of $\underline{Z}$ {\it etc.}, we deduce that
$B$ is actually defined over a Zariski open neighbourhood
$W\hookrightarrow U\ts_S U$ of $0\ts 0$, which, again,
is a classical topology, so without loss of generality
$B$ is defined everywhere. As such if $q_{ji}$ are the
various projections of $\cF_U\ts_\cX \cF_U\ts_\cX \cF_U$
to $\cF_U\ts_\cX \cF_U$ there is a unique element $\g\in\uG$
such that
\begin{equation}\label{Lef118}
p_{31}^*B^\g= p_{32}^*B p_{21}^* B
\end{equation}
Necessarily, however, $\g$ vanishes on the (proper) fibre
over $0\ts 0\ts 0$, so it vanishes on all the fibres
over a Zariski neighbourhood of the same, whence
without loss of generality on $U\ts_S U\ts_S U$, {\it i.e.}
$B$ defines a descent datum for the smooth map
$q:\cF_U \ra \cX$.

The Zariski open image of $q$ contains a
Zariski open neighbourhood $\cU'$ of $\cH':=\cH\ts_S S'$
for $S'$ a Zariski open neighbourhood of our initial
$s$, so we've certainly lifted the
restriction, $\cE'=\cE\ts_\cH \cH'$, of our torsor
to a $\G$-torsor $\cG'$ over $\cU'$. In order to glue such torsors,
$\cG'$, $\cG''$ over Zariski
open neighbourhoods $S'$, $S''$, say, one 
proceeds exactly as in \eqref{l121}, {\it i.e.}
via the unique lifting of the given gluing
of $\cE'$ with $\cE''$. The existence of such 
a unique such isomorphism on Zariski neighbourhoods
of $\cH'\cap \cH''$ is assured by \ref{cor:l1new},
and, since every fibre is a fortiori everywhere
connected in dimension 1 in the \'etale sense, 
it extends, {\it cf.}
\ref{fact:l11},  to all of $\cU'\cap\cU''$. As such,
just as in \ref{fact:Lef02} there is even a unique
(depending on $\cE$)
maximal Zariski neighbourhood, $\cU_{\mathrm{max}}$, to
which
$\cE$ lifts. 
\end{proof}
In the presence of sufficient depth we therefore have
the following variations
\begin{cor}\label{prop:l11}
Let everything be as above, \ref{prop:Lef11},
and suppose further that
for all $s\in S$, the fibre
$\cX_s$ is everywhere of homotopy depth $3$
(in the Zariski sense over $k(s)$)
then the inclusion $i:\cH\hookrightarrow\cX$ 
of a hyperplane section
affords
an equivalence of categories
\begin{equation}\label{l115}
i^*:\et_1(\cX) \xrightarrow{\sim} \et_1(\cH)
\end{equation}
\end{cor}
\begin{proof}
Again, we require to
lift a $\G$-torsor, $\cE\ra\cH$,
to a $\G$-torsor over $\cX$ for any finite group $\G$,
and we already know, \ref{prop:Lef11},
or more accurately the end of the proof,
that this can be done over a unique 
maximal Zariski neighbourhood $\cU$ of $\cH$.
Now let $\cG$ be the lifting and $s\in S$, then
by hypothesis $\cG_s$ extends over the whole
fibre $\cX_s$, so by proper base change, \ref{cor:l13}, there
is a Nistnevich neighbourhood $V\ra S$ of $s$
such that $\cG_V$ extends over all of $\cX_V$
to some $\G$-torsor $\bar{\cG}$.
Each fibre of $\bar{\cG}$ is, however, a fortiori
connected in dimension 1 in the \'etale sense,
so 
for $(s,t):V\ts_S V\ra V$
the descent isomorphism between $s^*\cG_V$ and $t^*\cG_V$
extends uniquely to a descent isomorphism between
$s^*\bar{\cG}$ and $t^*\bar{\cG}$.
As such, the fibre of $\cU_{\mathrm{max}}$ contains
$\cX_s$, and $s$ was arbitrary.
\end{proof}
\begin{cor}\label{cor:Lef11}
Let everything be as in \ref{prop:Lef11},
and suppose further that
$\cX/S$ itself is the universal family
of sufficiently ample hyperplanes,
i.e. $\cF_U/U$, $U\hookrightarrow P$
Zariski open,  in the notation of \eqref{l118},
of some $\cY/B$ such that
for all $b\in B$, the fibre
$\cY_b$ is everywhere locally simply connected
(in the Zariski sense over $k(b)$) in dimension 2, i.e.
\begin{equation}\label{CorLef111}
\pi_1(U\bsh y)\xrightarrow{\sim} \pi_1(U),\quad
\mathrm{Trdeg}_{k(b)}k(y) < 2
\end{equation}
for  $U$ the Zariski local ring at $y$,
then the inclusion $i:\cH\hookrightarrow\cX$ 
of a hyperplane section
affords
an equivalence of categories
\begin{equation}\label{CorLef115}
i^*:\et_1(\cX) \xrightarrow{\sim} \et_1(\cH)
\end{equation}
\end{cor}
\begin{proof} Let everything be as in the proof
of \ref{prop:l11}, then as before it suffices
to extend the torsor $\cG_s$ over $\cX_s$,
{\it i.e.} in a neighbourhood of $\cZ_s$, 
where $\cZ$ is the Zariski closed complement
of $\cU$. Now the universal family of hyperplanes
has the particular property that 
for $z\in \cZ_s$ one can find
a constructible $T\ni s$ (in fact a 
Zariski open in a smooth curve over $B$)
such that $\cX_T\ra\cY$ is a Zariski local
isomorphism at $z$ and $\cZ_T$ is of relative
dimension at most 1, and whence we can apply 
\eqref{CorLef111}.
\end{proof}
In all of which we can usefully observe that the
systematic use of proper base change has
a certain optimality, {\it i.e.}
\begin{rmk}\label{Notube}
There is no tubular neighbourhood theorem for
the \'etale topology. Specifically the usual
example for the necessity of the dimension
condition in \ref{prop:l11} is to take $H$
to be a smooth curve of positive genus in
$X=\bp^2_\bc$. Of course, \ref{fact:411}, 
the completion of $X$ in $H$ has the same
homotopy type as $H$. Given, however, a finite 
\'etale cover $H'\ra H$ it is not true that
there exists a (not necessarily proper over
its image) \'etale neighbourhood $p:V\ra H$
together with an embedding $H'\hookrightarrow V$.
Indeed such a $V$ would be normal, and without
loss of generality connected, whence it's irreducible,
so that we can just compactify everything to
$\bar{p}:\bar{V}\ra X$, $\bar{V}$ normal. Now $H'$ is a connected
component of $p^{-1}(H)$, so it's equally a 
connected component of $\bar{p}^{-1}(H)$, but
the latter is connected by \ref{fact:l02}. Consequently,
$H'=\bar{p}^{-1}(H)$, so $\bar{p}$ is an \'etale
covering in a neighbourhood of $H$,  and whence
everywhere by purity, which is absurd.
\end{rmk}
Similarly let us consider the necessity
of the depth conditions by way of,
\begin{Cscholion}\label{cschol:Lef11}\ref{schol:l01}/\ref{cschol:l11}
The necessity of the depth conditions \eqref{l110}
at closed points is equally obvious- {\it e.g.} join two
projective spaces in a point, take a variety with
an isolated quotient singularity, or, whatever. In
the situation of algebraic depth, however, the
analogous
conditions to \eqref{l110} are referred to as rectified
depth, \cite[10.8.1]{egaIV2}, and is wholly implied,
{\it op. cit.} 10.8.6, by the corresponding condition
at closed points. This is, however, false for
homotopy depth (whence we've
dropped `rectified' of 
\cite[Expos\'e XIII.6, D\'efinition 2]{sga2}) {\it i.e.}
\begin{equation}\label{CschLef1}
\text{\eqref{l110} at closed points for $\pi_0$
and $\pi_1$ does not imply connected in dimension 2}
\end{equation}
and, as it happens, the same example will equally prove
\begin{equation}\label{CschLef2}
\text{Connected in dimension 2 is necessary in the
Lefschetz theorem, \ref{prop:Lef11}, for $\pi_1$.}
\end{equation}
Specifically, therefore, let $X/\bc$ be the join
of two copies $A$, $B$ of $\bp^3_\bc$ in a line
$L$. This is manifestly connected in dimension 1
in the \'etale sense but not connected in dimension
2 in the Zariski sense. Similarly off $L$, \eqref{l111},
is trivially true at closed points because
$5>1$, and only slightly
less trivial ($5>1$ plus Van Kampen) on $L$ itself.
This shows \eqref{CschLef1}, and, again Van Kampen
shows that $X$ is simply connected. Now
take $a$, $b$ to be smooth quadrics in $A$,
respectively $B$, meeting the line $L$ in the same
pair, $p$, $q$, of distinct points, then 
({\it e.g.} Mayer-Vietoris) the
join, $H$, of $a$ to $b$ in $\{p\}\cup\{q\}$ isn't
simply connected. If, however, a cover $E\ra H$
could be extended to a Zariski neighbourhood of $H$
then it would extend everywhere because the 
depth conditions hold at closed points, and
whence \eqref{CschLef2}. 
\end{Cscholion}
\subsection{Lefschetz for \texorpdfstring{$\Pi_2$}{Pi\_2}}\label{SS:L2}
At which point the general inductive schema should
be clear since we now have
\begin{cor}\label{cor:l21}
Let everything be as in \ref{prop:l11}, respectively
\ref{cor:Lef11}, 
albeit with the fibres of $\cX$ of homotopy depth 3, respectively
those of $\cY$ simply connected in dimension 2, in the \'etale
rather than the Zariski sense,
then the 2-functor
\begin{equation}\label{l21}
i^*:\et_2(\cX)\ra\et_2(\cH)
\end{equation}
is fully faithful, {\it i.e.} for any pair of $0$-cells
$q,q'$ on the left, 
\begin{equation}\label{l122}
i^*:\Hom(q,q')\ra\Hom(i^* q, i^* q')
\end{equation}
is an equivalence of categories.
\end{cor}
\begin{proof} By the 2-Galois correspondence \ref{prop:376},
the respective sides of \eqref{l21} are equivalent to the 2-category
of representations of $\Pi_2(\cX_*)$, respectively $\Pi_2(\cH_*)$,
while 
under the conditions of \ref{prop:l11} 
or \ref{cor:Lef11}
the  
fundamental groups are isomorphic. As such the assertion
is equivalent to
\begin{equation}\label{l23}
\pi_2(\cH_*)\twoheadrightarrow \pi_2(\cX_*)
\end{equation}
Indeed,  without loss of generality, the representation
on $q$, $q'$ are transitive in the sense immediately following
\eqref{eq:groupPlus3}, which in turn can be replaced,
\ref{fact:groupPlus1}, by their pointed stabiliser 
representations. The only part of the explicit
description, \ref{sum:group1}, of the same 
which depends on $\pi_2$ among the data defining
$q$, respectively $q'$, is the representation of
$\pi_2$ in the centre of the stabiliser of $q$,
respectively $q'$, and the only part of the data in
the description of a 1-cell $q\ra q'$ which depends
on $\pi_2$ is that these representations should map
to themselves under the implied map from the stabiliser
of $q$ to $q'$. Consequently, the easy direction is
that \eqref{l23} implies that \eqref{l122} is an
equivalence of categories. Conversely, in the 
notation of \ref{sum:group1}, given a finite quotient
$\G$ of $\pi_2(\cX_*)$, one should take $q$ to be
defined by the triple
$(\underline{\pi}', \underline{\g'}, \underline{A}')$,
where $\underline{\pi}'=(\pi'_1, \G, k'_3)$ for some
sufficiently small open sub-group $\pi'_1$ acting
trivially on $\G$; $\underline{\g'}$ the topological
2-type of the universal fibration in $K(\G,1)$'s; and
$\underline{A'}$ the natural forgetful map between them.
Similarly, given a character $\chi$ of $\G$, one takes
$q'$ to be defined by a triple 
$(\underline{\pi}'', \underline{\g}'', \underline{A}'')$,
where $\underline{\pi}''$ is the previous $\underline{\pi}'$;
$\underline{\g}''$ the topological
2-type of the universal fibration in $K(\mathrm{Im}(\chi),1)$'s;
and again $\underline{A}''$ the natural forgetful map.
As such 
$\chi$ affords a
natural map from $q\ra q'$ and
the isomorphism class of objects of the category
$\Hom(q, q')$ are a principal homogeneous space
under $\rH^1(\pi_1', \mathrm{Im}(\chi))$. If, however,
$\D$ is the image of $\pi_2(\cH_*)$ in $\G$ then 
the orbit under $\rH^1(\pi_1', \mathrm{Im}(\chi))$
of the isomorphism classes in $\Hom(i^* q, i^* q')$ 
is isomorphic to characters $F$ such that,
\begin{equation}\label{l24}
\Hom(\G, \bq/\bz)\ni F \mpo \chi\vert_\D \in \Hom(\D, \bq,\bz)
\end{equation}
and whence \eqref{l122} an equivalence implies
\eqref{l23} surjective.

This said, it therefore suffices by \ref{factdef:plagAct1}
to prove that if $q$ is a 2-Galois cell over $\cX$, then
$i^* q$ is 2-Galois over $\cH$, and, even slightly more
conveniently, \ref{fact:noplag3} \& \ref{fact:noplag5noplag5},
we can replace 2-Galois by quasi-minimal. To this end 
let $\cY\xrightarrow{p} \cY'\xrightarrow{r} \cX$ be a
factorisation of $q$ into a locally constant gerbe followed
by a representable cover, then we have a commuting diagram
of functors
\begin{equation}\label{l25}
\begin{CD}
\et_1(\cY') @>>{p^*}> \et_1(\cY)\\
@V{i^*}VV @VV{i^*}V\\
\et_1(i^*\cY') @>{(p\vert_\cH)^*}>> \et_1(i^* \cY)
\end{CD}
\end{equation}
where by the very definition of quasi-minimal the top horizontal
is an equivalence of categories, while the verticals are
equivalences by   \ref{prop:l11} or \ref{cor:Lef11}, so the bottom horizontal
is an equivalence.
\end{proof}
Alternatively, we have the weaker but more generally valid
\begin{cor}\label{cor:Lef21}
Suppose only
that everything is as in \ref{prop:Lef11}
and define $\et^{\mathrm{Zar}}_{2,\cH}$ in
the obvious way, {\it i.e.} $0$-cells are champs
which are proper and \'etale over a Zariski
neighbourhood of $\cH$ etc.,
then restriction affords a fully faithful 2-functor
\begin{equation}\label{ll21}
i^*:\et_{2,\cH}^{\mathrm{Zar}}\ra\et_2(\cH)
\end{equation}
\end{cor}
\begin{proof} As per \ref{prop:Lef11} for the 1-Galois
correspondence in Zariski neighbourhoods of $\cH$, the
2-Galois correspondence is valid in the obvious way, so
everything is exactly as in the proof of \ref{cor:l21}
modulo replacing \ref{prop:l11}, respectively
\ref{cor:Lef11}, by 
\ref{prop:Lef11}
\end{proof}
While the 
rationality issues encountered in
\eqref{l0p1}, \ref{prop:l11}
won't completely disappear, \ref{cor:LLef21}
\& \ref{rmk:LLef21},
we can usefully observe
\begin{fact}\label{lem:l21}
Let everything be as in the set up \ref{setup:l01} but
suppose $S$ is a field $k$ with separable closure
$i:k\hookrightarrow \bar{k}$ then 
$i$ 
and a choice of $\rp\ra \mathrm{Spec}(\bar{k})$ 
determines a unique (up to equivalence) lifting 
(up to equivalence) of the base
point $*:\rp\ra\cX$ and on understanding
$\pi_2(\cX_*)$ as $\pi_2$ of the connected
component of the base point,
\begin{equation}\label{l26}
\pi_2(\cX_{\bar{k}*})\xrightarrow{\sim}
\pi_2(\cX_*)
\end{equation}
\end{fact}
\begin{proof} Lifting $*$ is just the definition of
fibre products \eqref{FibreA3}. Similarly if,  again,
$\cY\xrightarrow{p} \cY'\xrightarrow{r} \cX$  is the 
factorisation into a local constant gerbe 
followed by a representable cover of a quasi-minimal
$0$-cell $q$, then, more or less by definition, every
base change of $p$ by a representable \'etale cover 
of $\cY'$ is again quasi-minimal. Of course for a 
finite separable extension $K/k$ identified with a
subfield of $\bar{k}$ it could  happen that
$\cY'_K$ is disconnected, but it is, nevertheless, a
finite sum of connected representable covers of $\cY'$,
amongst which pointing $\cY'$ a priori picks out exactly 1,
say $\cY''$, so that $\cY\ts_{\cY'}\cY''$ over the
corresponding connected component of $\cX_K$ is certainly
quasi-minimal, and whence \eqref{l26} is surjective.
Finally, for any such $K$, $\et_2(\cX_K)$ is a full
sub-2-category of $\et_2(\cX)$ so \eqref{l26} is injective.
\end{proof}
In a similar vein we can refine
\ref{cor:l13}, {\it i.e.}
\begin{fact}\label{fact:Lef13}
Let $\cX/V$ be a proper Deligne-Mumford champ over
the spectrum of a Henselian local ring, with $0\hookrightarrow V$
the closed point, then there is a 2-equivalence
\begin{equation}\label{llef18}
\et_2(\cX) \xrightarrow{\sim} \et_2 (\cX_0)
\end{equation}
\end{fact}
\begin{proof}
By the Whitehead theorem \ref{Whitehead}, or more
accurately a minor variant of it
as in the proof of \ref{cor:l1new}
to cope with non-trivial $\pi_0$, it suffices to
prove that the homotopy groups of either side of
\eqref{Lef70} coincide, where, without loss of
generality $\cX_V$ and $\cX_0$ are connected. 
Arguing however as above, \eqref{l25}, 
we already know from \ref{cor:l13} that we
have an isomorphism on $\pi_1$, and a surjection on $\pi_2$.
It therefore suffices to lift a 2-Galois cell $q$ over
$\cX_0$ to the same over $\cX_V$. 
On the other hand if
$\cE\xrightarrow{q'}\cX_0'\xrightarrow{p}\cX_0$ 
is the factorisation of $q$ into a locally
constant gerbe in $\rB_Z$'s, and a representable cover
for some locally
constant sheaf $\uZ/\cX_0$, then since $p$ can
be lifted to some \'etale cover $\cX'_V\ra\cX_V$, we 
can on replacing $\cX_V$ by $\cX'_V$,
suppose that the 2-Galois cell is actually a locally constant
gerbe. Similarly, another application of \ref{cor:l13}
implies that $\uZ$ is the restriction of a locally
constant sheaf on $\cX_V$, so we're done by proper base
change, \ref{cor:l12}, and
\ref{fact:384}.
\end{proof}

Irrespectively, the main proposition should,
at this juncture, be clear 
\begin{prop}\label{prop:Lef21}
Let everything be as in the set up \ref{setup:l01},
with 
$\cX/S$ proper enjoying 
for all $s\in S$
fibres, $\cX_s$,
everywhere of dimension 4 while being connected
in dimension 3, and
simply connected in dimension 2 (both in the \'etale sense), then for
$i:\cH\hookrightarrow \cX$  a hyperplane section 
there is an equivalence of
categories
\begin{equation}\label{Lef70}
i^*:\et_{2,\cH}^{\mathrm{Zar}} \xrightarrow{\sim} \et_2(\cH)
\end{equation}
\end{prop}
To this end we assert
\begin{claim}\label{claim:Lef21}
Denoting by $\uZ$ the constant sheaf with values in
the finite abelian group $Z$ it will suffice, 
under exactly the hypothesis of \ref{prop:Lef21} to
prove
\begin{equation}\label{Lef71}
\varinjlim_{\cU\supset \cH} \rH^2(\cU, \uZ)
\twoheadrightarrow \rH^2(\cH,\uZ)
\end{equation} 
where the limit is taken over Zariski open
neighbourhoods of $\cH$.
\end{claim}
\begin{proof}
Arguing as in the proof of \ref{fact:Lef13}
proves the assertion
under the weaker hypothesis that $\uZ$ is an arbitrary
locally constant sheaf in finite abelian groups. By
definition, however, as a group $\pi_2$ is unchanged
by \'etale coverings so by lifting a further representable
cover of $\cH$ to $\cX$ we can suppose that $\uZ$ is
functions to $Z$.
\end{proof}
Now let us return to
\begin{proof}[proof of \ref{prop:Lef21}] Fix a constant
sheaf, $\uZ$, of finite abelian groups; 
$e\in\rH^2(\cH,\uZ)$; $s\in S$; retake the
notation of \eqref{l118} for the universal hyperplane; and,
again, identify $\cH_s$ with the $k(s)$-point $0\in P$.
As such, following the 
notation and proof of \ref{prop:Lef11} in the
obvious manner there is an \'etale neighbourhood $V$ of
$0$ and a class
\begin{equation}\label{Lef72}
\rH^2(\cF_V,\uZ)\ni g \mpo q^*e_V \in \rH^2(\cD_V, \uZ)
\end{equation}
Now for $U$ the Zariski open image of $V$ in $P$ we
have a simplicial space
\begin{equation}\label{Lef73}
\cdots V^{(2)}:= V\times_U V\times_U V {\build\ra_{\ra}^{\ra}}
Y^{(1)}:= V\times_U V \, {\build\rras_{s}^{t}} V^{(0)}:= V
\end{equation}
whose pull-back to $\cF$ affords the descent
spectral sequence
\begin{equation}\label{Lef74}
E_2^{i,j}(\cF_U, \uZ) = \check{\rH}^i ( \rH^j(\cF_{V^{(i)}}, \uZ)
\Rightarrow \rH^{i+j} (\cF_U, \uZ)
\end{equation}
and similarly for the co-homology of the 
hyperplane section $\cD_U\hookrightarrow \cF_U$,
where by our current  hypothesis of
\ref{prop:Lef21} and \ref{cor:Lef11}  
\begin{equation}\label{Lef75}
\rE_2^{i,j}(\cF_U, \uZ)\xrightarrow{\sim}
\rE_2^{i,j}(\cD_U, \uZ),\quad j=0,\,\,\text{or} 1 
\end{equation}
while the $\rE_2^{0,2}$ term mimics \eqref{l0p4},
{\it i.e.} a commutative diagram
\begin{equation}\label{Lef76}
\begin{CD}
0@>>> \rE_2^{0,2} (\cD_U, \uZ) @>>> \rH^2 (\cD_V, \uZ) @>>{s^*-t^*}> 
\rH^2(\cD_{V\ts_U V},\uZ)\\
@. @AAA @AAA @AAA \\
0@>>> \rE_2^{0,2} (\cF_U, \uZ) @>>> 
\rH^2 (\cF_V, \uZ) @>{s^*-t^*}>> \rH^2(\cF_{V\ts_U V},\uZ)
\end{CD}
\end{equation}
Here our hypothesis allow us to apply the respective
part of \ref{cor:l21}, so that by \eqref{3129} all
the verticals in \eqref{Lef76} are injective. 
Consequently $g$ of \eqref{Lef72}, can be identified
with a class in the bottom left group mapping to 
the image of $q^*e_U$ in the top left, and a simple
diagram chase using \eqref{Lef75} shows that we
can replace $V$ with the Zariski open $U$ in \eqref{Lef72}.

At this point we 
change notation in order
to compare the descent spectral
sequences,
say
\begin{equation}\label{Lef77}
\rH^{i+j}(\cU',\uZ) \Leftarrow
\rE^{i,j}_r(\cX)\xrightarrow{\e_r^{i.j}}\rE^{i,j}_r(\cH)
\Rightarrow \rH^{i+j}(\cH, \uZ)
\end{equation}
for $q:\cF_U\ra \cU'$- $\cU'\supset \cH$ the Zariski
open image- and $\cD_U\ra \cH$. Now as we've seen
in \eqref{l0p44}, \eqref{Lef117}, and \eqref{Lef118}
the trick is to shrink $U\ni 0$, or products thereof, as
necessary in order to extract some extra surjectivity
or injectivity that isn't quite guaranteed by the
Lefschetz theorems that we already know. In order therefore
to lighten the exposition from here to the end of this
comparison of spectral sequences
\begin{equation}\label{Lef78}
\text{$P\supseteq U\ni 0$ means a sufficiently small
Zariski neighbourhood}
\end{equation}
Similarly, there will  appear intermediate \'etale
neighbourhoods, where the words `surjective', `injective',
or isomorphism,
may only be valid after Henselisation, but since
we're performing a finite diagram change there will be
no loss of generality in confusing the notions
`sufficiently small
\'etale neighbourhood' with Henselisation. With this
in mind what we've so far established is that $\e_1^{0,2}$
is an isomorphism, so the obvious next step is
\begin{claim}\label{claim:Lef22}
Notation and conventions as above, $\e_2^{0,2}$ is
an isomorphism, {\it i.e.} $\e_1^{1,2}$ is injective.
\end{claim}
\begin{proof}[sub-proof]
Let $V\ra U\ts_S U$ be a sufficiently small \'etale
neighbourhood of $0\ts 0$ and consider the resulting map of
descent spectral sequences
\begin{equation}\label{Lef79}
\rH^{i+j}(\cF_{U}\ts_\cX\cF_U,\uZ) \Leftarrow
\rF^{i,j}_r(\cX)\xrightarrow{\phi_r^{i.j}}\rF^{i,j}_r(\cH)
\Rightarrow \rH^{i+j}(\cD_U\ts_\cH\cD_U, \uZ)
\end{equation}
then the following facts about $\phi_1^{i,j}$
formally imply the assertion.
\begin{itemize}
\item[(a)] If $i=0$, it's an isomorphism $\forall j$-
\ref{cor:l12}.
\item[(b)] If $j=0$, it's an isomorphism $\forall i$-
\ref{fact:l02}.
\item[(c)] If $j=1$, it's injective $\forall i$-
\ref{cor:l11}.\qedhere
\end{itemize}
\end{proof}
\begin{claim}\label{claim:Lef23}
Notation and conventions as above, $\e_3^{0,2}$ is
an isomorphism, {\it i.e.} $\e_2^{2,1}$ is injective.
\end{claim}
\begin{proof}[sub-proof] This breaks up into 
\begin{equation}\label{Lef80}
\text{$\e_1^{1,1}$
is an isomorphism,} 
\end{equation}
which as in the
proof of \ref{claim:Lef22} is a formal consequence
of (a)-(c) of {\it op. cit.}, and
\begin{subclaim}\label{subclaim:lef24}
$\e_1^{2,1}$ is injective
\end{subclaim}
\begin{proof}[sub-sub-proof]
Let $V\ra U\ts_S U\ts_S U$ be a sufficiently small \'etale
neighbourhood of $0\ts 0\ts 0$ and consider the resulting map of
descent spectral sequences
\begin{equation}\label{Lef81}
\begin{split}
& \rH^{i+j}(\cF_{U}\ts_\cX\cF_U\ts_\cX\cF_U,\uZ) \Leftarrow
\rF^{i,j}_r(\cX)\xrightarrow{\phi_r^{i.j}}\\
&\rF^{i,j}_r(\cH)
\Rightarrow \rH^{i+j}(\cD_U\ts_\cH\cD_U\ts_\cH \cD_U, \uZ)
\end{split}
\end{equation}
then the following facts about $\phi_1^{i,j}$
formally imply the assertion.
\begin{itemize}
\item[(a)] If $i=0$, it's an isomorphism $\forall j$-
\ref{cor:l12}.
\item[(b)] If $j=0$, it's injective $\forall i$-
\ref{fact:l01}.\qedhere
\end{itemize}
\end{proof}
Which in turn completes the proof of \ref{claim:Lef23}
\end{proof}
\begin{claim}\label{claim:Lef25}
Notation and conventions as above, $\e_4^{0,2}$ is
an isomorphism, {\it i.e.} $\e_3^{3,0}$ is injective.
\end{claim}
\begin{proof}[sub-proof] Given \eqref{Lef80} this amounts to
\begin{equation}\label{Lef83}
\text{$\e_2^{3,0}$ is injective}
\end{equation}
Now, bearing in mind \eqref{Lef78}, 
by the obvious variant of 
\eqref{l0p44}, $\e_1^{p,0}$ is injective for all $p$,
so the assertion follows a fortiori from
\begin{equation}\label{Lef84}
\text{$\e_1^{2,0}$ is an isomorphism}
\end{equation}
which in turn is a formal consequence of (a) and (b)
in the proof of \ref{subclaim:lef24}.
\end{proof}
A priori such tedia only give an isomorphism 
(Zariski locally around $s\in S$) of the
$\mathrm{gr}^0$ parts of either side of \eqref{Lef71}
implied by \eqref{Lef77}, but
\begin{equation}\label{Lef85}
\begin{split}
\text{$\mathrm{gr}^1$ is an isomorphism because
$\e_2^{1,1}$ is- \eqref{Lef80}, \ref{subclaim:lef24},
\& \ref{cor:Lef11}- and \eqref{Lef83}}\\
\text{$\mathrm{gr}^2$'s an iso since $\e_2^{2,0}$
-\eqref{Lef84}, \ref{fact:l02}, \& $\e_1^{3,0}$ inj-
and $\e_2^{1,0}$ are- \ref{cor:Lef11} \& \ref{cor:l11}.}
\end{split}
\end{equation}

As such, we've solved the problem everywhere Zariski
locally on $S$, and we have to patch liftings $e'$, $e''$
of $e$ to Zariski neighbourhoods $\cU'$, $\cU''$ of
$\cH'=\cH\ts_S S'$, $\cH''=\cH\ts_S S''$ for $S', S''$
themselves Zariski open open in $S$. Now by \ref{cor:Lef21}
$e'$ and $e''$ restrict to the same class on some 
intermediary Zariski neighbourhood, $\cU'''$, 
of 
$\cH'\cap\cH''$, whose complement in 
$\cU'\cap\cU''$ we denote by $\cZ$. What we therefore
require to prove is that our depth hypothesis imply,
\begin{equation}\label{Lef86}
\rH^2_{\cZ}(\cX\ts_S S'\cap S'', \uZ)=0
\end{equation}
We already know, however, the corresponding proposition
for $\rH^i$, $i\leq 1$, {\it cf.} \ref{fact:l11}, and
even in a strong local sense since globally disconnected
implies locally disconnected for generic $s\in S$. 
Consequently, it will, by the local global spectral
sequence for co-homology with support suffice to prove
\eqref{Lef86} with $S'\cap S''$ replaced by the 
Nistnevich local neighbourhood $N\ra S'\cap S''$ of
an arbitrary $s\in S'\cap S''$. We have, however, 
a
commutative diagram
\begin{equation}\label{Lef87}
\begin{CD}
\rH^1(\cX_N,\uZ)
@>>> \rH^1 (\cX_N\bsh\cZ, \uZ) 
@>>> \rH^2_\cZ (\cX_N, \uZ) @>>> 
\rH^2(\cX_N,\uZ)\\
@VVV @VVV @VVV @VVV \\
\rH^1(\cX_s,\uZ)
@>>> \rH^1 (\cX_s\bsh\cZ, \uZ) 
@>>> \rH^2_\cZ (\cX_s, \uZ) @>>> 
\rH^2(\cX_s,\uZ)
\end{CD}
\end{equation}
with exact rows wherein the leftmost horizontals
are isomorphisms by \ref{prop:l11} (or, slightly
more accurately its proof), whence the rightmost
horizontals are injections, so that, by Nistnevich
proper base change for the ultimate vertical, the
vertical between the local co-homology groups is
injective, while the vanishing of the bottom $\rH^2_\cZ$
follows a fortiori from simply connected fibres
in dimension 1 in the \'etale sense. As such   
we get \eqref{Lef86}, and, just as in \ref{fact:Lef02},
we even find a maximal 
(depending on $e$)
Zariski open neighbourhood of $\cH$ to which the
corresponding gerbe $\cE\ra\cX$ may be extended
as a locally constant fibration in $\rB_Z$'s.
\end{proof}
In the presence of further
depth hypothesis we therefore obtain
\begin{cor}\label{cor:LLef21}
Let everything be as above, \ref{prop:Lef21},
and suppose further that
for all $s\in S$, the fibre
$\cX_s$ is everywhere of homotopy depth $4$
(in the Zariski sense over $k(s)$)
then the inclusion $i:\cH\hookrightarrow\cX$ 
of a hyperplane section
affords
an equivalence of 2-categories
\begin{equation}\label{Lef2115}
i^*:\et_2(\cX) \xrightarrow{\sim} \et_2(\cH)
\end{equation}
\end{cor}
\begin{proof}
By the Whitehead theorem, \ref{Whitehead},
and \ref{cor:l21}, it suffices to lift a
2-Galois cell $q:\cE\ra\cH$ to all of $\cX$.
As such let $\cE\xrightarrow{q'}\cH'\xrightarrow{p}\cH$
be its factorisation into a locally constant gerbe
in $\rB_Z$'s for some locally constant (but not
necessarily constant) sheaf of abelian groups
$\uZ$ followed by a representable cover $p$. 
An immediate application of \ref{prop:l11} implies
that $\uZ$ is defined on all of $\cX$ and $p$
extends to $\cX'\ra\cX$. 
Now by \ref{prop:Lef21}, $q'$ extends to a Zariski
neighbourhood $\cU'$ of $\cH':=p^*\cH$ as say,
$q':\cE'\ra\cU'$, and better, by the conclusion
of the proof of {\it op. cit.} we can take $\cU'$
to be maximal with respect to the properties that
the extension $q'$ exists and $p\mid_{\cU'}$ is 
proper over its image. Plainly we're done if 
the Zariski closed complement of $p(\cU')$ is empty,
so, suppose otherwise and let $z$ be a generic
point of the same lying over $s\in S$ with $N\ni s$
the Nistnevich local neighbourhood of the same.
By the definition, \ref{def:l11}, of Zariski
homotopy depth and the Whitehead theorem
(or more correctly another minor variant
since there's need for a little bit of
care about the meaning of a base point,
{\it cf.} \cite[Expos\'e XIII, pg. 15]{sga2})
there
is a local extension, $f_s:\cF_s\ra \cV_s$ of $p_sq'_s$
to a Zariski open neighbourhood of $z$ in the
fibre $\cX_s$. On the  other hand, $f$ has it's
own factorisation say,
$\cF_s\xrightarrow{f'_s}\cV'_s\xrightarrow{f_s}\cV_s$,
into a locally constant gerbe followed by a
representable cover, which, in turn is unique,
so $\cV'_s$ can be identified with a Zariski
open neighbourhood of $\cX'_s$, and $f'_s$ with
an extension of the locally constant gerbe $q'_s$
over a Zariski neighbourhood of $p_s^{-1}(z)$.
Arguing similarly at any other closed points
in the complement of $\cU'_s$ implies by 
Mayer-Vietoris that the locally constant
gerbe defined by $q'$ extends to all of $\cX'_s$,
and whence to $\cX_N$ by Nistnevich proper base
change. Such an extension is, however, unique
by \ref{cor:l21} and \ref{3129} so it glues
exactly as in \eqref{Lef74}-\eqref{Lef75}
to an extension of $q'$ 
over some $\cX'_W$ for $W$
a Zariski neighbourhood
of $s$ in $S$, which, plainly,  
is contrary to the maximality of $\cU'$.
\end{proof}
Since the above is slightly longer than is
strictly desirable, let's note
\begin{rmk}\label{rmk:LLef21}
If we were to make the stronger hypothesis of
homotopy depth 4 in the \'etale sense, then
for every 
locally constant sheaf $\uZ$ of finite abelian
groups, and every
closed point, $x$, in every fibre
$\cX_s$ we have $\rH^q_x(\cX_s, \uZ)=0$ for
$q<4$, so that going from \ref{prop:Lef21}
to \eqref{Lef2115} would be immediate by proper
base change.
\end{rmk}

\newpage
\bibliography{elvis}{}
\bibliographystyle{Gamsalpha}

\newpage

\appendix
\renewcommand{\thesection}{\Alph{section}}
\renewcommand{\thesubsection}{\Alph{section}.\roman{subsection}}
\section{Some 2-Category stuff for which Mr. Google doesn't
give a useful answer}
{\bf Convention.} Throughout this
appendix the words `2-category', 
respectively `strict 2-category', will
be used as a shorthand for: weak, respectively
strict, 2-category in which all
2-cells are invertible, {\it e.g.} the 2-category of groupoids
but not the 2-category of categories. 
\subsection{Fibre products}
Although logically subordinate to \S.\ref{SS:DIL}
on limits, fibre products so permeate the manuscript
as to merit a specific treatment. This said,
the basic rule for constructing a universal property
in a 2-category is to take its analogue in a 1-category,
and replace all commutative diagrams with 2-commutative
diagrams. Applying this rule of thumb to fibre products
yields
\begin{defn}\label{def:Fibre2}
Let a pair, $F_i:\cX\ra\cS$, $i=1$ or $2$, of 1-cells in
a 2-category be given, then a {\it fibre product} consists of
a 0-cell, $\cX_1\ts_\cS\cX_2$, 1-cells, $P_i:\cX_1\ts_\cS\cX_2\ra\cX_i$,
a 2-cell $\s:F_1 P_1\Rightarrow F_2 P_2$, or, equivalently
a diagram,
\begin{equation}\label{FibreA1}
 \xy 
 (0,0)*+{\cX_1\ts_\cS \cX_2}="A";
 (30,0)*+{\cX_1}="B";
 (30,-18)*+{\cS}="C";
 (0,-18)*+{\cX_2}="D";
 (22, -6)*+{}="E";
 (8,-12)*+{}="F";
    {\ar_{P_{1}} "A";"B"};
    {\ar_{P_{2}} "A";"D"};
    {\ar^{F_{1}} "B";"C"};
    {\ar_{ F_{2}} "D";"C"};
 {\ar@{=>}_{ \s} "E";"F"}; 
 \endxy 
\end{equation}
with the following universal property: given any other such square
\begin{equation}\label{FibreA2}
 \xy  
 (0,0)*+{\cY}="A";
 (30,0)*+{\cX_1}="B"; 
 (30,-18)*+{\cS}="C"; 
 (0,-18)*+{\cX_2}="D"; 
 (22, -6)*+{}="E"; 
 (8,-12)*+{}="F"; 
    {\ar_{Q_{1}} "A";"B"}; 
    {\ar_{Q_{2}} "A";"D"}; 
    {\ar^{F_{1}} "B";"C"}; 
    {\ar_{ F_{2}} "D";"C"}; 
 {\ar@{=>}_{ \tau} "E";"F"};  
 \endxy  
\end{equation}
there is a 2-commutative diagram
\begin{equation}\label{FibreA3}
 \xy 
 (0,0)*+{\cX_1\ts_\cS\cX_2}="A";
 (30,0)*+{\cX_1}="B";
 (30,-18)*+{\cS}="C";
 (0,-18)*+{\cX_2}="D";
 (22, -6)*+{}="E";
 (8,-12)*+{}="F";
 (-20,15)*+{\cY}="H";
 (-2,4)*+{}="a";
 (2,8)*+{}="b";
 (-3,-3)*+{}="c";
 (-7,-7)*+{}="d";
    {\ar_{P_{1}} "A";"B"};
    {\ar^{P_{2}} "A";"D"};
    {\ar^{F_{1}} "B";"C"};
    {\ar_{ F_{2}} "D";"C"};
 {\ar@{=>}_{ \s} "E";"F"}; 
 {\ar^{Q_1} "H";"B"};
 {\ar_{Q_2} "H";"D"};
  {\ar^{Q} "H";"A"};
 {\ar@{=>}^{ \phi_{1}} "b";"a"}; 
 {\ar@{=>}_{ \phi_{2}} "d";"c"};
(-10,14)*+{}="e";
(-18,6)*+{}="f";
(-10,15)*+{}="g";
(-22,8)*+{\tau}="h";
"e";"f" **\crv{(-35,34)} ?(.99)*\dir{>}; 
"g";"f" **\crv{(-35,35)} ?(.99)*\dir{};
\endxy 
\end{equation}
which 
for $K(-)$ the associator
is equivalent to saying that there is a
commutative diagram
\begin{equation}\label{FibreA4}
 \xy 
 (0,0)*+{F_1Q_1}="A";
 (30,0)*+{F_1 (P_1 Q) }="B";
 (30,-18)*+{F_2 (P_2 Q)}="C";
 (0,-18)*+{F_2 Q_2}="D";
 (60,0)*+{(F_1P_1)Q}="E";
 (60,-18)*+{(F_2P_2)Q}="F";
    {\ar@{=>}_{(F_1)_*\phi_1} "A";"B"};
    {\ar@{=>}_{K(F_1, P_1, Q)} "B";"E"};
    {\ar@{=>}^{Q^*\s} "E";"F"};
    {\ar@{=>}^{(F_2)_*\phi_2} "D";"C"};
{\ar@{=>}^{K(F_2, P_2, Q)} "C";"F"};
    {\ar@{=>}_{\tau } "A";"D"};
 \endxy 
\end{equation} 
with the following uniqueness property: if $(Q', \phi'_i)$
is another such diagram then there is a unique 2-cell
$\psi:Q\Rightarrow Q'$ such that 
for $i=1$ or $2$ the following diagrams
commute
\begin{equation}\label{FibreA5}
 \xy
 (-18,0)*+{Q_i}="L";
 (18,0)*+{P_i Q'}="R";
 (0,16)*+{P_i Q}="T";
    {\ar@{=>}^{\phi_i} "L";"T"};
    {\ar@{=>}^{(P_i)_*\psi} "T";"R"};
    {\ar@{=>}^{\phi'_i} "L";"R"};
 \endxy
\end{equation}
\end{defn}
It follows that a fibre product is only unique
up to unique equivalence between 1-cells which
themselves are far from unique. Now while this
is the nature of the world, 
it can be
awkward, {\it cf.} \ref{faq8}, so
we distinguish
\begin{defn}\label{def:FibreS}
Let the givens be as in \ref{def:Fibre2} then a
{\it strict fibre product}, $\cX_1\uts_\cS\cX_2$,
is a fibre product, $(\cX_1\uts_\cS\cX_2, P_i, \s)$,
with the further proviso that \eqref{FibreA3},
equivalently, 
\eqref{FibreA4}, holds with $\phi_i=\mathbf{1}$,
$i=1$ or $2$, for a unique 1-cell, $Q:\cY\ra \cX_1\uts_\cS\cX_2$.
In particular, therefore, strict fibre products
are uniquely unique, {\it i.e.} unique up to
unique strict equivalence.
\end{defn}
We ignore the question of when the existence
of fibre products and strict fibre products are
equivalent since this holds in all cases relevant
to this manuscript, {\it e.g.} \cite[2.2.2]{L-MB}
provides a rather general construction of strict
fibre products, which generalises to inverse limits,
\ref{exInv}.
\subsection{Direct and inverse limits}\label{SS:DIL}
Before applying the rule of thumb for universal
properties to limits, 
we need to define what we're taking limits of
{\it i.e.}
\begin{defn}\label{def:LimitA}
A directed, respectively inverse, system in
a 2-category $\gC$ is a 2-functor $I\ra\gC$,
respectively $I\op\ra\gC$, for some partially
ordered
set,
$I$, viewed as a 1-category in the usual way.
\end{defn}
The content of the notion of 2-functor from
a 1-category to a 2-category
may usefully
be noted for future reference, {\it i.e.}
\begin{equation}\label{LimitA1}
\begin{split}
& \text{(a). A map $i\mpo \cX_i$ from objects of $I$
to $0$-cells}\\
& \text{(b). A map, $f\mpo F_f$ from arrows to $1$-cells,
written $F_{ij}$, respectively $F^{ji}$, in the}\\
&\text{directed, respectively inverse, context of \ref{def:LimitA}}\\
& \text{(c). For every pair of compossible arrows, $g,f$,
a 2-cell, $\g_{gf}: F_{gf}\Rightarrow F_g F_f$, written}\\
&\text{$\g_{ijk}$, respectively $\g^{kji}$,
in the directed, respectively inverse, context such that}\\
&\quad\quad\quad\quad \xy 
 (0,0)*+{F_h(F_g F_f)}="A";
 (60,-9)*+{F_{hgf}}="B";
 (0,-18)*+{(F_h F_g) F_f }="D";
 (30,0)*+{F_h F_{gf}}="E";
 (30,-18)*+{F_{hg}F_f}="F";
    {\ar@{=>}_{K(F_h, F_g, F_f) } "A";"D"};
    {\ar@{=>}_{\g_{hg,f}} "B";"F"};
    {\ar@{=>}^{(F_h)_*\g_{g,f}} "E";"A"};
    {\ar@{=>}^{\g_{h,gf}} "B";"E"};
    {\ar@{=>}_{(F_{f})^* \g_{h,g}} "F";"D"};
 \endxy 
\end{split}
\end{equation}
commutes for $K(-)$ the associator of $\gC$, which,
in the directed, respectively inverse, context,
is just the tetrahedron condition of \ref{def:PosSeq}.(e)
if $\gC$ is strict. This said, we have
\begin{defn}\label{def:LimitB}
For $\cX_i$ a direct, respectively inverse, system in
$\gC$ a direct, respectively inverse, limit consists
in a 0-cell, $\varinjlim_i\cX_i$, respectively 
$\varprojlim_i\cX_i$, 1-cells $X_i:\cX_i\ra \varinjlim_i\cX_i$,
respectively $X^i:\varprojlim_i\cX_i\ra\cX_i$, 
2-cells, $\xi_{ij}:X_j\Rightarrow X_i F_{ij}$,
respectively $\xi^{ji}:X^j\Rightarrow F^{ji}X^i$,
whenever $j\ra i$, for which whenever $k\ra j\ra i$
the tetrahedron
\begin{equation}\label{LimitA2}
 \xy
 (20,-24)*+{\varinjlim_i \cX_i}="A";
 (38,6)*+{\cX_k}="B";
 (0,0)*+{\cX_j}="C";
(15,12)*+{\cX_i}="D";
{\ar_{}^{X_{k} } "B";"A"};
    {\ar^{F_{jk} } "B";"C"};
    {\ar_{X_{j} } "C";"A"};
    {\ar^{F_{ij} } "C";"D"};
    {\ar_{F_{ik} } "B";"D"};
{\ar@{-->}^{X_{i}} "D";"A"}
\endxy
\end{equation}
(with the obvious faces)
2-commutes, respectively the opposite diagram for
inverse limits such that given any other 0-cell,
$\cY$, 1-cells $Y_i$, respectively $Y^i$, and
2-cells $\eta_{ij}$, respectively $\eta^{ji}$ 
forming such a diagram there is a 1-cell $Z:\varinjlim_i\cX_i\ra\cY$,
respectively $Z:\cY\ra\varprojlim_i\cX_i$ along
with 2-cells $\z_i:Y_i\ra Z X_i$, respectively
$\z^i:Y^i\ra X^i Z$, such that for $j\ra i$ the
tetrahedron
\begin{equation}\label{LimitA3}
 \xy
 (20,-24)*+{\cY}="A";
 (38,6)*+{\cX_j}="B";
 (0,0)*+{\cX_i}="C";
(15,12)*+{\varinjlim_i \cX_i}="D";
{\ar_{}^{Y_{j} } "B";"A"};
    {\ar^{F_{ij} } "B";"C"};
    {\ar_{Y_{i} } "C";"A"};
    {\ar^{X_{i} } "C";"D"};
    {\ar_{X_{j} } "B";"D"};
{\ar@{-->}^{Z} "D";"A"}
\endxy
\end{equation}
(with the obvious faces)
2-commutes, respectively the opposite diagram for
inverse limits while enjoying the following
uniqueness property: if $Z'$, $\z'_i$ {\it etc.}
also yields a series of such 2-commutative diagrams
then there is a unique $\psi:Z\Rightarrow Z'$ such that
\begin{equation}\label{LimitA4}
 \xy
 (-18,0)*+{Y_i}="L";
 (18,0)*+{Z'X_i}="R";
 (0,16)*+{Z X_i}="T";
    {\ar@{=>}^{\z_i} "L";"T"};
    {\ar@{=>}^{(X_i)^*\psi} "T";"R"};
    {\ar@{=>}^{\z'_i} "L";"R"};
 \endxy
\end{equation}
respectively (modulo notation) \eqref{FibreA5},
commutes for all $i$.
\end{defn}
Again, if limits exists they're unique up to
equivalence, while pertinent existence results
are
\begin{ex}\label{exDirect} {\it Direct limits of groupoids}
In this case we can define a groupoid, $\cX$, whose 
objects are
\begin{equation}\label{LimitA5}
\mathrm{ob}(\cX):= \coprod_i \mathrm{ob}(\cX_i)
\end{equation}
Now suppose that we have arrows $k\ra i\ra h$, $j\ra i\ra h$,
and objects $x_k\in \mathrm{ob}(\cX_i)$, 
$x_j\in \mathrm{ob}(\cX_j)$ then we have Hom-sets,
\begin{equation}\label{LimitA6}
\Hom^i(x_j, x_k) := \Hom_{\cX_i} (F_{ij}(x_j), F_{ik}(x_k))
\end{equation}
with transition functions
\begin{equation}\label{LimitA7}
T_{hi}:\Hom^i(x_j, x_k)\ra \Hom^h(x_j, x_k):
f\mpo \g_{hik}^{-1}(x_k) F_{hi}(f) \g_{hij}(x_j)
\end{equation} 
satisfying $T_{gi}=T_{gh}T_{hi}$ for any further arrow $h\ra g$
by \eqref{LimitA1}, and
\begin{equation}\label{LimitA8}
\Hom_\cX (x_j, x_k) := \varinjlim_i \Hom^i (x_j, x_k)
\end{equation}
where the latter is just directed limit in $\Ens$ 
with respect to the transition functions $T_{hi}$
does
the job.
\end{ex}
\begin{ex}\label{exInv}{\it Inverse limits of groupoids}
Here we define a groupoid, $\cX$, with objects a subset of
\begin{equation}\label{LimitA9}
\prod_{i\geq j} \mathrm{ob}(\cX^i)\ts \mathrm{Ar}(\cX^j) 
\ni (x^i, \xi^{ji})
\end{equation}
where we further insist that
\begin{equation}\label{LimitA10}
\xi^{ji}:F^{ji} (x^i) \ra x^j,\,\,\text{and}\,\, 
\xi^{ki}=\xi^{kj}F^{kj}(\xi^{ji})\g^{kji},\,\,
\text{whenever}\,\, i\geq j\geq k
\end{equation} 
while an arrow between two such objects, 
$x^i\ts \xi^{ji}$, $y^i\ts \eta^{ji}$ is an element
$(f^i)\in \prod_i \mathrm{Ar}(\cX_i)$ with
source $x^i$, sink $y^i$ for which the following
diagram commutes 
\begin{equation}\label{LimitA11}
\begin{CD}
F^{ji}(x^i)@>>{F^{ji}(f^i)}> y^i\\
@V{\xi^{ji}}VV @VV{\eta^{ji}}V\\
x^j @>f^j>> y^j
\end{CD}
\end{equation} 
\end{ex}
In the particular case of groups, {\it i.e.} groupoids
with 1-object, one can usefully observe the contrasting
behaviour of the direct and inverse limits, to wit:
\begin{rmk}\label{rmk:Limit}
If $\cX_i$ are groups, $\G_i\rras \rp$, 
with $I$ filtered on the left and the directed,
respectively inverse, system is in groups, {\it i.e.}
$\g_{ijk}$, respectively $\g^{kji}$, is the identity, then:

(a) The directed limit is again a group, and coincides
with the group direct limit.

(b) The inverse limit need not be a group. In fact the
isomorphism classes of objects in $\varprojlim_i \cX_i$-
which is an invariant of groupoids up to equivalence- is
${\varprojlim_i}^{(1)} \G_i$, while the stabaliser of
the object defined by the points of the $\cX_i$ is
the inverse limit in groups. Consequently,
\cite[Th\'eor\`eme 7.1]{jensen}, inverse limits of
finite groups are pro-finite groups, but otherwise
anything can happen.
\end{rmk} 
\subsection{Champ means 2-sheaf}\label{SS:TwoSheaf}
For $E$ a site, the fact that an $E$-groupoid in the
sense of \cite[2.1]{L-MB}- replace $\mathrm{Aff}/S$ in
{\it op. cit.} by $E$- or an $E$ fibred groupoid in
the sense of \cite[I.1.0.2]{giraud} is equivalent to
giving a normalised ({\it i.e.} trivial
on identities) 2-functor
\begin{equation}\label{2Sh1}
\cX: E\op \ra \Grpd: U\mpo \cX(U);\,\, 
\{V\xrightarrow{f} U\} \mpo X_f;\,\,
(g,f)\mpo \xi_{g,f}
\end{equation}
as defined in \eqref{LimitA1}
is \cite[Expos\'e VI.7-8]{sga1}, where we do not adopt the
shorthand $X_f= f^*$ to avoid confusion with pull-backs
of natural transformations, or the actual 1-functor
of \eqref{2Sh3}. The only reason, therefore, why one
doesn't call this a pre-champ (or functorially with respect
to the ideas pre-2-sheaf in English) is because the
gluing/sheaf condition is itself double barrelled, 
\eqref{2Sh4},
and
one wants 
(perhaps not necessarily correctly)
to reserve `pre-champ' for the 2-functors of
\eqref{2Sh1} which satisfies some but not all of the
gluing condition. More precisely for $x,y\in\cX(U)$ one
has by \eqref{LimitA1} an actual (and not just a
pseudo) functor
\begin{equation}\label{2Sh2}
\underline{\Hom}_U (x,y):E/U\op\ra\Ens: \{V\xrightarrow{v} U\}
\mpo \Hom_{\cX(V)} ( X_v(x), X_v(y))
\end{equation}
where for $W\xrightarrow{f} V$ the transition maps are given 
by, 
\begin{equation}\label{2Sh3} 
\underline{\Hom}_U(x,y)(V) \xrightarrow{f^*}
\underline{\Hom}_U(x,y)(W): \rho \mpo \xi_{f,v}^{-1} X_f (\rho) \xi_{f,v}
\end{equation}
and one asks that
for every cover $V\ra U\in E$,
\begin{equation}\label{2Sh4}
1\ra\cX(U)\ra\cX(V){\build\rras_{X_s}^{X_t}}\cX(V\ts_U V) 
{\build\rightarrow_{\rightarrow}^{\rightarrow}}
\cX(V\ts_U V\ts_U V)
\end{equation}
is ``exact'', 
according to the only possible meaning that
exact could have,
{\it i.e.}, 
\cite[3.1]{L-MB} or \cite[II.1.2.1]{giraud},
for 
every pair of objects \eqref{2Sh2}-\eqref{2Sh3}
defines a sheaf (a.k.a. the pre-champ condition),  while for
$p_{ij}$, respectively
$p_i$, the projections from $V\ts_U V\ts_U V$ to
$V\ts_U V$, respectively $V$, the objects of $\cX(U)$
are, up to
isomorphism, precisely objects of $x\in \cX(V)$ together 
with an arrow $\phi: X_s(x)\ra X_t(x)$ such that
in the sense of \eqref{2Sh3}- albeit relative to
$V$ rather than $U$- $p_{13}^*\phi=p_{23}^*\phi p_{12}^*\phi$.
Evidently, 
therefore, one thinks of $x\in\cX(U)$ as a global
section, and
we may usefully observe,

a) In principle there is a competing definition of
global section over $U$, {\it i.e.} a transformation
between the identity 2-functor on $E/U$ and $\cX$
restricted to the same. However, essentially for
the same reason that \eqref{2Sh2}-\eqref{2Sh3}
define a sheaf, the objects $V\mpo X_v (x)$,
and arrows $f\mpo \xi_{f,v}$, in the notation of
{\it op. cit.} define such a transformation, and,
better still, any other transformation admits a
modification into one enjoying this special form,
{\it i.e. $\cX(U)$ is all possible global sections
up to isomorphism}.

b) There is a minor, but informative, generalisation,
to wit: let $U\in E$, $U_i\in E/U$, $i=1$ or $2$ be given, 
along with objects $x_i\in \cX(U_i)$ then,
\begin{equation}\label{2Sh5} 
E/U\ni V \mpo (v_i:V\ra U_i, \, \phi: X_{v_1}(x_1) \ra X_{v_2}(x_2))
\end{equation}
is
for $j:U_1\ts_U U_2\ra U$ the fibre product,
with projections $p_i$,
representable by the sheaf,
\begin{equation}\label{2Sh6}
j_{!}\,\bigl(
\underline{\Hom}_{U_1\ts_U U_2} (X_{p_1}(x_1), X_{p_2}(x_2))
\bigr)
\end{equation}
Consequently the espace \'etal\'e generalises 
as it should, {\it i.e.}
\begin{cor}\label{cor:2Sh} Let $E$ be the \'etale site of
a 
topological, respectively algebraic, space $S$, and $\cX/E$ a
champ 
(as defined in \eqref{2Sh1} and \eqref{2Sh4}) 
then there is a cover $U\ra S$ which viewed as a sheaf
factors through $\cX$, and, \ref{def:Fibre2}, 
$(s,t):R:=U\ts_\cX U\rras U$,
$U$, but 
not necessarily $R$, separated,
is an \'etale groupoid such that (over $E$) $\cX$ is
equivalently to the classifying champ $[U/R]$.
\end{cor}
\begin{proof} Following the case of the analogous 
proposition for sheaves- \cite[V.1.5]{milne}- 
we just put $U$ to be the co-product over all pairs
$(x,V)$ consisting in a
sufficiently small open (so by the respective definitions
this will be supposed even if it isn't standard usage
in the topological case to imply separated) $V\ra S$
and a section $x\in \cX(V)$. By item.(a) above, the
$x$, define a 2-functor from $U$ (identified with the
identity on $(E/U)\op$) to $\cX$ while by
item.(b) the resulting fibre product $U\ts_\cX U$ is
represented by a sheaf, of which $R$ is 
the espace \'etal\'e.
\end{proof}
Plainly, \ref{cor:2Sh}, invites the reading ``every 2-sheaf
is representable'', which, although correct, 
like the espace \'etal\'e of a sheaf
requires
a little care if things aren't locally constant.
For example, whether algebraically or topologically,
consider the classifying champ $\cY$ of the action,
$\ba^1_\bc\ts\mu_2 \rras \ba^1_\bc$, by $\pm 1$ on
the affine line. The coarse quotient, 
$Y=\ba^1_\bc/\pm\xrightarrow{\sim}\ba^1_\bc$ is a
moduli space for $\cY$, but the ``2-sheaf of sections
of $\cY\ra Y$'' is just the sheaf $j_! Y^*$ for
$j:Y^*\hookrightarrow Y$ the complement of the origin,
and one fails to recover $\cY$.
Similarly, the space of arrows $R$-
even in cases such as 
$S$ a  manifold where the topology of $S$
poses no obstruction-  needn't
be separated, and while $R$ (as opposed to $s\ts t$)
separated isn't, {\it cf.} \cite[4.1]{L-MB}, always
required in the definition of an algebraic champ,
quasi-compactness of $s\ts t$ usually is,
which in turn, needs further,
if wholly reasonable hypothesis, {\it e.g.} the sheaves defined by
\eqref{2Sh3}-\eqref{2Sh4} are constructible.

\end{document}